\documentclass[a4paper,12pt,amsart,frenchb]{article}

\usepackage{amsmath,amsbsy,amsfonts,amssymb}
\usepackage[french]{babel}
\newcommand{\ass}[2]{\vskip0.3cm\noindent
{\bf {#1}}. { \sl {#2}}\vskip0.3cm\noindent
}
  
\begin{document}

    \title{ Stabilisation de la formule des traces tordue I: endoscopie tordue sur un corps local}
\author{J.-L. Waldspurger}
\date{18 janvier 2014}
\maketitle

{\bf Introduction}

  Ceci est le premier d'une s\'erie d'articles, en collaboration avec C. Moeglin, visant \`a  stabiliser la formule des traces tordue. L'essentiel du travail consiste \`a reprendre dans ce cadre tordu la d\'emonstration (colossale) qu'Arthur a mise au point dans le cas non tordu. Mais auparavant, un certain nombre de travaux pr\'eparatoires sont n\'ecessaires. Le texte qui suit est l'un d'eux. On y pr\'esente les d\'efinitions et propri\'et\'es de base de la th\'eorie de l'endoscopie tordue sur un corps local de caract\'eristique nulle, du c\^ot\'e "g\'eom\'etrique", c'est-\`a-dire du c\^ot\'e des int\'egrales orbitales. Cela fournira, on l'esp\`ere, un socle pour la suite de nos travaux. 

Ainsi, ce texte ne contient gu\`ere de r\'esultats originaux. Il reprend largement les travaux fondamentaux de Kottwitz-Shelstad, Labesse et Shelstad sur la question.  On a toutefois modifi\'e sur certains points la pr\'esentation de ces auteurs. Donnons un peu plus de d\'etails. La premi\`ere section donne les d\'efinitions de base des espaces tordus et de leurs donn\'ees endoscopiques. Les espaces tordus ont \'et\'e introduits par Labesse et remplacent, avantageusement nous semble-t-il, les couples form\'es d'un groupe connexe et d'un automorphisme de celui-ci. Notons que, dans le cadre le plus g\'en\'eral, on doit aussi associer aux donn\'ees endoscopiques des espaces tordus. On en donne en 1.7 une d\'efinition parfaitement canonique, ce qui est l'un des points nouveaux de notre pr\'esentation. Un autre point nouveau est que l'on a fait dispara\^{\i}tre le traditionnel groupe quasi-d\'eploy\'e $G^*$. A notre avis, ce groupe est mal adapt\'e \`a l'endoscopie tordue, parce qu'il n'y a pas d'espace tordu $\tilde{G}^*$. Plus exactement, on peut d\'efinir un tel espace tordu, mais il n'y a pas de correspondance canonique entre les classes de conjugaison stable dans l'espace de d\'epart $\tilde{G}$ et  les classes de conjugaison stable dans cet espace $\tilde{G}^*$. Pour \'etudier la correspondance entre classes de conjugaison stable dans $\tilde{G}$ et dans un espace endoscopique $\tilde{G}'$, correspondance qui est parfaitement canonique et \'equivariante pour les actions galoisiennes, ce n'est pas un bon point de d\'epart de la d\'ecomposer en deux correspondances entre $\tilde{G}$ et $\tilde{G}^*$ d'une part, entre $\tilde{G}^*$ et $\tilde{G}'$ d'autre part, qui ne sont ni canoniques, ni \'equivariantes pour les actions galoisiennes. En fait, le groupe $G^*$ sert rarement. Ce qui sert, c'est son tore maximal $T^*$. Mais ce tore se r\'ecup\`ere facilement en utilisant la m\'ethode qu'on a apprise de Deligne: c'est $\underline{le}$ tore maximal de $G$ muni de son action galoisienne canonique, cf. 1.2. Dans la section 2, on r\'ecrit la d\'efinition des facteurs de transfert d'apr\`es Kottwitz et Shelstad, puis celle du transfert des int\'egrales orbitales. Une donn\'ee endoscopique ${\bf G}'=(G',{\cal G}',\tilde{s})$ \'etant fix\'ee, pour d\'efinir ce transfert d'int\'egrales orbitales, on doit fixer des donn\'ees auxiliaires, en particulier un groupe $G'_{1}$ au-dessus de $G'$, et un facteur de transfert pour ces donn\'ees. Malheureusement,  la stabilisation de la formule des traces tordue n\'ecessite de pouvoir changer de donn\'ees auxiliaires. La raison en est que si $\tilde{M}$ est un espace de Levi de $\tilde{G}$ et ${\bf M}'$ est une donn\'ee endoscopique de $\tilde{M}$, ${\bf M}'$ peut appara\^{\i}tre comme "donn\'ee de Levi" de plusieurs donn\'ees endoscopiques de $\tilde{G}$ et qu'on ne peut pas assurer que les restrictions \`a ${\bf M}'$ des donn\'ees auxiliaires affect\'ees \`a ces diverses donn\'ees co\"{\i}ncident. Il convient donc de savoir ce qui se passe quand on change de donn\'ees auxiliaires. Il s'av\`ere que les objets construits \`a l'aide de deux s\'eries de donn\'ees auxiliaires sont canoniquement isomorphes. Mais alors, il est aussi simple d'\'eliminer formellement les donn\'ees auxiliaires en rempla\c{c}ant ces objets par leur limite inductive (par ces isomorphismes canoniques) sur toutes les donn\'ees auxiliaires possibles. C'est ce que l'on fait en 2.5. Cette pr\'esentation permet ensuite de d\'efinir naturellement  sur ces objets une action du groupe d'automorphismes de la donn\'ee endoscopique ${\bf G}'$, cf. 2.6. Cette action est assez subtile car, dans la situation tordue, ce groupe d'automorphismes contient un sous-groupe qui agit trivialement sur le groupe $G'$. Mais il agit sur l'espace des fonctions sur ce groupe par multiplication par des caract\`eres. La section 3 compare les donn\'ees endoscopiques d'espaces de Levi avec les Levi de donn\'ees endoscopiques.  La section 4 d\'ecrit exactement l'image du transfert des int\'egrales orbitales. La section 5 introduit ce que l'on appelle les distributions $\omega$-\'equivariantes "g\'eom\'etriques", qui sont  celles dont le support est r\'eduit \`a une r\'eunion finie de classes de conjugaison. On a dualement un transfert entre de telles distributions et on d\'etermine son noyau. On examine aussi le comportement de ces distributions par descente d'Harish-Chandra. 
Signalons que, dans le cas d'un corps archim\'edien, les r\'esultats des sections 4 et 5 reposent essentiellement sur ceux de Renard et Shelstad. Enfin, on traite dans la section 6 le cas "non ramifi\'e", o\`u l'on peut d\'efinir des facteurs de transfert canoniques, modulo le choix "d'espaces hypersp\'eciaux".

\section{Les d\'efinitions de base}

\subsection{Groupes et espaces tordus}
Soit $F$ un corps de caract\'eristique nulle, dont on fixe une cl\^oture alg\'ebrique $\bar{F}$. Posons $\Gamma_{F}=Gal(\bar{F}/F)$. Soit $G$ un groupe alg\'ebrique d\'efini sur $F$, r\'eductif et connexe.  
On l'identifie au groupe de ses points sur $\bar{F}$. Le groupe $\Gamma_{F}$ agit sur $G$. Pour $\sigma\in \Gamma_{F}$, on note encore $\sigma$ son action sur $G$, ou $\sigma_{G}$
s'il semble bon de pr\'eciser.   Pour $g\in G$, on note $ad_{g}$ l'automorphisme int\'erieur $x\mapsto gxg^{-1}$ de $G$. On note $Z(G)$ le centre de $G$ et $A_{G}$ le plus grand sous-tore de $Z(G)$ qui soit d\'eploy\'e sur $F$ (remarquons que $A_{G}$ d\'epend du corps $F$). On pose ${\cal A}_{G}=X_{*}(A_{G})\otimes_{{\mathbb Z}}{\mathbb R}$, avec la notation $X_{*}$ usuelle.  On note $G_{AD}$ le groupe adjoint de $G$ et $G_{SC}$ le rev\^etement simplement connexe du groupe d\'eriv\'e de $G$. On notera souvent de la m\^eme fa\c{c}on un \'el\'ement, ou un sous-ensemble, de $G_{SC}$ et son image dans $G$. N\'eanmoins, si besoin est, on notera $\pi:G_{SC}\to G$ l'homomorphisme naturel. Si $X$ est un sous-ensemble de $G$, on note $X_{ad}$ son image dans $G_{AD}$ et $X_{sc}$ l'image r\'eciproque de $X_{ad}$ dans $G_{SC}$ (ce qui n'est pas forc\'ement l'image r\'eciproque de $X$). On aura tendance \`a noter de la m\^eme fa\c{c}on deux objets qui se d\'eduisent l'un de l'autre par fonctorialit\'e. Par exemple, pour $g\in G$, on note encore $ad_{g}$ les automorphismes de $G_{AD}$ ou de $G_{SC}$ qui se d\'eduisent de l'automorphisme $ad_{g}$ de $G$. 

Soit $\tilde{G}$ un espace tordu sous $G$. C'est une vari\'et\'e alg\'ebrique sur $F$. Le groupe $G$ agit \`a droite et \`a gauche sur $\tilde{G}$ et, pour chaque action, $\tilde{G}$ est un espace principal homog\`ene sous $G$.  Il y a une application $\gamma\mapsto ad_{\gamma}$ de $\tilde{G}$ dans le groupe des automorphismes de $G$ telle que $\gamma g=ad_{\gamma}(g)\gamma$ pour tout $g\in G$. On a l'\'egalit\'e $ad_{g\gamma g'}=ad_{g}\circ ad_{\gamma}\circ ad_{g'}$ pour tous $g,g'\in G$ et $\gamma\in \tilde{G}$. Les actions et applications ci-dessus sont toutes alg\'ebriques et d\'efinies sur $F$. Pour $\gamma\in \tilde{G}$, on note $Z_{G}(\gamma)$ son commutant dans $G$ (c'est-\`a-dire l'ensemble des points fixes de $ad_{\gamma}$). On note $G_{\gamma}=Z_{G}(\gamma)^0$ la composante neutre de ce groupe. L'image de $ad_{\gamma}$ dans le groupe des automorphismes ext\'erieurs  de $G$ ne d\'epend pas de $\gamma$. D'autre part, l'automorphisme $ad_{\gamma}$ d\'efinit par fonctorialit\'e des automorphismes de divers objets. Quand ils ne d\'ependent pas de $\gamma$ (ou m\^eme de $\gamma$ dans un sous-ensemble indiqu\'e), on note ces automorphismes $\theta$. Ainsi, il y a un automorphisme $\theta$ du centre $Z(G)$. On note $A_{\tilde{G}}$ le plus grand tore d\'eploy\'e sur $F$ contenu dans $Z(G)^{\theta}$.  On pose ${\cal A}_{\tilde{G}}=X_{*}(A_{\tilde{G}})\otimes{\mathbb R}$.

On dira que $\tilde{G}$ est \`a torsion int\'erieure si, pour $\gamma\in \tilde{G}$, l'automorphisme $ad_{\gamma}$ de $G$ est int\'erieur. En fixant $\gamma$ et en le multipliant  par un \'el\'ement convenable de $G$, on obtient un \'el\'ement tel que $ad_{\gamma}$ soit l'identit\'e. Alors l'application $g\gamma \mapsto g$ identifie $\tilde{G}$ \`a $G$ muni de ses actions par multiplication \`a droite et \`a gauche. Mais cet isomorphisme n'est en g\'en\'eral d\'efini que sur $\bar{F}$, car on ne peut pas toujours trouver de $\gamma$ comme ci-dessus qui appartienne \`a $ \tilde{G}(F)$. 

{\bf Exemple.} On fixe un entier $n\geq1$ et un \'el\'ement $d\in F^{\times}$. On prend $G=SL(n)$ et $\tilde{G}=\{g\in GL(n); det(g)=d\}$. Cet espace tordu est trivial sur $F$ si et seulement si $d$ appartient au groupe $F^{\times,n}$ des puissances $n$-i\`emes dans $F^{\times}$.

\bigskip

\subsection{Paires de Borel}

On appelle paire de Borel de $G$ un couple $(B,T)$ form\'e d'un sous-groupe de Borel $B$ et d'un sous-tore maximal $T$ de $B$. On ne suppose pas que $B$ ou $T$ soient d\'efinis sur $F$. On appelle paire de Borel \'epingl\'ee un triplet ${\cal E}=(B,T,(E_{\alpha})_{\alpha\in \Delta})$ o\`u $(B,T)$ est une paire de Borel et $(E_{\alpha})_{\alpha\in \Delta}$ est un \'epinglage relatif \`a cette paire. C'est-\`a-dire que $\Delta$ est l'ensemble des racines simples de $T$ agissant dans l'alg\`ebre de Lie $\mathfrak{u}$ du radical unipotent de $B$ et, pour tout $\alpha\in \Delta$, $E_{\alpha}$ est un \'el\'ement non nul de la droite radicielle $\mathfrak{u}_{\alpha}\subset \mathfrak{u}$ associ\'ee \`a $\alpha$. Pour deux paires de Borel \'epingl\'ees ${\cal E}=(B,T,(E_{\alpha})_{\alpha\in \Delta})$ et ${\cal E}'=(B',T',(E'_{\alpha'})_{\alpha'\in \Delta'})$, il existe $g\in G_{SC}$ tel que $ad_{g}$ transporte ${\cal E}$ sur ${\cal E}'$. Cet \'el\'ement $g$ n'est pas unique mais sa classe $gZ(G_{SC})$ l'est. Les restrictions de $ad_{g}$ \`a $B$ et $T$ sont uniquement d\'etermin\'ees. Cela autorise \`a d\'efinir $\underline{la}$ paire de Borel \'epingl\'ee ${\cal E}^*=(B^*,T^*,(E^*_{\alpha})_{\alpha\in \Delta})$ comme la limite inductive de toutes les paires de Borel \'epingl\'ees, les applications de transition \'etant celles ci-dessus. Par un m\^eme proc\'ed\'e de limite inductive, on d\'efinit l'ensemble $\Sigma$ des racines de $T^*$ dans l'alg\`ebre de Lie de $G$, l'ensemble $\check{\Sigma}$ des coracines et le groupe de Weyl $W$. Pour une paire de Borel \'epingl\'ee ${\cal E}$, ces ensembles s'identifient \'evidemment aux m\^emes ensembles relatifs \`a cette paire.

 Le groupe $\Gamma_{F}$ agit naturellement sur l'ensemble des paires de Borel ou des paires de Borel \'epingl\'ees. On en d\'eduit une action de $\Gamma_{F}$ sur ${\cal E}^*$, not\'ee $\sigma\mapsto \sigma_{G^*}$. Pour n'importe quelle paire de Borel \'epingl\'ee ${\cal E}$, $\sigma_{G^*}$ est la compos\'ee des isomorphismes
 $${\cal E}^*\simeq {\cal E}\stackrel{\sigma_{G}}{\to}\sigma_{G}({\cal E})\simeq {\cal E}^*.$$
 On en d\'eduit une action de $\Gamma_{F}$ sur $\Delta$, $\Sigma$, $\check{\Sigma}$ et $W$. 
 
 Pour une paire de Borel \'epingl\'ee ${\cal E}$ et pour $\sigma\in \Gamma_{F}$, choisissons $u_{{\cal E}}(\sigma)\in G_{SC}$ tel que $ad_{u_{{\cal E}}(\sigma)}\circ \sigma_{G}({\cal E})={\cal E}$. Alors l'isomorphisme de ${\cal E}$ sur ${\cal E}^*$ transporte l'action $\sigma\mapsto ad_{u_{{\cal E}}(\sigma)}\circ \sigma_{G}$ sur $\sigma\mapsto \sigma_{G^*}$. L'application $\sigma\mapsto( u_{{\cal E}}(\sigma))_{ad}$ est un cocycle \`a valeurs dans $G_{AD}$ dont la classe ne d\'epend pas de la paire ${\cal E}$. On dit qu'une paire de Borel ou une paire de Borel \'epingl\'ee est d\'efinie sur $F$ si et seulement si elle est fixe par  l'action naturelle $\sigma\mapsto \sigma_{G}$. Dans le cas d'une paire de Borel \'epingl\'ee ${\cal E}$, cela revient \`a dire que l'on peut choisir $u_{{\cal E}}(\sigma)=1$ pour tout $\sigma$ (mais, bien s\^ur, $\sigma_{G}$ peut agir sur $\Delta$ par une permutation non triviale). Dans ce cas,  on peut identifier ${\cal E}^*$ \`a ${\cal E}$ et l'action $\sigma\mapsto \sigma_{G^*}$ \`a l'action naturelle $\sigma\mapsto \sigma_{G}$. On dit que $G$ est quasi-d\'eploy\'e si et seulement s'il existe une paire de Borel \'epingl\'ee d\'efinie sur $F$ (il suffit d'ailleurs qu'il existe une paire de Borel tout court d\'efinie sur $F$).
 
 Pour toute paire de Borel \'epingl\'ee ${\cal E}$, notons $Z(\tilde{G},{\cal E})$ l'ensemble des $e\in \tilde{G}$ tels que $ad_{e}$ conserve ${\cal E}$. C'est un espace principal homog\`ene sous $Z(G)$, \`a droite comme \`a gauche. Notons ${\cal Z}(\tilde{G},{\cal E})$ le quotient de $Z(\tilde{G},{\cal E})$ par l'action par conjugaison de $Z(G)$. Alors ${\cal Z}(\tilde{G},{\cal E})$ est un espace principal homog\`ene, \`a droite comme \`a gauche, sous ${\cal Z}(G):=Z(G)/(1-\theta)(Z(G))$ (on note $1-\theta$ l'homomorphisme $z\mapsto z\theta(z)^{-1}$). Si ${\cal E}'$ est une autre paire de Borel \'epingl\'ee, on choisit comme ci-dessus $g\in G_{SC}$ tel que $ad_{g}({\cal E})={\cal E}'$. Alors $ad_{g}:Z(\tilde{G},{\cal E})\to Z(\tilde{G},{\cal E}')$ est un isomorphisme. Il n'est pas uniquement d\'efini car $g$ n'est  pas unique. Mais, par passage aux quotients, $ad_{g}$ d\'efinit un isomorphisme de ${\cal Z}(\tilde{G},{\cal E})$ sur ${\cal Z}(\tilde{G},{\cal E}')$ qui est uniquement d\'efini. On note ${\cal Z}(\tilde{G})$ la limite inductive des ${\cal Z}(\tilde{G},{\cal E})$ sur les paires de Borel \'epingl\'ees, les applications de transition \'etant les isomorphismes canoniques que l'on vient de d\'efinir. Alors ${\cal Z}(\tilde{G})$ est un espace tordu sous le groupe ${\cal Z}(G)$.  On d\'efinit une action $\sigma\mapsto \sigma_{G^*}$ de $\Gamma_F$ sur ${\cal Z}(\tilde{G})$  comme on a d\'efini l'action sur ${\cal E}^*$. On voit que ${\cal Z}(\tilde{G})$ est un espace tordu sous ${\cal Z}(G)$, d\'efini sur $F$. Remarquons que ${\cal Z}(\tilde{G})(F)$ peut \^etre vide.

Soit ${\cal E}= (B,T,(E_{\alpha})_{\alpha\in \Delta})$ une paire de Borel \'epingl\'ee. Pour $e\in Z(\tilde{G},{\cal E})$, l'automorphisme $ad_{e}$ de $G$ ne d\'epend pas du choix de $e$. On le note $\theta_{{\cal E}}$ ou simplement $\theta$. Remarquons que, si $\gamma\in \tilde{G}$ est tel que $ad_{\gamma}$ conserve seulement $(B,T)$, la restriction de $ad_{\gamma}$ \`a $T$ co\"{\i}ncide avec  celle de $\theta$. Par restriction puis passage \`a la limite, on obtient un automorphisme de ${\cal E}^*$ que l'on note $\theta^*$. Il commute \`a l'action galoisienne sur ${\cal E}^*$.  Rappelons deux propri\'et\'es cruciales du sous-groupe $W^{\theta^*}$ (avec la notation usuelle: c'est le sous-groupe des points fixes de $\theta^*$ agissant dans $W$):

(1) un \'el\'ement $\omega\in W$ appartient \`a $W^{\theta^*}$ si et seulement s'il conserve $(T^*)^{\theta^*}$ ou $(T^*)^{\theta^*,0}$;

(2) pour ${\cal E}$ et $e\in Z(\tilde{G},{\cal E})$ comme ci-dessus, $W^{\theta^*}$ s'identifie au groupe de Weyl de $G_{e}$ relatif \`a son sous-tore maximal $T^{\theta,0}$.

\bigskip

\subsection{El\'ements semi-simples}

 Un \'el\'ement $\gamma\in \tilde{G}$ est dit semi-simple si et seulement s'il existe une paire de Borel de $G$ qui est conserv\'ee par $ad_{\gamma}$ (la terminologie plus correcte est "quasi-semi-simple"; en vertu de l'hypoth\`ese "$\theta^*$ est d'ordre fini" que l'on imposera d\`es 1.5, on peut aussi bien abandonner le "quasi"). Supposons $\gamma$ semi-simple. On dit qu'il est fortement r\'egulier si et seulement si $Z_{G}(\gamma)$ est ab\'elien et  la composante neutre $G_{\gamma}$ est un tore. On note $\tilde{G}_{ss}$ l'ensemble des \'el\'ements semi-simples et $\tilde{G}_{reg}$ l'ensemble des \'el\'ements semi-simples et fortement r\'eguliers.

 Soient ${\cal E}= (B,T,(E_{\alpha})_{\alpha\in \Delta})$ une paire de Borel \'epingl\'ee et $\gamma\in \tilde{G}$ tel que $ad_{\gamma}$ conserve $(B,T)$. On pose $\theta=\theta_{{\cal E}}$. On a

(1) pour tout $e\in Z(\tilde{G},{\cal E})$, il existe $t\in T$ tel que $\gamma=te$;

(2) une paire de Borel $(B',T')$ de $G$ est conserv\'ee par $ad_{\gamma}$ si et seulement s'il existe $\omega\in W^{\theta}$  et $x\in G_{\gamma}$ tels que $(B',T')=ad_{x}\circ\omega(B,T)$.

Preuve. Il existe $t\in G$ tel que $\gamma=te$. Puisque $ad_{\gamma}$ et $ad_{e}$ conservent $(B,T)$, $ad_{t}$ aussi donc $t$ appartient \`a $T$. Pour $\omega\in W^{\theta}$, on rel\`eve  $\omega$ gr\^ace \`a 1.2(2)  en un \'el\'ement $n\in G_{e}$ qui normalise $T^{\theta,0}$, donc aussi son commutant $T$. La paire $\omega(B,T)=ad_{n}(B,T)$ est conserv\'ee par $ad_{e}$. Elle l'est aussi par $t\in T=ad_{n}(T)$, donc elle est conserv\'ee par $ad_{\gamma}$. Pour $x\in G_{\gamma}$, la paire $ad_{x}\circ\omega(B,T)$ l'est aussi. Inversement, soit $(B',T')$ une paire conserv\'ee par $ad_{\gamma}$. D'apr\`es [KS1] th\'eor\`eme 1.1.A, le couple $(B'\cap G_{\gamma},T'\cap G_{\gamma})$ est une paire de Borel de $G_{\gamma}$. Il existe donc $x\in G_{\gamma}$ tel que l'image de cette paire par $ad_{x}$ ait pour tore maximal $T^{\theta,0}$. Quitte \`a remplacer $(B',T')$ par $ad_{x}(B',T')$, on peut supposer $T'=T$. Cette paire est alors conserv\'ee par $ad_{t}$, donc aussi par $ad_{e}$. Par le m\^eme argument, le couple $(B'\cap G_{e},T^{\theta,0})$ est une paire de Borel de $G_{e}$. Gr\^ace \`a 1.2(2), il existe  $\omega\in W^{\theta}$ tel que $(B'\cap G_{e},T^{\theta,0})$ se d\'eduise de $(B\cap G_{e},T^{\theta,0})$ par l'action de $\omega$. Autrement dit, $(B',T)$ et $(\omega(B),T)$ ont m\^eme intersection avec $G_{e}$. Or, parce que $ad_{e}$ conserve un \'epinglage, cette op\'eration d'intersection avec $G_{e}$ est une bijection entre les paires de Borel de $G$ conserv\'ees par $ad_{e}$ et les paires de Borel de $G_{e}$, cf. [KS1] p.14.  Donc  $(B',T)=(\omega(B),T)$. $\square$

Notons   $p:T\to T/(1-\theta)(T)$ l'homomorphisme naturel. Le groupe $W^{\theta}$ agit sur sur le quotient $T/(1-\theta)(T)$.  Supposons $T$ d\'efini sur $F$ et  $\gamma\in \tilde{G}(F)$. Alors $\theta$ est d\'efini sur $F$. Ecrivons $\gamma=te$ comme en (1). Pour tout $\sigma\in \Gamma_{F}$, on introduit un \'el\'ement $u_{{\cal E}}(\sigma)\in G_{SC}$ comme en 1.2. On a

(3) $u_{{\cal E}}(\sigma)$ normalise $T$ et son image dans $W$ appartient \`a $W^{\theta}$;

(4) il existe $z(\sigma)\in Z(G)$ tel que $u_{{\cal E}}(\sigma)\sigma(e)u_{{\cal E}}(\sigma)^{-1} =z(\sigma)^{-1}e$ et $\sigma\circ p(t)=p(z(\sigma)t)$.

Preuve. La paire $(\sigma(B),T)$ est conserv\'ee par $ad_{\gamma}$, donc aussi par $ad_{e}$. Cela entra\^{\i}ne comme ci-dessus qu'elle se d\'eduit de $(B,T)$ par l'action d'un \'el\'ement de $W^{\theta}$. Or $(B,T)=ad_{u_{{\cal E}}(\sigma)}(\sigma(B),T)$, d'o\`u (3).   On peut \'ecrire $u_{{\cal E}}(\sigma)=n(\sigma)t(\sigma)$ o\`u $t(\sigma)\in T$ et $n(\sigma)\in G_{e}$. L'\'el\'ement $u_{{\cal E}}(\sigma)\sigma(e)u_{{\cal E}}(\sigma)^{-1}$ appartient encore \`a $Z(\tilde{G},{\cal E})$, donc est de la forme $z(\sigma)^{-1}e$, avec $z(\sigma)\in Z(G)$. On obtient l'\'egalit\'e $\sigma(e)=(\theta-1)(t(\sigma))z(\sigma)^{-1}e$. Puisque $\gamma=te$ et $\sigma(\gamma)=\gamma$, on a aussi $\sigma(t)=z(\sigma)(1-\theta)(t(\sigma))t$, donc $\sigma\circ p(t)=p(z(\sigma)t)$. $\square$

Levons les hypoth\`eses pr\'ec\'edentes et supposons $\gamma$ fortement r\'egulier. Alors

(5) $p(t)$ est r\'egulier au sens que son fixateur dans $W^{\theta}$ est r\'eduit \`a l'unit\'e.

 Preuve. Soit $\omega\in W^{\theta}$ qui fixe $p(t)$. On peut relever $\omega$ en un \'el\'ement $n\in G_{e}$. L'\'egalit\'e $\omega\circ p(t)=p(t)$ signifie qu'il existe $t'\in T$ tel que $t'ntn^{-1}\theta(t')^{-1}=t$. Mais alors $t'n\gamma(t'n)^{-1}=\gamma$ donc $t'n\in Z_{G}(\gamma)$. Puisque $\gamma$ est fortement  r\'egulier,  $Z_{G}(\gamma)=T^{\theta}$ et cela entra\^{\i}ne $\omega=1$. $\square$
 
 Remarquons que si $\gamma\in \tilde{G}_{reg}(F)$,  $T$ est  uniquement d\'etermin\'e par $\gamma$ et est d\'efini sur $F$: c'est le commutant dans $G$ de $G_{\gamma}$.  
 
 Soit $(B,T)$ une paire de Borel de $G$.   Soit $\tilde{T}$ le normalisateur commun de $B$ et $T$. Nous dirons que $\tilde{T}$ est un tore tordu maximal de $\tilde{G}$ si $T$ est d\'efini sur $F$ (mais pas forc\'ement $B$) et $\tilde{T}\cap \tilde{G}(F)$ est non vide. Dans ce cas, $\tilde{T}$ est aussi d\'efini sur $F$. Pour un tel tore tordu, notons $\theta$ l'automorphisme $ad_{\gamma}$ de $T$ pour un \'el\'ement quelconque $\gamma\in \tilde{T}$. On dit que $\tilde{T}$ est elliptique si et seulement si le plus grand sous-tore d\'eploy\'e de $T^{\theta,0}$ est $A_{\tilde{G}}$. 

 \bigskip
 
 \subsection{$L$-groupes}
 D\'esormais, $F$ sera soit un corps local, soit un corps de nombres. On note $W_{F}$ son groupe de Weil. Via l'homomorphisme naturel de $W_{F}$ dans $\Gamma_{F}$, le groupe $W_{F}$ agit sur tout ensemble sur lequel agit $\Gamma_{F}$.

 Soit $\hat{G}$ le groupe dual de $G$. Rappelons ce que cela signifie. C'est un groupe r\'eductif connexe d\'efini sur ${\mathbb C}$. On d\'efinit comme en 1.2 $ \underline{sa}$ paire de Borel \'epingl\'ee $\hat{\boldsymbol{{\cal E}}}=(\hat{\bf B},\hat{\bf T},(\hat{{\bf E}}_{\alpha})_{\alpha\in \hat{\Delta}})$. Des isomorphismes en dualit\'e $X_{*}(T^*)\simeq X^*(\hat{\bf T})$, $X^*(T^*)\simeq X_{*}(\hat{\bf T})$ sont donn\'es, qui \'echangent ensembles de racines et ensembles de coracines et respectent les ordres d\'efinis par $B^*$ et $\hat{\bf B}$. Le groupe $\hat{G}$ est muni d'une action alg\'ebrique de $\Gamma_{F}$ not\'ee $w\mapsto w_{G}$.   Il en r\'esulte une action sur $\hat{\boldsymbol{{\cal E}}}$. On suppose que les isomorphismes ci-dessus sont \'equivariants pour les actions galoisiennes. On suppose de plus que $\hat{G}$ poss\`ede une paire de Borel \'epingl\'ee qui est conserv\'ee par l'action galoisienne.   
On note $^LG$ le produit semi-direct $\hat{G}\rtimes W_{F}$.
 
 Par dualit\'e, il se d\'eduit de $\theta^*$ un automorphisme $\hat{\theta}$ de $\hat{\bf T}$. Soulignons que $\theta^*\mapsto \hat{\theta}$ est bien une dualit\'e, c'est-\`a-dire est contravariante.  Identifions $\hat{\boldsymbol{\cal E}}$ \`a une paire de Borel \'epingl\'ee de $\hat{G}$ conserv\'ee par l'action galoisienne. Alors $\hat{\theta}$ se prolonge de fa\c{c}on unique en un automorphisme $\hat{\boldsymbol{\theta}}$ de $\hat{G}$ qui pr\'eserve cette paire. L'automorphisme $\hat{\boldsymbol{\theta}}$ commute \`a l'action de $\Gamma_{F}$.  Remarquons que l'ensemble $\hat{G}\hat{\boldsymbol{\theta}}$ est naturellement un espace tordu sous $\hat{G}$, d\'efini sur ${\mathbb C}$. Cela nous permet d'utiliser pour lui les notations et terminologie introduites pour $\tilde{G}$. On peut aussi introduire l'espace $^L\tilde{G}={^LG}\hat{\boldsymbol{\theta}}$ qui est, en un sens convenable, un espace tordu sous $^LG$.
 
  Il est g\^enant de se limiter aux paires de Borel \'epingl\'ees de $\hat{G}$ conserv\'ees par l'action galoisienne, l'ensemble de ces paires n'\'etant pas invariant par conjugaison. On peut s'affranchir de cette limitation de la fa\c{c}on suivante. Soit $\hat{\cal E}$ une paire de Borel \'epingl\'ee quelconque de $\hat{G}$. On choisit $y\in \hat{G}_{SC}$ (le rev\^etement simplement connexe de $\hat{G}$) tel que $ad_{y^{-1}}(\hat{\cal E})$ soit la paire que l'on a fix\'ee ci-dessus. On d\'efinit une nouvelle action de $\Gamma_{F}$ sur $\hat{G}$ par $w\mapsto ad_{y}w_{G}ad_{y^{-1}}$. Elle conserve $\hat{\cal E}$. Le groupe $^LG$ est encore le produit semi-direct $\hat{G}\rtimes W_{F}$ pour cette nouvelle action: on envoie $(g,w)$ sur $(gw_{ G}(y)y^{-1},w)$. On pose $\hat{\theta}=y\hat{\boldsymbol{\theta}}y^{-1}\in {^L\tilde{G}}$. L'automorphisme d\'eduit de $\hat{\theta}$ (que l'on note encore $\hat{\theta}$) conserve $\hat{\cal E}$, commute \`a la nouvelle action galoisienne et on a l'\'egalit\'e $\hat{G}\hat{\boldsymbol{\theta}}=\hat{G}\hat{\theta}$. Ces d\'efinitions d\'ependent du choix de $y$ qui n'est d\'etermin\'e que modulo $Z(\hat{G}_{SC})$, mais ce choix s'av\'erera sans importance. Ainsi, pour une paire $\hat{\cal E}$ fix\'ee, on choisira $y$, on d\'efinira $\hat{\theta}$ comme ci-dessus et une action galoisienne, que l'on notera encore $w\mapsto w_{G}$ en esp\'erant que cela ne cr\'ee pas d'ambigu\"{\i}t\'e.

 \bigskip
 \subsection{Donn\'ees endoscopiques}
Pour la suite de l'article, $F$ est un corps local de caract\'eristique nulle, $G$ est un groupe r\'eductif connexe   et $\tilde{G}$ est  un espace tordu sous $G$, tous deux d\'efinis sur $F$.   On fixe de plus une classe de cohomologie ${\bf a}\in H^1(W_{F},Z(\hat{G}))$. D'apr\`es un th\'eor\`eme de Langlands, ce groupe de cohomologie s'envoie surjectivement, et m\^eme bijectivement si $F\not={\mathbb R}$, sur le groupe des caract\`eres continus de $G(F)$ (on rappellera cette correspondance en 1.13).  On note $\omega$ le caract\`ere de $G(F)$ associ\'e \`a ${\bf a}$. On impose les hypoth\`eses suivantes:

$\bullet$ $\tilde{G}(F)\not=\emptyset$;

$\bullet$ $\theta^*$ est d'ordre fini.

On peut aussi imposer l'hypoth\`ese

$\bullet$ $\omega$ est trivial sur $Z(G;F)^{\theta}$,

\noindent sinon toute la th\'eorie est vide. Mais, parce que cette hypoth\`ese n'est pas stable par passage \`a un groupe de Levi, il vaut mieux ne pas l'imposer.

Une donn\'ee endoscopique pour $(G,\tilde{G},{\bf a})$ est un triplet ${\bf G}'=(G',{\cal G}',\tilde{s})$ v\'erifiant les conditions qui suivent. Le terme $G'$ est un groupe r\'eductif connexe d\'efini et quasi-d\'eploy\'e sur $F$. Le terme $\tilde{s}$ est un \'el\'ement semi-simple de $\hat{G}\hat{\boldsymbol{\theta}}$.
 Le terme ${\cal G}'$ est un sous-groupe ferm\'e de $^LG$. On suppose que ${\cal G}'\cap \hat{G}=\hat{G}_{\tilde{s}}$ (composante neutre du commutant de $\tilde{s}$). On a donc une suite:
$$1\to \hat{G}_{\tilde{s}}\to {\cal G}'\to W_{F}\to 1,$$
o\`u la troisi\`eme fl\`eche est la restriction de la projection naturelle de $^LG$ sur $W_{F}$. On suppose que cette suite est exacte et scind\'ee, c'est-\`a-dire qu'il existe une section $W_{F}\to {\cal G}'$ qui soit  un homomorphisme continu. Fixons une paire de Borel \'epingl\'ee $\hat{\cal E}'=(\hat{ B}',\hat{T}',(\hat{E}'_{\alpha'})_{\alpha'\in \Delta'})$ de $\hat{G}_{\tilde{s}}$. Pour $w\in W_{F}$, on peut choisir  $g_{w}=(g(w),w)\in {\cal G}'$ tel que $ad_{g_{w}}$ conserve cette paire. L'application $w\mapsto w_{G'}=ad_{g_{w}}$ s'\'etend en une action galoisienne de $\Gamma_{F}$ sur $\hat{G}_{\tilde{s}}$. On suppose que $\hat{G}_{\tilde{s}}$ muni de cette action est un groupe dual de $G'$. Cela nous autorise \`a noter $\hat{G}_{\tilde{s}}=\hat{G}'$. On suppose enfin qu'il existe un cocycle $a:W_{F}\to Z(\hat{G})$, dont la classe est ${\bf a}$, tel que pour tout $(g,w)\in {\cal G}'$, on ait l'\'egalit\'e
$$ad_{\tilde{s}}(g,w)=(a(w)g,w).$$

 Soient ${\bf G}'_{1}=(G'_{1},{\cal G}'_{1},\tilde{s}_{1})$ et  ${\bf G}'_{2}=(G'_{2},{\cal G}'_{2},\tilde{s}_{2})$ deux donn\'ees comme ci-dessus. Une \'equivalence entre ces donn\'ees est un \'el\'ement $x\in \hat{G}$ tel que $x{\cal G}'_{1}x^{-1}={\cal G}'_{2}$ et $x \tilde{s}_{1}x^{-1}\in Z(\hat{G})\tilde{s}_{2}$. De $ad_{x}^{-1}:\hat{G}'_{2}\to \hat{G}'_{1}$ se d\'eduit par dualit\'e un automorphisme $\alpha_{x}:G'_{1}\to G'_{2}$ d\'efini sur $F$, ou plus exactement une classe de tels isomorphismes  modulo l'action de l'un ou l'autre des groupes $G'_{1,AD}(F)$ ou $G'_{2,AD}(F)$. En particulier, pour une seule donn\'ee ${\bf G}'$, on note $Aut({\bf G}')$ le groupe de ses automorphismes, c'est-\`a-dire des \'equivalences entre cette donn\'ee et elle-m\^eme. Ce groupe contient $\hat{G}'$. Notons $Out({\bf G}')$ le sous-groupe form\'e des $\alpha_{x}$ dans le groupe $Out(G')$ des automorphismes ext\'erieurs de $G'$. On a une suite exacte ([KS1] p.19)
 $$1\to (Z(\hat{G})/(Z(\hat{G})\cap \hat{G}'))^{\Gamma_{F}}\to Aut({\bf G}')/\hat{G}'\to Out({\bf G}')\to 1.$$
 
 Soit ${\bf G}'=(G',{\cal G}',\tilde{s})$ une donn\'ee endoscopique pour $(G,\tilde{G},{\bf a})$. Fixons une paire de Borel \'epingl\'ee $\hat{{\cal E}}=(\hat{B},\hat{T},(\hat{E}_{\alpha})_{\alpha\in \Delta})$ de $\hat{G}$ telle que $ad_{\tilde{s}}$ conserve $\hat{B}$ et $\hat{T}$. Posons $\hat{B}'=\hat{B}\cap \hat{G}'$, $\hat{T}'=\hat{T}\cap \hat{G}'$ et compl\'etons $(\hat{B}',\hat{T}')$ en une paire de Borel \'epingl\'ee $\hat{\cal E}'=(\hat{ B}',\hat{T}',(\hat{E}'_{\alpha'})_{\alpha'\in \Delta'})$ de $\hat{G}'$. Ainsi qu'on l'a expliqu\'e en 1.4,   en r\'ef\'erence \`a la paire $\hat{\cal E}$, on modifie l'action $\sigma\mapsto \sigma_{G}$ de $\Gamma_{F}$ sur $\hat{G}$, on modifie l'isomorphisme $^LG\simeq \hat{G}\rtimes W_{F}$ et on d\'efinit l'\'el\'ement $\hat{\theta}\in \hat{G}\hat{\boldsymbol{\theta}}$. On peut \'ecrire $\tilde{s}=s\hat{\theta}$, avec $s\in \hat{T}$. On construit comme ci-dessus l'action galoisienne $\sigma\mapsto \sigma_{G'}$ qui conserve $\hat{{\cal E}}'$. On a l'\'egalit\'e $\hat{T}'=\hat{T}^{\hat{\theta},0}$. Cette \'egalit\'e identifie le groupe de Weyl $W'$ de $\hat{G}'$ (ou $G'$) \`a un sous-groupe des \'el\'ements invariants par $\hat{\theta}$ du groupe de Weyl de $\hat{G}$, lequel s'identifie par dualit\'e \`a $W^{\theta^*}$. Le plongement $\hat{\xi}:\hat{T}'\subset \hat{T}$ n'est pas \'equivariant pour les actions galoisiennes. Il existe un cocycle $\omega_{G'}:\Gamma_{F}\to W^{\theta^*}$ tel que $\omega_{G'}(\sigma)\circ\sigma(\hat{\xi})=\hat{\xi}$. Remarquons que le groupe $Z(\hat{G})\cap \hat{G}'$ qui intervient dans la suite exacte ci-dessus est \'egal \`a $Z(\hat{G})\cap T^{\hat{\theta},0}$. Introduisons $\underline{la}$ paire de Borel \'epingl\'ee ${\cal E}^{_{'}*}=(B^{_{'}*},T^{_{'}*},(E^{_{'}*}_{\alpha'})_{\alpha'\in \Delta'})$ de $G'$. Les tores $\hat{ T}$ et $\hat{T}'$ sont duaux de $T^*$ et $T^{_{'}*}$. Le tore $\hat{ T}^{\hat{\theta},0}$ est dual de $T^*/(1-\theta^*)(T^*)$. Du plongement $\hat{\xi}$ se d\'eduit par dualit\'e un homomorphisme
 $$\xi:T^*\to T^*/(1-\theta^*)(T^*)\simeq T^{_{'}*}.$$
 Pour $\sigma\in \Gamma_{F}$, on a l'\'egalit\'e $\sigma(\xi)=\xi\circ \omega_{G'}(\sigma)$.
 
 Les constructions ci-dessus d\'ependent du choix de la paire $\hat{\cal E}$. La plupart du temps, pour une donn\'ee endoscopique ${\bf G}'$ fix\'ee, on supposera choisie une telle paire et on utilisera ces constructions sans plus de commentaires.

\bigskip
\subsection{Syst\`emes de racines}
Notons $\Sigma(T^*)$ l'ensemble des racines de $T^*$ dans l'alg\`ebre de Lie de $G$, $\Sigma(\hat{T})$ celui des racines de $\hat{T}$ dans l'alg\`ebre de Lie de $\hat{G}$ et $\check{\Sigma}(T^*)$, $\check{\Sigma}(\hat{T})$ les ensembles de coracines.  Par les isomorphismes $X_{*}(T^*)\simeq X^*(\hat{T})$, $X^*(T^*)\simeq X_{*}(\hat{T})$, l'ensemble $\Sigma(T^*)$ s'identifie \`a $\check{\Sigma}(\hat{T})$ et $\check{\Sigma}(T^*)$ s'identifie \`a $\Sigma(\hat{T})$. On note $\alpha\mapsto \hat{\alpha}$ la bijection de $\Sigma(T^*)$ sur $\Sigma(\hat{T})$ telle que, par les identifications pr\'ec\'edentes, $\hat{\alpha}$ s'identifie \`a la coracine $\check{\alpha}$. Pour $\alpha\in \Sigma(T^*)$, on note $N\alpha$  la somme des \'el\'ements de l'orbite de $\alpha$  sous l'action du groupe d'automorphismes engendr\'e par $\theta^*$. On note $\alpha_{res}$ la restriction de $\alpha$ \`a $T^{*,\theta^*,0}$. On pose $\Sigma(T^*)_{res}=\{\alpha_{res};\alpha\in \Sigma(T^*)\}$. De m\^eme, pour $\alpha\in \Sigma(\hat{T})$, on note $N\alpha$ la somme des \'el\'ements de l'orbite de $\alpha$  sous l'action du groupe d'automorphismes engendr\'e par $\hat{\theta}$. On note $\alpha_{res}$ la restriction de $\alpha$ \`a $\hat{T}^{\hat{\theta},0}$. On pose $\Sigma(\hat{T})_{res}=\{\alpha_{res};\alpha\in \Sigma(\hat{T})\}$.  Les ensembles $\Sigma(T^*)_{res} $ et $\Sigma(\hat{T})_{res} $ sont des syst\`emes de racines non r\'eduits en g\'en\'eral. On dit que $\alpha\in \Sigma(T^*)$ est de type 1 si ni $\alpha_{res}/2$, ni  $2\alpha_{res}$ n'appartiennent  \`a $\Sigma(T^*)_{res}$, de type 2 si $2\alpha_{res}\in \Sigma(T^*)_{res}$ et de type 3 si $\alpha_{res}/2
\in \Sigma(T^*)_{res}$. On d\'efinit de m\^eme le type d'une racine $\alpha\in \Sigma(\hat{T})$. Pour $\alpha\in \Sigma(T^*)$, l'\'el\'ement $\hat{\alpha}\in \Sigma(\hat{T})$ est de m\^eme type que $\alpha$.

Soit ${\bf G}'=(G',{\cal G}',\tilde{s})$ une donn\'ee endoscopique pour $(G,\tilde{G},{\bf a})$. L'ensemble $\Sigma(\hat{T}')$ des racines de $\hat{ T}'$ dans l'alg\`ebre de Lie de $\hat{G}'$ est form\'e des $\alpha_{res}$ pour $\alpha\in \Sigma(\hat{T})$ telles que
$$N\alpha(s)=\left\lbrace\begin{array}{cc}1,&\text{ si }\alpha\text{ est de type 1 ou 2 }\\ -1,&\text{ si }\alpha\text{ est de type 3.}\\ \end{array}\right.$$ 
 (on rappelle que $\tilde{s}=s\hat{\theta}$). Par composition avec l'homomorphisme $\xi$, l'ensemble $\Sigma(T^{_{'}*})$ des racines de $T^{_{'}*}$ dans l'alg\`ebre de Lie de $G'$ s'identifie \`a  un ensemble de caract\`eres de $T^*$. D'apr\`es [KS1] 1.3.9, c'est l'ensemble suivant:
$$\{N\alpha; \alpha\in \Sigma(T^*) \text{ de type 1 },N\hat{\alpha}(s)=1\}$$
$$\cup \{2N\alpha; \alpha\in \Sigma(T^*) \text{ de type 2 },N\hat{\alpha}(s)=1\}$$
$$\cup \{N\alpha; \alpha\in \Sigma(T^*) \text{ de type 3 },N\hat{\alpha}(s)=-1\}.$$

 \bigskip
 \subsection{Espace endoscopique tordu}
   Soit ${\bf G}'=(G',{\cal G}',\tilde{s})$ une donn\'ee endoscopique.      On a
   
   (1) $\xi(Z(G))\subset Z(G')$.
   
   Preuve. Pour $z\in Z(G)$, on a $\alpha(z)=1$ pour tout $\alpha\in \Sigma(T^*)$. A fortiori $N\alpha(z)=1$. Pour toute racine $\alpha' \in \Sigma(T^{_{'}*})$, il existe $\alpha\in \Sigma(T^*)$ telle que $\alpha'\circ\xi=N\alpha$ ou $2N\alpha$.  Donc $\alpha'(\xi(z))=1$ pour tout $\alpha'\in \Sigma(T^{_{'}*})$ et cela \'equivaut \`a $\xi(z)\in Z(G')$. $\square$

 La restriction de $\xi$ \`a $Z(G)$ se quotiente \'evidemment en un homomorphisme $\xi_{{\cal Z}}:{\cal Z}(G)\to Z(G')$. On v\'erifie que celui-ci est \'equivariant pour les actions galoisiennes. On pose $\tilde{G}'=G'\times_{{\cal Z}(G)}{\cal Z}(\tilde{G})$, c'est-\`a-dire le quotient de $G'\times {\cal Z}(\tilde{G})$ par la relation d'\'equivalence $(g'\xi_{{\cal Z}}(z), \tilde{z})\equiv (g',z\tilde{z})$ pour $z\in {\cal Z}(G)$. Les actions \`a droite et \`a gauche de $G'$ sur $G'\times {\cal Z}(\tilde{G})$ se descendent en des actions \`a droite et \`a gauche sur $\tilde{G}'$. L'action galoisienne sur $G'\times {\cal Z}(\tilde{G})$ se descend aussi en une action sur $\tilde{G}'$. On voit que $\tilde{G}'$ est un espace tordu sur $G'$, d\'efini sur $F$.
 
{\bf Remarques.} (2) L'ensemble $\tilde{G}'(F)$ peut \^etre vide. Par exemple, soient $d\in F^{\times}$, $G=SL(2)$, $\tilde{G}=\{\gamma\in GL(2); det(\gamma)=d\}$ et ${\bf a}=1$. Pour toute extension quadratique  $E$ de $F$, il y a une donn\'ee endoscopique ${\bf G}'$ telle que $G'(F)$ est le groupe des \'el\'ements de $E$ de norme $1$. Alors $\tilde{G}'(F)$ est l'ensemble des \'el\'ements de $E$ de norme $d$. On peut trouver choisir $d$ et $E$ de sorte que cet ensemble soit vide.
 
 (3) $\tilde{G}'$ est \`a torsion int\'erieure.
 
  {\bf Cas particulier.} On dira que $(G,\tilde{G},{\bf a})$ est quasi-d\'eploy\'e et \`a torsion int\'erieure si $G$ est quasi-d\'eploy\'e sur $F$, $\tilde{G}$ est \`a torsion int\'erieure et ${\bf a}=1$. Dans ce cas, on a $\hat{\boldsymbol{\theta}}=1$ et la donn\'ee ${\bf G}=(G,{^LG},\tilde{s}=1)$ est une donn\'ee endoscopique "maximale". L'espace endoscopique que l'on en d\'eduit est bien s\^ur l'espace $\tilde{G}$ lui-m\^eme. Remarquons que, pour toute donn\'ee endoscopique ${\bf G}'=(G',{\cal G}',\tilde{s})$, le couple $(G',\tilde{G}')$ compl\'et\'e par le cocycle trivial est quasi-d\'eploy\'e et \`a torsion int\'erieure.
  
\bigskip

 \subsection{Correspondance entre classes de conjugaison semi-simples}
 Soit $\gamma\in \tilde{G}_{ss}$. Par d\'efinition des \'el\'ements semi-simples, on peut fixer une paire de Borel $(B,T)$ de $G$ qui est conserv\'ee par $ad_{\gamma}$. On la compl\`ete en une paire de Borel \'epingl\'ee ${\cal E}$. On identifie cette paire \`a ${\cal E}^*$. D'apr\`es 1.3(1), on peut \'ecrire $\gamma=te$, avec $t\in T$ et $e\in Z(\tilde{G},{\cal E})$.   Soit $\bar{t}$ l'image de $t$ dans $(T^*/(1-\theta^*)(T^*))/W^{\theta^*}$, $\bar{e}$ l'image de $e$ dans ${\cal Z}(\tilde{G})$ et $\bar{\gamma}$ l'image de $(\bar{t},\bar{e})$ dans $((T^*/(1-\theta^*)(T^*))/W^{\theta})\times_{{\cal Z}(G)}{\cal Z}(\tilde{G})$. Montrons que:
 
 (1) l'\'el\'ement $\bar{\gamma}$ ne d\'epend pas des choix; l'application
$\gamma\mapsto \bar{\gamma}$ se quotiente en une bijection de l'ensemble des classes de conjugaison semi-simples dans $\tilde{G}$ sur $((T^*/(1-\theta^*)(T^*))/W^{\theta^*})\times_{{\cal Z}(G)}{\cal Z}(\tilde{G})$; cette bijection est d\'efinie sur $F$.

Preuve. Pour ${\cal E}$ fix\'ee, on peut remplacer $(t,e)$ par $(tz,z^{-1}e)$, avec $z\in Z(G)$. Cela remplace $(\bar{t},\bar{e})$ par $(\bar{t}\bar{z},\bar{z}^{-1}\bar{e})$, o\`u $\bar{z}$ est l'image de $z$ dans ${\cal Z}(G)$, et cela ne change pas  $\bar{\gamma}$. Laissons fix\'ee $(B,T)$, mais changeons d'\'epinglage. La nouvelle paire de Borel \'epingl\'ee ${\cal E}'$ se d\'eduit de ${\cal E}$ par $ad_{y}$ pour un $y\in T$. Posons $e'=ad_{y}(e)$. On a $e'\in Z(\tilde{G},{\cal E}')$ et $e'=(1-\theta)(y)e$ o\`u $\theta=\theta_{{\cal E}}=\theta_{{\cal E}'}$. On peut \'ecrire $\gamma=t'e'$ avec $t'=(\theta-1)(y)t$. On voit que $\bar{t}'=\bar{t}$ et $\bar{e}'=\bar{e}$. Donc $\bar{\gamma}$ ne change pas. Rempla\c{c}ons $(B,T)$ par une autre paire $(B',T)$ de m\^eme tore.  Comme on l'a vu dans la preuve de 1.3(2), la paire $(B',T)$ se d\'eduit de $(B,T)$ par l'action d'un \'el\'ement de $W^{\theta}$, que l'on peut repr\'esenter par un \'el\'ement $n\in G_{e}$. Posons ${\cal E}'=ad_{n}({\cal E})$. Alors $e$ appartient \`a $Z(\tilde{G},{\cal E}')$ et on peut changer ${\cal E}$ en ${\cal E}'$ tout en conservant la d\'ecomposition $\gamma=te$.  Parce que $e$ est fixe par $ad_{n}$, son image dans ${\cal Z}(\tilde{G})$ est la m\^eme, que la paire de r\'ef\'erence soit ${\cal E}$ ou ${\cal E}'$. Les identifications de $T$ \`a $T^*$ relatives aux deux paires ${\cal E}$ et ${\cal E}'$ diff\`erent par l'action d'un \'el\'ement de $W^{\theta^*}$, donc les applications compos\'ees $T\to (T^*/(1-\theta^*)(T^*))/W^{\theta^*}$ sont les m\^emes et $\bar{t}$ ne change pas quand on remplace ${\cal E}$ par ${\cal E}'$.  Donc $\bar{\gamma}$ ne change pas non plus. Rempla\c{c}ons maintenant $(B,T)$ par une paire quelconque $(B',T')$.  D'apr\`es la preuve de 1.3(2), il existe $g\in G_{\gamma}$ tel que $ad_{g}(T)=T'$. L'\'etape pr\'ec\'edente nous permet de changer $B$ de sorte que l'on ait aussi $ad_{g}(B)=B'$. On choisit alors ${\cal E}'=ad_{g}({\cal E})$ et pour d\'ecomposition $\gamma=t'e'$, avec $t'=ad_{g}(t)$ et $e'=ad_{g}(e)$. Les diverses applications relatives \`a ${\cal E}'$ sont  les compos\'ees des applications relatives \`a  ${\cal E}$ avec $ad_{g}^{-1}$. Donc $\bar{\gamma}$ ne change pas. Cela prouve la premi\`ere assertion. La deuxi\`eme est facile. Soit $\sigma\in \Gamma_{F}$.  On utilise une paire ${\cal E}$ pour calculer $\bar{\gamma}$ et la paire $\sigma({\cal E})$ pour calculer $\overline{\sigma(\gamma)}$. D'une d\'ecomposition $\gamma=te$ se d\'eduit la d\'ecomposition $\sigma(\gamma)=\sigma(t)\sigma(e)$. On a $\overline{\sigma(t)}=\sigma_{G^*}(\bar{t})$ et $\overline{\sigma(e)}=\sigma_{G^*}(\bar{e})$ par d\'efinition des actions galoisiennes sur $T^*$ et  ${\cal Z}(\tilde{G})$.   Donc  $\overline{\sigma(\gamma)}$ est bien l'image de $\bar{\gamma}$ par l'action  $\sigma_{G^*}$. $\square$

Soit ${\bf G}'=(G',{\cal G}',\tilde{s})$ une donn\'ee endoscopique pour $(G,\tilde{G},{\bf a})$.  Les classes de conjugaison semi-simples dans $\tilde{G}'$ sont de m\^eme param\'etr\'ees par $(T^{_{'}*}/W^{G'})\times_{{\cal Z}(G')}{\cal Z}(\tilde{G}')$. On a ${\cal Z}(G')=Z(G')$ et, par construction, ${\cal Z}(\tilde{G}')=Z(G')\times_{{\cal Z}(G)}{\cal Z}(\tilde{G})$. Donc 
$$(T^{_{'}*}/W^{G'})\times_{{\cal Z}(G')}{\cal Z}(\tilde{G}')=(T^{_{'}*}/W^{G'})\times_{ {\cal Z}(G)}{\cal Z}(\tilde{G}).$$
  En utilisant  l'isomorphisme $T'\simeq T^*/(1-\theta^*)(T^*)$  par lequel   $W^{G'}$ s'identifie \`a un sous-groupe de $ W^{\theta^*}$, on obtient une surjection
$$(T^{_{'}*}/W^{G'})\times_{ {\cal Z}(G)}{\cal Z}(\tilde{G})\to ((T^*/(1-\theta^*)(T^*))/W^{\theta^*})\times_{{\cal Z}(G)}{\cal Z}(\tilde{G}),$$
c'est-\`a-dire une surjection de l'ensemble des classes de conjugaison semi-simples dans $\tilde{G}'$ sur l'ensemble des classes de conjugaison semi-simples dans $\tilde{G}$. Cette application est d\'efinie sur $F$. 

{\bf Remarque.} Restreinte aux \'el\'ements invariants par $\Gamma_F$, l'application n'est plus surjective en g\'en\'eral. D'autre part, une classe de conjugaison semi-simple dans $\tilde{G}$ peut \^etre d\'efinie sur $F$ sans contenir d'\'el\'ement de $\tilde{G}(F)$.

\bigskip
On dit qu'un \'el\'ement de $\tilde{G}'_{ss}$ est $\tilde{G}$-fortement r\'egulier si et seulement si l'image de sa classe de conjugaison par l'application ci-dessus est une classe de conjugaison dans $\tilde{G}$ form\'ee d'\'el\'ements fortement r\'eguliers.

On note ${\cal D}({\bf G}')$ l'ensemble des couples $(\delta,\gamma)\in \tilde{G}'(F)\times \tilde{G}(F)$ form\'es d'\'el\'ements semi-simples dont les classes de conjugaison (sur $\bar{F}$) se correspondent et tels que $\gamma$ est fortement r\'egulier dans $\tilde{G}$. On dit que ${\bf G}'$ est "relevant" si ${\cal D}({\bf G}')$ n'est pas vide.

\bigskip

\subsection{Remarques sur le cas quasi-d\'eploy\'e et \`a torsion int\'erieure}
On suppose $(G,\tilde{G},{\bf a})$ quasi-d\'eploy\'e et \`a torsion int\'erieure. L'ensemble $Z(\tilde{G},{\cal E})$ attach\'e \`a une paire de Borel \'epingl\'ee ${\cal E}$ est en fait ind\'ependant de ${\cal E}$: c'est l'ensemble des $e\in \tilde{G}$ tels que $ad_{e}$ soit l'identit\'e. L'ensemble $Z(\tilde{G})$ s'identifie donc \`a ce m\^eme ensemble.

Soit ${\bf G}'=(G',{\cal G}',\tilde{s})$ une donn\'ee endoscopique de $(G,\tilde{G},{\bf a})$.

\ass{Lemme}{  Supposons $\tilde{G}'(F)\not=\emptyset$. Alors l'ensemble des \'el\'ements $\tilde{G}$-fortement r\'eguliers de $\tilde{G}'(F)$ n'est pas vide et, pour tout  \'el\'ement $\delta$ de cet ensemble, il existe $\gamma\in \tilde{G}_{reg}(F)$ tel que $(\delta,\gamma)\in {\cal D}({\bf G}')$. A fortiori, ${\bf G}'$ est relevant.  }

 Preuve.  Puisque $\tilde{G}'(F)$ n'est pas vide, le sous-ensemble $\tilde{G}'_{ss}(F)$ ne l'est pas non plus: la partie semi-simple d'un \'el\'ement de $\tilde{G}'(F)$ appartient \`a cet ensemble. Soit $\epsilon\in \tilde{G}'_{ss}(F)$. Fixons un tore maximal $T'$ de $G'_{\epsilon}$ d\'efini sur $F$.  Pour $t'\in T'(F)$ en position g\'en\'erale, $t'\epsilon$ est $\tilde{G}$-fortement r\'egulier. D'o\`u la premi\`ere assertion. Fixons maintenant un \'el\'ement $\delta\in \tilde{G}'(F)$ qui soit $\tilde{G}$-fortement r\'egulier. Fixons une paire de Borel $(B',T')$ de $G'$ qui soit conserv\'ee par $ad_{\delta}$. On a $T'=G'_{\delta}$, donc $T'$ est d\'efini sur $F$. Soit $(B^*,T^*)$ une paire de Borel de $G$ d\'efinie sur $F$. Des deux paires de Borel se d\'eduit un isomorphisme $\xi_{T^*,T'}:T^*\to T'$. Il existe un cocycle $\omega_{T'}:\Gamma_{F}\to W$ tel que $\xi\circ \omega_{T'}(\sigma)\circ\sigma=\sigma\circ\xi$ pour tout $\sigma\in \Gamma_{F}$. Puisque $G$ est quasi-d\'eploy\'e, on peut appliquer le corollaire 2.2 de [K1]: il existe $g\in G(\bar{F})$ tel que $ad_{g^{-1}}(T^*)$ soit d\'efini sur $F$ et que, pour tout $\sigma\in \Gamma_{F}$,  on ait l'\'egalit\'e suivante sur $T$: $\omega_{T'}(\sigma)\circ \sigma\circ ad_{g}=ad_{g}\circ \sigma$. Posons $(B,T)=ad_{g^{-1}}(B^*,T^*)$.  De $(B,T)$ et $(B',T')$ se d\'eduit un isomorphisme $\xi_{T,T'}:T\to T'$ qui est maintenant \'equivariant pour les actions galoisiennes. On v\'erifie  que $\xi_{T,T'}$ s'\'etend en un isomorphisme $\tilde{\xi}_{T,T'}:T\times_{Z(G)}{\cal Z}(\tilde{G})\to T'\times_{Z(G')}{\cal Z}(\tilde{G}')$ qui est encore \'equivariant pour les actions galoisiennes. L'\'el\'ement $\delta$ appartient \`a l'ensemble d'arriv\'ee. Soit $\gamma$ son image r\'eciproque par $\tilde{\xi}_{T,T'}$. Puisque $\tilde{\xi}_{T,T'}$ est \'equivariant pour les actions galoisiennes, $\gamma$ appartient \`a $\tilde{G}(F)$ et il est clair que $(\delta,\gamma)$ appartient \`a ${\cal D}({\bf G}')$. $\square$

\bigskip

\subsection{Correspondance entre \'el\'ements semi-simples}
Soit ${\bf G}'=(G',{\cal G}',\tilde{s})$ une donn\'ee endoscopique pour $(G,\tilde{G},{\bf a})$. Appelons diagramme un sextuplet $(\epsilon,B',T',B,T,\eta)$ v\'erifiant les conditions (1) \`a (6) suivantes:  

(1) $\epsilon\in \tilde{G}'_{ss}(F)$ et $\eta\in \tilde{G}_{ss}(F)$;

(2) $(B',T')$ est une paire de Borel de $G'$ et $(B,T)$ est une paire de Borel  de $G$;

(3) $ad_{\epsilon}$ conserve $(B',T')$ et $ad_{\eta}$ conserve $(B,T)$;

(4) $T$ et $T'$ sont d\'efinis sur $F$.

A l'aide de $(B',T')$, resp. $(B,T)$, on identifie $T'$ \`a $T^{_{'}*}$ et $T$ \`a $T^*$. L'homomorphisme $\xi$ se transforme en un homomorphisme $\xi_{T,T'}:T\to T'$.

(5)  L'homomorphisme $\xi_{T,T'}$ est d\'efini sur $F$.

 Compl\'etons $(B,T)$ en une paire de Borel \'epingl\'ee ${\cal E}$, \'ecrivons $\eta=te$, avec $e\in Z(\tilde{G},{\cal E})$ et $t\in T$, cf. 1.3(1). Notons $e'$ l'image de $e$ dans ${\cal Z}(\tilde{G}')$. L'\'el\'ement $\xi_{T,T'}(t)e'$ de $\tilde{G}'$ ne d\'epend pas de ces choix: la preuve de cette assertion est contenue dans celle de 1.8(1). Alors

(6)   pour de quelconques choix comme ci-dessus, $\epsilon=\xi_{T,T'}(t)e'$.

{\bf Remarque.} Soit un diagramme $(\epsilon,B',T',B,T,\eta)$ et soit $B'_{1}$ un sous-groupe de Borel de $G'$ contenant $T'$. Il existe un unique \'el\'ement $w$ du groupe de Weyl $W^{G'}$ de $G'$ relativement \`a $T'$ tel que $B'_{1}=w(B')$. Cet \'el\'ement s'identifie \`a un \'el\'ement de $W$ (le groupe de Weyl de $G$ relativement \`a $T)$ qui est invariant par $\theta=\theta_{e}$ pour $e$ comme ci-dessus. Posons $B_{1}=w(B)$. Alors $(\epsilon,B_{1}',T',B_{1},T,\eta)$ est encore un diagramme. 
\bigskip

Pour $\epsilon$ et $\eta$ v\'erifiant (1), on dit que $\epsilon$ et $\eta$ se correspondent si et seulement s'il existe un diagramme joignant $\epsilon$ \`a $\eta$. Il est clair que si $\epsilon$ et $\eta$ se correspondent, les classes de conjugaison sur $\bar{F}$ de $\epsilon$ et $\eta$ se correspondent. La r\'eciproque est fausse en g\'en\'eral, c'est-\`a-dire que, si les classes de conjugaison sur $\bar{F}$ de $\epsilon$ et $\eta$ se correspondent, il n'existe pas toujours de diagramme joignant $\epsilon$ et $\eta$. Le lemme suivant  pr\'ecise ce point.

\ass{Lemme}{(i) Soit $(\delta,\gamma)\in {\cal D}({\bf G}')$. Alors il existe un diagramme $(\delta,B',T',B,T,\gamma)$.

(ii) Soient  $\epsilon\in \tilde{G}'_{ss}(F)$ et $\eta\in \tilde{G}_{ss}(F)$. Alors ces deux \'el\'ements se correspondent si et seulement si $(\epsilon,\eta)$ appartient \`a l'adh\'erence de ${\cal D}({\bf G}')$.}

Preuve. (i) On fixe $(B',T')$ et $(B,T)$ tels que (3) soit v\'erifi\'ee (pour $\epsilon=\delta$, $\eta=\gamma$). Les tores $T$ et $T'$ sont uniquement d\'etermin\'es puisque nos \'el\'ements sont fortement r\'eguliers, donc (4) est v\'erifi\'ee.  On compl\`ete $(B,T)$ en une paire de Borel \'epingl\'ee ${\cal E}$. Il existe un cocyle $\omega_{T',T}:\Gamma_{F}\to W^{\theta}$ (o\`u $\theta=\theta_{{\cal E}}$) tel que $\sigma_{T'}\circ\xi_{T,T'}=\xi_{T,T'}\circ\omega_{T',T}(\sigma)\circ\sigma_{T}$. On \'ecrit $\gamma=te$, avec $t\in T$ et $e\in Z(\tilde{G},{\cal E})$. On peut aussi \'ecrire $\delta=t'e'$ o\`u $t'\in T'$ et $e'$ est l'image de $e$ dans ${\cal Z}(\tilde{G}')$. L'hypoth\`ese que les classes de conjugaison de $\delta$ et $\gamma$ se correspondent signifie qu'il existe $w\in W^{\theta}$ tel que $\xi_{T,T'}\circ w(t)=t'$. On peut relever $w$ en un \'el\'ement $n$ de $G_{SC,e}$ qui normalise $T$. Rempla\c{c}ons  ${\cal E}$ par ${\cal E}_{1}=ad_{n^{-1}}({\cal E})$. Cela remplace $\xi_{T,T'}$   par $\xi_{T,T',1}=\xi_{T,T'}\circ w$. On a alors $\xi_{T,T',1}(t)=t'$. En oubiant cette construction, on suppose $\xi_{T,T'}(t)=t'$. Soit $\sigma\in \Gamma_{F}$. D'apr\`es 1.3(4), il existe $z(\sigma)\in Z(G)$ tel que $ad_{u_{{\cal E}}(\sigma)}
\circ \sigma(e)=z(\sigma)^{-1}e$ et l'image de $\sigma(t)$ dans $T/(1-\theta)(T)$ soit \'egale \`a celle de $t$ multipli\'ee par $z(\sigma)$ (en notant encore $z(\sigma)$ l'image de cet \'el\'ement dans les divers quotients de $Z(G)$). La premi\`ere relation entra\^{\i}ne $\sigma_{G^*}(\bar{e})=z(\sigma)^{-1}\bar{e}$  (o\`u $\bar{e}$ est l'image de $e$ dans ${\cal Z}(\tilde{G})$) puis $\sigma_{G'}(e')=z(\sigma)^{-1}e'$.  La seconde relation entra\^{\i}ne $\xi_{T,T'}\circ\sigma (t)=z(\sigma)\xi_{T,T'}(t)=z(\sigma)t'$. On a aussi
$$t'e=\delta=\sigma(\delta)=\sigma(t')\sigma(e')=\sigma(t')z(\sigma)^{-1}e',$$
d'o\`u $\sigma(t')=z(\sigma)t'$. Alors $\xi_{T,T'}\circ\sigma(t)=\sigma\circ \xi_{T,T'}(t)$. Mais ce terme est aussi \'egal \`a $\xi_{T,T'}\circ\omega_{T',T}(\sigma)\circ\sigma(t)$. D'o\`u $\omega_{T',T}(\sigma)=1$ puisque $\gamma$ est fortement r\'egulier, cf. 1.3(5). Cela prouve (i).

(ii) Supposons que $\epsilon$ et $\eta$ se correspondent. Fixons un diagramme $(\epsilon,B',T',B,T,\eta)$. Soit $t\in T(F)$, posons $t'=\xi_{T,T'}(t)$. Alors $(t'\epsilon,B',T',B,T,t\eta)$ est encore un diagramme.  Si $t$ est en position g\'en\'erale, $t\eta$ est fortement r\'egulier. Donc $(t'\epsilon,t\eta)\in {\cal D}({\bf G}')$. On peut choisir $t$ aussi proche de $1$ que l'on veut.  Donc $(\epsilon,\eta)$ appartient \`a l'adh\'erence de ${\cal D}({\bf G}')$. Inversement, supposons cette condition v\'erifi\'ee. On fixe une suite d'\'el\'ements $(\delta_{n},\gamma_{n})\in {\cal D}({\bf G}')$, pour $n\in {\mathbb N}$, qui tend vers $(\epsilon,\eta)$. Les r\'esultats usuels de la th\'eorie de la descente valent dans le cas tordu. En notant par des lettres gothiques les alg\`ebres de Lie, on peut fixer un voisinage $\mathfrak{u}_{\eta}$ de $0$ dans $\mathfrak{g}_{\eta}(F)$ de sorte que tout point assez voisin de $\eta$ soit conjugu\'e par un  \'el\'ement de $G(F)$ \`a un \'el\'ement $exp(X)\eta$ o\`u $X\in \mathfrak{u}_{\eta}$. On peut fixer un voisinage similaire $\mathfrak{u}_{\epsilon}$ de $0$ dans $\mathfrak{g}'_{\epsilon}(F)$. Quitte \`a conjuguer nos \'el\'ements $\delta_{n}$ et $\gamma_{n}$ et \`a supprimer un nombre fini de termes de la suite, on peut donc \'ecrire $\delta_{n}=exp(Y_{n})\epsilon$, $\gamma_{n}=exp(X_{n})\eta$. Puisqu'il s'agit d'\'el\'ements semi-simples, les $X_{n}$ et $Y_{n}$ le sont aussi. Puisqu'il n'y a qu'un nombre fini de classes de conjugaison par $G_{\eta}(F)$ de sous-tores maximaux de $G_{\eta}$ d\'efinis sur $F$ (et de m\^eme pour $G'_{\epsilon}$), on peut, quitte \`a extraire une sous-suite, fixer de tels sous-tores maximaux $T^{\natural}\subset G_{\eta}$ et $T'\subset G'_{\epsilon}$ et supposer $X_{n}\in \mathfrak{t}^{\natural}(F)$, $Y_{n}\in \mathfrak{t}'(F)$. D'apr\`es (i), on peut fixer des diagrammes $(\delta_{n},B'_{n},T'_{n},B_{n},T_{n},\gamma_{n})$. Il n'y a pas le choix pour les tores:  on a n\'ecessairement $T'_{n}=T'$ tandis que $T_{n}$ est  le commutant de $T^{\natural}$ dans $G$. Puisque ces tores n'appartiennent qu'\`a un nombre fini de paires de Borel, on peut, quitte \`a extraire une sous-suite, fixer $B$ contenant $T$ et $B'$ contenant $T'$ et supposer que $B_{n}=B$ et $B'_{n}=B'$ pour tout $n$.   Puisque $\gamma_{n}\in T(F)\eta$  et que $ad_{\gamma_{n}}$ conserve $(B,T)$, $ad_{\eta}$ conserve aussi cette paire. On \'ecrit $\eta=te$ comme au d\'ebut du paragraphe, avec $t\in T$. De m\^eme, on peut \'ecrire $\epsilon=t'e'$, o\`u $e'$ est l'image de $e$ dans ${\cal Z}(\tilde{G}')$ et $t'\in T'$. On a alors $\gamma_{n}=exp(X_{n})te$, $\delta_{n}=exp(Y_{n})t'e'$. D'apr\`es (6) appliqu\'e au diagramme joignant $\delta_{n}$ et $\gamma_{n}$, on a $\xi_{T,T'}(exp(X_{n})t)=exp(Y_{n})t'$. Quand $n$ tend vers l'infini, $X_{n}$ et $Y_{n}$ tendent vers $0$. D'o\`u $\xi_{T,T'}(t)=t'$. Mais alors $(\epsilon,B',T',B,T,\eta)$ est un diagramme. Cela ach\`eve la preuve. $\square$

\bigskip

\subsection{$K$-espaces}
 On suppose dans ce paragraphe $F={\mathbb R}$. Consid\'erons une famille finie $(G_{p},\tilde{G}_{p})_{p\in \Pi}$, o\`u, pour tout $p$, $G_{p}$ est un groupe r\'eductif connexe sur ${\mathbb R}$ et $\tilde{G}_{p}$ est un espace tordu sur $G_{p}$. On  suppose donn\'ees des familles $(\phi_{p,q})_{p,q\in \Pi}$, $(\tilde{\phi}_{p,q})_{p,q\in \Pi}$ et $(\nabla_{p,q})_{p,q\in \Pi}$. Pour $p,q\in \Pi$, $\phi_{p,q}:G_{q}\to G_{p}$ et $\tilde{\phi}_{p,q}:\tilde{G}_{q}\to \tilde{G}_{p}$ sont des isomorphismes compatibles d\'efinis sur ${\mathbb C}$  et $\nabla_{p,q}:\Gamma_{{\mathbb R}}\to G_{p,SC}$ est un cocycle. On suppose les hypoth\`eses (1) \`a (5) v\'erifi\'ees pour tous $p,q,r \in \Pi$ et $\sigma\in \Gamma_{{\mathbb R}}$:
 
 (1) $\phi_{p,q}\circ\sigma(\phi_{p,q})^{-1}=ad_{\nabla_{p,q}(\sigma)}$ et 
 $\tilde{\phi}_{p,q}\circ\sigma(\tilde{\phi}_{p,q})^{-1}=ad_{\nabla_{p,q}(\sigma)}$ (ce dernier automorphisme est l'action par conjugaison de $\nabla_{p,q}(\sigma)$ sur $\tilde{G}_{p}$);
 
 (2) $\phi_{p,q}\circ\phi_{q,r}=\phi_{p,r}$ et $\tilde{\phi}_{p,q}\circ\tilde{\phi}_{q,r}=\tilde{\phi}_{p,r}$;
 
 (3) $\nabla_{p,r}(\sigma)=\phi_{p,q}(\nabla_{q,r}(\sigma))\nabla_{p,q}(\sigma)$;
 
 (4) $\tilde{G}_{p}({\mathbb R})\not=\emptyset$.
 
 Pour $x\in \tilde{G}_{p}({{\mathbb R}})$, $ad_{x}$ d\'efinit naturellement un automorphisme de $H^1(\Gamma_{{\mathbb R}},G_{p})$ qui ne d\'epend pas du choix de $x$. Conform\'ement \`a nos conventions, on note cet automorphisme $\theta$.  Alors
 
 (5) la famille $(\nabla_{p,q})_{q\in \Pi}$ s'envoie bijectivement sur $\pi(H^1(\Gamma_{{\mathbb R}};G_{p,SC}))\cap H^1(\Gamma_{{\mathbb R}};G_{p})^{\theta}$.
 
 Dans une telle situation, on d\'efinit le $K$-groupe $KG$ comme la r\'eunion disjointe des $G_{p}$ pour $p\in \Pi$ et le $K$-espace tordu $K\tilde{G}$ comme la r\'eunion disjointe des $\tilde{G}_{p}$. On introduit les sous-ensembles \'evidents $K\tilde{G}_{ss}$ et $K\tilde{G}_{reg}$. Pour $\gamma_{p}\in \tilde{G}_{p}$ et $\gamma_{q}\in \tilde{G}_{q}$, on dit que $\gamma_{p}$ et $\gamma_{q}$ sont conjugu\'es si $\tilde{\phi}_{p,q}(\gamma_{q})$ est conjugu\'e \`a $\gamma_{p}$  dans $\tilde{G}_{p}$.

 {\bf Remarque.} On adopte la terminologie $K$-groupe par commodit\'e. Telle qu'on l'a d\'efinie, cette notion n'est pas intrins\`eque aux groupes puisque la condition (5) d\'epend de l'espace tordu.
 \bigskip
 
 De $\phi_{p,q}$ se d\'eduit une bijection ${\cal E}_{q}\mapsto \phi_{p,q}({\cal E}_{q})$ entre paires de Borel \'epingl\'ees de $G_{q}$ et $G_{p}$. Il s'en d\'eduit une identification ${\cal E}_{q}^*\simeq {\cal E}_{p}^*$ \'equivariante pour les actions galoisiennes. Elle transporte l'automorphisme $\theta^*_{q}$ sur $\theta^*_{p}$. On peut noter simplement ${\cal E}^*$ et $\theta^*$ ces objets. On supposera comme en 1.5 que $\theta^*$ est d'ordre fini. Les groupes $G_{p}$ ont  un $L$-groupe $^LG$ commun et un $L$-espace $^L\tilde{G}$ commun. La donn\'ee d'un ${\bf a}\in H^1(W_{{\mathbb R}};Z(\hat{G}))$ d\'etermine des caract\`eres $\omega_{p}$ de chaque $G_{p}({\mathbb R})$. L'application $\tilde{\phi}_{p,q}$ se restreint en une bijection de $Z(\tilde{G}_{q},{\cal E}_{q})$ sur $Z(\tilde{G}_{p},\phi_{p,q}({\cal E}_{q}))$. Il s'en d\'eduit une bijection ${\cal Z}(\tilde{G}_{q})\simeq {\cal Z}(\tilde{G}_{p})$ elle-aussi \'equivariante pour les actions galoisiennes.
 
 Une donn\'ee endoscopique ${\bf G}'=(G',{\cal G}',\tilde{s})$ pour $(\tilde{G}_{p},{\bf a})$ est aussi une donn\'ee endoscopique pour $(\tilde{G}_{q},{\bf a})$ pour tout $q$.  Changer $(\tilde{G}_{p},{\bf a})$ en $(\tilde{G}_{q},{\bf a})$ ne change pas l'espace endoscopique $\tilde{G}'$. On peut donc consid\'erer ${\bf G}'$ comme une donn\'ee endoscopique pour $(K\tilde{G},{\bf a})$. Pour chaque $p\in \Pi$, notons plus pr\'ecis\'ement ${\cal D}_{\tilde{G}_{p}}({\bf G}')$ l'ensemble d\'efini en 1.8 quand on consid\`ere ${\bf G}'$ comme une donn\'ee endoscopique de $(\tilde{G}_{p},{\bf a})$. On pose ${\cal D}_{K\tilde{G}}({\bf G}')=\sqcup_{p\in {\Pi}}{\cal D}_{\tilde{G}_{p}}({\bf G}')$.
 
 Montrons qu'\`a partir d'un couple $(G,\tilde{G})$ v\'erifiant les conditions de 1.5, on peut construire un $K$-espace comme ci-dessus. On fixe un ensemble $\Pi$ de cocycles $p:\Gamma_{{\mathbb R}}\to G_{SC}$ qui s'envoie bijectivement sur $\pi(H^1(\Gamma_{{\mathbb R}},G_{SC}))\cap H^1(\Gamma_{{\mathbb R}},G)^{\theta}$. Pour $p\in {\Pi}$, fixons un groupe $G_{p}$ et un espace tordu $\tilde{G}_{p}$ sous ce groupe, tous deux d\'efinis sur ${\mathbb R}$, munis d'isomorphismes compatibles $\phi_{p}:G_{p}\to G$ et $\tilde{\phi}_{p}:\tilde{G}_{p}\to \tilde{G}$, d\'efinis sur ${\mathbb C}$, de sorte que, pour tout $\sigma\in \Gamma_{{\mathbb R}}$, on ait les \'egalit\'es $\phi_{p}\circ\sigma(\phi_{p})^{-1}=ad_{p(\sigma)}$ et  $\tilde{\phi}_{p}\circ\sigma(\tilde{\phi}_{p})^{-1}=ad_{p(\sigma)}$. De tels objets existent: il suffit de poser $G_{p}=G$, $\tilde{G}_{p}=\tilde{G}$, de  prendre pour $\phi_{p}$ et $\tilde{\phi}_{p}$ les identit\'es et de d\'efinir les actions galoisiennes sur $G_{p}$ et $\tilde{G}_{p}$ par les \'egalit\'es pr\'ec\'edentes. Pour $p,q\in \Pi$ et $\sigma\in G_{{\mathbb R}}$, on d\'efinit $\phi_{p,q}=\phi_{p}^{-1}\circ \phi_{q}$ et $\nabla_{p,q}(\sigma)=\phi_{p}^{-1}(q(\sigma)p(\sigma)^{-1})$. La v\'erification des propri\'et\'es (1) \`a (5) est routini\`ere. Indiquons simplement la peuve de (4), qui justifie la condition d'invariance par $\theta$ impos\'ee aux cocycles. Fixons $\gamma\in \tilde{G}({\mathbb R})$. L'image de $p$ dans $H^1(\Gamma_{{\mathbb R}},G)$ est invariante par $ad_{\gamma}$. On peut donc fixer $g\in G$ tel que $ad_{\gamma}(p(\sigma))=g^{-1}p(\sigma)\sigma(g)$ pour tout $\sigma$. Cela implique 
$$\sigma(g\gamma)=\sigma(g)\gamma=p(\sigma)^{-1}g\gamma p(\sigma)=ad_{p(\sigma)^{-1}}(g\gamma).$$
Posons $\gamma_{p}=\phi_{p}^{-1}(g\gamma)$. Alors 
$$\sigma(\gamma_{p})=\sigma(\phi_{p})^{-1}(\sigma(g\gamma))=\sigma(\phi_{p})^{-1}\circ ad_{p(\sigma)^{-1}}(g\gamma)=\phi_{p}^{-1}(g\gamma)=\gamma_{p}. $$
Donc $\gamma_{p}\in \tilde{G}_{p}({\mathbb R})$.

Inversement, si on part de donn\'ees comme ci-dessus et si on fixe un $p_{0}\in \Pi$, on peut identifier $K\tilde{G}$ \`a un $K$-espace tordu d\'efini comme on vient de le faire \`a partir du couple $(G,\tilde{G})=(G_{p_{0}},\tilde{G}_{p_{0}})$.

\bigskip
\subsection{L'ensemble $\tilde{G}_{ab}(F)$}
Le corps $F$ est de nouveau un corps local de caract\'eristique nulle.   Soit $A$ un groupe  et $B$ un ensemble muni d'une action {\bf \`a droite} de $A$. On suppose $A$ et $B$ munis d'actions de $\Gamma_{F}$ compatibles  \`a cette action. Notons $Z^{1,0}(\Gamma_{F};A\circlearrowleft B)$ l'ensemble des couples $(\alpha,b)$ o\`u $b\in B$ et $\alpha:\Gamma_{F}\to A$ est un cocycle tels que $\sigma(b)=b\alpha(\sigma)$ pour tout $\sigma\in \Gamma_{F}$. On introduit la relation d'\'equivalence $(\alpha,b)\equiv (\alpha',b')$ si et seulement s'il existe $a\in A$ tel que $\alpha'(\sigma)=a^{-1}\alpha(\sigma)\sigma(a)$ et $b'=ba$. On note $H^{1,0}(\Gamma_{F};A\circlearrowleft B) $ le quotient de $Z^{1,0}(\Gamma_{F};A\circlearrowleft B)$ par cette relation d'\'equivalence.

Il y a un cas particulier important de la construction pr\'ec\'edente. Consid\'erons deux groupes $A$ et $B$ munis d'actions de $\Gamma_{F}$ et un homomorphisme de groupes $f:A\to B$ \'equivariant pour cette action. On peut consid\'erer que $A$ agit sur $B$ par $(a,b)\mapsto bf(a)$. On note alors $H^{1,0}(\Gamma_{F};A\stackrel{f}{\to }B)$ l'ensemble $H^{1,0}(\Gamma_{F};A\circlearrowleft B)$ pr\'ec\'edent. Si $A$ et $B$ sont ab\'eliens, c'est  aussi un groupe ab\'elien.

{\bf Remarque.} Ces ensembles ont \'et\'e d\'efinis par divers auteurs. F\^acheusement, les uns les notent $H^0$, les autres $H^1$ et les d\'efinitions varient par des signes. Nous avons adopt\'e la notation $H^{1,0}$ qui est lourde mais a l'avantage de m\'econtenter tout le monde. Labesse utilise la notation $H^0$ et sa d\'efinition diff\`ere de la n\^otre car il consid\`ere une action \`a gauche de $A$ sur $B$. Kottwitz et Shelstad  ne consid\`erent que des groupes ab\'eliens et utilisent la notation $H^1$. A cette diff\'erence de notation pr\`es, notre d\'efinition est la m\^eme que la leur. Signalons que, sous certaines hypoth\`eses topologiques suppl\'ementaires, on peut d\'efinir comme ci-dessus des ensembles $H^{1,0}(W_{F};A\stackrel{f}{\to}B)$, cf.  [KS1] A.3.
 
\bigskip

 Ainsi, on d\'efinit l'ensemble $G_{ab}(F)=H^{1,0}(\Gamma_{F};G_{SC}\stackrel{\pi}{\to}G)$ (pour nous, $G_{SC}$ agit \`a droite sur $G$), cf. [Lab1] 1.6. L'application naturelle  de $H^{1,0}(\Gamma_{F};Z(G_{SC})\stackrel{\pi}{\to}Z(G))$ dans cet ensemble $G_{ab}(F)$ est bijective, ce qui munit $G_{ab}(F)$ d'une structure de groupe. Il y a un homomorphisme naturel injectif
$$G(F)/\pi(G_{SC}(F))\to G_{ab}(F),$$
qui est surjectif si $F\not={\mathbb R}$.

 Ainsi, on d\'efinit l'ensemble $H^{1,0}(\Gamma_F;G_{SC}\circlearrowleft\tilde{G})$, que l'on peut noter $\tilde{G}_{ab}(F)$. On a une application:
  $$\begin{array}{ccc} Z^{1,0}(\Gamma_{F};G_{SC}\circlearrowleft\tilde{G})\times Z^{1,0}(\Gamma_{F};Z(G_{SC})\stackrel{\pi}{\to}Z(G)) &\to&Z^{1,0}(\Gamma_{F};G_{SC}\circlearrowleft\tilde{G})\\ ((\mu,\gamma),(\zeta,z))&\mapsto&(\mu\zeta,\gamma z).\\ \end{array}$$
  Elle se quotiente en une action \`a droite du groupe $G_{ab}(F)\simeq H^{1,0}(\Gamma_{F};Z(G_{SC})\stackrel{\pi}{\to}Z(G))$ sur $\tilde{G}_{ab}(F)$. On a:
  
  (1) $\tilde{G}_{ab}(F)$ est un espace principal homog\`ene sous $G_{ab}(F)$.
  
  Preuve. Soient $(\zeta,z)$, $(\zeta',z')$ deux \'el\'ements de $Z^{1,0}(\Gamma_{F};Z(G_{SC})\stackrel{\pi}{\to}Z(G))$ et soit $(\mu,\gamma)\in Z^{1,0}(\Gamma_{F};G_{SC}\circlearrowleft\tilde{G})$. Supposons $(\mu\zeta,\gamma z)$ cohomologue \`a $(\mu\zeta',\gamma z')$. Alors il existe $x\in G_{SC}$ tel que $\mu(\sigma)\zeta'(\sigma)=x^{-1}\mu(\sigma)\zeta(\sigma)\sigma(x)$ et $\gamma z'=\gamma z \pi(x)$. Cette derni\`ere relation implique que $z'=z\pi(x)$ et que $x$ appartient \`a $Z(G_{SC})$. La premi\`ere relation implique alors que $\zeta'(\sigma)=x^{-1}\zeta(\sigma)\sigma(x)$, donc les couples $(\zeta,z)$ et $(\zeta',z')$ sont cohomologues. Cela prouve que l'action de $G_{ab}(F)$ sur $\tilde{G}_{ab}(F)$ est libre. Soient maintenant $(\mu,\gamma)$ et $(\mu',\gamma')$ deux \'el\'ements de $Z^{1,0}(\Gamma_{F};G_{SC}\circlearrowleft\tilde{G})$. Soit $g\in G$ l'\'el\'ement tel que $\gamma'=\gamma g$, \'ecrivons $g=\pi(x)z$ avec $x\in G_{SC}$ et $z\in Z$. Le couple $(\mu',\gamma')$ est cohomologue \`a $(\mu'',\gamma z)$, o\`u $\mu''(\sigma)=x\mu'(\sigma)\sigma(x)^{-1}$. Posons $\zeta(\sigma)=\mu(\sigma)^{-1}\mu''(\sigma) $. Les \'egalit\'es $\sigma(\gamma)=\gamma\pi(\mu(\sigma))$ et $\sigma(\gamma z)=\gamma z\pi(\mu''(\sigma))$ entra\^{\i}nent que $\sigma(z)=z\pi(\zeta(\sigma))$. Cela implique que $\zeta(\sigma)$ appartient \`a $Z(G_{SC})$. Cette propri\'et\'e et le fait que $\mu$ et $\mu''$ sont des cocycles implique que $\zeta$ est aussi un cocycle. Alors $(\zeta,z)$ appartient \`a $Z^{1,0}(\Gamma_{F};Z(G_{SC})\stackrel{\pi}{\to}Z(G))$. Le couple $(\mu',\gamma')$ est cohomologue au produit  de $(\mu,\gamma)$ et de $(\zeta,z)$. Cela prouve que l'action de $G_{ab}(F)$ sur $\tilde{G}_{ab}(F)$ est transitive. $\square$
  
  Remarquons que l'on pourrait aussi bien d\'efinir une action \`a gauche de $G_{ab}(F)$ sur $\tilde{G}_{ab}(F)$, jouissant des m\^emes propri\'et\'es.
  
  Il y a une application naturelle $\tilde{G}(F)\to \tilde{G}_{ab}(F)$: \`a $\gamma\in \tilde{G}(F)$, on associe l'image dans $\tilde{G}_{ab}(F)$ de $(\mu=1,\gamma)\in Z^{1,0}(\Gamma_{F};G_{SC}\circlearrowleft\tilde{G})$.
  
  On va d\'efinir une application
  $$(2) \qquad \tilde{G}_{ab}(F)\to H^{1,0}(\Gamma_{F};{\cal Z}(G_{SC})\circlearrowleft {\cal Z}(\tilde{G})).$$ 
  
  Soit $(\mu,\gamma)\in Z^{1,0}(\Gamma_{F};G_{SC}\circlearrowleft\tilde{G})$. Fixons une paire de Borel \'epingl\'ee ${\cal E}$ et une cocha\^{\i}ne $u_{{\cal E}}$ comme en 1.2.  On peut choisir, et on choisit, $x\in G_{SC}$ et $e\in Z(\tilde{G},{\cal E})$ tels que $\gamma=e\pi(x)$. Posons  $\mu'(\sigma)=x\mu(\sigma)\sigma(x)^{-1}$, puis $\nu(\sigma)=ad_{e}^{-1}(u_{{\cal E}}(\sigma))\mu'(\sigma)u_{{\cal E}}(\sigma)^{-1} $. L'\'egalit\'e $\sigma(\gamma)=\gamma\pi(\mu(\sigma))$ entra\^{\i}ne $\sigma(e)=e\pi(\mu'(\sigma))$, puis 
  
  (3) $ad_{u_{{\cal E}}(\sigma)}(\sigma(e))= e\pi(\nu(\sigma)) $. 
  
  Or $ad_{u_{{\cal E}}(\sigma)}\circ\sigma$ conserve ${\cal E}$, donc aussi $Z(\tilde{G},{\cal E})=eZ(G)$. Donc $ad_{u_{{\cal E}}(\sigma)}(\sigma(e))\in eZ(G)$ et l'\'egalit\'e (3) implique que $\nu(\sigma)$ appartient \`a $Z(G_{SC})$. Rappelons que le cobord $du_{{\cal E}}$ prend ses valeurs dans $Z(G_{SC})$. Montrons que
  
  (4) $d\nu=(\theta^{-1}-1)(du_{{\cal E}})$.
  
  Pour cela, d\'efinissons un espace tordu $\tilde{G}_{\star}$ sur le groupe $G_{SC}$ de la fa\c{c}on suivante. Il est \'egal \`a $e_{\star}G_{SC} $, o\`u $e_{\star}$ est un point fix\'e.  L'action de $G_{SC}$ \`a droite est l'action naturelle, celle \`a gauche est d\'efinie par $ge_{\star}=e_{\star}ad_{e}^{-1}(g)$. La structure galoisienne est $(\sigma,e_{\star}g)\mapsto e_{\star}\mu'(\sigma)\sigma(g)$. On v\'erifie que cette d\'efinition est loisible. On a la relation analogue \`a (3):
  $$(5) \qquad ad_{u_{{\cal E}}(\sigma)}(\sigma(e_{\star}))=e_{\star}\nu(\sigma).$$
 Soient $\sigma_{1},\sigma_{2}\in \Gamma_{F}$.  En  rempla\c{c}ant dans (5) $\sigma$ par $\sigma_{1}$ et en multipliant \`a droite l'\'egalit\'e obtenue par $\sigma_{1}(\nu(\sigma_{2}))$, on obtient
 $$ad_{u_{{\cal E}}(\sigma_{1})}(\sigma_{1}(e_{\star}\nu(\sigma_{2})))=e_{\star}\nu(\sigma_{1})\sigma_{1}(\nu(\sigma_{2})),$$
 puisque $\nu(\sigma_{2})$ est central. On remplace le terme $e_{\star}\nu(\sigma_{2})$ du membre de gauche par sa valeur donn\'ee par (5) et on obtient
 $$ad_{u_{{\cal E}}(\sigma_{1})\sigma_{1}(u_{{\cal E}}(\sigma_{2}))}(\sigma_{1}\sigma_{2}(e_{\star}))=e_{\star}\nu(\sigma_{1})\sigma_{1}(\nu(\sigma_{2})),$$
 ou encore
 $$ad_{du_{{\cal E}}(\sigma_{1},\sigma_{2})u_{{\cal E}}(\sigma_{1}\sigma_{2})}(\sigma_{1}\sigma_{2}(e_{\star}))=e_{\star}\nu(\sigma_{1},\sigma_{2})d\nu(\sigma_{1},\sigma_{2}).$$
 On exprime le membre de gauche gr\^ace \`a l'\'egalit\'e (5) pour $\sigma=\sigma_{1}\sigma_{2}$. On obtient
 $$ad_{du_{{\cal E}}(\sigma_{1},\sigma_{2})}(e_{\star}\nu(\sigma_{1}\sigma_{2}))=e_{\star}\nu(\sigma_{1},\sigma_{2})d\nu(\sigma_{1},\sigma_{2}).$$
 Cela entra\^{\i}ne la relation (4).

Notons $z\mapsto \bar{z}$ les applications naturelles de $Z(G_{SC})$ dans ${\cal Z}(G_{SC})$ ou de $Z(\tilde{G},{\cal E})$ dans ${\cal Z}(\tilde{G})$.   La relation (4) entra\^{\i}ne que $\bar{\nu}$ est un cocycle. La relation (3) et la d\'efinition de l'action galoisienne sur ${\cal Z}(\tilde{G})$ entra\^{\i}nent que $\sigma(\bar{e})=\bar{e}\pi(\bar{\nu}(\sigma))$. Donc $(\bar{\nu},\bar{e})$ appartient \`a $Z^{1,0}(\Gamma_{F};{\cal Z}(G_{SC})\circlearrowleft {\cal Z}(\tilde{G}))$. Montrons que

(6) la classe de cohomologie de $(\bar{\nu},\bar{e})$ ne d\'epend pas des choix effectu\'es et ne d\'epend que de la classe de cohomologie de $(\mu,\gamma)$.

On a choisi ${\cal E}$,  $u_{{\cal E}}$, $x$ et $e$. L'ind\'ependance de $u_{{\cal E}}$ est claire: on ne peut modifier $u_{{\cal E}}(\sigma)$ que par un \'el\'ement de $Z(G_{SC})$, ce qui ne change pas l'image $\bar{\nu}(\sigma)$ dans ${\cal Z}(G_{SC})$.  Supposons d'abord ${\cal E}$ et $(\mu,\gamma)$ fix\'es. On ne peut modifier $x$ et $e$ qu'en rempla\c{c}ant $x$ par $z^{-1}x$ et $e$ par $e\pi(z)$ pour un \'el\'ement $z\in Z(G_{SC})$. On voit que cela remplace $ \bar{\nu}(\sigma)$ par $\bar{\nu}_{1}(\sigma)=\bar{z}^{-1}\bar{\nu}(\sigma)\sigma(\bar{z})$ et $\bar{e}$ par $\bar{e}_{1}=\bar{e}\bar{z}$. Or $(\bar{\nu}_{1},\bar{e}_{1})$ est cohomologue \`a $(\bar{\nu},\bar{e})$. Supposons maintenant ${\cal E}$ fix\'e et rempla\c{c}ons $(\mu,\gamma)$ par $(\mu_{1},\gamma_{1})$ cohomologue \`a $(\mu,\gamma)$. Soit $v\in G_{SC}$ tel que $\mu_{1}(\sigma)=v^{-1}\mu(\sigma)\sigma(v)$ et $\gamma_{1}=\gamma v$. Pour le couple $(\mu_{1},\gamma_{1})$, on peut choisir $e_{1}=e$ et $x_{1}=xv$. Alors $\mu'_{1}=\mu'$ et le couple $(\bar{\nu},\bar{e})$ ne change pas. Il reste \`a remplacer ${\cal E}$   par une autre paire de Borel \'epingl\'ee ${\cal E}_{1}$, $(\mu,\gamma)$ \'etant fix\'e.  On fixe  $r\in G_{SC}$ tel que  $ad_{r}({\cal E})={\cal E}_{1}$. On peut choisir $u_{{\cal E}_{1}}(\sigma)=ru_{{\cal E}}(\sigma)\sigma(r)^{-1}$,  $e_{1}=ad_{r}(e)=e\pi(s)$, o\`u $s=ad_{e}^{-1}(r)r^{-1}$, et $x_{1}=s^{-1}x$. On a $\bar{e}_{1}=\bar{e}$ par d\'efinition de l'ensemble ${\cal Z}(\tilde{G})$. On a $\mu'_{1}(\sigma)=s^{-1}\mu'(\sigma)\sigma(s)$, puis 
$$\nu_{1}(\sigma)=ad_{e_{1}}^{-1}(u_{{\cal E}_{1}}(\sigma))\mu'_{1}(\sigma)u_{{\cal E}_{1}}(\sigma)^{-1}$$
$$=ad_{r}\circ ad_{e}^{-1}\circ ad_{r}^{-1}(ru_{{\cal E}}(\sigma)\sigma(r)^{-1})s^{-1}\mu'(\sigma)\sigma(s)\sigma(r)u_{{\cal E}}(\sigma)^{-1}r^{-1}$$
$$=r ad_{e}^{-1}(u_{{\cal E}}(\sigma)\sigma(r)^{-1})\mu'(\sigma)\sigma(ad^{-1}_{e}(r))u_{{\cal E}}(\sigma)^{-1}r^{-1}=ra\nu(\sigma)br^{-1},$$
o\`u $a=ad_{e}^{-1}(u_{{\cal E}}(\sigma)\sigma(r)^{-1}u_{{\cal E}}(\sigma)^{-1})$ et $b=u_{{\cal E}}(\sigma)\sigma(ad_{e}^{-1}(r))u_{{\cal E}}(\sigma)^{-1}$. Puisqu'on sait que $\nu_{1}(\sigma)$ est central, on peut aussi bien conjuguer par $ra$ et on obtient $\nu_{1}(\sigma)=\nu(\sigma)ba$. Introduisons l'action $\sigma\mapsto \sigma_{G^*}$ de $\Gamma_{F}$ sur $G$ d\'efinie par $\sigma_{G^*}=ad_{u_{{\cal E}}(\sigma)}\circ \sigma_{G}$. Le fait que $ad_{u_{{\cal E}}(\sigma)}(\sigma(e))\in Z(G)e$ entra\^{\i}ne que $ad_{e}$ commute \`a cette action. Or $a=ad_{e}^{-1}\circ \sigma_{G^*}(r)^{-1}$ et $b=\sigma_{G^*}\circ ad_{e}^{-1}(r)$. Donc $a=b^{-1}$ et $\nu_{1}(\sigma)=\nu(\sigma)$. Cela prouve (6).

D'apr\`es (6), on a d\'efini l'application cherch\'ee
$$ \tilde{G}_{ab}(F)\to H^{1,0}(\Gamma_{F}; {\cal Z}(G_{SC})\circlearrowleft {\cal Z}(\tilde{G})).$$
Il est facile de voir comme en (1) que l'ensemble d'arriv\'ee est un espace principal homog\`ene sous $H^{1,0}(\Gamma_{F};{\cal Z}(G_{SC})\stackrel{\pi}{\to}{\cal Z}(G))$. 

{\bf Cas particulier.} Dans le cas o\`u $\tilde{G}$ est \`a torsion int\'erieure, ce dernier ensemble n'est autre que $G_{ab}(F)$. La fl\`eche (2) \'etant bien s\^ur \'equivariante pour les actions de $G_{ab}(F)$ et les ensembles de d\'epart et d'arriv\'ee \'etant tous deux des espaces principaux homog\`enes sous ce groupe, la fl\`eche est bijective.

\bigskip

 Le groupe $Z(G)$ est naturellement un sous-groupe de $T^*$. On pose ${\cal Z}_{0}(G)=Z(G)/(Z(G)\cap(1-\theta^*)(T^*))$.   Il y a un homomorphisme surjectif ${\cal Z}(G)\to {\cal Z}_{0}(G)$. On pose ${\cal Z}_{0}(\tilde{G})={\cal Z}_{0}(G)\times_{{\cal Z}(G)}{\cal Z}(\tilde{G})$, la notation ayant le m\^eme sens qu'en 1.7. L'application (2) se pousse en une application que nous notons
$$ N^{\tilde{G}}:\tilde{G}_{ab}(F)\to H^{1,0}(\Gamma_{F}; {\cal Z}_{0}(G_{SC})\circlearrowleft {\cal Z}_{0}(\tilde{G})).$$

Soit ${\bf G}'=(G ',{\cal G}',\tilde{s})$ une donn\'ee endoscopique pour $(G,\tilde{G},{\bf a})$. Rappelons que l'on a un homomorphisme $Z(G)\to Z(G')$. Il se factorise en une suite
$$Z(G)\to {\cal Z}_{0}(G)\stackrel{\xi_{0}}{\to}Z(G')$$
et $\xi_{0}$ est injectif. On a de m\^eme une suite
$${\cal Z}(\tilde{G})\to {\cal Z}_{0}(\tilde{G})\stackrel{\tilde{\xi}_{0}}{\to }{\cal Z}(\tilde{G}'),$$
et $\tilde{\xi}_{0}$ est injectif. On a une suite d'extensions
$$\hat{G}'\to \hat{G}'_{ad}=\hat{G}'/(\hat{G}'\cap Z(\hat{G}))\to \hat{G}'_{AD}=\hat{G}'/Z(\hat{G}'),$$
dont on d\'eduit une suite duale
$$G'\leftarrow G'_{sc}\leftarrow G'_{SC}.$$
Il y a donc une application naturelle
$$(7)\qquad H^{1,0}(\Gamma_{F};Z(G'_{SC})\circlearrowleft {\cal Z}(\tilde{G}'))\to H^{1,0}(\Gamma_{F}; Z(G'_{sc})\circlearrowleft {\cal Z}(\tilde{G}')).$$
Un tore maximal de $\hat{G}'_{ad}$ est naturellement isomorphe \`a ${\hat{T}}^{\hat{\theta},0}/(Z(\hat{G})\cap {\hat{T}}^{\hat{\theta},0})$, qui n'est autre que ${\hat{T}}_{ad}^{\hat{\theta}}$, o\`u ${\hat{T}}_{ad}$ est l'image de ${\hat{T}}$ dans $\hat{G}_{AD}$ (on rappelle que ${\hat{T}}_{ad}^{\hat{\theta}}$ est connexe). Dualement, un tore maximal de $G'_{sc}$ est donc isomorphe \`a $T_{sc}^*/(1-\theta^*)(T_{sc}^*)$.   On en d\'eduit une suite analogue \`a celle ci-dessus:
$$Z(G_{SC})\to {\cal Z}_{0}(G_{SC})\stackrel{\xi_{0,sc}}{\to}Z(G'_{sc}),$$
o\`u $\xi_{0,sc}$ est injectif. D'o\`u une application naturelle
$$(8)\qquad H^{1,0}(\Gamma_{F}; {\cal Z}_{0}(G_{SC})\circlearrowleft {\cal Z}_{0}(\tilde{G}))\to H^{1,0}(\Gamma_{F}; Z(G'_{sc})\circlearrowleft {\cal Z}(\tilde{G}')).$$
Montrons qu'elle est bijective. Consid\'erons le diagramme
$$\begin{array}{ccc} {\cal Z}_{0}(G_{SC})&\to&{\cal Z}_{0}(G)\\ \xi_{0,sc}\,\downarrow\,\,&&\xi_{0}\,\downarrow\,\,\\ Z(G'_{sc})&\stackrel{\pi'}{\to}&Z(G')\\ \end{array}$$
Alors

(9) $Z(G')$ est engendr\'e par  les images de $\pi'$ et de $\xi_{0}$;

(10) l'image r\'eciproque par $\pi'$ de l'image de $\xi_{0}$ est l'image de $\xi_{0,sc}$.

Le tore $T^*$ est engendr\'e par $Z(G)$ et par l'image de $T^*_{sc}$ et (9) en r\'esulte. Soit $x\in Z(G'_{sc})$ tel que $\pi'(x)$ appartient \`a l'image de $\xi_{0}$. Choisissons un \'el\'ement $t_{sc}\in T^*_{sc}$ dont $x$ soit l'image dans $T^*_{sc}/(1-\theta^*)(T^*_{sc})$. L'hypoth\`ese signifie que $\pi(t_{sc})\in Z(G)(1-\theta^*)(T^*)$. Ecrivons $\pi(t_{sc})=z(1-\theta^*)(t)$, avec $z\in Z(G)$ et $t\in T^*$. Ecrivons $t=z'\pi(t'_{sc})$, avec $z'\in Z(G)$ et $t'_{sc}\in T^*_{sc}$. Alors $\pi(t_{sc}(\theta^*-1)(t'_{sc}))=z(1-\theta^*)(z')$. Cela entra\^{\i}ne que $t_{sc}(\theta^*-1)(t'_{sc})$ appartient \`a $Z(G_{SC})$. Puisque $t_{sc}(\theta^*-1)(t'_{sc})$ a aussi pour image $x$ dans $T^*_{sc}/(1-\theta^*)(T^*_{sc})$, cela montre que $x$ appartient \`a l'image de $Z(G_{SC})$, qui n'est autre que celle de l'application $\xi_{0,sc}$. Cela prouve (10).

 Soit $(\zeta',e')\in Z^{1,0}( \Gamma_{F}; Z(G'_{sc})\circlearrowleft {\cal Z}(\tilde{G}'))$.  La relation (9) entra\^{\i}ne que l'on peut \'ecrire $e'=\tilde{\xi}_{0}(e)\pi'(z'_{sc})$, avec $z'_{sc}\in Z(G'_{sc})$ et $e\in {\cal Z}_{0}(\tilde{G})$. Alors $(\zeta',e')$ est cohomologue \`a $(\zeta'_{1},\tilde{\xi}_{0}(e))$, o\`u $\zeta'_{1}(\sigma)=z_{sc}'\zeta'(\sigma)\sigma(z'_{sc})^{-1}$. La relation $\sigma\circ\tilde{\xi}_{0}(e)=\tilde{\xi}_{0}(e)\pi'(\zeta'_{1}(\sigma))$ entra\^{\i}ne que $\pi'\circ\zeta'_{1}$ prend ses valeurs dans l'image de $\xi_{0}$. D'apr\`es (10), on peut \'ecrire $\zeta'_{1}=\xi_{0,sc}(\zeta)$, o\`u $\zeta$ est \`a valeurs dans ${\cal Z}_{0}(G_{SC})$. Puisque $\xi_{0,sc}$ et $\tilde{\xi}_{0}$ sont injectifs, le couple $(\zeta,e)$ v\'erifie les conditions requises pour appartenir \`a $Z^{1,0}(\Gamma_{F}; {\cal Z}_{0}(G_{SC})\circlearrowleft {\cal Z}_{0}(\tilde{G}))$. La classe de cohomologie de $(\zeta',e')$ est l'image par l'application (8) de celle de $(\zeta,e)$. Cela prouve la surjectivit\'e de (8). Inversement, soient $(\zeta_{1},e_{1})$ et $(\zeta_{2},e_{2})$ deux \'el\'ements de $Z^{1,0}(\Gamma_{F}; {\cal Z}_{0}(G_{SC})\circlearrowleft {\cal Z}_{0}(\tilde{G}))$ qui ont m\^eme image dans 
 $H^{1,0}(\Gamma_{F}; Z(G'_{sc})\circlearrowleft {\cal Z}(\tilde{G}'))$. Il existe $z'_{sc}\in Z(G'_{sc})$ tel que $\xi_{0,sc}(\zeta_{1}(\sigma))=\xi_{0,sc}(\zeta_{2}(\sigma))(z'_{sc})^{-1}\sigma(z'_{sc})$ et $\tilde{\xi}_{0}(e_{1})=\tilde{\xi}_{0}(e_{2})\pi'(z'_{sc})$. Cette deuxi\`eme relation entra\^{\i}ne que $\pi'(z'_{sc})$ appartient \`a l'image de $\xi_{0}$. D'apr\`es (10), il existe $z_{sc}\in {\cal Z}_{0}(G_{SC})$ tel que $z'_{sc}=\xi_{0,sc}(z_{sc})$. D'apr\`es l'injectivit\'e de $\xi_{0,sc}$ et $\tilde{\xi}_{0}$, on a alors $\zeta_{1}(\sigma)=\zeta_{2}(\sigma)(z_{sc})^{-1}\sigma(z_{sc})$ et $e_{1}=e_{2}\pi(z_{sc})$. Donc les couples $(\zeta_{1},e_{1})$ et $(\zeta_{2},e_{2})$ sont cohomologues, ce qui prouve l'injectivit\'e de (8).
 \bigskip

L'ensemble de d\'epart de (8) n'est autre que $\tilde{G}'_{ab}(F)$, puisque $\tilde{G}'$ est \`a torsion int\'erieure. Par composition de (7) et de l'inverse de (8),  on obtient une application que nous notons
$$N^{\tilde{G}',\tilde{G}}:\tilde{G}'_{ab}(F)\to H^{1,0}(\Gamma_{F}; {\cal Z}_{0}(G_{SC})\circlearrowleft {\cal Z}_{0}(\tilde{G})).$$

{\bf Remarque.} On note aussi $N^{\tilde{G}',\tilde{G}}$ la compos\'ee de cette application avec l'application $\tilde{G}'(F)\to \tilde{G}'_{ab}(F)$.

Il est plus parlant d'identifier l'ensemble d'arriv\'ee de cette application. Introduisons le groupe $G_{0}$ quasi-d\'eploy\'e sur $F$ dual du groupe $\hat{G}_{0}=\hat{G}^{\hat{\theta},0}$, muni de l'action galoisienne provenant de celle sur $\hat{G}$. Notons ${\cal G}'_{0}$ le sous-groupe $\hat{G}_{0}\rtimes W_{F}$ de $^LG$. Le cocycle ${\bf a}$ ne joue ici aucun r\^ole. On peut  remplacer ${\bf a}$ par le caract\`ere trivial ${\bf 1}$. Alors le triplet ${\bf G}_{0}=(G_{0},{\cal G}'_{0},\hat{\theta})$ est une donn\'ee endoscopique  pour $(G,\tilde{G},{\bf 1})$ \`a laquelle on applique les constructions ci-dessus. Pour cette donn\'ee, on a $Z(\hat{G}_{0})=Z(\hat{G})\cap \hat{T}^{\hat{\theta},0}$. Cela r\'esulte du fait que les racines simples pour la paire de Borel $({\hat{B}}\cap \hat{G}_{0},{\hat{T}}\cap\hat{G}_{0}={\hat{T}}^{\hat{\theta},0})$  de $\hat{G}_{0}$ sont exactement les restrictions \`a ${\hat{T}}^{\hat{\theta},0}$ des racines simples pour la paire de Borel $({\hat{B}},{\hat{T}})$ de $\hat{G}$, cf. 1.6. Il en r\'esulte que $\hat{G}_{0,ad}=\hat{G}_{0,AD}$, puis $G_{0,sc}=G_{0,SC}$. Donc, pour cette donn\'ee ${\bf G}_{0}$, l'application (7) est l'identit\'e. Donc  l'application $N^{\tilde{G}_{0},\tilde{G}}$ est bijective, ce qui nous permet d'identifier $H^{1,0}(\Gamma_{F}; {\cal Z}_{0}(G_{SC})\circlearrowleft {\cal Z}_{0}(\tilde{G}))$ \`a $\tilde{G}_{0,ab}(F)$.

Revenons \`a notre donn\'ee ${\bf G}'$. On a construit des applications
$$(11)\left\lbrace\begin{array}{ccccc}\tilde{G}(F)&\to&\tilde{G}_{ab}(F)&&\\&&&\searrow\,N^{\tilde{G}}&\\&&&&\tilde{G}_{0,ab}(F)\\&&&\nearrow\,N^{\tilde{G}',\tilde{G}}&\\\tilde{G}'(F)&\to&\tilde{G}'_{ab}(F)&&\\ \end{array}\right.$$
Les termes extr\^emes sont des espaces principaux homog\`enes sous respectivement $G_{ab}(F)$, $G'_{ab}(F)$ et $G_{0,ab}(F)$. Il est clair qu'il y a des homomorphismes similaires
$$G_{ab}(F)\stackrel{N^G}{\to} G_{0,ab}(F)\stackrel{N^{G',G}}{\leftarrow} G'_{ab}(F)$$
compatibles avec les applications ci-dessus.

Supposons un instant que $F={\mathbb R}$. On a introduit en 1.11 un $K$-espace $K\tilde{G}$. On d\'efinit $K\tilde{G}_{ab}({\mathbb R})$ comme la r\'eunion disjointe des $\tilde{G}_{p,ab}({\mathbb R})$ pour $p\in{\Pi}$ et on obtient un diagramme similaire au pr\'ec\'edent o\`u $\tilde{G}({\mathbb R})$ et $\tilde{G}_{ab}({\mathbb R})$ sont remplac\'es par  $K\tilde{G}({\mathbb R})$ et $K\tilde{G}_{ab}({\mathbb R})$.

\bigskip
 \subsection{Caract\`eres de $G(F)$, $G_{0,ab}(F)$, $G_{0,ab}(F)/N^G(G_{ab}(F))$}
 Comme on l'a dit dans le paragraphe pr\'ec\'edent, on a l'\'egalit\'e
 $$G_{ab}(F)= H^{1,0}(\Gamma_{F};Z(G_{SC})\to Z(G)).$$
  Fixons un tore maximal $T$ de $G$ d\'efini sur $F$. On  introduit le tore dual $\hat{T}$ muni de l'action galoisienne  duale de celle de $T$.  L'homomorphisme naturel
  $$H^{1,0}(\Gamma_{F};Z(G_{SC})\to Z(G))\to H^{1,0}(\Gamma_{F};T_{sc}\to T)$$
  est bijectif. D'apr\`es [KS1] lemme A.3.B,  le groupe de caract\`eres continus du dernier groupe est le quotient de $H^{1,0}(W_{F};\hat{T}\to \hat{T}_{ad})$ par l'image naturelle de $\hat{T}_{ad}^{\Gamma_{F},0}$. On v\'erifie que cette image est nulle et que l'homomorphisme naturel
  $$H^1(W_{F};Z(\hat{G}))\to H^{1,0}(W_{F};\hat{T}\to \hat{T}_{ad})$$ 
  est bijectif. On en d\'eduit que le groupe des caract\`eres continus de $G_{ab}(F)$ est isomorphe \`a $H^1(W_{F};Z(\hat{G}))$.
  
  Cela nous permet de pr\'eciser la correspondance qui, \`a ${\bf a}\in H^1(W_{F};Z(\hat{G}))$, associe le caract\`ere $\omega$ de $G(F)$. On a un homomorphisme
  $$G(F)\to G_{ab}(F)= H^{1,0}(\Gamma_{F};Z(G_{SC})\to Z(G)).$$
  Concr\`etement, pour $g\in G(F)$, on \'ecrit $g=\pi(g_{sc})z$, avec $g_{sc}\in G_{SC}$ et $z\in Z(G)$. L'image de $g$ par l'application ci-dessus est repr\'esent\'ee par le couple $(\mu,z)$, o\`u $\mu(\sigma)=g_{sc}\sigma(g_{sc})^{-1}$.  Alors $\omega(g)$ est le produit par l'accouplement
  $$ H^{1,0}(\Gamma_{F};T_{sc}\to T)\times H^{1,0}(W_{F};\hat{T}\to \hat{T}_{ad})\to {\mathbb C}^{\times}$$
  des images de $g$ dans le premier groupe et de ${\bf a}$ dans le second.
  
  On v\'erifie sur les constructions que le dual de l'homomorphisme
  $$G_{ab}(F)\stackrel{N^G}{\to}G_{0,ab}(F)$$
  est l'homomorphisme naturel
  $$(1) \qquad H^1(W_{F};Z(\hat{G}_{0}))\to H^1(W_{F};Z(\hat{G})).$$
  On a vu que $Z(\hat{G}_{0})=Z(\hat{G})\cap \hat{T}^{\hat{\theta},0}$. Notons $Z(\hat{G})_{*}$ le groupe des $x\in Z(\hat{G})$ tels que $\sigma(x)x^{-1}\in Z(\hat{G})\cap \hat{T}^{\hat{\theta},0}$ pour tout $\sigma\in \Gamma_{F}$. Le quotient $Z(\hat{G})_{*}/(Z(\hat{G})\cap \hat{T}^{\hat{\theta},0})$ n'est autre que le groupe des invariants $(Z(\hat{G})/(Z(\hat{G})\cap \hat{T}^{\hat{\theta},0}))^{\Gamma_{F}}$. On a un homomorphisme
  $$Z(\hat{G})_{*}/(Z(\hat{G})\cap \hat{T}^{\hat{\theta},0})Z(\hat{G})^{\Gamma_{F}}\to H^1(W_{F};Z(\hat{G}_{0}))$$
  qui, \`a $x\in Z(\hat{G})_{*}$, associe le cocycle $w\mapsto w(x)x^{-1}$. On v\'erifie qu'il se quotiente en un isomorphisme de $Z(\hat{G})_{*}/(Z(\hat{G})\cap \hat{T}^{\hat{\theta},0})Z(\hat{G})^{\Gamma_{F}}$ sur le noyau de l'homomorphisme (1). Le groupe $Z(\hat{G})_{*}/(Z(\hat{G})\cap \hat{T}^{\hat{\theta},0})Z(\hat{G})^{\Gamma_{F}}$ s'identifie ainsi au groupe dual de $G_{0,ab}(F)/N^G(G_{ab}(F))$. Pour $x\in Z(\hat{G})_{*}/(Z(\hat{G})\cap \hat{T}^{\hat{\theta},0})Z(\hat{G})^{\Gamma_{F}}$, on note $\mu_{x}$ le caract\`ere  associ\'e de $G_{0,ab}(F)/N^G(G_{ab}(F))$.
  
  L'application $N^{\tilde{G}}:\tilde{G}_{ab}(F)\to \tilde{G}_{0,ab}(F)$ \'etant compatible \`a $N^G$, on voit qu'\`a tout $x\in Z(\hat{G})_{*}/(Z(\hat{G})\cap \hat{T}^{\hat{\theta},0})Z(\hat{G})^{\Gamma_{F}}$, on peut aussi associer une fonction $\tilde{\mu}_{x}$ sur $\tilde{G}_{0,ab}(F)$ telle que
  
  (2) $\tilde{\mu}_{x}$ vaut $1$ sur $N^{\tilde{G}}(\tilde{G}_{ab}(F))$;
  
  (3) $\tilde{\mu}_{x}(g_{0}\gamma_{0})=\mu_{x}(g_{0})\tilde{\mu}_{x}(\gamma_{0})$ pour tous $g_{0}\in G_{0,ab}(F)$ et tout $\gamma_{0}\in \tilde{G}_{0,ab}(F)$.
  
  Pour $\gamma_{0}\in \tilde{G}_{0,ab}(F)$, la somme 
  $$\vert Z(\hat{G})_{*}/(Z(\hat{G})\cap \hat{T}^{\hat{\theta},0})Z(\hat{G})^{\Gamma_{F}}\vert ^{-1}\sum_{x\in Z(\hat{G})_{*}/(Z(\hat{G})\cap \hat{T}^{\hat{\theta},0})Z(\hat{G})^{\Gamma_{F}}}\tilde{\mu}_{x}(\gamma_{0})$$
  vaut $1$ si $\gamma_{0}\in N^{\tilde{G}}(\tilde{G}_{ab}(F))$, $0$ sinon.

\bigskip

\subsection{Image de la correspondance}
Soit ${\bf G}'=(G',{\cal G}',\tilde{s})$ une donn\'ee endoscopique pour $(G,\tilde{G},{\bf a})$.   Rappelons que ${\bf G}'$ est dit elliptique si et seulement si $Z(\hat{G}')^{\Gamma_{F},0}=Z(\hat{G})^{\Gamma_{F},\hat{\theta}, 0}$.

{\bf D\'efinition.} Nous dirons qu'un \'el\'ement semi-simple $\gamma\in \tilde{G}(F)$ est elliptique si et seulement s'il existe un tore tordu maximal elliptique $\tilde{T}$ de $\tilde{G}$ tel que $\gamma\in \tilde{T}(F)$.

  Si $F$ est non-archim\'edien, cette condition \'equivaut \`a l'\'egalit\'e  $A_{G_{\gamma}}=A_{\tilde{G}}$. Si $F$ est archim\'edien, la condition d'ellipticit\'e entra\^{\i}ne cette \'egalit\'e $A_{G_{\gamma}}=A_{\tilde{G}}$, mais la r\'eciproque n'est pas toujours vraie.

\ass{Proposition}{(i) Soit $(\delta,\gamma)\in {\cal D}({\bf G}')$. Alors les images de $\delta$ et $\gamma$ dans $\tilde{G}_{0,ab}(F)$ par le diagramme 1.12(11) sont \'egales.

(ii) Supposons ${\bf G}'$ elliptique et $F\not={\mathbb R}$. Soit $\delta\in \tilde{G}'_{ss}(F)$. On suppose que $\delta$ est elliptique et $\tilde{G}$-r\'egulier, et que l'image de $\delta$ dans $\tilde{G}_{0,ab}(F)$ appartient \`a l'image de $\tilde{G}_{ab}(F)$ par l'application $N^{\tilde{G}}$. Alors il existe $\gamma\in \tilde{G}(F)$ tel que $(\delta,\gamma)$ appartienne \`a ${\cal D}({\bf G}')$.

(iii) Supposons $F={\mathbb R}$. L'assertion (ii) devient vraie si l'on remplace $\tilde{G}_{ab}(F)$ et ${\cal D}({\bf G}')$ par $K\tilde{G}_{ab}({\mathbb R})$ et ${\cal D}_{K\tilde{G}}({\bf G}')$.}

Preuve.  Soit $(\delta,\gamma)\in {\cal D}({\bf G}')$. Gr\^ace au lemme 1.10, on choisit un diagramme $(\delta,B',T',B,T,\gamma)$ et on utilise les notations de 1.10 pour celui-ci. On note  $\xi_{sc}:T_{sc}\to T'_{sc}$ l'homomorphisme relevant $\xi_{T,T'}$, o\`u $T'_{sc}$ est l'image r\'eciproque de $T'$ dans $G'_{sc}$. Cet homomorphisme est \'equivariant pour les actions galoisiennes. On n'a aucun mal \`a relever 1.10(6) sous la forme: on peut \'ecrire $\gamma=e\pi(t)$, $\delta=e'\pi(t')$, avec $t\in T_{sc}$, $e\in Z(\tilde{G},{\cal E})$ et $t'=\xi_{sc}(t)$. D'apr\`es les d\'efinitions, les images de $\delta$ et $\gamma$ dans $\tilde{G}_{0,ab}(F)$ sont repr\'esent\'es respectivement par les couples $(\nu',e_{0})$ et $(\nu,e_{0})$, o\`u $\nu'(\sigma)=t'\sigma(t')^{-1}$, $\nu(\sigma)=ad_{e}^{-1}(u_{{\cal E}}(\sigma))t\sigma(t)^{-1}u_{{\cal E}}(\sigma)^{-1}$  et $e_{0}$ est  l'image de $e$ dans ${\cal Z}_{0}(\tilde{G})$. Pour prouver (i), il suffit de prouver l'\'egalit\'e $\xi_{sc}(\nu(\sigma))=\nu'(\sigma)$. Puisque $\nu(\sigma)$ est central, on a aussi bien $\nu(\sigma)=u_{{\cal E}}(\sigma)^{-1}ad_{e}(u_{{\cal E}}(\sigma))t^{-1}\sigma(t)$. On sait que $u_{{\cal E}}(\sigma)$ d\'efinit un \'el\'ement de $W^{\theta}$ que l'on peut relever en un \'el\'ement de $G_{e}$. On peut donc \'ecrire $u_{{\cal E}}(\sigma)=n(\sigma)t(\sigma)$, o\`u $n(\sigma)\in G_{e}$ et $t(\sigma)\in T_{sc}$. Alors $\nu(\sigma)=(\theta^{-1}-1)(t(\sigma))t\sigma(t)^{-1}$, d'o\`u $\xi_{sc}(\nu(\sigma))=\xi_{sc}(t\sigma(t)^{-1})$. Puisque $\xi_{sc}$ est \'equivariant pour les actions galoisiennes, on en d\'eduit l'\'egalit\'e cherch\'ee $\xi_{sc}(\nu(\sigma))=\nu'(\sigma)$.

Pla\c{c}ons-nous sous les hypoth\`eses de (ii). On choisit une paire de Borel $(B',T')$ de $G'$ conserv\'ee par $ad_{\delta}$ et on identifie $\underline{la}$ paire de Borel \'epingl\'ee ${\cal E}^*$ de $G$ \`a une paire particuli\`ere. On choisit une cocha\^{\i}ne $u_{{\cal E}^*}$ pour cette paire, on la note simplement $u^*$. Munissons $G$ de l'action galoisienne $\sigma\mapsto \sigma_{G^*}=ad_{u^*(\sigma)}\circ\sigma$. Sa restriction \`a $T^*$ est l'action d\'ej\`a introduite sur ce tore et $G$ est quasi-d\'eploy\'e pour cette action. Posons $\theta=\theta_{{\cal E}^*}$. Les deux paires $(B',T')$ et $(B^*,T^*)$ d\'eterminent un homomorphisme $\xi_{T^*,T'}:T^*\to T'$. Il y a un cocycle $\omega_{T'}:\Gamma_{F}\to W^{\theta}$ tel que $\sigma_{G'}\circ\xi_{T^*,T'}\circ\sigma_{G^*}^{-1}=\xi_{T^*,T'}\circ \omega_{T'}(\sigma)$. Le groupe $G_{SC}^{\theta}$ est lui-aussi quasi-d\'eploy\'e. D'apr\`es [K1] corollaire 2.2, on peut fixer $g\in G_{SC}^{\theta}$ tel qu'en posant $T=ad_{g^{-1}}(T^*)$, le tore $T$ soit d\'efini sur $F$ pour l'action $\sigma\mapsto \sigma_{G^*}$ et $\xi_{T,T'}=\xi_{T^*,T'}\circ ad_{g}$ v\'erifie $\sigma_{G'}\circ\xi_{T,T'}=\xi_{T,T'}\circ \sigma_{G^*}$. Remarquons qu'en posant ${\cal E}=ad_{g}^{-1}({\cal E}^*)$ et $B=ad_{g}^{-1}(B^*)$, l'homomorphisme $\xi_{T,T'}$ est celui associ\'e aux deux paires $(B',T')$ et $(B,T)$. D'autre part, puisque $g$ est fixe par $\theta$, on a $Z(\tilde{G},{\cal E})=Z(\tilde{G},{\cal E}^*)$ et $\theta=\theta_{{\cal E}}$.
 
  Par hypoth\`ese, l'image de $\delta$ dans $\tilde{G}_{0,ab}(F)$ est aussi l'image d'un \'el\'ement de $\tilde{G}_{ab}(F)$. On peut repr\'esenter ce dernier par un \'el\'ement $(\mu,e)\in Z^{1,0}(\Gamma_{F}; G_{SC}\circlearrowleft \tilde{G})$, o\`u $e$ appartient \`a $Z(\tilde{G},{\cal E})$. Son image dans $H^{1,0}(\Gamma_{F};{\cal Z}_{0}(G_{SC})\circlearrowleft {\cal Z}_{0}(\tilde{G}))$ est repr\'esent\'ee par le couple $(\nu_{0},e_{0})$ suivant:  $e_{0}$ est l'image de $e$ dans ${\cal Z}_{0}$ et  $\nu_{0}(\sigma)$ est l'image  de $\nu(\sigma)=\theta^{-1}(u^*(\sigma))\mu(\sigma)u^*(\sigma)^{-1}$  dans ${\cal Z}_{0}(G_{SC})$.  D'apr\`es la preuve de la bijectivit\'e de l'application 1.12(8), on peut \'ecrire $\delta=f'\pi(t')$, o\`u $t'\in T'_{sc}$, $f\in {\cal Z}(\tilde{G})$ et $f'$ est l'image de $f$ dans ${\cal Z}(\tilde{G}')$. L'image de $\delta$ dans $H^{1,0}(\Gamma_{F};{\cal Z}_{0}(G_{SC})\circlearrowleft {\cal Z}_{0}(\tilde{G}))$ est repr\'esent\'ee par le couple $(\nu',f_{0})$, o\`u $\nu'(\sigma)=t'\sigma(t')^{-1}$ et $f_{0}$ est l'image de $f$ dans ${\cal Z}_{0}(\tilde{G})$. L'\'egalit\'e des images de $\delta$ et $(\mu,e)$ signifie que les couples $(\nu_{0},e_{0})$ et $(\nu',f_{0})$ sont cohomologues, c'est-\`a-dire qu'il existe $z\in Z(G_{SC})$ tel que $\nu'(\sigma)=z^{-1}\nu_{0}(\sigma)\sigma(z)$ et $f_{0}=e_{0}z$ (pour simplifier, on note encore $z$ l'image de cet \'el\'ement dans divers quotients de $Z(G_{SC})$). Quitte \`a remplacer le couple $(\mu,e)$ par le couple cohomologue $(\mu',ez)$, o\`u $\mu'(\sigma)=z^{-1}\mu(\sigma)\sigma(z)$, on se ram\`ene \`a la situation o\`u $f_{0}=e_{0}$, donc $f'=e'$, et $\nu'=\nu_{0}$.  Rappelons que $\nu$ est \`a valeurs dans $Z(G_{SC})\subset T$. L'\'egalit\'e $\nu'=\nu_{0}$ signifie que
   signifie que  $\xi_{sc}(\nu(\sigma))=\nu'(\sigma)$ pour tout $\sigma\in \Gamma_{F}$, o\`u $\xi_{sc}:T_{sc}\to T'_{sc}$ rel\`eve $\xi_{T,T'}$. Soit $t\in T_{sc}$ tel que $\xi_{sc}(t)=t'$. D'apr\`es l'\'equivariance de $\xi_{sc}$, l'\'egalit\'e pr\'ec\'edente signifie que $\nu(\sigma)$ et $t\sigma_{G^*}(t)^{-1}$ ont m\^eme image dans $T_{sc}/(1-\theta)(T_{sc})$. On peut choisir une cocha\^{\i}ne $y:\Gamma_{F}\to T_{sc}$ telle que 
 $$(1) \qquad \nu(\sigma)=(1-\theta^{-1})(y(\sigma))t\sigma_{G^*}(t)^{-1}.$$
 On note $d$ la diff\'erentielle pour l'action naturelle $\sigma\mapsto \sigma_{G}$ et $d^*$ celle pour l'action $\sigma\mapsto \sigma_{G^*}$. Puisque $\nu$ est \`a valeurs centrales, on a $d\nu=d^*\nu$. D'autre part, $\theta$ commute \`a l'action $\sigma\mapsto \sigma_{G^*}$. De l'\'egalit\'e ci-dessus se d\'eduit la relation $d\nu=(1-\theta^{-1})(d^*y)$ puis $(1-\theta^{-1})(du^*d^*y)=1$ gr\^ace \`a 1.12(4). Puisque $du^*$ est \`a valeurs centrales, c'est un cocycle pour chacune des actions galoisiennes. Donc $du^*d^*y$ est un cocycle pour l'action $\sigma\mapsto \sigma_{G^*}$ et l'\'egalit\'e pr\'ec\'edente montre qu'il prend ses valeurs dans $T_{sc}^{\theta}$. 
 
 {\bf Remarque.} La notation $T_{sc}^{\theta}$ d\'esigne l'ensemble des points fixes par $\theta$ dans $T_{sc}$, et non pas l'image r\'eciproque dans $G_{SC}$ de $T^{\theta}$. L'ensemble $T_{sc}^{\theta}$ est connexe, donc est un tore. 
 
   Les hypoth\`eses d'ellipticit\'e de ${\bf G}'$ et de $\delta$  et l'\'equivariance de $\xi_{sc}$ entra\^{\i}nent que ce tore $T_{sc}^{\theta}$, muni de l'action $\sigma\mapsto \sigma_{G^*}$, est elliptique. Donc $H^2(\Gamma_{F},T_{sc}^{\theta})=0$ et $du^*d^*y$ est le cobord d'une cocha\^{\i}ne \`a valeurs dans $T_{sc}^{\theta}$. Quitte \`a multiplier $y$ par l'inverse de  cette cocha\^{\i}ne, on peut supposer $du^*d^*y=1$. Posons $Y(\sigma)=y(\sigma)u^*(\sigma)$. L'\'egalit\'e pr\'ec\'edente et un calcul standard montrent que $Y$ est un cocycle pour l'action naturelle $\sigma\mapsto \sigma_{G}=ad_{u^*(\sigma)^{-1}}\circ \sigma_{G^*}$. Posons $\gamma_{1}=et$ (ou plus exactement $\gamma_{1}=e\pi(t)$). Puisque $(\mu,e)$ appartient \`a $Z^{1,0}(\Gamma_{F}; G_{SC}\circlearrowleft \tilde{G})$, on a $\sigma(e)=e\mu(\sigma)$, d'o\`u $\sigma(\gamma_{1})=e\mu(\sigma)\sigma(t)$. On a
   $$\mu(\sigma)\sigma(t)=\theta^{-1}(u^*(\sigma)^{-1})\nu(\sigma)u^*(\sigma)u^*(\sigma)^{-1}\sigma_{G^*}(t)u^*(\sigma).$$
   En utilisant (1), on obtient $\mu(\sigma)\sigma(t)=\theta^{-1}(Y(\sigma)^{-1})tY(\sigma)$, d'o\`u 
   $$(2) \qquad \sigma(\gamma_{1})=\pi(Y(\sigma)^{-1})\gamma_{1}\pi(Y(\sigma)),$$
 o\`u on a r\'etabli l'homomorphisme $\pi$ pour plus de pr\'ecision. Jusque-l\`a, nous n'avons pas utilis\'e l'hypoth\`ese que $F$ est non archim\'edien. Utilisons-la. Le cocycle $Y$ est \`a valeurs dans $G_{SC}$. Or $H^1(\Gamma_{F},G_{SC})=0$. Donc on peut choisir $g_{1}\in G_{SC}$ tel que $Y(\sigma)=g_{1}^{-1}\sigma(g_{1})$. Posons $\gamma=g_{1}\gamma_{1}g_{1}^{-1}$. La relation (2) implique que $\gamma$ appartient \`a $\tilde{G}(F)$.  La classe de conjugaison sur $\bar{F}$ de $\gamma$ est la m\^eme que celle de $\gamma_{1}$. En appliquant les d\'efinitions de 1.8, la d\'efinition $\gamma_{1}=et$ montre que sa classe correspond \`a celle de $\delta$. Cela prouve (ii).
 
 Supposons maintenant $F={\mathbb R}$ et consid\'erons un $K$-espace tordu. On peut supposer qu'il est issu d'un couple $(G,\tilde{G})$ comme en 1.11. On a encore (2). Fixons $\gamma_{2}\in \tilde{G}(F)$, \'ecrivons $\gamma_{1}=x\gamma_{2}$, avec $x\in G$. La relation (2) entra\^{\i}ne
$$ad_{\gamma_{2}}\circ\pi(Y(\sigma))=x^{-1}\pi(Y(\sigma))\sigma(x).$$
Donc la classe du cocycle $\pi(Y)$ est fixe par $\theta$. Il existe $p\in {\Pi}$ et $g_{1}\in G$ tels que $\pi(Y(\sigma))=g_{1}^{-1}\pi(p(\sigma))\sigma(g_{1})$. La relation (2) se r\'ecrit
$$\sigma(g_{1}\gamma_{1}g_{1}^{-1})=ad_{p(\sigma)^{-1}}(g_{1}\gamma_{1}g_{1}^{-1}).$$
Posons $\gamma=\tilde{\phi}_{p}^{-1}(g_{1}\gamma_{1}g_{1}^{-1})$. Alors $\gamma$ appartient \`a $\tilde{G}_{p}({\mathbb R})$ et, de nouveau, les classes de conjugaison de $\gamma$ et $\delta$ se correspondent. Cela prouve (iii). $\square$

 \bigskip
 
 \section{Transfert}
 
 \bigskip
\subsection{Facteurs de transfert}
La situation est la m\^eme qu'en 1.5. Soit ${\bf G}'=(G',{\cal G}',\tilde{s})$ une donn\'ee endoscopique relevante pour $(G,\tilde{G},{\bf a})$. On introduit des donn\'ees auxiliaires $G'_{1}$, $\tilde{G}'_{1}$, $C_{1}$, $\hat{\xi}_{1}$. Le terme $G'_{1}$ est un groupe r\'eductif connexe d\'efini et quasi-d\'eploy\'e sur $F$, $C_{1}\subset G'_{1}$ est un tore central d\'efini sur $F$ et induit (c'est-\`a-dire que $X_{*}(C_{1})$ poss\`ede une base conserv\'ee par l'action de $\Gamma_F$). Il y a une suite exacte
$$1\to C_{1}\to G'_{1}\to G' \to 1.$$
Le terme $\tilde{G}'_{1}$ est un espace tordu sur $G'_{1}$, d\'efini sur $F$, \`a torsion int\'erieure, tel que $\tilde{G}'_{1}(F)\not=\emptyset$. Il y a une surjection $\tilde{G}'_{1}\to \tilde{G}'$ compatible avec la surjection $G'_{1}\to G'$. Le terme $\hat{\xi}_{1}:{\cal G}'\to {^LG}'_{1}$ est un plongement  compatible aux projections sur $W_{F}$ dont la restriction \`a $\hat{G}'$ est un homomorphisme $\hat{G}'\to \hat{G}'_{1}$ dual de $G'_{1}\to G'$. Il existe de telles donn\'ees auxiliaires, cf. [KS1] paragraphe 2.2. Fixons-en.

 Pour $w\in W_{F}$, soit $g_{w}=(g(w),w)\in {\cal G}'$. Ecrivons $\hat{\xi}_{1}(g_{w})=(g'_{1}(w),w)$. L'image $z_{C_{1}}(w)$ de $g'_{1}(w)$ dans $\hat{G}'_{1}/\hat{G}'=\hat{C}_{1}$ ne d\'epend pas du choix de $g_{w}$. L'application $w\mapsto z_{C_{1}}(w)$ est un cocycle, qui d\'etermine un caract\`ere $\lambda_{1}$ de $C_{1}(F)$.
 
 Notons ${\cal D}_{1}$ l'ensemble des $(\delta_{1},\gamma)\in \tilde{G}'_{1}(F)\times \tilde{G}(F)$ tels que $(\delta,\gamma)\in {\cal D}({\bf G}')$, o\`u $\delta$ est l'image de $\delta_{1}$ dans $\tilde{G}'(F)$. Kottwitz et Shelstad d\'efinissent ce que l'on peut appeler un bifacteur de transfert, que l'on note  $\boldsymbol{\Delta}_{1}:{\cal D}_{1}\times {\cal D}_{1}\to {\mathbb C}^{\times}$. On rappelle sa d\'efinition (l\'eg\`erement modifi\'ee: on  supprime les termes $\Delta_{IV}$) au paragraphe suivant. Il ne d\'epend que des donn\'ees d\'ej\`a fix\'ees. Un facteur de transfert est une application $\Delta_{1}:{\cal D}_{1}\to {\mathbb C}^{\times}$ telle que 
 $$\Delta_{1}(\delta_{1},\gamma)\Delta_{1}(\underline{\delta}_{1},\underline{\gamma})^{-1}=\boldsymbol{\Delta}_{1}(\delta_{1},\gamma;\underline{\delta}_{1},\underline{\gamma}).$$
 Il  existe un tel facteur. Il est unique \`a homoth\'etie pr\`es. La valeur $\Delta_{1}(\delta_{1},\gamma)$ ne d\'epend que de la classe de conjugaison stable de $\delta_{1}$ (on rappelle que, $\delta_{1}$ \'etant fortement r\'egulier, sa classe de conjugaison stable est l'intersection de $\tilde{G}'_{1}(F)$ avec la classe de conjugaison g\'eom\'etrique de $\delta_{1}$, c'est-\`a-dire sa classe de conjugaison par $G'_{1}=G'_{1}(\bar{F})$). Pour $c_{1}\in C_{1}(F)$ et $g\in G(F)$, on a l'\'egalit\'e
 $$\Delta_{1}(c_{1}\delta_{1},g^{-1}\gamma g)=\lambda_{1}(c_{1})^{-1}\omega(g)\Delta_{1}(\delta_{1},\gamma).$$
 
 Supposons $F={\mathbb R}$ et consid\'erons un $K$-espace $K\tilde{G}$. En utilisant \'evidemment les m\^emes donn\'ees auxiliaires pour chaque espace $\tilde{G}_{p}$, on d\'efinit l'ensemble ${\cal D}_{K\tilde{G},1}$ r\'eunion disjointe des ${\cal D}_{\tilde{G}_{p},1}$ relatifs \`a chaque $\tilde{G}_{p}$. Comme l'a remarqu\'e Kottwitz, on peut d\'efinir un bifacteur de transfert $\boldsymbol{\Delta}_{1}:{\cal D}_{K\tilde{G},1}\times {\cal D}_{K\tilde{G},1}\to {\mathbb C}^{\times}$, cf. 2.3.

  \bigskip
 \subsection{D\'efinition du bifacteur de transfert}
On conserve la situation du paragraphe pr\'ec\'edent. On fixe des paires de Borel \'epingl\'ees $\hat{\cal E}$ et $\hat{\cal E}'$ comme en 1.5 et on utilise les constructions de ce paragraphe relatives \`a ces paires.  On fixe deux \'el\'ements  $(\delta_{1},\gamma)$ et $(\underline{\delta}_{1},\underline{\gamma})$ de ${\cal D}_{1}$.

 On fixe un diagramme $(\delta,B',T',B,T,\gamma)$ et on utilise pour celui-ci les notations de 1.10.  On compl\`ete $(B,T)$ en une paire de Borel \'epingl\'ee ${\cal E}$. On fixe $e\in Z(\tilde{G},{\cal E})$ et on pose $\theta=\theta_{e}$. On note $\Sigma(T)$ l'ensemble des racines de $T$ dans l'alg\`ebre de Lie de $G$. Il s'identifie  \`a $\Sigma(T^*)$ par l'identification ${\cal E}\simeq {\cal E}^*$. Mais il est muni  d'une action galoisienne naturelle du fait que $T$ est d\'efini sur $F$ et c'est cette action que l'on consid\`ere dans la suite.  L'automorphisme $\theta$ agit sur $\Sigma(T)$. Comme en 1.6, on note $\Sigma(T)_{res}$ l'ensemble des restrictions $\alpha_{res}$ d'\'el\'ements $\alpha\in \Sigma(T)$ \`a $T^{\theta,0}$. On note $\Sigma_{res,ind}$ le sous-ensemble des \'el\'ements indivisibles de $\Sigma(T)_{res}$. On fixe des $a$-data $(a_{\alpha})_{\alpha\in \Sigma(T)_{res,ind}}$   pour l'ensemble $\Sigma(T)_{res,ind}$ muni de son action galoisienne, cf. [LS] paragraphe 2.2.  On les rel\`eve en des $a$-data pour $\Sigma(T)$ en posant $a_{\alpha}=a_{\alpha_{res}}$ si $\alpha_{res}$ est indivisible, $a_{\alpha}=a_{\alpha_{res}/2}$ sinon.   On d\'efinit une fonction $r_{T}:\Gamma_{F}\to T_{sc}^{\theta}$ par
$$r_{T}(\sigma)=\prod_{\alpha\in \Sigma(T),\alpha>0,\sigma^{-1}(\alpha)<0}\check{\alpha}(a_{\alpha}),$$
o\`u la positivit\'e est relative \`a $B$ et o\`u on consid\`ere que les coracines prennent leurs valeurs dans $G_{SC}$.
Comme en 1.2, on fixe pour tout $\sigma\in \Gamma_{F}$ un \'el\'ement $u_{{\cal E}}(\sigma)\in G_{SC}$ tel que $ad_{u_{{\cal E}}(\sigma)}\circ\sigma$ conserve ${\cal E}$.  L'\'el\'ement $u_{{\cal E}}(\sigma)^{-1}$ d\'efinit un \'el\'ement de $W^{\theta}$ que nous notons $\omega_{T}(\sigma)$. D'autre part, \`a la paire de Borel \'epingl\'ee ${\cal E}$ est associ\'ee une section de Springer $n_{{\cal E}}:W\to G_{SC}$, cf. [LS] 2.1.   On d\'efinit une cocha\^{\i}ne $V_{T}:\Gamma_{F}\to T_{sc}$ par
$$V_{T}(\sigma)=r_{T}(\sigma)n_{{\cal E}}(\omega_{T}(\sigma))u_{{\cal E}}(\sigma).$$
Notons que  $n_{{\cal E}}(\omega_{T}(\sigma))\in G_{SC,e}$ car $n_{{\cal E}}$ est \'equivariante pour l'action de $\theta$. On v\'erifie que $dV_{T}=du_{{\cal E}}$. Notons   $T'_{1}$ le commutant de $\delta_{1}$ dans $G'_{1}$. On a deux homomorphismes \'equivariants pour les actions galoisiennes
$$T'_{1}\stackrel{\xi_{T'_{1},T'}}{\to} T'\stackrel{\xi_{T,T'}}{ \leftarrow }T.$$
Notons $\mathfrak{T}_{1}$ le produit fibr\'e de $T'_{1}$ et $T$ au-dessus de $T'$, c'est-\`a-dire 
$$\mathfrak{T}_{1}=\{(t_{1},t)\in T'_{1}\times T; \xi_{T'_{1},T'}(t_{1})=\xi_{T,T'}(t)\}.$$
Notons $e'$ l'image naturelle de $e$ dans ${\cal Z}(\tilde{G}')$. Relevons-le en un \'el\'ement $e'_{1}\in {\cal Z}(\tilde{G}_{1})$. Ecrivons $\gamma=\nu e$ et $\delta_{1}=\mu_{1}e'_{1}$ puis posons $\nu_{1}=(\mu_{1},\nu)$. Alors $ \nu_{1}$ appartient \`a $\mathfrak{T}_{1}$: l'image commune de $\nu$ et $\mu_{1}$ dans $T'$ est l'\'el\'ement $\mu$ tel que $\delta=\mu e'$. Remarquons qu'il y a un homomorphisme naturel $1-\theta:T_{sc}\to \mathfrak{T}_{1}$: \`a $t_{sc}\in T_{sc}$, il associe le couple $(1,(1-\theta)\circ\pi(t_{sc}))\in \mathfrak{T}_{1}$. On v\'erifie l'\'egalit\'e
$$(1-\theta)(V_{T}(\sigma))=(z_{1}(\sigma),z(\sigma))\sigma(\nu_{1})\nu_{1}^{-1},$$
o\`u $z(\sigma)$ et $z_{1}(\sigma)$ sont les \'el\'ements de $Z(G)$, resp. $Z(G'_{1})$, tels que $u_{{\cal E}}(\sigma)\sigma(e)u_{{\cal E}}(\sigma)^{-1}=z(\sigma)^{-1}e$, resp. $\sigma(e'_{1})=z_{1}(\sigma)^{-1}e'_{1}$.

On effectue les m\^emes constructions pour la paire $(\underline{\delta}_{1},\underline{\gamma})$. On utilise les m\^emes notations, en les soulignant. Il est essentiel d'effectuer pour ces donn\'ees des choix coh\'erents avec ceux faits pour la premi\`ere paire. Pour cela, on fixe $r\in G_{SC}$ tel que $ad_{r}({\cal E})=\underline{\cal E}$. On choisit $\underline{e}=ad_{r}(e)$, $u_{\underline{\cal E}}(\sigma)=ru_{{\cal E}}(\sigma)\sigma(r)^{-1}$ et $\underline{e}'_{1}=e'_{1}$ (ce dernier choix est loisible puisque $e$ et $\underline{e}$ ont m\^eme image $e'$ dans ${\cal Z}(\tilde{G}')$.  D\'efinissons le tore $U=(T_{sc}\times\underline{T}_{sc})/diag_{-}(Z(G_{SC}))$, o\`u $diag_{-}$ est le plongement antidiagonal.   On d\'efinit une cocha\^{\i}ne $V:\Gamma_{F}\to U$: $V(\sigma)$ est l'image dans $U$ de $(V_{T}(\sigma),V_{\underline{T}}(\sigma)^{-1})$. C'est un cocycle. Introduisons le groupe $\mathfrak{Z}_{1}$ form\'e des couples $(z_{1},z)\in Z(G'_{1})\times Z(G)$ qui ont m\^eme image dans $Z(G')$. D\'efinissons le tore $S_{1}=(\mathfrak{T}_{1}\times\underline{\mathfrak{T}}_{1})/diag_{-}(\mathfrak{Z}_{1})$. Notons $\boldsymbol{\nu}_{1}=(\nu_{1},\underline{\nu}_{1}^{-1})$. Des homomorphismes $1-\theta$ d\'efinis ci-dessus se d\'eduit un autre homomorphisme $1-\theta:U\to S_{1}$. On v\'erifie que le couple $(V,\boldsymbol{\nu}_{1})$ appartient \`a $Z^{1,0}(\Gamma_{F};U\stackrel{1-\theta}{\to}S_{1})$.

On va effectuer des constructions similaires du c\^ot\'e dual. Des deux paires  de Borel ${\cal E}$ et $\hat{\cal E}$ se d\'eduisent  des isomorphismes en dualit\'e $X_{*}(T)\simeq X^* (\hat{ T})$ et $X^*(T)\simeq X_{*}(\hat{ T})$.  Pour $\sigma\in \Gamma_{F}$, on a d\'efini plus haut l'\'el\'ement $\omega_{T}(\sigma)\in W^{\theta}$. On peut munir le tore $\hat{ T}$ d'une nouvelle action galoisienne de sorte que $\sigma$ agisse par $\sigma_{T}=\omega_{T}(\sigma)\sigma_{G}$ (o\`u $\sigma_{G}$ est l'action qui conserve $\hat{\cal E}$, cf. 1.5).  On v\'erifie que, pour cette action, les isomorphismes ci-dessus deviennent \'equivariants, autrement dit $\hat{T}$, muni de cette action, est le tore dual de $T$. C'est cette action que l'on utilise dans la suite. On note $\Sigma(\hat{T})_{res,ind}$ l'ensemble des racines de $\hat{T}^{\hat{\theta},0}$ dans l'alg\`ebre de Lie de $\hat{G}^{\hat{\theta},0}$.  Il s'identifie \`a l'ensemble   des \'el\'ements indivisibles dans $\Sigma(\hat{T})_{res}$, cf. 1.6. Il est de plus muni de l'action galoisienne provenant de celle sur $\hat{T}$. 
On fixe des $\chi$-data $(\chi_{\alpha})_{\alpha\in \Sigma(\hat{T})_{res,ind}}$ pour cette action, cf. 
 [LS] paragraphe 2.5.

Consid\'erons l'ensemble des orbites de  l'action galoisienne dans  $\Sigma(\hat{T})_{res,ind}$. Disons qu'une orbite ${\cal O}$ est sym\'etrique si ${\cal O}=-{\cal O}$ (ou, ce qui revient au m\^eme, si ${\cal O}\cap (-{\cal O})\not=\emptyset$) et qu'elle est asym\'etrique sinon. Consid\'erons un couple $({\cal O},-{\cal O})$ d'orbites asym\'etriques. Fixons $\alpha\in {\cal O}$, notons $F_{\alpha}$ l'extension de $F$ telle que $\Gamma_{F_{\alpha}}$ soit le fixateur de $\alpha$ dans $\Gamma_{F}$. Fixons un ensemble de repr\'esentants $w_{1},...,w_{n}$ du quotient $W_{F_{\alpha}}\backslash W_{F}$. Soit $w\in W_{F}$. Pour tout $i=1,...,n$, il y a un unique couple $(j,v_{i}(w))\in \{1,...,n\}\times W_{F_{\alpha}}$ tel que $w_{i}w=v_{i}(w)w_{j}$. On pose
$$\hat{r}^{{\cal O},-{\cal O}}_{T}(w)=\left(\prod_{\beta\in {\cal O}; \beta>0,w^{-1}\beta<0}\check{\beta}(-1)\right)\left(\prod_{i=1,...,n}(w_{i}^{-1}\check{\alpha})(\chi_{\alpha}(v_{i}(w)))\right).$$
  La positivit\'e est relative \`a $\hat{ B}\cap\hat{G}^{\hat{\theta},0}$. Gr\^ace \`a l'isomorphisme du corps de classes, on a identifi\'e $\chi_{\alpha}$ \`a un caract\`ere de $W_{F_{\alpha}}$. Consid\'erons maintenant une orbite sym\'etrique ${\cal O}$. On fixe $\alpha\in {\cal O}$ et des \'el\'ements $w_{0},w_{1},...,w_{n}\in W_{F}$ de sorte que $w_{0}^{-1}\alpha=-\alpha$ et $i\mapsto w_{i}^{-1}\alpha$ soit une bijection de $\{1,...,n\}$ sur l'ensemble des \'el\'ements positifs de ${\cal O}$. Pour $i=1,...,n$, on pose $w_{-i}=w_{0}w_{i}$. Soit $w\in W_{F}$. Pour tout $i=1,...,n$, il y a un unique couple $(j,v_{i}(w))\in \{\pm 1,...,\pm n\}\times W_{F_{\alpha}}$ (avec la d\'efinition ci-dessus de $F_{\alpha}$) de sorte que $w_{i}w=v_{i}(w)w_{j}$. On pose
$$\hat{r}^{{\cal O}}_{T}(w)=\prod_{i=1,...,n}(w_{i}^{-1}\check{\alpha})(\chi_{\alpha}(v_{i}(w))).$$

On note $\hat{r}_{T}(w)$ le produit des $\hat{r}^{{\cal O},-{\cal O}}(w)$ sur les paires $({\cal O},-{\cal O})$ d'orbites asym\'etriques et des $\hat{r}_{T}^{{\cal O}}(w)$ sur les orbites ${\cal O}$ sym\'etriques. Cela d\'efinit une cocha\^{\i}ne $\hat{r}_{T}:W_{F}\to \hat{T}_{sc}^{\hat{\theta}}$. On a effectu\'e de nombreux choix, mais on montre qu'ils n'affectent cette cocha\^{\i}ne que par multiplication par un cobord, ce qui est sans importance pour la suite.

  On peut effectuer des constructions analogues dans le groupe $\hat{G}'$. Il existe un cocycle $\omega_{T,G'}:\Gamma_{F}\to W^{G'}$ de sorte qu'en munissant le tore $\hat{T}'$  de l'action $(\sigma,t)\mapsto \sigma_{T'}(t)=\omega_{T,G'}(\sigma)\sigma_{G'}(t)$, ce tore  s'identifie au tore dual de $T'$. En fait l'\'egalit\'e  $\hat{T}'=\hat{ T}^{\hat{\theta},0} $ est compatible aux  actions  que l'on vient de d\'efinir sur $\hat{T}'$ et $\hat{T}$.  C'est une cons\'equence du fait que l'application $\xi_{T,T'}:T\to T'$ est \'equivariante pour les actions galoisiennes. On munit l'ensemble $\Sigma(\hat{T}')$ des racines de $\hat{T}'$ dans l'alg\`ebre de Lie de $\hat{G}'$   de l'action galoisienne provenant de celle que l'on vient de d\'efinir sur $\hat{T}'$ . Cet ensemble n'est pas forc\'ement inclus dans $\Sigma(\hat{T})_{res,ind}$, mais il y a n\'eanmoins une injection naturelle du premier dans le second: l'image de $\alpha\in\Sigma(\hat{T}')$ est le seul \'el\'ement de $\Sigma(\hat{T})_{res,ind}\cap\{\alpha/2,\alpha\}$, cf. 1.6.  Cette injection est \'equivariante pour les actions galoisiennes. De nos $\chi$-data se d\'eduisent des $\chi$-data pour l'ensemble $\Sigma(\hat{T}')$. On d\'efinit alors une cocha\^{\i}ne $\hat{r}_{T,G'}:W_{F}\to \hat{T}'_{sc}$ o\`u $\hat{T}'_{sc}$ est l'image r\'eciproque de $\hat{T}'$ dans $\hat{G}'_{SC}$. Sa d\'efinition est copi\'ee sur celle de $\hat{r}_{T}$.

On introduit les sections de Springer $\hat{n}:W^{\theta}\to \hat{G}_{SC}^{\hat{\theta}}$ et $\hat{n}_{G'}:W^{G'}\to \hat{G}'_{SC}$ associ\'ees aux paires de Borel \'epingl\'ees $\hat{\cal E}$ et $\hat{\cal E}'$. Plus exactement, dans le cas de $\hat{n}$, \`a la paire de Borel \'epingl\'ee de $\hat{G}_{SC}^{\hat{\theta}}$ qui se d\'eduit naturellement de $\hat{\cal E}$. Celle-ci a pour paire de Borel sous-jacente la paire $(\hat{B}_{sc}\cap \hat{G}_{SC}^{\hat{\theta}},\hat{T}_{sc}^{\hat{\theta}})$ et les \'el\'ements de l'\'epinglage sont les $\hat{E}_{\alpha}+\hat{E}_{\hat{\theta}\alpha}+...+\hat{E}_{\hat{\theta}^{n_{\alpha}-1}\alpha}$ pour $\alpha\in \Delta$, o\`u $n_{\alpha}\geq1$ est le plus petit entier $n\geq1$ tel que $\hat{\theta}^{n}\alpha=\alpha$.  Rappelons que l'on a modifi\'e l'isomorphisme $^LG\simeq \hat{G}\rtimes W_{F}$, cf. 1.5. On fixe une application
$$\begin{array}{ccc}W_{F}&\to&{\cal G}'\\w&\mapsto &g_{w}=(g(w),w)\\ \end{array}$$
de sorte que $ad_{g(w)}w_{G}$ agisse comme $w_{G'}$ sur $\hat{G}'$. 
Pour $w\in W_{F}$, posons
$$t_{T}(w)= \hat{r}_{T}(w)\hat{n}(\omega_{T}(w))g(w)^{-1}\hat{n}_{G'}(\omega_{T,G'}(w))^{-1}\hat{r}_{T,G'}(w)^{-1} .$$
   L'action galoisienne sur $\hat{T}$, relev\'ee en une action de $W_{F}$, est $w\mapsto  \hat{n}(\omega_{T}(w))w_{G}=ad_{\hat{n}(\omega_{T}(w))g(w)^{-1}}w_{G'}$. Restreinte \`a $\hat{T}'$, elle est \'egale \`a $ad_{\hat{n}_{G'}(\omega_{T,G'}(w))}w_{G'}$. Donc l'\'el\'ement $\hat{n}(\omega_{T}(w))g(w)^{-1}\hat{n}_{G'}(\omega_{T,G'}(w))^{-1}$ appartient \`a $\hat{T}$.  Il en r\'esulte que $t_{T}(w)\in \hat{T}$. On montre que le cobord $dt_{T}$ de la cocha\^{\i}ne $t_{T}$ est \'egal \`a  celui de la cocha\^{\i}ne $w\mapsto g(w)^{-1}$, qui prend ses valeurs dans $Z(\hat{G}')$. Rappelons que l'on a un plongement $\hat{\xi}_{1}:{\cal G}'\to {^LG}'_{1}$. Notons $\hat{T}'_{1}$ le commutant dans $\hat{G}'_{1}$ de $\hat{\xi}_{1}(\hat{T}')$ et $\hat{B}'_{1}$ le groupe engendr\'e par $\hat{T}'_{1}$ et $\hat{\xi}_{1}(\hat{B}')$. Le triplet $(\hat{B}'_{1},\hat{T}'_{1},(\hat{\xi}_{1}(\hat{E}'_{\alpha}))_{\alpha\in \Delta'})$ est une paire de Borel \'epingl\'ee de $\hat{G}'_{1}$. Comme en 1.5, on modifie l'isomorphisme $^LG'_{1}\simeq \hat{G}'_{1}\rtimes W_{F}$ de sorte que l'action d'un \'el\'ement de $W_{F}$ conserve cette paire. On munit $\hat{T}'_{1}$ de la nouvelle action galoisienne $(\sigma,t_{1})\mapsto \sigma_{T}(t_{1})=\omega_{T,G'}(\sigma)\sigma_{G'_{1}}(t_{1})$. Muni de cette action, $\hat{T}'_{1}$ est le tore dual de $T'_{1}$. Posons $\hat{\xi}_{1}(g_{w})=(\zeta_{1}(w),w)$. D'apr\`es la d\'efinition de $g_{w}$, $\zeta_{1}(w)$ appartient au centre de $\hat{G}_{1}$, a fortiori \`a $\hat{T}'_{1}$. Notons $\hat{\mathfrak{T}}_{1}$ le quotient de $\hat{T}'_{1}\times \hat{T}$ par la relation d'\'equivalence $(t_{1}\hat{\xi}(t'),t)=(t_{1},t't)$ pour tout $t'\in \hat{T}'$. C'est le tore dual de $\mathfrak{T}_{1}$. On d\'efinit une cocha\^{\i}ne $\hat{V}_{\mathfrak{T}_{1}}:W_{F}\to \hat{\mathfrak{T}}_{1}$: $\hat{V}_{\mathfrak{ T}_{1}}(w)$ est l'image de $(\zeta_{1}(w),t_{T}(w))$ dans $\hat{\mathfrak{T}}_{1}$. C'est un cocycle.  
   
   On d\'efinit les objets similaires relatifs au tore $\underline{T}$. Remarquons que, quand on oublie les actions galoisiennes, on a l'\'egalit\'e $\hat{\mathfrak{T}}_{1}=\underline{\hat{\mathfrak{T}}}_{1}$ et qu'il y a un homomorphisme naturel $j:\hat{T}_{sc}\to\hat{\mathfrak{T}}_{1}=\underline{\hat{\mathfrak{T}}}_{1}$. On v\'erifie que, pour tout $\omega\in W^{G'}$, l'application  $\omega-1$ de $\hat{\mathfrak{T}}_{1}$ dans lui-m\^eme se rel\`eve en une application naturelle encore not\'ee $\omega-1:\mathfrak{T}_{1}\to \hat{T}_{sc}$. Autrement dit, on a un diagramme commutatif
   $$\begin{array}{ccccc}\hat{\mathfrak{T}}_{1}&&\stackrel{\omega-1}{\to}&&\hat{\mathfrak{T}}_{1}\\&\searrow\,\omega-1&&\nearrow\,j&\\ &&\hat{T}_{sc}&&\\ \end{array}$$
   Notons $\hat{S}_{1}$ le sous-tore de $\hat{\mathfrak{T}}_{1}\times \underline{\hat{\mathfrak{T}}}_{1}\times \hat{T}_{sc}$ form\'e des $(t,\underline{t},t_{sc})$ tels que $j(t_{sc})=t\underline{t}^{-1}$. On le munit de l'action de $\Gamma_{F}$ d\'efinie par
   $$(\sigma,(t,\underline{t},t_{sc}))\mapsto (\sigma_{T}(t),\sigma_{\underline{T}}(\underline{t}),\sigma_{T}(t_{sc})(\omega_{T,G'}(\sigma)\omega_{\underline{T},G'}(\sigma)^{-1}-1)\sigma_{\underline{T}}(\underline{t}))$$
   $$=(\sigma_{T}(t),\sigma_{\underline{T}}(\underline{t}),\sigma_{\underline{T}}(t_{sc})(1-\omega_{\underline{T},G'}(\sigma)\omega_{T,G'}(\sigma)^{-1})\sigma_{T}(t)).$$
   On v\'erifie que $\hat{S}_{1}$ est le tore dual de $S_{1}$. Pour $w\in W_{F}$, on fixe un \'el\'ement $g_{sc}(w)\in \hat{G}_{SC}$ qui ait m\^eme image que $g(w)$ dans $\hat{G}_{AD}$.  On d\'efinit une cocha\^{\i}ne $t_{T,sc}:W_{F}\to \hat{T}_{sc}$ par 
$$t_{T,sc}(w)=\hat{r}_{T}(w)\hat{n}(\omega_{T}(w))g_{sc}(w)^{-1}\hat{n}_{G'}(\omega_{T,G'}(w))^{-1}\hat{r}_{T,G'}(w)^{-1} ,$$
puis la cocha\^{\i}ne $t_{sc}=t_{T,sc}t_{\underline{T},sc}^{-1}$. On d\'efinit ensuite une cocha\^{\i}ne
$\hat{V}_{1}:W_{F}\to \hat{\mathfrak{T}}_{1}\times \underline{\hat{\mathfrak{T}}}_{1}\times \hat{T}_{sc}$  par $\hat{V}_{1}(w)=(\hat{V}_{\mathfrak{ T}_{1}}(w),\hat{V}_{\underline{\mathfrak{ T}}_{1}}(w),t_{sc}(w))$. Elle prend ses valeurs dans $\hat{S}_{1}$ et c'est un cocycle.

Le tore dual de $U$ est $\hat{U}=(\hat{T}_{sc}\times \underline{\hat{T}}_{sc})/diag(Z(\hat{G}_{SC}))$, o\`u $diag$ est le plongement diagonal.   On fixe un \'el\'ement $s_{sc}$ ayant m\^eme image que $s$ dans $\hat{G}_{AD}$ (rappelons que $\tilde{s}=s\hat{\theta}$). On d\'efinit l'\'el\'ement ${\bf s}=(s_{sc},s_{sc})$ de $\hat{U}$. On dispose de l'homomorphisme $1-\hat{\theta}:\hat{S}_{1}\to \hat{U}$ dual de l'homomorphisme $1-\theta:U\to S_{1}$. On v\'erifie que le couple $(\hat{V}_{1},{\bf s})$ appartient \`a $Z^{1,0}(W_{F};\hat{S}_{1}\stackrel{1-\hat{\theta}}{\to}\hat{U})$.

D'apr\`es [KS1] A.3, on dispose d'un produit
$$<.,.>:H^{1,0}(\Gamma_{F};U\stackrel{1-\theta}{\to}S_{1})\times H^{1,0}(W_{F};\hat{S}_{1}\stackrel{1-\hat{\theta}}{\to}\hat{U})\to {\mathbb C}^{\times}.$$
On pose
$$\boldsymbol{\Delta}_{imp}(\delta_{1},\gamma;\underline{\delta}_{1},\underline{\gamma})=<(V,\boldsymbol{\nu}_{1}),(\hat{V}_{1},{\bf s})>^{-1},$$
en notant de la m\^eme fa\c{c}on les \'el\'ements de $Z^{1,0}$ et leurs images dans $H^{1,0}$.

La bijection $\alpha\mapsto \hat{\alpha}$ de $\Sigma(T)$ sur $\Sigma(\hat{T})$ induit une bijection $\alpha_{res}\mapsto \hat{\alpha}_{res}$ de $\Sigma(T)_{res,ind}$ sur $\Sigma(\hat{T})_{res,ind}$. On peut donc consid\'erer nos $\chi$-data comme des $\chi$-data pour l'ensemble $\Sigma(T)_{res,ind}$. Consid\'erons un \'el\'ement de $ \Sigma(T)_{res,ind}$ que l'on \'ecrit $\alpha_{res}$, avec $\alpha\in \Sigma(T)$.   Puisque $\alpha_{res}$ est indivisible, $\alpha$ est du type 1 ou 2. On distingue les cas suivants:

(a) $\alpha$ est de type 1 et $(N\hat{\alpha})(s)\not=1$, autrement dit $(\hat{\alpha})_{res}\not\in \Sigma(\hat{T}')$;

(b) $\alpha$ est de type 2 et $(N\hat{\alpha})(s)\not=\pm 1$, autrement dit ni $(\hat{\alpha})_{res}$, ni $2(\hat{\alpha})_{res}$ n'appartiennent \`a $\Sigma(\hat{T}')$;

(c) $\alpha$ est de type 2 et $(N\hat{\alpha})(s)=-1$, autrement dit $2(\hat{\alpha})_{res}\in \Sigma(\hat{T}')$;

(d) $\alpha$ est de type 1 ou 2 et $(N\hat{\alpha})(s)=1$.

On pose
$$\Delta_{II,\alpha_{res}}(\delta,\gamma)=\left\lbrace\begin{array}{cc}\chi_{\alpha_{res}}(\frac{(N\alpha)(\nu)-1}{a_{\alpha_{res}}}),& \text{ dans le cas (a),}\\ \chi_{\alpha_{res}}(\frac{(N\alpha)(\nu)^2-1}{a_{\alpha_{res}}}),& \text{ dans le cas (b),}\\ \chi_{\alpha_{res}}((N\alpha)(\nu)+1),&\text{ dans le cas (c),}\\ 1,& \text{ dans le cas (d) }\\ \end{array}\right.$$

Ce terme ne d\'epend que de l'orbite de $\alpha_{res}$ pour l'action de $\Gamma_{F}$. On pose
$$\Delta_{II}(\delta,\gamma)=\prod_{\alpha_{res}}\Delta_{II,\alpha_{res}}(\delta,\gamma),$$
o\`u le produit porte sur les orbites de l'action de $\Gamma_{F}$ dans $\Sigma(T)_{res,ind}$.

On d\'efinit alors le bifacteur de transfert
$$\boldsymbol{\Delta}_{1}(\delta_{1},\gamma;\underline{\delta}_{1},\underline{\gamma})=\Delta_{II}(\delta,\gamma)\Delta_{II}(\underline{\delta},\underline{\gamma})^{-1}\boldsymbol{\Delta}_{imp}(\delta_{1},\gamma;\underline{\delta}_{1},\underline{\gamma}).$$

  {\bf Remarques.} (1) Ce terme est ind\'ependant de tous les choix de donn\'ees auxiliaires.
  
  (2) On a rassembl\'e dans le facteur $\boldsymbol{\Delta}_{imp}$ les facteurs plus habituels $\Delta_{I}$ et $\Delta_{III}$. Cela parce que l'on a fait dispara\^{\i}tre le traditionnel groupe $G^*$ qui  nous semble inadapt\'e \`a l'endoscopie tordue.  
  
  (3) On a tent\'e d'incorporer dans les d\'efinitions les changements de signes introduits dans [KS2] 5.4. On n'est pas s\^ur d'avoir r\'eussi.

   \bigskip
  \subsection{Bifacteur de transfert et $K$-groupes}
  On suppose ici $F={\mathbb R}$, on consid\`ere un $K$-espace tordu comme en 1.11. On veut d\'efinir le bifacteur de transfert sur ${\cal D}_{K\tilde{G},1}\times {\cal D}_{K\tilde{G},1}$. On reprend les constructions pr\'ec\'edentes. Du c\^ot\'e dual, il n'y a rien de chang\'e, l'espace $K\tilde{G}$ n'intervenant pas. Du  c\^ot\'e des groupes sur ${\mathbb R}$, les tores $U$ et $S_{1}$ se d\'efinissent aussi bien si $\gamma$ et $\underline{\gamma}$ appartiennent \`a des composantes connexes diff\'erentes de $K\tilde{G}({\mathbb R})$ (il suffit pour les d\'efinir  d'identifier les centres des diff\'erents groupes $G_{p}$).   La seule chose \`a changer est la condition de coh\'erence impos\'ee aux choix de $e$, $u_{{\cal E}}(\sigma)$, $\underline{e}$ et $u_{\underline{\cal E}}(\sigma)$. Dans le paragraphe pr\'ec\'edent, on avait choisi $r\in G_{SC}$ tel que $ad_{r}({\cal E})=\underline{\cal E}$. Supposons maintenant que $\gamma\in \tilde{G}_{p}({\mathbb R})$ et $\underline{\gamma}\in \tilde{G}_{\underline{p}}({\mathbb R})$. On choisit $r\in G_{\underline{p},SC}$ tel que $ ad_{r}\circ \phi_{\underline{p},p}({\cal E})=\underline{\cal E}$. On impose $\underline{e}=ad_{r}\circ\tilde{\phi}_{\underline{p},p}(e)$ et
  $$u_{\underline{\cal E}}(\sigma)=r\phi_{\underline{p},p}(u_{{\cal E}}(\sigma))\nabla_{\underline{p},p}(\sigma)\sigma(r)^{-1}.$$

     \bigskip
  \subsection{Transfert}
 Les donn\'ees sont les m\^emes qu'en 2.1. On fixe une mesure de Haar sur $G(F)$. Soit  $\gamma\in \tilde{G}(F)$.   On pose
 $$D^{\tilde{G}}(\gamma)=\vert det(1-ad_{\gamma_{ss}})_{\vert  \mathfrak{g}/\mathfrak{g}_{\gamma_{ss}}}\vert_{F}, $$
 o\`u $\gamma_{ss}$ est  la partie semi-simple de $\gamma$ et $\vert .\vert_{F} $ la valeur absolue usuelle de $F$. On   fixe une mesure de Haar sur $G_{\gamma}(F)$.  Soit $f\in C_{c}^{\infty}(\tilde{G}(F))$. Dans le cas o\`u  $\omega$ est trivial sur $G_{\gamma}(F)$, on pose
 $$I^{\tilde{G}}(\gamma,\omega,f)=D^{\tilde{G}}(\gamma)^{1/2} \int_{G_{\gamma}(F)\backslash G(F)}\omega(g)f(g^{-1}\gamma g)dg.$$
 Dans le cas o\`u $\omega$ n'est pas trivial sur $ G_{\gamma}(F)$, on pose $I^{\tilde{G}}(\gamma,\omega,f)=0$.
 
 {\bf Remarque.} Il n'est pas  clair que la normalisation que l'on a choisie soit la plus simple. On aurait pu int\'egrer sur $Z_{G}(\gamma;F)\backslash G(F)$  au lieu de $G_{\gamma}(F)\backslash G(F)$.  Auquel cas, la condition sur $\omega$ serait d'\^etre trivial sur $Z_{G}(\gamma;F)$.
 Notons que cela ne cr\'ee pas d'ambigu\"{\i}t\'e: si $\omega$ est trivial sur $G_{\gamma}(F)$ mais pas sur $Z_{G}(\gamma;F)$, l'int\'egrale sur $G_{\gamma}(F)\backslash G(F)$ est nulle.
 \bigskip
 
   On note $I(\tilde{G}(F),\omega)$ le quotient de $C_{c}^{\infty}(\tilde{G}(F))$ par le sous-espace annul\'e par toutes les $I^{\tilde{G}}(\gamma,\omega,.)$, $\gamma$  tr\`es r\'egulier. 
   
   {\bf Remarque.}  Dans le cas o\`u $\omega$ est trivial, on supprime $\omega$ de la notation: $I^{\tilde{G}}(\gamma,f)$ et $I(\tilde{G}(F))$ au lieu de $I^{\tilde{G}}(\gamma,\omega,f)$ et $I(\tilde{G}(F),\omega)$. D'autres simplifications similaires seront utilis\'ees dans la suite.
   \bigskip
       
 On note $C_{c,\lambda_{1}}^{\infty}(\tilde{G}'_{1}(F))$ l'espace des fonctions $f_{1}:\tilde{G}'_{1}(F)\to {\mathbb C}$ telles que $f_{1}(c_{1}\delta_{1})=\lambda_{1}(c_{1})^{-1}f_{1}(\delta_{1})$ pour $c_{1}\in C_{1}(F)$ et $f_{1}$ est lisse et \`a support compact modulo $C_{1}(F)$. On fixe une mesure de Haar sur $G'(F)$. Pour $\delta_{1}\in \tilde{G}'_{1}(F)$, on fixe une mesure de Haar sur $G'_{\delta}(F)$ et, pour $f_{1}\in C_{c,\lambda_{1}}^{\infty}(\tilde{G}'_{1}(F))$, on pose:
 $$I^{\tilde{G}'}(\delta_{1},f_{1})=D^{\tilde{G}'}(\delta)^{1/2}\int_{G'_{\delta}(F)\backslash G'(F)}f_{1}(x^{-1}\delta_{1}x)dx.$$
 Si $\delta_{1}$ est semi-simple fortement r\'egulier, on pose
 $$S^{\tilde{G}'}(\delta_{1},f_{1})=\sum_{\delta'_{1}}I^{\tilde{G}'}(\delta'_{1},f_{1}),$$
 o\`u $\delta'_{1}$ parcourt la classe de conjugaison stable de $\delta_{1}$ modulo conjugaison par $G'(F)$.    On note $SI_{\lambda_{1}}(\tilde{G}'_{1}(F))$ le quotient de $C_{c,\lambda_{1}}^{\infty}(\tilde{G}'_{1}(F))$ par le sous-espace annul\'e par toutes les $S^{\tilde{G}'}(\delta_{1},.)$ pour $\delta_{1}$ fortement r\'egulier.

 On fixe un facteur de transfert $\Delta_{1}$.  Soit  $\delta_{1}\in \tilde{G}'_{1}(F)$, semi-simple et fortement $\tilde{G}$-r\'egulier. Pour $\gamma\in \tilde{G}(F)$ tel que $(\delta_{1},\gamma)\in {\cal D}_{1}$,  il y a un homomorphisme naturel $G_{\gamma}(F)\to G'_{\delta}(F)$, qui est un rev\^etement sur son image. En choisissant un diagramme $(\delta,B',T',B,T,\gamma)$ comme en 1.10, c'est la restriction de $\xi_{T,T'}$ \`a $G_{\gamma}(F)=T^{\theta,0}(F)$.   On fixe les mesures de Haar sur ces deux groupes, de sorte qu'elles se correspondent localement par cet isomorphisme. On pose
 $$d(\theta^*)=\vert det(1-\theta^*)_{\vert \mathfrak{t}^*/(\mathfrak{t}^{*})^{\theta^*}}\vert _{F}.$$
  Pour $f\in C_{c}^{\infty}(\tilde{G}(F))$, on pose
   $$I^{\tilde{G}}(\delta_{1},f)=d(\theta^*)^{1/2}\sum_{\gamma}\Delta_{1}(\delta_{1},\gamma) [Z_{G}(\gamma;F):G_{\gamma}(F)]^{-1}I^{\tilde{G}}(\gamma,\omega,f),$$
 o\`u $\gamma$ parcourt les \'el\'ements de $\tilde{G}(F)$ tels que $(\delta_{1},\gamma)\in {\cal D}_{1}$, modulo conjugaison par $G(F)$.  On montre ([KS1] lemme 4.4.C) que pour tous ces $\gamma$, $\omega$ est trivial sur $Z_{G}(\gamma;F)$, les termes $I^{\tilde{G}}(\gamma,\omega,f)$ sont donc de v\'eritables int\'egrales orbitales. Pour $f_{1}\in C_{c,\lambda_{1}}^{\infty}(\tilde{G}'_{1}(F))$, on dit que $f_{1}$ est un transfert de $f$ si et seulement si $S^{\tilde{G}'}(\delta_{1},f_{1})=I^{\tilde{G}}(\delta_{1},f)$ pour tout $\delta_{1}$ fortement $\tilde{G}$-r\'egulier. On peut d'ailleurs aussi bien demander que cette \'egalit\'e ne soit v\'erifi\'ee que pour un sous-ensemble topologiquement dense. La conjecture de transfert est maintenant prouv\'ee:

 \ass{Th\'eor\`eme}{Tout \'el\'ement de $ C_{c}^{\infty}(\tilde{G}(F))$ admet un transfert dans $C_{c,\lambda_{1}}^{\infty}(\tilde{G}'_{1}(F))$. }
  
  Par passage aux quotients, le transfert appara\^{\i}t comme une application lin\'eaire $I(\tilde{G}(F),\omega)\to SI_{\lambda_{1}}(\tilde{G}'_{1}(F))$. Il d\'epend des choix des donn\'ees auxiliaires, du facteur de transfert et des mesures de Haar. On peut s'affranchir de ce dernier choix en notant $Mes(G(F))$ la droite complexe port\'ee par une mesure de Haar sur $G(F)$. On peut voir le transfert comme une application lin\'eaire 
  $$I(\tilde{G}(F),\omega)\otimes Mes(G(F))\to SI_{\lambda_{1}}(\tilde{G}'_{1}(F))\otimes Mes(G'(F)).$$

  \bigskip
 \subsection{Recollement de donn\'ees auxiliaires}
  Soit ${\bf G}'=(G',{\cal G}',\tilde{s})$ une donn\'ee endoscopique relevante pour $(G,\tilde{G},{\bf a})$. Consid\'erons des donn\'ees $G'_{1}$, $\tilde{G}'_{1}$, $C_{1}$, $\hat{\xi}_{1}$ comme en 2.1, plus un facteur de transfert $\Delta_{1}$. On consid\`ere une autre s\'erie de donn\'ees $G'_{2}$, $\tilde{G}'_{2}$, $C_{2}$, $\hat{\xi}_{2}$, $\Delta_{2}$. On introduit le produit fibr\'e $G'_{12}$ de $G'_{1}$ et $G'_{2}$ au-dessus de $G'$. On a $Z(\hat{G}'_{12})=(Z(\hat{G}'_{1})\times Z(\hat{G}'_{2}))/diag_{-}(Z(\hat{G}'))$. Pour $w\in W_{F}$, soit $g_{w}=(g(w),w)\in {\cal G}'$ tel que $ad_{g_{w}}$ agisse par $w_{G'}$ sur $\hat{G}'$ (on a modifi\'e l'isomorphisme $^LG\simeq \hat{G}\rtimes W_{F}$ comme en 1.5; pour $i=1,2$, on modifie de m\^eme les isomorphismes $^LG'_{i}\simeq \hat{G}'_{i}\rtimes W_{F}$ comme en 2.2). Pour $i=1,2$, on a $\hat{\xi}_{i}(g_{w})=(\zeta_{i}(w),w)$, avec $\zeta_{i}(w)\in Z(\hat{G}'_{i})$. Soit $\zeta_{12}(w)$ l'image de $(\zeta_{1}(w),\zeta_{2}(w)^{-1})$ dans $Z(\hat{G}'_{12})$. Ce terme est bien d\'efini et $\zeta_{12}$ est un cocycle de $W_{F}$ dans $Z(\hat{G}'_{12})$, qui d\'etermine un caract\`ere $\lambda_{12}$ de $G'_{12}(F)$. La restriction de ce caract\`ere \`a $C_{1}(F)\times C_{2}(F)$ est $\lambda_{1}\times \lambda_{2}^{-1}$. Introduisons le produit fibr\'e $\tilde{G}'_{12}$ de $\tilde{G}'_{1}$ et $\tilde{G}'_{2}$ au-dessus de $\tilde{G}'$. Soient $(\delta_{1},\gamma)$ et $(\underline{\delta}_{1},\underline{\gamma})$ deux \'el\'ements de ${\cal D}_{1}$. Soient $\delta_{2}, \underline{\delta}_{2}\in \tilde{G}'_{2}(F)$ tels que $(\delta_{1},\delta_{2}) $ et $(\underline{\delta}_{1},\underline{\delta}_{2})$ appartiennent \`a $\tilde{G}'_{12}(F)$. Alors $(\delta_{2},\gamma)$ et $(\underline{\delta}_{2},\underline{\gamma})$ appartiennent \`a ${\cal D}_{2}$.
  
  \ass{Lemme}{Sous ces hypoth\`eses, on a l'\'egalit\'e
  $$\boldsymbol{\Delta}_{2}(\delta_{2},\gamma;\underline{\delta}_{2},\underline{\gamma})=\lambda_{12}(x_{1},x_{2})\boldsymbol{\Delta}_{1}(\delta_{1},\gamma;\underline{\delta}_{1},\underline{\gamma}),$$
  o\`u $(x_{1},x_{2})\in G_{12}(F)$ est l'\'el\'ement tel que $(\delta_{1},\delta_{2})=(x_{1},x_{2})(\underline{\delta}_{1},\underline{\delta}_{2})$.}
  
  Preuve. On calcule les bifacteurs de transfert en utilisant la d\'efinition de 2.2, en affectant d'un indice $2$ les termes relatifs \`a la deuxi\`eme famille de donn\'ees auxiliaires. Quand on remplace une famille par l'autre, les termes $\Delta_{II}$ ne changent pas et les termes $V$ et ${\bf s}$ non plus. De m\^eme que l'on a d\'efini les tores $\mathfrak{T}_{1}$ et $\mathfrak{T}_{2}$, on introduit le tore $\mathfrak{T}_{12}$ qui est le produit fibr\'e de $T'_{1}$, $T'_{2}$ et $T$ au-dessus de $T'$. On note $\nu_{12}$ l'\'el\'ement $(\mu_{1},\mu_{2},\nu)$ de ce tore. On introduit le groupe $\mathfrak{Z}_{12}$ form\'e des $(z_{1},z_{2},z)\in Z(G'_{1})\times Z(G'_{2})\times Z(G)$ qui ont m\^eme image dans $Z(G')$ puis le tore $S_{12}=(\mathfrak{T}_{12}\times \underline{\mathfrak{T}}_{12})/diag_{-}(\mathfrak{Z}_{12})$. Notons $\boldsymbol{\nu}_{12}$ l'image de $( \nu_{12},\underline{\nu}_{12}^{-1})$ dans $S_{12}$. L'oubli d'une variable d\'efinit des homomorphismes
  $$\begin{array}{ccccc}&&S_{12}&&\\ &\swarrow&&\searrow&\\ S_{1}&&&&S_{2}\\ \end{array}$$
  qui envoient $\boldsymbol{\nu}_{12}$ respectivement sur $\boldsymbol{\nu}_{1}$ et $\boldsymbol{\nu}_{2}$. D'o\`u des homomorphismes
  $$\begin{array}{ccccc}&&H^{1,0}(\Gamma_{F},U\stackrel{1-\theta}{\to}S_{12})&&\\ &\swarrow\,p_{1}&&\searrow\,p_{2}&\\ H^{1,0}(\Gamma_{F},U\stackrel{1-\theta}{\to}S_{1})&&&&H^{1,0}(\Gamma_{F},U\stackrel{1-\theta}{\to}S_{2})\\ \end{array}$$
  qui envoient $(V,\boldsymbol{\nu}_{12})$ respectivement sur $(V,\boldsymbol{\nu}_{1})$ et $(V,\boldsymbol{\nu}_{2})$. Il y a des homomorphismes duaux
  $$\begin{array}{ccccc}H^{1,0}(W_{F};\hat{S}_{1}\stackrel{1-\hat{\theta}}{\to}\hat{U})&&&&H^{1,0}(W_{F};\hat{S}_{2}\stackrel{1-\hat{\theta}}{\to}\hat{U})\\ &\nwarrow\,\hat{p}_{1}&&\nearrow\,\hat{p}_{2}&\\ &&H^{1,0}(W_{F};\hat{S}_{12}\stackrel{1-\hat{\theta}}{\to}\hat{U})&&\\ \end{array}$$
   D'apr\`es les propri\'et\'es de compatibilit\'e des produits de groupes de cohomologie, on a les \'egalit\'es
   $$<(V,\boldsymbol{\nu}_{1}),(\hat{V}_{1},{\bf s})>=<(V,\boldsymbol{\nu}_{12}),\hat{p}_{1}(\hat{V}_{1},{\bf s})>,$$
  
  $$<(V,\boldsymbol{\nu}_{2}),(\hat{V}_{2},{\bf s})>=<(V,\boldsymbol{\nu}_{12}),\hat{p}_{2}(\hat{V}_{2},{\bf s})>.$$
  En posant
  $$X=\boldsymbol{\Delta}_{2}(\delta_{2},\gamma;\underline{\delta}_{2},\underline{\gamma})\boldsymbol{\Delta}_{1}(\delta_{1},\gamma;\underline{\delta}_{1},\underline{\gamma})^{-1},$$
  on obtient
  $$X=<(V,\boldsymbol{\nu}_{12}),\hat{p}_{1}(\hat{V}_{1},{\bf s})\hat{p}_{2}(\hat{V}_{2},{\bf s})^{-1}>.$$
  Le tore $\hat{\mathfrak{T}}_{12}$ dual de $\mathfrak{T}_{12}$ est le quotient de $\hat{T}'_{1}\times \hat{T}'_{2}\times \hat{T}$ par le sous-groupe 
  $$\{(\hat{\xi}_{1}(t'_{1}),\hat{\xi}_{2}(t'_{2}),t'); t'_{1},t'_{2},t'\in \hat{T}', t'_{1}t'_{2}t'=1\}.$$
   Pour $w\in W_{F}$, notons encore $\zeta_{12}(w)$ l'image de $(\zeta_{1}(w),\zeta_{2}(w)^{-1},1)$ dans ce tore. Alors $\zeta_{12}$ est un cocycle. Le tore dual $\hat{S}_{12}$ de $S_{12}$ est le groupe des $(t,\underline{t},t_{sc})\in \hat{\mathfrak{T}}_{1}\times \hat{\mathfrak{T}}_{2}\times \hat{T}_{sc}$ tels que $t\underline{t}^{-1}=j(t_{sc})$, en g\'en\'eralisant la notation $j$ de 2.2. Notons $\hat{V}_{12}$ le cocycle $w\mapsto (\zeta_{12}(w),\zeta_{12}(w),1)$ de $W_{F}$ dans $\hat{S}_{12}$. On calcule $\hat{p}_{1}(\hat{V}_{1},{\bf s})\hat{p}_{2}(\hat{V}_{2},{\bf s})^{-1}$: c'est la classe de l'\'el\'ement $(\hat{V}_{12},1) \in Z^{1,0}(W_{F};\hat{S}_{12}\stackrel{1-\theta}{\to}\hat{U})$.  D'o\`u
  $$X=<(V,\boldsymbol{\nu}_{12}),(\hat{V}_{12},1)>.$$

  Introduisons le produit fibr\'e $T'_{12} $ de $T'_{1}$ et $T'_{2}$ au-dessus de $T'$, qui n'est autre que le commutant de $(\delta_{1},\delta_{2})$ dans $G'_{12}$. Introduisons le tore $\Sigma_{12}=(T'_{12}\times \underline{T}'_{12})/diag_{-}(Z(G'_{12}))$. Il y a un homomorphisme naturel $q:S_{12}\to \Sigma_{12}$. Dualement, on a $\hat{T}'_{12}=(\hat{T}'_{1}\times \hat{T}'_{2})/diag_{-}(\hat{T}')$ et 
  $$\hat{\Sigma}_{12}=\{(t,\underline{t},t'_{sc})\in \hat{T}'_{12}\times \underline{\hat{T}}'_{12}\times \hat{T}'_{sc}; j(t'_{sc} )=t\underline{t}^{-1}\},$$
  o\`u on note encore $j$ l'homomorphisme naturel et o\`u $\hat{T}'_{sc}$ est l'image r\'eciproque de $\hat{T}'$ dans $\hat{G}'_{SC}$. On a une suite d'homomorphismes
  $$Z(\hat{G}'_{12})\stackrel{diag}{\to}\hat{\Sigma}_{12}\stackrel{\hat{q}}{\to} \hat{S}_{12}.$$
  L'homomorphisme $\hat{q}$ prend ses valeurs dans le noyau de $1-\hat{\theta}$. Il y a donc un homomorphisme naturel
  $$H^1(W_{F},\hat{\Sigma}_{12})\to H^{1,0}(W_{F};\hat{S}_{12}\stackrel{1-\hat{\theta}}{\to}\hat{U}).$$
  L'\'el\'ement $\hat{V}_{12}$ est l'image  par cet homomorphisme de $diag(\zeta_{12})$. En vertu de la relation de compatibilit\'e [KS1] (A.3.13) (o\`u le signe n\'egatif dispara\^{\i}t d'apr\`es la correction [KS2] 4.3), on obtient
  $$X=<q(\boldsymbol{\nu}_{12}),diag(\zeta_{12})>,$$
  o\`u le produit est celui sur $H^0(\Gamma_{F};\Sigma_{12})\times H^1(W_{F};\hat{\Sigma}_{12})$.
  Le tore $\Sigma_{12}$ est un sous-tore maximal du groupe $\mathfrak{G}'_{12}=(G'_{12}\times G'_{12})/diag_{-}(Z(G'_{12}))$. L'homomorphisme $diag:Z(\hat{G}'_{12})\to \hat{\Sigma}_{12}$  se factorise en  
  $$Z(\hat{G}'_{12})\stackrel{\iota}{\to} Z(\hat{\mathfrak{G}}'_{12})\to \hat{\Sigma}_{12}.$$
  On se rappelle que tout \'el\'ement de $H^1(W_{F};Z(\hat{\mathfrak{G}}'_{12}))$ d\'efinit un caract\`ere de $\mathfrak{G}'_{12}(F)$. Donc 
  $$X=\omega_{12}(q(\boldsymbol{\nu}_{12})),$$
   o\`u $\omega_{12}$ est le caract\`ere de $\mathfrak{G}'_{12}(F)$ d\'efini par $\iota(\zeta_{12})$. Remarquons que 
   $$q(\boldsymbol{\nu}_{12})=((\delta_{1},\delta_{2}),(\underline{\delta}_{1}^{-1},\underline{\delta}_{2}^{-1})),$$
    en identifiant ce quadruplet \`a son image naturelle dans $\mathfrak{G}'_{12}(F)$. On peut d\'ecomposer 
  $$q(\boldsymbol{\nu}_{12})=((x_{1},x_{2}),(1,1))diag_{-}(\underline{\delta}_{1},\underline{\delta}_{2}).$$
  On a un homomorphisme
  $$G'_{12}\times G'_{12}\to \mathfrak{G}'_{12}.$$
  Par composition avec cet homomorphisme, $\omega_{12}$ d\'efinit un caract\`ere de $G'_{12}(F)\times G'_{12}(F)$. D'apr\`es les propri\'et\'es usuelles de compatibilit\'e, ce dernier caract\`ere est \'egal \`a   $\lambda_{12}\times \lambda_{12}$. D'o\`u $\omega_{12}((x_{1},x_{2}),(1,1))=\lambda_{12}(x_{1},x_{2})$. Pour achever la preuve du lemme, il reste \`a prouver que $\omega_{12}(diag_{-}(\underline{\delta}_{1},\underline{\delta}_{2}))=1$.  On peut dire que $\omega_{12}(diag_{-}(\underline{\delta}_{1},\underline{\delta}_{2}))$ est la valeur de notre quotient $X$ quand les triplets $(\delta_{1},\delta_{2},\gamma)$ et $(\underline{\delta}_{1},\underline{\delta}_{2},\underline{\gamma})$ sont \'egaux et qu'alors ce quotient vaut $1$ car, d'apr\`es [KS1] lemme 5.1.A, les deux termes $\boldsymbol{\Delta}_{1}(\delta_{1},\gamma;\underline{\delta}_{1},\underline{\gamma})$ et $\boldsymbol{\Delta}_{2}(\delta_{2},\gamma;\underline{\delta}_{2},\underline{\gamma})$ valent $1$. On peut dire aussi que $diag_{-}(\underline{\delta}_{1},\underline{\delta}_{2})$ appartient \`a l'image de l'homomorphisme naturel
  $$\underline{T}'_{12}/Z(G'_{12})\stackrel{diag_{-}}{\to} (\underline{T}'_{12}\times \underline{T}'_{12})/diag_{-}(Z(G'_{12}))$$
  Or, d'apr\`es sa construction, $\iota(\zeta_{12})$ est annul\'e par l'homomorphisme dual. $\square$

   Gr\^ace \`a ce lemme, il existe une unique fonction $\tilde{\lambda}_{12}$ sur $\tilde{G}'_{12}(F)$ telle que
  
  (i) pour $(\delta_{1},\delta_{2})\in \tilde{G}'_{12}(F)$ et $(x_{1},x_{2})\in G'_{12}(F)$, $\tilde{\lambda}_{12}(x_{1}\delta_{1},x_{2}\delta_{2})=\lambda_{12}(x_{1},x_{2})\tilde{\lambda}_{12}(\delta_{1},\delta_{2})$ (on abr\'egera cette propri\'et\'e en disant que $\tilde{\lambda}_{12}$ se transforme selon le caract\`ere $\lambda_{12}$); 
  
  (ii) pour $(\delta_{1},\gamma)\in {\cal D}_{1}$ et $\delta_{2}\in \tilde{G}'_{2}(F)$ tel que $(\delta_{1},\delta_{2})\in \tilde{G}'_{12}(F)$, $\Delta_{2}(\delta_{2},\gamma)=\tilde{\lambda}_{12}(\delta_{1},\delta_{2})\Delta_{1}(\delta_{1},\gamma)$.
  
  On d\'efinit une application lin\'eaire
  $$\begin{array}{ccc}C_{c,\lambda_{1}}^{\infty}(\tilde{G}'_{1}(F))&\to&C_{c,\lambda_{2}}^{\infty}(\tilde{G}'_{2}(F))\\ f_{1}&\mapsto&f_{2}\\ \end{array}$$
  par $f_{2}(\delta_{2})=\tilde{\lambda}_{12}(\delta_{1},\delta_{2})f_{1}(\delta_{1})$, o\`u $\delta_{1}$ est n'importe quel \'el\'ement tel que $(\delta_{1},\delta_{2})\in \tilde{G}'_{12}(F)$. C'est un isomorphisme qui se descend en un isomorphisme de $SI_{\lambda_{1}}(\tilde{G}'_{1}(F))$ sur $SI_{\lambda_{2}}(\tilde{G}'_{2}(F))$. Le diagramme 
  $$\begin{array}{ccccc}&& I(\tilde{G}(F),\omega)&&\\ &\swarrow&&\searrow&\\SI_{\lambda_{1}}(\tilde{G}'_{1}(F))&&\simeq&&SI_{\lambda_{2}}(\tilde{G}'_{2}(F))\\ \end{array}$$
  est commutatif, o\`u les deux fl\`eches descendantes sont les transferts.
  
   On a envie de d\'efinir $C_{c}^{\infty}({\bf G}')$ et $SI({\bf G}')$ comme les limites inductives des $C_{c,\lambda_{1}}^{\infty}(\tilde{G}'_{1}(F))$, resp. $SI_{\lambda_{1}}(\tilde{G}'_{1}(F))$, la limite \'etant prise sur toutes les donn\'ees $G'_{1}$, ...,$\Delta_{1}$, les applications de transition \'etant celles que l'on vient de d\'efinir. Alors le transfert devient une application lin\'eaire
  $$I(\tilde{G}(F),\omega)\otimes Mes(G(F))\to SI({\bf G}')\otimes Mes(G'(F)),$$
  qui ne d\'epend plus d'aucune donn\'ee auxiliaire. La construction pose un probl\`eme de logique car nos donn\'ees auxiliaires ne forment pas un ensemble: l'ensemble des groupes n'existe pas. Il y a plusieurs moyens de r\'esoudre cette difficult\'e. L'un, que l'on se contentera d'esquisser,  consiste \`a fixer un ensemble de couples $(G,\tilde{G})$ v\'erifiant les hypoth\`eses de 1.5, stable par quelques op\'erations \'el\'ementaires (le produit de deux couples de l'ensemble appartient \`a l'ensemble, un sous-objet d'un \'el\'ement de l'ensemble appartient \`a l'ensemble...) et tel que, pour  tout couple v\'erifiant les hypoth\`eses de 1.5, il existe un couple isomorphe appartenant \`a l'ensemble. Un tel ensemble existe puisque pour tout entier $n$, il n'y a qu'un nombre fini de classes d'isomorphisme de couples tels que $dim(G)=n$. On se limite ensuite \`a ne consid\'erer que des couples appartenant \`a l'ensemble fix\'e. Un autre moyen plus simple pour r\'esoudre le probl\`eme est de dire qu'une fois fix\'e le groupe $G$ et l'espace tordu $\tilde{G}$,  les donn\'ees ${\bf G}'$ que l'on rencontrera au cours de notre travail seront sinon en nombre fini, du moins d\'eduites des donn\'ees initiales par un nombre fini d'op\'erations. Elles restent dans un ensemble. On peut donc pour chacune d'elles fixer arbitrairement des donn\'ees auxiliaires $G'_{1}$,...,$\Delta_{1}$ et d\'efinir $C_{c}^{\infty}({\bf G}')$ et $SI({\bf G}')$  comme \'etant les espaces $C_{c,\lambda_{1}}^{\infty}(\tilde{G}'_{1}(F))$, resp. $SI_{\lambda_{1}}(\tilde{G}'_{1}(F))$ pour ces donn\'ees particuli\`eres. L'important est que, quand interviendront d'autres donn\'ees auxiliaires, on identifiera les espaces associ\'es \`a ces donn\'ees \`a $C_{c}^{\infty}({\bf G}')$ et $SI({\bf G}')$   par les isomorphismes d\'efinis ci-dessus. 
  
 Remarquons que les notions suivantes ont un sens:
 
 - le support dans $\tilde{G}'(F)$ d'un \'el\'ement de $C_{c}^{\infty}({\bf G}')$: on r\'ealise cet \'el\'ement dans un espace $C_{c,\lambda_{1}}^{\infty}(\tilde{G}'_{1})$; la projection dans $\tilde{G}'(F)$ de son support ne d\'epend pas des donn\'ees auxiliaires;
 
 - la multiplication d'un \'el\'ement de $C_{c}^{\infty}({\bf G}')$ par une fonction lisse sur $\tilde{G}'(F)$ (par le m\^eme argument).
 
 {\bf Cas particulier.} Supposons $(G,\tilde{G},{\bf a})$ quasi-d\'eploy\'e et \`a torsion int\'erieure, cf. 1.7. On dispose de la donn\'ee endoscopique maximale ${\bf G}=(G,{^LG},\tilde{s}=1)$. Pour cette donn\'ee, on peut choisir pour donn\'ees auxiliaires $G'_{1}=G$, $\tilde{G}'_{1}=\tilde{G}$ et $\Delta_{1}$ valant $1$ sur les couples qui se correspondent. Les espaces $C_{c}^{\infty}({\bf G})$ et $SI({\bf G})$ sont simplement $C_{c}^{\infty}(\tilde{G}(F))$ et $SI(\tilde{G}(F))$. 
 
   \bigskip
    \subsection{Action de groupes d'automorphismes}
  Soient ${\bf G}'=(G',{\cal G}',\tilde{s})$ une donn\'ee endoscopique relevante, $\underline{{\bf G}}'=(\underline{G}',\underline{{\cal G}}',\underline{\tilde{s}})$ une donn\'ee \'equivalente et $x\in \hat{G}$ d\'efinissant  l'\'equivalence.    Soit $\alpha_{x}:G'\to \underline{G}'$ un isomorphisme associ\'e \`a $x$, cf. 1.5. Remarquons que le diagramme
  $$\begin{array}{ccccc}&& Z(G)&&\\ &\swarrow&&\searrow&\\ Z(G')&&\stackrel{\alpha_{x}}{\to}&&Z(\underline{G}')\\ \end{array}$$
  est commutatif, donc de $\alpha_{x}$ se d\'eduit un isomorphisme $\tilde{\alpha}_{x}:\tilde{G}'=G'\times_{{\cal {Z}}(G)}{\cal Z}(\tilde{G})\to \underline{\tilde{G}}'=\underline{G}'\times_{{\cal Z}(G)}{\cal Z}(\tilde{G})$. 
  
  Fixons des donn\'ees auxiliaires $G'_{1}$,...,$\Delta_{1}$ relatives \`a la premi\`ere donn\'ee. On pose $\underline{G}'_{1}=G'_{1}$,   $\underline{C}_{1}=C_{1}$, avec pour application $\underline{G}'_{1}\to \underline{G}'$ la compos\'ee de $G'_{1}\to G'$ et de $\alpha_{x}:G'\to \underline{G}'$. On pose $\underline{\tilde{G}}'_{1}=\tilde{G}'_{1}$, avec pour application $\underline{\tilde{G}}'_{1}\to \underline{\tilde{G}}'$ la compos\'ee de $\tilde{G}'_{1}\to \tilde{G}'$ et de $\tilde{\alpha}_{x}:\tilde{G}'\to \underline{\tilde{G}}'$. On pose $\underline{\hat{\xi}}_{1}=\hat{\xi}_{1}\circ ad_{x^{-1}}:\underline{{\cal G}}'\to {^L\underline{G}}'_{1}={^LG}'_{1}$. Ces donn\'ees v\'erifient les conditions requises relativement \`a la donn\'ee $\underline{{\bf G}}'$. On v\'erifie que les bifacteurs de transfert d\'eduits de ces deux s\'eries de donn\'ees co\"{\i}ncident.  Donc  la fonction $\underline{\Delta}_{1}=\Delta_{1}$ est encore un facteur de transfert pour ces donn\'ees auxiliaires. On a alors un isomorphisme
  $$C_{c}^{\infty}({\bf G}')\simeq C_{c,\lambda_{1}}^{\infty}(\tilde{G}'_{1}(F))=C_{c,\lambda_{1}}^{\infty}(\underline{\tilde{G}}'_{1}(F))\simeq C_{c}^{\infty}(\underline{{\bf G}}').$$
  On en d\'eduit un isomorphisme $SI({\bf G}')\simeq SI(\underline{{\bf G}}')$. Par construction, il est compatible au transfert, c'est-\`a-dire que le diagramme suivant est commutatif:
  $$\begin{array}{ccccc}&&I(\tilde{G}(F),\omega)&&\\ &\text{transfert}\swarrow&&\searrow\,\text{transfert}&\\ SI({\bf G}')&&\simeq&&SI(\underline{{\bf G}}')\\ \end{array}$$

  Dans le cas particulier o\`u ${\bf G}'=\underline{{\bf G}}'$, on peut identifier ${\cal E}^{_{'}*}$ \`a une paire de Borel \'epingl\'ee d\'efinie sur $F$ (puisque $G'$ est quasi-d\'eploy\'ee) puis  pr\'eciser $\alpha_{x}$ en imposant que cet automorphisme pr\'eserve  cette paire de Borel \'epingl\'ee. On obtient une action du groupe $Aut({\bf G}')$ sur $C_{c}^{\infty}({\bf G}')$. 
  
  {\bf Remarque.} Comme me l'a fait remarquer Chaudouard, cette action d\'epend du choix de   la paire de Borel \'epingl\'ee,  qui n'est d\'etermin\'e que modulo l'action de $G'_{AD}(F)$. L'action devient canonique dans deux cas:
  
  - quand on passe \`a un quotient o\`u cette action dispara\^{\i}t, par exemple l'action sur l'espace $SI({\bf G}')$ est canonique;
  
  - si on se restreint aux $x$ pour lesquels $\alpha_{x}=1$.  
      
  On v\'erifie que le sous-groupe $Z(\hat{G})^{\Gamma_F}\hat{G}'$ de $Aut({\bf G}')$ agit trivialement. On a donc une action de $Aut({\bf G}')/\hat{G}'$ et en particulier de son sous-groupe $(Z(\hat{G})/(Z(\hat{G})\cap \hat{T}^{\hat{\theta},0}))^{\Gamma_F}$.  On a vu en 1.13 comment associer \`a un \'el\'ement $x$ de ce groupe un caract\`ere $\mu_{x}$ de $G_{0,ab}(F)$ et une fonction $\tilde{\mu}_{x}$ sur  $\tilde{G}_{0,ab}(F)$.

  \ass{Lemme}{Pour $x\in Z(\hat{G})$ tel que  $x(Z(\hat{G})\cap \hat{T}^{\hat{\theta},0})$ soit fixe par $\Gamma_{F}$, l'action de $x$ sur $C_{c}^{\infty}({\bf G}')$ est la multiplication par la fonction $\tilde{\mu}_{x}\circ N^{\tilde{G}',\tilde{G}}$.}

 Preuve. Fixons des donn\'ees auxiliaires $G'_{1}$,...,$\Delta_{1}$ dont on d\'eduit, \`a l'aide de $x$, de nouvelles donn\'ees comme ci-dessus. Mais, au lieu de les souligner, on note ces nouvelles donn\'ees $G'_{2}$,...,$\Delta_{2}$. En fait, ces donn\'ees sont les m\^emes que les premi\`eres, sauf $\hat{\xi}_{1}$ qui est remplac\'e par $\hat{\xi}_{2}=\hat{\xi}_{1}\circ ad_{x^{-1}}$. L'action de $x$ sur $C_{c,\lambda_{1}}^{\infty}(\tilde{G}'_{1}(F))$ est la compos\'ee de l'identit\'e de cet espace sur $C_{c}^{\infty}(\tilde{G}'_{2}(F),\lambda_{2})$ et de l'application de transition de ce deuxi\`eme espace sur le premier d\'efinie au paragraphe pr\'ec\'edent. Autrement dit, c'est la multiplication par la fonction $\delta_{1}\mapsto \tilde{\lambda}_{21}(\delta_{1},\delta_{1})$. Cette fonction se transforme selon le caract\`ere $g_{1}\mapsto \lambda_{21}(g_{1},g_{1})$ de $G'_{1}(F)$. Celui-ci est associ\'e au cocycle $w\mapsto (\zeta_{2}(w),\zeta_{1}(w)^{-1})$ de $W_{F}$ dans $Z(\hat{G}'_{12})$. Avec les notations de 2.2, on a $(\zeta_{2}(w),w)=\hat{\xi}_{2}(g(w),w)=\hat{\xi}_{1}(x^{-1}w(x)g(w),w)$, d'o\`u $\zeta_{2}(w)=\hat{\xi}_{1}(w(x)x^{-1})\zeta_{1}(w)$. Notre cocycle est donc le produit des deux cocycles $w\mapsto (\zeta_{1}(w),\zeta_{1}(w)^{-1})$ et $w\mapsto (\hat{\xi}_{1}(w(x)x^{-1}),1)$. On voit comme dans la preuve de 2.5 que le premier vaut $1$ sur la diagonale de $G'_{21}(F)$. Le deuxi\`eme d\'efinit le caract\`ere compos\'e de la projection de $G'_{21}(F)$ sur $G'(F)$ et du caract\`ere de ce dernier groupe associ\'e au cocycle $w\mapsto w(x)x^{-1}$. Ce dernier caract\`ere est le compos\'e de $\mu_{x}$ et de l'homomorphisme $G'(F)\to G'_{ab}(F)\to G'_{0,ab}(F)$. Cela d\'emontre que notre fonction $\delta_{1}\mapsto \tilde{\lambda}_{21}(\delta_{1},\delta_{1})$ se transforme selon le m\^eme caract\`ere que la fonction $\tilde{\mu}_{x}\circ N^{\tilde{G}',\tilde{G}}$ (ou plus exactement que cette fonction compos\'ee avec la projection $G'_{1}(F)\to G'(F)$). Pour  que ces deux fonctions soient \'egales, il suffit qu'elles le soient en un point. Puisque $\Delta_{2}=\Delta_{1}$ et que la multiplication par la fonction de transition envoie $\Delta_{2}$ sur $\Delta_{1}$, on a $\tilde{\lambda}_{21}(\delta_{1},\delta_{1})=1$ pour tout $\delta_{1}$ qui est composante d'un couple  $(\delta_{1},\gamma)\in {\cal D}_{1}$.  Pour un tel $\delta_{1}$, on a aussi $\tilde{\mu}_{x}\circ N^{\tilde{G},\tilde{G}}(\delta)=1$ d'apr\`es la d\'efinition de $\tilde{\mu}_{x}$ et la proposition 1.14(i). Or un tel $\delta_{1}$ existe puisque ${\bf G}'$ est relevante. $\square$
 
 \ass{Corollaire}{Un \'el\'ement de $C_{c}^{\infty}({\bf G}')$ est invariant par l'action de $(Z(\hat{G})/(Z(\hat{G})\cap\hat{T}^{\hat{\theta},0})^{\Gamma_{F}}$ si et seulement si son support est contenu dans l'ensemble des $\delta\in \tilde{G}'(F)$ tels que $N^{\tilde{G}',\tilde{G}}(\delta)$ appartienne \`a $N^{\tilde{G}}(\tilde{G}_{ab}(F))$.}

   \bigskip
  \subsection{Une propri\'et\'e de transformation du facteur de transfert} 
 Posons $G_{\sharp}=G/Z(G)^{\theta}$. Le groupe $G_{\sharp}(F)$ agit par conjugaison sur $\tilde{G}(F)$.  Soit $(B,T)$ une paire de Borel de $G$. On a l'\'egalit\'e $Z(G)^{\theta}=Z(G)\cap T^{\theta}$, o\`u  $\theta$ d\'esigne  la restriction de $ad_{\gamma}$ \`a $T$ pour n'importe quel $\gamma\in \tilde{G}$ tel que $ad_{\gamma}$ conserve $(B,T)$. D'o\`u une suite exacte
 $$(1) \qquad 1\to T/Z(G)^{\theta}\to T/T^{\theta}\times T_{ad}\to T_{ad}/T_{ad}^{\theta}\to 1.$$
 La deuxi\`eme fl\`eche  est le produit des applications naturelles. La premi\`ere  est le produit de l'application naturelle $T/Z(G)^{\theta}\to T_{ad}$ et de l'inverse de l'application naturelle  $T/Z(G)^{\theta}\to T/T^{\theta}$. En identifiant $T/T^{\theta}$ \`a $(1-\theta)(T)$ par l'homomorphisme $1-\theta$ et en identifiant de m\^eme $T_{ad}/T_{ad}^{\theta}$ \`a $(1-\theta)(T_{ad})$, on obtient une suite exacte
 $$1\to T/Z(G)^{\theta}\to (1-\theta)(T)\times T_{ad}\to (1-\theta)(T_{ad})\to 1.$$
 Dualement, en fixant une paire de Borel \'epingl\'ee $\hat{\cal E}$ de $\hat{G}$ et en utilisant les notations de 1.4, un tore maximal $\hat{T}_{\sharp }$ de $\hat{G}_{\sharp}$ s'ins\`ere dans une suite exacte
 $$(2) \qquad 1\to \hat{T}_{sc}/\hat{T}_{sc}^{\hat{\theta}}\stackrel{(\pi,1-\hat{\theta})}{\to}\hat{T}/\hat{T}^{\hat{\theta},0}\times  \hat{T}_{sc}\to \hat{T}_{\sharp}\to 1.$$
Dualement \`a l'homomorphisme $T_{sc}\to T/Z(G)^{\theta}$, on dispose d'un homomorphisme $\hat{T}_{\sharp}\to \hat{T}_{ad}$. Puisque $Z(\hat{G}_{\sharp})$ est le noyau de  cet homomorphisme, on d\'eduit ais\'ement de la suite ci-dessus la suite exacte
$$1\to Z(\hat{G}_{SC})/Z(\hat{G}_{SC})^{\hat{\theta}}\stackrel{(\pi,1-\hat{\theta})}{\to}Z(\hat{G})/(Z(\hat{G})\cap \hat{ T}^{\hat{\theta},0})\times Z(\hat{G}_{SC})\to Z(\hat{G}_{\sharp})\to 1.$$

 Soit ${\bf G}'=(G',{\cal G}',\tilde{s})$ une donn\'ee endoscopique relevante pour $(G,\tilde{G},{\bf a})$. On suppose que $\hat{{\cal E}}$ est adapt\'ee \`a cette donn\'ee, cf. 1.5. En particulier $\tilde{s}=s\hat{\theta}$, avec $s\in \hat{T}$.  Pour $w\in W_{F}$, soit $g_{w}=(g(w),w)\in {\cal G}'$ tel que $ad_{g_{w}}$ agisse par $w_{G'}$ sur $\hat{G}'$. Choisissons $z(w)\in Z(\hat{G})$ et $g_{sc}(w)\in \hat{G}_{SC}$ tels que $g(w)=z(w)\pi(g_{sc}(w))$. Choisissons $s_{sc}\in \hat{G}_{SC}$ qui a m\^eme image que $s$ dans $\hat{G}_{AD}$. On d\'efinit $a_{sc}(w)\in \hat{G}_{SC}$ par
 $$s_{sc}\hat{\theta}(g_{sc}(w))w(s_{sc})^{-1}=a_{sc}(w)g_{sc}(w).$$
 On note $z_{\sharp}(w)$ l'image de $(z(w),a_{sc}(w))$ dans $Z(\hat{G}_{\sharp})$ par l'application de la suite ci-dessus. Ce terme est bien d\'efini et $z_{\sharp}$ est un cocycle, qui d\'efinit un caract\`ere $\omega_{\sharp}$ de $G_{\sharp}(F)$. 
 
 {\bf Attention:} le caract\`ere $\omega_{\sharp}$ d\'epend de la donn\'ee endoscopique.
 
 Soient  $G'_{1}$,...,$\Delta_{1}$ des donn\'ees auxiliaires. 
 
 \ass{Lemme}{Pour $(\delta_{1},\gamma)\in {\cal D}_{1}$ et $x\in G_{\sharp}(F)$, on a
 $$\Delta_{1}(\delta_{1},x^{-1}\gamma x)=\omega_{\sharp}(x)\Delta_{1}(\delta_{1},\gamma).$$}
 
 Preuve. Il s'agit de calculer $\boldsymbol{\Delta}_{1}(\delta_{1},x^{-1}\gamma x;\delta_{1},\gamma )$. Choisissons une d\'ecomposition $x=z\pi(x_{sc})$, avec $z\in Z(G)$ et $x_{sc}\in G_{SC}$. Reprenons les constructions de 2.2. Si $(\delta,B',T',B,T,x^{-1}\gamma x)$ est le diagramme relatif \`a $(\delta,x^{-1}\gamma x)$, on prend $(\delta,B',T',ad_{x}(B),ad_{x}(T),\gamma)$ pour diagramme relatif \`a $(\delta,\gamma)$ et $r=x_{sc}$. D'o\`u $u_{\cal E}(\sigma)=x_{sc}^{-1}u_{\underline{\cal E}}(\sigma)\sigma(x_{sc})$. On en d\'eduit
 $$V_{T}(\sigma)=x_{sc}^{-1}V_{\underline{T}}(\sigma)\sigma(x_{sc})=x_{sc}^{-1}V_{\underline{T}}(\sigma)x_{sc}{\bf x}_{sc}(\sigma)^{-1} ,$$
 o\`u on a pos\'e ${\bf x}_{sc}(\sigma)=\sigma(x_{sc})^{-1}x_{sc}$. On a aussi $x^{-1}\gamma x=z^{-1}x_{sc}^{-1}\underline{\nu} \underline{e}x_{sc}z=z^{-1}\theta(z)x_{sc}^{-1}\underline{\nu} x_{sc}e$, d'o\`u $\nu=z^{-1}\theta(z)x_{sc}^{-1}\underline{\nu} x_{sc}$.  Le couple $({\bf x}_{sc},z)$ appartient \`a $Z^{1,0}(\Gamma_{F};T_{sc}\to T/Z(G)^{\theta})$. On a le diagramme commutatif
 $$\begin{array}{ccc}T_{sc}&\to&T/Z(G)^{\theta}\\ \downarrow&&\downarrow\,1-\theta\\   U&\stackrel{1-\theta}{\to}&S_{1}\\ \end{array}$$
 d'o\`u un homomorphisme
 $$(3) \qquad H^{1,0}(\Gamma_{F};T_{sc}\to T/Z(G)^{\theta})\to H^{1,0}(\Gamma_{F}; U\stackrel{1-\theta}{\to}S_{1}).$$
 Le terme $(V,\boldsymbol{\nu}_{1})$ est le produit de l'inverse de l'image de $({\bf x}_{sc},z)$ par cet homomorphisme et d'un \'el\'ement $(V',\boldsymbol{\nu}'_{1})$ qu'il est facile de reconna\^{\i}tre: en identifiant $T$ \`a $\underline{T}$ par $ad_{x_{sc}}$, $(V',\boldsymbol{\nu}'_{1})$ est le cocycle associ\'e au quadruplet diagonal $(\delta_{1},\gamma; \delta_{1},\gamma)$. Du c\^ot\'e dual, la conjugaison par $x$ ne se voit pas et le cocycle $(\hat{V}_{1},{\bf s})$ est le m\^eme que celui associ\'e \`a cette paire diagonale. On a donc
 $$<(V',\boldsymbol{\nu}'_{1}),(\hat{V}_{1},{\bf s})>=\boldsymbol{\Delta}_{1}(\delta_{1},\gamma;\delta_{1},\gamma)^{-1}=1.$$
 Donc $\boldsymbol{\Delta}_{1}(\delta_{1},x^{-1}\gamma x;\delta_{1},\gamma )$ est le produit de $({\bf x}_{sc},z)$ et de l'image   dans $H^{1,0}(W_{F};\hat{T}_{\sharp}\to \hat{T}_{ad})$ de $(\hat{V}_{1},{\bf s})$ par l'homomorphisme dual de (3). Par l'homomorphisme
 $$\hat{T}\to \hat{T}/\hat{T}^{\hat{\theta},0},$$
 le cocycle $t_{T}$ d\'efinit un cocycle \`a valeurs dans $\hat{T}/\hat{T}^{\hat{\theta},0}$, que nous notons  $t'_{T}$. Le cocycle $((t'_{T})^{-1},1)$ \`a valeurs dans $\hat{T}/\hat{T}^{\hat{\theta},0}\times \hat{T}_{sc}$ se descend par la suite (2) en un cocycle $\hat{V}'_{1}$ \`a valeurs dans $\hat{T}_{\sharp}$. Notons $s_{ad}$ l'image de $s$ dans $\hat{T}_{ad}$. On voit que l'image de $(\hat{V}_{1},{\bf s})$ dans $H^{1,0}(W_{F}; \hat{T}_{\sharp}\to \hat{T}_{ad})$ est la classe du couple $(\hat{V}'_{1},s_{ad})$. 
 
 {\bf Remarque.} L'inversion de $t'_{T}$ provient du fait que, dans la suite (1), l'homomorphisme $T/Z(G)^{\theta}\to T/T^{\theta}$ est l'inverse de l'homomorphisme naturel.
 
 Notons $t'_{T_{sc}}$ l'image de la cocha\^{\i}ne $t_{T_{sc}}$ par l'homomorphisme
 $$\hat{T}_{sc}\to\hat{T}_{sc}/\hat{T}_{sc}^{\hat{\theta}}.$$
 On ne change pas $\hat{V}'_{1}$ en multipliant la cocha\^{\i}ne $((t'_{T})^{-1},1)$ par l'image par le premier homomorphisme de la suite (2) de la cocha\^{\i}ne $t'_{T_{sc}}$, autrement dit en rempla\c{c}ant $((t'_{T})^{-1},1)$ par $((t'_{T})^{-1}\pi(t'_{T_{sc}}),(1-\hat{\theta})(t'_{T_{sc}}))$. On a
 $$t'_{T}(w)^{-1}\pi(t'_{T_{sc}}(w))=z(w).$$
 Les termes $\hat{r}_{T}(w)$, $\hat{n}(\omega_{T}(w)) $ et $\hat{r}_{T,G'}(w)$  sont invariants par $\hat{\theta}$. D'o\`u
 $$(1-\hat{\theta})(t'_{T_{sc}}(w))=\hat{\theta}(\hat{n}_{G'}(\omega_{T,G'}(w))g_{sc}(w))g_{sc}(w)^{-1}\hat{n}_{G'}(\omega_{T,G'}(w))^{-1} .$$
 On peut remplacer $\hat{\theta}$ par $ad_{s_{sc}^{-1}}\circ ad_{s_{sc}}\circ\hat{\theta}$. Or $ad_{s_{sc}}\circ\hat{\theta}$ fixe $\hat{n}_{G'}(\omega_{T,G'}(w))$ (par d\'efinition de $\hat{G}'$) et envoie $g_{sc}(w)$ sur $a_{sc}(w)g_{sc}(w)w_{G}(s_{sc})s_{sc}^{-1}$. D'o\`u
 $$(1-\hat{\theta})(t'_{T_{sc}}(w))=s_{sc}^{-1}\hat{n}_{G'}(\omega_{T,G'}(w))a_{sc}(w)g_{sc}(w)w_{G}(s_{sc})g_{sc}(w)^{-1}\hat{n}_{G'}(\omega_{T,G'}(w))^{-1}  .$$
 Le terme $a_{sc}(w)$ est central. Le compos\'e de la conjugaison par $\hat{n}_{G'}(\omega_{T,G'}(w))g_{sc}(w)$ et de l'op\'erateur $w_{G}$ n'est autre que l'op\'erateur $w_{T}$. On obtient
 $$(1-\hat{\theta})(t'_{T_{sc}}(w))=s_{sc}^{-1}w_{T}(s_{sc})a_{sc}(w).$$
Le cocycle $\hat{V}'_{1}$ est donc l'image naturelle de $z_{\sharp}$ par l'homomorphisme $Z(\hat{G}_{\sharp})\to \hat{T}_{\sharp}$ et de l'image naturelle du cocycle $w\mapsto s_{sc}^{-1}w_{T}(s_{sc})\in \hat{T}_{sc}$. Or le couple form\'e de cette image et de $s_{ad}$ est un cobord. Donc la classe de $(\hat{V}'_{1},s_{ad})$ est \'egale \`a l'image de $z_{\sharp}$ par l'homomorphisme 
$$H^1(W_{F};Z(\hat{G}_{\sharp}))\to H^{1,0}(W_{F}; \hat{T}_{\sharp}\to \hat{T}_{ad}).$$
D'autre part, le couple $({\bf x}_{sc},z)$ est  l'image naturelle de $x\in G_{\sharp}(F)$ par la suite d'applications
$$G_{\sharp}(F)\to G_{\sharp,ab}(F)\simeq H^{1,0}(\Gamma_{F};T_{sc}\to T/Z(G)^{\theta}).$$
D'apr\`es 1.13, le produit de $({\bf x}_{sc},z)$ et de  $(\hat{V}'_{1},s_{ad})$  est \'egal \`a $\omega_{\sharp}(x)$. $\square$

\bigskip
\subsection{Le cas $F={\mathbb R}$}
  On suppose $F={\mathbb R}$. Soit ${\bf G}'=(G',{\cal G}',\tilde{s})$ une donn\'ee endoscopique relevante de $(G,\tilde{G},{\bf a})$.  On fixe des donn\'ees auxiliaires $G'_{1}$,...,$\Delta_{1}$. Le groupe $W_{{\mathbb R}}$ contient $W_{{\mathbb C}}={\mathbb C}^{\times}$ comme sous-groupe d'indice $2$. Pour $w\in W_{{\mathbb C}}$, soit $g_{w}=(g(w),w)\in {\cal G}'$ tel que $ad_{g_{w}}$ agisse sur $\hat{G}'$ comme $w_{G'}$, c'est-\`a-dire par l'identit\'e. N\'ecessairement, $g(w)$ appartient \`a $\hat{ T}$. On a aussi $\hat{\xi}_{1}(g_{w})=(\zeta_{1}(w),w)$, avec $\zeta_{1}(w)\in Z(\hat{G}'_{1})\subset \hat{ T}'_{1}$. Notons $\hat{\mathfrak{T}}$ le quotient de $\hat{ T}'_{1}\times\hat{ T}$ par l'image de  $\hat{ T}^{\hat{\theta},0}$ plong\'e par $t\mapsto (\xi_{1}(t),t^{-1})$. On note $\rho(w)$ l'image de $(\zeta_{1}(w)^{-1},g(w))$ dans $\hat{\mathfrak{ T}}$. Cette image ne d\'epend pas du choix de $g_{w}$ et l'application $\rho$ ainsi d\'efinie est un homomorphisme continu de ${\mathbb C}^{\times}$ \`a valeurs dans $\hat{\mathfrak{ T}}$. Rappelons les propri\'et\'es suivantes, valables pour tout tore complexe $\hat{T}$. A tout \'el\'ement $b\in X_{*}(\hat{T})\otimes{\mathbb C}$ est associ\'e un homomorphisme du groupe multiplicatif ${\mathbb R}_{>0}$ dans $\hat{T}$: on \'ecrit $b=\sum_{i=1,...,n}s_{i}b_{i}$ avec des $b_{i}\in X_{*}(\hat{T})$ et des $s_{i}$ complexes; pour $x\in {\mathbb R}_{>0}$, on pose $b(x)=\prod_{i=1,...,n}b_{i}(x^{s_{i}})$. Si $\lambda$ est un homomorphisme continu de ${\mathbb C}^{\times}$ dans un tore complexe $\hat{T}$, il existe d'uniques $b_{\lambda}$, $b'_{\lambda}\in X_{*}(\hat{T})\otimes {\mathbb C}$ de sorte que $b_{\lambda}-b'_{\lambda}\in X_{*}(\hat{T})$ et $\lambda(w)=(b_{\lambda}-b'_{\lambda})(w)b'_{\lambda}(w\bar{w})$ pour tout $w\in {\mathbb C}^{\times}$. A notre homomorphisme $\rho$ sont ainsi associ\'es $b_{\rho}$ et $b'_{\rho}\in X_{*}(\hat{\mathfrak{ T}})\otimes {\mathbb C}$.  On a une suite exacte
$$0\to X_{*}(\hat{T}^{\hat{\theta},0})\otimes {\mathbb C}\stackrel{x\mapsto (\hat{\xi}_{1}(x),-x)}{\to}(X_{*}(\hat{ T}'_{1})\otimes{\mathbb C})\oplus(X_{*}(\hat{ T})\otimes{\mathbb C})\stackrel{\hat{p}}{\to} X_{*}(\hat{\mathfrak{T}})\otimes{\mathbb C}\to 0$$
L'espace $(X_{*}(\hat{ T}'_{1})\otimes{\mathbb C})\oplus(1-\hat{\theta})(X_{*}(\hat{ T})\otimes{\mathbb C})$ est un suppl\'ementaire du noyau de $\hat{p}$ et s'identifie donc \`a $X_{*}(\hat{\mathfrak{ T}})\otimes{\mathbb C}$. On peut consid\'erer que $b_{\rho}$ et $b'_{\rho}$ appartiennent \`a ce suppl\'ementaire et on pose simplement $b=b_{\rho}$. Montrons que
  
 (1) $b$ appartient \`a $(X_{*}(Z(\hat{ G}'_{1})^0)\otimes{\mathbb C})\oplus(1-\hat{\theta})(X_{*}(Z(\hat{G})^0)\otimes{\mathbb C}) $.

Preuve. Notons $b_{1}$ et $b_{2}$ les deux composantes de $b$. Soit $\alpha$ une racine de $\hat{ T}'_{1}$ dans $\hat{G}'_{1}$. On veut montrer que $<\alpha,b_{1}>=0$. La racine $\alpha$ se restreint (via $\hat{\xi}_{1}$) en une racine  de $\hat{ T}^{\hat{\theta},0}$ dans $\hat{G}'$, qui est la restriction d'une racine $\beta$ de $\hat{ T}$ dans $\hat{G}$. On d\'efinit $N\beta$ comme en 1.6 et on note $n(\beta)$ l'entier positif tel que la restriction de $N\beta$ \`a $\hat{T}^{\hat{\theta},0}$ co\"{\i}ncide avec celle de $n(\beta)\alpha$. L'\'el\'ement $(n(\beta)\alpha,N\beta)$ appartient \`a $X^*(\hat{\mathfrak{ T}})$. Parce que $N\beta$ est invariant par $\hat{\theta}$, on a  $<N\beta,b_{2}>=0$, d'o\`u l'\'egalit\'e $n(\beta)<\alpha,b_{1}>=<(n(\beta)\alpha,N\beta),b>$. Pour prouver que ce terme est nul, il suffit de prouver que $(n(\beta)\alpha,N\beta)\circ\rho(w)=1$ pour tout $w\in {\mathbb C}^{\times}$. Mais $\alpha(\zeta_{1}(w))=1$ parce que $\zeta_{1}(w)$ est central dans $\hat{G}'_{1}$ et $(N\beta)(g(w))=1$ parce que $\beta$ se restreint en une racine de $\hat{G}'$ et que $g_{w}$ agit par l'identit\'e sur ce groupe.  Cela prouve que $b_{1}$ est central. Notons $\rho'$ l'homomorphisme de ${\mathbb C}^{\times}$ dans $(1-\hat{\theta})(\hat{ T})$ d\'efini par $\rho'(w)=(1-\hat{\theta})(g(w))$. On a $(1-\hat{\theta})(b_{2})=b_{\rho'}$. On a la relation $s\hat{\theta}(g(w))w_{G}(s)^{-1}=a(w)g(w)$, o\`u $a$ est un  cocycle de $W_{{\mathbb R}}$ dans $Z(\hat{G})$, de classe ${\bf a}$. Ici, on se restreint \`a $w\in {\mathbb C}^{\times}$ donc $w_{G}=1$. De plus, $s$ commute \`a $g(w)\in \hat{T}$. L'\'egalit\'e pr\'ec\'edente se simplifie en $\rho'(w)=a(w)^{-1}$. L'application $a$, restreinte \`a ${\mathbb C}^{\times}$, est un homomorphisme continu  dont l'image est connexe, donc contenue dans $Z(\hat{G})^0$. On obtient $b_{a}=-b_{\rho'}=(\hat{\theta}-1)(b_{2})$. D'o\`u $b_{a}\in (1-\hat{\theta})(X_{*}(\hat{ T})\otimes{\mathbb C})\cap X{*}(Z(\hat{G})^0)\otimes{\mathbb C}$. La d\'ecomposition
$$X_{*}(\hat{ T})\otimes{\mathbb C}= (X_{*}(Z(\hat{G})^0)\otimes{\mathbb C})\oplus (X_{*}(\hat{T}_{sc})\otimes{\mathbb C})$$
est stable par $1-\hat{\theta}$ et cela entra\^{\i}ne que l'intersection pr\'ec\'edente est \'egale \`a $(1-\hat{\theta})(X_{*}(Z(\hat{G})^0)\otimes{\mathbb C})$ ou encore \`a $(1-\hat{\theta})^2(X_{*}(Z(\hat{G})^0)\otimes{\mathbb C})$. L'\'egalit\'e $b_{a}=(\hat{\theta}-1)(b_{2})$ et l'injectivit\'e de $(1-\hat{\theta})$ sur $(1-\hat{\theta})(X_{*}(\hat{ T})\otimes{\mathbb C})$ entra\^{\i}nent alors que $b_{2}$ appartient \`a $(1-\hat{\theta})(X_{*}(Z(\hat{G})^0)\otimes{\mathbb C}) $. $\square$

Soit $(\delta_{1},\gamma)\in {\cal D}_{1}$. On note $T'_{1}$ et $T'$ les commutants de $\delta_{1}$ dans $G'_{1}$ et de $\delta$ dans $G'$ et on note $T$ le commutant de $G_{\gamma}$ dans $G$. On a des projections
$$T'_{1}\to T'\leftarrow T$$
d\'efinies sur ${\mathbb R}$, d'o\`u des projections
$$ \mathfrak{t}'_{1}\to \mathfrak{t}'\leftarrow \mathfrak{t}$$
au niveau des alg\`ebres de Lie. L'\'el\'ement $b$ s'identifie \`a un \'el\'ement de  $\mathfrak{t}'_{1}({\mathbb C})^*\oplus \mathfrak{t}({\mathbb C})^*$. Soient $Y_{1}\in \mathfrak{t}'_{1}({\mathbb R})$ et $X\in\mathfrak{t}({\mathbb R})$ ayant m\^eme image $Y$ dans $\mathfrak{t}'({\mathbb R})$. Pour $\lambda\in {\mathbb R}$ assez proche de $0$, le couple $(exp(\lambda Y_{1})\delta_{1},exp(\lambda X)\gamma)$ appartient \`a ${\cal D}_{1}$. On dispose donc du facteur de transfert $\Delta_{1}( exp(\lambda Y_{1})\delta_{1},exp(\lambda X)\gamma)$.

\ass{Lemme}{La fonction
$$\lambda\mapsto \Delta_{1}( exp(\lambda Y_{1})\delta_{1},exp(\lambda X)\gamma)$$
est $C^{\infty}$ au voisinage de $0$. On a l'\'egalit\'e
$$\frac{d}{d\lambda}\Delta_{1}( exp(\lambda Y_{1})\delta_{1},exp(\lambda X)\gamma)_{\vert \lambda=0}=<b,Y_{1}\oplus X>\Delta_{1}( \delta_{1},\gamma).$$}

Preuve. Dans ces assertions, on peut remplacer $\Delta_{1}( exp(\lambda Y_{1})\delta_{1},exp(\lambda X)\gamma)$ par
 $$\boldsymbol{\Delta}_{1}(exp(\lambda Y_{1})\delta_{1},exp(\lambda X)\gamma;\delta_{1},\gamma).$$
  Reprenons les constructions de 2.2 pour calculer ce bifacteur. On ajoute un $\lambda$ dans les notations et on le supprime de nouveau pour noter les valeurs en $\lambda=0$. Par exemple, on note $\boldsymbol{\nu}_{1}(\lambda)$ le terme not\'e $\boldsymbol{\nu}_{1}$ en 2.2 et on pose $\boldsymbol{\nu}_{1}=\boldsymbol{\nu}_{1}(0)$. Dans la d\'efinition de $\boldsymbol{\Delta}_{imp}(exp(\lambda Y_{1})\delta_{1},exp(\lambda X)\gamma;\delta_{1},\gamma)$, le seul terme qui d\'epend vraiment de  
$\lambda$ est $\boldsymbol{\nu}_{1}(\lambda)$. Ce terme est le produit de $\boldsymbol{\nu}_{1}$ et de l'image de $(exp(\lambda Y_{1}),exp(\lambda X))\in \mathfrak{T}_{1}({\mathbb R})$ par l'homomorphisme naturel $\mathfrak{ T}_{1}({\mathbb R})\to S_{1}({\mathbb R})$.Posons simplement $Z=Y_{1}\oplus X\in X_{*}(\mathfrak{ T}_{1})\otimes{\mathbb C}$. Une propri\'et\'e de compatibilit\'e d\'ej\`a utilis\'ee entra\^{\i}ne alors l'\'egalit\'e
$$\boldsymbol{\Delta}_{imp}(exp(\lambda Y_{1})\delta_{1},exp(\lambda X)\gamma;\delta_{1},\gamma)=<\hat{V}_{\mathfrak{ T}_{1}},exp(\lambda Z)>^{-1}\boldsymbol{\Delta}_{imp}(\delta_{1},\gamma;\delta_{1},\gamma).$$
En fait, le dernier terme vaut $1$. Le cocycle $\hat{V}_{\mathfrak{ T}_{1}}$ d\'efinit un caract\`ere disons $\omega_{\mathfrak{T}_{1}}$ de $\mathfrak{ T}_{1}({\mathbb R})$. Par une propri\'et\'e g\'en\'erale, la restriction $\hat{V}_{\mathfrak{T}_{1},{\mathbb C}}$ de $\hat{V}_{\mathfrak{ T}_{1}}$ \`a $W_{{\mathbb C}}={\mathbb C}^{\times}$ d\'efinit le caract\`ere $\omega_{\mathfrak{T}_{1}}\circ Norm$ de $\mathfrak{ T}_{1} ({\mathbb C})$, o\`u $Norm:\mathfrak{ T}_{1}({\mathbb C})\to \mathfrak{ T}_{1}({\mathbb R})$ est la norme habituelle. On a $exp(\lambda  Z)=Norm(exp( \lambda Z/2))$ d'o\`u
 $$<\hat{V}_{\mathfrak{ T}_{1}},exp(\lambda  Z)>=<\hat{V}_{\mathfrak{T}_{1},{\mathbb C}},exp(\lambda Z/2) )>.$$
 Ce dernier terme est calcul\'e dans [Bor] 9.1. En notant simplement $b_{1}$ et $b'_{1}$ les \'el\'ements de $X_{*}(\hat{\mathfrak{ T}}_{1})\otimes{\mathbb C}=X^*(\mathfrak{ T}_{1})\otimes {\mathbb C}$ associ\'es \`a   $\hat{V}_{\mathfrak{ T}_{1},{\mathbb C}}$, on a
 $$<\hat{V}_{\mathfrak{ T}_{1},{\mathbb C}},exp(\lambda Z/2)>= exp(\lambda(<b_{1},Z>+<b'_{1},\bar{Z}>)/2),$$
 o\`u $Z\mapsto \bar{Z}$ est l'identit\'e sur $X_{*}(\mathfrak{ T}_{1})$ et la conjugaison complexe sur ${\mathbb C}$. Parce que $\hat{V}_{\mathfrak{ T}_{1},{\mathbb C}}$ est la restriction d'un cocycle d\'efini sur $W_{{\mathbb R}}$, on a $<b'_{1},\bar{Z}>=<b_{1},\sigma(Z)>$ o\`u $\sigma$ est le produit des deux conjugaisons complexes sur $X_{*}(\mathfrak{ T}_{1})$ et ${\mathbb C}$. Mais $Z$ est d\'efini sur ${\mathbb R}$ donc $\sigma(Z)=Z$ et le terme ci-dessus vaut simplement $exp(\lambda<b_{1},Z>)$. Calculons $b_{1}$. Pour $w\in W_{{\mathbb C}}$, les formules de 2.2 se simplifient: $\omega_{T}(w)=1$ et $\omega_{T,G'}(w)=1$. D'o\`u
 $$\hat{V}_{\mathfrak{T}_{1}}(w)=(\zeta_{1}(w),\hat{r}_{T}(w)g(w)^{-1}\hat{r}_{T,G'}(w)^{-1}).$$
 Cet homomorphisme est le produit de $\rho^{-1}$ et de l'image naturelle de l'homomorphisme $\rho'$ de ${\mathbb C}^{\times}$ dans $\hat{T}$ d\'efini par
 $$\rho'(w)=\hat{r}_{T}(w)\hat{r}_{T,G'}(w)^{-1}.$$
 On obtient ainsi
 $$(2) \qquad \boldsymbol{\Delta}_{imp}(exp(\lambda Y_{1})\delta_{1},exp(\lambda X)\gamma;\delta_{1},\gamma)=<\hat{V}_{\mathfrak{ T}_{1}},exp(\lambda Z)>^{-1}$$
 $$=exp(\lambda<b,Y_{1}\oplus X>)exp(-\lambda<b_{\rho'},X>).$$
 
 On va calculer $b_{\rho'}$. Pour d\'efinir  le bifacteur de transfert, on a d\^u fixer un sous-groupe de Borel $B$ contenant $T$, qui d\'etermine une positivit\'e sur $\Sigma(T)_{res,ind}$. Notons $\sigma$ la conjugaison complexe et notons $C$ le caract\`ere de ${\mathbb C}^{\times}$ d\'efini par $C(w)=\frac{w}{ \vert w\vert }$, o\`u ici $\vert w\vert =(\bar{w}w)^{1/2}$. On peut choisir nos $\chi$-data telles que, pour $\alpha_{res}\in \Sigma(T)_{res,ind}$, 
$$\chi_{\alpha_{res}}=\left\lbrace\begin{array}{cc}1,&\text{ si }\sigma\alpha_{res}\not=-\alpha_{res},\\ C,&\text{ si }\sigma\alpha_{res}=-\alpha_{res}\text{ et }\alpha_{res}>0,\\ C^{-1},&\text{ si }\sigma\alpha_{res}=-\alpha_{res}\text{ et }\alpha_{res}<0.\\ \end{array}\right.$$
Avec ces d\'efinitions, on voit que, pour $w\in W_{{\mathbb C}}$, on a
$$\hat{r}_{T}(w)=\prod_{\beta\in \Sigma(\hat{T})_{res,ind}; \sigma\beta=-\beta,\beta>0}\check{\beta}\circ C(w),$$
$$\hat{r}_{T,G'}(w)=\prod_{\beta\in \Sigma(\hat{T}');\sigma\beta=-\beta,\beta>0}\check{\beta}'\circ C(w).$$
Attention:  on a not\'e $\check{\beta}$ la coracine pour le groupe $\hat{G}^{\hat{\theta},0}$ associ\'ee \`a $\beta\in \Sigma(\hat{T})_{res,ind}$ et $\check{\beta}'$ celle pour le groupe $\hat{G}'$ associ\'ee \`a $\beta\in \Sigma(\hat{T}')$.   On d\'eduit de ces formules l'\'egalit\'e
$$b_{\rho'}=\frac{1}{2}\left(\sum_{\beta\in \Sigma(\hat{T})_{res,ind}; \sigma\beta=-\beta,\beta>0}\check{\beta}\right)-\frac{1}{2}\left(\sum_{\beta\in \Sigma(\hat{T}');\sigma\beta=-\beta,\beta>0}\check{\beta}'\right).$$
On doit identifier toutes ces coracines \`a des caract\`eres de $T$. Pour cela, on utilise 1.6. Soit $\alpha_{res}$ un \'el\'ement de $\Sigma(T)_{res,ind}$, qui est la restriction d'un \'el\'ement $\alpha\in \Sigma$ de type 1 ou 2 (puisque $\alpha_{res}$ est indivisible). Il lui est toujours associ\'e un \'el\'ement $\beta=(\hat{\alpha})_{res}\in \Sigma(\hat{T})_{res,ind}$ et $\check{\beta}$ s'identifie \`a $N\alpha$ si $\alpha$ est de type 1, \`a $2N\alpha$ si $\alpha$ est de type 2.
 Il est associ\'e \`a $\alpha_{res}$ un \'el\'ement $\beta\in \Sigma(\hat{T}')$ si $\alpha$ est de type 1 et $N\hat{\alpha}(s)=1$ ou si $\alpha$ est de type 2 et $N\hat{\alpha}(s)=\pm 1$ (si $N\hat{\alpha}(s)=-1$, $\beta$ est plus exactement associ\'e \`a $2\alpha_{res}$ qui est la restriction d'une racine de type 3). Alors $\check{\beta}$ s'identifie \`a
$$\left\lbrace\begin{array}{cc}N\alpha,&\text{ si }\alpha\text{ est de type 1 et }N\hat{\alpha}(s)=1;\\2N\alpha,&\text{ si }\alpha\text{ est de type 2 et }N\hat{\alpha}(s)=1;\\N\alpha,&\text{ si }\alpha\text{ est de type 2 et }N\hat{\alpha}(s)=-1.\\ \end{array}\right.$$

 Reprenons la classification en types (a), (b), (c) et (d) de la fin du paragraphe 2.2. Les formules ci-dessus conduisent \`a l'\'egalit\'e
$$(3) \qquad b_{\rho'}=\frac{1}{2}\left(\sum_{\alpha_{res}\in \Sigma_{\star} \text{ de type (a) ou (c)}}N\alpha\right)+\left(\sum_{\alpha_{res}\in \Sigma_{\star} \text{ de type (b)}}N\alpha\right),$$
o\`u on a not\'e $\Sigma_{\star}$ l'ensemble des $\alpha_{res}\in \Sigma(T)_{res,ind}$ tels que $\sigma\alpha_{res}=-\alpha_{res}$ et $\alpha_{res}>0$.

 D'apr\`es les d\'efinitions et notre choix de $\chi$-data, on a
$$\Delta_{II}(exp(\lambda Y)\delta,exp(\lambda X)\gamma)\Delta_{II}(\delta,\gamma)^{-1}=\prod_{\alpha_{res}\in \Sigma_{\star}}C(z_{\alpha_{res}}(\lambda)),$$
o\`u
$$z_{\alpha_{res}}(\lambda)=\left\lbrace\begin{array}{cc}\frac{(N\alpha)(\nu(\lambda))-1}{(N\alpha)(\nu)-1},&\text{ dans le cas (a),}\\ \frac{(N\alpha)(\nu(\lambda))^2-1}{(N\alpha)(\nu)^2-1},&\text{ dans le cas (b),}\\\frac{(N\alpha)(\nu(\lambda))+1}{(N\alpha)(\nu)+1},&\text{ dans le cas (c),}\\ 1,&\text{ dans le cas (d).} \end{array}\right.$$
Parce que $\gamma$ appartient \`a $\tilde{G}({\mathbb R})$, il r\'esulte de 1.3(4) que l'image de $\nu$ dans $T/(1-\theta)(T)Z(G)$ est fixe par $\sigma$. Pour $\alpha\in \Sigma_{\star}$, on a $\sigma(N\alpha)=-N\alpha$. Ces deux propri\'et\'es entra\^{\i}nent que $(N\alpha)(\nu)$ est un nombre complexe de module $1$. De m\^eme pour $(N\alpha)(\nu(\lambda))$.  Remarquons que $(N\alpha)(\nu(\lambda))=exp(\lambda<N\alpha,X>)(N\alpha)(\nu)$. Un calcul montre alors que pour $\lambda$ proche de $0$, on a
$$z_{\alpha_{res}}(\lambda)\in\left\lbrace\begin{array}{cc}exp(\lambda<N\alpha,X>/2){\mathbb R}_{>0},&\text{ dans les cas (a) et (c),}\\exp(\lambda<N\alpha,X>){\mathbb R}_{>0},&\text{ dans le cas (b),}\\ {\mathbb R}_{>0},&\text{ dans le cas (d).}\\ \end{array}\right.$$
En comparant avec (3), on en d\'eduit
$$\Delta_{II}(exp(\lambda Y)\delta,exp(\lambda X)\gamma)\Delta_{II}(\delta,\gamma)^{-1}=exp(\lambda<b_{\rho'},X>),$$
puis, gr\^ace \`a (2)
$$\boldsymbol{\Delta}_{1}(exp(\lambda Y_{1})\delta_{1},exp(\lambda X)\gamma;\delta_{1},\gamma)=exp(\lambda<b,X_{1}\oplus X>).$$
Le lemme r\'esulte de cette formule. $\square$

Soit $f\in C_{c}^{\infty}(\tilde{G}({\mathbb R}))$. On introduit la fonction $F_{f,\gamma}$ d\'efinie au voisinage de $\gamma$ dans $T({\mathbb R})\gamma$ par $F_{f,\gamma}( exp(X)\gamma)=[T^{\theta}({\mathbb R}):T^{\theta,0}({\mathbb R})]^{-1}I^{\tilde{G}}( exp(X)\gamma,\omega,f)$ et la fonction $F'_{f,\gamma}$ d\'efinie au voisinage de $\delta_{1}$ dans $T'_{1}({\mathbb R})\delta_{1}$ par 
$$F'_{f,\gamma}(exp(Y_{1})\delta_{1})=\Delta_{1}(exp(Y_{1})\delta_{1}, exp(X)\gamma)[T^{\theta}({\mathbb R}):T^{\theta,0}({\mathbb R})]^{-1}I^{\tilde{G}}( exp(X)\gamma,\omega,f),$$
 o\`u $X$ est n'importe quel \'el\'ement de $\mathfrak{t}({\mathbb R})$, assez petit, dont l'image dans $\mathfrak{t}'({\mathbb R})$ co\"{\i}ncide avec celle de $Y_{1}$ (l'expression ne d\'epend pas de $X$:  pour $X$ assez proche de $0$, la classe de conjugaison de $ exp(X)\gamma$  est d\'etermin\'ee par l'image de $X$ dans $\mathfrak{t}'({\mathbb R})$). La preuve du lemme montre que ces deux fonctions sont $C^{\infty}$ (rappelons que $\gamma$ est fortement r\'egulier). Le tore $T({\mathbb R})$ agit \`a gauche sur l'espace des fonctions sur $T({\mathbb R})\gamma$. Il s'en d\'eduit une action par op\'erateurs diff\'erentiels de l'alg\`ebre $Sym(\mathfrak{t}({\mathbb C}))$  sur l'espace des fonctions $C^{\infty}$ d\'efinies au voisinage de $\gamma$ dans $ T({\mathbb R})\gamma$. De m\^eme, on a une action de $Sym(\mathfrak{t}'_{1}({\mathbb C}))$ sur l'espace des fonctions $C^{\infty}$ d\'efinies au voisinage de $\delta_{1}$ dans $T'_{1}({\mathbb R})\delta_{1}$.
 
 {\bf Remarque.} Une abondante litt\'erature concernant les groupes r\'eels privil\'egie les actions \`a droite.  On pr\'ef\`ere les actions \`a gauche. On esp\`ere que cela ne cr\'eera pas trop de perturbations.
 
 \bigskip
 
  On a des homomorphismes
$$Sym(\mathfrak{t}'_{1}({\mathbb C}))\to Sym(\mathfrak{t}'({\mathbb C}))\leftarrow Sym(\mathfrak{t}({\mathbb C})).$$
On d\'efinit un automorphisme ${\bf b}$ de $Sym(\mathfrak{t}({\mathbb C}))$: c'est l'unique automorphisme tel que, pour $X\in \mathfrak{t}({\mathbb C})$, on ait ${\bf b}(X)=X+<b,X>$. On d\'efinit    un automorphisme ${\bf b}'_{1}$ de $Sym(\mathfrak{t}'_{1}({\mathbb C})) $: c'est l'unique automorphisme tel que, pour $Y_{1}\in \mathfrak{t}'_{1}({\mathbb C})$, on ait ${\bf b}'_{1}(Y_{1})=Y_{1}+<b,Y_{1}>$. Montrons que

(4) soient $U\in Sym(\mathfrak{t}({\mathbb C}))$ et $U'_{1}\in Sym(\mathfrak{t}'_{1}({\mathbb C}))$; supposons que $({\bf b}'_{1})^{-1}(U'_{1})$ et ${\bf b}(U)$ aient m\^eme image dans $Sym(\mathfrak{t}'({\mathbb C}))$; alors
$$U'_{1}F'_{f,\gamma}(\delta_{1})=\Delta_{1}(\delta_{1},\gamma)UF_{f,\gamma}(\gamma).$$

Preuve. Consid\'erons d'abord le cas o\`u $U'_{1}=Y_{1}+<b,Y_{1}>$ et $U=X-<b,X>$, avec $Y_{1}\in \mathfrak{t}'_{1}({\mathbb R})$ et $X\in \mathfrak{t}({\mathbb R})$ ayant m\^eme image dans $\mathfrak{t}'({\mathbb R})$. Dans ce cas, la relation cherch\'ee r\'esulte d'un simple calcul et du lemme pr\'ec\'edent. En fait, on obtient une relation plus g\'en\'erale: la fonction $U'_{1}F'_{f,\gamma}$ se d\'eduit de $UF_{f,\gamma}$ comme $F'_{f,\gamma}$ se d\'eduit de $F_{f,\gamma}$. Par r\'ecurrence, on obtient la m\^eme relation dans le cas o\`u $U'_{1}=U_{1}^{'(1)}...U_{1}^{'(n)}$ et $U=U^{(1)}...U^{(n)}$, si chaque couple $(U_{1}^{'(i)},U^{(i)})$ v\'erifie les conditions ci-dessus. En g\'en\'eral, on peut \'ecrire $(U'_{1},U)$ comme combinaison lin\'eaire de tels couples $(U_{1}^{'(1)}...U_{1}^{'(n)},U^{(1)}...U^{(n)})$ et d'un couple $(U'_{1},0)$. Il nous reste \`a traiter ce  cas. Supposons donc $U=0$. Alors $U'_{1}$ appartient \`a l'id\'eal engendr\'e par les ${\bf b}_{1}(Y_{1})$ o\`u $Y_{1}$ appartient au noyau de la projection $\mathfrak{t}'_{1}({\mathbb R})\to \mathfrak{t}'({\mathbb R})$. Il suffit de prouver que pour un tel $U'_{1}$, on a $U'_{1}F'_{f,\gamma}=0$. Or cela r\'esulte du premier cas trait\'e: il suffit de compl\'eter $Y_{1}$ en le couple $(Y_{1},X=0)$.   $\square$

Notons $\mathfrak{ Z}(G)$ le centre de l'alg\`ebre enveloppante de l'alg\`ebre de Lie de $G$. D'apr\`es Harish-Chandra, on a l'isomorphisme $\mathfrak{ Z}(G)\simeq Sym({\mathfrak t}({\mathbb C}))^{W}$.  On en d\'eduit des homomorphismes
$$(5) \qquad \mathfrak{ Z}(G'_{1})\simeq Sym({\mathfrak t}'_{1}({\mathbb C}))^{W'}\to \mathfrak{ Z}(G')\simeq Sym({\mathfrak t}'({\mathbb C}))^{W'}\leftarrow \mathfrak{ Z}(G)\simeq Sym({\mathfrak t}({\mathbb C}))^W.$$
Les automorphismes ${\bf b}$ et ${\bf b}'_{1}$ d\'efinis plus haut se restreignent en des automorphismes de $\mathfrak{ Z}(G)$ et $\mathfrak{ Z}(G'_{1})$: cela r\'esulte de (1). L'alg\`ebre $\mathfrak{ Z}(G)$ agit \`a gauche et \`a droite sur $C_{c}^{\infty}(\tilde{G}({\mathbb R}))$.  L'alg\`ebre $\mathfrak{ Z}(G'_{1})$ agit  \`a gauche et \`a droite sur $C_{c,\lambda_{1}}^{\infty}(\tilde{G}'_{1}({\mathbb R}))$. On consid\`ere les actions \`a gauche.  

\ass{Corollaire}{Soient $U'_{1}\in \mathfrak{ Z}(G'_{1})$, $U\in \mathfrak{Z}(G)$, $f\in C_{c}^{\infty}(\tilde{G}({\mathbb R}))$ et $f_{1}\in C_{c,\lambda_{1}}^{\infty}(\tilde{G}'_{1}({\mathbb R}))$. Supposons que $f_{1}$ soit un transfert de $f$ et que $({\bf b}'_{1})^{-1}(U'_{1})$ et ${\bf b}(U)$ aient m\^eme image dans $\mathfrak{ Z}(G')$ par les homomorphismes (4). Alors $U'_{1}f_{1}$ est un transfert de $Uf$.}

Preuve.  Soit $\delta_{1}\in \tilde{G}'_{1}({\mathbb R})$ un \'el\'ement fortement $\tilde{G}$-r\'egulier. On a 
$$I^{\tilde{G}}(\delta_{1},f)=d(\theta^*)^{1/2}\sum_{\gamma}\Delta_{1}(\delta_{1},\gamma)[Z_{G}(\gamma,{\mathbb R}):G_{\gamma}({\mathbb R})]^{-1}I^{\tilde{G}}(\gamma,\omega,f),$$
cf. 2.5. Pour chaque $\gamma$ dans l'ensemble de sommation, introduisons les fonctions $F_{f,\gamma}$ et $F'_{f,\gamma}$ comme plus haut. La formule ci-dessus se r\'ecrit
$$I^{\tilde{G}}(\delta_{1},f)=d(\theta^*)^{1/2}\sum_{\gamma}F'_{f,\gamma}(\delta_{1}).$$
 Pour tout $\gamma$, on a l'\'egalit\'e
$$F_{Uf,\gamma}=UF_{f,\gamma}.$$
 Ceci est un th\'eor\`eme d'Harish-Chandra dans le cas non tordu et on v\'erifie que la preuve s'\'etend dans notre cas. Cette relation jointe \`a (4) entra\^{\i}ne
 $$F'_{Uf,\gamma}=U'_{1}F'_{f,\gamma}.$$
 On en d\'eduit
 $$(6) \qquad I^{\tilde{G}}(\delta_{1},Uf)=U'_{1}I^{\tilde{G}}(\delta_{1},f),$$
 o\`u, \`a droite, on consid\`ere $I^{\tilde{G}}(\delta_{1},f)$ comme une fonction d\'efinie au voisinage de $\delta_{1}$ dans $T'_{1}({\mathbb R})\delta_{1}$, $T'_{1}$ \'etant comme pr\'ec\'edemment le commutant de $\delta_{1}$. 
 
  Une m\^eme relation vaut pour l'int\'egrale orbitale stable $S^{\tilde{G}'}(\delta_{1},f_{1})$. C'est en fait essentiellement le cas particulier o\`u $\tilde{G}=\tilde{G}'$. On obtient
 $$(7) \qquad S^{\tilde{G}'}(\delta_{1},U'_{1}f_{1})=U'_{1}S^{\tilde{G}'}(\delta_{1},f_{1}).$$
 Puisque $f_{1}$ est un transfert de $f$, les deux membres de droite de (6) et (7) sont \'egaux. Donc aussi les deux membres de gauche. Cette derni\`ere \'egalit\'e signifie que $U'_{1}f_{1}$ est un transfert de $Uf$ $\square$
 
 \bigskip
 
 \section{Levi et image du transfert}
 
 \bigskip
 
  \subsection{Espaces paraboliques, espaces de Levi}
 Appelons paire parabolique un couple $(P,M)$  form\'e d'un sous-groupe parabolique $P$ de $G$ et d'une composante de Levi $M$ de $P$. Provisoirement, on ne suppose pas que $P$ ou $M$ sont d\'efinis sur $F$. On note $\tilde{P}$ le normalisateur de $P$ dans $\tilde{G}$ ($\tilde{P}=\{\gamma\in \tilde{G}; ad_{\gamma}(P)=P\}$) et $\tilde{M}$ le normalisateur commun de $P$ et $M$. Si $\tilde{P}$ n'est pas vide, $\tilde{M}$ ne l'est pas non plus (si $P$ et $M$ sont d\'efinis sur $F$, on a mieux: $\tilde{P}(F)$ et $\tilde{M}(F)$ sont tous deux non vides). On dit alors que $\tilde{P}$ est un espace parabolique de $\tilde{G}$, que $\tilde{M}$ est un espace de Levi de $\tilde{G}$ et que $(\tilde{P},\tilde{M})$ est une paire parabolique de $\tilde{G}$. Remarquons que $\tilde{P}$ est uniquement d\'etermin\'e par $P$, mais $\tilde{M}$ n'est pas uniquement d\'etermin\'e par $M$. Toutefois, dans le cas particulier o\`u $\tilde{G}$ est \`a torsion int\'erieure, $\tilde{P}$ est toujours non vide et $\tilde{M}$ est uniquement d\'etermin\'e par $M$: c'est l'ensemble des $\gamma\in \tilde{G}$ tels que $ad_{\gamma}\in M/Z(G)$.

 {\bf Exemples.} Supposons $G=GL(3)$. Posons 
 $$J=\left(\begin{array}{ccc}0&0&1\\0&-1&0\\1&0&0\\ \end{array}\right).$$
 Notons $\theta^*$ l'automorphisme $g\mapsto J{^tg}^{-1}J$ de $G$ et posons  $\tilde{G}=G\theta^*$. 
 
 (1) Soit $P$ le sous-groupe parabolique triangulaire sup\'erieur \`a deux blocs de longueurs $2$ et $1$. Alors $\tilde{P}$ est vide. 
 
 (2) Soit $P$ le sous-groupe de Borel triangulaire sup\'erieur et $M$ le sous-groupe diagonal. Alors $\tilde{P}=P\theta^*$, $\tilde{M}=M\theta^*$. Soit $s$ un \'el\'ement du groupe de Weyl. Posons $P'=sPs^{-1}$ et $M'=sMs^{-1}=M$. Alors $\tilde{P}'=s\tilde{P}s^{-1}=sP\theta^*(s)^{-1}\theta^*$ et $\tilde{M}'=s\tilde{M}s^{-1}=sM\theta^*(s)^{-1}\theta^*$.  Si $\theta^*(s)\not=s$, on a $\tilde{M}'\not=\tilde{M}$. 
 
 (3)  Consid\'erons le groupe
 $$M=\{\left(\begin{array}{ccc}\star&0&\star\\ 0&\star&0\\ \star&0&\star\\ \end{array}\right)\}.$$
 C'est un Levi de $G$ qui est stable par $\theta^*$. Mais il n'y a aucun sous-groupe parabolique $P$ de $G$, de composante de Levi $M$, pour lequel $\tilde{P}$ soit non vide. 
 
 \bigskip
 
 Fixons une paire parabolique $(P_{0},M_{0})$ de $G$ d\'efinie sur $F$ et minimale. On d\'efinit comme ci-dessus les normalisateurs $\tilde{P}_{0}$ et $\tilde{M}_{0}$. Fixons une paire de Borel \'epingl\'ee ${\cal E}=(B,T,(E_{\alpha})_{\alpha\in \Delta})$ de $G$ telle que $T\subset M_{0}$ et $B\subset P_{0}$. Fixons $e\in Z(\tilde{G},{\cal E})$ (un tel \'el\'ement n'a pas de raison d'appartenir \`a $\tilde{G}(F)$). On a $\tilde{M}_{0}=M_{0}e$, $\tilde{P}_{0}=P_{0}e$, et $\tilde{M}_{0}(F)\not=\emptyset$. On introduit l'action galoisienne $\sigma\mapsto \sigma_{G^*}$ qui pr\'eserve la paire ${\cal E}$, pour laquelle $G$ devient quasi-d\'eploy\'e, cf. 1.2. Fixons une paire de Borel \'epingl\'ee $\hat{{\cal E}}=(\hat{B},\hat{T},(\hat{E}_{\alpha})_{\alpha\in \Delta})$ de $\hat{G}$. On modifie l'isomorphisme $^LG\simeq \hat{G}\rtimes W_{F}$ de sorte qu'elle devienne stable par l'action galoisienne et on fixe un \'el\'ement $\hat{\theta}$ relatif  \`a cette paire, cf. 1.4.

Rappelons qu'il y a des bijections naturelles entre les divers ensembles suivants:

- les classes de conjugaison de paires paraboliques de $G$;

- les  paires paraboliques de $G$  qui sont standard,  c'est-\`a-dire qu'elles contiennent $(B,T)$;

- les classes de conjugaison de paires paraboliques de $\hat{G}$;

- les  paires paraboliques de $\hat{G}$  qui sont standard,  c'est-\`a-dire qu'elles contiennent $(\hat{ B},\hat{ T})$.

Ces ensembles sont munis d'actions galoisiennes (sur le deuxi\`eme, c'est celle provenant de l'action quasi-d\'eploy\'ee $\sigma\mapsto \sigma_{G^*}$). Les bijections   sont \'equivariantes pour les actions galoisiennes. Celle entre paires standard transporte l'action de $\theta$ sur celle de $\hat{\theta}^{-1}$. Ainsi, la paire $(P_{0},M_{0})$ correspond \`a une paire standard $(\hat{P}_{0},\hat{M}_{0})$ qui est fix\'ee par les actions de $\Gamma_{F}$ et $\hat{\theta}$. Alors les bijections pr\'ec\'edentes induisent des bijections entre
  
  - les classes de conjugaison de paires paraboliques de $\tilde{G}$ d\'efinies sur $F$;
  
  - les  paires paraboliques standard de $G$ fix\'ees par $\Gamma_{F}$ et $\theta$ et qui contiennent $(P_{0},M_{0})$;
  
  - les paires paraboliques standard de $\hat{G}$ fix\'ees par $\Gamma_{F}$ et $\hat{\theta}$ et qui contiennent $(\hat{P}_{0},\hat{M}_{0})$.

Appelons sous-groupe parabolique de $^LG$ un sous-groupe ${\cal P}\subset {^LG}$ pour lequel la projection sur $W_{F}$ induit une suite exacte
  $$1\to \hat{P} \to {\cal P}\to W_{F}\to 1,$$
  o\`u $\hat{P}$ est un sous-groupe parabolique de $\hat{G}$.   Appelons composante de Levi d'un tel sous-groupe un sous-groupe ${\cal M}\subset {\cal P}$  pour lequel la projection sur $W_{F}$ induit une suite exacte
  $$1\to \hat{M} \to {\cal M}\to W_{F}\to 1,$$
  o\`u $\hat{M}$ est une composante de Levi de $\hat{P}$.  Remarquons que ${\cal P}$ est d\'etermin\'e par $\hat{P}$: c'est le normalisateur de $\hat{P}$ dans $^LG$. De m\^eme, ${\cal M}$ est d\'etermin\'e par $\hat{P}$ et $\hat{M}$. Pour de tels ${\cal P}$ et ${\cal M}$, notons $\tilde{{\cal P}}$ le normalisateur de ${\cal P}$ dans $^L\tilde{G}={^LG}\hat{\theta}$ et $\tilde{{\cal M}}$ le normalisateur commun de ${\cal P}$ et ${\cal M}$. Si $\tilde{{\cal P}}$ n'est pas vide, $\tilde{{\cal M}}$ ne l'est pas non plus et on appelle $\tilde{{\cal P}}$ un espace parabolique de $^L\tilde{G}$, $\tilde{{\cal M}}$ un espace de Levi de $^L\tilde{G}$ et $(\tilde{{\cal P}},\tilde{{\cal M}})$ une paire parabolique de $^L\tilde{G}$. Le groupe $\hat{G}$ agit par conjugaison sur l'ensemble de ces paires paraboliques.  Montrons que
  
  (4) l'ensemble des classes de conjugaison de paires paraboliques de $^L\tilde{G}$ est en bijection avec l'ensemble des paires paraboliques standard de $\hat{G}$ qui sont invariantes par $\Gamma_{F}$ et $\hat{\theta}$.
  
  Preuve. Soit $(\tilde{{\cal P}},\tilde{{\cal M}})$ une paire parabolique de $^L\tilde{G}$. Le groupe ${\cal P}$ est bien d\'etermin\'e:  c'est le sous-groupe des $x\in {^LG}$ tels que $x\tilde{{\cal P}}=\tilde{{\cal P}}$ . De m\^eme, le groupe ${\cal M}$ est bien d\'etermin\'e. Les groupes $\hat{P}$ et $\hat{M}$ sont bien d\'etermin\'es: ce sont les intersections de ${\cal P}$ et ${\cal M}$ avec $\hat{G}$. Il existe une unique paire parabolique standard $(\hat{P}',\hat{M}')$ de $\hat{G}$ qui est conjugu\'ee \`a $(\hat{P},\hat{M})$. Quitte \`a effectuer une conjugaison, on se ram\`ene au cas o\`u $(\hat{P},\hat{M})$ est elle-m\^eme standard. Soit $(g,w)\in {\cal P}$. Puisque ${\cal P}$ est un groupe, la conjugaison par $(g,w)$ conserve $\hat{P}$, autrement dit $gw(\hat{P})g^{-1}=\hat{P}$. Puisque $(\hat{ B},\hat{ T})$ est conserv\'e par $\Gamma_{F}$, $w(\hat{P})$ est encore standard. Deux sous-groupes paraboliques standard ne sont conjugu\'es que s'ils sont \'egaux. Donc $w(\hat{P})=\hat{P}$. L'\'egalit\'e $gw(\hat{P})g^{-1}=\hat{P}$ entra\^{\i}ne alors que $g\in \hat{P}$. Cela d\'emontre que $\hat{P}$ est conserv\'e par $\Gamma_{F}$ et que ${\cal P}=\hat{P}\rtimes W_{F}$. Fixons un \'el\'ement de l'ensemble $\tilde{{\cal P}}$, qui n'est pas vide. Quitte \`a le multiplier par un \'el\'ement de ${\cal P}$, on peut le supposer de la forme $g\hat{\theta}$, avec $g\in \hat{G}$. Cet \'el\'ement normalise ${\cal P}$, donc aussi son intersection $\hat{P}$ avec $\hat{G}$. De nouveau, parce que $\hat{\theta}(\hat{P})$ est standard, cela entra\^{\i}ne que $\hat{\theta}(\hat{P})=\hat{P}$, puis que $g\in \hat{P}$. Donc $\tilde{{\cal P}}=(\hat{P}\rtimes W_{F})\hat{\theta}$. Un raisonnement analogue vaut pour les composantes de Levi: $\hat{M}$ est n\'ecessairement stable par $\Gamma_{F}$ et $\hat{\theta}$ et on a $\tilde{{\cal M}}=(\hat{M}\rtimes W_{F})\hat{\theta}$. L'assertion (4) s'ensuit. $\square$ 

 Ainsi, les bijections pr\'ec\'edentes se prolongent en une injection de l'ensemble des classes  de conjugaison de paires paraboliques de $\tilde{G}$ d\'efinies sur $F$ dans celui des  classes de conjugaison de paires paraboliques de $^L\tilde{G}$. C'est une bijection si et seulement si $G$ est quasi-d\'eploy\'e.  Remarquons que, si la classe de $(\tilde{{\cal P}},\tilde{{\cal M}})$ correspond \`a celle de $(P,M)$ par cette application, le groupe ${\cal M}$ s'identifie \`a $^LM$ et $\tilde{{\cal M}}$ \`a $^L\tilde{M}$. Mais une telle identification n'est pas intrins\`eque aux deux ensembles $\tilde{M}$ et $\tilde{{\cal M}}$, elle d\'epend des paraboliques.

 On aura aussi besoin de consid\'erer des Levi ou sous-groupes paraboliques semi-standard. Pour un sous-groupe parabolique $P$ de $G$ semi-standard, c'est-\`a-dire contenant $T$, notons $\Sigma^P(T)$ l'ensemble des racines de $T$ dans l'alg\`ebre de Lie de $P$. De m\^eme, pour un Levi semi-standard $M$ de $G$, on d\'efinit l'ensemble $\Sigma^M(T)$. Pour un sous-groupe parabolique semi-standard $\hat{P}$ de $\hat{G}$, ou pour un Levi semi-standard $\hat{M}$, on d\'efinit de m\^eme les ensembles de racines $\Sigma^{\hat{P}}(\hat{T})$ et $\Sigma^{\hat{M}}(\hat{T})$. Montrons que
 
 (5) il y a une bijection $P\mapsto \hat{P}$ entre l'ensemble des sous-groupes paraboliques semi-standard de $G$ et celui des sous-groupes paraboliques semi-standard de $\hat{G}$ caract\'eris\'ee par l'\'egalit\'e $\Sigma^{\hat{P}}(\hat{T})=\{\hat{\alpha}; \alpha\in \Sigma^P(T)\}$; 
 
  (6) il y a une bijection $M\mapsto \hat{M}$ entre l'ensemble des  Levi semi-standard de $G$ et  celui des Levi semi-standard de $\hat{G}$ caract\'eris\'ee par l'\'egalit\'e $\Sigma^{\hat{M}}(\hat{T})=\{\hat{\alpha}; \alpha\in \Sigma^M(T)\}$.
  
  Preuve. L'application $P\mapsto \Sigma^P(T)$ est une bijection entre l'ensemble des sous-groupes paraboliques semi-standard de $G$ et  l'ensemble des sous-ensembles $\Pi\subset \Sigma(T)$ v\'erifiant les deux propri\'et\'es:
  
  (7) $\Pi\cup(-\Pi)=\Sigma$;
  
  (8) si $\alpha, \beta\in \Pi$ sont tels que $\alpha+\beta\in \Sigma(T)$, alors $\alpha+\beta\in \Pi$.
  
  On a une assertion analogue du c\^ot\'e dual. D'autre part, on peut identifier $\Sigma(\hat{T})$ \`a l'ensemble de coracines $\check{\Sigma}(T)$. Pour prouver (5), il suffit de prouver que l'application $\Pi\mapsto \check{\Pi}=\{\check{\alpha}; \alpha\in \Pi\}$ \'echange les conditions (7) et (8) avec les analogues pour l'ensemble $\check{\Sigma}(T)$. Evidemment, si $\Pi$ v\'erifie (7), $\check{\Pi}$ v\'erifie la condition analogue. Soient   $\check{\alpha}$, $\check{\beta}\in \check{\Pi}$, supposons $\check{\alpha}+\check{\beta}\in \check{\Sigma}(T)$. Alors $\alpha$ et $\beta$ forment une base d'un syst\`eme de racines irr\'eductible de rang $2$. D'apr\`es (8), les \'el\'ements positifs de ce syst\`eme appartiennent \`a $\Pi$. En inspectant les trois syst\`emes de racines possibles de rang $2$, on v\'erifie que $\check{\alpha}+\check{\beta}$ est toujours la coracine d'un \'el\'ement positif de ce syst\`eme. Donc $\check{\alpha}+\check{\beta}$ appartient \`a $\check{\Pi}$. Cela prouve (5). Les ensembles $\Sigma^M(T)$ pour $M$ semi-standard sont exactement ceux de la forme $\Pi\cap (-\Pi)$, pour $\Pi$ v\'erifiant (7) et (8). Alors la m\^eme preuve s'applique \`a (6). $\square$  
  
   {\bf Changement de terminologie.}  Dor\'enavant, on appellera "sous-groupe parabolique"  de $G$ ou "espace parabolique" de $\tilde{G}$ de tels objets d\'efinis sur $F$. On appellera "Levi" de $G$  une composante de Levi  d\'efinie sur $F$ d'un sous-groupe parabolique d\'efini sur $F$ et on appellera "espace de Levi" de $\tilde{G}$ une composante de Levi d\'efinie sur $F$ d'un espace parabolique de $\tilde{G}$ d\'efini sur $F$. 
   
   On utilisera les notations d'Arthur concernant ces objets. Par exemple, pour un espace de Levi $\tilde{M}$ de $\tilde{G}$, on note ${\cal L}(\tilde{M})$ l'ensemble des espaces de Levi $\tilde{L}$ contenant $\tilde{M}$. On utilise des notations analogues pour les groupes et espaces duaux. Notons $({\cal P}_{0},{\cal M}_{0}) $ la paire parabolique de $^L\tilde{G}$ issue de $(\hat{P}_{0},\hat{M}_{0})$, c'est-\`a-dire ${\cal P}_{0}=(\hat{P}_{0}\rtimes W_{F})\hat{\theta}$, ${\cal M}_{0}=(\hat{M}_{0}\rtimes W_{F})\hat{\theta}$. Alors
   
   (9) il y a une bijection $\tilde{M}\mapsto {\cal M}$ de ${\cal L}(\tilde{M}_{0})$ sur ${\cal L}({\cal M}_{0})$ caract\'eris\'ee ainsi: si $M$ est le Levi sous-jacent \`a $\tilde{M}$ et $\hat{M}$ le Levi sous-jacent \`a ${\cal M}$, $M$ s'envoie sur $\hat{M}$ par la bijection (6).
   
   C'est \'evident puisque la bijection (6) est \'equivariante pour les actions galoisiennes et \'echange les actions de $\theta$ et de $\hat{\theta}^{-1}$.
   
    Soient $\tilde{M}$ et $\underline{\tilde{M}}$ deux espaces de Levi de $\tilde{G}$ et soient $\tilde{{\cal M}}$ et $\underline{\tilde{\cal M}}$ deux Levi de $^L\tilde{G}$. On suppose que $\tilde{{\cal M}}$ et  $\underline{\tilde{\cal M}}$ s'identifient \`a $^L\tilde{M}$ et $^L\tilde{\underline{M}}$ gr\^ace \`a des choix de paraboliques comme plus haut et on fixe de telles identifications. Notons 
    $$W(\tilde{M},\underline{\tilde{M}})=\{g\in G(F); ad_{g}(\tilde{M})=\underline{\tilde{M}}\}/M(F),$$
    $$W(\tilde{{\cal M}},\underline{\tilde{{\cal M}}})=\{x\in \hat{G}; ad_{x}(\tilde{{\cal M}})=\underline{\tilde{{\cal M}}}\}/\hat{M}.$$
    Alors
    
    (10) il y a une bijection naturelle entre   $W(\tilde{M},\underline{\tilde{M}})$ et $W(\tilde{{\cal M}},\underline{\tilde{{\cal M}}})$.
    
    Preuve. En oubliant les choix faits pr\'ec\'edemment, on fixe maintenant des paires de Borel \'epingl\'ees dans chacun des groupes intervenant, dont on note les tores $T$, $\underline{T}$, $\hat{T}$, $\underline{\hat{T}}$. On normalise les actions de $\Gamma_{F}$ sur $\hat{M}$ et $\underline{\hat{M}}$ de sorte qu'elles pr\'eservent les paires de Borel \'epingl\'ees. On choisit $\hat{\theta}\in \tilde{{\cal M}}\cap \hat{G}\boldsymbol{\hat{\theta}}$ et $\underline{\hat{\theta}}\in \underline{\tilde{{\cal M}}}\cap\hat{G}\boldsymbol{\hat{\theta}}$ qui pr\'eservent aussi ces paires. De m\^eme, on choisit $e\in \tilde{M}$ et $\underline{e}\in \underline{\tilde{M}}$ pr\'eservant les paires de Borel \'epingl\'ees et on introduit les actions galoisiennes quasi-d\'eploy\'ees $\sigma\mapsto \sigma_{M^*}$ et $\sigma\mapsto \sigma_{\underline{M}^*}$ qui prolongent aux groupes $M$ et $\underline{M}$ celles de 1.2. Puisque l'on a fix\'e des identifications de $\tilde{{\cal M}}$ et  $\underline{\tilde{\cal M}}$  \`a $^L\tilde{M}$ et $^L\tilde{\underline{M}}$, les tores $\hat{T}$ et $\underline{\hat{T}}$ s'identifient aux duaux de $T$ et $\underline{T}$. Ces identifications sont \'equivariantes pour les actions galoisiennes et transportent $\hat{\theta}$ et $\underline{\hat{\theta}}$ en les inverses de $\theta=ad_{e}$ et $\underline{\theta}=ad_{\underline{e}}$. Soit $x\in\hat{G}$ tel que $ad_{x}(\tilde{{\cal M}})=\underline{\tilde{{\cal M}}}$. Quitte \`a multiplier $x$ \`a droite par un \'el\'ement de $\hat{M}$, on peut supposer que $ad_{x}$ transporte la paire de $\hat{M}$ sur celle de $\hat{\underline{M}}$, donc $\hat{T}$ sur $\hat{\underline{T}}$. Par dualit\'e puis inversion, il s'en d\'eduit un isomorphisme $\iota:T\to \underline{T}$. Celui-ci est la restriction d'un automorphisme $ad_{g}$ pour un $g\in G$. En effet nos identifications sont issues de choix de paraboliques. A conjugaison pr\`es, on peut les supposer tous standard, pour des paires de Borel fix\'ees de $G$ et $\hat{G}$. Alors  $\underline{\hat{T}}$ devient \'egal \`a $\hat{T}$, l'isomorphisme $ad_{x}$ de ce tore est un \'el\'ement du groupe de Weyl de $\hat{G}$ et $\iota$ est l'\'el\'ement du groupe de Weyl de $G$ qui lui correspond. Soit donc $g\in G$ tel  que $\iota$ soit la restriction de $ad_{g}$ \`a $T$. La d\'efinition de $x$ entra\^{\i}ne que $ad_{x}$ envoie $\hat{M}$ sur $\underline{\hat{M}}$, qu'il est \'equivariant pour les actions galoisiennes et transporte $\hat{\theta}$ sur $\underline{\hat{\theta}}$ (ces \'el\'ements \'etant vus ici comme des automorphismes de $\hat{G}$). Par dualit\'e, $ad_{g}$ envoie $M$ sur $\underline{M}$, est \'equivariant pour les actions galoisiennes quasi-d\'eploy\'ees et transporte $\theta$ sur $\underline{\theta}$. Cette derni\`ere condition implique que $ad_{g}(\tilde{M})=\underline{\tilde{M}}$. Parce que les actions galoisiennes naturelles ne diff\`erent des actions quasi-d\'eploy\'ees que par des automorphismes int\'erieurs, la condition d'\'equivariance entra\^{\i}ne que la classe $gM$ est fixe par $\Gamma_{F}$ dans $G/M$. Or $(G/M)(F)=G(F)/M(F)$.  Quitte \`a multiplier $g$ \`a droite par un \'el\'ement de $M$, on peut supposer que $g\in G(F)$. Alors $gM(F)\in W(\tilde{M},\underline{\tilde{M}})$. Evidemment, cette classe ne d\'epend que de la classe $x\hat{M}$ et on a ainsi d\'efini une application de $W(\tilde{{\cal M}},\underline{\tilde{{\cal M}}})$ dans $W(\tilde{M},\underline{\tilde{M}})$. On v\'erifie qu'elle ne d\'epend pas des choix de paires de Borel \'epingl\'ees. On d\'efinit l'application r\'eciproque de fa\c{c}on analogue. Cela prouve (10). $\square$
   
   Les propri\'et\'es suivantes sont utiles:
  
  (11) soit $T\subset G$ un tore d\'efini et d\'eploy\'e sur $F$; notons $Z_{\tilde{G}}(T)$ l'ensemble des $\gamma\in \tilde{G}$ tels que $ad_{\gamma}(t)=t$ pour tout $t\in T$; si cet ensemble n'est pas vide, c'est un espace de Levi de $\tilde{G}$;
  
  (12) soit $\tilde{M}$ un espace de Levi de $\tilde{G}$; alors $\tilde{M}=Z_{\tilde{G}}(A_{\tilde{M}})$.
  
  Preuve. Cela est bien connu dans le cas non tordu o\`u $\tilde{G}=G$. Dans la situation de (11), le commutant $M$ de $T$ dans $G$ est un Levi. Soit $x_{*}\in X_{*}(T)$ en position g\'en\'erale. Il d\'etermine un sous-groupe parabolique $P$ de composante de Levi $M$: $P$ est engendr\'e par $M$ et les sous-espaces radiciels pour l'action de $T$ dans l'alg\`ebre de Lie de $G$ associ\'es aux racines $\alpha$ telles que $<\alpha,x_{*}>>0$. Le normalisateur commun $\tilde{M}$ de $P$ et $M$ dans $\tilde{G}$ est un espace de Levi, s'il est non vide. Mais $Z_{\tilde{G}}(T)$ est inclus dans $\tilde{M}$ et est non vide par hypoth\`ese. Donc $\tilde{M}$ est un espace de Levi. C'est un espace principal homog\`ene pour l'action disons \`a gauche de $M$. Or $Z_{\tilde{G}}(T)$ est stable par cette action. L'inclusion $Z_{\tilde{G}}(T)\subset \tilde{M}$ est donc une \'egalit\'e. Dans la situation de (12), on a l'inclusion $\tilde{M}\subset Z_{\tilde{G}}(A_{\tilde{M}})$ et ce deuxi\`eme ensemble est un espace de Levi comme on vient de le prouver. Il suffit de prouver que les Levi associ\'es dans $G$ sont \'egaux, autrement dit que $M=Z_{G}(A_{\tilde{M}})$. Soit $\tilde{P}$ un sous-espace parabolique de $\tilde{G}$ de sous-espace de Levi $\tilde{M}$. Soit $y_{*}\in X_{*}(A_{M})$ d\'eterminant $P$ par la construction ci-dessus. Notons $x_{*}$ la somme des \'el\'ements de l'orbite de $y_{*}$ pour l'action du groupe d'automorphismes de $X_{*}(A_{M})$ engendr\'e par $\theta$, o\`u $\theta=ad_{\gamma}$ pour un \'el\'ement quelconque $\gamma\in \tilde{M}$. Alors $x_{*}\in X_{*}(A_{\tilde{M}})$. Comme $\theta$ pr\'eserve les racines de $A_{M}$ positives pour $P$, on voit que le couple $(P,M)$ co\"{\i}ncide avec celui construit dans la preuve de (11). Donc $Z_{G}(A_{\tilde{M}})\subset M$ et la conclusion. $\square$
  
   Soit $\tilde{M}$ un espace de Levi de $\tilde{G}$. Consid\'erons un espace parabolique $\tilde{P}$ de composante de Levi $\tilde{M}$ et une paire de Borel \'epingl\'ee ${\cal E}=(B,T,(E_{\alpha})_{\alpha\in \Delta})$ de $G$ telle que $B\subset P$ et $T\subset M$. Alors ${\cal E}^M=(B\cap M, T, (E_{\alpha})_{\alpha\in \Delta^M})$ est une paire de Borel \'epingl\'ee de $M$. On a une injection $Z(\tilde{G},{\cal E})\subset Z(\tilde{M},{\cal E}^M)$. On d\'eduit par passage aux quotients une application ${\cal Z}(\tilde{G},{\cal E})\to {\cal Z}(\tilde{M},{\cal E})$ qui s'identifie \`a une application ${\cal Z}(\tilde{G})\to {\cal Z}(\tilde{M})$. On laisse le lecteur v\'erifier que
  
 (14) cette application ${\cal Z}(\tilde{G})\to {\cal Z}(\tilde{M})$ ne d\'epend pas des choix de $\tilde{P}$ et de ${\cal E}$.
 
 Soit $\tilde{M}$ un espace de Levi de $\tilde{G}$. Fixons un espace parabolique $\tilde{P}$ de composante $\tilde{M}$ et un sous-groupe compact maximal $K$ de $G(F)$, en bonne position relativement \`a $M$ et sp\'ecial si $F$ est non-archim\'edien. On note $U$ le radical unipotent de $P$. Fixons des mesures de Haar sur $G(F)$ et $M(F)$. On en d\'eduit une mesure sur $U(F)\times K$ de sorte que  l'\'egalit\'e suivante soit v\'erifi\'ee
 $$\int_{G(F)}f(g)\,dg\,=\int_{M(F)\times U(F)\times K}f(muk)\,dk\,du\,dm$$
 pour toute $f\in C_{c}^{\infty}(G(F))$. On d\'efinit un homomorphisme
 $$\begin{array}{ccc}C_{c}^{\infty}(\tilde{G}(F))&\to&C_{c}^{\infty}(\tilde{M}(F))\\ f&\mapsto&f_{\tilde{M},\omega}\\ \end{array}$$
 par la formule 
 $$f_{\tilde{M},\omega}(\gamma)=\int_{U(F)\times K}f(k^{-1}u^{-1}\gamma uk)\omega^{-1}(k)\,du\,dk.$$
 Cet homomorphisme d\'epend des choix de $K$ et $\tilde{P}$. Mais il s'en d\'eduit un homomorphisme $I(\tilde{G}(F),\omega)\to I(\tilde{M}(F),\omega)$ qui n'en d\'epend plus. Pour $\gamma\in \tilde{M}(F)\cap \tilde{G}_{reg}(F)$, on a simplement $I^{\tilde{M}}(\gamma,\omega,f_{\tilde{M},\omega})=I^{\tilde{G}}(\gamma,\omega,f)$ pourvu bien s\^ur que l'on choisisse une mesure unique sur le groupe $M_{\gamma}(F)=G_{\gamma}(F)$. L'homomorphisme ci-dessus d\'epend encore des choix de mesures de Haar, mais on le rend canonique en le consid\'erant comme un homomorphisme
 $$\begin{array}{ccc}I(\tilde{G}(F),\omega)\otimes Mes(G(F))&\to&I(\tilde{M}(F),\omega)\otimes Mes(M(F))\\ {\bf f}&\mapsto&{\bf f}_{\tilde{M},\omega}\\ \end{array}$$ 
Notons $Norm_{G(F)}(\tilde{M})$ le normalisateur de $\tilde{M}$ dans $G(F)$ et posons $W(\tilde{M})=Norm_{G(F)}(\tilde{M})/M(F)$. Le groupe $Norm_{G(F)}(\tilde{M})$ agit sur $C_{c}^{\infty}(\tilde{M}(F))$ par $(x,f)\mapsto xf$, o\`u
$$(xf)(m)=\omega(x)f(x^{-1}mx).$$
Cette action se descend en une action de $W(\tilde{M})$ sur $I(\tilde{M}(F),\omega)$, donc aussi sur $I(\tilde{M}(F),\omega)\otimes Mes(M(F))$. L'image de l'homomorphisme ci-dessus est contenu dans le  sous-espace des invariants par cette action. On d\'ecrira cette image en 4.3.

On note $I_{cusp}(\tilde{G}(F),\omega)$ l'espace des $f\in I(\tilde{G}(F),\omega)$ tels que $f_{\tilde{M},\omega}=0$ pour tout espace de Levi propre $\tilde{M}$ de $\tilde{G}$. On note $C_{cusp}^{\infty}(\tilde{G}(F),\omega)$ l'espace des $f\in C_{c}^{\infty}(\tilde{G}(F),\omega)$ dont l'image dans $I(\tilde{G}(F),\omega)$ appartient \`a $I_{cusp}(\tilde{G}(F),\omega)$.

Consid\'erons le cas o\`u $(G,\tilde{G},{\bf a})$ est quasi-d\'eploy\'e et \`a torsion int\'erieure. Pour $\gamma\in \tilde{M}(F)\cap \tilde{G}_{reg}(F)$, on sait qu'un ensemble de repr\'esentants des classes de conjugaison par $M(F)$ dans la classe de conjugaison stable de $\gamma$ dans $\tilde{M}(F)$ est aussi un tel ensemble de repr\'esentants des classes de conjugaison par $G(F)$ dans la classe de conjugaison stable de $\gamma$ dans $\tilde{G}(F)$. Pour $f\in I(\tilde{G}(F),\omega)$, on a donc l'\'egalit\'e $S^{\tilde{M}}(\gamma,f_{\tilde{M}})=S^{\tilde{G}}(\gamma,f)$. Il en r\'esulte que l'homomorphisme compos\'e
$$I(\tilde{G}(F))\to I(\tilde{M}(F))\to SI(\tilde{M}(F))$$
se factorise en un homomorphisme
$$SI(\tilde{G}(F))\to SI(\tilde{M}(F))$$
que nous noterons aussi $f\mapsto f_{\tilde{M}}$. On note $SI_{cusp}(\tilde{G}(F))$ l'espace des $f\in SI(\tilde{G}(F))$ tels que $f_{\tilde{M}}=0$ (dans $SI(\tilde{M}(F))$) pour tout espace de Levi propre $\tilde{M}$ de $\tilde{G}$. Ces d\'efinitions s'adaptent au cas o\`u on consid\`ere une extension
$$1\to C_{1}\to G_{1}\to G\to 1$$
o\`u $C_{1}$ est un tore central induit, une extension compatible
$$\tilde{G}_{1}\to \tilde{G}$$
et un caract\`ere $\lambda_{1}$ de $C_{1}(F)$, et o\`u on remplace l'espace $C_{c}^{\infty}(\tilde{G}(F))$ par $C_{c,\lambda_{1}}^{\infty}(\tilde{G}_{1}(F)$.

\bigskip
 
 \subsection{ Donn\'ees endoscopiques d'espace de Levi}
 Consid\'erons un espace de Levi $\tilde{M}$ de $\tilde{G}$. Comme on l'a expliqu\'e, on peut r\'ealiser le $L$-groupe $^LM$ comme un sous-groupe de $^LG$. Pr\'ecis\'ement, apr\`es avoir fix\'e comme en 1.4 une paire de Borel \'epingl\'ee $\hat{\cal E}=(\hat{B},\hat{T},(\hat{E}_{\alpha})_{\alpha\in \Delta})$ de $\hat{G}$, on peut fixer une paire parabolique standard $(\hat{P},\hat{M})$ fixe par $\Gamma_{F}$ et $\hat{\theta}$ de sorte que $\hat{M}\rtimes W_{F}$ soit le $L$-groupe de $M$ et $(\hat{M}\rtimes W_{F})\hat{\theta}$ soit le $L$-espace $^L\tilde{M}$. On a alors un homomorphisme $H^1(W_{F};Z(\hat{G}))\to H^1(W_{F};Z(\hat{M}))$. En fait, il ne d\'epend pas des choix. On note ${\bf a}_{M}$ l'image de ${\bf a}$ dans $H^1(W_{F};Z(\hat{M}))$. Consid\'erons une donn\'ee endoscopique ${\bf M}'=(M',{\cal M}',\tilde{\zeta})$ pour $(\tilde{M},{\bf a}_{M})$.  Quitte \`a conjuguer $\hat{\cal E}$ par un \'el\'ement de $\hat{M}$, on suppose que $\tilde{\zeta}$ fixe $(\hat{B},\hat{T})$. Dans la d\'efinition d'une telle donn\'ee intervient un cocycle $a_{M}$ tel que $ad_{\tilde{\zeta}}(m,w)=(a_{M}(w)m,w)$ pour tout $w\in W_{F}$. Sa classe est ${\bf a}_{M}$. Si on remplace $\tilde{\zeta}$ par un \'el\'ement de $Z(\hat{M})\tilde{\zeta}$, ce cocycle est modifi\'e par un cobord. Pour quelques instants, notons plus pr\'ecis\'ement $a_{M,\tilde{\zeta}}$ le cocycle associ\'e \`a $\tilde{\zeta}$. On a
 
 (1) dans l'ensemble $Z(\hat{M})\tilde{\zeta}$, il existe une unique classe modulo $Z(\hat{M})^{\Gamma_{F}}Z(\hat{G})$ telle que, pour $\tilde{\zeta}'$ dans cette classe, $a_{M,\tilde{\zeta}'}$ prenne ses valeurs dans $Z(\hat{G})$.
 
 Preuve. L'hypoth\`ese que ${\bf a}_{M}$ provient d'un \'el\'ement de $H^1(W_{F};Z(\hat{G}))$ entra\^{\i}ne qu'il existe au moins un $\tilde{\zeta}'\in Z(\hat{M})\tilde{\zeta}$ tel que $a_{M,\tilde{\zeta}'}$ prenne ses valeurs dans $Z(\hat{G})$. Fixons-en un et pour simplifier les notations, supposons que ce soit $\tilde{\zeta}$ lui-m\^eme. Pour $z\in Z(\hat{M})$, on calcule $a_{M,z\tilde{\zeta}}(w)=zw(z)^{-1}a_{M,\tilde{\zeta}}(w)$.  Ce terme appartient \`a $Z(\hat{G})$ pour tout $w$ si et seulement si l'image $z_{ad}$ de $z$ dans $Z(\hat{M}_{ad})$ est fixe par $\Gamma_{F}$ (o\`u $\hat{M}_{ad}=\hat{M}/Z(\hat{G})$). Or $Z(\hat{M}_{ad})^{\Gamma_{F}}$ est connexe (c'est bien connu; on rappelle la preuve dans celle de  3.3(2) ci-dessous) donc est l'image naturelle de $Z(\hat{M})^{\Gamma_{F}}$. La condition \'equivaut donc \`a $z\in Z(\hat{M})^{\Gamma_{F}}Z(\hat{G})$. $\square$
 
 Quitte \`a remplacer $\tilde{\zeta}$ par un \'el\'ement convenable de $Z(\hat{M})\tilde{\zeta}$, on peut supposer que $\tilde{\zeta}$ appartient \`a l'unique classe d\'etermin\'ee par (1). C'est ce que l'on supposera  toujours, pour simplifier les notations. Autrement dit, on suppose que $a_{M}$ prend ses valeurs dans $Z(\hat{G})$.  Remarquons qu'alors, la classe de $a_{M}$ dans $H^1(W_{F};Z(\hat{G}))$ est \'egale \`a ${\bf a}$, d'apr\`es:
 
 (2) l'homomorphisme $H^1(W_{F};Z(\hat{G}))\to H^1(W_{F};Z(\hat{M}))$ est injectif. 
 
 Par la suite longue de cohomologie, cela r\'esulte de la surjectivit\'e remarqu\'ee ci-dessus de l'application $Z(\hat{M})^{\Gamma_{F}}\to Z(\hat{M}_{ad})^{\Gamma_{F}}$.
  
 Au lieu d'un espace de Levi et d'une donn\'ee endoscopique de cet espace, consid\'erons deux telles paires $(\tilde{M},{\bf M}')$ et $(\underline{\tilde{M}},\underline{{\bf M}}')$, soumises aux m\^emes hypoth\`eses que ci-dessus. On r\'ealise $^LM$ et $^L\underline{M}$ comme sous-groupes de $^LG$  et $^L\tilde{M}$ et $^L\underline{\tilde{M}}$ comme sous-ensembles de $^L\tilde{G}$ (il n'est pas n\'ecessaire d'utiliser une paire de Borel commune). Appelons \'equivalence entre ces donn\'ees un \'el\'ement $x\in \hat{G}$ tel que $ad_{x}(\hat{M})=\hat{\underline{M}}$, $ad_{x}({\cal M}')=\underline{{\cal M}}'$, $ad_{x}(\tilde{\zeta})\in Z(\underline{\hat{M}})\underline{\tilde{\zeta}}$. Remarquons que les ensembles $\tilde{{\cal M}}={^L\tilde{M}}$ et $\underline{\tilde{{\cal M}}}={^L\underline{\tilde{M}}}$, r\'ealis\'es comme sous-ensembles de $^L\tilde{G}$, sont des espaces de Levi et que les conditions impos\'ees \`a $x$ entra\^{\i}nent que $ad_{x}(\tilde{{\cal M}})=\underline{\tilde{{\cal M}}}$. D'apr\`es 3.1(10), \`a $x$ est donc associ\'e une classe $gM(F)$ dans $W(\tilde{M},\underline{\tilde{M}})$. 
 
 Fixons un isomorphisme $\iota':M'\to \underline{M}'$ d\'efini sur $F$ dual \`a la restriction de $ad_{x}^{-1}$ \`a $\hat{\underline{M}}'$. Remarquons que $ad_{g}$ d\'efinit un isomorphisme de ${\cal Z}(\tilde{M})$ sur ${\cal Z}(\underline{\tilde{M}})$. De ces deux isomorphismes r\'esulte un isomorphisme  $\tilde{\iota}':\tilde{M}'\to \underline{\tilde{M}}'$. Supposons ${\bf M}'$ relevant et fixons des donn\'ees suppl\'ementaires $M'_{1}$,...,$\Delta_{1}$. Posons $\underline{M}'_{1}=M'_{1}$, $\underline{C}_{1}=C_{1}$, avec pour homomorphisme
$\underline{M}'_{1}\to \underline{M}'$ le compos\'e de $M'_{1}\to M'$ et de $\iota'$.  On pose $\underline{\tilde{M}}'_{1}=\tilde{M}'_{1}$ muni de l'application $\underline{\tilde{M}}'_{1}\to \underline{\tilde{M}}'$ compos\'ee de $\tilde{M}'_{1}\to \tilde{M}'$ et de  $\tilde{\iota}'$. On pose $\underline{\hat{\xi}}_{1}=\hat{\xi}_{1}\circ ad_{x}^{-1}:\underline{{\cal M}}\to{^L\underline{M}}'_{1}={^LM}'_{1}$. Ces donn\'ees v\'erifient les conditions requises relativement \`a la donn\'ee $\underline{{\bf M}}'$. Pour $(\delta_{1},\gamma)\in {\cal D}_{1}$ (l'ensemble relatif aux premi\`eres donn\'ees), on a $(\delta_{1},g\gamma g^{-1})\in \underline{{\cal D}}_{1}$ (l'ensemble relatif aux secondes). On v\'erifie l'\'egalit\'e
$$\underline{\Delta}_{1}(\delta_{1},g\gamma g^{-1}; \delta'_{1},g\gamma' g^{-1})=\Delta_{1}(\delta_{1},\gamma;\delta'_{1},\gamma').$$
On choisit alors pour facteur de transfert pour les secondes donn\'ees le facteur 
$$\underline{\Delta}_{1}(\delta_{1},gÊ\gamma g^{-1})=\omega(g)\Delta_{1}(\delta_{1},\gamma).$$
Cette d\'efinition ne d\'epend que de la classe $gM(F)$. Ces choix fournissent les isomorphismes extr\^emes de la suite
$$C_{c}^{\infty}({\bf M}')\simeq C_{c,\lambda_{1}}^{\infty}(\tilde{M}'_{1}(F))=C_{c,\lambda_{1}}^{\infty}(\underline{\tilde{M}}'_{1}(F))\simeq C_{c}^{\infty}(\underline{{\bf M}}').$$
Ici encore, l'isomorphisme obtenu d\'epend du choix de $\iota'$. Mais il devient ind\'ependant de ce choix si on se limite \`a des fonctions invariantes par l'action des groupes adjoints. Comme en 2.6, dans le cas particulier o\`u $\underline{\tilde{M}}=\tilde{M}$ et $\underline{{\bf M}}'={\bf M}'$, on note $Aut(\tilde{M},{\bf M}')$ le groupe des automorphismes de la paire $(\tilde{M},{\bf M}')$ (c'est-\`a-dire de ses \'equivalences avec elle-m\^eme). On obtient une action de ce groupe sur
  $C_{c}^{\infty}({\bf M}')$. il y a une suite exacte
 $$1\to Aut({\bf M}')\to Aut(\tilde{M},{\bf M}')\to W(\tilde{M},{\bf M}')\to 1$$
 o\`u $W(\tilde{M},{\bf M}')$ est un sous-groupe de $W(\tilde{M})$. En particulier, on a une \'egalit\'e d'espaces invariants
 $$SI({\bf M}')^{Aut(\tilde{M},{\bf M}')}=(SI({\bf M}')^{Aut({\bf M}')})^{W(\tilde{M},{\bf M}')}.$$

 \bigskip

\subsection{Donn\'ees endoscopiques de $\tilde{G}$ associ\'ees \`a une donn\'ee endoscopique d'un espace de Levi} 
Soient $\tilde{M} $ un espace de Levi de $\tilde{G}$ et ${\bf M}'=(M',{\cal M}',\tilde{\zeta})$ une donn\'ee endoscopique de $(\tilde{M},{\bf a}_{M})$. On reprend la situation du d\'ebut du paragraphe pr\'ec\'edent et on note $\hat{P}$ le sous-groupe parabolique standard dont $\hat{M}$ est la composante de Levi standard. Pour $\tilde{s}\in Z(\hat{M})^{\Gamma_{F}}\tilde{\zeta}$, posons $\hat{G}'(\tilde{s})=Z_{\hat{G}}(\tilde{s})^0$ et ${\cal G}'(\tilde{s})=\hat{G}'(\tilde{s}){\cal M}'$. On v\'erifie que ${\cal G}'(\tilde{s})$ est un groupe. Remarquons que:
 
 (1) $\hat{M}'$ est un Levi de $\hat{G}'(\tilde{s})$.
 
 En effet, d'apr\`es les d\'efinitions, $\hat{M}'$ est \'egal \`a $(\hat{M}\cap \hat{G}'(\tilde{s}))^0$. La m\^eme preuve qu'en 3.1(11) montre que $\hat{M}$ est le commutant de $Z(\hat{M})^{\hat{\theta},0}$ dans $\hat{G}$. Donc $\hat{M}\cap \hat{G}'(\tilde{s})$ est le commutant de $Z(\hat{M})^{\hat{\theta},0}$ dans $\hat{G}'(\tilde{s})$. Remarquons que $Z(\hat{M})^{\hat{\theta},0}$ est un tore dans $\hat{G}'(\tilde{s})$. Donc $\hat{M}\cap \hat{G}'(\tilde{s})$ est un Levi de ce groupe. Un Levi est connexe et (1) s'ensuit.
 
  Fixons une paire de Borel  \'epingl\'ee de $\hat{G}'(\tilde{s})$ pour laquelle $(\hat{P'}(\tilde{s}),\hat{M}')$ est standard, o\`u $\hat{P'}(\tilde{s})=\hat{G'}(\tilde{s})\cap \hat{P}$.    On munit $\hat{G}'(\tilde{s})$ de l'unique action $\sigma\mapsto \sigma_{G'(\tilde{s})}$ de $\Gamma_{F}$ conservant cette paire de Borel \'epingl\'ee et telle que, pour tout $(m,w)\in {\cal M}'$, l'action par conjugaison de $(m,w)$ sur $\hat{G}'(\tilde{s})$ soit \'egale \`a $w_{G'(\tilde{s})}$ compos\'e avec un automorphisme int\'erieur.  Cette action conserve la paire $(\hat{P}'(\tilde{s}),\hat{M}')$. On introduit un groupe dual $G'(\tilde{s})$ r\'eductif connexe d\'efini et quasi-d\'eploy\'e sur $F$. Alors ${\bf G}'(\tilde{s})=(G'(\tilde{s}),{\cal G}'(\tilde{s}),\tilde{s})$ est une donn\'ee endoscopique pour $(G,\tilde{G},{\bf a})$. En particulier, il y a un espace endoscopique $\tilde{G}'(\tilde{s})$. Puisque la paire $(\hat{P}'(\tilde{s}),\hat{M}')$ est invariante par $\Gamma_{F}$, $M'$ s'identifie \`a un Levi de $G'(\tilde{s})$ et on v\'erifie que l'espace endoscopique $\tilde{M}'$ s'identifie conform\'ement \`a un espace de Levi de $\tilde{G}'(\tilde{s})$.
 
 Soient $\delta\in \tilde{M}'_{reg}(F)$ et $\gamma\in \tilde{M}_{reg}(F)$. Si la classe de conjugaison par $M'$ de $\delta$ correspond \`a la classe de conjugaison par $M$ de $\gamma$, alors la classe de conjugaison par $G'(\tilde{s})$ de $\delta$ correspond \`a la classe de conjugaison par $G$ de $\gamma$. Autrement dit ${\cal D}({\bf M}')\subset {\cal D}({\bf G}'(\tilde{s}))$. Inversement, pour $(\delta,\gamma)\in {\cal D}({\bf G}'(\tilde{s}))\cap (\tilde{M}'(F)\times \tilde{M}(F))$, il existe un \'el\'ement $n\in Norm_{G(F)}(\tilde{M})$   tel que  $(\delta,n\gamma n^{-1})$ appartienne \`a ${\cal D}({\bf M}')$.  Supposons ${\bf M}'$ relevant. Alors ${\bf G}'(\tilde{s})$ l'est aussi. On voit que le bifacteur de transfert pour la donn\'ee ${\bf M}'$ co\"{\i}ncide avec la restriction \`a ${\cal D}({\bf M}')\times {\cal D}({\bf M}')$ du bifacteur de transfert pour la donn\'ee ${\bf G}'(\tilde{s})$. Fixons des donn\'ees auxiliaires $G'_{1}(\tilde{s})$, $\tilde{G}'_{1}(\tilde{s})$, $C_{1}(\tilde{s})$, $\hat{\xi}_{1}(\tilde{s})$, $\Delta_{1}(\tilde{s})$. On note $\lambda_{1}(\tilde{s})$ le caract\`ere de $C_{1}(\tilde{s})$ associ\'e \`a ces donn\'ees. On note $M'_{1}(\tilde{s})$ et $\tilde{M}'_{1}(\tilde{s})$ les images r\'eciproques de $M'$ et $\tilde{M}'$ dans $G'_{1}(\tilde{s})$ et $\tilde{G}'_{1}(\tilde{s})$. On  note $\hat{\xi}_{1,M'}(\tilde{s})$ la restriction de $\hat{\xi}_{1}(\tilde{s})$ \`a ${\cal M}$, ${\cal D}_{1,M'}$ l'image r\'eciproque de ${\cal D}({\bf M}')$ dans ${\cal D}_{1}$ et $\Delta_{1, M'}(\tilde{s})$ la restriction de $\Delta_{1}(\tilde{s})$ \`a ${\cal D}_{1,M'}$. Alors $(M'_{1}(\tilde{s}),\tilde{M}'_{1}(\tilde{s}),C_{1}(\tilde{s}),\hat{\xi}_{1,M'}(\tilde{s}),\Delta_{1,M'}(\tilde{s}))$ sont des donn\'ees auxiliaires pour ${\bf M}'$.  Par une variante de la construction de 3.1, on a un homomorphisme
  $$\begin{array}{ccc}I_{\lambda_{1}(\tilde{s})}(\tilde{G}'_{1}(\tilde{s};F))\otimes Mes(G'(\tilde{s};F))&\to&I_{\lambda_{1}(\tilde{s})}(\tilde{M}'_{1}(\tilde{s};F))\otimes Mes(M'(F))\\ f&\mapsto&f_{\tilde{M}'}\\ \end{array}$$
  On v\'erifie que, quand on change de donn\'ees auxiliaires, ces homomorphismes sont compatibles aux applications de recollement de 2.5. On obtient un homomorphisme
 $$\begin{array}{ccc}I({\bf G}'(\tilde{s}))\otimes Mes(G'(\tilde{s};F))&\to&I({\bf M}')\otimes Mes(M'(F))\\ f&\mapsto&f_{\tilde{M}'}\\ \end{array}$$

 Pour $\lambda\in Z(\hat{M} )^{\Gamma_{F}}$ et $\nu\in Z(\hat{G})^{\Gamma_{F}}$, posons $\tilde{s}'=\nu \lambda\tilde{s}\lambda^{-1}$. Alors la donn\'ee ${\bf G}'(\tilde{s}')$ est \'equivalente \`a ${\bf G}'(\tilde{s})$, l'\'equivalence \'etant d\'efinie par $\lambda$. Dans les constructions o\`u seule la classe d'\'equivalence de ${\bf G}'(\tilde{s})$ importe, on pourra consid\'erer que $\tilde{s}$ parcourt l'ensemble des classes de conjugaison par $Z(\hat{M})^{\Gamma_{F}}$ dans $\tilde{\zeta}Z(\hat{M})^{\Gamma_{F}}/Z(\hat{G})^{\Gamma_{F}} $.  Par l'application $z\mapsto \tilde{\zeta}z$, cet ensemble de classes de conjugaison s'identifie \`a  $Z(\hat{M})^{\Gamma_{F}}/(Z(\hat{G})^{\Gamma_{F}}(1-\hat{\theta})(Z(\hat{M})^{\Gamma_{F}}))$.
 On le remplacera souvent  par $\tilde{\zeta}Z(\hat{M})^{\Gamma_{F},\hat{\theta}}/Z(\hat{G})^{\Gamma_{F},\hat{\theta}}$ gr\^ace \`a l'assertion suivante. On y note $\theta^{\tilde{M}}$ l'automorphisme de ${\cal A}_{M}$ induit par $ad_{\gamma}$ pour n'importe quel $\gamma\in \tilde{M}$. On a

(2) l'homomorphisme naturel
$$Z(\hat{M})^{\Gamma_{F},\hat{\theta}}/Z(\hat{G})^{\Gamma_{F},\hat{\theta}}\to Z(\hat{M})^{\Gamma_{F}}/(Z(\hat{G})^{\Gamma_{F}}(1-\hat{\theta})(Z(\hat{M})^{\Gamma_{F}}))$$
est surjectif; son noyau a pour nombre d'\'el\'ements $\vert det((1-\theta^{\tilde{M}})_{\vert {\cal A}_{M}/({\cal A}_{\tilde{M}}+{\cal A}_{G}})\vert $.

Preuve.   Introduisons l'ensemble des racines simples $\Delta$ de $\hat{T}$, le sous-ensemble $\Delta^M$ associ\'e \`a $\hat{M}$ et celui des copoids fondamentaux $\{\check{\varpi}_{\alpha}; \alpha\in \Delta\}\subset X_{*}(\hat{ T}_{ad})$. Le groupe $Z(\hat{M}_{ad})^{\Gamma_{F}}$ est le sous-groupe des \'el\'ements $\prod_{\alpha\in \Delta- \Delta^M}\check{\varpi}_{\alpha}(t_{\alpha})\in \hat{ T}_{ad}$ avec $t_{\alpha}\in {\mathbb C}^{\times}$ et $\alpha\mapsto t_{\alpha}$ est constante sur les orbites de $\Gamma_{F}$ dans $\Delta- \Delta^M$. Donc $Z(\hat{M}_{ad})^{\Gamma_{F}}$ est connexe.  Il en r\'esulte que l'homomorphisme
$$Z(\hat{M})^{\Gamma_{F}}\to Z(\hat{M}_{ad})^{\Gamma_{F}}$$
est surjectif. Le m\^eme calcul montre que $Z(\hat{M}_{ad})^{\Gamma_{F},\hat{\theta}}$ est connexe et que  l'homomorphisme
$$Z(\hat{M})^{\Gamma_{F}, \hat{\theta}}\to Z(\hat{M}_{ad})^{\Gamma_{F},\hat{\theta}}$$
est surjectif. Les ensembles de d\'epart et d'arriv\'ee de  l'homomorphisme (2) s'identifient respectivement \`a $Z(\hat{M}_{ad})^{\Gamma_{F},\hat{\theta}}$ et $Z(\hat{M}_{ad})^{\Gamma_{F}}/(1-\hat{\theta})(Z(\hat{M}_{ad})^{\Gamma_{F}})$. Tout se d\'ecompose selon les orbites dans $\Delta- \Delta^M$ de l'action du groupe engendr\'e par $\Gamma_{F}$ et $\hat{\theta}$, ce qui nous ram\`ene au cas o\`u il n'y a qu'une seule orbite. Fixons un \'el\'ement $\alpha\in \Delta- \Delta^M$, notons $[\alpha]$ son orbite sous l'action de $\Gamma_{F}$, $n$ le plus petit entier $\geq1$ tel que $\hat{\theta}^{n}(\alpha)\in[\alpha]$ et posons $\check{\varpi}_{[\alpha]}=\sum_{\beta\in [\alpha]}\check{\varpi_{\beta}}$. Un \'el\'ement de $Z(\hat{M}_{ad})^{\Gamma_{F}}$ s'\'ecrit $\prod_{i=0,...,n-1}\check{\varpi}_{\hat{\theta}^{i}[\alpha]}(t_{i})$, avec des $t_{i}\in {\mathbb C}^{\times}$. Il appartient \`a $(1-\hat{\theta})(Z(\hat{M}_{ad})^{\Gamma_{F}})$ si et seulement si $\prod_{i}t_{i}=1$. Il appartient \`a $Z(\hat{M}_{ad})^{\Gamma_{F},\hat{\theta}}$ si et seulement si les $t_{i}$ sont tous \'egaux. Il r\'esulte de cette description que notre homomorphisme est surjectif et que son noyau a $n$ \'el\'ements. Or $n$ est \'egal au d\'eterminant figurant dans l'assertion (2). $\square$

On a aussi

(3) l'ensemble des $\tilde{s}\in \tilde{\zeta}Z(\hat{M})^{\Gamma_{F},\hat{\theta}}/Z(\hat{G})^{\Gamma_{F},\hat{\theta}}$ tels que ${\bf G}'(\tilde{s})$ soit une donn\'ee endoscopique elliptique de $\tilde{G}$ est fini; si ${\bf M}'$ est une donn\'ee endoscopique elliptique de $\tilde{M}$, cet ensemble n'est pas vide.

Preuve. Cf. [W2] 3.2(1) pour la finitude. Pour la deuxi\`eme assertion,   utilisons les m\^emes notations que dans la preuve pr\'ec\'edente. Ecrivons $\tilde{\zeta}=\zeta\hat{\theta}$. Soit $\Delta_{0}$ un ensemble de repr\'esentants dans $\Delta- \Delta^M$ des orbites pour l'action du groupe engendr\'e par $\Gamma_{F}$ et $\hat{\theta}$. L'homomorphisme
$$\begin{array}{ccc}Z(\hat{M})^{\Gamma_{F},\hat{\theta},0}/Z(\hat{G})^{\Gamma_{F},\hat{\theta},0}&\to&({\mathbb C}^{\times})^{\Delta_{0}}\\ x&\mapsto&(\alpha(x))_{\alpha\in \Delta_{0}}\\ \end{array}$$
est surjective \`a noyau fini. Il existe donc $x\in Z(\hat{M})^{\Gamma_{F},\hat{\theta},0}$ tel que $(N\alpha)(x\zeta)=1$ pour tout $\alpha\in \Delta_{0}$. Pour un tel \'el\'ement, posons $\tilde{s}=x\tilde{\zeta}$.  L'alg\`ebre de Lie de $\hat{G}'(\tilde{s})$ contient $\sum_{i=0,...,n_{\alpha}-1}(ad_{\tilde{s}})^{i}(\hat{E}_{\alpha})$ pour tout $\alpha\in \Delta_{0}$, o\`u  $n_{\alpha}$ est le plus petit entier $i\geq1$ tel que $\hat{\theta}^{i}(\alpha)=\alpha$.  Un \'el\'ement de $Z(\hat{G}'(\tilde{s}))$ fixe cet \'el\'ement donc aussi chaque composante $\hat{E}_{\hat{\theta}^{i}\alpha}$. Remarquons que les actions galoisiennes relatives \`a $\hat{G}$ et \`a $\hat{G}'(\tilde{s})$ co\"{\i}ncident sur $Z(\hat{G}'(\tilde{s}))\cap Z(\hat{M})$. Un \'el\'ement de $Z(\hat{G}'(\tilde{s}))^{\Gamma_{F}}\cap Z(\hat{M}) $ fixe donc $\hat{E}_{\sigma\hat{\theta}^{i}\alpha}$ pour tous $\alpha\in \Delta_{0}$, $i\in {\mathbb N}$ et $\sigma\in \Gamma_{F}$. Donc il fixe  $\hat{E}_{\alpha}$ pour tout $\alpha\in \Delta-\Delta^M$. Appartenant de plus \`a $Z(\hat{M})$, il fixe tout $\hat{G}$. Donc $Z(\hat{G}'(\tilde{s}))^{\Gamma_{F}}\cap Z(\hat{M})\subset Z(\hat{G})$. Or 
$$Z(\hat{G}'(s))^{\Gamma_{F},0}\subset Z(\hat{M'})^{\Gamma_{F},0}=Z(\hat{M})^{\Gamma_{F},\hat{\theta},0}$$
par l'hypoth\`ese d'ellipticit\'e de ${\bf M}'$. Donc $Z(\hat{G}'(s))^{\Gamma_{F},0}\subset Z(\hat{G})$  et forc\'ement $Z(\hat{G}'(s))^{\Gamma_{F},0}\subset Z(\hat{G})^{\Gamma_{F},\hat{\theta},0}$. $\square$

\bigskip

\subsection{Levi de donn\'ees endoscopiques}
Soient ${\bf G}'=(G',{\cal G}',\tilde{s})$ une donn\'ee endoscopique de $(G,\tilde{G},{\bf a})$ et $M'\subset G'$ un  Levi, auquel est associ\'e un espace de Levi $\tilde{M}'$ (puisque $\tilde{G}'$ est \`a torsion int\'erieure). On fixe une paire de Borel \'epingl\'ee de $\hat{G}'$ et on normalise l'action galoisienne sur ce groupe de sorte qu'elle conserve cette paire. Le choix d'un parabolique $P'$ de $G'$ de composante de Levi $M'$ permet d'identifier $\hat{M}'$ \`a un Levi standard de $\hat{G}'$, donc \`a un sous-groupe de $\hat{G}$. Notons $\hat{M}$, ${\cal M}$, $\tilde{{\cal M}}$ les commutants de $Z(\hat{M}')^{\Gamma_{F},0}$ dans $\hat{G}$, $^LG$, $^L\tilde{G}$.   Fixons $x_{*}\in X_{*}(Z(\hat{M}')^{\Gamma_{F},0})$ en position g\'en\'erale. Il d\'etermine un sous-groupe parabolique $\hat{P}$ de $\hat{G}$, engendr\'e par $\hat{M}$ et les sous-groupes radiciels associ\'es aux racines $\alpha$ de $\hat{T}$   telles que $<\alpha,x_{*}>>0$ ($\hat{T}$ \'etant choisi comme en 1.5). On pose ${\cal P}=\hat{P}{\cal M}$, $\tilde{{\cal P}}=\hat{P}\tilde{{\cal M}}$.   Le couple $(\tilde{{\cal P}},\tilde{{\cal M}})$ est une paire parabolique  de $^L\tilde{G}$. Les seuls points non \'evidents \`a v\'erifier sont que la projection de ${\cal P}$ sur $W_{F}$ est surjective et que $\tilde{\cal{P}}$ est non vide. Mais $\tilde{s}$ appartient \`a $\tilde{{\cal P}}$, ce qui v\'erifie ce deuxi\`eme point. Pour $w\in W_{F}$, il existe $g_{w}=(g(w),w)\in {\cal G}'$ tel que $ad_{g_{w}}$ agisse sur $\hat{G}'$ comme $w_{G'}$. Alors $ad_{g_{w}}$ fixe $x_{*}$ donc aussi $\hat{P}$. Donc $g_{w}\in {\cal P}$, ce qui v\'erifie le premier point. On pose ${\cal M}'={\cal G}'\cap {\cal M}$.  On se rappelle qu'il y a une injection de l'ensemble des paires paraboliques de $\tilde{G}$ dans celui des paires paraboliques de $^L\tilde{G}$. Si $G$ n'est pas quasi-d\'eploy\'e, $(\tilde{{\cal P}},\tilde{{\cal M}})$ peut ne pas appartenir \`a l'image: c'est le cas si et seulement si $(\hat{P},\hat{M})$ ne contient pas de conjugu\'e d'une paire $(\hat{P}_{0},\hat{M}_{0})$ comme en 3.1. On sait que les Levi $\hat{M}_{0}$ ont une propri\'et\'e particuli\`ere: tous les paraboliques ayant $\hat{M}_{0}$ comme composante de Levi sont conjugu\'es. Cela entra\^{\i}ne que la condition pr\'ec\'edente ne d\'epend que de $\hat{M}$ et pas du choix de $\hat{P}$. Supposons que $(\tilde{{\cal P}},\tilde{{\cal M}})$ soit l'image d'une paire parabolique $(\tilde{P},\tilde{M})$ de $\tilde{G}$. On dira dans ce cas que $\hat{M}$ correspond \`a l'espace de Levi $\tilde{M}$. Alors ${\bf M}'=(M',{\cal M}',\tilde{s})$ est une donn\'ee endoscopique pour $(\tilde{M},{\bf a}_{M})$. Cette donn\'ee est elliptique par construction. M\^eme si ${\bf G}'$ est relevant, il peut se produire que ${\bf M}'$ ne le soit pas. On dira que $M'$ est relevant si d'une part, $\hat{M}$ correspond \`a un espace de Levi $\tilde{M}$, d'autre part ${\bf M}'$ est relevant. Dans ce cas, comme dans le paragraphe pr\'ec\'edent, des donn\'ees auxiliaires pour ${\bf G}'$ se restreigent en des donn\'ees auxiliaires pour ${\bf M}'$ et on d\'efinit un homomorphisme
 $$\begin{array}{ccc}I({\bf G}')\otimes Mes(G'(F))&\to&I({\bf M}')\otimes Mes(M'(F))\\ f&\mapsto&f_{\tilde{M}'}\\ \end{array}$$

  En fait, seule la classe d'\'equivalence des donn\'ees $(\tilde{M},{\bf M}')$ est bien d\'etermin\'ee car on a effectu\'e divers choix.  Changer ces choix compose l'homomorphisme ci-dessus par des \'el\'ements de $Aut(\tilde{M},{\bf M}')$. Cela entra\^{\i}ne la propri\'et\'e suivante: si $f$ est un \'el\'ement de $ I({\bf G}')\otimes Mes(G'(F))$ et $\varphi$ est un \'el\'ement de $I({\bf M}')\otimes Mes(M'(F))$ invariant par l'action de $Aut(\tilde{M},{\bf M}')$, alors la relation $f_{\tilde{M}'}=\varphi$ est ind\'ependante des choix. De m\^eme, levons l'hypoth\`ese que $M'$ est relevant, supposons seulement que ${\bf G}'$ le soit. On ne peut plus d\'efinir d'espace $I({\bf M}')$. Mais, pour $f\in  I({\bf G}')\otimes Mes(G'(F))$, la relation $f_{\tilde{M}'}=0$ a un sens: elle signifie que si, par le choix de donn\'ees auxiliaires, on identifie $f$ \`a un \'el\'ement $f_{1}\in C_{c,\lambda_{1}}^{\infty}(\tilde{G}'_{1}(F))\otimes Mes(G'(F))$, alors $(f_{1})_{\tilde{M}'_{1}}=0$. Ceci est ind\'ependant du choix des donn\'ees auxiliaires.

On peut remplacer dans les constructions ci-dessus les espaces $I({\bf G}')$ par $SI({\bf G}')$.

\bigskip
\subsection{$K$-espaces}
Supposons $F={\mathbb R}$ et consid\'erons un $K$-espace $K\tilde{G}$ sur un $K$-groupe $KG$ comme en 1.11. Les constructions des quatre paragraphes pr\'ec\'edents valent pour chaque composante $\tilde{G}_{p}$. Mais en travaillant composante par composante, on perd la notion de $K$-espace. Pour la retrouver, il faut d\'efinir correctement les notions d'espace parabolique et d'espace de Levi d'un $K$-espace. Sur ${\mathbb C}$, tous les groupes $G_{p}$ ou espaces $\tilde{G}_{p}$ sont isomorphes, d'o\`u une correspondance bijective entre leurs classes de conjugaison de paires paraboliques. On d\'efinit une paire parabolique $(KP,KM)$ de $KG$  sur ${\mathbb C}$ comme une famille $(P_{p},M_{p})_{p\in {\Pi}}$,  o\`u $(P_{p},M_{p})$ est une paire parabolique (sur ${\mathbb C}$) de $G_{p}$ de sorte que, pour $p,p'\in {\Pi}$, les classes de conjugaison de $(P_{p},M_{p})$ et $(P_{p'},M_{p'})$ se correspondent. On d\'efinit de m\^eme une paire parabolique de $K\tilde{G}$. La d\'efinition est plus subtile sur ${\mathbb R}$. On d\'efinit une paire parabolique $(KP,KM)$ (sur ${\mathbb R}$, pr\'ecision que l'on omettra dans la suite) comme une famille $(P_{p},M_{p})_{\alpha\in{\Pi}'}$ o\`u

- ${\Pi}'$ est un sous-ensemble non vide de ${\Pi}$;

- pour tout $p\in {\Pi}'$, $(P_{p},M_{p})$ est une paire parabolique (sur ${\mathbb R}$) de $G_{p}$;

- pour $p,p'\in {\Pi}'$, les classes de conjugaison de $(P_{p},M_{p})$ et $(P_{p'},M_{p'})$ se correspondent;

- pour $p\in {\Pi}-{\Pi}'$, la classe de conjugaison de paires paraboliques de $G_{p}$  correspondant \`a celles des $(P_{p'},M_{p'})$ pour $p'\in {\Pi}'$ ne contient aucun \'el\'ement d\'efini sur ${\mathbb R}$.

En particulier, si ${\Pi}'\not={\Pi}$, une telle paire n'est pas une paire parabolique sur ${\mathbb C}$. On d\'efinit un Levi de $KG$ comme une famille $KM$ intervenant dans une paire parabolique $(KP,KM)$. On d\'efinit de m\^eme les paires paraboliques et les espaces de Levi de $K\tilde{G}$. On appellera plut\^ot ces derniers des $K$-espaces de Levi.  Si $(K\tilde{P},K\tilde{M})$ est une paire parabolique de $K\tilde{G}$, la paire sous-jacente $(KP,KM)$ est une paire parabolique de $KG$. On a

(1) tout espace de Levi $K\tilde{M}$ s'identifie \`a un $K$-espace tordu sur le $K$-groupe $KM$.

   Preuve. On compl\`ete $K\tilde{M}$ en une paire parabolique $(K\tilde{P},K\tilde{M})$. On fixe $p_{0}$ dans l'ensemble d'indices ${\Pi}'$ relatif \`a cette paire, on pose $G=G_{p_{0}}$,   $M=M_{p_{0}}$  etc... Pour $p\in {\Pi}'$, on choisit $x_{p}\in G_{SC}$ tel que $ad_{x_{p}}\circ\phi_{p_{0},p}$ envoie $(P_{p},M_{p})$ sur $(P,M)$. On note $\phi_{p}^M$ la restriction de $ad_{x_{p}}\circ\phi_{p_{0},p}$ \`a $M_{p}$ et $\tilde{\phi}^M_{p}$ celle de $ad_{x_{p}}\circ\tilde{\phi}_{p_{0},p}$ \`a $\tilde{M}_{p}$. Pour $\sigma\in\Gamma_{{\mathbb R}}$, on pose $\nabla_{p}^M(\sigma)=x_{p} \nabla_{p_{0},p}(\sigma)\sigma(x_{p})^{-1}$. On v\'erifie que $\nabla_{p}^M$ est un cocycle, \`a valeurs dans $G_{SC}$. On a $\phi^M_{p}\circ\sigma(\phi^M_{p})^{-1}=ad_{\nabla_{p}^M(\sigma)}$, $\tilde{\phi}^M_{p}\circ\sigma(\tilde{\phi}^M_{p})^{-1}=ad_{\nabla_{p}^M(\sigma)}$. Puisque $\phi^M_{p}\circ\sigma(\phi^M_{p})^{-1}$ pr\'eserve $(P,M)$, on en d\'eduit $\nabla_{p}^M(\sigma)\in M_{sc}$. D'apr\`es le  th\'eor\`eme 1.2 de [K2], l'image de l'application
  $$(2) \qquad H^1(\Gamma_{{\mathbb R}};M_{SC})\to H^1(\Gamma_{{\mathbb R}};M_{sc})$$
   est le noyau d'une application $H^1(\Gamma_{{\mathbb R}};M_{sc})\to \pi_{0}(Z(\hat{M}_{ad})^{\Gamma_{{\mathbb R}}})$. Or
   $Z(\hat{M}_{ad})^{\Gamma_{{\mathbb R}}}$ est connexe car $Z(\hat{M}_{ad})$ est un tore induit. Donc l'application (2) est surjective et, quitte \`a modifier l'\'el\'ement $x_{p}$, on peut relever $\nabla_{p}^M$ en un cocycle $\nabla_{p}^{M_{SC}}$ \`a valeurs dans $M_{SC}$.
    Pour prouver que $K\tilde{M}$ est un $K$-espace tordu issu de $\tilde{M}$ comme en 1.11, il reste \`a prouver que la famille $(\nabla_{p}^{M_{SC}})_{p\in {\Pi}'}$ s'envoie bijectivement sur $\pi(H^1(\Gamma_{{\mathbb R}};M_{SC}))\cap H^1(\Gamma_{{\mathbb R}};M)^{\theta}$ (o\`u $\theta$ est d\'etermin\'e par $\tilde{M}$).  Puisque $M$ est un Levi de $G$, l'application $H^1(\Gamma_{{\mathbb R}};M)\to H^1(\Gamma_{{\mathbb R}};G)$ est injective. Elle est \'equivariante pour l'action de $\theta$. Il en r\'esulte qu'un \'el\'ement de $H^1(\Gamma_{{\mathbb R}};M)$ est invariant par $\theta$ si et seulement si son image dans $H^1(\Gamma_{{\mathbb R}};G)$ l'est. L'image de $\nabla_{p}^{M_{SC}}$ dans $H^1(\Gamma_{{\mathbb R}};G)$ est \'egale \`a celle de $\nabla_{p_{0},p}$, donc est invariante par $\theta$. Donc l'image de $\nabla_{p}^{M_{SC}}$ dans $H^1(\Gamma_{{\mathbb R}};M)$ est invariante par $\theta$.  De m\^eme, pour $p,q\in \Pi'$ avec $p\not=q$,  les images de $\nabla_{p}^{M_{SC}}$ et $\nabla_{q}^{M_{SC}}$ dans cet ensemble sont distinctes car leurs images dans $H^1(\Gamma_{{\mathbb R}};G)$ le sont. Soit enfin $\nabla^{M}:\Gamma_{{\mathbb R}}\to M$ un cocycle dont la classe appartient \`a $\pi(H^1(\Gamma_{{\mathbb R}},M_{SC}))\cap H^1(\Gamma_{{\mathbb R}};M)^{\theta}$. Son image $\nabla^G$ dans $H^1(\Gamma_{{\mathbb R}};G)$ appartient \`a $\pi(H^1(\Gamma_{{\mathbb R}},G_{SC}))\cap H^1(\Gamma_{{\mathbb R}};G)^{\theta}$. Il existe donc $p\in \Pi$ tel que $\nabla^G$ soit cohomologue \`a $\nabla_{p_{0},p}$. Fixons $y\in G$ tel que $\nabla^M(\sigma)=y\nabla_{p_{0},p}(\sigma)\sigma(y)^{-1}$ pour tout $\sigma\in \Gamma_{{\mathbb R}}$. Puisque $\nabla^M$ prend ses valeurs dans $M$, cette relation implique que l'image r\'eciproque $(P'_{p},M'_{p})$ de $(P,M)$ par l'application $ad_{y}\circ \phi_{p_{0},p}$ est une paire de Borel de $G_{p}$ qui est d\'efinie sur ${\mathbb R}$. Cette paire est conjugu\'ee par un \'el\'ement de $G_{p}({\mathbb C})$ \`a l'image r\'eciproque de $(P,M)$ par l'application $\phi_{p_{0},p}$. Il en r\'esulte que $p\in \Pi'$ et que les paires de Borel $(P'_{p},M'_{p})$ et $(P_{p},M_{p})$ sont conjugu\'ees par un \'el\'ement de $G_{p}({\mathbb C})$. Etant toutes deux d\'efinies sur ${\mathbb R}$, elles sont conjugu\'ees par un \'el\'ement de $G_{p}({\mathbb R})$. On peut  donc fixer un \'el\'ement $g_{p}\in G_{p}({\mathbb R})$ tel que $ad_{y}\circ\phi_{p_{0},p}\circ ad_{g_{p}}(P_{p},M_{p})=(P,M)$. En posant $g=\phi_{p_{0},p}(g_{p})$, cela \'equivaut \`a $ad_{yg}\circ \phi_{p_{0},p}(P_{p},M_{p})=(P,M)$. 
    Cela entra\^{\i}ne que l'\'el\'ement $m=ygx_{p}^{-1}$ appartient \`a $M$. Parce que $g_{p}\in G_{p}({\mathbb R})$, on v\'erifie que la multiplication de $y$ par $g$ ne modifie pas l'\'egalit\'e de cocycles ci-dessus, c'est-\`a-dire que l'on a  $\nabla^M(\sigma)=yg\nabla_{p_{0},p}(\sigma)\sigma(yg)^{-1}$ pour tout $\sigma\in \Gamma_{{\mathbb R}}$.   Ou encore $\nabla^M(\sigma)=m\nabla_{p}^M(\sigma)\sigma(m)^{-1}$. Donc $\nabla^M$ a m\^eme classe dans $H^1(\Gamma_{{\mathbb R}};M)$ que $\nabla^M_{p}$. Cela ach\`eve la preuve de (1). $\square$

On doit d\'ecrire comme en 3.1 la correspondance entre classes de conjugaison de paires paraboliques de $K\tilde{G}$ et classes de conjugaison de paires paraboliques de $\hat{G}$. Dans le cas non tordu, cette correspondance est d\'ecrite par le lemme 2.1 de [A1]. A priori, celui-ci ne s'applique pas dans le cas g\'en\'eral car, comme on l'a dit en 1.11, notre notion de $K$-groupes est plus restrictive que celle d'Arthur. Nous allons prouver que ce lemme reste malgr\'e tout valable. Fixons une paire de Borel \'epingl\'ee $\hat{{\cal E}}=(\hat{B},\hat{T},(\hat{E}_{\hat{\alpha}})_{\hat{\alpha}\in \hat{\Delta}})$ de $\hat{G}$. On suppose qu'elle est stable par l'action galoisienne et on fixe un \'el\'ement $\hat{\theta}$ relatif \`a cette paire. On note $\sigma\mapsto \sigma_{G^*}$ l'action galoisienne. Les sous-groupes paraboliques standard $\hat{P}=\hat{M}\hat{U}$ qui sont stables par $\hat{\theta}$ et par l'action galoisienne sont en bijection avec les sous-ensembles $\hat{\Delta}^{\hat{M}}$ de $\hat{\Delta}$ qui v\'erifient les m\^emes propri\'et\'es de stabilit\'e ($\hat{\Delta}^{\hat{M}}$ est l'ensemble des racines de $\hat{T}$ dans $\hat{M}$). 

D'autre part, fixons une composante de notre $K$-espace $K\tilde{G}$, que l'on note simplement $\tilde{G}$. Fixons une paire de Borel \'epingl\'ee ${\cal E}=(B,T,(E_{\alpha})_{\alpha\in \Delta})$ de $G$ et fixons une cocha\^{\i}ne $\sigma\mapsto u(\sigma)$ de $\Gamma_{{\mathbb R}}$ dans $G_{SC}$ de sorte que $ad_{u(\sigma)}\circ\sigma_{G}({\cal E})={\cal E}$ (o\`u $\sigma\mapsto \sigma_{G}$ est l'action naturelle). On d\'efinit l'action quasi-d\'eploy\'ee $\sigma\mapsto \sigma_{G^*}=ad_{u(\sigma)}\circ\sigma_{G}$ de $\Gamma_{{\mathbb R}}$ sur $G$ et, pour simplifier, on note $G^*$ le groupe $G$ muni de cette action. On note $\theta^*$ l'automorphisme $ad_{e}$ pour un \'el\'ement $e\in Z(\tilde{G},{\cal E})$ quelconque. Cet automorphisme pr\'eserve ${\cal E}$ et l'action galoisienne quasi-d\'eploy\'ee. La bijection naturelle $\alpha\mapsto \hat{\alpha}$ de $\Delta$ sur $\hat{\Delta}$ est \'equivariante pour les actions galoisiennes et \'echange l'action de $\theta^*$ avec celle de $\hat{\theta}$. Posons $u^*(\sigma)=u(\sigma)^{-1}$ et notons $u^*_{ad}(\sigma)$ l'image de $u^*(\sigma)$ dans $G^*_{AD}$. On v\'erifie que $u^*_{ad}$ est un cocycle, qui d\'efinit un \'el\'ement de $H^1(\Gamma_{{\mathbb R}};G^*_{AD})$ not\'e encore $u^*_{ad}$. On a une application naturelle
$$H^1(\Gamma_{{\mathbb R}};G^*_{AD})\to H^2(\Gamma_{{\mathbb R}};Z(G^*_{SC})).$$
Ce dernier groupe s'identifie facilement au groupe des caract\`eres de $Z(\hat{G}_{SC})^{\Gamma_{{\mathbb R}}}$ qui sont triviaux sur l'image de la norme
$$Z(\hat{G}_{SC})\to Z(\hat{G}_{SC})^{\Gamma_{R}}.$$
On renvoie pour cela \`a [K2], th\'eor\`eme 1.2. Ainsi, $u^*_{ad}$ d\'efinit un caract\`ere $\chi_{K\tilde{G}}$ de $Z(\hat{G}_{SC})^{\Gamma_{{\mathbb R}}}$.  On a fait divers choix, qui affectent m\^eme notre construction de $G^*$. Quand on change de choix, on voit que les deux groupes $G^*$ construits s'identifient naturellement et que le caract\`ere $\chi_{K\tilde{G}}$ obtenu est le m\^eme. C'est facile \`a voir pourvu que l'on conserve la m\^eme composante connexe  $\tilde{G}$. Consid\'erons une autre composante $\tilde{G}'$. Par d\'efinition, il y a un isomorphisme $\phi:G'\to G$ et un cocycle $\nabla\in H^1(\Gamma_{{\mathbb R}};G_{SC})$ tel que $\phi\circ \sigma(\phi)^{-1}=ad_{\nabla(\sigma)}$ pour tout $\sigma\in\Gamma_{{\mathbb R}}$. On prend pour paire de Borel \'eping\'ee ${\cal E}'=\phi^{-1}({\cal E})$. On v\'erifie que l'on  peut choisir $u'(\sigma)= \phi^{-1}(u(\sigma)\nabla(\sigma))$. Il est clair que $\phi$ d\'efinit un isomorphisme d\'efini sur ${\mathbb R}$ de  $G^{'*}$ sur $G^*$. Via cet isomorphisme, $u^{'*}(\sigma)$ s'identifie \`a $\nabla(\sigma) u^*(\sigma)$. Le calcul montre que la condition que $\nabla$ est un cocycle (pour l'action naturelle sur $G$) \'equivaut \`a ce que $d(\nabla u^*)=d(u^*)$, o\`u $d$ est la diff\'erentielle sur $G^*$. Les images de $u^*$ et $\nabla u^*$ dans $H^2(\Gamma_{{\mathbb R}};Z(G^*_{SC}))$ sont donc les m\^emes et on r\'ecup\`ere ainsi le m\^eme caract\`ere $\chi_{K\tilde{G}}$. Remarquons que, par hypoth\`ese, $\tilde{G}({\mathbb R})$ est non vide. On peut donc fixer $\gamma\in \tilde{G}({\mathbb R})$. Ecrivons $\gamma=ge$, avec $g\in G$ et $e\in Z(\tilde{G},{\cal E})$. Pour tout $\sigma\in \Gamma_{{\mathbb R}}$, on a encore $ad_{u(\sigma)}\circ \sigma_{G}(e)\in Z(\tilde{G},{\cal E})$, donc il existe $z(\sigma)\in Z(G)$ tel que $ ad_{u(\sigma)}\circ \sigma_{G}(e)=z(\sigma)^{-1}e$. La condition $\gamma\in \tilde{G}({\mathbb R})$ \'equivaut \`a ce que, pour tout $\sigma\in \Gamma_{{\mathbb R}}$, on ait l'\'egalit\'e
$\sigma_{G}(\gamma)=\gamma$. Or on a les \'equivalences suivantes
$$\sigma_{G}(\gamma)=\gamma \iff \sigma_{G}(g)\sigma_{G}(e)=ge \iff g^{-1}\sigma_{G}(g)ad_{u(\sigma)^{-1}}(z(\sigma)^{-1}e)=e $$
$$\iff g^{-1}\sigma_{G}(g)u(\sigma)^{-1}\theta^*(u(\sigma))z(\sigma)^{-1}e=e \iff g^{-1}\sigma_{G}(g)u(\sigma)^{-1}\theta^*(u(\sigma))z(\sigma)^{-1}=1$$
$$  \iff g^{-1}u^*(\sigma)\sigma_{G^*}(g)
=z(\sigma)\theta^*(u^*(\sigma)).$$
Il en r\'esulte que la classe du cocycle $u^*_{ad}$ est invariante par $\theta^*$, donc $\chi_{K\tilde{G}}$ est invariant par $\hat{\theta}$.

 Pour $x_{*}\in X_{*}(\hat{T}_{ad})$, choisissons un entier $N\geq1$ tel que $Nx_{*}\in X_{*}(\hat{T}_{sc})$. Alors l'\'el\'ement $Nx_{*}(e^{2\pi i/N})$ appartient \`a $Z(\hat{G}_{SC})$ et ne d\'epend pas du choix de $N$. L'application $x_{*}\mapsto Nx_{*}(e^{2\pi i/N})$ se quotient en un isomorphisme
 $$ X_{*}(\hat{T}_{ad})/X_{*}(\hat{T}_{sc})\simeq Z(\hat{G}_{SC}).$$ 
 
A tout \'el\'ement $\alpha\in\Delta$ est naturellement associ\'e un copoids $\varpi_{\alpha}\in X_{*}(\hat{T}_{ad})$. On note $\varpi_{\alpha}^{\Gamma_{{\mathbb R}}}$ la somme des \'el\'ements $\varpi_{\alpha'}$ pour les $\alpha'$ dans l'orbite de $\alpha$ sous l'action de $\Gamma_{{\mathbb R}}$ (puisque ce groupe a deux \'el\'ements, les orbites ont au plus deux \'el\'ements). L'\'el\'ement $\varpi_{\alpha}^{\Gamma_{{\mathbb R}}}$ s'envoie sur un \'el\'ement de $Z(\hat{G}_{SC})^{\Gamma_{{\mathbb R}}}$. On note $\Delta_{min}$ l'ensemble des $\alpha\in \Delta$ tels que $\chi_{K\tilde{G}}(\varpi_{\alpha}^{\Gamma_{{\mathbb R}}})\not=1$.  On note $\hat{\Delta}_{min} $ l'ensemble des $\hat{\alpha}$ pour $\alpha\in \Delta_{min}$. Cet ensemble est stable par l'action galoisienne  et aussi par $\hat{\theta}$ puisque $\chi_{\tilde{G}}$ l'est. 

\ass{Lemme}{ Soit $\hat{P}=\hat{M}\hat{U}$ un sous-groupe parabolique standard de $\hat{G}$ stable par $\hat{\theta}$ et par l'action galoisienne. Alors $\hat{P}$ correspond \`a une classe de conjugaison de sous-$K$-espaces paraboliques de $K\tilde{G}$ si et seulement si $\hat{\Delta}^{\hat{M}}$ contient $\hat{\Delta}_{min}$.} 

C'est exactement l'\'enonc\'e du lemme 2.1 de [A1]. Nous le prouverons dans le paragraphe suivant.

Il r\'esulte de ce lemme que
 
 (3) parmi les classes de conjugaison par $KG({\mathbb R})$ de paires paraboliques de $KG$, il y a une unique classe minimale.
 
Une propri\'et\'e \'equivalente est qu'il y a au moins un $p\in {\Pi}$ tel que $G_{p}$ soit "plus quasi-d\'eploy\'e" que les autres composantes. 

On doit d\'efinir correctement les espaces ${\cal L}(K\tilde{M})$, ${\cal P}(K\tilde{M})$ et ${\cal F}(K\tilde{M})$ pour un $K$-espace de Levi $K\tilde{M}$. Si l'on d\'efinit ${\cal L}(K\tilde{L})$ comme l'ensemble des $K$-espaces de Levi de $K\tilde{G}$ contenant $K\tilde{M}$, il y en a beaucoup trop. Pour cela, on fixe pour tout $p\in \Pi$ une paire parabolique minimale $(P_{p,0},M_{p,0})$, qui donne naissance \`a une paire d'espaces tordus $(KP_{p,0},KM_{p,0})$. Le r\'esultat pr\'ec\'edent entra\^{\i}ne qu'il existe un unique sous-ensemble non vide $\Pi^{M_{0}}$ de $\Pi$ v\'erifiant les deux conditions suivantes:

- la famille $K\tilde{M}_{0}=(\tilde{M}_{p,0})_{p\in \Pi^{M_{0}}}$ est un $K$-espace de Levi de $K\tilde{G}$;

- pour tous $p\in \Pi$, $p'\in \Pi^{M_{0}}$, il existe $x_{p',p}\in G_{p'}$ tel que $ad_{x_{p',p}}\circ \tilde{\phi}_{p',p}(\tilde{P}_{p,0},\tilde{M}_{p,0})$ contienne $(\tilde{P}_{p',0},\tilde{M}_{p',0})$.

On fixe de tels \'el\'ements $x_{p',p}$. La construction suivante ne d\'ependra pas de leur choix. Il est facile de montrer que, pour tout $K$-espace de Levi $K\tilde{L}=(\tilde{L}_{p})_{p\in \Pi^L}$ de $K\tilde{G}$, l'ensemble d'indices $\Pi^L$ contient $\Pi^{M_{0}}$. On note ${\cal L}(K\tilde{M}_{0})$ l'ensemble des $K$-espaces de Levi $K\tilde{L}=(\tilde{L}_{p})_{p\in \Pi^L}$ de $K\tilde{G}$ v\'erifiant les deux conditions suivantes:

- $\tilde{L}_{p}\supset \tilde{M}_{p,0}$ pour tout $p\in \Pi^L$;

- $ad_{x_{p',p}}\circ \tilde{\phi}_{p',p}(\tilde{L}_{p})=\tilde{L}_{p'}$ pour tous $p\in \Pi^L$, $p'\in \Pi^{M_{0}}$.

Pour $K\tilde{M}=(\tilde{M}_{p})_{p\in \Pi^M}\in {\cal L}(K\tilde{M}_{0})$, on note ${\cal L}(K\tilde{M})$ l'ensemble des $K\tilde{L}=(\tilde{L}_{p})_{p\in \Pi^L}\in {\cal L}(K\tilde{M}_{0})$ tels que $\Pi^M\subset \Pi^L$ et $\tilde{M}_{p}\subset \tilde{L}_{p}$ pour tout $p\in \Pi^M$. On d\'efinit de fa\c{c}on similaire les ensembles ${\cal P}(K\tilde{M})$ et ${\cal F}(K\tilde{M})$.

Les consid\'erations des quatre paragraphes pr\'ec\'edents s'adaptent aux objets d\'efinis ci-dessus. Du c\^ot\'e dual, il faut bien s\^ur prendre pour paire $(\hat{P}_{0},\hat{M}_{0})$ une paire qui correspond \`a $(KP_{0},KM_{0})$. 

\bigskip
\subsection{Preuve du lemme 3.5}
La n\'ecessit\'e de la condition r\'esulte du lemme d'Arthur. Nos $K$-groupes peuvent se compl\'eter en $K$-groupes au sens d'Arthur. Si un sous-groupe parabolique $\hat{P}=\hat{M}\hat{U}$ (standard, invariant par $\hat{\theta}$ et par l'action galoisienne) correspond \`a une classe de conjugaison de sous-$K$-espaces paraboliques de $K\tilde{G}$, il correspond a fortiori \`a une classe de conjugaison de sous-$K$-groupes paraboliques de ce $K$-groupe \'etendu, donc v\'erifie l'inclusion $\hat{\Delta}^{\hat{M}}\supset \hat{\Delta}_{min}$. 

Pour la r\'eciproque, il suffit de traiter l'unique sous-groupe parabolique $\hat{P}=\hat{M}\hat{U}$ tel que $\hat{\Delta}^{\hat{M}}=\hat{\Delta}_{min}$. En effet, si celui-ci correspond bien \`a une classe de conjugaison de sous-$K$-espaces paraboliques de $K\tilde{G}$, on peut fixer une composante $\tilde{G}$ de $K\tilde{G}$ et un sous-espace parabolique $\tilde{P}$ de $\tilde{G}$ correspondant \`a $\hat{P}$. Les consid\'erations de 3.1 s'appliquent \`a cette composante. En particulier, tout sous-groupe parabolique $\hat{P}' $ contenant $\hat{P}$ et invariant par $\hat{\theta}$ et par l'action galoisienne correspond \`a un sous-espace parabolique $\tilde{P}'$ de $\tilde{G}$ contenant $\tilde{P}$. Dor\'enavant, on note $\hat{P}$ le sous-groupe "minimal" d\'efini ci-dessus.

Montrons que l'on peut se ramener au cas o\`u $K\tilde{G}$ n'a pas d'autre espace de Levi que lui-m\^eme. En effet, supposons qu'il existe un espace parabolique propre $K\tilde{Q}$, de Levi $K\tilde{L}$. Il correspond \`a $K\tilde{L}$ un sous-ensemble $\Delta^L$ de $\Delta$, d'o\`u un sous-ensemble $\hat{\Delta}^L$ de $\hat{\Delta}$. Rempla\c{c}ant dans les constructions $K\tilde{G}$ par $K\tilde{L}$, on d\'efinit un sous-ensemble $\hat{\Delta}^L_{min}$ de $\hat{\Delta}^L$. Si on suppose l'assertion prouv\'ee pour $K\tilde{L}$, il correspond \`a ce sous-ensemble $\hat{\Delta}^L_{min}$ un sous-espace parabolique de $K\tilde{L}$, d'o\`u aussi un sous-espace parabolique de $K\tilde{G}$. Pour obtenir l'assertion cherch\'ee pour $\tilde{G}$, il suffit de prouver l'\'egalit\'e

(1) $\Delta_{min}=\Delta^L_{min}$.

Par le sens d\'ej\`a prouv\'e du lemme, on a en tout cas $\Delta_{min}\subset \Delta^L$. En affectant des exposants $L$ aux termes construits \`a l'aide de $K\tilde{L}$, les d\'efinitions nous ram\`enent \`a prouver l'\'egalit\'e 

(2) $\chi_{K\tilde{G}}(\varpi_{\alpha}^{\Gamma_{{\mathbb R}}})=\chi_{K\tilde{L}}(\varpi_{\alpha}^{L,\Gamma_{{\mathbb R}}})$ pour tout $\alpha\in \Delta^L$. 

 Fixons une composante $\tilde{L}$ de $K\tilde{L}$, qui est incluse dans une composante $\tilde{G}$ de $K\tilde{G}$. On utilise ces composantes pour effectuer les constructions du paragraphe pr\'ec\'edent, en les affectant d'exposants $G$ ou $L$. On suppose que $\tilde{L}$ est standard pour la paire de Borel \'epingl\'ee ${\cal E}$ et on prend pour paire de Borel \'epingl\'ee ${\cal E}^L$ la restriction de ${\cal E}$. On peut alors supposer que $u^*(\sigma)$ est le produit d'un \'el\'ement de $Z(L_{sc})$ et de l'image de $u^{*L}(\sigma)\in L^*_{SC}$  dans $G^*_{SC}$. Alors $u^*$ est une cocha\^{\i}ne \`a valeurs dans $L_{sc}$, qui d\'efinit un \'el\'ement de $H^1(\Gamma_{{\mathbb R}};L_{ad})$ que l'on note $v^*$. On a des applications naturelles
 $$\begin{array}{ccc}&&H^1(\Gamma_{{\mathbb R}};G_{AD})\\ &\nearrow&\\ H^1(\Gamma_{{\mathbb R}};L_{ad})\\&\searrow&\\ &&H^1(\Gamma_{{\mathbb R}};L_{AD}).\\ \end{array}$$
 L'\'el\'ement $v^*$ s'envoie sur $u^*_{ad}$ par la fl\`eche du haut et sur $u^{*L}_{ad}$ par celle du bas. D'apr\`es [K2] th\'eor\`eme 1.2, $v^*$ d\'efinit un caract\`ere $\chi$ de $Z(\hat{L}_{sc})^{\Gamma_{{\mathbb R}}}/Z(\hat{L}_{sc})^{\Gamma_{{\mathbb R}},0}$. On a un diagramme dual
 $$\begin{array}{ccc}Z(\hat{G}_{SC})^{\Gamma_{{\mathbb R}}}&&\\ &\searrow&\\ && Z(\hat{L}_{sc})^{\Gamma_{{\mathbb R}}}/Z(\hat{L}_{sc})^{\Gamma_{{\mathbb R}},0}\\&\nearrow&\\ Z(\hat{L}_{SC})^{\Gamma_{{\mathbb R}}}&&\\ \end{array}$$
 Le caract\`ere $\chi_{K\tilde{G}}$ est compos\'e de $\chi$ et de la fl\`eche du haut tandis que $\chi_{K\tilde{L}}$ est compos\'e de $\chi$ et de la fl\`eche du bas. Cela nous ram\`ene \`a prouver que, pour $\alpha\in \Delta^L$, les images dans $Z(\hat{L}_{sc})^{\Gamma_{{\mathbb R}}}/Z(\hat{L}_{sc})^{\Gamma_{{\mathbb R}},0}$ de $\varpi_{\alpha}^{\Gamma_{{\mathbb R}}}$ et de $\varpi_{\alpha}^{L,\Gamma_{{\mathbb R}}}$ sont \'egales. Ecrivons $\varpi_{\alpha}^{\Gamma_{{\mathbb R}}}\in X_{*}(\hat{T}_{ad})$ sous la forme $\frac{1}{N}(x_{*}+y_{*})$, o\`u $N$ est un entier strictement positif, $x_{*}\in X_{*}(Z(\hat{L}_{sc})^0)$ et $y_{*}\in X_{*}(\hat{T}^L_{sc})$. Ici $ \hat{T}^L_{sc}$ est l'image r\'eciproque de $\hat{T}$ dans $\hat{L}_{SC}$. Le groupe $X_{*}(\hat{T}^L_{sc})$ est engendr\'e par les \'el\'ements de $\Delta^L$ (un \'el\'ement $\beta\in \Delta$ \'etant identifi\'e \`a la coracine associ\'ee \`a $\hat{\beta}\in \hat{\Delta}$). Il r\'esulte des d\'efinitions que $\varpi_{\alpha}^{L,\Gamma_{{\mathbb R}}}=\frac{1}{N}y_{*}$ et que $x_{*}$ est invariant par $\Gamma_{{\mathbb R}}$. Par d\'efinition, l'\'el\'ement de $Z(\hat{G}_{SC})$ correspondant \`a $\varpi_{\alpha}^{\Gamma_{{\mathbb R}}}$ est $x_{*}(\zeta)y_{*}(\zeta)$, o\`u $\zeta=e^{2\pi i/N}$ tandis que l'\'el\'ement de $Z(\hat{L}_{SC})$ correspondant \`a $\varpi_{\alpha}^{L,\Gamma_{{\mathbb R}}}$ est $y^*(\zeta)$. Quand on pousse ces \'el\'ements dans $Z(\hat{L}_{sc})$, ces deux \'el\'ements diff\`erent par $x_{*}(\zeta)$, qui appartient \`a $Z(\hat{L}_{sc})^{\Gamma_{{\mathbb R}},0}$. Cela prouve (2) et (1).
 
 On suppose d\'esormais que $K\tilde{G}$ n'a pas d'autre espace de Levi que lui-m\^eme. Remarquons qu'il revient au m\^eme de supposer que, pour chaque composante $\tilde{G}$, le groupe $G$ lui-m\^eme n'a pas de groupe de Levi propre. On a vu en effet qu'un groupe de Levi minimal donnait naissance \`a un espace de Levi. Remarquons aussi que, sous notre hypoth\`ese, la propri\'et\'e \`a prouver est l'\'egalit\'e  $\Delta_{min}=\Delta$.
 
 Montrons maintenant que l'on peut supposer que $G$ est simplement connexe. En effet, fixons une composante $\tilde{G}$ de $K\tilde{G}$ et un \'el\'ement $\gamma\in \tilde{G}({\mathbb R})$. L'automorphisme $ad_{\gamma}$ se rel\`eve en un automorphisme de $G_{SC}$. On peut introduire un espace tordu $\tilde{G}_{SC}$ sur $G_{SC}$, que l'on note formellement $G_{SC}\gamma_{sc}$, de la fa\c{c}on suivante. La multiplication \`a gauche est \'evidente. Celle de droite est d\'efinie par $g_{sc}\gamma_{sc}x_{sc}=g_{sc}ad_{\gamma}(x_{sc})\gamma_{sc}$. Enfin l'action galoisienne est $\sigma(g_{sc}\gamma_{sc})=\sigma(g_{sc})\gamma_{sc}$. L'application $\tilde{G}_{SC}\to \tilde{G}$ d\'efinie par $g_{sc}\gamma_{sc}\mapsto \pi(g_{sc})\gamma$ est un homomorphisme d'espaces tordus en un sens \'evident. On peut compl\'eter $\tilde{G}_{SC}$ en un $K$-espace $K\tilde{G}_{SC}$ et on v\'erifie que l'application pr\'ec\'edente s'\'etend en un homomorphisme $K\tilde{G}_{SC}\to K\tilde{G}$ (remarquons toutefois que l'application qui s'en d\'eduit entre les ensembles de composantes connexes de ces espaces n'est en g\'en\'eral ni injective, ni surjective).  Il est clair que l'hypoth\`ese sur $K\tilde{G}$ est  aussi v\'erifi\'ee pour $K\tilde{G}_{SC}$:  $G_{SC}$  et les autres groupes de $K\tilde{G}_{SC}$ n'ont pas d'autres groupes de  Levi qu'eux-m\^emes. L'ensemble $\Delta_{min}$ ne change pas puisque n'interviennent dans sa d\'efinition que les groupes $\hat{G}_{SC}$ et $\hat{G}_{AD}$ qui n'ont pas chang\'e. Si on suppose d\'emontr\'ee  l'assertion pour $K\tilde{G}_{SC}$, on conclut $\Delta_{min}=\Delta$, ce qui est la m\^eme assertion que pour $K\tilde{G}$.
 
 On suppose d\'esormais que $G$ est simplement connexe. On conserve toutefois la notation $G_{SC}$ quand elle est  plus suggestive. On fixe une composante $\tilde{G}$ de $K\tilde{G}$ et un \'el\'ement $\gamma\in \tilde{G}({\mathbb R})$ fortement r\'egulier. On choisit une paire de Borel \'epingl\'ee ${\cal E}$ de $G=G_{SC}$ dont le tore sous-jacent $T=T_{sc}$ est conserv\'e par $ad_{\gamma}$. On utilise cette paire de Borel \'epingl\'ee dans les constructions du paragraphe pr\'ec\'edent. Le tore est  d\'efini sur ${\mathbb R}$ pour l'action naturelle comme pour l'action quasi-d\'eploy\'ee. Il en r\'esulte que $u^*(\sigma)$ normalise $T_{sc}$ pour tout $\sigma\in \Gamma_{{\mathbb R}}$.  N\'ecessairement, son image $w(\sigma)$ dans le groupe de Weyl $W$ est invariante par $\theta^*$. L'hypoth\`ese que $G$ n'a pas d'espace de Levi propre entra\^{\i}ne que $T=T_{sc}$ est elliptique. En notant $\boldsymbol{\sigma}$ l'unique \'el\'ement non trivial de $\Gamma_{{\mathbb R}}$, $w(\boldsymbol{\sigma})\circ\boldsymbol{\sigma}_{G^*}$ agit donc par $-1$ sur $X_{*}(T_{sc})$. Il en r\'esulte que $w(\boldsymbol{\sigma})$ envoie toute racine positive sur une racine n\'egative. C'est donc l'\'el\'ement de $W$ de plus grande longueur, que l'on note ${\bf w}$. Introduisons la section de Springer $n:W\to G_{SC}$, cf. [LS] 2.1. A ce point, on a prouv\'e que l'on pouvait supposer supposer
 $$u^*(1)=1, \,\, u^*(\boldsymbol{\sigma})=tn({\bf w}),$$
 pour un \'el\'ement $t\in T_{sc}$. Soit $\alpha\in \Delta$. On dispose d\'ej\`a de l'\'el\'ement $E_{\alpha}$ de l'\'epinglage. On introduit l'\'el\'ement $E_{-\alpha}$ de l'espace radiciel de $\mathfrak{g}$ associ\'e \`a $-\alpha$, normalis\'e de sorte que $[E_{\alpha},E_{-\alpha}]=\check{\alpha}$, en identifiant la coracine $\check{\alpha}$ \`a un \'el\'ement de $\mathfrak{t}$. Notons $G_{\alpha}$ le sous-groupe de $G$ engendr\'e par $T$ et les sous-groupes radiciels associ\'es \`a $\alpha$ et $-\alpha$. Puisque l'action galoisienne naturelle \'echange $\alpha$ et $-\alpha$, ce groupe est d\'efini sur ${\mathbb R}$. Puisque $G$ est semi-simple et n'a pas de Levi propre, $G({\mathbb R})$ est compact, donc aussi $G_{\alpha}({\mathbb R})$. Comme on le sait ([S2] paragraphe 2), cela implique qu'il existe des \'el\'ements $c_{\alpha},c_{-\alpha}\in {\mathbb C}^{\times}$ tels que $[c_{-\alpha}E_{\alpha},c_{\alpha}E_{\alpha}]=\check{\alpha}$ et $\boldsymbol{\sigma}_{G}(c_{\alpha}E_{\alpha})=-c_{-\alpha}E_{-\alpha}$. La premi\`ere relation dit que $c_{-\alpha}=c_{\alpha}^{-1}$. 
 
 Montrons que l'on a
 $$(3) \qquad ad_{n({\bf w})}\circ\boldsymbol{\sigma}_{G^*}(E_{\alpha})=-E_{-\alpha}.$$
 On a $ad_{n({\bf w})}\circ\boldsymbol{\sigma}_{G^*}(\check{\alpha})=-\check{\alpha}$ et il existe des nombres complexes non nuls $x$ et $y$ de sorte que $ad_{n({\bf w})}\circ\boldsymbol{\sigma}_{G^*}(E_{\alpha})=xE_{-\alpha}$, $ad_{n({\bf w})}\circ\boldsymbol{\sigma}_{G^*}(E_{-\alpha})=yE_{\alpha}$. Ces trois relations entra\^{\i}nent $xy=1$. Notons $s_{\alpha}$ la sym\'etrie relative \`a $\alpha$. ¬Par d\'efinition,
 $$n(s_{\alpha})=exp(X_{\alpha})exp(-X_{-\alpha})exp(X_{\alpha}).$$
 Un calcul matriciel entra\^{\i}ne l'\'egalit\'e $n({\bf w})\boldsymbol{\sigma}_{G^*}(n(s_{\alpha}))n({\bf w})^{-1}=\check{\alpha}(-x^{-1})n(s_{\alpha})$. Mais le lemme 2.1.A de [LS] entra\^{\i}ne $n({\bf w})n(s_{\alpha})n({\bf w})^{-1}=n(s_{\alpha})$. D'o\`u $\check{\alpha}(-x^{-1})=1$ et $x=-1$ puisque notre groupe est simplement connexe. Cela prouve (3). 
 
Il r\'esulte de (3) que $\boldsymbol{\sigma}_{G}(c_{\alpha}E_{\alpha})=-\alpha(t)^{-1}\overline{c_{\alpha}}E_{-\alpha}=-\alpha(t)^{-1}(c_{\alpha}\overline{c_{\alpha}})c_{-\alpha}E_{-\alpha}$. La condition de compacit\'e nous dit donc que $\alpha(t)$ est un r\'eel positif, et cela pour tout $\alpha\in \Delta$. Cette propri\'et\'e implique que l'on peut trouver un \'el\'ement $t'=\prod_{\alpha\in \Delta}\check{\alpha}(t_{\alpha})$, avec des  $t_{\alpha}$ r\'eels positifs, tel que $(t')^2$  a m\^eme image que $t$ dans $T_{ad}$. Notons que $\boldsymbol{\sigma}_{G}(t')=(t')^{-1}$.  Donc $t=\zeta t'\boldsymbol{\sigma}_{G}(t')^{-1}$, avec $\zeta\in Z(G)=Z(G_{SC})$. Alors 
$$u^*(\boldsymbol{\sigma})=\zeta t'\boldsymbol{\sigma}_{G}(t')^{-1}n({\bf w})= \zeta t'n({\bf w})\boldsymbol{\sigma}_{G^*}(t')^{-1}.$$
En rempla\c{c}ant ${\cal E}$ par $ad_{t'}^{-1}({\cal E})$, on fait dispara\^{\i}tre le cobord et on obtient
$$u^*(\boldsymbol{\sigma})=\zeta n({\bf w}),$$
avec $\zeta\in Z(G_{SC})$. Mais on peut toujours multiplier notre cocha\^{\i}ne par une cocha\^{\i}ne \`a valeurs dans $Z(G_{SC})$. Cela nous ram\`ene au cas o\`u
$$u^*(1)=1,\,\, u^*(\boldsymbol{\sigma})=n({\bf w}).$$
Calculons le cobord $du^*$. On a $du^*(1,1)=du^*(\boldsymbol{\sigma},1)=du^*(1,\boldsymbol{\sigma})=1$ et $du^*(\boldsymbol{\sigma},\boldsymbol{\sigma})=n({\bf w})\boldsymbol{\sigma}_{G^*}(n({\bf w}))$. L'\'el\'ement ${\bf w}$ est invariant par l'action galoisienne et $n$ est \'equivariant pour cette action. Donc $\boldsymbol{\sigma}_{G^*}(n({\bf w}))=n({\bf w})$. En appliquant de nouveau le lemme 2.1.A de [LS], on obtient
$$du^*(\boldsymbol{\sigma},\boldsymbol{\sigma})=\prod_{\alpha>0}\check{\alpha}(-1),$$
o\`u le produit est pris sur toutes les racines de $T$ dans $G$ qui sont positives pour $B$. Il est d'usage de noter $2\check{\rho}$ la somme $\sum_{\alpha>0}\check{\alpha}$. On prendra garde \`a cette notation: $\check{\rho}$ n'est pas forc\'ement une somme de coracines \`a coefficients entiers, mais seulement \`a coefficients demi-entiers.    En tout cas, $\check{\rho}$ appartient \`a $X_{*}(T_{ad})$ car on sait que $<\alpha,\check{\rho}>=1$ pour tout $\alpha\in \Delta$. On peut \'ecrire de fa\c{c}on unique $2\check{\rho}$ comme somme d'un \'el\'ement de $2X_{*}(T_{sc})$ et d'un \'el\'ement
$$\check{\epsilon}=\sum_{\alpha\in \Delta}\epsilon_{\alpha}\check{\alpha},$$
avec des coefficients $\epsilon_{\alpha}$ \'egaux \`a $0$ ou $1$. On obtient $du^*(\boldsymbol{\sigma},\boldsymbol{\sigma})=(2\check{\rho})(-1)=\check{\epsilon}(-1)$. 

Rappelons comment on identifie un \'el\'ement de $H^2(\Gamma_{{\mathbb R}};Z(G_{SC}))$ \`a un caract\`ere de $Z(\hat{G}_{SC})^{\Gamma_{{\mathbb R}}}$. Tout d'abord, fixons un entier $N\geq1$ tel que $NX_{*}(T_{ad})\subset X_{*}(T_{sc})$ et une racine primitive d'ordre $N$ de l'unit\'e $\zeta\in {\mathbb C}^{\times}$. L'application $x_{*}\mapsto (Nx_{*})(\zeta)$ d\'efinie sur $X_{*}(T_{ad})$ se quotiente en un isomorphisme
$$X_{*}(T_{ad})/X_{*}(T_{sc})\simeq Z(G_{SC}).$$
Il n'est pas \'equivariant par l'action galoisienne: puisque $\boldsymbol{\sigma}(\zeta)=\zeta^{-1}$, l'isomorphisme transporte l'action de $\boldsymbol{\sigma}$ en l'oppos\'e de cette action. Un \'el\'ement  de $H^2(\Gamma_{{\mathbb R}};Z(G_{SC}))$ peut toujours se repr\'esenter par une cocha\^{\i}ne $v$ v\'erifiant comme ci-dessus $v(1,1)=v(1,\boldsymbol{\sigma})=v(\boldsymbol{\sigma},1)=1$. L'\'el\'ement ${\bf v}=v(\boldsymbol{\sigma},\boldsymbol{\sigma})$ v\'erifie ${\bf v}=\boldsymbol{\sigma}({\bf v})=1$ (par la condition de cocycle) et s'identifie donc \`a un \'el\'ement $x\in X_{*}(T_{ad})/X_{*}(T_{sc})$ tel que $x\boldsymbol{\sigma}(x)=1$. On voit  que $x$ est uniquement d\'etermin\'e par la classe de $v$ modulo un \'el\'ement de la forme $y\boldsymbol{\sigma}(y)^{-1}$. Puisque $X_{*}(T_{ad})$ est le dual de $X_{*}(\hat{T}_{sc})$ et $X_{*}(T_{sc})$ est le dual de $X_{*}(\hat{T}_{ad})$, les deux groupes $X_{*}(T_{ad})/X_{*}(T_{sc})$ et $X_{*}(\hat{T}_{ad})/X_{*}(\hat{T}_{sc})\simeq Z(\hat{G}_{SC})$ sont duaux. Donc $x$ d\'efinit un caract\`ere de $Z(\hat{G}_{SC})$. La restriction de ce caract\`ere au sous-groupe $Z(\hat{G}_{SC})^{\Gamma_{{\mathbb R}}}$ ne change pas si on multiplie $x$ par un \'el\'ement de la forme $y\boldsymbol{\sigma}(y)^{-1}$. Cette restriction ne d\'epend donc que de $v$. C'est le caract\`ere associ\'e \`a $v$. 

Appliqu\'ee \`a $du^*$, cette construction nous dit que le caract\`ere $\chi_{K\tilde{G}}$ s'identifie au caract\`ere de $X_{*}(\hat{T}_{ad})/X_{*}(\hat{T}_{sc})$ associ\'e \`a l'\'el\'ement $\check{\rho}\in X_{*}(T_{ad})$.  Par d\'efinition, l'ensemble $\Delta_{min}$ est alors la r\'eunion de

- l'ensemble des $\alpha\in \Delta$ tels que $\boldsymbol{\sigma}_{G^*}(\alpha)=\alpha$ et $<\varpi_{\alpha},\check{\rho}>\not\in {\mathbb Z}$;

- l'ensemble des $\alpha\in \Delta$ tels que $\boldsymbol{\sigma}_{G^*}(\alpha)\not=\alpha$ et $<\varpi_{\alpha},\check{\rho}>+<\boldsymbol{\sigma}_{G^*}(\varpi_{\alpha}),\check{\rho}>\not\in {\mathbb Z}$.

L'\'el\'ement $\check{\rho}$ est invariant par l'action galoisienne et son produit avec tout \'el\'ement $\varpi_{\alpha}$ appartient \`a $\frac{1}{2}{\mathbb Z}$. Le second ensemble ci-dessus est donc vide. D'autre part la condition  $<\varpi_{\alpha},\check{\rho}>\not\in {\mathbb Z}$ \'equivaut \`a $\epsilon_{\alpha}=1$. On obtient que $\Delta_{min}$ est form\'e d'\'el\'ements fixes par l'action galoisienne et que l'on a une \'egalit\'e
$$\check{\epsilon}=(\sum_{\alpha\in \Delta_{min}}\check{\alpha})+(\sum_{\alpha\in \Delta'}\check{\alpha}+\boldsymbol{\sigma}_{G^*}(\check{\alpha})),$$
o\`u $\Delta'$ est un certain sous-ensemble de $\Delta-\Delta_{min}$ form\'e d'\'el\'ements $\alpha$ tels que $\boldsymbol{\sigma}_{G^*}(\alpha)\not=\alpha$. Remarquons que, puisque $\check{\rho}$ est invariant par $\theta^*$, $\check{\epsilon}$ l'est aussi. Donc $\Delta_{min}$ l'est (ce qui \'etait d\'ej\`a \'evident) ainsi que l'ensemble $\Delta'\sqcup \boldsymbol{\sigma}_{G^*}(\Delta')$.

Reprenons les calculs  effectu\'es dans le paragraphe pr\'ec\'edent. On peut  \'ecrire $\gamma=te$, avec $t\in T$ et $e\in Z(\tilde{G},{\cal E})$. Comme on l'a dit, on a pour tout $\sigma\in \Gamma_{{\mathbb R}}$ une \'egalit\'e
$$ad_{u(\sigma)}\circ\sigma_{G}(e)=z(\sigma)^{-1}e,$$
avec $z(\sigma)\in Z(G)=Z(G_{SC})$. Ou encore 
$$\sigma_{G}(e)=z(\sigma)^{-1}ad_{u^*(\sigma)}(e)=z(\sigma)^{-1}u^*(\sigma)\theta^*(u^*(\sigma))^{-1}e.$$
Mais $\theta^*(u^*(\sigma))=u^*(\sigma)$. La condition devient simplement $\sigma_{G}(e)=z(\sigma)^{-1}e$. Puisque $\gamma\in \tilde{G}({\mathbb R})$, on a $\sigma_{G}(te)=te$, ou encore
$\sigma_{G}(t)=z(\sigma)t$. Cela entra\^{\i}ne que l'image $t_{ad}$ de $t$ dans $T_{ad}$ appartient \`a $T_{ad}({\mathbb R})$. Mais $T_{ad}$ est elliptique. Donc $T_{ad}({\mathbb R})$ est connexe et l'application $\pi:T_{sc}({\mathbb R})\to T_{ad}({\mathbb R})$ est surjective. On peut donc \'ecrire $t=t_{0}\zeta$, avec $t_{0}\in T_{sc}({\mathbb R})$ et $\zeta\in Z(G_{SC})$.    Alors $\zeta e=t_{0}^{-1}\gamma\in \tilde{G}({\mathbb R})$. Quitte \`a remplacer $e$ par $\zeta e$, on a construit un \'el\'ement $e\in Z(\tilde{G},{\cal E})$ qui appartient \`a $\tilde{G}({\mathbb R})$. 

Traduisons maintenant ce que l'on cherche. On veut trouver un cocycle $\nabla:\Gamma_{{\mathbb R}}\to G_{SC}$ tel que sa classe dans $H^1(\Gamma_{{\mathbb R}};G_{SC})$ soit invariante par $\theta$ et tel que la condition suivante soit v\'erifi\'ee. Introduisons un groupe $G'$ sur ${\mathbb R}$ muni d'un isomorphisme $\phi:G'\to G$ de sorte que $\phi\circ\sigma(\phi)^{-1}=ad_{\nabla(\sigma)}$ pour tout $\sigma\in \Gamma_{{\mathbb R}}$. Notons $P^*=M^*U^*$ le sous-groupe parabolique standard de $G^*$ tel que l'ensemble de racines simples associ\'e \`a $M^*$ soit $\Delta_{min}$. On veut que $P^*$ se transf\`ere \`a $G'$. Comme on l'a dit dans le paragraphe pr\'ec\'edent,  que $\nabla$ soit un cocycle \`a valeurs dans $G_{SC}$ revient \`a dire que $d(\nabla u^*)=du^*$.  De plus, quand on remplace $G$ par $G'$, on remplace $u^*$ par $\nabla u^*$. La derni\`ere condition ci-dessus signifie que l'image de $\nabla u^*$ dans $G^*_{ad}$ est cohomologue \`a une cocha\^{\i}ne \`a valeurs dans $M^*_{ad}$. Traduisons la condition d'invariance par $\theta$. On se rappelle que cette action $\theta$ est l'action $ad_{\gamma}$ pour un \'el\'ement $\gamma\in \tilde{G}({\mathbb R})$. On peut choisir pour $\gamma$  l'\'el\'ement $e$ fix\'e ci-dessus. Alors $\theta=\theta^*$ et la condition signifie qu'il existe $g\in G_{SC}$ tel que $\theta^*(\nabla(\sigma))=g\nabla(\sigma)\sigma_{G}(g)^{-1}$ pour tout $\sigma\in \Gamma_{{\mathbb R}}$. Puisque $u^*(\sigma)$ est fixe par $\theta^*$, cette relation \'equivaut \`a 
$$\theta^*(\nabla(\sigma)u^*(\sigma))=g\nabla(\sigma)\sigma_{G}(g)^{-1}u^*(\sigma)=g\nabla(\sigma)u^*(\sigma)\sigma_{G^*}(g)^{-1}.$$
Supposons trouv\'e une cocha\^{\i}ne  $v^*:\Gamma_{{\mathbb R}}\to M^*_{sc}=M^*$ telle que

(4) $dv^*=du^*$;

(5) il existe $t\in T_{sc}$ tel que $\theta^*(v^*(\sigma))=tv^*(\sigma)\sigma_{G^*}(t)^{-1}$ pour tout $\sigma\in \Gamma_{{\mathbb R}}$.

Alors le cocycle $\nabla=v^*(u^*)^{-1}$ r\'epond \`a la question.  

Pour construire $v^*$, on a besoin de quelques remarques pr\'eliminaires concernant les ensembles $\Delta_{min}$ et $\Delta'$. Rappelons que $\Delta_{min}\sqcup \Delta'\sqcup \boldsymbol{\sigma}_{G^*}(\Delta')$ est l'ensemble des $\alpha\in \Delta$ tels que, quand on \'ecrit $2\check{\rho}=\sum_{\beta\in \Delta}c_{\beta}\beta$, le coefficient $c_{\alpha}$ soit impair. Or on sait calculer $2\check{\rho}$ pour chaque syst\`eme de racines irr\'eductible. On renvoie aux tables de Bourbaki ([Bour]). On s'aper\c{c}oit en consultant ces tables que $\Delta_{min}\sqcup \Delta'\sqcup \boldsymbol{\sigma}_{G^*}(\Delta')$ est form\'e de racines deux \`a deux orthogonales. Puisque de plus, $\boldsymbol{\sigma}_{G^*}$ fixe tout \'el\'ement de $\Delta_{min}$, il en r\'esulte que $M^*_{SC}$ est un produit de groupes $SL(2)$ index\'es par les racines $\alpha\in \Delta_{min}$. 
Introduisons l'\'el\'ement de plus grande longueur du groupe de Weyl de $M^*$, que l'on note $\boldsymbol{\omega}$. C'est simplement le produit des sym\'etries $s_{\alpha}$ associ\'ees aux $\alpha\in \Delta_{min}$ et  on a $\boldsymbol{\omega}(\alpha)=-\alpha$ pour tout $\alpha\in \Delta_{min}$. Enfin, puisque $\Delta'\sqcup \boldsymbol{\sigma}_{G^*}(\Delta')$ est orthogonal \`a $\Delta_{min}$, on a $\boldsymbol{\omega}(\alpha)=\alpha$ pour tout $\alpha\in \Delta'\sqcup \boldsymbol{\sigma}_{G^*}(\Delta')$. 

Introduisons l'\'el\'ement
$$x=\prod_{\alpha\in \Delta'}\check{\alpha}(-1)\in T_{sc}.$$
D'efinissons la cocha\^{\i}ne $v^*$ par $v^*(1)=1$ et $v^*(\boldsymbol{\sigma})=xn(\boldsymbol{\omega})$. Elle prend ses valeurs dans $M^*_{sc}$.  On va montrer qu'elle v\'erifie les conditions (4) et (5).

On a
$$dv^*(\boldsymbol{\sigma},\boldsymbol{\sigma})=xad_{n(\boldsymbol{\omega})}\circ\boldsymbol{\sigma}_{G^*}(x)^{-1}n(\boldsymbol{\omega})\boldsymbol{\sigma}_{G^*}(n(\boldsymbol{\omega})).$$
On a $\sigma_{G^*}(x)=\prod_{\alpha\in \Delta'}\boldsymbol{\sigma}(\check{\alpha})(-1)$. On a vu plus haut que toutes les coracines intervenant ici sont fixes par $\boldsymbol{\omega}$. D'o\`u
$$xad_{n(\boldsymbol{\omega})}\circ\boldsymbol{\sigma}_{G^*}(x)^{-1}=\prod_{\alpha\in \Delta'}\check{\alpha}(-1)\boldsymbol{\sigma}_{G^*}(\check{\alpha})(-1).$$
On calcule $n(\boldsymbol{\omega})\boldsymbol{\sigma}_{G^*}(n(\boldsymbol{\omega}))$ comme on a calcul\'e plus haut $n({\bf w})\boldsymbol{\sigma}_{G^*}(n({\bf w}))$. Ce terme vaut $(2\check{\rho}^{M^*})(-1)$, o\`u $2\check{\rho}^{M^*}$ est la somme des racines positives dans $M^*$. Puisque $M^*_{SC}$ est un produit de groupes $SL(2)$, on a simplement $2\check{\rho}^{M^*}=\sum_{\alpha\in \Delta_{min}}\check{\alpha}$. Cela conduit \`a l'\'egalit\'e
$$dv^*(\boldsymbol{\sigma},\boldsymbol{\sigma})=\check{\epsilon}(-1),$$
autrement dit
$$dv^*(\boldsymbol{\sigma},\boldsymbol{\sigma})=du^*(\boldsymbol{\sigma},\boldsymbol{\sigma}).$$
Cela v\'erifie la condition (4).

On a $x\boldsymbol{\sigma}_{G^*}(x)=\prod_{\alpha\in \Delta'\sqcup \boldsymbol{\sigma}_{G^*}(\Delta')}\check{\alpha}(-1)$. Or l'ensemble $\Delta'\sqcup \boldsymbol{\sigma}_{G^*}(\Delta')$ est invariant par $\theta^*$. Donc $x\boldsymbol{\sigma}_{G^*}(x)$ est invariant par $\theta^*$. Autrement dit, l'\'el\'ement $y=\theta^*(x)x^{-1}$ v\'erifie $y\boldsymbol{\sigma}_{G^*}(y)=1$. Consid\'erons le sous-tore $T''$ de $T_{sc}$ tel que $X_{*}(T'')=\Delta'\sqcup \boldsymbol{\sigma}_{G^*}(\Delta')$, muni de l'action $\sigma\mapsto \sigma_{G^*}$. C'est un tore induit donc $H^1(\Gamma_{{\mathbb R}};T'')=0$.  L'application $1\mapsto 1$, $\boldsymbol{\sigma}\mapsto y$ est un cocycle \`a valeurs dans ce tore, donc est un cobord. Il existe donc $t\in T''$ tel que $y=t\sigma_{G^*}(t)^{-1}$. Parce que $\boldsymbol{\omega}$ op\`ere trivialement sur $\Delta'\sqcup \boldsymbol{\sigma}_{G^*}(\Delta')$, $ad_{v^*(\boldsymbol{\sigma})}$ fixe $T''$. On a aussi bien $y=tad_{v^*(\boldsymbol{\sigma})}\circ\boldsymbol{\sigma}_{G^*}(t)^{-1}$. Autrement dit
$$\theta^*(x)x^{-1}=tv^*(\boldsymbol{\sigma})\sigma_{G^*}(t)^{-1}v^*(\boldsymbol{\sigma})^{-1},$$
ou encore
$$\theta^*(x)x^{-1}v^*(\boldsymbol{\sigma})=tv^*(\boldsymbol{\sigma})\sigma_{G^*}(t)^{-1},$$
ou encore
$$\theta^*(x)n(\boldsymbol{\omega})=tv^*(\boldsymbol{\sigma})\sigma_{G^*}(t)^{-1},$$
ou encore
$$\theta^*(v^*(\boldsymbol{\sigma}))=tv^*(\boldsymbol{\sigma})\sigma_{G^*}(t)^{-1},$$
puisque $n(\boldsymbol{\omega})$ est fixe par $\theta^*$. La relation pr\'ec\'edente \'equivaut \`a (5). Cela ach\`eve la d\'emonstration. $\square$

 \bigskip
 
 \section{Stabilit\'e et image du transfert}

\bigskip

\subsection{Rappels sur la descente d'Harish-Chandra et la transformation de Fourier}

  Le corps $F$ est de nouveau un corps local quelconque de caract\'eristique nulle. Dans les premiers paragraphes, on fixe des mesures de Haar pour se d\'ebarrasser des espaces de mesures.
  
   Oublions pour un temps les espaces tordus, c'est-\`a-dire supposons $\tilde{G}=G$, mais conservons le caract\`ere $\omega$.   Un certain nombre de d\'efinitions se descendent aux alg\`ebres de Lie, par exemple les int\'egrales orbitales. On utilise pour ces alg\`ebres des notations analogues \`a celles pour les groupes. 
  
      On introduit une transformation de Fourier $f\mapsto \hat{f}$ dans l'espace $C_{c}^{\infty}(\mathfrak{g}(F))$ relative \`a un bicaract\`ere invariant par conjugaison par $G(F)$ (en appelant conjugaison l'action adjointe). Cette transformation de Fourier conserve le noyau de l'homomorphisme $C_{c}^{\infty}(\mathfrak{g}(F))\to I(\mathfrak{g}(F),\omega)$, donc passe au quotient en une transformation $f\mapsto \hat{f}$ dans $I(\mathfrak{g}(F),\omega)$. D'autre part, pour tout Levi $M$ de $G$, on a une \'egalit\'e $(\hat{f})_{M,\omega}=(f_{M,\omega})\hat{}$. Cela entra\^{\i}ne que la transformation de Fourier conserve le sous-espace $C_{cusp}^{\infty}(\mathfrak{g}(F),\omega)\subset C_{c}^{\infty}(\mathfrak{g}(F))$ des fonctions $f$ telles que $f_{M,\omega}=0$ dans $I(\mathfrak{m}(F),\omega)$ pour tout Levi propre $M$ de $G$.
 
 Les propri\'et\'es suivantes r\'esultent d'une part de la conjecture de Howe (qui n'est plus une conjecture depuis longtemps), ou plut\^ot de sa variante concernant les int\'egrales orbitales tordues par $\omega$, d'autre part de l'int\'egrabilit\'e des transform\'ees  de Fourier d'int\'egrales orbitales.
        
  Soit $\mathfrak{u}$ un ouvert de $\mathfrak{g}_{reg}(F)$ dont l'adh\'erence contienne un voisinage de $0$. Alors
  
  (1) si $F$ est non-archim\'edien, pour tout $f\in C_{c}^{\infty}(\mathfrak{g}(F))$, il existe $f'\in C_{c}^{\infty}(\mathfrak{u})$ telle que les int\'egrales orbitales de $f$ et de $\hat{f}'$ co\"{\i}ncident dans un voisinage de $0$.
  
  Notons $\mathfrak{g}(F)_{ell}$ le sous-ensemble des \'el\'ements semi-simples r\'eguliers et elliptiques dans $\mathfrak{g}(F)$. Alors
  
  (2) si $F$ est non-archim\'edien,  pour tout $f\in C_{cusp}^{\infty}(\mathfrak{g}(F),\omega)$, il existe $f'\in C_{c}^{\infty}(\mathfrak{u}\cap \mathfrak{g}(F)_{ell})$  telle que les int\'egrales orbitales de $f$ et de $\hat{f}'$ co\"{\i}ncident dans un voisinage de $0$.
  
    Supposons donn\'e un groupe  $\Xi$ d'automorphismes de $G$, d\'efinis sur $F$ et conservant le caract\`ere $\omega$. Supposons que  l'image de $\Xi$ dans le groupe d'automorphismes ext\'erieurs de $G$ soit finie. On peut supposer que le bicaract\`ere utilis\'e pour d\'efinir la transformation de Fourier est invariant par $\Xi$. Alors la transformation de Fourier est \'equivariante pour l'action de $\Xi$. Dans les assertions pr\'ec\'edentes, si l'on suppose que $\mathfrak{u}$ est invariant par $\Xi$ et que l'image de $f$ dans $I(\mathfrak{g}(F),\omega)$ est fixe par $\Xi$, on peut imposer qu'il en est de m\^eme de celle de $f'$.

  Revenons au cas g\'en\'eral (on ne suppose plus $\tilde{G}=G$).   Soient $\eta\in \tilde{G}_{ss}(F)$ et $\mathfrak{u}$ un voisinage ouvert  de $0$ dans $\mathfrak{g}_{\eta}(F)$ v\'erifiant les deux conditions suivantes
  
  - $\mathfrak{u}$ est invariant par conjugaison par $Z_{G}(\eta,F)$;
  
  - si $X\in \mathfrak{u}$, alors sa partie semi-simple $X_{ss}$ appartient \`a $\mathfrak{u}$.
  
  On va \'enoncer des propri\'et\'es qui sont vraies pourvu que $\mathfrak{u}$ soit assez petit. En particulier, on suppose $\mathfrak{u}$  assez petit pour que l'exponentielle y soit d\'efinie.  On pose $U_{\eta}=exp(\mathfrak{u})\subset G_{\eta}(F)$. 
  Notons $\tilde{U}$ l'ensemble des \'el\'ements de $\tilde{G}(F)$ qui sont conjugu\'es par un \'el\'ement de $G(F)$ \`a un \'el\'ement  de $ U_{\eta}\eta$. C'est un ouvert de $\tilde{G}(F)$. Notons $I(\tilde{U},\omega)$ l'image  de $C_{c}^{\infty}(\tilde{U})$ dans $I(\tilde{G}(F),\omega)$, $I(U_{\eta},\omega)$ celle de $C_{c}^{\infty}(U_{\eta})$ dans $I(G_{\eta}(F),\omega)$ et $I(\mathfrak{u},\omega)$ celle de $C_{c}^{\infty}(\mathfrak{u})$ dans $I(\mathfrak{g}_{\eta}(F),\omega)$. L'exponentielle \'etablit un isomorphisme entre $I(U_{\eta},\omega)$ et $I(\mathfrak{u},\omega)$. 
  Remarquons que le groupe $Z_{G}(\eta;F)$ agit naturellement sur  $I(G_{\eta}(F),\omega)$ et $I(\mathfrak{g}_{\eta}(F),\omega)$. D\'efinissons une correspondance entre $C_{c}^{\infty}(\tilde{U})$ et $C_{c}^{\infty}(U_{\eta})$ par: $f\in C_{c}^{\infty}(\tilde{U})$ et $\phi\in C_{c}^{\infty}(U_{\eta})$ se correspondent si et seulement si on a l'\'egalit\'e $I^{\tilde{G}}(x\eta,\omega,f)=I^{G_{\eta}}(x,\omega,\phi)$ pour tout \'el\'ement r\'egulier $x\in U_{\eta}$ tel que $x\eta$ soit fortement r\'egulier dans $\tilde{G}$ (il est sous-entendu que les mesures sur $G_{x\eta}(F)=(G_{\eta})_{x}(F)$ qui interviennent dans la d\'efinition de ces int\'egrales orbitales sont  les m\^emes pour les deux int\'egrales). La th\'eorie de la descente affirme que cette correspondance se quotiente en un isomorphisme 
  $$desc_{\eta}^{\tilde{G}}:I(\tilde{U},\omega)\to I(U_{\eta},\omega)^{Z_{G}(\eta;F)},$$
  o\`u, selon l'usage, l'exposant $Z_{G}(\eta;F)$ d\'esigne le sous-espace d'invariants par ce groupe. Via l'exponentielle, on peut aussi consid\'erer que $desc_{\eta}^{\tilde{G}}$ prend ses valeurs dans $I(\mathfrak{u},\omega)^{Z_{G}(\eta,F)}$.
  
  Supposons $\eta$ elliptique. Alors le m\^eme r\'esultat vaut pour les fonctions cuspidales. C'est-\`a-dire, d\'efinissons $C_{cusp}^{\infty}(\tilde{U})=C_{cusp}^{\infty}(\tilde{G}(F))\cap C_{c}^{\infty}(\tilde{U})$, notons $I_{cusp}(\tilde{U},\omega)$ son image dans $I_{cusp}(\tilde{G}(F),\omega)$ et d\'efinissons de m\^eme $C_{cusp}^{\infty}(U_{\eta})$ et $I_{cusp}(U_{\eta},\omega)$. L'application pr\'ec\'edente se restreint en un isomorphisme
  $$desc_{\eta}^{\tilde{G}}:I_{cusp}(\tilde{U},\omega)\to I_{cusp}(U_{\eta},\omega)^{Z_{G}(\eta;F)}.$$

  \bigskip
  
  \subsection{Filtration de $I(\tilde{G}(F),\omega)$}
  L'espace $\tilde{G}$ et le corps $F$ sont quelconques. Pour un entier $n\geq-1$, notons ${\cal F}^nI(\tilde{G}(F),\omega)$ l'espace des $f\in I(\tilde{G}(F),\omega)$ tels que $f_{\tilde{M},\omega}=0$ pour tout espace de Levi $\tilde{M}$ tel que $a_{\tilde{M}}>n$.  C'est aussi l'espace des $f\in I(\tilde{G}(F),\omega)$ qui v\'erifient la condition
  
  (1) pour tout $\gamma\in \tilde{G}_{reg}(F)$ tel que $dim(A_{G_{\gamma}})>n$, on a $I^{\tilde{G}}(\gamma,\omega,f)=0$.
   
  Ces espaces forment une filtration
  $$\{0\}={\cal F}^{a_{\tilde{G}}-1}I(\tilde{G}(F),\omega)\subset I_{cusp}(\tilde{G}(F),\omega)={\cal F}^{a_{\tilde{G}}}(\tilde{G}(F),\omega)\subset {\cal F}^{a_{\tilde{G}}+1}(\tilde{G}(F),\omega)\subset... $$
  $$\subset I(\tilde{G}(F),\omega)= {\cal F}^{a_{\tilde{M}_{0}}}(\tilde{G}(F),\omega),$$
  o\`u $\tilde{M}_{0}$ est un espace de Levi minimal. On note $GrI(\tilde{G}(F),\omega)$ l'espace gradu\'e associ\'e \`a cette filtration. Fixons un ensemble de repr\'esentants $\underline{{\cal L}}$ des classes de conjugaison par $G(F)$ d'espaces de Levi de $\tilde{G}$. Notons $\underline{{\cal L}}^n$ le sous-ensemble des $\tilde{M}\in \underline{{\cal L}}$ tels que $a_{\tilde{M}}=n$. L'application 
  $$\begin{array}{ccc}{\cal F}^nI(\tilde{G}(F),\omega)&\to&\oplus_{\tilde{M}\in \underline{{\cal L}}^n}I(\tilde{M}(F),\omega)^{W(\tilde{M})}\\ f&\mapsto& (f_{\tilde{M},\omega})_{\tilde{M}\in \underline{{\cal L}}^n}\\ \end{array}$$
  se quotiente en un homomorphisme injectif
  $$Gr^nI(\tilde{G}(F),\omega)={\cal F}^nI(\tilde{G}(F),\omega)/{\cal F}^{n-1}I(\tilde{G}(F),\omega)\to \oplus_{\tilde{M}\in \underline{{\cal L}}^n}I_{cusp}(\tilde{M}(F),\omega)^{W(\tilde{M})}.$$
  
  \ass{Lemme}{Cet homomorphisme est bijectif.}
  
  Preuve. Dans le cas o\`u $F$ est r\'eel, l'assertion est prouv\'ee par Bouaziz ([Boua], th\'eor\`eme 3.3.1)  dans le cadre non tordu et par Renard ([R1] th\'eor\`eme 11.2) dans le cadre tordu mais pour $\omega=1$. La preuve de Renard s'\'etend au cas $\omega$ quelconque. En effet, un argument de descente nous ram\`ene \`a une question analogue pour l'alg\`ebre de Lie.  Introduisons le groupe $G_{\natural}=Z(G)^0\times G_{SC}$ et l'espace $I(\mathfrak{g}_{\natural}(F))$ des int\'egrales orbitales relatives \`a ce groupe et \`a son caract\`ere trivial. Il y a un homomorphisme $\pi_{\natural}:G_{\natural}(F)\to G(F)$ de conoyau fini et $\omega$ se factorise par ce conoyau. D'autre part, $G_{\natural}$ et $G$ ont m\^eme alg\`ebre de Lie. Le conoyau $G(F)/\pi_{\natural}(G_{\natural}(F))$ agit naturellement sur $I(\mathfrak{g}_{\natural}(F))$. Alors notre espace $I(\mathfrak{g}(F),\omega)$ d'int\'egrales orbitales tordues par $\omega$ s'identifie au sous-espace de $I(\mathfrak{g}_{\natural}(F))$ o\`u ce conoyau agit par le caract\`ere $\omega$.  Passer \`a un tel sous-espace est une op\'eration \`a peu pr\`es triviale et tous les r\'esultats voulus pour $I(\mathfrak{g}(F),\omega)$ se d\'eduisent ainsi de ceux concernant $I(\mathfrak{g}_{\natural}(F))$.
  
  Le cas o\`u $F={\mathbb C}$ se ram\`ene au cas $F={\mathbb R}$ en rempla\c{c}ant les groupes et les espaces par leurs images par restriction des scalaires de ${\mathbb C}$ \`a ${\mathbb R}$.
  
  On suppose maintenant $F$ non-archim\'edien. On doit prouver la surjectivit\'e de l'homomorphisme.  On va d'abord prouver un analogue partiel pour les alg\`ebres de Lie. Supposons pour un moment que $\tilde{G}=G$. On a de m\^eme une filtration sur $I(\mathfrak{g}(F),\omega)$ et un homomorphisme injectif
   $$Gr^nI(\mathfrak{g}(F),\omega)={\cal F}^nI(\mathfrak{g}(F),\omega)/{\cal F}^{n-1}I(\mathfrak{g}(F),\omega)\to \oplus_{M\in \underline{{\cal L}}^n}I_{cusp}(\mathfrak{m}(F),\omega)^{W(\tilde{M})}.$$
   Montrons que:
   
   (2) pour tout \'el\'ement $(f^{\mathfrak{m}})_{M\in \underline{{\cal L}}^n}\in\oplus_{M\in \underline{{\cal L}}^n}I_{cusp}(\mathfrak{m}(F),\omega)^{W(M)}$, il existe $\varphi\in {\cal F}^nI(\mathfrak{g}(F),\omega)$ tel que, pour tout $M\in \underline{{\cal L}}^n$, les int\'egrales orbitales de $\varphi$ et de $f^{\mathfrak{m}}$ co\"{\i}ncident dans un voisinage de $0$ dans $\mathfrak{m}(F)$.
   
   On peut fixer $M\in \underline{{\cal L}}^n$ et supposer $f^{\mathfrak{m}'}=0$ pour tout $M'\in \underline{{\cal L}}^n$ diff\'erent de $M$. En fixant un bicaract\`ere invariant par conjugaison de $\mathfrak{g}(F)$, on introduit une  transformation de Fourier dans $C_{c}^{\infty}(\mathfrak{g}(F))$, cf. 4.1. On a de m\^eme des transformations de Fourier dans $C_{c}^{\infty}(\mathfrak{l}(F))$ pour tout Levi $L$ de $G$. D'apr\`es 4.1(2), on peut fixer $f'\in C_{c}^{\infty}(\mathfrak{m}(F))$ telle que 
   
   - son support est form\'e d'\'el\'ements elliptiques dans $\mathfrak{m}(F)$ et r\'eguliers dans $\mathfrak{g}(F)$;
   
   - les int\'egrales orbitales de $f^{\mathfrak{m}}$ et de $\hat{f}'$ co\"{\i}ncident dans un voisinage de $0$. 
   
   En rempla\c{c}ant $f'$ par la moyenne de ses conjugu\'es par un ensemble de repr\'esentants de $W(M)$, on peut supposer l'image de $f'$ dans $I(\mathfrak{m}(F),\omega)$ invariante par $W(M)$. Parce que le support de $f'$ est form\'e d'\'el\'ements r\'eguliers, on n'a aucun mal \`a trouver une fonction $\varphi'\in C_{c}^{\infty}(\mathfrak{g}(F))$  telle que
   
   - $\varphi'_{M,\omega}=f'$ dans $I(\mathfrak{m}(F),\omega)$;
   
   - le support de $\varphi'$ est un  voisinage assez petit dans $\mathfrak{g}(F)$ de celui de $f'$.
   
   Cette deuxi\`eme condition implique que le support de $\varphi'$ est form\'e d'\'el\'ements r\'eguliers dans $\mathfrak{g}(F)$ et conjugu\'es par $G(F)$ \`a des \'el\'ements elliptiques de $\mathfrak{m}(F)$.  Si $M'$ est un Levi de $G$, un tel \'el\'ement ne peut appartenir \`a $\mathfrak{m}'(F)$ que si $M'$ contient un conjugu\'e de $M$. A fortiori $\varphi'_{M',\omega}=0$ si $M'$ ne v\'erifie pas cette condition
   
   Posons $\varphi=\hat{\varphi}'$. On a $\varphi_{M,\omega}=\hat{f}'$, donc les int\'egrales orbitales de $\varphi$ et de $f^{\mathfrak{m}}$ co\"{\i}ncident dans un voisinage de $0$ dans $\mathfrak{m}(F)$. Soit $M'$ un Levi de $G$  qui v\'erifie soit $a_{M'}>n$, soit $a_{M'}=n$ et $M'$ n'est pas conjugu\'e \`a $M$. Alors $\varphi_{M',\omega}=(\varphi'_{M',\omega})\hat{}=0$. Cela entra\^{\i}ne que $\varphi\in {\cal F}^nI(\mathfrak{g}(F),\omega)$ et que $\varphi_{M',\omega}=0$ pour tout $M'\in \underline{{\cal L}}$ diff\'erent de $M$. Alors $\varphi$ satisfait les conditions de (2).
   
   Supposons de plus qu'un groupe  $\Xi$ agit sur $G$ par automorphismes d\'efinis sur $F$ en conservant le caract\`ere $\omega$. Supposons que l'image de $\Xi$ dans le groupe d'automorphismes ext\'erieurs de $G$ est fini.  Supposons les transformations de Fourier \'equivariantes pour cette action. L'action du  groupe $\Xi$ conserve la filtration $({\cal F}^nI(\mathfrak{g}(F),\omega))_{n\in {\mathbb N}}$. Il agit naturellement sur l'espace $\oplus_{M\in \underline{{\cal L}}^n}I_{cusp}(\mathfrak{m}(F),\omega)^{W(\tilde{M})}$ (un \'el\'ement $\xi\in \Xi$ envoie le terme index\'e par $M$ sur celui index\'e par l'unique \'el\'ement de $\underline{{\cal L}}^n$ conjugu\'e \`a $\xi(M)$). En prenant les invariants par $\Xi$, on obtient un homomorphisme
    $$Gr^nI(\mathfrak{g}(F),\omega)^{\Xi}={\cal F}^nI(\mathfrak{g}(F),\omega)^{\Xi}/{\cal F}^{n-1}I(\mathfrak{g}(F),\omega)^{\Xi}\to( \oplus_{M\in \underline{{\cal L}}^n}I_{cusp}(\mathfrak{m}(F),\omega)^{W(\tilde{M})})^{\Xi}.$$
   On peut aussi bien remplacer ici $\Xi$ par son image finie dans le groupe des automorphismes de $G$ quotient\'e par celui des automorphismes int\'erieurs d\'efinis par des \'el\'ements de $G(F)$. En moyennant sur ce groupe fini, on obtient pour cet homomorphisme une assertion analogue \`a (2).

 Revenons \`a l'assertion du lemme. Un argument familier de partition de l'unit\'e nous ram\`ene \`a prouver l'assertion suivante:
  
  (3) soient $\tilde{M}\in \underline{{\cal L}}^n$, $f\in I_{cusp}(\tilde{M}(F),\omega)^{W(\tilde{M})}$ et $\eta\in \tilde{M}_{ss}(F)$; alors il existe $\varphi\in {\cal F}^nI(\tilde{G}(F),\omega)$ tel que
  
  - $\varphi_{\tilde{M}',\omega}=0$ pour tout $\tilde{M}'\in \underline{{\cal L}}^n$ diff\'erent de $\tilde{M}$;
  
  - les int\'egrales orbitales de $f$ et $\varphi$ co\^{\i}ncident dans un voisinage de $\eta$ dans $\tilde{M}(F)$.
  
  Fixons donc $\tilde{M}\in \underline{{\cal L}}^n$,  $f\in I_{cusp}(\tilde{M}(F),\omega)^{W(\tilde{M})}$ et $\eta\in \tilde{M}_{ss}(F)$. Si $\eta $ n'est pas elliptique dans $\tilde{M}$, les int\'egrales orbitales de $f$ sont nulles au voisinage de $\eta$  par cuspidalit\'e de $f$ et la fonction $\varphi=0$ r\'esout la question. On suppose maintenant $\eta$ elliptique dans $\tilde{M}(F)$. Fixons un voisinage $\mathfrak{u}$ de $0$ dans $\mathfrak{g}_{\eta}(F)$, ouvert et ferm\'e et v\'erifiant les conditions de 4.1. Posons $\mathfrak{u}_{M}=\mathfrak{u}\cap \mathfrak{m}_{\eta}(F)$. On d\'eduit de $\mathfrak{u}$ et $\mathfrak{u}_{M}$ des ouverts $\tilde{U}\subset \tilde{G}(F)$ et $\tilde{U}_{M}\subset \tilde{M}(F)$. Posons ${\cal F}^nI(\tilde{U},\omega)=I(\tilde{U},\omega)\cap {\cal F}^nI(\tilde{G}(F),\omega)$, ${\cal F}^nI(\mathfrak{u},F)=I(\mathfrak{u},\omega)\cap {\cal F}^nI(\mathfrak{g}_{\eta}(F),\omega)$. La descente nous fournit un isomorphisme
  $I(\tilde{U},\omega)\simeq I(\mathfrak{u},\omega)^{Z_{G}(\eta;F)}$. Celui-ci se restreint en un isomorphisme
   $$(4) \qquad {\cal F}^n(\tilde{U},\omega)\simeq {\cal F}^n(\mathfrak{u},\omega)^{Z_{G}(\eta;F)}.$$
   C'est clair en utilisant la caract\'erisation (1) des filtrations. Par ailleurs, la descente nous fournit un isomorphisme
   $$I_{cusp}(\tilde{U}_{M},\omega)\simeq I_{cusp}(\mathfrak{u}_{M},\omega)^{Z_{M}(\eta;F)}.$$
   Notons $f_{loc}$ l'image par cet isomorphisme de la restriction de $f$ \`a $\tilde{U}_{M}$. Soit $Norm(\tilde{M},\eta;F)$ l'intersection de $Z_{G}(\eta;F)$ avec le normalisateur de $\tilde{M}$ dans $G$. Ce groupe est \'egal au normalisateur de $M_{\eta}$ dans $Z_{G}(\eta;F)$: un \'el\'ement de $Z_{G}(\eta;F)$ normalise $\tilde{M}$ ou $M_{\eta}$ si et seulement s'il normalise $A_{\tilde{M}}=A_{M_{\eta}}$.  Parce que $f$ est invariante par $W(\tilde{M})$,  $f_{loc}$ est invariante par $Norm(\tilde{M},\eta;F)$. Notons $\underline{{\cal L}}^n_{\eta}$ l'analogue de $\underline{{\cal L}}^n$ pour le groupe $G_{\eta}$. Pour $R\in \underline{{\cal L}}^n_{\eta}$, on d\'efinit un \'el\'ement $f^{\mathfrak{r}}\in I_{cusp}(\mathfrak{l}(F),\omega)$ de la fa\c{c}on suivante. Si $R$ n'est pas conjugu\'e \`a $M_{\eta}$ par un \'el\'ement de $Z_{G}(\eta;F)$, on pose $f^{\mathfrak{r}}=0$.  Si $R$ est conjugu\'e \`a $M_{\eta}$ par un \'el\'ement de $Z_{G}(\eta;F)$, on fixe un tel \'el\'ement $x$. L'automorphisme $ad_{x}$ d\'efinit un isomorphisme de $I_{cusp}(\mathfrak{m}_{\eta}(F),\omega)$ sur $I_{cusp}(\mathfrak{r}(F),\omega)$ et $f^{\mathfrak{l}^{\eta}}$ est l'image de $f_{loc}$ par cet isomormophisme. La propri\'et\'e d'invariance ci-dessus montre que cette d\'efinition ne d\'epend pas du choix de $x$. La famille $(f^{\mathfrak{r}})_{R\in \underline{{\cal L}}^n_{\eta}}$ appartient \`a $\oplus _{R\in \underline{{\cal L}}^n_{\eta}}I_{cusp}(\mathfrak{r}(F),\omega)^{W^{G_{\eta}}(R)}$ et, par construction, elle est invariante par l'action de $Z_{G}(\eta;F)$. En appliquant l'assertion (2) renforc\'ee comme on l'a expliqu\'e ci-dessus, on choisit un \'el\'ement $\varphi_{loc}\in {\cal F}^nI(\mathfrak{g}(F),\omega)^{Z_{G}(\eta;F)}$ satisfaisant la conclusion de (2).    En utilisant (4), on rel\`eve $\varphi_{loc}$ en un \'el\'ement $\varphi'$ de ${\cal F}^n(\tilde{U},\omega)$. Consid\'erons un voisinage $\mathfrak{u}'$ de $0$ dans $\mathfrak{g}_{\eta}(F)$ v\'erifiant les m\^emes conditions que $\mathfrak{u}$. On en d\'eduit un voisinage $\tilde{U}'$ de $\eta$ dans $\tilde{G}(F)$. Notons $\varphi$ le produit de $\varphi'$ et de la fonction caract\'eristique de $\tilde{U}'$. On va montrer que, si $\mathfrak{u}'$ est assez petit, $\varphi$ v\'erifie (3).  Cette fonction appartient \`a ${\cal F}^n(\tilde{U},\omega)$, cet espace \'etant \'evidemment stable par multiplication par la fonction caract\'eristique d'un ensemble ouvert et ferm\'e et invariant par conjugaison par $G(F)$. Pour $X\in \mathfrak{m}_{\eta}(F)$ assez proche de $0$, on a
   $$I^{\tilde{G}}(exp(X)\eta,\omega,\varphi)=I^{\tilde{G}}(exp(X)\eta,\omega,\varphi')=I^{G_{\eta}}(X,\omega,\varphi_{loc})$$
   $$=I^{M_{\eta}}(X,\omega,f_{loc})=I^{\tilde{M}}(exp(X)\eta,\omega,f),$$
   ce qui est la derni\`ere condition requise. Soient $\tilde{M}'\in \underline{{\cal L}}^n$ diff\'erent de $\tilde{M}$ et $\gamma$ un \'el\'ement $\tilde{G}$-r\'egulier de $\tilde{M}'(F)$. On doit montrer que $I^{\tilde{G}}(\gamma,\omega,\varphi)=0$. C'est clair si $\gamma\not\in \tilde{U}'$. Supposons $\gamma\in \tilde{U}'$. On peut alors \'ecrire $\gamma=g^{-1}exp(X)\eta g$, avec $g\in G(F)$ et $X\in \mathfrak{u}'$. Quitte \`a changer $g$, on peut conjuguer $X$ par un \'el\'ement de $G_{\eta}(F)$ et supposer $X$ assez proche de $0$. Posons $\tilde{M}''=g\tilde{M}'g^{-1}$. Puisque $\gamma\in \tilde{M}'(F)$, on a $exp(X)\eta\in \tilde{M}''(F)$. Donc $A_{\tilde{M}''}\subset Z_{G}(exp(X)\eta)$. Pour $X$ assez petit, ce commutant est inclus dans $Z_{G}(\eta)$. Alors $\eta\in \tilde{M}''(F)$, puis $X\in \mathfrak{m}''_{\eta}(F)$.  On a comme ci-dessus
$$\omega(g) I^{\tilde{G}}( g^{-1} exp(X)\eta g,\omega,\varphi)=I^{\tilde{G}}(exp(X)\eta,\omega,\varphi')=I^{G_{\eta}}(X,\omega,\varphi_{loc})=I^{M''_{\eta}}(X,\omega,\varphi_{loc,M''_{\eta},\omega}).$$
 On a $A_{\tilde{M}''}\subset A_{M''_{\eta}}$. Si cette inclusion est stricte, $dim(A_{M''_{\eta}})>n$ et les int\'egrales orbitales ci-dessus sont nulles puisque $\varphi_{loc}\in {\cal F}^nI(\mathfrak{g}_{\eta}(F),\omega)$. Si l'inclusion ci-dessus est une \'egalit\'e, $M''_{\eta}$ est conjugu\'e par $G_{\eta}(F)$ \`a un \'el\'ement de $\underline{{\cal L}}^n_{\eta}$ et il r\'esulte de notre construction  que les int\'egrales ci-dessus sont encore nulles sauf si $M''_{\eta}$ est conjugu\'e \`a $M_{\eta}$ par un \'el\'ement de $Z_{G}(\eta;F)$. Il reste \`a exclure cette possibilit\'e. Mais, parce que l'on a \`a la fois $A_{\tilde{M}''}=A_{M''_{\eta}}$ et $A_{\tilde{M}}=A_{M_{\eta}}$, dire que $M''_{\eta}$ et $M_{\eta}$ sont conjugu\'es par un \'el\'ement de $Z_{G}(\eta;F)$ revient \`a dire que $\tilde{M}''$ et $\tilde{M}$ le sont. Puisque $\tilde{M}''$ et $\tilde{M}'$ sont conjugu\'es par $g$, cela est exclu par notre hypoth\`ese que $\tilde{M}'$ n'est conjugu\'e \`a $\tilde{M}$ par   aucun \'el\'ement de $G(F)$. $\square$

\bigskip

\subsection{Image de la restriction}

Pour un espace de Levi $\tilde{M}$ de $\tilde{G}$, on note  $res_{\tilde{M}}$ l'homomorphisme
$$\begin{array}{ccc}I(\tilde{G}(F),\omega)&\to&I(\tilde{M}(F),\omega)\\ f&\mapsto&f_{\tilde{M},\omega},\\ \end{array}$$
ou sa variante envoyant $I(\tilde{G}(F),\omega)\otimes Mes(G(F))$ dans $I(\tilde{M}(F),\omega)\otimes Mes(M(F))$. Soit $(\tilde{M}_{j})_{j=1,...,k}$ une famille finie d'espaces de Levi de $\tilde{G}$. Consid\'erons l'application lin\'eaire
$$res=\oplus_{j=1,...,k}res_{\tilde{M}_{j}}:I(\tilde{G}(F),\omega)\otimes Mes(G(F))\to \oplus_{j=1,...,k}I(\tilde{M}_{j}(F),\omega)\otimes Mes(M_{j}(F)).$$

\ass{Lemme}{L'image de $res$ est l'espace des $(\boldsymbol{\varphi}_{j})_{j=1,...,k}\in\oplus_{j=1,...,k} I(\tilde{M}_{j}(F),\omega)\otimes Mes(M_{j}(F))$ qui v\'erifient les conditions \'equivalentes suivantes:

(i) soient $j,j'\in \{1,...,k\}$,  $\gamma \in \tilde{M}_{j}(F)$ et $\gamma'\in \tilde{M}_{j'}(F)$ deux \'el\'ements $\tilde{G}$-r\'eguliers et soit $g\in G(F)$ tel que $\gamma'=g\gamma g^{-1}$;  munissons $G_{\gamma}(F)$ et $G_{\gamma'}(F)$ de mesures de Haar se correspondant par $ad_{g}$; alors $I^{\tilde{M}_{j'}}(\gamma',\omega,\boldsymbol{\varphi}_{j'})=\omega(g)I^{\tilde{M}_{j}}(\gamma,\omega,\boldsymbol{\varphi}_{j})$;

(ii) soient $j,j'\in \{1,...,k\}$, $\tilde{R}$ un espace de Levi de $\tilde{M}_{j}$ et $\tilde{R}'$ un espace de Levi de $\tilde{M}_{j'}$ et soit $g\in G(F)$ tel que $\tilde{R}'=ad_{g}(\tilde{R})$; alors $\boldsymbol{\varphi}_{\tilde{R}',\omega}$ est l'image de $\boldsymbol{\varphi}_{\tilde{R},\omega}$ par l'isomorphisme $I(\tilde{R},\omega)\otimes Mes(R(F))\to I(\tilde{R}',\omega)\otimes Mes(R'(F))$ d\'eduit de $ad_{g}$.}

{\bf Remarque.}  Dans (i), la donn\'ee de $\gamma$ et d'une mesure de Haar sur $G_{\gamma}(F)$ d\'efinit une int\'egrale orbitale qui est naturellement une forme lin\'eaire sur $ I(\tilde{M}_{j}(F),\omega)\otimes Mes(M_{j}(F))$. 

\bigskip

Preuve. Pour simplifier les notations, on oublie les espaces de mesures. Il est clair que les deux conditions de l'\'enonc\'e sont \'equivalentes et qu'elles sont v\'erfi\'ees sur les \'el\'ements de l'image de $res $. Posons
$$I=\oplus_{j=1,...,k} I(\tilde{M}_{j}(F),\omega)$$
et, pour tout $n$,
$${\cal F}^nI=\oplus_{j=1,...,k} {\cal F}^nI(\tilde{M}_{j}(F),\omega).$$
Notons $J$ le sous-espace des $(\varphi_{j})_{j=1,...,k}\in I$ satisfaisant les conditions (i) ou (ii). Il est clair que $res $ envoie ${\cal F}^nI(\tilde{G}(F),\omega)$ dans ${\cal F}^nI $, donc aussi dans $J\cap {\cal F}^nI $.  Donc $res $ d\'efinit une application
$$(1) \qquad Gr^nI(\tilde{G}(F),\omega)\to (J\cap {\cal F}^nI )/(J\cap {\cal F}^{n-1}I) .$$
On va montrer qu'elle est surjective. L'espace de d\'epart est isomorphe \`a 
$$(2) \qquad \oplus_{\tilde{L}\in \underline{{\cal L}}^n}I_{cusp}(\tilde{L}(F),\omega)^{W(\tilde{L})}$$
tandis que l'espace d'arriv\'ee est inclus dans
$$(3)\qquad  Gr^nI \simeq \oplus_{j=1,...,k} \oplus_{\tilde{R}\in \underline{{\cal L}}^{\tilde{M}_{j},n}}I_{cusp}(\tilde{R}(F),\omega)^{W^{M_{j}}(\tilde{R})}.$$
L'image dans l'espace (2) de $(J\cap {\cal F}^nI )/(J\cap {\cal F}^{n-1}I)$ est contenu dans le sous-espace des \'el\'ements v\'erifiant la condition (ii) restreinte aux espaces de Levi $\tilde{R} \in \underline{{\cal L}}^{\tilde{M}_{j},n}$ et $\tilde{R}' \in \underline{{\cal L}}^{\tilde{M}_{j'},n}$. Pour un \'el\'ement $(\varphi_{j}^{\tilde{R}})_{j=1,...k,\tilde{R}\in \underline{{\cal L}}^{\tilde{M}_{j},n}}$ v\'erifiant cette condition, on d\'efinit un \'el\'ement $(f^{\tilde{L}})_{\tilde{L}\in \underline{{\cal L}}^n}$ de l'espace (2) de la fa\c{c}on suivante. Soit $\tilde{L}\in \underline{{\cal L}}^n$. S'il n'existe pas de $j\in \{1,...,k\}$ et de $\tilde{R}\in \underline{{\cal L}}^{\tilde{M}_{j},n}$ tel suq $\tilde{L}$ soit conjugu\'e \`a $\tilde{R}$ par un \'el\'ement de $G(F)$, on pose
  $f^{\tilde{L}}=0$. Si  au contraire il existe un tel couple $(j,\tilde{R})$, on en fixe un et   on choisit   un \'el\'ement $g$ tel que $ad_{g}(\tilde{R})=\tilde{L}$. Alors $f^{\tilde{L}}$ est l'image de $\varphi^{\tilde{R}}$ par l'isomorphisme d\'eduit de $ad_{g}$. La condition (ii) entra\^{\i}ne que cela ne d\'epend pas des choix et que la fonction $f^{\tilde{L}}$ est bien invariante par $W(\tilde{L})$. Il est clair que  $(\varphi_{j}^{\tilde{R}})_{j=1,...,k,\tilde{R}\in \underline{{\cal L}}^{\tilde{M}_{j},n}}$ est l'image de $(f^{\tilde{L}})_{\tilde{L}\in \underline{{\cal L}}^n}$ par la compos\'ee de l'application (1) et de l'inclusion de son espace d'arriv\'ee dans l'espace (3). Cela d\'emontre la surjectivit\'e de l'application (1)

Par r\'ecurrence sur $n$, on en d\'eduit que l'application
$$res:{\cal F}^nI(\tilde{G}(F),\omega)\to J\cap {\cal F}^nI$$
est surjective. Pour $n$ grand, cela signifie que $J$ est bien l'image de l'application $res$. $\square$ 
         
\bigskip

\subsection{Conjugaison stable}
On a d\'ej\`a rappel\'e la notion de conjugaison stable pour les \'el\'ements de $\tilde{G}_{reg}(F)$: deux \'el\'ements de cet ensemble sont stablement conjugu\'es si et seulement s'ils sont conjugu\'es par un \'el\'ement de $G=G(\bar{F})$. Pour un \'el\'ement $\eta\in \tilde{G}_{ss}(F)$, on note $I_{\eta}=G_{\eta}Z(G)^{\theta}$ et on pose
$${\cal Y}(\eta)=\{y\in G; \forall \sigma\in \Gamma_{F}, y\sigma(y)^{-1}\in I_{\eta}\}.$$
Pour deux \'el\'ements $\eta,\eta'\in \tilde{G}_{ss}(F)$, on appelle diagramme joignant $\eta$ et $\eta'$ un sextuplet $(\eta,B,T,B',T',\eta')$ tel que

(1) $(B,T)$ et $(B',T')$  sont des paires de Borel de $G$;

(2) $ad_{\eta}$ conserve $(B,T)$ et $ad_{\eta'}$ conserve $(B',T')$;

(3) $T$ et $T'$ sont d\'efinis sur $F$ et l'isomorphisme $\xi_{T,T'}:T\to T'$ issu des deux paires est \'equivariant pour les actions galoisiennes;

compl\'etons les deux paires en des paires de Borel \'epingl\'ees ${\cal E}$ et ${\cal E}'$, \'ecrivons
$\eta=te$, avec $t\in T$ et $e\in Z(\tilde{G},{\cal E})$ et \'ecrivons de m\^eme $\eta'=t'e'$; on impose que $e$ et $e'$ aient m\^eme image dans ${\cal Z}(\tilde{G})$; alors

(4) $\xi_{T,T'}(t)\in t'(1-\theta')(T')$, o\`u $\theta'$ est l'automorphisme de $T'$ d\'etermin\'e par ${\cal E}'$.

On voit que la condition (4) ne d\'epend pas des choix auxiliaires.

Dans le cas o\`u $\eta$ et $\eta'$ sont fortement r\'eguliers, on montre comme au lemme 1.10(i) qu'il existe un diagramme joignant $\eta$ et $\eta'$ si et seulement si ces deux \'el\'ements sont stablement conjugu\'es. En g\'en\'eral, consid\'erons les conditions suivantes:

(st1) il existe $y\in {\cal Y}(\eta)$ tel que $\eta'=y^{-1}\eta y$;

(st2) il existe un diagramme   $(\eta,B,T,B',T',\eta')$;

(st3)   il existe un diagramme   $(\eta,B,T,B',T',\eta')$ tel que

(st3)(a) si $F$ est non archim\'edien, $T^{\theta,0}$ est elliptique dans $G_{\eta}$ (c'est-\`a-dire $T^{\theta,0}/Z(G_{\eta})$ ne contient pas de sous-tore d\'eploy\'e non trivial) et $(T')^{\theta',0}$ est elliptique dans $G_{\eta'}$;

(st3)(b) si $F$ est r\'eel, $T^{\theta,0}$ est fondamental dans $G_{\eta}$ et $(T')^{\theta',0}$ est fondamental dans $G_{\eta'}$;

(st4) $(\eta,\eta')$ appartient \`a l'adh\'erence dans $\tilde{G}(F)\times \tilde{G}(F)$ de l'ensemble des couples $(\gamma,\gamma')\in \tilde{G}_{reg}(F)\times \tilde{G}_{reg}(F)$ tels que $\gamma$ et $\gamma'$ sont stablement conjugu\'es.

\ass{Lemme }{Les conditions (st1) \`a (st4) ci-dessus sont \'equivalentes.}

Preuve. La m\^eme preuve qu'au lemme 1.10(ii) montre l'\'equivalence de (st2) et (st4). 

Supposons (st2) v\'erifi\'ee et fixons un diagramme $(\eta,B,T,B',T',\eta')$. On compl\`ete les paires de Borel en des paires \'epingl\'ees et on \'ecrit $\eta$ et $\eta'$ comme en (4). Soit $x\in G$ tel que $ad_{x}$ envoie ${\cal E}$ sur ${\cal E}'$. Les \'el\'ements $e$ et $ad_{x}(e)$ ont par d\'efinition m\^eme image dans ${\cal Z}(\tilde{G})$. L'hypoth\`ese de (4) est que $e'$ et $e$ ont m\^eme image dans ${\cal Z}(\tilde{G})$. Cela signifie que, quitte \`a multiplier $x$ par un \'el\'ement de $Z(G)$, on peut supposer $ad_{x}(e)=e'$. L'isomorphisme $\xi_{T,T'}$ n'est autre que la restriction \`a $T$ de $ad_{x}$. D'apr\`es (4), on peut donc \'ecrire $ad_{x}(t)=t'(1-\theta')(t'')$, avec un $t''\in T'$. Alors $x\eta x^{-1}=t''\eta'(t'')^{-1}$. Posons $y=x^{-1}t''$. On a $y^{-1}\eta y=\eta'$. L'isomorphisme $\xi_{T,T'}$ est encore la restriction de $ad_{y}$. Puisqu'il est d\'efini sur $F$, $y\sigma(y)^{-1}$ commute \`a $T$, donc appartient \`a $T$, pour tout $\sigma\in \Gamma_{F}$. L'\'egalit\'e $y^{-1}\eta y=\eta'$ et le fait que $\eta$ et $\eta'$ appartiennent \`a $\tilde{G}(F)$ entra\^{\i}nent que $y\sigma(y)^{-1}$ appartient aussi \`a $Z_{G}(\eta)$. Or  $T\cap Z_{G}(\eta)\subset I_{\eta}$ ([W1] 3.1(1)). Donc $y\in {\cal Y}(\eta)$ et (st1) est v\'erifi\'ee. 

Supposons (st1) v\'erifi\'ee, fixons $y\in {\cal Y}(\eta)$ tel que $y^{-1}\eta y=\eta'$. Fixons, ainsi qu'il est loisible, une paire de Borel $(B,T)$ conserv\'ee par $ad_{\eta}$, telle que $T$ soit d\'efini sur $F$ et $T^{\theta,0}$ soit elliptique dans $G_{\eta}$ si $F$ est non archim\'edien, ou fondamental si $F={\mathbb R}$. L'automorphisme $ad_{y^{-1}}$ envoie $G_{\eta}$ sur $G_{\eta'}$ et l'hypoth\`ese que $y$ appartient \`a ${\cal Y}(\eta)$ entra\^{\i}ne que sa restriction \`a $G_{\eta}$ est un torseur int\'erieur entre ces deux groupes. On sait qu'un tore elliptique, ou fondamental, se transf\`ere \`a toute forme int\'erieure (et son transfert est encore elliptique ou fondamental). Quitte \`a multiplier $y$ \`a droite par un \'el\'ement de $G_{\eta'}$, on peut donc supposer que $ad_{y^{-1}}(T^{\theta,0})$ est d\'efini sur $F$ et que la restriction de $ad_{y^{-1}}:T^{\theta,0}\to ad_{y^{-1}}(T^{\theta,0})$ est d\'efini sur $F$. Posons $B'=ad_{y^{-1}}(B)$, $T'=ad_{y^{-1}}(T)$. Puisque $T$ est le commutant de $T^{\theta,0}$, les propri\'et\'es pr\'ec\'edentes impliquent que $T$ est d\'efini sur $F$ et que $ad_{y^{-1}}:T\to T'$ l'est aussi. Evidemment, $ad_{\eta'}$ conserve $(B',T')$. On compl\`ete nos paires en des paires \'epingl\'ees et on \'ecrit $\eta$ et $\eta'$ comme en (4).  Soit $x\in G$ qui envoie ${\cal E}'$ sur ${\cal E}$. Comme ci-dessus, on peut imposer que $ad_{x}(e')=e$. Puisque $ad_{y^{-1}}$ et $ad_{x^{-1}}$ envoient tous deux $(B,T)$ sur $(B',T')$, on peut \'ecrire $y=xt''$, avec un $t''\in T'$. L'\'egalit\'e $ad_{y^{-1}}(\eta)=\eta'$ entra\^{\i}ne que $ad_{x^{-1}}(t)=t'(1-\theta')(t'')$. Puisque $\xi_{T,T'}$ est la restriction de $ad_{x^{-1}}$ \`a $T$, on obtient (4). Donc $(\eta,B,T,B',T',\eta')$ est un diagramme v\'erifiant les conditions suppl\'ementaires de (st3).

Enfin, (st3) implique \'evidemment (st2). $\square$

\ass{D\'efinition}{On dit  que $\eta$ et $\eta'$ sont stablement conjugu\'es si et seulement si les  conditions  (st1),...,(st4) sont v\'erifi\'ees.}

\bigskip

\subsection{Conjugaison stable et application $N^{\tilde{G}}$}

\ass{Lemme }{Soient $\eta$, $\eta'$ deux \'el\'ements stablement conjugu\'es de $\tilde{G}_{ss}(F)$. Alors on a l'\'egalit\'e $N^{\tilde{G}}(\eta)=N^{\tilde{G}}(\eta')$ dans $\tilde{G}_{0,ab}(F)$.}

Preuve. On fixe une paire de Borel ${\cal E}$, on \'ecrit $\eta=\pi(x)e$, avec $x\in G_{SC}$ et $e\in Z(\tilde{G},{\cal E})$. On pose $\theta=ad_{e}$. L'\'el\'ement $N^{\tilde{G}}(\eta)$ est l'image dans $\tilde{G}_{0,ab}(F)$ du cocycle $(\bar{\nu},\bar{e})\in Z^{1,0}(\Gamma_{F};{\cal Z}(G_{SC})\circlearrowleft {\cal Z}(\tilde{G}))$, o\`u $\nu(\sigma)=\theta(u_{{\cal E}}(\sigma))x^{-1}\sigma(x)u_{{\cal E}}(\sigma)^{-1}$ (les $\bar{}$ d\'esignent les images dans ${\cal Z}(G_{SC})$ ou ${\cal Z}(\tilde{G})$). Soit $y\in {\cal Y}(\eta)$ tel que $\eta'=y^{-1}\eta y$. Ecrivons $y=z \pi(v)$, avec $z\in Z(G)$ et $v\in G_{SC}$. Alors $\eta'=\pi(x')e'$, avec $x'=v^{-1}x\theta(v)$, $e'=z^{-1}\theta(z)e$. L'\'el\'ement $N^{\tilde{G}}(\eta')$ est l'image du cocycle $(\bar{\nu}',\bar{e}')$, o\`u 
$$\nu'(\sigma)=\theta(u_{{\cal E}}(\sigma))\theta(v)^{-1}x^{-1}v\sigma(v)^{-1}\sigma(x)\sigma(\theta(v))u_{{\cal E}}(\sigma)^{-1}.$$
Introduisons l'action quasi-d\'eploy\'ee $\sigma\mapsto \sigma_{G^*}=ad_{u_{{\cal E}}(\sigma)}\circ\sigma$ qui pr\'eserve ${\cal E}$. Puisque $e\in Z(\tilde{G},{\cal E})$, $\theta=ad_{e}$ est fixe pour cette action. Donc 
$$u_{{\cal E}}(\sigma)\sigma(\theta(v))u_{{\cal E}}(\sigma)^{-1}=\theta(u_{{\cal E}}(\sigma))\theta(\sigma(v))\theta(u_{{\cal E}}(\sigma))^{-1}.$$
Puisque $\nu'$ est \`a valeurs centrales, on peut aussi bien  conjuguer $\nu'(\sigma)$ par cette expression et on obtient
$$\nu'(\sigma)=\theta(u_{{\cal E}}(\sigma)\theta(\sigma(v)v^{-1})x^{-1}v\sigma(v)^{-1}\sigma(x)u_{{\cal E}}(\sigma)^{-1}.$$
L'hypoth\`ese $y\in {\cal Y}(\eta)$ entra\^{\i}ne que $\pi(v\sigma(v)^{-1})\in Z(G)G_{\eta}$, a fortiori $v_{ad}\sigma(v_{ad})^{-1}\in G_{AD,\eta}$. Mais $G_{SC,\eta}$ s'envoie surjectivement sur $G_{AD,\eta}$. Donc $v\sigma(v)^{-1}\in Z(G_{SC})G_{SC,\eta}$. Ecrivons $v\sigma(v)^{-1}=\zeta(\sigma)g(\sigma)$, avec $\zeta(\sigma)\in Z(G_{SC})$ et $g(\sigma)\in G_{SC,\eta}$. Cette derni\`ere relation signifie que $x\theta(g(\sigma))x^{-1}=g(\sigma)$. On calcule alors
$$\nu'(\sigma)=\theta(\zeta(\sigma))^{-1}\zeta(\sigma)\nu(\sigma).$$
Donc $\bar{\nu}'=\bar{\nu}$. On  a aussi $\bar{e}'=\bar{e}$  et le lemme s'ensuit. $\square$

\bigskip

\subsection{Description locale des classes de conjugaison stable}

Pour $\eta\in \tilde{G}_{ss}(F)$, fixons un ensemble de repr\'esentants $\dot{\cal Y}(\eta)$ de l'ensemble de doubles classes $I_{\eta}\backslash{\cal Y}(\eta)/ G(F)$. L'application qui \`a $y\in \dot{\cal Y}(\eta)$ associe la classe de conjugaison par $G(F)$ de $\eta[y]=y^{-1}\eta y$ est une surjection de $\dot{\cal Y}(\eta)$ sur l'ensemble des classes de conjugaison par $G(F)$ contenues dans la classe de conjugaison stable de $\eta$. En g\'en\'eral, elle n'est pas injective. C'est toutefois le cas si $\eta$ est fortement r\'egulier.

Soit $\eta\in \tilde{G}_{ss}(F)$. Fixons une forme quasi-d\'eploy\'ee $\bar{G}$ de $G_{\eta}$. On peut, si on veut, fixer un torseur int\'erieur entre ces deux groupes. Nous pr\'ef\'erons dire que nous fixons une identification entre $\underline{la}$ paire de Borel \'epingl\'ee de $\bar{G}$ et celle de $G_{\eta}$. Pour tout $y\in {\cal Y}(\eta)$, l'automorphisme $ad_{y^{-1}}$ permet d'identifier $\underline{la}$ paire de Borel \'epingl\'ee de $G_{\eta}$ et celle de $G_{\eta[ y]}$, d'o\`u une identification de cette derni\`ere avec celle de $\bar{G}$.  Il y a donc une correspondance entre classes de conjugaison stable semi-simples dans $G_{\eta[ y]}(F)$ et classes de conjugaison stable semi-simples dans $\bar{G}(F)$. D'autre part, les groupes $Z_{G}(\eta [y])/I_{\eta [y]}$ s'identifient de fa\c{c}on \'equivariante pour les actions. On note $\Xi$ ce groupe commun. On le fait agir sur $\bar{G}$ de sorte que cette action conserve  une paire de Borel  \'epingl\'ee d\'efinie sur $F$ fix\'ee. Cette action est fid\`ele (seul l'\'el\'ement neutre de $\Xi$ agit par l'identit\'e). Les actions galoisiennes sur $\Xi$ et $\bar{G}$ sont compatibles. En particulier, $\Xi^{\Gamma_{F}}$ agit par automorphismes d\'efinis sur $F$.

Fixons un ouvert $\mathfrak{\bar{u}}$ de $\mathfrak{\bar{g}}(F)$ contenant $0$, tel que 

- $\bar{X}\in \bar{\mathfrak{u}}$ si et seulement si $\bar{X}_{ss}\in \bar{\mathfrak{u}}$, o\`u $\bar{X}_{ss}$ est la partie semi-simple de $\bar{X}$;

- si $\bar{X}\in\bar{\mathfrak{u}}$ et $\bar{X}'\in \mathfrak{\bar{g}}(F)$ sont conjugu\'es par un \'el\'ement de $\bar{G}(\bar{F})$, alors $\bar{X}'\in \bar{\mathfrak{u}}$;

- $\mathfrak{\bar{u}}$ est invariant par $\Xi^{\Gamma_{F}}$.

 Pour tout $y$, il lui correspond un tel voisinage $\mathfrak{u}_{\eta[ y]}\subset \mathfrak{g}_{\eta[ y]}(F)$, form\'e des $X$ tels que la classe de conjugaison stable de $X_{ss}$ corresponde \`a celle d'un \'el\'ement de $\bar{\mathfrak{u}}$. Soit $\bar{X}\in \bar{\mathfrak{u}}\cap \bar{\mathfrak{g}}_{reg}(F)$. Pour tout $y\in \dot{\cal Y}(\eta)$, fixons un ensemble $\dot{\cal X}(\bar{X},y)\subset \mathfrak{u}_{\eta[ y]}$ de repr\'esentants des classes de conjugaison par $I_{\eta[ y]}(F)$ dans la classe de conjugaison stable de $\mathfrak{g}_{\eta[ y]}(F)$ correspondant \`a celle de $\bar{X}$, si cette classe existe. Sinon, on pose  $\dot{\cal X}(\bar{X},y)=\emptyset$. Notons $C(\bar{X})$ la classe de conjugaison stable commune dans $\tilde{G}(F)$ des $exp(X)\eta[ y]$, pour $y\in \dot{\cal Y}(\eta)$ et $X\in \dot{\cal X}(\bar{X},y)$. Notons $\bar{\mathfrak{u}}_{\tilde{G}-reg}$ le sous-ensemble des $\bar{X}$ tels que $C(\bar{X})$ soit form\'e d'\'el\'ements fortement r\'eguliers dans $\tilde{G}$.

  Notons $\tilde{U}$ l'ensemble des \'el\'ements $\gamma\in \tilde{G}(F)$ tels que la partie semi-simple de $\gamma$ soit  stablement conjugu\'ee \`a un \'el\'ement $exp(X)\eta[ y]$ pour un $y\in \dot{\cal Y}(\eta)$ et un $X\in \mathfrak{u}_{\eta[ y]}$ (en supposant $\bar{\mathfrak{u}}$ assez petit pour que ces exponentielles soient d\'efinies). Notons $\tilde{U}'$ des \'el\'ements $\gamma\in \tilde{G}(F)$ tels que la partie semi-simple de $\gamma$ soit   conjugu\'ee par un \'el\'ement de $G(F)$ \`a un \'el\'ement $exp(X)\eta [y]$ pour un $y\in \dot{\cal Y}(\eta)$ et un $X\in \mathfrak{u}_{\eta[ y]}$. 
  
      \ass{Lemme}{Si $\bar{\mathfrak{u}}$ est assez petit, les propri\'et\'es suivantes sont v\'erifi\'ees.
      
 (i) L'ensemble $\tilde{U}$ est ouvert et \'egal \`a $\tilde{U}'$.
 
 (ii) L'application $\bar{X}\mapsto C(\bar{X})$ est une surjection  de $ \bar{\mathfrak{u}}_{\tilde{G}-reg}$   sur l'ensemble des classes de conjugaison stable contenues dans $\tilde{U}\cap \tilde{G}_{reg}(F)$.
 
 (iii) On a $C(\bar{X})=C(\bar{X}')$ si et seulement s'il existe $\xi\in \Xi^{\Gamma_{F}}$ tel que $\xi(\bar{X})$ soit stablement conjugu\'e \`a $\bar{X}'$.
 
 (iv) Pour tout $\bar{X}\in \bar{\mathfrak{u}}_{\tilde{G}-reg}$, l'ensemble $\{exp(X)\eta[ y]; y\in \dot{\cal Y}(\eta), X\in \dot{\cal X}(\bar{X},y)\}$ est un ensemble de repr\'esentants des classes de conjugaison par $G(F)$ dans $C(\bar{X})$. }

 Preuve.  On a \'evidemment $\tilde{U}'\subset \tilde{U}$. Pour d\'emontrer l'inclusion oppos\'ee, on peut se limiter aux \'el\'ements semi-simples. Soit $\gamma\in \tilde{U}$ un tel \'el\'ement. On peut fixer $y\in \dot{\cal Y}(\eta)$,  $X\in \mathfrak{u}_{\eta[ y]}$ et un diagramme $(\gamma,B,T,B',T',\gamma')$, o\`u $\gamma'=exp(X)\eta[ y]$. Posons $\theta'=ad_{\gamma'}$. Le tore $(T')^{\theta',0}$ est un sous-tore maximal de $G_{\gamma'}$. Si $\bar{\mathfrak{u}}$ est assez petit,  $G_{\gamma'}$ est le commutant de $X$ dans $G_{\eta [y]}$. Donc $X$ appartient au centre de $\mathfrak{g}_{\gamma'}$, a fortiori \`a $(\mathfrak{t}')^{\theta'}(F)$. Soit $Y$ l'image de $X$ par l'application $\xi_{T',T}:\mathfrak{t}'(F)\to \mathfrak{t}(F)$. Alors $Y$ est fixe par $ad_{\gamma}$ et on v\'erifie que $(exp(-Y)\gamma,B,T,B',T',\eta[ y])$ est un diagramme. Donc $exp(-Y)\gamma$ est stablement conjugu\'e \`a $\eta[ y]$. Il existe donc $y_{1}\in \dot{\cal Y}(\eta)$ tel que $exp(-Y)\gamma$ soit conjugu\'e \`a $\eta [y_{1}]$ par un \'el\'ement de $G(F)$. Quitte \`a effectuer une telle conjugaison, on peut supposer que ces deux \'el\'ements sont \'egaux. Alors $\gamma=exp(Y)\eta [y_{1}]$, avec $Y\in \mathfrak{g}_{\eta[ y_{1}]}(F)$ (parce que $Y$ commute \`a $\gamma$ et \`a $exp(Y)$). Il r\'esulte des d\'efinitions des voisinages et de l'hypoth\`ese  $X\in \mathfrak{u}_{ \eta[ y]}$ que $Y$ appartient \`a $\mathfrak{u}_{\eta [y_{1}]}$. Cela prouve l'\'egalit\'e $\tilde{U}=\tilde{U}'$. L'ensemble $\tilde{U}'$ \'etant clairement ouvert, cela prouve (i).

  Le (ii) est \'evident. Le (iv) est le lemme 3.8 de [W1] (dans cette r\'ef\'erence, le corps $F$ est non-archim\'edien, mais la preuve vaut aussi bien pour $F$ archim\'edien). Pour le (iii), on  peut identifier $\bar{G}$ \`a $G_{\eta}$ muni d'une action galoisienne de la forme $\sigma\mapsto \sigma_{\bar{G}}=ad_{u(\sigma)}\circ\sigma$, o\`u $u(\sigma)\in G_{\eta,SC}$. Posons $\gamma=exp(\bar{X})\eta$, $\gamma'=exp(\bar{X}')\eta$. Dire que $C(\bar{X})=C(\bar{X}')$ revient \`a dire qu'il existe $g\in G$ tel que $g\gamma g^{-1}=\gamma'$. Si $\bar{\mathfrak{u}}$ est assez petit, cela entra\^{\i}ne $g\in Z_{G}(\eta)$. Pour $\sigma\in \Gamma_{F}$, on a $\sigma(g)\sigma(\gamma)\sigma(g)^{-1}=\sigma(\gamma')$. Puisque $\bar{X}\in \bar{\mathfrak{g}}(F)$, on a $\sigma(\gamma)=u(\sigma)^{-1}\gamma u(\sigma)$. De m\^eme, $\sigma(\gamma')=u(\sigma)^{-1}\gamma'u(\sigma)$.  D'o\`u $u(\sigma)\sigma(g)u(\sigma)^{-1}\gamma u(\sigma)g^{-1}u(\sigma)^{-1}=\gamma'$. Alors $g^{-1}u(\sigma)\sigma(g)u(\sigma)^{-1}$ fixe $\gamma$, donc est contenu dans $I_{\gamma}$, lui-m\^eme contenu dans $I_{\eta}$. Donc l'image de $g$ dans $\Xi$ est fixe par $\Gamma_{F}$ et la conclusion de (iii) s'ensuit. La r\'eciproque est claire. $\square$
  
  \bigskip
  
  \subsection{Conjugaison stable et $K$-espaces tordus}
 
 Dans le cas o\`u $F={\mathbb R}$, les d\'efinitions et r\'esultats des trois paragraphes pr\'ec\'edents  s'adaptent aux $K$-espaces tordus. Il suffit de d\'efinir correctement la notion de conjugaison stable et les ensembles ${\cal Y}(\eta)$ et $\dot{\cal Y}(\eta)$. Pour des \'el\'ements $\gamma\in \tilde{G}_{p,reg}({\mathbb R})$ et $\gamma' \in \tilde{G}_{p',reg}({\mathbb R})$, on dit simplement qu'ils sont stablement conjugu\'es si $\gamma$ est conjugu\'e \`a $\tilde{\phi}_{p,p'}(\gamma')$ par un \'el\'ement de $G_{p}$. Soit $\eta\in \tilde{G}_{p,ss}({\mathbb R})$. Pour $p'\in \Pi$, on note ${\cal Y}_{p'}(\eta)$ l'ensemble des $y\in G_{p'}$ tels que $y\sigma(y)^{-1}\nabla_{p',p}(\sigma)^{-1}\in I_{\tilde{\phi}_{p',p}(\eta)}$ pour tout $\sigma\in \Gamma_{F}$.  Pour $y\in {\cal Y}_{p'}(\eta)$, on pose $\eta[ y]=y^{-1}\tilde{\phi}_{p',p}(\eta)y$. On note $\dot{\cal Y}_{p'}(\eta)$ un ensemble de repr\'esentants des doubles classes $ I_{\tilde{\phi}_{p',p}(\eta)}\backslash {\cal Y}_{p'}(\eta)/G_{p'}(F) $. On pose ${\cal Y}(\eta)=\sqcup_{p'\in \Pi}{\cal Y}_{p'}(\eta)$, $\dot{\cal Y}(\eta)=\sqcup_{p'\in \Pi}\dot{\cal Y}_{p'}(\eta)$.  Remarquons que, puisque les paires de Borel des diff\'erents groupes $G_{p}$ s'identifient, on peut d\'efinir sans changement la notion de diagramme joignant deux \'el\'ements semi-simples de $K\tilde{G}({\mathbb R})$. Avec les d\'efinitions ci-dessus, les propri\'et\'es (st1) \`a (st4) de 4.4 restent \'equivalentes pour $\eta,\eta'\in K\tilde{G}_{ss}({\mathbb R})$. On dit que $\eta$ et $\eta'$ sont stablement conjugu\'es si et seulement si ces conditions sont v\'erifi\'ees.
 
 \bigskip
 
 \subsection{Descente d'Harish-Chandra et stabilit\'e}
 
  Supposons   $(G,\tilde{G},{\bf a})$ quasi-d\'eploy\'e et \`a torsion int\'erieure. On sait que tout \'el\'ement semi-simple de $\tilde{G}(F)$ est stablement conjugu\'e \`a un \'el\'ement $\epsilon$ pour lequel $G_{\epsilon}$ est quasi-d\'eploy\'e. Soit $\epsilon$ v\'erifiant ces conditions. Posons
 $\Xi_{\epsilon}=Z_{G}(\epsilon)/G_{\epsilon}$. C'est le m\^eme groupe qu'en 4.6 compte tenu du fait que $G_{\epsilon}=I_{\epsilon}$ puisque la torsion est int\'erieure.  On a vu que le groupe $\Xi_{\epsilon}^{\Gamma_{F}}$ agissait sur $G_{\epsilon}$ par automorphismes d\'efinis sur $F$. Pour simplifier, on note cette action comme une conjugaison. Soit $\mathfrak{u}$ un voisinage ouvert de $0$ dans $\mathfrak{g}_{\epsilon}(F)$ v\'erifiant les conditions suivantes
 
 - $X\in \mathfrak{u}$ si et seulement si sa partie semi-simple $X_{ss}$ appartient \`a $\mathfrak{u}$;
 
 - si $X\in \mathfrak{u}$ et $X'\in \mathfrak{g}_{\epsilon}(F)$ sont conjugu\'es par un \'el\'ement de $G_{\epsilon}(\bar{F})$, alors $X'\in \mathfrak{u}$;
 
 - $\mathfrak{u}$ est invariant par l'action de $\Xi_{\epsilon}^{\Gamma_{F}}$.
 
 On suppose $\mathfrak{u}$ assez petit, en particulier l'exponentielle y est d\'efinie.
 
   Pour tout $y\in {\cal Y}(\epsilon)$, on d\'efinit $\mathfrak{u}_{\epsilon[y]}$ comme en 4.6 et on pose $U_{\epsilon[y]}=exp(\mathfrak{u}_{\epsilon[y]})$ (simplement $U_{\epsilon}=exp(\mathfrak{u})$). On note $\tilde{U}$ l'ensemble des \'el\'ements de $\tilde{G}(F)$ dont la partie semi-simple est stablement  conjugu\'ee \`a un \'el\'ement  de $U_{\epsilon}\epsilon$. C'est l'ensemble du (i) du lemme 4.6. En effet,  pour $y\in \dot{\cal Y}(\epsilon)$, tout \'el\'ement  semi-simple de $U_{\epsilon[y]}\epsilon[ y]$ est stablement conjugu\'e \`a un \'el\'ement de $ U_{\epsilon}\epsilon$, cela parce que $G_{\epsilon}$ est quasi-d\'eploy\'e.  On d\'efinit une correspondance entre $C_{c}^{\infty}(\tilde{U})$ et $C_{c}^{\infty}((U_{\epsilon})\simeq C_{c}^{\infty}(\mathfrak{u})$ par: $f\in C_{c}^{\infty}(\tilde{U})$ et $\phi\in C_{c}^{\infty}(U_{\epsilon})$ se correspondent si et seulement si on a l'\'egalit\'e $S^{\tilde{G}}(x\epsilon,f)=S^{G_{\epsilon}}(x,\phi)$ pour tout \'el\'ement $x\in U_{\epsilon}$ tel que $x\epsilon$ soit fortement r\'egulier dans $\tilde{G}$. Avec des notations \'evidentes, on a le r\'esultat suivant.

 \ass{Lemme}{Cette correspondance se quotiente en un isomorphisme 
 $$desc_{\epsilon}^{st}:SI(\tilde{U}) \to SI(U_{\epsilon})^{\Xi_{\epsilon}^{\Gamma_{F}}}\simeq SI(\mathfrak{u})^{\Xi_{\epsilon}^{\Gamma_{F}}}.$$
  Si $\epsilon$ est elliptique dans $\tilde{G}(F)$, cet isomorphisme se restreint en un isomorphisme   
  $$SI_{cusp}(\tilde{U}) \to SI_{cusp}(U_{\epsilon})^{\Xi_{\epsilon}^{\Gamma_{F}}}\simeq SI_{cusp}(\mathfrak{u})^{\Xi_{\epsilon}^{\Gamma_{F}}}.$$}

 Preuve. Notons $C_{c}^{\infty}(\tilde{U})'$ et $C_{c}^{\infty}(\mathfrak{u})'$ les projections dans $C_{c}^{\infty}(\tilde{U})$ et $C_{c}^{\infty}(\mathfrak{u})$ du graphe de la correspondance. Notons $SI(\tilde{U})'$ et $SI(\mathfrak{u})'$ leurs images dans $SI(\tilde{U})$ et $SI(\mathfrak{u})$.  
   Puisque toute classe de conjugaison stable dans $\tilde{U}$ contient un \'el\'ement $exp(X)\epsilon$ avec $X\in \mathfrak{u}$, la correspondance se quotiente  alors en un isomorphisme entre   $SI(\tilde{U})'$ et  $SI(\mathfrak{u})'$. Ce dernier espace est inclus dans $SI(\mathfrak{u})^{\Xi_{\epsilon}^{\Gamma_{F}}}$: cela r\'esulte du fait que, pour $g\in \Xi_{\epsilon}^{\Gamma_{F}}$ et $X\in \mathfrak{u}$, l'\'el\'ement $exp(g^{-1}Xg)\epsilon$ est stablement conjugu\'e \`a $exp(X)\epsilon$. On va montrer que $C_{c}^{\infty}(\tilde{U})'=C_{c}^{\infty}(\tilde{U})$ tandis que $C_{c}^{\infty}(\mathfrak{u})'$ est l'espace des \'el\'ements de $C_{c}^{\infty}(\mathfrak{u})$ dont l'image dans $SI(\mathfrak{u})$ est invariante par $\Xi_{\epsilon}^{\Gamma_{F}}$.

 Soient $f\in C_{c}^{\infty}(\tilde{U})$ et $X$ un \'el\'ement r\'egulier de $ \mathfrak{u}$. Le lemme 4.6(iv) d\'ecrit un ensemble de repr\'esentants des classes de conjugaison par $G(F)$ dans la classe stable de $exp(X)\epsilon$. En appliquant les d\'efinitions, on obtient
 $$S^{\tilde{G}}(exp(X)\epsilon,f)=\sum_{y\in \dot{\cal Y}(\epsilon)}\sum_{X'\in\dot{\cal X}(X,y)}I^{\tilde{G}}(exp(X')\epsilon[ y],f).$$
On effectue la descente d'Harish-Chandra au voisinage de chaque point $\epsilon[ y]$. La fonction $f$ correspond ainsi \`a une fonction disons $\phi'_{y}\in C_{c}^{\infty}(\mathfrak{ u}_{\epsilon [y]})$. Le groupe quasi-d\'eploy\'e $G_{\epsilon}$ se compl\`ete de la fa\c{c}on habituelle en une donn\'ee endoscopique de $G_{\epsilon[ y]}$ et la fonction $\phi'_{y}$ se transf\`ere en une fonction $\phi_{y}\in C_{c}^{\infty}(\mathfrak{u})$. La formule pr\'ec\'edente devient
$$S^{\tilde{G}}(exp(X)\epsilon,f)=\sum_{y\in \dot{\cal Y}(\epsilon)}S^{G_{\epsilon}}(X,\phi_{y}).$$
Donc la fonction $\phi_{f}=\sum_{y\in \dot{\cal Y}(\epsilon)}\phi_{y}$ correspond \`a $f$. Cela prouve l'\'egalit\'e $C_{c}^{\infty}(\tilde{U})'=C_{c}^{\infty}(\tilde{U})$.

 Inversement, soit $\phi\in C_{c}^{\infty}(\mathfrak{u})$ dont l'image dans $SI(\mathfrak{u})$ est invariante par $\Xi_{ \epsilon}^{\Gamma_{F}}$. On a une inclusion $Z_{G}(\epsilon;F)/G_{\epsilon}(F)\subset \Xi_{\epsilon}^{\Gamma_{F}}$. Sans changer l'image de $\phi$ dans $SI(\mathfrak{u})$, on peut remplacer $\phi$ par la fonction
$$X\mapsto\vert Z_{G}(\epsilon;F)/G_{\epsilon}(F)\vert ^{-1}\sum_{g}\phi(g^{-1}Xg),$$
o\`u $g$ parcourt un ensemble de repr\'esentants de $Z_{G}(\epsilon;F)/G_{\epsilon}(F)$. On peut ainsi supposer que l'image de $\phi$ dans $I(\mathfrak{u})$ est invariante par $Z_{G}(\epsilon;F)$. Appliquant la descente d'Harish-Chandra, on peut trouver $f\in C_{c}^{\infty}(\tilde{U})$ qui correspond \`a $\phi$ et dont les int\'egrales orbitales sont nulles en tout point qui n'est pas conjugu\'e par un \'el\'ement de $G(F)$ \`a un \'el\'ement de  $exp(\mathfrak{u})\epsilon$. Appliquant la premi\`ere partie du raisonnement \`a cette fonction, on construit une fonction  $\phi_{f}\in C_{c}^{\infty}(\mathfrak{u})$ qui correspond \`a $f$. On va montrer que l'image de $\phi_{f}$ dans $SI(\mathfrak{u})$ est \'egale \`a celle de $N\phi$, o\`u $N$ est un entier non nul, ce qui ach\`evera la preuve de la premi\`ere assertion du lemme. On a une inclusion naturelle
$$  \Xi_{\epsilon}^{\Gamma_{F}}/Z_{G}(\epsilon;F)\to G_{\epsilon}\backslash {\cal Y}(\epsilon)/G(F).$$
Notons $\dot{\cal Y}_{0}(\epsilon)$ le sous-ensemble de $\dot{\cal Y}(\epsilon)$ repr\'esentant l'image de cette inclusion. On peut supposer que, pour $y\in \dot{\cal Y}_{0}(\epsilon)$, $\epsilon[ y]=\epsilon$ et l'automorphisme $ad_{y}$  de $G_{\epsilon}$ est un \'el\'ement de $\Xi_{\epsilon}^{\Gamma_{F}}$. On peut aussi supposer que $y=1$ appartient \`a $\dot{\cal Y}_{0}(\epsilon)$. Pour $y=1$, $\phi_{1}=\phi'_{1}$ a par d\'efinition m\^eme image que $\phi$ dans $I(\mathfrak{u})$, a fortiori dans $SI(\mathfrak{u})$. Pour $y\in \dot{\cal Y}_{0}(\epsilon)$, $\phi'_{y}=\phi'_{1}$ puisque $y^{-1}\epsilon y=\epsilon$. D'apr\`es la propri\'et\'e ci-dessus de $ad_{y}$, le transfert $\phi_{y}$ de $\phi'_{y}$ a m\^eme image dans $SI(\mathfrak{u})$ que l'image de $\phi$ par l'action d'un \'el\'ement de $\Xi_{\epsilon}^{\Gamma_{F}}$. Puisque cette derni\`ere image est invariante par ce groupe, $\phi_{y}$ a m\^eme image que $\phi$ dans $SI(\mathfrak{u})$. Pour $y\in \dot{\cal Y}(\epsilon)- \dot{\cal Y}_{0}(\epsilon)$, aucun \'el\'ement de $ U_{\epsilon[y]}\epsilon [y]$ n'est conjugu\'e par un \'el\'ement de $G(F)$ \`a un \'el\'ement de $ U_{\epsilon}\epsilon$. Sinon, en supposant $\mathfrak{u}$ assez petit, cela entra\^{\i}nerait que $\epsilon[ y]$ serait conjugu\'e \`a $\epsilon$ par un \'el\'ement de $G(F)$ et on voit que cela contredirait l'hypoth\`ese que $y\not\in \dot{\cal Y}_{0}(\epsilon)$. On peut donc supposer $\phi'_{y}=0$ pour ces $y$ et on conclut comme on le voulait que l'image de $\phi_{f}$ dans $SI(\mathfrak{u})$ est \'egale \`a celle de $\vert \dot{\cal Y}_{0}(\epsilon)\vert \phi$. Cela ach\`eve la preuve de la premi\`ere assertion de l'\'enonc\'e.
           
    Si $\epsilon$ est elliptique, pour $X\in \mathfrak{u}$ r\'egulier, $X$ est elliptique dans $\mathfrak{g}_{\epsilon}(F)$ si et seulement si $exp(X)\epsilon$ est elliptique dans $\tilde{G}(F)$.   Il en r\'esulte que l'isomorphisme de la premi\`ere assertion conserve la cuspidalit\'e. $\square$

    {\bf Variante.} Supposons donn\'ee une extension
    $$1\to C_{1}\to G_{1}\to G\to 1$$
    o\`u $C_{1}$ est un tore central induit, une extension compatible
    $$\tilde{G}_{1}\to \tilde{G}$$
    avec $\tilde{G}_{1}$ \`a torsion int\'erieure et un caract\`ere $\lambda_{1}$ de $C_{1}(F)$. Soit $\epsilon$  comme pr\'ec\'edemment. Fixons $\epsilon_{1}\in \tilde{G}_{1}(F)$ se projetant sur $\epsilon$. On a une suite exacte
    $$0\to \mathfrak{c}_{1}\to \mathfrak{g}_{1,\epsilon_{1}}\to \mathfrak{g}_{\epsilon}\to 0$$
    On a besoin de scinder convenablement cette suite. La partie semi-simple de $\mathfrak{g}_{\epsilon}$ se scinde canoniquement par le diagramme
    $$\begin{array}{ccc}\mathfrak{g}_{1,\epsilon_{1},SC}&\simeq&\mathfrak{g}_{\epsilon,SC}\\ \downarrow&&\downarrow\\ \mathfrak{g}_{1,\epsilon_{1}}&\to&\mathfrak{g}_{\epsilon}\\ \end{array}$$
   Notons $Z_{\epsilon_{1}}$ et $Z_{\epsilon}$ les centres de $G_{1,\epsilon_{1}}$ et $G_{\epsilon}$. Le groupe $Z_{G}(\epsilon)$ agit par conjugaison sur $\tilde{G}_{1}$.  Cette action conserve $G_{1,\epsilon_{1}}$. En effet, un \'el\'ement $g\in Z_{G}(\epsilon)$ envoie $\epsilon_{1}$ sur $c(g)\epsilon_{1}$ pour un unique $c(g)\in C_{1}$, donc envoie $G_{1,\epsilon_{1}}$ sur $G_{1,c(g)\epsilon_{1}}=G_{1,\epsilon_{1}}$. L'action de $Z_{G}(\epsilon)$ se restreint en une action sur $Z_{\epsilon_{1}}$, qui est l'identit\'e sur $C_{1}$. On peut alors trouver une d\'ecomposition
   $$\mathfrak{z}_{\epsilon_{1}}=\mathfrak{c}_{1}\oplus \mathfrak{s}$$
   stable pour les actions de $\Gamma_{F}$ et de $Z_{G}(\epsilon)$. On fixe une telle d\'ecomposition. La projection $\mathfrak{g}_{1,\epsilon_{1}}\to \mathfrak{g}_{\epsilon}$ se restreint en un isomorphisme
   $$\mathfrak{s}\oplus \mathfrak{g}_{1,\epsilon_{1},SC}\to \mathfrak{g}_{\epsilon}$$
   et on prend pour section l'isomorphisme r\'eciproque. Soit $\mathfrak{u}$ un voisinage comme pr\'ec\'edemment, que l'on identifie par la section \`a un sous-ensemble de $\mathfrak{g}_{1,\epsilon_{1}}(F)$. On note $\tilde{U}_{1}$ l'image r\'eciproque de $\tilde{U}$ dans $\tilde{G}_{1}(F)$ et on d\'efinit l'espace $SI_{\lambda_{1}}(\tilde{U}_{1})$, quotient de $C_{c,\lambda_{1}}^{\infty}(\tilde{U}_{1})$ par le sous-espace des fonctions dont les int\'egrales orbitales stables sont nulles. On d\'efinit comme pr\'ec\'edemment  une correspondance naturelle entre  $C_{c,\lambda_{1}}^{\infty}(\tilde{U}_{1})$  et $C_{c}^{\infty}(\mathfrak{u})$. Comme on l'a dit ci-dessus, un \'el\'ement de $Z_{G}(\epsilon)$ envoie $\epsilon_{1}$ sur $c(g)\epsilon_{1}$ pour un unique $c(g)\in C_{1}$. On a $c(g)=1$ pour $g\in G_{\epsilon}$. D'autre part, si l'image de $g$ dans $Z_{G}(\epsilon)/G_{\epsilon}=\Xi_{\epsilon}$ est fixe par $\Gamma_{F}$, $c(g)$ appartient \`a $C_{1}(F)$. On obtient un caract\`ere $g\mapsto \lambda_{1}(c(g)^{-1})$ du groupe $\Xi_{\epsilon}^{\Gamma_{F}}$. Alors
   
   (1) la correspondance ci-dessus se quotiente en un isomorphisme entre $SI_{\lambda_{1}}(\tilde{U}_{1})$ et le sous-espace des \'el\'ements de $SI(\mathfrak{u})$ qui se transforment selon ce caract\`ere de $\Xi_{\epsilon}^{\Gamma_{F}}$.
   
   Consid\'erons maintenant d'autres extensions
   $$1\to C_{2}\to G_{2}\to G\to 1,\,\, \tilde{G}_{2}\to \tilde{G}$$
   et un caract\`ere $\lambda_{2}$ de $C_{2}(F)$, v\'erifiant les m\^emes conditions que ci-dessus. Introduisons comme en 2.5 les produits fibr\'es $G_{12}$ et $\tilde{G}_{12}$ et supposons donn\'es un caract\`ere $\lambda_{12}$ de $G_{12}(F)$ et une fonction non nulle $\tilde{\lambda}_{12}$ sur $\tilde{G}_{12}(F)$ v\'erifiant les conditions de ce paragraphe, c'est-\`a-dire
   
   - la restriction de $\lambda_{12}$ \`a $C_{1}(F)\times C_{2}(F)$ est $\lambda_{1}\times \lambda_{2}^{-1}$;
   
   - pour $(\gamma_{1},\gamma_{2})\in \tilde{G}_{12}(F)$ et $(x_{1},x_{2})\in G_{12}(F)$, on a l'\'egalit\'e $\tilde{\lambda}_{12}(x_{1}\gamma_{1},x_{2}\gamma_{2})=\lambda_{12}(x_{1},x_{2})\tilde{\lambda}_{12}(\gamma_{1},\gamma_{2})$. 
   
   Par la construction ci-dessus, chaque s\'erie de donn\'ees d\'efinit un caract\`ere de $\Xi_{\epsilon}^{\Gamma_{F}}$. On a
   
   (2) ces caract\`eres sont \'egaux. 
   
   Fixons $\epsilon_{1}$ comme plus haut et $\epsilon_{2}$ de fa\c{c}on similaire. Soit $g\in Z_{G}(\epsilon)$ s'envoyant sur un \'el\'ement de $\Xi_{\epsilon}^{\Gamma_{F}}$. Pour $i=1,2$, on a $ad_{g}(\epsilon_{i})=c_{i}(g)\epsilon_{i}$ avec $c_{i}(g)\in C_{i}(F)$. Il s'agit de prouver que $\lambda_{1}(c_{1}(g))=\lambda_{2}(c_{2}(g))$. En posant $\epsilon_{12}=(\epsilon_{1},\epsilon_{2})$ et $\epsilon'_{12}=(ad_{g}(\epsilon_{1}),ad_{g}(\epsilon_{2}))$, 
   il revient au m\^eme de prouver que $\tilde{\lambda}_{12}(\epsilon_{12})=\tilde{\lambda}_{12}(\epsilon'_{12})$. Puisque $G_{12}$ est quasi-d\'eploy\'e, il co\"{\i}ncide avec le groupe $G_{12,0}$  qu'on lui a associ\'e en 1.12. Il en r\'esulte que l'application $N^{\tilde{G}_{12}}$ se quotiente en l'injection $\pi(G_{12,SC}(F))\backslash \tilde{G}_{12}(F)\to \tilde{G}_{12,ab}(F)$. Par construction, les \'el\'ements $\epsilon_{12}$ et $\epsilon'_{12}$ sont stablement conjugu\'es. D'apr\`es le lemme 4.5, on a $N^{\tilde{G}_{12}}(\epsilon_{12})=N^{\tilde{G}_{12}}(\epsilon'_{12})$, donc $\epsilon'_{12}\in \pi(G_{12,SC}(F))\epsilon_{12}$. Le caract\`ere $\lambda_{12}$ est forc\'ement trivial sur $\pi(G_{12,SC}(F))$. Donc $\tilde{\lambda}_{12}(\epsilon_{12})=\tilde{\lambda}_{12}(\epsilon'_{12})$ comme on le voulait. Cela prouve (2). 

 \bigskip
 
 \subsection{Conjugaison stable et endoscopie}
 Soit ${\bf G}'$ une donn\'ee endoscopique relevante pour $(G,\tilde{G},{\bf a})$. Fixons  un diagramme $(\epsilon,B',T',B,T,\eta)$. On fixe une forme quasi-d\'eploy\'ee $\bar{G}$ de $G_{\eta}$. On fixe de m\^eme une forme quasi-d\'eploy\'ee $G^{_{'}*}_{\epsilon}$ de $G'_{\epsilon}$. A l'aide du diagramme, on a construit en [W1] 3.5 une donn\'ee endoscopique ${\bf \bar{G}'}=(\bar{G}',\bar{{\cal G}}',\bar{s})$ de $\bar{G}_{SC}$. Il s'agit d'endoscopie usuelle, il n'y a ici ni torsion, ni caract\`ere. Les deux groupes $G^{_{'}*}_{\epsilon,SC}$ et $\bar{G}'_{SC}$ forment une paire endoscopique non standard ([W1] 1.7). Pr\'ecisons les correspondances de tores. Fixons des paires de Borel dans chacun des groupes, dont on note les tores $\bar{T}$ pour $\bar{G}$, $\bar{T}'$ pour $\bar{G}'$ et $T^{_{'}*}$ pour $G^{_{'}*}_{\epsilon}$. Si on oublie les actions galoisiennes, on peut identifier $\bar{T}$ \`a $T^{\theta,0}$, o\`u $\theta=ad_{\eta}$, et $T^{_{'}*}$ \`a $T'$. De l'homomorphisme $\xi_{T,T'}$ se d\'eduit un isomorphisme
 $$X_{*}(\bar{T})\otimes {\mathbb Q}\to X_{*}(T^{_{'}*})\otimes {\mathbb Q}.$$
 De m\^eme, on peut choisir un homomorphisme $\xi_{\bar{T}_{sc},\bar{T}' }$ (qui est un isomorphisme puisque la situation n'est pas tordue), d'o\`u un isomorphisme
 $$X_{*}(\bar{T}_{sc})\otimes {\mathbb Q}\to X_{*}(\bar{T}')\otimes {\mathbb Q}.$$
 Enfin, sous-jacent \`a la notion d'endoscopie non standard, il y a un isomorphisme
 $$X_{*}(T^{_{'}*}_{sc})\otimes {\mathbb Q}\to X_{*}(\bar{T}'_{sc})\otimes {\mathbb Q},$$
 qui, lui, est \'equivariant pour les actions galoisiennes. Ces homomorphismes sont compatibles. De plus, il s'en d\'eduit un isomorphisme
 $$X_{*}(Z(G'_{\epsilon})^0)\otimes {\mathbb Q}\to (X_{*}(Z(\bar{G})^0)\otimes {\mathbb Q})\oplus (X_{*}(Z(\bar{G}')^0)\otimes {\mathbb Q})$$
 qui est compatible aux actions galoisiennes. 
 
 Ces isomorphismes induisent des correspondances compatibles entre classes de conjugaison stable d'\'el\'ements semi-simples r\'eguliers dans les alg\`ebres de Lie des diff\'erents groupes.
 
 Rappelons que l'on dit que $\epsilon$ et $\eta$ se correspondent s'il existe un diagramme les joignant.
 
 {\bf Remarque.} Si $\epsilon$ et $\eta$ se correspondent, il existe un diagramme $(\epsilon,B',T',B,T,\eta)$ tel que $T'$ est un tore elliptique de $G'_{\epsilon}$ si $F$ est non-archim\'edien, resp. est un tore fondamental de $G'_{\epsilon}$ si $F$ est archim\'edien. A l'aide des rappels ci-dessus, cela r\'esulte que, puisque ${\bf \bar{G}}'$ est une donn\'ee endoscopique relevante de $G_{\eta,SC}$, tout sous-tore maximal elliptique, resp. fondamental, de $\bar{G}'$ se transf\`ere \`a $G_{\eta,SC}$.
 \bigskip
 
  Cette correspondance induit une correspondance entre classes de conjugaison stable d'\'el\'ements semi-simples dans $\tilde{G}'(F)$ et $\tilde{G}(F)$. Pr\'ecis\'ement, pour de tels \'el\'ements
 
 (1) si $\epsilon$ correspond \`a $\eta$ et $\eta'$, alors $\eta$ et $\eta'$ sont stablement conjugu\'es;
 
 (2) si $\epsilon$ correspond \`a $\eta$ et $\epsilon'$ est stablement conjugu\'e \`a $\epsilon$, alors $\epsilon'$ correspond \`a $\eta$;
 
 (3) si $\epsilon$ correspond \`a $\eta$ et si $\tilde{\alpha}_{x}$ est un automorphisme d\'efini sur $F$ de $\tilde{G}'$ provenant d'un \'el\'ement $x\in Aut({\bf G}')$, alors $\tilde{\alpha}_{x}(\epsilon)$ correspond \`a $\eta$.
 
 Le (1) est le lemme 3.4 de [W1]. Pour (2),  d'apr\`es la remarque ci-dessus,  s'il existe un diagramme $(\epsilon,B',T',B,T,\eta)$, on peut le remplacer par un autre o\`u $T'$ est elliptique ou fondamental dans $G'_{\epsilon}$. Un tel tore se transf\'erant \`a toute forme int\'erieure, (2) s'ensuit. Le (3) r\'esulte des d\'efinitions.
 
 Remarquons que les assertions r\'eciproques de (1) et (2) sont fausses en g\'en\'eral. La r\'eciproque de (1) devient toutefois vraie si ${\bf G}'$ est elliptique ainsi que $\epsilon$ (avec notre d\'efinition: $\epsilon$ est elliptique si il appartient \`a un sous-tore maximal elliptique de $G'$). D'autre part, parce que l'on sait que dans la classe de conjugaison stable de $\epsilon$, il y a toujours un \'el\'ement dont le commutant connexe est quasi-d\'eploy\'e, (2) nous permet de nous limiter \`a consid\'erer des $\epsilon$ v\'erifiant cette propri\'et\'e.
 
 Restreignons-nous maintenant aux \'el\'ements elliptiques. Pour un \'el\'ement semi-simple elliptique $\eta\in \tilde{G}(F)$, consid\'erons les couples $({\bf G}',\epsilon)$ o\`u ${\bf G}'$ est une donn\'ee endoscopique elliptique de $(G,\tilde{G},{\bf a})$ et $\epsilon\in \tilde{G}'(F)$ est un \'el\'ement semi-simple elliptique qui correspond \`a $\eta$ et dont le commutant connexe $G'_{\epsilon}$ est quasi-d\'eploy\'e. Disons que deux couples $({\bf G}'_{1},\epsilon_{1})$ et $({\bf G}_{2}',\epsilon_{2})$ sont \'equivalents si et seulement s'il existe un isomorphisme $\tilde{\alpha}:\tilde{G}'_{1}\to \tilde{G}'_{2}$ d\'efini sur $F$ et provenant d'une \'equivalence entre ${\bf G}'_{1}$ et ${\bf G}'_{2}$ de sorte que $\epsilon_{2}$ soit stablement conjugu\'e \`a $\tilde{\alpha}(\epsilon_{1})$. On fixe un ensemble $\dot{\cal X}^{{\cal E}}(\eta)$ de repr\'esentants des classes d'\'equivalence de ces couples. Pour tout $({\bf G}',\epsilon)\in \dot{\cal X}^{{\cal E}}(\eta)$, on fixe des donn\'ees auxiliaires $G'_{1},...,\Delta_{1}$ (notons que ${\bf G}'$ est forc\'ement relevant) et un \'el\'ement $\epsilon_{1}\in \tilde{G}'_{1}(F)$ qui rel\`eve $\epsilon$.
 
 Consid\'erons d'abord le cas o\`u $\eta$ est fortement r\'egulier et $F$ est non archim\'edien. On a d'abord
 
 - si $\omega$ n'est pas trivial sur $Z_{G}(\eta;F)$, alors $\dot{\cal X}^{{\cal E}}(\eta)=\emptyset$, cf. [KS1] lemme 4.4.C.
 
Supposons $\omega$ trivial sur $Z_{G}(\eta;F)$. Fixons un ensemble de repr\'esentants $\dot{\cal X}(\eta)$ des classes de conjugaison par $G(F)$ dans la classe de conjugaison stable de $\eta$. On d\'efinit les deux applications lin\'eaires
 $$(4) \qquad\begin{array}{ccc}{\mathbb C}^{\dot{\cal X}(\eta)}&\to&{\mathbb C}^{\dot{\cal X}^{{\cal E}}(\eta)}\\ (x_{\eta'})_{\eta'\in \dot{\cal X}(\eta)}&\mapsto &(y_{({\bf G}',\epsilon)})_{({\bf G}',\epsilon)\in \dot{\cal X}^{{\cal E}}(\eta)}\\ \end{array}$$
 o\`u
 $$y_{({\bf G}',\epsilon)}=d(\theta^*)^{1/2} \sum_{\eta'\in \dot{\cal X}(\eta)}\Delta_{1}(\epsilon_{1},\eta')[Z_{G}(\eta';F):G_{\eta'}(F)]^{-1}x_{\eta'}$$
 (le $\Delta_{1}$ est bien s\^ur celui de ${\bf G}'$);
 $$(5) \qquad\begin{array}{ccc}{\mathbb C}^{\dot{\cal X}^{\cal E}(\eta)}&\to&{\mathbb C}^{\dot{\cal X}^(\eta)}\\ (y_{({\bf G}',\epsilon)})_{({\bf G}',\epsilon)\in \dot{\cal X}^{{\cal E}}(\eta)}&\mapsto&(x_{\eta'})_{\eta'\in \dot{\cal X}(\eta)} \\ \end{array}$$
 o\`u
 $$ x_{\eta'}=[Z_{G}(\eta';F):G_{\eta'}(F)]\vert \dot{\cal X}(\eta)\vert ^{-1}d(\theta^*)^{-1/2} \sum_{({\bf G}',\epsilon)\in \dot{\cal X}^{{\cal E}}(\eta)}\Delta_{1}(\epsilon_{1},\eta')^{-1}y_{({\bf G}',\epsilon)}.$$
 L'assertion fondatrice de la th\'eorie de l'endoscopie tordue est que ces deux applications lin\'eaires sont inverses l'une de l'autre. On renvoie pour cette assertion \`a Kottwitz-Shelstad ([KS]) et \`a Labesse ([Lab2]), bien que ces auteurs d\'etaillent plut\^ot le cas o\`u le corps de base est un corps de nombres. 
 
 Dans le cas o\`u $F={\mathbb R}$, on doit consid\'erer un $K$-espace tordu. Pour $\eta\in K\tilde{G}_{reg}({\mathbb R})$, on d\'efinit sans changement l'ensemble $\dot{\cal X}^{{\cal E}}(\eta)$. On fixe  pour tout $p\in \Pi$ un ensemble de repr\'esentants $\dot{\cal X}_{p}(\eta)$ des classes de conjugaison par $G_{p}({\mathbb R})$ dans l'intersection de $\tilde{G}_{p}({\mathbb R})$ avec la classe de conjugaison stable de $\eta$. On pose $\dot{\cal X}(\eta)=\sqcup_{p\in \Pi}\dot{\cal X}_{p}(\eta)$. Avec ces d\'efinitions, les applications (4) et (5) sont encore inverses l'une de l'autre. C'est la raison d'\^etre des $K$-espaces tordus. 
 
 La correspondance entre \'el\'ements semi-simples elliptiques non fortement r\'eguliers est plus compliqu\'ee. L'important pour nous est qu'elle forme un "bord" satisfaisant \`a celle des \'el\'ements fortement r\'eguliers. Notons $\tilde{G}_{ss}(F)_{ell}$ l'ensemble des \'el\'ements semi-simples elliptiques de $\tilde{G}(F)$, pas forc\'ement r\'eguliers. Notons $\tilde{G}_{ss}(F)_{ell}/st-conj$ l'ensemble des classes de conjugaison stable dans $\tilde{G}_{ss}(F)_{ell}$.   Soit ${\bf G}'$ une donn\'ee endoscopique elliptique pour $(G,\tilde{G},{\bf a})$. On d\'efinit de m\^eme l'espace $\tilde{G}'_{ss}(F)_{ell}/st-conj$. D'apr\`es le lemme 4.5, l'application $N^{\tilde{G}'}$ restreinte \`a $\tilde{G}'_{ss}(F)_{ell}$ se factorise par cet ensemble de classes de conjugaison stable. A fortiori, l'application $N^{\tilde{G}',\tilde{G}} $ se factorise de m\^eme. Dans le cas o\`u $F$ est non-archim\'edien, on note $\tilde{G}'_{ss}(F)_{ell}^{\tilde{G}}/st-conj$ l'ensemble des \'el\'ements de $\tilde{G}'_{ss}(F)_{ell}/st-conj$ dont l'image par cette application appartient \`a l'image de $\tilde{G}_{ab}(F)$ par $N^{\tilde{G}}$. Dans le cas o\`u $F={\mathbb R}$ et o\`u on travaille avec des $K$-espaces tordus, on pose  la m\^eme d\'efinition en rempla\c{c}ant $\tilde{G}_{ab}({\mathbb R})$ par $K\tilde{G}_{ab}({\mathbb R})$. Montrons que
 
 (6) un \'el\'ement $\epsilon\in \tilde{G}'_{ss}(F)_{ell}$ correspond \`a un \'el\'ement semi-simple de $\tilde{G}(F)$ (ou de $K\tilde{G}({\mathbb R})$) si et seulement si sa classe de conjugaison stable appartient \`a $\tilde{G}'_{ss}(F)_{ell}^{\tilde{G}}/st-conj$.
 
 Preuve. Supposons qu'il existe un diagramme $(\epsilon,B',T',B,T,\eta)$.
   Pour $X\in \mathfrak{t}^{\theta}(F)$ assez petit et en position g\'en\'erale  et pour  $Y=\xi_{T,T'}(X)$, les \'el\'ements $exp(Y)\epsilon$ et $exp(X)\eta$ sont fortement r\'eguliers et se correspondent. D'apr\`es la proposition 1.14(i), l'image par $N^{\tilde{G}',\tilde{G}}$ de $exp(Y)\epsilon$ appartient \`a l'image de $\tilde{G}_{ab}(F)$ par $N^{\tilde{G}}$. Cette image \'etant ferm\'ee, $N^{\tilde{G}',\tilde{G}}(\epsilon)$ lui appartient aussi et l'image de $\epsilon$ dans $\tilde{G}'_{ss}(F)_{ell}/st-conj$ appartient \`a $\tilde{G}'_{ss}(F)_{ell}^{\tilde{G}}/st-conj$. Inversement, supposons que cette condition soit v\'erifi\'ee. Supposons pour simplifier $F$ non archim\'edien, l'extension aux $K$-espaces \'etant similaire.  Puisque $\epsilon$ est elliptique, on peut fixer un sous-tore maximal $T'$ de $G'_{\epsilon}$, d\'efini sur $F$ et elliptique dans $G'$. Pour $Y\in \mathfrak{t}'(F)$ assez petit et en position g\'en\'erale, $exp(Y)\epsilon$ est elliptique r\'egulier et son image par $N^{\tilde{G}',\tilde{G}}$ appartient \`a l'image de $\tilde{G}_{ab}(F)$ par $N^{\tilde{G}}$. Par la proposition 1.14(ii), il existe $\gamma\in \tilde{G}(F)_{reg}$ tel que $(exp(Y)\epsilon,\gamma)\in {\cal D}$. On peut fixer un diagramme joignant $exp(Y)\epsilon$ et $\gamma$. Le tore $T'$ de ce diagramme est impos\'e: c'est le commutant de $exp(Y)\epsilon$, donc c'est le tore $T'$ d\'ej\`a introduit. Notons $(exp(Y)\epsilon,B',T',B,T,\gamma)$ ce diagramme. Puisque $exp(Y)\epsilon$ conserve $(B',T')$ et $Y\in \mathfrak{t}'(F)$, $\epsilon$ conserve lui-aussi $(B',T')$. De l'application $\xi_{T,T'}$ r\'esulte un isomorphisme $\mathfrak{t}^{\theta}(F)\to \mathfrak{t}'(F)$. Soit $X\in \mathfrak{t}^{\theta}(F)$ correspondant \`a $Y$, posons $\eta=exp(-X)\gamma$. Par le m\^eme argument, $\eta$ conserve $(B,T)$. Alors $(\epsilon,B',T',B,T,\eta)$ est un diagramme. $\square$

 D'apr\`es (1) ci-dessus, et en remarquant qu'un \'el\'ement elliptique de $\tilde{G}'(F)$ ne peut correspondre qu'\`a un \'el\'ement elliptique de $\tilde{G}(F)$, on a une application  
 $$(7) \qquad\tilde{G}'_{ss}(F)_{ell}^{\tilde{G}}/st-conj\to \tilde{G}_{ss}(F)_{ell}/st-conj.$$
 Munissons $\tilde{G}_{ss}(F)_{ell}$ de la topologie induite par celle de $\tilde{G}(F)$ et $\tilde{G}_{ss}(F)_{ell}/st-conj$ de la topologie la moins fine pour laquelle la projection $\tilde{G}_{ss}(F)_{ell}\to \tilde{G}_{ss}(F)_{ell}/st-conj$ est continue. On munit de m\^eme $\tilde{G}'_{ss}(F)_{ell}/st-conj$ d'une topologie et $\tilde{G}'_{ss}(F)_{ell}^{\tilde{G}}/st-conj$ de la topologie induite.
 
 \ass{Lemme}{L'espace $\tilde{G}_{ss}(F)_{ell}/st-conj$ est s\'epar\'e et localement compact. La projection $\tilde{G}_{ss}(F)_{ell}\to \tilde{G}_{ss}(F)_{ell}/st-conj$ est ouverte. L'application (7) est continue  et propre.}
 
 Preuve. Soient $\eta_{1}$ et $\eta_{2}$ deux \'el\'ements de $\tilde{G}_{ss}(F)_{ell}$ qui ne sont pas stablement conjugu\'es. On  construit comme en 4.6 des voisinages $\tilde{U}_{1}$ et $\tilde{U}_{2}$ de $\eta_{1}$ et $\eta_{2}$. La caract\'erisation  du lemme 4.6(i) montre que l'on peut les construire disjoints. Ils sont invariants par conjugaison stable. Alors leurs images dans $\tilde{G}_{ss}(F)_{ell}/st-conj$ sont des voisinages disjoints des images de $\eta_{1}$ et $\eta_{2}$. Pour un seul \'el\'ement $\eta$, construisons un voisinage $\tilde{U}$ comme en 4.6 issu d'un voisinage $\mathfrak{u}$ de $0$ dans $\mathfrak{g}_{\eta}(F)$ qui est compact modulo conjugaison par $G_{\eta}(F)$. Alors son image dans $\tilde{G}_{ss}(F)_{ell}/st-conj$ est un voisinage compact de l'image de $\eta$. Cela prouve les deux premi\`eres assertions de l'\'enonc\'e. Par ailleurs, l'image de $\tilde{U}$ est \'egale \`a celle de $exp(\mathfrak{u})\eta$. En effet, un \'el\'ement  de $\tilde{U}$ est conjugu\'e par $G(F)$ \`a un \'el\'ement $exp(X)\eta [y]$ pour un $y\in \dot{{\cal Y}}(\eta)$ et $X\in \mathfrak{u}_{\eta[ y]}$. Si l'\'el\'ement est semi-simple elliptique, il existe un sous-tore maximal elliptique $T_{\natural}$ de $G_{\eta[ y]}$ tel que $X\in \mathfrak{t}_{\natural}(F)$.  Parce que ce tore est elliptique, il se transf\`ere par le torseur $ad_{y}$ en un sous-tore elliptique de $G_{\eta}$ et notre \'el\'ement est stablement conjugu\'e \`a un \'el\'ement de $exp(\mathfrak{u})\eta$. Cela prouve l'assertion. Mais alors l'image de  $exp(\mathfrak{u})\eta$  dans $\tilde{G}_{ss}(F)_{ell}/st-conj$ est un voisinage de celle de $\eta$. Puisqu'on peut prendre $\mathfrak{u}$ aussi petit que l'on veut, modulo conjugaison par $G_{\eta}(F)$, la projection $\tilde{G}_{ss}(F)_{ell}\to S\tilde{G}_{ss}(F)_{ell}$ est ouverte.   Puisque l'application $\tilde{G}'_{ss}(F)_{ell}\to \tilde{G}'_{ss}(F)_{ell}/st-conj$ est ouverte, il suffit, pour prouver la continuit\'e de (7), de prouver que l'application compos\'ee $\tilde{G}'_{ss}(F)_{ell}\to \tilde{G}_{ss}(F)_{ell}/st-conj$ l'est. Soient $\epsilon\in \tilde{G}'_{ss}(F)_{ell}$ et $\eta\in \tilde{G}_{ss}(F)_{ell}$ qui se correspondent. Pour tout \'el\'ement $\epsilon'$ de $\tilde{G}'_{ss}(F)_{ell}$ assez proche de $\epsilon$, il y a un sous-tore elliptique $T'$ de $G'_{\epsilon}$ tel que $\epsilon'=exp(Y)\epsilon$, avec $Y\in \mathfrak{t}'(F)$ et $Y$ proche de $0$. Puisqu'il n'y a \`a conjugaison pr\`es qu'un nombre fini de tores elliptiques $T'$, on peut fixer celui-ci.
  On voit en pr\'ecisant ce que l'on a dit plus haut que  l'on peut fixer un diagramme $(\epsilon,B',T',B,T,\eta)$ o\`u $T'$ est le tore fix\'e. En fixant une section de l'homormophisme $\xi_{T,T'}:\mathfrak{t}(F)\to \mathfrak{t}'(F)$, on voit que, quand $Y$ tend vers $0$ dans $\mathfrak{t}'(F)$, l'\'el\'ement $exp(Y)\epsilon$ correspond \`a un \'el\'ement $exp(X)\eta$ avec $X\in \mathfrak{t}(F)$ tendant vers $0$. Cela prouve la continuit\'e de (7). Soit maintenant   $\eta\in \tilde{G}_{ss}(F)_{ell}$. Fixons un ensemble de repr\'esentants $\underline{{\cal X}}$ des classes de conjugaison par $G'(F)$ dans l'ensemble des \'el\'ements de $\tilde{G}'_{ss}(F)_{ell}$ qui correspondent \`a $\eta$. C'est un ensemble fini puisqu'il est en tout cas inclus dans un ensemble fini de classes de conjugaison par $G(\bar{F})$. Soit $(\epsilon_{n},\eta_{n})_{n\in {\mathbb N}}$ une suite de couples qui se correspondent dans $\tilde{G}'_{ss}(F)_{ell}\times\tilde{G}_{ss}(F)_{ell}$ telle que $\eta_{n}$ tend vers $\eta$. Un raisonnement similaire \`a celui de la preuve du lemme 1.10(ii) montre que, quitte \`a  remplacer $\epsilon_{n}$ par un \'el\'ement stablement conjugu\'e, on peut supposer que $\epsilon_{n}$ appartient \`a un voisinage arbitraire de $\underline{{\cal X}}$ quand $n$ est assez grand. Autrement dit, l'image dans $\tilde{G}'_{ss}(F)_{ell}/st-conj$ d'un voisinage de $\underline{{\cal X}}$ contient l'image r\'eciproque par (7) d'un voisinage assez petit de l'image de $\eta$ dans  $\tilde{G}_{ss}(F)_{ell}/st-conj$. Cela entra\^{\i}ne que (7) est propre. $\square$
  
  \bigskip
  
On peut pr\'eciser la derni\`ere assertion de la fa\c{c}on suivante. Soit $\eta\in \tilde{G}_{ss}(F)_{ell}$.  On fixe comme plus haut un ensemble $\dot{\cal X}^{\cal E}(\eta)$. Pour tout $({\bf G}',\epsilon)\in \dot{\cal X}^{\cal E}(\eta)$, fixons un voisinage $U'_{\epsilon}$ de $\epsilon$ dans $\tilde{G}'(F)$. Alors il existe un voisinage $U$ de $\eta$ dans $\tilde{G}(F)$ tel que, pour tout $\gamma\in U$ elliptique r\'egulier, on peut choisir pour $\dot{\cal X}^{\cal E}(\gamma)$ un ensemble tel que, pour tout \'el\'ement $({\bf G}',\delta)$ de cet ensemble, il existe $\epsilon$ tel que $({\bf G}',\epsilon)\in \dot{\cal X}^{\cal E}(\eta)$ et $\delta\in U'_{\epsilon}$.
 
 \bigskip
 
 \subsection{Rappels sur la transformation de Fourier et l'endoscopie}  
 
    Supposons $F$ non-archim\'edien, $\tilde{G}=G$ et $\omega=1$.  La th\'eorie de l'endoscopie vaut aussi pour les alg\`ebres de Lie, avec quelques simplifications. Par exemple, pour une donn\'ee endoscopique ${\bf G}'$, les donn\'ees auxiliaires $G'_{1}$, $C_{1}$ et $\hat{\xi}_{1}$ ne servent plus \`a rien.  Modulo le choix d'un facteur de transfert, on peut poser $C_{c}^{\infty}(\boldsymbol{\mathfrak{g}}')=C_{c}^{\infty}(\mathfrak{g}'(F))$. Fixons une  transformation de Fourier dans $C_{c}^{\infty}(\mathfrak{g}(F))$ comme en 4.1. Elle  en d\'etermine une dans $C_{c}^{\infty}(\mathfrak{g}'(F))$, cf. [W1]. Elle se quotiente en une transformation de $SI(\mathfrak{g}'(F))$. On a
  
  (1) il existe un nombre complexe non nul $\gamma(\mathfrak{g})$ tel que, pour toute donn\'ee endoscopique ${\bf G}'$  et toutes $f\in I(\mathfrak{g}(F))$, $f'\in SI(\mathfrak{g}'(F))$, l'\'egalit\'e $f'=transfert(f)$ \'equivaut \`a $\gamma(\mathfrak{g}')\hat{f}'=transfert(\gamma(\mathfrak{g})\hat{f})$.
  
  Arthur a prouv\'e en [A2] lemme 3.4 que
  
  (2) l'homomorphisme de transfert
  $$I(\mathfrak{g}(F))\to \oplus_{{\bf G}'\in {\cal E}(G)}SI(\mathfrak{g}'(F))$$
  se restreint en un isomorphisme
   $$I_{cusp}(\mathfrak{g}(F))\simeq \oplus_{{\bf G}'\in {\cal E}(G)}SI_{cusp}(\mathfrak{g}'(F))^{Aut({\bf G}')}.$$
   
   {\bf Remarques.} (3) L'action de $Aut({\bf G}')$ est d\'efinie comme en 2.6.   On peut d\'efinir une action intrins\`eque de $Aut({\bf G}')$ dans $\mathfrak{g}'(F)$ mais  l'action que l'on consid\`ere est cette action intrins\`eque tordue par un caract\`ere qui tient compte du facteur de transfert.
   
   (4) Supposons $G$ quasi-d\'eploy\'e. Par d\'efinition, $SI_{cusp}(\mathfrak{g}(F))$ est le sous-espace de $SI(\mathfrak{g}(F))$ annul\'e par les applications $f\mapsto f_{M}$ pour tout Levi propre. C'est donc l'image du sous-espace des $f\in C_{c}^{\infty}(\mathfrak{g}(F))$ telles que $S^{G}(X,f)=0$ pour tout $X$ r\'egulier dans une sous-alg\`ebre de Levi propre. Ce sous-espace contient \'evidemment $C_{cusp}^{\infty}(\mathfrak{g}(F))$ mais ne lui est pas \'egal. En fait, l'assertion (2) montre que $SI_{cusp}(\mathfrak{g}(F))$ est bien l'image de $C_{cusp}^{\infty}(\mathfrak{g}(F))$. On reviendra sur ce point en 4.15.

  Soient maintenant $G$ et $G'$ deux groupes en situation d'endoscopie non standard, cf. [W1] 1.7. Rappelons que $G$ et $G'$ sont quasi-d\'eploy\'es et simplement connexes et qu'il y a une application de transfert entre $C_{c}^{\infty}(\mathfrak{g}(F))$ et $C_{c}^{\infty}(\mathfrak{g}'(F))$ (avec facteur de transfert \'egal \`a $1$ sur les couples qui se correspondent). On a
  
   (5) l'homomorphisme de transfert d\'efinit des isomorphismes 
 $$SI(\mathfrak{g}(F))\simeq SI(\mathfrak{g}'(F)),$$
   $$SI_{cusp}(\mathfrak{g}(F))\simeq SI_{cusp}(\mathfrak{g}'(F)),$$
   qui commutent \`a la transformation de Fourier.
 
    \bigskip
   
   \subsection{Image du transfert}
   On fixe un ensemble de repr\'esentants ${\cal E}(\tilde{G},{\bf a})$ de repr\'esentants des classes d'\'equivalence de donn\'ees endoscopiques elliptiques et relevantes de $(G,\tilde{G},{\bf a})$. On l'\'etend en un ensemble des repr\'esentants ${\cal E}_{+}(\tilde{G},{\bf a})$ de repr\'esentants des classes d'\'equivalence de couples $(\tilde{M},{\bf M}')$ o\`u $\tilde{M}$ est un espace de Levi de $\tilde{G}$ et ${\bf M}'$ est une donn\'ee endoscopique elliptique et relevante pour $(\tilde{M},{\bf a}_{M})$. On note $I^{{\cal E}}_{+}(\tilde{G}(F),\omega)$ le sous-espace des \'el\'ements $({\bf f}_{(\tilde{M},{\bf M}')})\in\oplus_{(\tilde{M},{\bf M}')\in {\cal E}_{+}(\tilde{G},{\bf a})}SI({\bf M}')\otimes Mes(M'(F))$  qui v\'erifient les conditions suivantes:
   
   (1) pour tout $(\tilde{M},{\bf M}')\in {\cal E}_{+}(\tilde{G},{\bf a})$, ${\bf f}_{(\tilde{M},{\bf M}')}$ est invariant par $Aut(\tilde{M},{\bf M}')$;
   
   (2) soit ${\bf G}'\in {\cal E}(\tilde{G},{\bf a})$ et  $M'$ un Levi de $G'$ qui est relevant; soit $(\tilde{M},{\bf M}')$ l'\'el\'ement de ${\cal E}_{+}(\tilde{G},{\bf a})$ qui lui est associ\'e par la construction de 3.4; alors $({\bf f}_{{\bf G}'})_{\tilde{M}'}={\bf f}_{(\tilde{M},{\bf M}')}$;
   
   (3) soit ${\bf G}'\in {\cal E}(\tilde{G},{\bf a})$ et  $M'$ un Levi de $G'$ qui n'est pas relevant; alors $({\bf f}_{{\bf G}'})_{\tilde{M}'}=0$.
   
   D'apr\`es (2) et 3.3(3), la projection naturelle de $I^{\cal E}_{+}(\tilde{G}(F),\omega)$ dans $\oplus_{{\bf G}'\in {\cal E}(\tilde{G},{\bf a})}SI({\bf G}')\otimes Mes(G'(F))$ est injective. On note $I^{\cal E}(\tilde{G}(F),\omega)$ l'image de cette projection.

 Dans le cas o\`u $F={\mathbb R}$, on travaille avec un $K$-espace tordu $K\tilde{G}$. Les espaces $I(\tilde{G}(F),\omega)$ et $I_{cusp}(\tilde{G}(F),\omega)$ ont des analogues \'evidents $I(K\tilde{G}({\mathbb R}),\omega)$ et $I_{cusp}(K\tilde{G}({\mathbb R}),\omega)$.   Il est peut-\^etre judicieux de noter $I_{+}^{\cal E}(K\tilde{G}({\mathbb R}),\omega)$ et  $I^{\cal E}(K\tilde{G}({\mathbb R}),\omega)$ les espaces $I_{+}^{\cal E}(\tilde{G}({\mathbb R}),\omega)$ et $I^{\cal E}(\tilde{G}({\mathbb R}),\omega)$, bien que leurs d\'efinitions ne fassent pas r\'ef\'erence au $K$-espace.
        
   \ass{Proposition}{(i) Supposons $F$ non archim\'edien. Alors l'application de transfert
   $$I(\tilde{G}(F),\omega)\otimes Mes(G(F))\to \oplus_{{\bf G}'\in {\cal E}(\tilde{G},{\bf a})}SI({\bf G}')\otimes Mes(G'(F))$$
  est injective et a pour image l'espace $I^{\cal E}(\tilde{G}(F),\omega)$. L'image de $I_{cusp}(\tilde{G}(F),\omega)\otimes Mes(G(F))$ est 
  $$ \oplus_{{\bf G}'\in {\cal E}(\tilde{G},{\bf a})}SI_{cusp}({\bf G}')^{Aut({\bf G}')}\otimes Mes(G'(F)).$$

   (ii) Supposons $F={\mathbb R}$. L'assertion devient vraie si on remplace $I(\tilde{G}(F),\omega)$, $I_{cusp}(\tilde{G}(F),\omega)$ et $I^{\cal E}(\tilde{G}(F),\omega)$ par $I(K\tilde{G}({\mathbb R}),\omega) $,   $I_{cusp}(K\tilde{G}({\mathbb R}),\omega)$ et $I^{\cal E}(K\tilde{G}({\mathbb R}),\omega)$.}
   
        La preuve occupe les  paragraphes 4.12 et 4.13. Remarquons que l'on peut d\'efinir une application de transfert
        $$(4)\qquad \begin{array}{ccc}I(\tilde{G}(F),\omega)\otimes Mes(G(F))&\to&\sum_{(\tilde{M},{\bf M}')\in {\cal E}_{+}(\tilde{G},{\bf a})}SI({\bf M}')\otimes Mes(M'(F))\\ {\bf f}&\mapsto&({\bf f}_{(\tilde{M},{\bf M}')})_{(\tilde{M},{\bf M}')\in {\cal E}_{+}(\tilde{G},{\bf a})}\\ \end{array}$$
 o\`u ${\bf f}_{(\tilde{M},{\bf M}')}$ est le transfert \`a ${\bf M}'$ de ${\bf f}_{\tilde{M},\omega}\in I(\tilde{M}(F),\omega)\otimes Mes(M(F))$ (on peut \'evidemment remplacer les $I(\tilde{G}(F),\omega)$ etc... par des $I(K\tilde{G}({\mathbb R}),\omega)$ etc... dans le cas r\'eel). L'application du (ii) de l'\'enonc\'e  est la compos\'ee de cette application et d'une projection naturelle. Or il est clair par construction et d'apr\`es 2.6 que l'image de l'application (4) est contenue dans l'espace $I^{{\cal E}}_{+}(\tilde{G}(F),\omega)$. Donc l'application de transfert de l'\'enonc\'e prend ses valeurs dans $I^{{\cal E}}(\tilde{G}(F),\omega)$. D'autre part, la premi\`ere assertion de l'\'enonc\'e \'equivaut \`a dire que l'image de l'application (4) est $I^{{\cal E}}_{+}(\tilde{G}(F),\omega)$.
        
        Dans les deux paragraphes suivants, on suppose fix\'ees des mesures de Haar sur tous les groupes intervenant, ce qui nous d\'ebarrasse des espaces de mesures.
        
    \bigskip
   \subsection{Preuve de la proposition 4.11 dans le cas non-archim\'edien}
   On a d\'efini en 4.2 la filtration $({\cal F}^nI(\tilde{G}(F),\omega))_{n\in {\mathbb N}}$. 
     Notons ${\cal F}^{n}I^{\cal E}_{+}(\tilde{G}(F),\omega)$ le sous-espace des \'el\'ements $(f_{(\tilde{L},{\bf L}')})\in I^{{\cal E}}_{+}(\tilde{G}(F),\omega)$ tels que $f_{(\tilde{L},{\bf L}')}=0$  pour tout espace de Levi $\tilde{L}$ tel  que $a_{\tilde{L}}>n$. Ces sous-espaces forment une filtration de $I^{{\cal E}}_{+}(\tilde{G}(F),\omega)$.    Notons $GrI(\tilde{G}(F),\omega)$ et $GrI^{\cal E}_{+}(\tilde{G}(F),\omega)$  les gradu\'es associ\'es \`a ces filtrations. Fixons un ensemble de repr\'esentants $\underline{\cal L}$ des classes de conjugaison par $G(F)$ d'espaces de Levi de $\tilde{G}$. D'apr\`es le lemme 4.2, on a l'isomorphisme
  $$(1)\qquad GrI(\tilde{G}(F),\omega)\simeq \oplus_{\tilde{M}\in \underline{\cal L}}I_{cusp}(\tilde{M}(F),\omega)^{W(\tilde{M})}.$$
  On a d'autre part une inclusion naturelle
  $$(2) \qquad GrI^{\cal E}_{+}(\tilde{G}(F),\omega)\subset  \oplus_{\tilde{M}\in \underline{\cal L}}(\oplus_{{\bf M}'\in {\cal E}(\tilde{M},{\bf a})}SI_{cusp}({\bf M}')^{Aut({\bf M}')})^{W(\tilde{M})}$$
  $$= \oplus_{(\tilde{M},{\bf M}')\in{\cal E}_{+}(\tilde{G},{\bf a})}  SI_{cusp}({\bf M}')^{Aut(\tilde{M},{\bf M}')}.$$
L' application  de transfert (4) de 4.11 est compatible aux filtrations et l'application qui en r\'esulte entre les gradu\'es n'est autre que la somme des applications naturelles de transfert.   Supposons  prouv\'e que le transfert induit un isomorphisme  
$$(3) \qquad I_{cusp}(\tilde{G}(F),\omega)\simeq \oplus_{{\bf G}'\in {\cal E}(\tilde{G},{\bf a})}SI_{cusp}({\bf G}')^{Aut({\bf G}')}.$$
On a alors un isomorphisme analogue
$$I_{cusp}(\tilde{M}(F),\omega)\simeq \oplus_{{\bf M}'\in {\cal E}(\tilde{M},{\bf a})}SI_{cusp}({\bf M}')^{Aut({\bf M}')}.$$
 pour chaque $\tilde{M}\in \underline{\cal L}$. Le transfert est compatible aux actions de $W(\tilde{M})$, on peut donc remplacer les deux membres ci-dessus par leurs sous-espaces d'invariants par $W(\tilde{M})$. On voit alors que l'inclusion (2) est elle-aussi une \'egalit\'e et que l'application gradu\'ee  $GrI(\tilde{G}(F),\omega)\to GrI^{\cal E}_{+}(\tilde{G}(F),\omega)$ est un isomorphisme. Le (i) de la proposition 4.11 en r\'esulte.

 Il faut montrer que (3) est un isomorphisme. Rappelons d'abord une propri\'et\'e fondamentale. Notons $C_{ell}^{\infty}(\tilde{G}(F))$ le sous-espace de $C_{c}^{\infty}(\tilde{G}(F))$ form\'e des \'el\'ements \`a support elliptique fortement r\'egulier et notons $I_{ell}(\tilde{G}(F),\omega)$ son image dans $I(\tilde{G}(F))$. On d\'efinit de fa\c{c}on similaire des espaces $SI_{\tilde{G}-ell}({\bf G}')$ en rempla\c{c}ant la condition fortement r\'egulier par fortement $\tilde{G}$-r\'egulier. Alors 
 
 (4) le transfert d\'efinit un isomorphisme
 $$I_{ell}(\tilde{G}(F),\omega)\simeq \oplus_{{\bf G}'\in {\cal E}(\tilde{G},{\bf a})}SI_{\tilde{G}-ell}({\bf G}')^{Aut({\bf  G}')}.$$
 
 Cela r\'esulte des faits suivants. D'abord,  les $\eta\in \tilde{G}(F)$ elliptiques et fortement r\'eguliers pour lesquels  $\omega$ est non trivial sur $Z_{G}(\eta;F)$ ne comptent pas: du c\^ot\'e de $I_{ell}(\tilde{G}(F),\omega)$, les int\'egrales orbitales sont toutes nulles au voisinage d'un tel point; et il ne leur correspond rien du c\^ot\'e droit de la formule ci-dessus. Fixons $\eta\in \tilde{G}(F)_{ell}$ tel que $\omega$ soit trivial sur $Z_{G}(\eta;F)$. Pour chaque ${\bf G}'\in {\cal E}(\tilde{G},{\bf a})$, on fixe des donn\'ees auxiliaires $G'_{1},...,\Delta_{1}$.   Soit $f\in C_{ell}^{\infty}(\tilde{G}(F))$. Alors les familles $(I^{\tilde{G}}(\eta',\omega,f))_{\eta'\in \dot{\cal X}(\eta)}$ et $(S^{\tilde{G}'_{1}}(\epsilon_{1},f^{G'_{1}}))_{({\bf G}',\epsilon)\in \dot{\cal X}^{\cal E}(\eta)}$ se d\'eduisent l'une de l'autre par les transformations bijectives (4) et (5) de 4.9.
 
  L'application (3) est injective: si $f\in I_{cusp}(\tilde{G}(F),\omega)$ a un transfert nul, il r\'esulte de (4) (ou plus exactement de sa preuve) que $I^{\tilde{G}}(\gamma,\omega,f)=0$ pour tout \'el\'ement fortement r\'egulier et elliptique; la cuspidalit\'e de $f$ entra\^{\i}ne alors $f=0$.

 La preuve de la surjectivit\'e n\'ecessite quelques pr\'eparatifs. Fixons $\eta\in \tilde{G}_{ss}(F)_{ell}$ et une forme int\'erieure quasi-d\'eploy\'ee $\bar{G}$ de $G_{\eta}$. On fixe un voisinage $\bar{\mathfrak{u}}$ de $0$ dans $\bar{\mathfrak{g}}(F)$ v\'erifiant les conditions de 4.6. et on utilise les constructions de ce paragraphe.   La descente d'Harish-Chandra nous fournit une application
  $$(5) \qquad \begin{array}{ccc}I_{cusp}(\tilde{U},\omega)&\to& \oplus_{y\in \dot{\cal Y}(\eta)}I_{cusp}(\mathfrak{u}_{\eta[y]},\omega)^{Z_{G}(\eta[y],F)}\\ f&\mapsto&(f_{y})_{y\in\dot{\cal Y}(\eta)}\\ \end{array}.$$
  Son image est form\'ee des familles $(f_{y})_{y\in \dot{\cal Y}(\eta)}$ telles que $f_{y}=f_{y'}$ si $\eta[y]=\eta[y']$. Fixons une transformation de Fourier sur $C_{c}^{\infty}(\bar{\mathfrak{g}}(F))$, dont on d\'eduit de telles transformations dans chaque $C_{c}^{\infty}(\mathfrak{g}_{\eta[y]}(F))$. On v\'erifie que ces transformations sont les m\^emes dans le cas o\`u $\eta[y]=\eta[y']$.

Pour tout $({\bf G}',\epsilon)\in \dot{\cal X}^{{\cal E}}(\eta)$, on fixe un diagramme joignant $\epsilon$ \`a un \'el\'ement $\eta[y]$ (on peut d'ailleurs supposer $\eta[y]=\eta$ mais peu  importe). On utilise les constructions de 4.9 pour ce diagramme, en les affectant au besoin d'indices $\epsilon$. C'est-\`a-dire que l'on introduit la donn\'ee endoscopique $\bar{{\bf G}}'_{\epsilon}=(\bar{G}'_{\epsilon},\bar{{\cal G}}'_{\epsilon},\bar{s}_{\epsilon})$ de $\bar{G}_{SC}$. Les isomorphismes d\'ecrits en 4.9 fournissent une correspondance entre classes de conjugaison stable semi-simples dans $\mathfrak{g}'_{\epsilon}(F)$ et dans $\bar{\mathfrak{g}}(F)$. On note $\mathfrak{u}'_{\epsilon}$ l'ensemble des \'el\'ements de $\mathfrak{g}'_{\epsilon}(F)$ dont la partie semi-simple a une classe de conjugaison stable qui correspond \`a celle d'un \'el\'ement de $\bar{\mathfrak{u}}$. En scindant la projection $\mathfrak{g}'_{1,\epsilon_{1}}(F)\to \mathfrak{g}'_{\epsilon}(F)$ comme en 4.8, on identifie $\mathfrak{u}'_{\epsilon}$ \`a un sous-ensemble de $\mathfrak{g}'_{1,\epsilon_{1}}(F)$. On note $\tilde{U}''_{1,\epsilon_{1}}$ l'ensemble des \'el\'ements de $\tilde{G}'_{1}(F)$ dont la partie semi-simple est stablement conjugu\'ee \`a un \'el\'ement de $C_{1}(F)exp(\mathfrak{u}_{\epsilon})\epsilon_{1}$. Rappelons qu'un \'el\'ement  de $Aut({\bf G}')$ est d\'efini par un \'el\'ement $x\in \hat{G}$, lequel d\'etermine un automorphisme $\tilde{\alpha}_{x}$ de $\tilde{G}'$. On note $\tilde{U}'_{1,\epsilon_{1}}$ la r\'eunion des $\tilde{\alpha}_{x}(\tilde{U}''_{1,\epsilon_{1}})$ pour tous les $x\in Aut({\bf G}')$. Puisque les fonctions que l'on consid\`ere sur $\tilde{G}'_{1}(F)$ se transforment selon le caract\`ere $\lambda_{1}$ de $C_{1}(F)$, la descente d\'efinit    une application 
$$ SI_{\lambda_{1},cusp}(\tilde{U}'_{1,\epsilon_{1}})\to SI_{cusp}(\mathfrak{u}'_{\epsilon}).$$
D'apr\`es 4.8, son image est le sous-espace des \'el\'ements de $SI_{cusp}(\mathfrak{u}'_{\epsilon})$ qui se transforment selon un certain carat\`ere de $\Xi_{\epsilon}^{\Gamma_{F}}$, o\`u $\Xi_{\epsilon}=Z_{G'}(\epsilon)/G'_{\epsilon}$.
  
L'espace  $SI_{\lambda_{1},cusp}(\tilde{U}'_{1,\epsilon_{1}})$ est stable par l'action de $Aut({\bf G}')$. Nous voulons d\'eterminer l'image de l'application
$$(6) \qquad   SI_{\lambda_{1},cusp}(\tilde{U}'_{1,\epsilon_{1}})^{Aut({\bf G}')}\to SI_{cusp}(\mathfrak{u}'_{\epsilon}).$$
Pour $x\in Aut({\bf G}')$,  l'action de $x$ n'impose une condition au voisinage de $\epsilon_{1}$ que si $\tilde{\alpha}_{x}(\epsilon)$ et $\epsilon$ sont stablement conjugu\'es. S'il en est ainsi, un \'el\'ement $g'\in G'$ qui \'etablit cette conjugaison stable d\'efinit un torseur int\'erieur entre les commutants connexes de ces \'el\'ements. Or on a suppos\'e ces groupes quasi-d\'eploy\'es. Quitte \`a modifier $g'$, on peut  donc supposer que ce torseur int\'erieur est un isomorphisme d\'efini sur $F$. Cela conduit \`a introduire l'ensemble $Aut_{\epsilon}$ des couples $(g', x)$ o\`u $x$ est comme ci-dessus et $g'\in G'$ est tel que $g'\tilde{\alpha}_{x}(\epsilon)(g')^{-1}=\epsilon$ et  que l'automorphisme $ad_{g'}\circ \alpha_{x}$ de $G'_{\epsilon}$ soit  d\'efini sur $F$. Soit $(g',x)\in Aut_{\epsilon}$. Consid\'erons les couples $(Y',Y)\in \mathfrak{u}'_{\epsilon}\times \mathfrak{u}'_{\epsilon}$ d'\'el\'ements  tels que  $Y'=ad_{g'}\circ\alpha_{x}(Y)$, avec $Y$ en position g\'en\'erale. D'apr\`es la construction de 2.6, il existe une fonction $(Y',Y)\mapsto\Lambda_{g',x}(Y',Y)$ sur cet ensemble de couples telle que pour $f'_{1}\in SI_{\lambda_{1},cusp}(\tilde{U}'_{1,\epsilon_{1}})$, la condition que $f'_{1}$ soit invariante par l'automorphisme d\'etermin\'e par $x$ se traduise par l'\'egalit\'e $S^{G'_{1}}(exp(Y')\epsilon_{1},f'_{1})=\Lambda_{g',x}(Y',Y)S^{G'_{1}}(exp(Y)\epsilon_{1},f'_{1})$ pour tout tel couple. En fait, la fonction $\Lambda_{g',x}$ est la restriction d'une fonction qui se transforme selon un caract\`ere du groupe $G'_{1}(F)\times G'_{1}(F)$. Pour $\mathfrak{u}'_{\epsilon}$ assez petit, elle est donc constante, de valeur disons $\Lambda(g',x)$. Par descente, la condition pr\'ec\'edente se traduit pour $f'\in SI_{cusp}(\mathfrak{u}'_{\epsilon})$ par l'\'egalit\'e $S^{G'_{\epsilon}}(ad_{g'}\circ\alpha_{x}(Y),f')=\Lambda(g',x)S^{G'_{\epsilon}}(Y,f')$ pour tout $Y\in \mathfrak{u}'_{\epsilon}$.  Notons que, dans le cas  $x=1$, $g'$ d\'efinit un \'el\'ement de $\Xi_{\epsilon}^{\Gamma_{F}}$ et cette \'egalit\'e n'est autre que la condition de transformation d\'ej\`a introduite sous l'action de ce groupe. La formule  pr\'ec\'edente d\'efinit une action du groupe $Aut_{\epsilon}$ sur $SI_{cusp}(\mathfrak{u}'_{\epsilon})$. On obtient

(7) l'image de l'application (6) est \'egale \`a $SI_{cusp}(\mathfrak{u}'_{\epsilon})^{Aut_{\epsilon}}$, l'invariance \'etant bien s\^ur relative l'action d\'efinie ci-dessus. 

Comme on l'a dit en 4.10, de la transformation de Fourier fix\'ee sur $C_{c}^{\infty}(\bar{\mathfrak{g}}(F))$ se d\'eduit une transformation de Fourier sur $C_{c}^{\infty}(\mathfrak{g}'_{\epsilon}(F))$. On peut supposer la premi\`ere invariante par toute action d'un \'el\'ement de $G$. La seconde l'est alors par l'action de $Aut_{\epsilon}$. Il en r\'esulte que 

(8) $SI_{cusp}(\mathfrak{g}'_{\epsilon}(F))^{Aut_{\epsilon}}$ est invariante par transformation de Fourier.

Pour $y\in \dot{\cal Y}(\eta)$ et $f_{y}\in C_{c}^{\infty}(\mathfrak{u}_{\eta[y]})$, nous allons construire une fonction $\varphi_{\epsilon,y}\in C_{c}^{\infty}(\mathfrak{u}'_{\epsilon})$. Par lin\'earit\'e, on peut supposer que $f_{y}=f_{y,Z}\otimes f_{y,sc}$, avec $f_{y,Z}\in C_{c}^{\infty}(\mathfrak{z}_{G_{\eta[y]}}(F))$ et $f_{y,sc}\in C_{c}^{\infty}(\mathfrak{g}_{\eta[y],SC}(F))$. Les centres $Z(\bar{G})$ et $Z(G_{\eta[y]})$ s'identifient. On peut donc identifier $f_{y,Z}$ \`a une fonction sur $\mathfrak{z}_{\bar{G}}(F)$. La donn\'ee ${\bf \bar{G}}'_{\epsilon}$ est aussi une donn\'ee endoscopique de $G_{\eta[y],SC}$ donc $f_{y,sc}$ se transf\`ere en une fonction disons $\phi_{y}$ sur $\mathfrak{\bar{g}}'_{\epsilon}(F)$. Par lin\'earit\'e, on peut supposer  $\phi_{y}=\phi_{y,Z}\otimes \phi_{y,sc}$, avec $\phi_{y,Z}\in C_{c}^{\infty}(\mathfrak{z}_{\bar{g}'_{\epsilon}}(F))$ et $\phi_{y,sc}\in C_{c}^{\infty}(\mathfrak{\bar{g}}'_{\epsilon,SC}(F))$. Par endoscopie non standard, $\phi_{y,sc}$ se transf\`ere en une fonction $\varphi_{\epsilon,y,sc}\in C_{c}^{\infty}(\mathfrak{g}'_{\epsilon,SC}(F))$. Par les isomorphismes de 3.7, on a l'identification $\mathfrak{z}_{G'_{\epsilon}}(F)=\mathfrak{z}_{\bar{G}}(F)\oplus \mathfrak{z}_{\bar{G}'_{\epsilon}}(F)$. La fonction $f_{y,Z}\otimes \phi_{y,Z}$ s'identifie \`a une fonction $\varphi_{\epsilon,y,Z}$ sur $\mathfrak{z}_{G'_{\epsilon}}(F)$. On pose $\varphi_{\epsilon,y}=\varphi_{\epsilon,y,Z}\otimes \varphi_{\epsilon,y,sc}$. Il est (plus ou moins) clair que l'on peut effectuer les choix de sorte que cette fonction soit \`a support dans $\mathfrak{u}'_{\epsilon}$. L'utilit\'e de cette construction est l'existence d'une famille $(c_{\epsilon,y})_{y\in \dot{\cal Y}(\eta)}$ de nombres complexes non nuls telle que la propri\'et\'e suivante soit v\'erifi\'ee. Soit $f\in I_{cusp}(\tilde{G}(F),\omega)$ (la condition de cuspidalit\'e ne sert ici \`a rien mais peu importe). Soit $(f_{y})_{y\in \dot{\cal Y}(\eta)}$ son image par l'application (5). Soit $\varphi\in SI_{\lambda_{1},cusp}(G'_{1}(F))$ le transfert de $f$ et soit $\varphi_{\epsilon}$ la fonction sur $\mathfrak{u}'_{\epsilon}$ qui se d\'eduit de $\varphi$ par (6). Alors 

(9) on a l'\'egalit\'e suivante dans $SI_{cusp}(\mathfrak{u}'_{\epsilon})$:
$$\varphi_{\epsilon}=\sum_{y\in \dot{\cal Y}(\eta)}c_{\epsilon,y}\varphi_{\epsilon,y}.$$

 Cela r\'esulte de la preuve de [W1] 3.11 (bien s\^ur, cela suppose que le voisinage $\bar{\mathfrak{u}}$ est assez petit).
 
 Prouvons maintenant la surjectivit\'e de (3). Le lemme 4.9 et un argument de partition de l'unit\'e sur l'espace $\tilde{G}_{ss}(F)_{ell}/st-conj$ montrent que, pour prouver cette surjectivit\'e, il suffit de prouver l'assertion suivante. Soient  $(f_{{\bf G}'})_{{\bf G}'\in {\cal E}(\tilde{G},{\bf a})}\in \oplus_{{\bf G}'\in {\cal E}(\tilde{G},{\bf a})}SI_{cusp}({\bf G}')^{Aut({\bf G}')}$  et $\eta\in \tilde{G}_{ss}(F)_{ell}$.  Alors il existe $f\in I_{cusp}(\tilde{G}(F),\omega)$ telle que pour tout  $({\bf G}',\epsilon)\in \dot{\cal X}^{\cal E}(\eta)$, les int\'egrales orbitales stables de $f_{{\bf G}'}$ et du  transfert $f^{{\bf G}'}$ de $f$ (ces fonctions \'etant identifi\'ees \`a des fonctions sur $\tilde{G}'_{1}(F)$) co\"{\i}ncident dans un voisinage de $\epsilon_{1}$. On fixe $\eta$ et on utilise les constructions ci-dessus. D'apr\`es les propri\'et\'es de l'application de descente (6), on peut aussi bien prouver l'assertion suivante. Soient $({\bf G}',\epsilon)\in \dot{\cal X}^{\cal E}(\eta)$ et $\phi\in SI_{cusp}(\mathfrak{u}'_{\epsilon})^{Aut_{\epsilon}}$. Alors il existe $f\in C_{c}^{\infty}(\tilde{U})$ dont les transferts $f^{\underline{{\bf G}}'}$ v\'erifient les deux conditions:
 
 (10) l'image de $f^{{\bf G}'}$ par descente au voisinage de $\epsilon_{1}$  a les m\^emes int\'egrales orbitales stables que $\phi$ dans un voisinage de $0$;
 
 (11) pour  $(\underline{{\bf G}}',\underline{\epsilon})\in \dot{\cal X}^{\cal E}(\eta)$ diff\'erent de $({\bf G}',\epsilon)\in \dot{\cal X}^{\cal E}(\eta)$, l'image de $f^{\underline{{\bf G}}'}$ par descente au voisinage de $\underline{\epsilon}_{1}$ a des int\'egrales orbitales stables nulles dans un voisinage de $0$.
 
 D'apr\`es 4.1(2), on peut trouver $\phi'\in SI_{cusp}(\mathfrak{u}'_{\epsilon})$ \`a support r\'egulier elliptique et  tel que $\phi$ et $\hat{\phi}'$ aient m\^emes int\'egrales orbitales stables au voisinage de $0$.  La propri\'et\'e (8) nous permet de supposer que $\phi'$ est invariante par le groupe $Aut_{\epsilon}$. On peut relever $\phi'$ en un \'el\'ement $\varphi'\in SI_{\lambda_{1},cusp}(\tilde{U}'_{1,\epsilon_{1}})^{Aut({\bf G}')}$ \`a support r\'egulier elliptique, et compl\'eter $\varphi'$ en un \'el\'ement de $\oplus_{\underline{\bf G}'\in {\cal E}(\tilde{G},{\bf a})}SI_{cusp}(\underline{\bf G}')^{Aut(\underline{\bf G}')}$, nul sur les autres composantes. D'apr\`es (4), c'est le transfert d'un \'el\'ement  $f'\in I_{ell}(\tilde{G}(F),\omega)$. Il est clair que l'on peut supposer $f'\in I_{cusp}(\tilde{U},\omega)$. Appliquons \`a $f'$ les constructions pr\'ec\'edant la formule (9), en les affectant d'un $'$.  On obtient les deux propri\'et\'es suivantes:
 
 - la fonction $\phi'$ a les m\^emes int\'egrales orbitales stables que  $ \sum_{y\in \dot{\cal Y}(\eta)}c_{\epsilon,y}\varphi'_{\epsilon,y}$; 
    
- pour  $(\underline{\bf G}',\underline{\epsilon})\in \dot{\cal Y}^{\cal E}(\xi)$ diff\'erent de $({\bf G}',\epsilon)$, la fonction $ \sum_{y\in \dot{\cal Y}(\eta)}c_{\underline{\epsilon},y}\varphi'_{\underline{\epsilon},y}$ a des int\'egrales orbitales stables nulles.

Pour tout $y\in \dot{\cal Y}(\eta)$, notons $f_{y}$ la fonction $\gamma(\mathfrak{g}'_{\epsilon})^{-1}\gamma(\mathfrak{g}_{\eta[y]})\hat{f}'_{y}$ restreinte \`a $\mathfrak{u}_{\eta[y]}$. D'apr\`es la description de l'image de (5), il existe $f\in I_{cusp}(\tilde{U},\omega)$ dont l'image par descente soit $(f_{y})_{y\in \dot{\cal Y}(\eta)}$. Soit $\varphi\in \oplus_{\underline{\bf G}'\in {\cal E}(\tilde{G},{\bf a})}SI_{cusp}(\underline{\bf G}')^{Aut(\underline{\bf G}')}$ le transfert de $f$. On applique \`a $f$ les constructions pr\'ec\'edant la formule (9). D'apr\`es 4.10, toutes les fonctions issues de $f$ se d\'eduisent de celles issues de $f'$ par transformation de Fourier et \'eventuellement multiplication par des constantes $\gamma$. On obtient que, pour $(\underline{\bf G}',\underline{\epsilon})\in \dot{\cal Y}^{\cal E}(\eta)$, l'image par descente de $f^{\underline{{\bf G}}'}$ a les m\^emes int\'egrales orbitales  stables que 
$$ \sum_{y\in \dot{\cal Y}(\eta)}c_{\underline{\epsilon},y}\varphi_{\underline{\epsilon},y},$$
 ou encore que 
 $$\gamma(\underline{\mathfrak{g}}'_{\underline{\epsilon}}) \gamma(\mathfrak{g}'_{\epsilon})^{-1}\sum_{y\in \dot{\cal Y}(\eta)}c_{\underline{\epsilon},y}\hat{\varphi}'_{\underline{\epsilon},y},$$
 ou encore que $\hat{\phi}'$ si $(\underline{\bf G}',\underline{\epsilon})=({\bf G}',\epsilon)$, $0$ sinon. D'apr\`es le choix de $\phi'$, $f$ satisfait (10) et (11), ce qui ach\`eve la d\'emonstration. 

\bigskip 
  \subsection{Preuve de la proposition  4.11 dans le cas r\'eel}
 On reprend la preuve du cas non-archim\'edien. Son d\'ebut reste pertinent.   En adaptant les notations aux $K$-espaces tordus,  il faut  prouver que le transfert induit un isomorphisme
 $$(1) \qquad I_{cusp}(K\tilde{G}({\mathbb R}),\omega)\simeq \oplus_{{\bf G}'\in {\cal E}(K\tilde{G},{\bf a})}SI_{cusp}({\bf G}')^{Aut({\bf G}')}.$$
  Commen\c{c}ons par d\'ecrire l'espace $I_{cusp}(\tilde{G}({\mathbb R}),\omega)$.  On a d\'efini en 1.3 la notion de tore tordu maximal elliptique dans $\tilde{G}$. Notons que, pour un tore tordu maximal $\tilde{T}$, la condition d'ellipticit\'e revient  \`a dire que  $ (T^{\theta,0}/A_{\tilde{G}})({\mathbb R})$ est compact. Il  y a au plus un nombre fini de classes de conjugaison par $G({\mathbb R})$ de tores tordus maximaux elliptiques (j'ignore s'il y en a au plus un comme dans le cas non tordu).  Fixons un ensemble de repr\'esentants $\tilde{{\cal T}}_{ell}$ des classes de conjugaison par $G({\mathbb R})$ parmi les tores tordus maximaux elliptiques $\tilde{T}$ tels que $\omega$ soit trivial sur $T^{\theta}({\mathbb R})$. Cet ensemble peut \^etre vide.
   Consid\'erons l'application qui \`a $f\in I_{cusp}(\tilde{G}({\mathbb R}),\omega)$ associe la  famille de fonctions $(\varphi_{\tilde{T}})_{\tilde{T}\in \tilde{{\cal T}}_{ell}}$, o\`u $\varphi_{\tilde{T}}$ est la fonction d\'efinie sur les \'el\'ements fortement r\'eguliers de $\tilde{T}({\mathbb R})$ par
  $$\varphi_{\tilde{T}}(\gamma)=I^{\tilde{G}}(\gamma,\omega,f).$$
  Elle est injective. Une famille  $(\varphi_{\tilde{T}})_{\tilde{T}\in \tilde{{\cal T}}_{ell}}$ dans l'image  v\'erifie la condition
  
  (2)  pour tout $\tilde{T}\in \tilde{{\cal T}}_{ell}$, tout \'el\'ement fortement r\'egulier $\gamma\in \tilde{T}({\mathbb R})$ et tout $g\in G({\mathbb R})$ tel que $g\gamma g^{-1}\in \tilde{T}({\mathbb R})$, on a $\varphi_{\tilde{T}}(g \gamma g^{-1})=\omega(g)\varphi_{\tilde{T}}(\gamma)$.
  
  Par descente d'Harish-Chandra,   nos fonctions v\'erifient localement les conditions de r\'egularit\'e ou de saut habituelles dans cette th\'eorie. Mais, parce que l'on consid\`ere ici des fonctions cuspidales, ces conditions se simplifient grandement. Soient $\tilde{T}\in \tilde{{\cal T}}_{ell}$ et  $\eta\in \tilde{T}({\mathbb R})$. Notons $\Sigma(T)_{\eta}$ l'ensemble des racines de $T^{\theta,0}$ dans $G_{\eta}$. Puisque $(T^{\theta,0}/A_{\tilde{G}})({\mathbb R})$ est compact, toutes ces racines sont imaginaires. Fixons un sous-ensemble de racines positives et d\'efinissons une fonction $\Delta_{\eta}$ sur le sous-ensemble des \'el\'ements de $\mathfrak{t}^{\theta}({\mathbb R})$ qui sont r\'eguliers dans $G_{\eta}$ par la formule 
  $$\Delta_{\eta}(X)=\prod_{\alpha\in \Sigma(T)_{\eta},\alpha>0}sgn(i\alpha(X)),$$
o\`u $sgn$ est le signe usuel d'un r\'eel non nul.  Cette fonction prend ses valeurs dans $\{\pm 1\}$. On a simplement
  
  (3) pour $\tilde{T}$ et $\eta$ comme ci-dessus, la fonction $X\mapsto \Delta_{\eta}(X)\varphi_{\tilde{T}}(exp(X)\eta)$ se prolonge en une fonction $C^{\infty}$  au voisinage de $0$ dans $\mathfrak{t}^{\theta}({\mathbb R})$.

 Inversement, la th\'eorie de la descente montre que toute famille $(\varphi_{\tilde{T}})_{\tilde{T}\in \tilde{{\cal T}}_{ell}}$ v\'erifiant (2) et (3) est l'image d'un \'el\'ement de $I_{cusp}(\tilde{G}({\mathbb R}),\omega)$. Ce r\'esultat se propage au $K$-espace $K\tilde{G}$. Pour $p\in \Pi$, on note plus pr\'ecis\'ement $\tilde{{\cal T}}_{ell,p}$ l'ensemble associ\'e \`a la composante $\tilde{G}_{p}$.   On pose $K\tilde{{\cal T}}_{ell}=\sqcup_{p\in \Pi}\tilde{{\cal T}}_{ell,p}$.  On obtient que l'application
  $$f\mapsto (\varphi_{\tilde{T}})_{\tilde{T}\in K\tilde{{\cal T}}_{ell}}$$
 est injective et que son image est form\'ee des familles v\'erifiant (2) et (3).  
 
 Soit ${\bf G}'\in {\cal E}(K\tilde{G},{\bf a})$. Fixons des donn\'ees suppl\'ementaires $G'_{1}$,...,$\Delta_{1}$ et identifions $C_{c}^{\infty}({\bf G}')$ \`a $C_{c,\lambda_{1}}^{\infty}(\tilde{G}'_{1}({\mathbb R}))$. Parce que $\tilde{G}'$ est \`a torsion int\'erieure, il y a au plus une  classe de conjugaison par $G'({\mathbb R})$ de tores tordus maximaux elliptiques dans $\tilde{G}'$. S'il n'y en a pas, il est clair que $SI_{cusp}({\bf G}')$ est nul. Supposons qu'il existe un tel tore tordu maximal elliptique, fixons-en un que l'on note $\tilde{T}'$. Notons $\tilde{T}'_{1}$ son image r\'eciproque dans $\tilde{G}'_{1}({\mathbb R})$. On consid\`ere l'application qui, \`a $f\in SI_{cusp}({\bf G}')$, associe la fonction $\varphi_{\tilde{T}'_{1}}$ sur $\tilde{T}'_{1}({\mathbb R})$ d\'efinie par $\varphi_{\tilde{T}'_{1}}(\delta_{1})=S^{\tilde{G}'_{1}}(\delta_{1},f)$ pour tout $\delta_{1}\in \tilde{T}'_{1}({\mathbb R})$ fortement r\'egulier. Cette application est injective. Un \'el\'ement de l'image v\'erifie les conditions
 
 (4) $\varphi_{\tilde{T}'_{1}}(c_{1}\delta_{1})=\lambda_{1}(c_{1})^{-1}\varphi_{\tilde{T}'_{1}}(\delta_{1})$ pour tout $\delta_{1}\in \tilde{T}'_{1}({\mathbb R})$ fortement r\'egulier et tout $c_{1}\in C_{1}({\mathbb R})$;
 
 (5) pour deux \'el\'ements $\delta_{1},\delta'_{1}\in \tilde{T}'_{1}({\mathbb R})$ fortement r\'eguliers et stablement conjugu\'es, $\varphi_{\tilde{T}'_{1}}(\delta'_{1})=\varphi_{\tilde{T}'_{1}}(\delta_{1})$.
 
 De nouveau, par descente, la fonction v\'erifie localement les conditions \'etablies par Shelstad. Puisqu'on travaille avec des fonctions cuspidales, ces conditions se simplifient. Soit $\epsilon\in \tilde{T}'({\mathbb R})$. On d\'efinit comme ci-dessus une fonction $\Delta_{\epsilon}$ sur  l'ensemble des \'el\'ements $\mathfrak{t}'({\mathbb R})$ qui sont r\'eguliers dans $G'_{\epsilon}$. On la remonte en une fonction d\'efinie presque partout sur $\mathfrak{t}'_{1}({\mathbb R})$. Alors  
 
 (6)  pour $\epsilon\in \tilde{T}'({\mathbb R})$ et $\epsilon_{1}\in \tilde{T}'_{1}({\mathbb R})$ au-dessus de $\epsilon$, la fonction $Y\mapsto \Delta_{\epsilon}(Y)\varphi_{\tilde{T}'_{1}}(exp(Y)\epsilon_{1})$ se prolonge en une fonction $C^{\infty}$ au voisinage de $0$ dans $\mathfrak{t}'_{1}({\mathbb R})$.
 
 Inversement, une fonction v\'erifiant les conditions (4), (5) et (6) est dans l'image de $SI_{cusp}({\bf G}')$, cf. [S1] th\'eor\`eme 12.1. On doit d\'eterminer l'image du sous-espace des invariants par $Aut({\bf G}')$. Notons $\tilde{T}'({\mathbb R})_{\sharp}$ l'ensemble des \'el\'ements $\delta\in \tilde{T}'({\mathbb R})$  tels que  $N^{\tilde{G}',K\tilde{G}}(\delta)$ appartient \`a l'image de $K\tilde{G}_{ab}({\mathbb R})$ par $N^{K\tilde{G}}$.  Cet ensemble est ouvert et ferm\'e (cela r\'esulte des d\'efinitions).   D'apr\`es le (iii) de la proposition 1.14, pour tout \'el\'ement $\tilde{G}$-r\'egulier $\delta\in \tilde{T}'({\mathbb R})_{\sharp}$, il existe $\gamma\in K\tilde{G}({\mathbb R})$ tel que $(\delta,\gamma)\in {\cal D}_{K\tilde{G}}$. Les d\'efinitions et le corollaire 2.6 entra\^{\i}nent que la condition d'invariance par $Aut({\bf G}')$ se traduit simplement par les deux conditions suivantes:
 
 (7) $\varphi_{\tilde{T}'_{1}}$ est nulle sur l'image r\'eciproque de $\tilde{T}'({\mathbb R})_{\sharp}$ dans $\tilde{T}'_{1}({\mathbb R})$;

 (8) pour deux \'el\'ements $\delta_{1},\delta'_{1}\in \tilde{T}'_{1}({\mathbb R})$ fortement r\'eguliers pour lesquels il existe $\gamma\in K\tilde{G}({\mathbb R})$ de sorte que $(\delta_{1},\gamma)$ et $(\delta'_{1},\gamma)$ appartiennent tous deux \`a ${\cal D}_{1,K\tilde{G}}$, on a l'\'egalit\'e
 $\Delta_{1}(\delta'_{1},\gamma)^{-1}\varphi_{\tilde{T}'_{1}}(\delta'_{1})=\Delta_{1}(\delta_{1},\gamma)^{-1}\varphi_{\tilde{T}'_{1}}(\delta_{1})$.
 
 Remarquons que cette condition implique (4) et (5).

 Quand on se limite \`a des fonctions \`a support r\'egulier elliptique, l'assertion 4.12(4) reste vraie sous la forme: le transfert d\'efinit un isomorphisme
 $$(9) \qquad I_{ell}(K\tilde{G}({\mathbb R}),\omega)\simeq \oplus_{{\bf G}'\in {\cal E}(K\tilde{G},{\bf a})}SI_{\tilde{G}-ell}({\bf G}')^{Aut({\bf  G}')}.$$
  Comme dans le cas non-archim\'edien, cela entra\^{\i}ne que le transfert est injectif sur $I_{cusp}(K\tilde{G}({\mathbb R}),\omega)$. 
 
 Notons ${\cal E}(K\tilde{G},{\bf a})_{0}$ l'ensemble des ${\bf G}'\in {\cal E}(K\tilde{G},{\bf a})$ tels que $\tilde{G}'$ poss\`ede un sous-tore tordu elliptique. Comme on l'a d\'ej\`a dit, il n'y  a qu'une classe de conjugaison de tels sous-tores et on en fixe un que l'on note $\tilde{T}[\tilde{G}']$. Consid\'erons une famille $(\varphi_{\tilde{T}[\tilde{G}']_{1}})_{{\bf G}'\in {\cal E}(K\tilde{G},{\bf a})_{0}}$, o\`u, pour tout ${\bf G}'\in {\cal E}(K\tilde{G},{\bf a})_{0}$, $\varphi_{\tilde{T}[\tilde{G}']_{1}}$ est une fonction sur $\tilde{T}[\tilde{G}']_{1}({\mathbb R})$ (d\'efinie presque partout) v\'erifiant (6), (7) et (8). Nous allons en d\'eduire une famille $(\varphi_{\tilde{T}})_{\tilde{T}\in K\tilde{{\cal T}}_{ell}}$ o\`u, pour tout $\tilde{T}\in K\tilde{{\cal T}}_{ell}$, $\varphi_{\tilde{T}}$ est une fonction d\'efinie presque partout sur $\tilde{T}({\mathbb R})$. Soient $\tilde{T}\in K\tilde{{\cal T}}_{ell}$ et $\gamma\in \tilde{T}({\mathbb R})\cap K\tilde{G}_{reg}(F)$. On peut supposer que chaque \'el\'ement de l'ensemble $\dot{\cal X}^{{\cal E}}(\gamma)$ de 4.9 est de la forme $({\bf G}',\delta)$ o\`u ${\bf G}'\in {\cal E}(K\tilde{G},{\bf a})_{0}$ et $\delta\in \tilde{T}[\tilde{G}']({\mathbb R})$. On pose alors
 $$\varphi_{\tilde{T}}(\gamma)= [T^{\theta}({\mathbb R}):T^{\theta,0}({\mathbb R})]\vert \dot{\cal X}(\gamma)\vert ^{-1}d(\theta^*)^{-1/2}\sum_{({\bf G}',\delta)\in \dot{\cal X}^{{\cal E}}(\gamma)}\Delta_{1}(\delta_{1},\gamma)^{-1}\varphi_{\tilde{T}[\tilde{G}']}(\delta_{1}),$$
 cf. 4.9(5). Dans le cas o\`u $(\varphi_{\tilde{T}[\tilde{G}']_{1}})_{{\bf G}'\in {\cal E}(K\tilde{G},{\bf a})_{0}}$ est \`a support r\'egulier, c'est-\`a-dire provient d'un \'el\'ement de $\oplus_{{\bf G}'\in {\cal E}(K\tilde{G},{\bf a})}SI_{\tilde{G}-ell}({\bf G}')^{Aut({\bf  G}')}$, la famille $(\varphi_{\tilde{T}})_{\tilde{T}\in K\tilde{{\cal T}}_{ell}}$ provient de l'\'el\'ement de $I_{ell}(K\tilde{G}({\mathbb R}),\omega)$ qui correspond \`a cet \'el\'ement par l'isomorphisme (9). Dans le cas g\'en\'eral, les \'el\'ements de la famille $(\varphi_{\tilde{T}})_{\tilde{T}\in K\tilde{{\cal T}}_{ell}}$ v\'erifient (2) par construction. Pour d\'emontrer la surjectivit\'e de l'application (1), il suffit de prouver qu'ils v\'erifient aussi la condition (3). Pour cela, fixons $\tilde{T}\in K\tilde{{\cal T}}_{ell}$ et $\eta\in \tilde{T}({\mathbb R})$.  Introduisons l'ensemble $\dot{\cal X}^{\cal E}(\eta)$. Comme ci-dessus, on peut supposer que tout \'el\'ement de cet ensemble est de la forme 
 $({\bf G}',\epsilon)$, o\`u ${\bf G}'\in {\cal E}(K\tilde{G},{\bf a})_{0}$ et $\epsilon\in \tilde{T}[\tilde{G}']({\mathbb R})$. On a m\^eme  $\epsilon\in \tilde{T}[\tilde{G}']({\mathbb R})_{\sharp}$ d'apr\`es 4.9(6). Soit $X_{0}\in \mathfrak{t}^{\theta}({\mathbb R})$ assez petit et r\'egulier dans $\mathfrak{g}_{\eta}$. L'\'el\'ement $\gamma_{0}=exp(X_{0})\eta$ est elliptique et fortement r\'egulier. Introduisons l'ensemble $\dot{\cal Y}^{\cal E}(\gamma_{0})$ et, pour simplifier, indexons-le par un ensemble $\{1,...,n\}$ d'entiers. D'apr\`es la remarque suivant le lemme 4.9, on peut supposer que, pour $k=1,...,n$, le $k$-i\`eme \'el\'ement de $\dot{\cal Y}^{\cal E}(\gamma_{0})$ est de la forme $({\bf G}_{k}',exp(Y_{k,0})\epsilon_{k})$, o\`u $({\bf G}'_{k},\epsilon_{k})\in \dot{\cal X}^{\cal E}(\eta)$ et $Y_{k,0}$ est un \'el\'ement r\'egulier de $\mathfrak{g}'_{k,\epsilon_{k}}({\mathbb R})$. Remarquons en passant que l'application $k\mapsto ({\bf G}'_{k},\epsilon_{k})$ n'est pas injective en g\'en\'eral. Notons $\tilde{T}'_{k}=\tilde{T}[\tilde{G}'_{k}]$. L'\'el\'ement $Y_{k,0}$ est elliptique. Puisque $T'_{k}$ est, \`a conjugaison pr\`es, l'unique sous-tore elliptique de $G'_{k,\epsilon_{k}}$, on peut supposer  $Y_{k,0}\in \mathfrak{t}_{k}'({\mathbb R})$. D'un diagramme reliant $exp(Y_{k,0})\epsilon_{k}$ \`a $exp(X_{0})\eta$ se d\'eduit alors un isomorphisme $\mathfrak{t}^{\theta}({\mathbb R})\simeq \mathfrak{t}_{k}'({\mathbb R})$ qui envoie $X_{0}$ sur $Y_{k,0}$.  En fixant une section $\mathfrak{t}_{k}'({\mathbb R})\to \mathfrak{t}'_{k,1}({\mathbb R})$ de la projection naturelle, on obtient un homomorphisme
$$\begin{array}{ccc}\mathfrak{t}^{\theta}({\mathbb R})&\to&\mathfrak{t}'_{k,1}({\mathbb R})\\ X&\mapsto&Y_{k}\\ \end{array}$$
Soit $X\in \mathfrak{t}^{\theta}({\mathbb R})$, assez petit et r\'egulier dans $\mathfrak{g}_{\eta}$, et posons $\gamma=exp(X)\eta$. Il est (plus ou moins) clair que l'on peut prendre pour ensemble $\dot{\cal X}^{\cal E}(\gamma)$ l'ensemble des $({\bf G}'_{k},exp(Y_{k})\epsilon_{k})$ pour $k=1,...,n$. En appliquant la d\'efinition ci-dessus, on obtient
$$\varphi_{\tilde{T}}(exp(X)\eta)=d(\theta^*)^{-1/2} [T^{\theta}({\mathbb R}):T^{\theta,0}({\mathbb R})]\vert \dot{\cal X}(\gamma_{0})\vert ^{-1}$$
$$ \sum_{ k=1,...,n }\Delta_{1}(exp(Y_{k})\epsilon_{k,1},exp(X)\eta)^{-1}\varphi_{\tilde{T}'_{k,1}}(exp(Y_{k})\epsilon_{k,1}).$$
On veut prouver que la fonction $X\mapsto \Delta_{\eta}(X)\varphi_{\tilde{T}}(exp(X)\eta)$ se prolonge en une fonction $C^{\infty}$ au voisinage de $0$. On sait d'apr\`es (6) que, pour tout $k$, la fonction $Y\mapsto\Delta_{\epsilon_{k}}(Y)\varphi_{\tilde{T}'_{k,1}}(exp(Y)\epsilon_{k,1})$ se prolonge en une telle fonction. Il suffit donc de prouver que, pour tout $k$, la fonction
$$X\mapsto \Delta_{\eta}(X)\Delta_{\epsilon_{k}}(Y_{k})^{-1}\Delta_{1}(exp(Y_{k})\epsilon_{k,1},exp(X)\eta)^{-1}$$
se prolonge en une fonction $C^{\infty}$ au voisinage de $0$. C'est ce que fait Shelstad dans [S1], dans une situation plus g\'en\'erale. Puisque l'on est ici dans un cas beaucoup plus simple, redonnons l'argument. Pour simplifier, on fixe $k$ et on abandonne les indices $k$. Il existe une constante $c\not=0$ telle que
$$\Delta_{1}(exp(Y)\epsilon_{1},exp(X)\eta)=c\boldsymbol{\Delta}_{1}(exp(Y)\epsilon_{1},exp(X)\eta;exp(Y_{0})\epsilon_{1},exp(X_{0})\eta).$$
Il est clair que le facteur $\Delta_{imp}(exp(Y)\epsilon_{1},exp(X)\eta;exp(Y_{0})\epsilon_{1},exp(X_{0})\eta)^{-1}$ est $C^{\infty}$ au voisinage de $0$.  Cela nous ram\`ene \`a consid\'erer la fonction
$$X\mapsto \Delta_{\eta}(X)\Delta_{\epsilon}(Y)^{-1}\Delta_{II}(exp(Y)\epsilon,exp(X)\eta)^{-1}.$$
Utilisons les notations de 1.6 et 2.2.    Le terme $\nu$ de 2.2 est de la forme $exp(X)\nu_{\eta}$. On a d\'ecrit en [W1] 3.3 l'ensemble de racines $\Sigma(T)_{\eta}$ du groupe $G_{\eta}$.  C'est
$$\Sigma(T)_{\eta}=\{\alpha_{res}; \alpha\in \Sigma(T), \,\,\alpha\text{ de type 1 ou 2 }, (N\alpha)(\nu_{\eta})=1\}$$
$$\cup\{2\alpha_{res}; \alpha\in \Sigma(T),\,\, \alpha\text{ de type 2 }, (N\alpha)(\nu_{\eta})=-1\}.$$
 On a aussi d\'ecrit l'ensemble de racines $\Sigma(T')_{\epsilon}$ du groupe $G'_{\epsilon}$.  C'est
$$\Sigma(T')_{\epsilon}=\{N\alpha; \alpha\in \Sigma(T),\,\,\alpha\text{ de type 1}, (N\hat{\alpha})(s)=1,\,(N\alpha)(\nu_{\eta})=1\}$$
$$
\cup\{2N\alpha; \alpha\in \Sigma(T),\,\,\alpha\text{ de type 2}, (N\hat{\alpha})(s)=1,\,(N\alpha)(\nu_{\eta})=\pm 1\}$$
$$\cup\{N\alpha; \alpha\in \Sigma(T),\,\,\alpha\text{ de type 2}, (N\hat{\alpha})(s)=-1,\,(N\alpha)(\nu_{\eta})=1\}.$$
Puisque $(T^{\theta,0}/A_{\tilde{G}})({\mathbb R})$ est elliptique, la conjugaison complexe agit sur $\Sigma(T)_{res,ind}$ par multiplication par $-1$. Fixons un ensemble $\Sigma_{\star}$ de repr\'esentants des orbites. Dans les d\'efinitions de $\Delta_{\eta}$ et $\Delta_{\epsilon}$, on peut remplacer les sous-ensembles de racines positives par des ensembles de repr\'esentants d'orbites pour la conjugaison complexe, cela ne change  ces fonctions que par des constantes. On peut supposer que ce sont les ensembles d\'eduits de ceux ci-dessus en ajoutant la condition $\alpha_{res}\in \Sigma_{\star}$.  Chacune des nos fonctions $\Delta_{\eta}(X)$, $\Delta_{\epsilon}(Y)^{-1}$ et $\Delta_{II}(exp(Y)\epsilon,exp(X)\eta)^{-1} $ est un produit index\'e par $\alpha_{res}\in \Sigma_{\star}$. Le terme index\'e par $\alpha_{res}$ est donn\'e par le tableau suivant
$$\begin{array}{cccccc}\text{ type de }\alpha&(N\alpha)(\nu_{\eta})&(N\hat{\alpha})(s)&\Delta_{\eta}(X)&\Delta_{\epsilon}(Y)^{-1}&\Delta_{II}(exp(Y)\epsilon,exp(X)\eta)^{-1} \\1&1&1&sgn(i\alpha_{res}(X))&sgn(i(N\alpha)(Y))&1\\1&1&\not=1&sgn(i\alpha_{res}(X))&1&\chi_{\alpha_{res}}(\frac{a_{\alpha_{res}}}{(N\alpha)(\nu)-1})\\1&\not=1&1&1&1&1\\1&\not=1&\not=1&1&1&\chi_{\alpha_{res}}(\frac{a_{\alpha_{res}}}{(N\alpha)(\nu)-1})\\2&1&1&sgn(i\alpha_{res}(X))&sgn(2i(N\alpha)(Y))&1\\2&1&-1&sgn(i\alpha_{res}(X))&sgn(i(N\alpha)(Y))&\chi_{\alpha_{res}}(\frac{1}{(N\alpha)(\nu)+1})\\2&1&\not=\pm 1&sgn(i\alpha_{res}(X))&1&\chi_{\alpha_{res}}(\frac{a_{\alpha_{res}}}{(N\alpha)(\nu)^2-1})\\2&-1&1&sgn(2i\alpha_{res}(X))&sgn(2i(N\alpha)(Y))&1\\2&-1&-1&sgn(2i\alpha_{res}(X))&1&\chi_{\alpha_{res}}(\frac{1}{(N\alpha)(\nu)+1})\\2&-1&\not=\pm 1&sgn(2i\alpha_{res}(X))&1&\chi_{\alpha_{res}}(\frac{a_{\alpha_{res}}}{(N\alpha)(\nu)^2-1})\\2&\not=\pm 1&1&1&1& 1\\2&\not=\pm 1&-1&1&1&\chi_{\alpha_{res}}(\frac{1}{(N\alpha)(\nu)+1})\\2&\not=\pm 1&\not=\pm 1&1&1&\chi_{\alpha_{res}}(\frac{a_{\alpha_{res}}}{(N\alpha)(\nu)^2-1})\\ \end{array}$$
 On peut choisir les $a$-data et les $\chi$-data de sorte que, pour tout $\alpha_{res}\in \Sigma_{\star}$, $a_{\alpha_{res}}=i$ et $\chi_{\alpha_{res}}(z)=z/\vert z\vert $.   On v\'erifie alors que, dans chaque cas, le produit des trois contributions ci-dessus est $C^{\infty}$ au voisinage de $X=0$. Par exemple, consid\'erons le cas $\alpha$ de type 2, $(N\alpha)(\nu_{\eta})=1$ et $(N\hat{\alpha})(s)=1$. L'homomorphisme $X\mapsto Y$ identifie $N\alpha$ \`a  $n_{\alpha}\alpha_{res}$, o\`u $n_{\alpha}$ est le plus petit entier $n\geq1$ tel que $\theta^{n}(\alpha)=\alpha$. Donc $sgn(2i(N\alpha)(Y))=sgn(i\alpha_{res}(X))$ et le produit de ces deux termes vaut $1$. Consid\'erons maintenant le cas  $\alpha$ de type 2, $(N\alpha)(\nu_{\eta})=1$ et $(N\hat{\alpha})(s)\not=\pm 1$. On a $(N\alpha)(\nu)^2=exp(2(N\alpha)(X))(N\alpha)(\nu_{\eta})^2 =exp(2n_{\alpha}\alpha_{res}(X))$ d'o\`u
 $$\chi_{\alpha_{res}}(\frac{a_{\alpha_{res}}}{(N\alpha)(\nu)^2-1})=i\vert exp(2n_{\alpha}\alpha_{res}(X))-1\vert (exp(2n_{\alpha}\alpha_{res}(X)-1)^{-1}.$$
 Le produit de cette expression avec $sgn(i\alpha_{res}(X))$ est $C^{\infty}$ au voisinage de $0$.
 On laisse les autres cas au lecteur. Cela ach\`eve la preuve.
 
 \bigskip
 
 \subsection{Un corollaire de la preuve dans le cas r\'eel}
 Le corps de base est ${\mathbb R}$. On suppose $(G,\tilde{G},\omega) $ quasi-d\'eploy\'e et \`a torsion int\'erieure. Notons $I_{cusp}^{st}(\tilde{G}({\mathbb R}))$ le sous-espace des  $f\in I_{cusp}(\tilde{G}({\mathbb R}))$ tels que la fonction $\gamma\mapsto I^{\tilde{G}}(\gamma,f)$ est constante sur les classes de conjugaison stable form\'ees d'\'el\'ements fortement r\'eguliers et elliptiques. 
 
 \ass{Lemme}{L'application naturelle $I^{st}_{cusp}(\tilde{G}({\mathbb R}))\to SI_{cusp}(\tilde{G}({\mathbb R}))$ est un isomorphisme.}
 
 {\bf Remarque.} Ce lemme vaut aussi sur un corps $F$ non-archim\'edien mais, dans ce cas, c'est une cons\'equence directe de la proposition 4.11. Dans le cas pr\'esent o\`u le corps de base est ${\mathbb R}$, cette proposition ne s'applique qu'\`a un $K$-espace. Ici, nous consid\'erons un seul espace $\tilde{G}$.
 \bigskip
 
 Preuve. On peut supposer que $\tilde{G}$ contient un tore tordu maximal elliptique, sinon les deux espaces sont nuls. Puisque $\tilde{G}$ est \`a torsion int\'erieure, il n'en contient qu'un \`a conjugaison pr\`es. On en fixe un, que l'on note $\tilde{T}$. L'espace $I_{cusp}(\tilde{G}({\mathbb R}))$, resp. $SI_{cusp}(\tilde{G}({\mathbb R}))$, s'identifie \`a celui  des fonctions $\varphi_{\tilde{T}}$ d\'efinies presque partout sur $\tilde{T}({\mathbb R})$ qui v\'erifient les conditions (2) et (3) du paragraphe pr\'ec\'edent, resp. (5) et (6) (la condition (4) est triviale  en identifiant $SI({\bf G})$ \`a $SI(\tilde{G}({\mathbb R}))$). On voit que ces deux derni\`eres conditions sont \'equivalentes \`a la r\'eunion des deux premi\`eres et de la condition: $\varphi_{\tilde{T}}$ est constante sur les classes de conjugaison stable form\'ees d'\'el\'ements fortement r\'eguliers et elliptiques. Il en r\'esulte que $\varphi_{\tilde{T}}\in SI_{cusp}(\tilde{G}({\mathbb R}))$ si et seulement si $\varphi_{\tilde{T}}\in I_{cusp}^{st}(\tilde{G}({\mathbb R}))$. On n'a pas tout-\`a-fait fini car l'application naturelle $I^{st}_{cusp}(\tilde{G}({\mathbb R}))\to SI_{cusp}(\tilde{G}({\mathbb R}))$ ne se traduit pas par l'identit\'e en termes de fonctions sur $\tilde{T}({\mathbb R})$, mais par l'application
 $\varphi_{\tilde{T}}\mapsto \varphi^{{\bf G}}_{\tilde{T}}$ d\'efinie par
 $$\varphi^{{\bf G}}_{\tilde{T}}(\delta)=\sum_{\gamma}\varphi_{\tilde{T}}(\gamma),$$
 o\`u on somme sur les $\gamma\in \tilde{T}({\mathbb R})$ stablement conjugu\'es \`a $\delta$, \`a conjugaison pr\`es par $G({\mathbb R})$. Il reste \`a voir que le nombre de ces \'el\'ements $\gamma$ ne d\'epend pas de $\delta$, pourvu que $\delta$ soit fortement r\'egulier. Mais ce nombre est \'egal au nombre d'\'el\'ements de l'ensemble
 $$T({\mathbb C})\backslash \{g\in G({\mathbb C}); g\sigma(g)^{-1}\in T({\mathbb C}) \text{ pour tout }\sigma\in \Gamma_{{\mathbb R}}\}/G({\mathbb R}).$$
 Cela ach\`eve la preuve. $\square$

 \bigskip
 
 \subsection{Filtration de l'espace $SI(\tilde{G}(F))$}
 
 On suppose $(G,\tilde{G},\omega)$ quasi-d\'eploy\'e et \`a torsion int\'erieure. On a filtr\'e en 4.2 l'espace $I(\tilde{G}(F))$. Il y a deux filtrations naturelles sur $SI(\tilde{G}(F))$. Pour un entier $n\geq-1$, notons ${\cal F}^nSI(\tilde{G}(F))$ le sous-espace des $f\in SI(\tilde{G}(F))$ tels que $f_{\tilde{M}}=0$ pour tout espace de Levi $\tilde{M}$ tel que $a_{\tilde{M}}>n$. Ces espaces forment l'une des filtrations. On note $GrSI(\tilde{G}(F))$ le gradu\'e associ\'e. On peut d'autre part consid\'erer l'image de la filtration de $I(\tilde{G}(F))$ par la projection naturelle de cet espace sur $SI(\tilde{G}(F))$. Autrement dit, si on note $I^{inst}(\tilde{G}(F))$ le noyau de cette projection, les termes de la filtration sont les espaces
 $$({\cal F}^nI(\tilde{G}(F))+I^{inst}(\tilde{G}(F)))/I^{inst}(\tilde{G}(F)).$$
 Il est clair que l'espace ci-dessus est inclus dans ${\cal F}^nSI(\tilde{G}(F))$.
 
\ass{Lemme}{Pour tout $n$, on a les \'egalit\'es:

$$({\cal F}^nI(\tilde{G}(F))+I^{inst}(\tilde{G}(F)))/I^{inst}(\tilde{G}(F))={\cal F}^nSI(\tilde{G}(F)$$
et
$$Gr^nSI(\tilde{G}(F)=\oplus_{\tilde{M}\in \underline{{\cal L}}^n}SI_{cusp}(\tilde{M}(F))^{W(\tilde{M})}.$$}

Preuve. Notons pour simplifier $E^n$ l'espace de gauche de la premi\`ere \'egalit\'e. On raisonne par r\'ecurrence et on suppose prouv\'e que $E^{n-1}={\cal F}^{n-1}SI(\tilde{G}(F))$. Puisque $E^n\subset {\cal F}^nSI(\tilde{G}(F))$, on a alors une injection
$$(1) \qquad E^n/E^{n-1}\subset Gr^nSI(\tilde{G}(F)).$$
Il s'agit de voir qu'elle est surjective. Le premier espace est quotient de $Gr^nI(\tilde{G}(F))$, ou encore, en utilisant le lemme 4.2, de
$$\oplus_{\tilde{M}\in \underline{{\cal L}}^n}I_{cusp}(\tilde{M}(F))^{W(\tilde{M})}.$$
Par d\'efinition, l'espace $Gr^nSI(\tilde{G}(F))$ s'envoie injectivement dans
$$\oplus_{\tilde{M}\in \underline{{\cal L}}^n}SI_{cusp}(\tilde{M}(F))^{W(\tilde{M})}.$$
L'homomorphisme (1)  compos\'e avec cette injection se quotiente en l'homomorphisme naturel
$$\oplus_{\tilde{M}\in \underline{{\cal L}}^n}I_{cusp}(\tilde{M}(F))^{W(\tilde{M})}\to \oplus_{\tilde{M}\in \underline{{\cal L}}^n}SI_{cusp}(\tilde{M}(F))^{W(\tilde{M})}.$$
Pour prouver les deux assertions de l'\'enonc\'e, il suffit de prouver que ce dernier est surjectif. Mais c'est un cas particulier de l'assertion 4.12(3) dans le cas non archim\'edien et c'est le lemme 4.14 dans le cas r\'eel (le cas complexe est trivial). $\square$

Comme toujours, il y a une variante de ce r\'esultat quand on consid\`ere des extensions centrales comme \`a la fin du paragraphe 4.8.

\bigskip

\subsection{Un corollaire  }

On suppose encore $(G,\tilde{G},{\bf a})$ quasi-d\'eploy\'e et \`a torsion int\'erieure. Soit $(\tilde{M}_{j})_{j=1,...,k}$ une famille finie d'espaces de Levi de $\tilde{G}$. Consid\'erons l'application lin\'eaire
$$res=\oplus_{j=1,...,k}res_{\tilde{M}_{j}}:I(\tilde{G}(F))\to \oplus_{j=1,...,k}I(\tilde{M}_{j}(F)).$$
\ass{Corollaire}{On a l'\'egalit\'e
$$res(I(\tilde{G}(F)))\cap \left(\oplus_{j=1,...,k}I^{inst}(\tilde{M}_{j}(F))\right)=res(I^{inst}(\tilde{G}(F))).$$}

Preuve. Posons
$$I=\oplus_{j=1,...,k}I(\tilde{M}_{j}(F)),\,\,I^{inst}=\oplus_{j=1,...,k}I^{inst}(\tilde{M}_{j}(F))$$
et, pour tout $n\in {\mathbb N}$, ${\cal F}^nI=\oplus_{j=1,...,k}{\cal F}^nI(\tilde{M}_{j}(F))$.
On va prouver  que, pour tout $n\in {\mathbb N}$,
$$(1)\qquad res({\cal F}^nI(\tilde{G}(F)))\cap I^{inst}\subset res(I^{inst}(\tilde{G}(F)))+(res({\cal F}^{n-1}I(\tilde{G}(F)))\cap I^{inst}).$$
Posons ${\cal F}^nI^{inst}(\tilde{G}(F))=I^{inst}(\tilde{G}(F))\cap {\cal F}^nI(\tilde{G}(F))$. On note $GrI^{inst}(\tilde{G}(F))$ le gradu\'e associ\'e \`a cette filtration.  En cons\'equence du lemme 4.15, la suite
$$0\to Gr^nI^{inst}(\tilde{G}(F))\to Gr^nI(\tilde{G}(F))\to Gr^nSI(\tilde{G}(F))\to 0$$
est exacte. Donc $Gr^nI^{inst}(\tilde{G}(F))$ est l'espace des $(f^{\tilde{L}})_{\tilde{L}\in \underline{{\cal L}}^n}\in\oplus_{\tilde{L}\in \underline{{\cal L}}^n}I_{cusp}(\tilde{L}(F))^{W(\tilde{L})}$ tels que les images de $f^{\tilde{L}}$ dans $SI_{cusp}(\tilde{L}(F))$ soient nulles pour tout $\tilde{L}$. Soit $f\in {\cal F}^nI(\tilde{G}(F))$ tel que $res(f)\in I^{inst}$. Soit $(f^{\tilde{L}})_{\tilde{ L}\in \underline{{\cal L}}^n}$ son image dans $\oplus_{\tilde{ L}\in \underline{{\cal L}}^n}I_{cusp}(\tilde{L}(F))^{W(\tilde{L})}$.  Notons $\underline{{\cal L}}^n_{\star}$ l'ensemble des $\tilde{L}\in \underline{{\cal L}}^n$ qui sont conjugu\'es par $G(F)$  \`a un espace inclus dans l'un des $\tilde{M}_{j}$. L'hypoth\`ese $res(f)\in I^{inst}$ entra\^{\i}ne que, si $\tilde{L}\in \underline{{\cal L}}^n_{\star}$, l'image de $f^{\tilde{L}}$ dans $SI_{cusp}(\tilde{L}(F))$ est nulle.
  Par le r\'esultat pr\'ec\'edent, on peut trouver $f_{0}\in {\cal F}^nI^{inst}(\tilde{G}(F))$ dont l'image $(f_{0}^{\tilde{L}})_{\tilde{ L}\in \underline{{\cal L}}^n} $ dans le gradu\'e v\'erifie $f_{0}^{\tilde{L}}=f^{\tilde{L}}$ si $\tilde{L}\in \underline{{\cal L}}^n_{\star}$, $f_{0}^{\tilde{L}}=0$ sinon. Alors, pour tout $j=1,...,k$, l'image de $res_{\tilde{M}_{j}}(f-f_{0})$ dans $Gr^nI(\tilde{M}_{j}(F))$ est nulle. Autrement dit $res(f-f_{0})\in {\cal F}^{n-1}I$. D'apr\`es la preuve du lemme 4.3, $res(I(\tilde{G}(F)))\cap {\cal F}^{n-1}I=res({\cal F}^{n-1}I(\tilde{G}(F)))$. Il existe donc $f'\in {\cal F}^{n-1}I(\tilde{G}(F))$ tel que $res(f-f_{0}-f')=0$.  On a encore $res(f')\in I^{inst}$. L'\'egalit\'e $res(f)=res(f_{0})+res(f')$ montre que $res(f)$ appartient au membre de droite de (1). Cela prouve cette relation.

Par r\'ecurrence sur $n$,  (1) implique que le membre de gauche de l'\'enonc\'e est inclus dans celui de droite. L'inclusion oppos\'ee \'etant \'evidente, cela d\'emontre le corollaire. $\square$

   \bigskip
\subsection{Produit scalaire}
Dans ce paragraphe, on suppose $\omega$ unitaire. On munit $G(F)$ d'une mesure de Haar. On doit aussi munir $A_{\tilde{G}}(F)$ d'une telle mesure. Par souci de coh\'erence avec [W3], on proc\`ede ainsi. On munit l'espace vectoriel r\'eel ${\cal A}_{\tilde{G}}$ d'une mesure de Haar. On dispose de l'homomorphisme habituel
$$H_{A_{\tilde{G}}}:A_{\tilde{G}}(F)\to {\cal A}_{\tilde{G}}.$$
Pour $a\in A_{\tilde{G}}(F)$ et $x^*\in X^*(A_{\tilde{G}})$, on a $\vert x^*(a)\vert_{F}=e^{<x^*,H_{A_{\tilde{G}}}(a)>}$. Notons $A_{\tilde{G}}(F)^{c}$ le noyau de $H_{A_{\tilde{G}}}$. C'est le sous-groupe compact maximal de $A_{\tilde{G}}(F)$. Si $F$ est non-archim\'edien, l'image $Im(H_{A_{\tilde{G}}})$ de l'homomorphisme $H_{A_{\tilde{G}}}$ est un r\'eseau de ${\cal A}_{\tilde{G}}$, tandis que $A_{\tilde{G}}(F)^{c}$ est un sous-groupe ouvert de $A_{\tilde{G}}(F)$. On munit $A_{\tilde{G}}(F)$ de la mesure de Haar telle que $mes(A_{\tilde{G}}(F)^{c})=mes({\cal A}_{\tilde{G}}/Im(H_{A_{\tilde{G}}}))$. Si $F$ est archim\'edien, on munit $A_{\tilde{G}}(F)^{c}$ de la mesure de Haar de masse totale $1$. La suite
$$1\to A_{\tilde{G}}(F)^{c}\to A_{\tilde{G}}(F)\to {\cal A}_{\tilde{G}}\to 0$$
est exacte et on munit $A_{\tilde{G}}(F)$ de la mesure compatible avec cette suite et avec les mesures d\'ej\`a fix\'ees sur les deux autres groupes.

 Commen\c{c}ons par supposer $F$ non-archim\'edien.    Pour tout sous-tore tordu maximal $\tilde{T}$ de $\tilde{G}$, munissons $T^{\theta,0}(F)$ d'une mesure de Haar. Notons $\tilde{G}_{reg}(F)/conj$ l'ensemble des classes de conjugaison par $G(F)$ dans l'ensemble $\tilde{G}_{reg}(F)$. Pour $\gamma\in \tilde{G}_{reg}(F)$, l'application
 $$\begin{array}{ccc} G_{\gamma}(F)&\to&\tilde{G}_{reg}(F)/conj\\t&\mapsto& classe(t\gamma)\\ \end{array}$$
 est injective dans un voisinage de $1$. On munit $\tilde{G}_{reg}(F)/conj$ de la topologie (ou de la structure de vari\'et\'e analytique sur $F$) et de la mesure telle que, pour tout $\gamma$, cette application soit, au voisinage de $1$, un isomorphisme pr\'eservant la mesure. On a alors la formule d'int\'egration, pour $f\in C_{c}^{\infty}(\tilde{G}(F))$:
 $$\int_{\tilde{G}(F)}f(\gamma) d\gamma=\int_{\tilde{G}_{reg}(F)/conj}\Phi(\gamma,f) D^{\tilde{G}}(\gamma)d\gamma,$$
 o\`u  
 $$\Phi(\gamma,f)=\int_{Z_{G}(\gamma;F)\backslash G(F)}f(g^{-1}\gamma g)dg.$$
 
 Soient $f_{1},f_{2}\in C_{c}^{\infty}(\tilde{G}(F))$, supposons les supports de $f_{1}$ et $f_{2}$ contenus dans l'ensemble $\tilde{G}(F)_{ell}$ des \'el\'ements elliptiques r\'eguliers de $\tilde{G}(F)$. On pose
 $$(1) \qquad J_{\tilde{G}}(\omega,f_{1},f_{2})=\int_{A_{\tilde{G}}(F)\backslash G(F)}\int_{\tilde{G}(F)} \overline{f_{1}(\gamma)}f_{2}(g^{-1} \gamma g)d\gamma \omega(g) dg.$$
 Cette int\'egrale est absolument convergente et on a
 $$(2) \qquad J_{\tilde{G}}(\omega,f_{1},f_{2})=\int_{\tilde{G}(F)_{ell}/conj}i(\gamma)^{-1}mes(A_{\tilde{G}}(F)\backslash G_{\gamma}(F))\overline{I^{\tilde{G}}(\gamma,\omega,f_{1})}I^{\tilde{G}}(\gamma,\omega,f_{2})d\gamma,$$
 o\`u on a pos\'e $i(\gamma)=[Z_{G}(\gamma;F):G_{\gamma}(F)]$ et o\`u on  rappelle la d\'efinition
 $$I^{\tilde{G}}(\gamma,\omega,f)=\left\lbrace\begin{array}{cc}D^{\tilde{G}}(\gamma)^{1/2}\int_{G_{\gamma}(F)\backslash G(F)}\omega(g)f(g^{-1}\gamma g)dg,& \text{ si }\omega\text{ est trivial sur }Z_{G}(\gamma;F)\\ 0,& \text{ sinon. }\\ \end{array}\right..$$
 Dans la formule (2), on peut consid\'erer que $f_{1}$ et $f_{2}$ ne sont plus des fonctions sur $\tilde{G}(F)_{ell}$ mais sont plut\^ot leurs images dans $I_{cusp}(\tilde{G}(F),\omega)$. Cela d\'efinit un produit hermitien sur un sous-espace de $I_{cusp}(\tilde{G}(F),\omega)$, \`a savoir l'image de l'espace des fonctions \`a support elliptique r\'egulier. Il r\'esulte de la formule des traces locale que la m\^eme formule (2) s'\'etend en un produit hermitien sur tout l'espace $I_{cusp}(\tilde{G}(F),\omega)$ (c'est-\`a-dire que cette formule reste absolument convergente), cf. [W3] 6.6(1).
 
 Consid\'erons le cas particulier o\`u $(G,\tilde{G},{\bf a})$ est quasi-d\'eploy\'e et \`a torsion int\'erieure. On dispose de la donn\'ee endoscopique maximale ${\bf G}$ pour laquelle $SI({\bf G})=SI(\tilde{G}(F))$. On a aussi $SI_{cusp}({\bf G})=SI_{cusp}(\tilde{G}(F))$. La proposition 4.11  identifie cet espace \`a un sous-espace de $I_{cusp}(\tilde{G}(F))$. C'est le sous-espace des $f\in I_{cusp}(\tilde{G}(F))$ dont les int\'egrales orbitales sont constantes sur toute classe de conjugaison stable fortement r\'eguli\`ere. Le produit hermitien ci-dessus se restreint en un tel produit sur ce sous-espace. 
 Notons  $\tilde{G}(F)_{ell}$ l'ensemble des \'el\'ements fortement r\'eguliers et elliptiques de $\tilde{G}(F)$ et $\tilde{G}(F)_{ell}/ st-conj$ l'ensemble des classes de conjugaison stable contenues dans $\tilde{G}(F)_{ell}$. Par le m\^eme proc\'ed\'e que ci-dessus, on le munit d'une topologie et d'une mesure. Pour $f_{1},f_{2}\in SI(\tilde{G}(F))$, on a l'\'egalit\'e
 $$(3) \qquad J_{\tilde{G}}(f_{1},f_{2})=\int_{\tilde{G}(F)_{ell}/st-conj}k(\delta)^{-1}mes(A_{\tilde{G}}(F)\backslash G_{\delta}(F))\overline{S^{\tilde{G}}(\delta,f_{1})}S^{\tilde{G}}(\delta,f_{2})d\delta$$
o\`u, pour toute classe de conjugaison stable $\delta$, on a not\'e $k(\delta) $ le nombre de classes de conjugaison par $G(F)$ contenues dans  $\delta$. Remarquons que  les centralisateurs sont connexes dans le cas o\`u la torsion est int\'erieure.

Revenons au cas g\'en\'eral, soit ${\bf G}'=(G',{\cal G}',\tilde{s})\in {\cal E}(\tilde{G},{\bf a})$. On peut choisir des donn\'ees auxiliaires $G'_{1}$,...,$\Delta_{1}$ de sorte que le caract\`ere $\lambda_{1}$ soit unitaire. Pour $f_{1},f_{2}\in SI_{\lambda_{1},cusp}(\tilde{G}'_{1}(F))$, la fonction
$$\delta_{1}\mapsto \overline{S^{\tilde{G}'_{1}}(\delta_{1},f_{1})}S^{\tilde{G}'_{1}}(\delta_{1},f_{2})$$ 
sur $\tilde{G}'_{1}(F)_{ell}$ se descend en une fonction de $\delta\in \tilde{G}'(F)_{ell}/st-conj$. Modulo les choix de mesures de Haar sur $G'(F)$ et ${\cal A}_{G'}$ (de cette derni\`ere se d\'eduisant une mesure sur $A_{G'}(F)$ comme plus haut), on peut donc d\'efinir le produit $J_{\tilde{G}'}(f_{1},f_{2})$ par la formule (3) o\`u $\tilde{G}$ est remplac\'e par $\tilde{G}'$. Quand on change de donn\'ees auxiliaires, ces formules se recollent et on obtient un produit hermitien $J_{{\bf G}'}$ sur l'espace $SI_{cusp}({\bf G}')$.

On suppose maintenant fix\'ees des mesures de Haar sur $G(F)$, sur $ {\cal A}_{\tilde{G}}$ et sur $G'(F)$ pour tout ${\bf G}'\in {\cal E}(\tilde{G},{\bf a})$. Pour tout tel ${\bf G}'$, on a un isomorphisme naturel $ {\cal A}_{\tilde{G}}\to {\cal A}_{G'}$.  On munit ${\cal A}_{G'}$ de la mesure telle que cet isomorphisme pr\'eserve  les mesures. 
On pose
$$c(\tilde{G},{\bf G}')= det((1-\theta)_{\vert {\cal A}_{G}/{\cal A}_{\tilde{G}}})\vert ^{-1}
\vert \pi_{0}(Z(\hat{G})^{\Gamma_{F}})\vert \vert Z(\hat{G}')^{\Gamma_{F}})\vert ^{-1}\vert $$
$$\vert Out({\bf G}')\vert ^{-1}\vert \pi_{0}(Z(\hat{G})^{\Gamma_{F},0}\cap  \hat{G}')\vert \vert \pi_{0}((Z(\hat{G})/(Z(\hat{G})\cap \hat{G}'))^{\Gamma_{F}})\vert ^{-1}.$$
La proposition 4.11 nous fournit un isomorphisme
$$\begin{array}{ccc}I_{cusp}(\tilde{G}(F),\omega)&\simeq &\oplus_{{\bf G}'\in {\cal E}(\tilde{G},{\bf a})}SI_{cusp}({\bf G}')^{Aut({\bf G}')}\\ f&\mapsto&(f^{{\bf G}'})_{{\bf G}'\in {\cal E}(\tilde{G},{\bf a})}\\ \end{array}.$$
Chaque espace est muni d'un produit hermitien.  

On a suppos\'e le corps $F$ non-archim\'edien. Dans le cas o\`u $F$ est r\'eel, toutes ces constructions s'adaptent aux $K$-espaces. Le produit hermitien sur $I_{cusp}(K\tilde{G}({\mathbb R}),\omega)$ est la somme directe des produits sur les diff\'erents $I_{cusp}(\tilde{G}_{p}({\mathbb R}),\omega)$.  Attention: dans la formule (3), $k(\delta)$ est un nombre de classes de conjugaison dans un $K$-espace  associ\'e \`a $\tilde{G}'$.

\ass{Proposition}{Soient $\underline{f},f\in I_{cusp}(\tilde{G}(F),\omega)$. Alors on a l'\'egalit\'e
$$J_{\tilde{G}}(\omega,\underline{f},f)=\sum_{{\bf G}'\in {\cal E}(\tilde{G},{\bf a})}c(\tilde{G},{\bf G}')J_{{\bf G}'}(\underline{f}^{{\bf G}'},f^{{\bf G}'}).$$}

{\bf Remarque.} La d\'emonstration s'inspire de celle du lemme 6.4.B de [KS1].
\bigskip

Preuve. Tous nos espaces d'int\'egration sont des rev\^etements de l'espace $\tilde{G}(F)_{ell}/st-conj$, cf. 4.9(7). Les mesures sur nos espaces d\'ependent de choix de mesures sur les tores. Si on impose \`a ces choix la m\^eme condition qu'en 2.4 (les mesures sur deux tores se correspondent localement quand il y a un isomorphisme naturel  entre ces deux tores), les rev\^etements pr\'eservent localement les mesures. L'\'egalit\'e de l'\'enonc\'e r\'esulte d'une \'egalit\'e plus forte: quand on consid\`ere les deux c\^ot\'es de la formule comme des int\'egrales sur $\tilde{G}(F)_{ell}/st-conj$, les fonctions que l'on int\`egre sont \'egales. C'est ce que l'on va prouver. Fixons $\gamma\in \tilde{G}(F)_{ell}$ et consid\'erons les valeurs de nos fonctions sur la classe de conjugaison stable de $\gamma$. Si $\omega$ n'est pas trivial sur $Z_{G}(\gamma;F)$, ces deux valeurs sont nulles. On suppose $\omega$ trivial sur $Z_{G}(\gamma;F)$. Pour ${\bf G}'\in {\cal E}(\tilde{G},{\bf a})$, le groupe $Out({\bf G}')$ agit librement sur l'ensemble des \'el\'ements de $\tilde{G}'(F)_{ell}/st-conj$ qui se projettent sur cette classe de conjugaison stable. L'ensemble  $\dot{\cal X}^{\cal E}(\gamma)$ est un ensemble de repr\'esentants de ces orbites. La fonction du membre de droite  vaut donc
$$(4) \qquad \sum_{({\bf G}',\delta)\in \dot{\cal X}^{\cal E}(\gamma)}c(\tilde{G},{\bf G}')\vert Out({\bf G}')\vert k(\delta)^{-1}mes(A_{G'}(F)\backslash G_{\delta}(F))\overline{S^{{\bf G}'}(\delta,\underline{f}^{{\bf G}'})}S^{{\bf G}'}(\delta,f^{{\bf G}'}).$$
Celle du membre de gauche vaut
$$ i(\gamma)^{-1}mes(A_{\tilde{G}}(F)\backslash G_{\gamma}(F))\sum_{\gamma'\in \dot{\cal X}(\gamma)}\overline{I^{\tilde{G}}(\gamma',\underline{f})}I^{\tilde{G}}(\gamma',f).$$
En utilisant la formule 4.9(5) qui exprime l'inverse du transfert et en se rappelant que $\vert \dot{\cal X}(\gamma)\vert=k(\gamma) $ on transforme cette expression en
$$d(\theta^*)^{-1}k(\gamma)^{-2}mes(A_{\tilde{G}}(F)\backslash Z_{G}(\gamma;F))\sum_{\gamma'\in \dot{\cal X}(\gamma)}\sum_{(\underline{{\bf G}'},\underline{\delta}),({\bf G}',\delta)\in \dot{\cal X}^{\cal E}(\gamma)}$$
$$\overline{\underline{\Delta}_{1}(\underline{\delta}_{1},\gamma')}^{-1}\Delta_{1}(\delta_{1},\gamma')^{-1}\overline{S^{\underline{\bf G}'}(\underline{\delta},\underline{f}^{\underline{\bf G}'})}S^{{\bf G}'}(\delta,f^{{\bf G}'}).$$
 Comme on le sait, la formule 4.9(5) exprime essentiellement une transformation de Fourier, les ensembles  $\dot{\cal X}(\gamma)$et $\dot{\cal X}^{\cal E}(\gamma)$ pouvant \^etre muni de structures de groupes ab\'eliens finis pour lesquelles ils sont duaux. La somme en $\gamma'$ des produits de facteurs de transfert vaut $\vert \dot{\cal X}(\gamma)\vert $, c'est-\`a-dire $k(\gamma)$, si $(\underline{{\bf G}'},\underline{\delta})=({\bf G}',\delta)$, $0$ sinon. On obtient
$$d(\theta^*)^{-1} k(\gamma)^{-1}mes(A_{\tilde{G}}(F)\backslash Z_{G}(\gamma;F)) \sum_{({\bf G}',\delta)\in \dot{\cal X}^{\cal E}(\gamma)}\overline{S^{{\bf G}'}(\delta,\underline{f}^{{\bf G}'})}S^{{\bf G}'}(\delta,f^{{\bf G}'}).$$
On veut prouver que cette expression est \'egale \`a (4). Il suffit de prouver que, pour tout $({\bf G}',\delta)\in \dot{\cal X}^{\cal E}(\gamma)$, on a l'\'egalit\'e
$$(5) \qquad c(\tilde{G},{\bf G}')=\vert Out({\bf G}')\vert ^{-1}k(\delta)mes(A_{G'}(F)\backslash G_{\delta}(F))^{-1}
d(\theta^*)^{-1} k(\gamma)^{-1}mes(A_{\tilde{G}}(F)\backslash Z_{G}(\gamma;F)) .$$
On note $c_{?}(\tilde{G},{\bf G}')$ le membre de droite de cette relation.   Notons $T$ le centralisateur de $G_{\gamma}$ dans $G$ et $T'=G_{\delta}$. On a $G_{\gamma}=T^{\theta,0}$,  $Z_{G}(\gamma)=T^{\theta}$ et $T'=T/(1-\theta)(T)$. De l'homomorphisme $\xi_{T,T'}$ se d\'eduit un homomorphisme 
$$a:A_{\tilde{G}}(F)\backslash T^{\theta}(F)\to A_{G'}(F)\backslash T'(F).$$
L'homomorphisme $\xi_{T,T'}:T^{\theta}(F)\to T'(F)$ conserve localement les mesures. Par contre, sa restriction $c:A_{\tilde{G}}(F)\to A_{G'}(F)$ ne les conserve pas. Notons $m' $ la mesure sur $A_{G'}(F)$ tel que $c$ conserve localement les mesures et $C$ la constante telle que notre mesure sur $A_{G'}(F)$ soit $Cm'$. On obtient alors
$$mes(A_{\tilde{G}}(F)\backslash T^{\theta}(F))=C\,mes(Im(a))\vert Ker(a)\vert .$$
On a aussi
$$mes(A_{G'}(F)\backslash T'(F))=mes(Im(a))\vert Coker(a)\vert .$$
D'o\`u
$$c_{?}(\tilde{G},{\bf G}')=C\vert Out({\bf G}')\vert ^{-1}d(\theta^*)^{-1}k(\delta)k(\gamma)^{-1}\vert Ker(a)\vert \vert Coker(a)\vert ^{-1}.$$
Consid\'erons le diagramme commutatif
$$\begin{array}{ccccccccc}1&\to&A_{\tilde{G}}(F)&\to&T^{\theta}(F)&\to&A_{\tilde{G}}(F)\backslash T^{\theta}(F)&\to&1\\&&\downarrow c&&\downarrow b&&\downarrow a&&\\1&\to&A_{G'}(F)&\to&T'(F)&\to&A_{G'}(F)\backslash T'(F)&\to&1\\ \end{array}$$
Ses lignes horizontales sont exactes. On en d\'eduit ais\'ement l'\'egalit\'e
$$\vert Ker(a)\vert \vert Coker(a)\vert ^{-1}=\vert Ker(c)\vert^{-1} \vert Coker(c)\vert \vert Ker(b)\vert \vert Coker(b)\vert^{-1} .$$
Montrons que

(6) $C=\vert Ker(c)\vert \vert Coker(c)\vert^{-1}$.

On peut identifier $A_{\tilde{G}}(F)$ \`a $A_{\tilde{G}}(F)^c\times Im(H_{A_{\tilde{G}}})$ et $A_{G'}(F)$ \`a $A_{G'}(F)^c\times Im(H_{A_{G'}})$ de sorte que $c$ se d\'ecompose conform\'ement en produit de deux homomorphismes. Le second homomorphisme est la restriction \`a $Im(H_{A_{\tilde{G}}})$ de l'isomorphisme de ${\cal A}_{\tilde{G}}$ sur ${\cal A}_{G'}$. D'apr\`es nos d\'efinitions,  il pr\'eserve les mesures (il s'agit des mesures de comptage dans le cas non-archim\'edien). Soit $V$ un ouvert compact de $Im(H_{A_{\tilde{G}}})$, posons $U=A_{\tilde{G}}(F)^c\times V$. 
 Si les mesures se correspondaient localement, on aurait l'\'egalit\'e $mes(c(U))=\vert Ker(c)\vert^{-1}  mes(U)$. Puisque ce n'est pas le cas, l'\'egalit\'e correcte est $mes(c(U))=C\vert Ker(c)\vert ^{-1}mes(U)$. On a $mes(U)=mes(A_{\tilde{G}}(F)^c)mes(V)$ et 
 $$mes(c(U))=mes(c(A_{\tilde{G}}(F)^c))mes(c(V))=[A_{G'}(F)^c:c(A_{\tilde{G}}(F)^c)]^{-1}mes(A_{G'}(F)^c)mes(V).$$
 On obtient
 $$C=\vert Ker(c)\vert  [A_{G'}(F)^c:c(A_{\tilde{G}}(F)^c)]^{-1}mes(A_{\tilde{G}}(F)^c)^{-1}mes(A_{G'}(F)^c).$$
Les mesures sur les groupes compacts sont d\'efinies de sorte que 
$$mes(A_{\tilde{G}}(F)^c)^{-1}mes(A_{G'}(F)^c)=[Im(H_{A_{G'}}):c(Im(H_{A_{\tilde{G}}}))]^{-1}.$$
On a aussi l'\'egalit\'e 
$$[Im(H_{A_{G'}}):c(Im(H_{A_{\tilde{G}}}))] [A_{G'}(F)^c:c(A_{\tilde{G}}(F)^c)]=\vert Coker(c)\vert .$$
Ces \'egalit\'es conduisent \`a (6).

Posons $V=(1-\theta)(T)$. Consid\'erons le diagramme commutatif
$$\begin{array}{cccccccc}(7)&&&&&1&&\\&&&&&\downarrow&&\\&&&&&V(F)&\stackrel{d}{\to}&V(F)\\&&&&&\downarrow&&\parallel\\&1&\to&T^{\theta}(F)&\to&T(F)&\stackrel{e}{\to}&V(F)\\&&&\parallel&&\downarrow f&&\\&&&T^{\theta}(F)&\stackrel{b}{\to}&T'(F)&&\\ \end{array}$$
o\`u $d$ et $e$ sont toutes deux \'egales \`a $1-\theta$. Les deuxi\`emes lignes horizontale et verticale sont exactes. On a $Ker(b)=T^{\theta}(F)\cap V(F)=Ker(d)$. On a aussi
$$\vert Coker(b)\vert =\vert Coker(f)\vert \vert T(F)/(T^{\theta}(F)V(F))\vert ,$$
$$\vert T(F)/(T^{\theta}(F)V(F))\vert=\vert e(T(F))/d(V(F))\vert =\vert Coker(d)\vert \vert Coker(e)\vert ^{-1}.$$
D'o\`u
$$\vert Ker(b)\vert \vert Coker(b)\vert^{-1} =\vert Ker(d)\vert \vert Coker(d)\vert^{-1} \vert Coker(e)\vert \vert Coker(f)\vert ^{-1}.$$
Consid\'erons un tore $D$ d\'efini sur $F$ et une isog\'enie $\varphi:D\to D$.   Notons ici $\varphi_{F}:D(F)\to D(F)$ l'homomorphisme qui s'en d\'eduit entre groupes de points sur $F$. Notons  $\mathfrak{d}$ l'alg\`ebre de Lie de $D$. On a
$$(8)\qquad \vert Ker(\varphi_{F})\vert \vert Coker(\varphi_{F})\vert ^{-1}=\vert X_{*}(D)^{\Gamma_{F}}/\varphi( X_{*}(D)^{\Gamma_{F}})\vert^{-1} \vert det(\varphi_{\vert \mathfrak{d}})\vert _{F}$$
$$=\vert det(\varphi_{\vert X_{*}(D)^{\Gamma_{F}}\otimes {\mathbb Q}})\vert^{-1}\vert det(\varphi_{\vert \mathfrak{d}})\vert _{F}.$$
Preuve de (8). Puisque $\varphi$ est injectif sur le ${\mathbb Z}$-module libre $X_{*}(D)^{\Gamma_{F}}$, on a l'\'egalit\'e
$$\vert X_{*}(D)^{\Gamma_{F}}/\varphi( X_{*}(D)^{\Gamma_{F}})\vert=\vert det(\varphi_{\vert X_{*}(D)^{\Gamma_{F}}\otimes {\mathbb Q}})\vert $$
et les deux derniers membres de (8) sont \'egaux.

Notons $D(F)^c$ le plus grand sous-groupe compact de $D(F)$ et $X=D(F)^c\backslash D(F)$. On utilise le diagramme commutatif:
$$\begin{array}{ccccccccc}1&\to&D(F)^c&\to&D(F)&\to&X&\to&1\\&&\downarrow\varphi_{F}^c&&\downarrow\varphi_{F}&&\downarrow\varphi_{X}&&\\1&\to&D(F)^c&\to&D(F)&\to&X&\to&1\\ \end{array}$$
Ses lignes \'etant exactes, on a
$$\vert Ker(\varphi_{F})\vert \vert Coker(\varphi_{F})\vert ^{-1}=\vert Ker(\varphi_{F}^c)\vert \vert Coker(\varphi_{F}^c)\vert ^{-1}\vert Ker(\varphi_{X})\vert \vert Coker(\varphi_{X})\vert ^{-1}.$$
Munissons $D(F)^c$ d'une mesure de Haar. On a
$$mes(D(F)^c)=\vert Coker(\varphi_{F}^c)\vert mes(Im(\varphi_{F}^c)),$$
$$mes(Im(\varphi_{F}^c))=j(\varphi_{F}^c)mes(D(F)^c)\vert Ker(\varphi_{F}^c)\vert ^{-1},$$
o\`u $j(\varphi_{F}^c)$ est le jacobien de $\varphi_{F}^c$. Si $F$ est non-archim\'edien, ce jacobien est la valeur absolue (au sens $\vert .\vert _{F}$) du d\'eterminant   de $\varphi$  agissant sur l'alg\`ebre de Lie de $D(F)^c$: $j(\varphi_{F}^c)=\vert det(\varphi_{\vert \mathfrak{d}})\vert _{F}$. Si $F$ est archim\'edien, le groupe $D(F)^c$ est un groupe de Lie r\'eel et $j(\varphi_{F}^c)$ est la valeur absolue r\'eelle du d\'eterminant de $\varphi$ agissant sur son alg\`ebre de Lie. Cette alg\`ebre de Lie est $\mathfrak{d}(F)/(X_{*}(D)^{\Gamma_{F}}\otimes {\mathbb R})$, d'o\`u 
$$j(\varphi_{F}^c)=\vert det(\varphi_{\vert \mathfrak{d}})\vert _{F}\vert det(\varphi_{\vert X_{*}(D)^{\Gamma_{F}}\otimes{\mathbb R}})\vert^{-1}=\vert det(\varphi_{\vert \mathfrak{d}})\vert _{F}\vert det(\varphi_{\vert X_{*}(D)^{\Gamma_{F}}\otimes{\mathbb Q}})\vert^{-1}.$$
 Si $F$ est archim\'edien, $X$ est un produit de groupes ${\mathbb R}_{+}^{\times}$ et $\varphi_{X}$ est bijectif. Si $F$ est non-archim\'edien, $\varphi_{X}$ est injectif et $\vert Coker(\varphi_{X})\vert =\vert det(\varphi_{\vert X\otimes {\mathbb Q}})\vert $. Fixons une uniformisante $\varpi_{F}$. L'application qui \`a $x_{*}\in X_{*}(D)^{\Gamma_{F}}$ associe l'image de $x_{*}(\varpi_{F})$ dans $X$ identifie $X_{*}(D)^{\Gamma_{F}}$ \`a un sous-groupe d'indice fini de $X$. Donc 
$$\vert det(\varphi_{\vert X\otimes {\mathbb Q}})\vert =\vert det(\varphi_{\vert X_{*}(D)^{\Gamma_{F}}\otimes {\mathbb Q}})\vert .$$
En mettant ces calculs bout \`a bout, on obtient (8) $\square$

On utilise (8) pour calculer
 $$\vert Ker(d)\vert \vert Coker(d)\vert ^{-1}=\vert det((1-\theta)_{\vert X_{*}(V)^{\Gamma_{F}}\otimes{\mathbb Q}})\vert ^{-1}\vert det((1-\theta)_{\vert (1-\theta)(\mathfrak{t})})\vert _{F}.$$
Le dernier terme n'est autre que $d(\theta^*)$.  En rassemblant les calculs pr\'ec\'edents, on obtient
$$(9) \qquad c_{?}(\tilde{G},{\bf G}')= \vert Out({\bf G}')\vert ^{-1}\vert det((1-\theta)_{\vert X_{*}(V)^{\Gamma_{F}}\otimes{\mathbb Q}})\vert ^{-1}$$
$$k(\delta)k(\gamma)^{-1}\vert Coker(e)\vert \vert Coker(f)\vert ^{-1}.$$
 
On consid\`ere la suite
$$H^1(\Gamma_{F};T^{\theta})\stackrel{g}{\to}H^1(\Gamma_{F};T)\stackrel{i}{\to}H^1(\Gamma_{F};G).$$
Supposons $F$ non archim\'edien. L'application qui \`a $y\in {\cal Y}(\gamma)$ (cf. 4.4) associe le cocycle $\sigma\mapsto y\sigma(y)^{-1}$ se quotiente en une bijection de $\dot{\cal Y}(\gamma)$ sur $Ker( i\circ g)$.  Donc
$$k(\gamma)=\vert \dot{\cal Y}(\gamma)\vert =\vert Ker(i\circ g)\vert.$$
De $g$ se d\'eduit une suite exacte
$$1\to Ker(g)\to Ker(i\circ g)\to Im(g)\cap Ker(i)\to 1.$$
Il est bien connu que $Ker(i)$ est \'egal \`a l'image de $j:H^1(\Gamma_{F};T_{sc})\to H^1(\Gamma_{F};T)$. La suite horizontale centrale de (7) se prolonge en une suite exacte de cohomologie
$$(10) \qquad T(F)\stackrel{e}{\to} V(F)\to H^1(\Gamma_{F};T^{\theta})\stackrel{g}{\to}H^1(\Gamma_{F};T)\stackrel{k}{\to}H^1(\Gamma_{F};V).$$
Donc $Im(g)=Ker(k)$ puis
$$k(\gamma)=\vert Ker(g)\vert \vert Ker(k)\cap Im(j)\vert .$$
Si $F={\mathbb R}$, parce que l'on travaille avec un $K$-espace, $\dot{\cal Y}(\gamma)$ s'identifie avec le sous-ensemble des \'el\'ements de $H^1(\Gamma_{{\mathbb R}};T^{\theta})$ dont l'image par $i\circ g$ appartient \`a l'image de l'application $H^1(\Gamma_{{\mathbb R}};G_{SC})\to H^1(\Gamma_{{\mathbb R}};G)$. La suite du calcul s'adapte et on obtient la m\^eme formule que ci-dessus. Revenons \`a $F$ quelconque. Consid\'erons la suite
$$H^1(\Gamma_{F} ;T_{sc})\stackrel{j}{\to}H^1(\Gamma_{F};T)\stackrel{k}{\to}H^1(\Gamma_{F};V).$$
Il s'en d\'eduit une suite exacte
$$1\to Ker(j)\to Ker(k\circ j)\to Ker(k)\cap Im(j)\to 1.$$
D'o\`u $\vert Ker(k)\cap Im(j)\vert =\vert Ker(k\circ j)\vert\vert Ker(j)\vert ^{-1}$ puis  
$$k(\gamma)=\vert Ker(g)\vert \vert Ker(k\circ j)\vert\vert Ker(j)\vert ^{-1}.$$
En utilisant la suite (10), on a
$$\vert Coker(e)\vert =\vert Ker(g)\vert .$$
La suite centrale verticale de (7) se prolonge elle-aussi en une suite exacte de cohomologie
$$T(F)\stackrel{f}{\to}T'(F)\to H^1(\Gamma_{F};V)\stackrel{l}{\to }H^1(\Gamma_{F};T).$$ 
 D'o\`u
$$\vert Coker(f)\vert =\vert Ker(l)\vert .$$
On obtient
$$(11) \qquad k(\gamma)^{-1}\vert Coker(e)\vert \vert Coker(f)\vert ^{-1}=\vert Ker(k\circ j)\vert^{-1}\vert Ker(j)\vert \vert Ker(l)\vert ^{-1},$$
o\`u on rappelle
$$j:H^1(\Gamma_{F};T_{sc})\to H^1(\Gamma_{F};T),\,\,k\circ j:H^1(\Gamma_{F},T_{sc})\to H^1(\Gamma_{F};V),\,\,l:H^1(\Gamma_{F};V)\to H^1(\Gamma_{F};T).$$
Tous ces groupes sont finis. On utilise l'\'egalit\'e
$$\vert Ker(j)\vert \vert H^1(\Gamma_{F};T)\vert =\vert Coker(j)\vert \vert H^1(\Gamma_{F};T_{sc})\vert $$
et les \'egalit\'es analogues pour $k\circ j$ et $l$.  On voit alors que, dans le membre de droite de (11), on peut remplacer les noyaux par les conoyaux. 

Le terme $k(\delta)$ se calcule comme $k(\gamma)$, le calcul \'etant beaucoup plus simple puisque la torsion est int\'erieure. On a $k(\delta)=\vert Im(m)\vert $, o\`u
$$m:H^1(\Gamma_{F};T'_{sc})\to H^1(\Gamma_{F};T').$$
On obtient
$$ k(\delta) k(\gamma)^{-1}\vert Coker(e)\vert \vert Coker(f)\vert ^{-1}=\vert Im(m)\vert \vert Coker(k\circ j)\vert^{-1}\vert Coker(j)\vert \vert Coker(l)\vert ^{-1}.$$

On utilise maintenant la dualit\'e. Par exemple $\pi_{0}(\hat{T}^{\Gamma_{F}})$ est le dual de $H^1(\Gamma_{F};T)$. On voit que $\vert Coker(j)\vert=\vert Ker(\hat{j})\vert $, o\`u $\hat{j}:\pi_{0}(\hat{T}^{\Gamma_{F}})\to \pi_{0}(\hat{T}_{ad}^{\Gamma_{F}})$ est dual de $j$. On calcule de m\^eme $\vert Coker(k\circ j)\vert $ et $\vert Coker(l)\vert $. On a aussi $\vert Im(m)\vert=\vert Im(\hat{m})\vert $ et la formule ci-dessus se transcrit en
$$(12) \qquad k(\delta) k(\gamma)^{-1}\vert Coker(e)\vert \vert Coker(f)\vert ^{-1}=\vert Im(\hat{m})\vert \vert Ker(\hat{j}\circ \hat{k})\vert^{-1}\vert Ker(\hat{j})\vert \vert Ker(\hat{l})\vert ^{-1}.$$
 Rappelons que $\hat{T}'=\hat{T}^{\hat{\theta},0}$. On a une suite exacte
$$\pi_{0}(Z(\hat{G}')^{\Gamma_{F}})\to\pi_{0}(\hat{T}^{\hat{\theta},0,\Gamma_{F}})\stackrel{\hat{m}}{\to} \pi_{0}((\hat{T}^{\hat{\theta},0}/Z(\hat{G}'))^{\Gamma_{F}}).$$
La donn\'ee ${\bf G}'$ est elliptique et $T'$ est un tore elliptique. Donc $Z(\hat{G}')^{\Gamma_{F},0}=\hat{T}^{\hat{\theta},\Gamma_{F},0}=Z(\hat{G})^{\hat{\theta},\Gamma_{F},0}$. La premi\`ere fl\`eche ci-dessus est injective, d'o\`u
$$(13)\qquad \vert Im(\hat{m})\vert =\vert \pi_{0}(\hat{T}^{\hat{\theta},0,\Gamma_{F}})\vert \vert \pi_{0}(Z(\hat{G}')^{\Gamma_{F}})\vert ^{-1}.$$
 Rappelons que $\hat{V}=\hat{T}/\hat{T}^{\hat{\theta},0}$. On a une suite exacte
$$\pi_{0}(\hat{T}^{\hat{\theta},0,\Gamma_{F}})\stackrel{\hat{n}}{\to}\pi_{0}(\hat{T}^{\Gamma_{F}})\stackrel{\hat{l}}{\to}\pi_{0}(\hat{V}^{\Gamma_{F}}).$$
Ici, la premi\`ere fl\`eche n'est pas injective. Son noyau est $(\hat{T}^{\hat{\theta},0}\cap\hat{T}^{\Gamma_{F},0})/\hat{T}^{\hat{\theta},\hat{\Gamma}_{F},0}=\pi_{0}(\hat{T}^{\hat{\theta},0}\cap \hat{T}^{\Gamma_{F},0})$. D'o\`u
$$(14)\qquad \vert Ker(\hat{l})\vert =\vert Coker(\hat{n})\vert =\vert \pi_{0}(\hat{T}^{\hat{\theta},0,\Gamma_{F}})\vert \vert \pi_{0}(\hat{T}^{\hat{\theta},0}\cap \hat{T}^{\Gamma_{F},0})\vert ^{-1}.$$
  De m\^eme, on a une suite exacte
  $$\pi_{0}(Z(\hat{G})^{\Gamma_{F}})\to \pi_{0}(\hat{T}^{\Gamma_{F}})\stackrel{\hat{j}}{\to}\pi_{0}(\hat{T}_{ad}^{\Gamma_{F}}).$$
  Le noyau de la premi\`ere fl\`eche est $\pi_{0}(Z(\hat{G})\cap\hat{T}^{\Gamma_{F},0})$ et on obtient
  $$(15)\qquad \vert Ker(\hat{j})\vert =\vert \pi_{0}(Z(\hat{G})^{\Gamma_{F}})\vert \vert \pi_{0}(Z(\hat{G})\cap\hat{T}^{\Gamma_{F},0})\vert ^{-1}.$$
  On a le diagramme commutatif
  $$\begin{array}{ccccc}\hat{T}/\hat{T}^{\hat{\theta},0}&& \to&&\hat{T}_{ad}\\ &\searrow&&\nearrow 1-\hat{\theta}&\\&&\hat{T}_{ad}/\hat{T}_{ad}^{\hat{\theta}}&&\\ \end{array}$$
  d'o\`u la factorisation
  $$\begin{array}{ccccc}\pi_{0}((\hat{T}/\hat{T}^{\hat{\theta},0})^{\Gamma_{F}})&&\stackrel{\hat{j}\circ\hat{k}}{\to}&&\pi_{0}(\hat{T}_{ad}^{\Gamma_{F}})\\ &\searrow&&\nearrow 1-\hat{\theta}&\\&&\pi_{0}((\hat{T}_{ad}/\hat{T}_{ad}^{\hat{\theta}})^{\Gamma_{F}})&&\\ \end{array}$$ 
  Puisque $\tilde{T}$ est elliptique, on a $X_{*}(\hat{T}_{ad})^{\Gamma_{F},\hat{\theta}}=0$. Donc $X_{*}(\hat{T}_{ad})^{\Gamma_{F}}\otimes{\mathbb Q}=(1-\hat{\theta})(X_{*}(\hat{T}_{ad})^{\Gamma_{F}})\otimes{\mathbb Q}$ puis $\hat{T}_{ad}^{\Gamma_{F},0}=(1-\hat{\theta})(\hat{T}_{ad}^{\Gamma_{F},0})$. Il en r\'esulte que l'homomorphisme
  $$1-\hat{\theta}:\pi_{0}((\hat{T}_{ad}/\hat{T}_{ad}^{\hat{\theta}})^{\Gamma_{F}})\to\pi_{0}(\hat{T}_{ad}^{\Gamma_{F}})$$ 
 est injectif. Le noyau de $\hat{j}\circ\hat{k}$ est donc \'egal \`a celui de l'homomorphisme
  $$\pi_{0}((\hat{T}/\hat{T}^{\hat{\theta},0})^{\Gamma_{F}})\to\pi_{0}((\hat{T}_{ad}/\hat{T}_{ad}^{\hat{\theta}})^{\Gamma_{F}}).$$
  Ce dernier se compl\`ete en la suite exacte
  $$ \pi_{0}((Z(\hat{G})/(Z(\hat{G})\cap\hat{T}^{\hat{\theta},0})^{\Gamma_{F}})\stackrel{\hat{p}}{\to}\pi_{0}((\hat{T}/\hat{T}^{\hat{\theta},0})^{\Gamma_{F}})\to\pi_{0}((\hat{T}_{ad}/\hat{T}_{ad}^{\hat{\theta}})^{\Gamma_{F}}).$$
  D'o\`u:
$$(16)\qquad \vert Ker(\hat{j}\circ\hat{k})\vert =\vert  \pi_{0}((Z(\hat{G})/(Z(\hat{G})\cap\hat{T}^{\hat{\theta},0}))^{\Gamma_{F}})\vert \vert Ker(\hat{p})\vert ^{-1}.$$
On calcule
$$Ker(\hat{p})=(Z(\hat{G})\cap \hat{T}^{\Gamma_{F},0}\hat{T}^{\hat{\theta},0})/(Z(\hat{G})^{\Gamma_{F},0}(Z(\hat{G})\cap\hat{T}^{\hat{\theta},0})).$$
La suite suivante est exacte:
$$1\to (\hat{T}^{\hat{\theta},0}\cap\hat{T}^{\Gamma_{F},0})/(Z(\hat{G})^{\Gamma_{F},0}\cap\hat{T}^{\hat{\theta},0})\to ((Z(\hat{G})\hat{T}^{\hat{\theta},0})\cap \hat{T}^{\Gamma_{F},0})/Z(\hat{G})^{\Gamma_{F},0}\to Ker(\hat{p})\to 1$$
Le premier terme de cette suite a pour nombre d'\'el\'ements
$$\vert \pi_{0}(\hat{T}^{\hat{\theta},0}\cap\hat{T}^{\Gamma_{F},0})\vert \vert \pi_{0}(Z(\hat{G})^{\Gamma_{F},0}\cap\hat{T}^{\hat{\theta},0})\vert ^{-1}.$$
Le deuxi\`eme terme s'ins\`ere dans la suite exacte
$$1\to (Z(\hat{G})\cap\hat{T}^{\Gamma_{F},0})/Z(\hat{G})^{\Gamma_{F},0}\to((Z(\hat{G})\hat{T}^{\hat{\theta},0})\cap \hat{T}^{\Gamma_{F},0})/Z(\hat{G})^{\Gamma_{F},0}$$
$$\to ((Z(\hat{G})\hat{T}^{\hat{\theta},0})\cap \hat{T}^{\Gamma_{F},0})/(Z(\hat{G})\cap\hat{T}^{\Gamma_{F},0})\to 1.$$
Le premier terme de cette suite n'est autre que $\pi_{0} (Z(\hat{G})\cap\hat{T}^{\Gamma_{F},0})$. On voit que le second n'est autre que $\hat{T}_{ad}^{\Gamma_{F},0,\hat{\theta}}$. Ce dernier groupe est fini puisque $\tilde{T}$ est elliptique. A ce point, on obtient
$$(17) \qquad \vert Ker(\hat{p})\vert =\vert \pi_{0}(\hat{T}^{\hat{\theta},0}\cap\hat{T}^{\Gamma_{F},0})\vert^{-1}\vert \pi_{0}(Z(\hat{G})^{\Gamma_{F},0}\cap\hat{T}^{\hat{\theta},0})\vert\vert \pi_{0} (Z(\hat{G})\cap\hat{T}^{\Gamma_{F},0})\vert \vert \hat{T}_{ad}^{\Gamma_{F},0,\hat{\theta}}\vert .$$
Puisque l'endomorphisme  $1-\hat{\theta}$ de $\hat{T}_{ad}^{\Gamma_{F},0}$ est une isog\'enie, son noyau $\hat{T}_{ad}^{\Gamma_{F},0,\hat{\theta}}$ a pour nombre d'\'el\'ements la valeur absolue du d\'eterminant de $1-\hat{\theta}$ agissant sur $X_{*}(\hat{T}_{ad})^{\Gamma_{F}}\otimes {\mathbb Q}$. Par dualit\'e, c'est aussi la valeur absolue du d\'eterminant de $1-\theta$ agissant sur $X_{*}(T_{sc})^{\Gamma_{F}}\otimes{\mathbb Q}$. On utilise les \'egalit\'es
$$X_{*}(V)^{\Gamma_{F}}\otimes {\mathbb Q}=(1-\theta)(X_{*}(T)^{\Gamma_{F}})\otimes{\mathbb Q}$$ 
$$=\left((1-\theta)(X_{*}(T_{sc})^{\Gamma_{F}})\otimes{\mathbb Q}\right)\oplus\left((1-\theta)(X_{*}(Z(G)^0)^{\Gamma_{F}})\otimes{\mathbb Q}\right).$$
 On a $(1-\theta)(X_{*}(T_{sc})^{\Gamma_{F}})\otimes{\mathbb Q}=X_{*}(T_{sc})^{\Gamma_{F}}\otimes{\mathbb Q}$ toujours parce que $\tilde{T}$ est elliptique. D'o\`u
 $$\vert \hat{T}_{ad}^{\Gamma_{F},0,\hat{\theta}}\vert =\vert det((1-\theta)_{\vert X_{*}(V)^{\Gamma_{F}}\otimes {\mathbb Q} })\vert \vert det((1-\theta)_{\vert (1-\theta)(X_{*}(Z(G)^0)^{\Gamma_{F}})\otimes{\mathbb Q}})\vert ^{-1}.$$
 Pour calculer ce dernier d\'eterminant, on peut  remplacer $(1-\theta)(X_{*}(Z(G)^0)^{\Gamma_{F}})\otimes{\mathbb Q}$ par $(1-\theta)(X_{*}(Z(G)^0)^{\Gamma_{F}})\otimes{\mathbb R}$. Cet espace est isomorphe \`a ${\cal A}_{G}/{\cal A}_{\tilde{G}}$. On obtient alors
 $$(18)\qquad  \vert \hat{T}_{ad}^{\Gamma_{F},0,\hat{\theta}}\vert =\vert det((1-\theta)_{\vert X_{*}(V)^{\Gamma_{F}}\otimes {\mathbb Q} })\vert \vert det((1-\theta)_{\vert {\cal A}_{G}/{\cal A}_{\tilde{G}}})\vert ^{-1}.$$
 Rassemblons les formules (9) et (12), (13),...,(18). On obtient
 $$c_{?}(\tilde{G},{\bf G}')= \vert det((1-\theta)_{\vert {\cal A}_{G}/{\cal A}_{\tilde{G}}})\vert ^{-1}
\vert \pi_{0}(Z(\hat{G})^{\Gamma_{F}})\vert \vert Z(\hat{G}')^{\Gamma_{F}})\vert ^{-1}\vert $$
$$\vert Out({\bf G}')\vert ^{-1}\vert \pi_{0}(Z(\hat{G})^{\Gamma_{F},0}\cap T^{\hat{\theta},0})\vert \vert \pi_{0}((Z(\hat{G})/(Z(\hat{G})\cap T^{\hat{\theta},0}))^{\Gamma_{F}})\vert ^{-1}.$$
 On a l'\'egalit\'e $Z(\hat{G})\cap T^{\hat{\theta},0}=Z(\hat{G})\cap \hat{G}'$. La formule ci-dessus est alors celle qui d\'efinit $c(\tilde{G},{\bf G}')$. Cela d\'emontre l'\'egalit\'e (5), ce qui ach\`eve la preuve. $\square$

 \bigskip
 
 \section{Distributions "g\'eom\'etriques"}

\bigskip

\subsection{Distributions "g\'eom\'etriques" dans le cas non-archim\'edien}
  On suppose $F$ non archim\'edien. On note $D_{g\acute{e}om}(\tilde{G}(F),\omega)$ l'espace des formes lin\'eaires sur $C_{c}^{\infty}(\tilde{G}(F))$ qui se factorisent en une forme lin\'eaire sur $I(\tilde{G}(F),\omega)$ et qui sont support\'ees par une r\'eunion finie de classes de conjugaison par $G(F)$. On a d\'ej\`a construit de telles formes lin\'eaires en 2.4: l'int\'egrale orbitale $f\mapsto I^{\tilde{G}}(\gamma,\omega,f)$ associ\'ee \`a un \'el\'ement $\gamma\in \tilde{G}(F)$ et aux choix de mesures sur $G(F)$ et $G_{\gamma}(F)$.   On se d\'ebarrasse du choix de la mesure sur $G(F)$ en consid\'erant cette forme lin\'eaire comme d\'efinie sur $C_{c}^{\infty}(\tilde{G}(F))\otimes Mes(G(F))$. On obtient donc un \'el\'ement de $D_{g\acute{e}om}(\tilde{G}(F),\omega)\otimes Mes(G(F))^*$. Il est commode de noter tout \'el\'ement de cet espace comme une int\'egrale orbitale. C'est-\`a-dire que, pour $\boldsymbol{\gamma}\in D_{g\acute{e}om}(\tilde{G}(F),\omega)\otimes Mes(G(F))^*$ et ${\bf f}\in C_{c}^{\infty}(\tilde{G}(F))\otimes Mes(G(F))$, on notera $I^{\tilde{G}}(\boldsymbol{\gamma},{\bf f})$ la valeur de $\boldsymbol{\gamma}$ sur ${\bf f}$. On utilisera diff\'erentes variantes de cette notation (pour les int\'egrales orbitales stables par exemple).

  Si ${\cal O}$ est une  r\'eunion finie de classes de conjugaison (par $G(F)$) semi-simples, on note $D_{g\acute{e}om}({\cal O},\omega)$ le sous-espace de ces distributions \`a support dans $\{\gamma\in \tilde{G}; \gamma_{ss}\in{\cal O}\}$, o\`u $\gamma_{ss}$ est la partie semi-simple de $\gamma$.  Notons  qu'un tel sous-espace  peut \^etre nul, \`a cause du caract\`ere $\omega$.
Plus concr\`etement, notons $I(\tilde{G}(F),\omega)_{{\cal O},0}$ le sous-espace des $f\in I(\tilde{G}(F),\omega)$ pour lesquels il existe un voisinage $\tilde{V}$ de ${\cal O}$ invariant par conjugaison par $G(F)$ tel que $I^{\tilde{G}}(\gamma,\omega,f)=0$ pour tout $\gamma\in \tilde{V}\cap \tilde{G}_{reg}(F)$. Posons $I(\tilde{G}(F),\omega)_{{\cal O},loc}=I(\tilde{G}(F),\omega)/I(\tilde{G}(F),\omega)_{{\cal O},0}$. La projection naturelle
  $$C_{c}^{\infty}(\tilde{G}(F))\to I(\tilde{G}(F),\omega)_{{\cal O},loc}$$
  est surjective et on a
  
  (1) $D_{g\acute{e}om}({\cal O},\omega)$ est l'espace des formes lin\'eaires sur $C_{c}^{\infty}(\tilde{G}(F))$ qui se factorisent par cette projection.
  
  Preuve. Notons $C_{c}^{\infty}(\tilde{G}(F))_{{\cal O},0}$ le sous-espace des \'el\'ements $ C_{c}^{\infty}(\tilde{G}(F))$ dont le support ne contient pas d'\'el\'ement de partie semi-simple dans ${\cal O}$. Par d\'efinition, $D_{g\acute{e}om}({\cal O},\omega)$ est l'espace des formes lin\'eaires sur $C_{c}^{\infty}(\tilde{G}(F))$ qui annulent $C_{c}^{\infty}(\tilde{G}(F))_{{\cal O},0}$ et qui se factorisent par $I(\tilde{G}(F),\omega)$. Il suffit donc de prouver que l'image de $C_{c}^{\infty}(\tilde{G}(F))_{{\cal O},0}$ dans $I(\tilde{G}(F),\omega)$ est \'egale \`a $I(\tilde{G}(F),\omega)_{{\cal O},0}$. Il est clair que cette image est contenue dans $I(\tilde{G}(F),\omega)_{{\cal O},0}$. Inversement, soit $f\in C_{c}^{\infty}(\tilde{G}(F))$ dont l'image dans $I(\tilde{G}(F),\omega)$ appartienne \`a ce sous-espace. On choisit un voisinage $\tilde{V}$ de ${\cal O}$ invariant par conjugaison tel que $I^{\tilde{G}}(\gamma,\omega,f)=0$ pour tout $\gamma\in \tilde{V}\cap \tilde{G}_{reg}(F)$. On peut supposer $\tilde{V}$ ouvert et ferm\'e. Notons ${\bf 1}_{\tilde{V}}$ sa fonction caract\'eristique. On a $f=f{\bf 1}_{\tilde{V}}+f(1-{\bf 1}_{\tilde{V}})$. Toutes les int\'egrales orbitales fortement r\'eguli\`eres de la fonction $f{\bf 1}_{\tilde{V}}$ sont nulles. Cela entra\^{\i}ne que l'image de cette fonction dans $I(\tilde{G}(F),\omega)$ est nulle.  La deuxi\`eme fonction $f(1-{\bf 1}_{\tilde{V}})$ appartient \`a $C_{c}^{\infty}(\tilde{G}(F))_{{\cal O},0}$. $\square$

   D'apr\`es (1),  $D_{g\acute{e}om}({\cal O},\omega)$  s'identifie au dual de $ I(\tilde{G}(F),\omega)_{{\cal O},loc}$.Il est bien connu que tout \'el\'ement de $D_{g\acute{e}om}(\tilde{G}(F),\omega)$ est combinaison lin\'eaire d'int\'egrales orbitales. Cela entra\^{\i}ne que $D_{g\acute{e}om}(\tilde{G}(F),\omega)$ est la somme directe de ses sous-espaces $D_{g\acute{e}om}({\cal O},\omega)$ quand ${\cal O}$ d\'ecrit  les classes de conjugaison semi-simples.   
   
   Soit $\tilde{M}$ un espace de Levi de $\tilde{G}$. Dualement \`a l'application
   $$\begin{array}{ccc}I(\tilde{G}(F),\omega)\otimes Mes(G(F))&\to&I(\tilde{M}(F),\omega)\otimes Mes(M(F))\\ {\bf f}&\mapsto &{\bf f}_{\tilde{M},\omega},\\ \end{array}$$
   on a un homomorphisme d'induction
   $$\begin{array}{ccc}D_{g\acute{e}om}(\tilde{M}(F),\omega)\otimes Mes(M(F))^*&\to &D_{g\acute{e}om}(\tilde{G}(F),\omega)\otimes Mes(G(F))^*\\ \boldsymbol{\gamma}&\mapsto&\boldsymbol{\gamma}^{\tilde{G}}\\ \end{array}$$
   Soit ${\cal O}$ une classe de conjugaison semi-simple contenant un \'el\'ement $\gamma$ tel que $\gamma\in \tilde{M}(F)$ et $G_{\gamma}\subset M$. Alors $D_{g\acute{e}om}({\cal O},\omega)$ est contenu dans l'image de cet homomorphisme d'induction.

   \bigskip
\subsection{Distributions "g\'eom\'etriques" dans le cas archim\'edien}
  On suppose $F={\mathbb R}$ ou ${\mathbb C}$. On munit $C_{c}^{\infty}(\tilde{G}(F))$ d'une topologie de la fa\c{c}on suivante. Notons $\mathfrak{U}(G)$ l'alg\`ebre enveloppante de l'alg\`ebre de Lie de $G$. Cette alg\`ebre agit sur $C_{c}^{\infty}(\tilde{G}(F))$ de deux fa\c{c}ons: via les translations \`a gauche ou \`a droite. Consid\'erons par exemple l'action via les translations \`a gauche. Pour $Y\in \mathfrak{U}(G)$, on d\'efinit la semi-norme $\nu_{Y}$ sur $C_{c}^{\infty}(\tilde{G}(F))$ par $\nu_{Y}(f)=sup\{\vert (Yf)(\gamma)\vert ; \gamma\in \tilde{G}(F)\}$. Soit $\tilde{H}$ un sous-ensemble compact de $\tilde{G}(F)$. Notons $C_{c}^{\infty}(\tilde{H})$ le sous-espace des \'el\'ements de $C_{c}^{\infty}(\tilde{G}(F))$ \`a support dans $\tilde{H}$.   On munit ce sous-espace de la topologie d\'efinie par les semi-normes $\nu_{Y}$ pour $Y\in \mathfrak{U}(G)$. L'espace $C_{c}^{\infty}(\tilde{G}(F))$ est     limite inductive des $C_{c}^{\infty}(\tilde{H})$ quand $\tilde{H}$ d\'ecrit les sous-ensembles compacts de $\tilde{G}(F)$ et on le munit de la topologie limite inductive des topologies sur ces sous-espaces. On appelle distribution sur $\tilde{G}(F)$ une forme lin\'eaire continue sur $C_{c}^{\infty}(\tilde{G}(F))$. Une distribution $\omega$-\'equivariante est une distribution qui se factorise par $I(\tilde{G}(F),\omega)$. En imitant Bouaziz, on munit l'espace $I(\tilde{G}(F),\omega)$ d'une topologie de la fa\c{c}on suivante. Fixons un ensemble  ${\tilde{\cal T}}$ de repr\'esentants des classes de conjugaison par $G(F)$ de tores tordus maximaux  $\tilde{T}$ tels que $\omega$ soit trivial sur $T^{\theta}(F)$. Un tel ensemble est fini. En fixant des mesures sur $G(F)$ et sur $T(F)$ pour tout $\tilde{T}\in {\tilde{\cal T}}$, on peut consid\'erer $I(\tilde{G}(F),\omega)$ comme un espace de familles $\varphi_{{\tilde{\cal T}}}=(\varphi_{\tilde{T}})_{\tilde{T}\in {\tilde{\cal T}}}$ o\`u $\varphi_{\tilde{T}}$ est une fonction $C^{\infty}$ sur $\tilde{T}(F)\cap \tilde{G}_{reg}(F)$ (l'int\'egrale orbitale sur $\tilde{T}(F)$). Dans la suite, on consid\'erera $I(\tilde{G}(F),\omega)$ soit comme un quotient de $C_{c}^{\infty}(\tilde{G}(F))$ (ses \'el\'ements seront alors not\'es $f$), soit comme un espace de telles familles (ses \'el\'ements seront alors not\'es $\varphi_{{\tilde{\cal T}}}$). On pose $\mathfrak{U}_{{\tilde{\cal T}}}=\prod_{\tilde{T}\in {\tilde{\cal T}}}\mathfrak{U}(T)$. Pour une famille $Y_{{\tilde{\cal T}}}=(Y_{\tilde{T}})_{\tilde{T}\in {\tilde{\cal T}}}\in \mathfrak{U}_{{\tilde{\cal T}}}$, on d\'efinit la semi-norme
  $$\nu_{Y_{{\tilde{\cal T}}}}(\varphi_{{\tilde{\cal T}}})=sup\{\vert (Y_{\tilde{T}}\varphi_{\tilde{T}})(\gamma)\vert ; \gamma\in \tilde{T}(F)\cap\tilde{G}_{reg}(F), \tilde{T}\in {\tilde{\cal T}}\}.$$
  Elle est bien d\'efinie c'est-\`a-dire que ce $sup$ est fini pour les \'el\'ements de $I(\tilde{G}(F))$. C'est un r\'esultat profond d'Harish-Chandra (sa g\'en\'eralisation au cas tordu par descente est imm\'ediate). Soit $\tilde{H}_{{\tilde{\cal T}}}=(\tilde{H}_{\tilde{T}})_{\tilde{T}\in {\tilde{\cal T}}}$ une famille telle que  pour tout $\tilde{T}$, $\tilde{H}_{\tilde{T}}$ est un sous-ensemble compact de $\tilde{T}(F)$. On note $I(\tilde{H}_{{\tilde{\cal T}}},\omega)$ le sous-espace des \'el\'ements $\varphi_{{\tilde{\cal T}}}=(\varphi_{\tilde{T}})_{\tilde{T}\in {\tilde{\cal T}}}\in I(\tilde{G}(F),\omega)$ tels que pour tout $\tilde{T}$, $\varphi_{\tilde{T}}$ est \`a support dans $\tilde{H}_{\tilde{T}}$. On munit ce sous-espace de la topologie d\'efinie par les semi-normes $\nu_{Y_{{\tilde{\cal T}}}}$ pour   $Y_{{\tilde{\cal T}}}\in \mathfrak{U}_{{\tilde{\cal T}}}$. Cela le munit   d'une topologie d'espace de Fr\'echet, c'est-\`a-dire que $I(\tilde{H}_{{\tilde{\cal T}}},\omega)$ est complet: les conditions de saut qui d\'efinissent l'espace des int\'egrales orbitales sont des conditions ferm\'ees.  On munit $I(\tilde{G}(F),\omega)$ de  la topologie limite inductive de celle sur les sous-espaces $I(\tilde{H}_{{\tilde{\cal T}}},\omega)$. On a
  
  (1)  l'homomorphisme $C_{c}^{\infty}(\tilde{G}(F))\to I(\tilde{G}(F),\omega)$ est une surjection continue et ouverte. 
  
  Cf. [R1] th\'eor\`eme 9.4. Renard suppose $\omega=1$ mais, ici encore, la preuve se g\'en\'eralise au cas $\omega$ quelconque.

 D'apr\`es (1), l'espace des distributions $\omega$-\'equivariantes s'identifie \`a celui des formes lin\'eaires continues sur $I(\tilde{G}(F),\omega)$. On note $D_{g\acute{e}om}(\tilde{G}(F),\omega)$ l'espace des distributions $\omega$-\'equivariantes qui sont support\'ees par un nombre fini de classes de conjugaison par $G(F)$.  Concr\`etement, consid\'erons un tore tordu $\tilde{T}\in \tilde{{\cal T}}$ et un \'el\'ement $\eta\in \tilde{T}(F)$. Fixons une composante connexe $\Omega$ de $\mathfrak{t}^{\theta}(F)\cap \mathfrak{g}_{\eta,reg}(F)$ et un op\'erateur diff\'erentiel $D$ sur $\mathfrak{t}^{\theta}(F)$ \`a coefficients constants. Pour $\varphi_{\tilde{{\cal T}}}=(\varphi_{\tilde{T}'})_{\tilde{T}'\in \tilde{{\cal T}}}\in I(\tilde{G}(F),\omega)$, la fonction $X\mapsto D\varphi_{\tilde{T}}(exp(X)\eta)$ est $C^{\infty}$ sur  $\Omega$ et a une limite quand $X$ tend vers $0$ dans $\Omega$. Notons $\boldsymbol{\gamma}_{\eta, \tilde{T},\Omega,D}(\varphi_{\tilde{{\cal T}}})$ cette limite. La forme lin\'eaire $\boldsymbol{\gamma}_{\eta, \tilde{T},\Omega,D}$ ainsi construite appartient \`a $D_{g\acute{e}om}(\tilde{G}(F),\omega)$ et cet espace est engendr\'e lin\'eairement par de telles formes lin\'eaires.

 Si ${\cal O}$ est une r\'eunion finie de classes de conjugaison (par $G(F)$) semi-simples, on d\'efinit  le sous-espace $D_{g\acute{e}om}({\cal O},\omega) $ comme dans le cas non-archim\'edien. Notons $I(\tilde{G}(F),\omega)_{{\cal O},0}$  le sous-espace des $f\in I(\tilde{G}(F),\omega)$ pour lesquels il existe un voisinage $\tilde{V}$ de ${\cal O}$ invariant par conjugaison par $G(F)$ tel que $I^{\tilde{G}}(\gamma,\omega,f)=0$ pour tout $\gamma\in \tilde{V}\cap \tilde{G}_{reg}(F)$.   Notons $C\ell I(\tilde{G}(F),\omega)_{{\cal O},0}$ sa cl\^oture dans $I(\tilde{G}(F))$. C'est le sous-espace des $\varphi_{{\tilde{\cal T}}}\in I(\tilde{G}(F),\omega)$ v\'erifiant la condition suivante. Soient $\tilde{T}\in {\tilde{\cal T}}$, $\eta\in \tilde{T}(F)\cap {\cal O}$ et $Y\in \mathfrak{U}(T)$. Alors la fonction $Y\varphi_{\tilde{T}}$ bien d\'efinie sur $\tilde{T}(F)\cap \tilde{G}_{reg}(F)$ a une limite nulle en $\eta$. On pose $I(\tilde{G}(F),\omega)_{{\cal O},loc}=I(\tilde{G}(F),\omega)/C\ell I(\tilde{G}(F),\omega)_{{\cal O},0}$ et on munit cet espace de la topologie quotient. Il y a un homomorphisme surjectif, continu et ouvert
  $$C_{c}^{\infty}(\tilde{G}(F))\to I(\tilde{G}(F),\omega)_{{\cal O},loc}.$$
  On a
  
  (2)   $D_{g\acute{e}om}({\cal O},\omega) $ est l'image par l'homomorphisme dual de l'espace des formes lin\'eaires continues sur $I(\tilde{G}(F),\omega)_{{\cal O},loc}$.
  
  Preuve. On note $C_{c}^{\infty}(\tilde{G}(F))_{{\cal O},0}$ le sous-espace des $f\in C_{c}^{\infty}(\tilde{G}(F))$ dont le support ne contient pas d'\'el\'ement de partie semi-simple dans ${\cal O}$. Son image dans $I(\tilde{G}(F),\omega)$ est \'evidemment contenue dans $I(\tilde{G}(F),\omega)_{{\cal O},0}$. En fait, cette image est \'egale \`a $I(\tilde{G}(F),\omega)_{{\cal O},0}$. La preuve est essentiellement la m\^eme que celle de 5.1(1). Il suffit d'y remplacer la fonction ${\bf 1}_{\tilde{V}}$ par une fonction $C^{\infty}$, invariante par conjugaison, \`a support dans $\tilde{V}$ et valant $1$ au voisinage des \'el\'ements de partie semi-simple dans ${\cal O}$. D'apr\`es (1) et la d\'efinition,  $D_{g\acute{e}om}({\cal O},\omega) $ est l'espace des formes lin\'eaires continues sur $I(\tilde{G}(F),\omega)$ qui annulent l'image de $C_{c}^{\infty}(\tilde{G}(F))_{{\cal O},0}$. Autrement dit qui annulent $I(\tilde{G}(F),\omega)_{{\cal O},0}$. Puisqu'il s'agit de formes continues, cela \'equivaut \`a annuler $C\ell I(\tilde{G}(F),\omega)_{{\cal O},0}$ ou encore \`a se factoriser en une forme lin\'eaire continue sur $I(\tilde{G}(F),\omega)_{{\cal O},loc}$. $\square$
  
  Remarquons que si $\tilde{M}$ est un espace de Levi de $\tilde{G}$, l'homomorphisme $f\mapsto f_{\tilde{M},\omega}$ de $I(\tilde{G}(F),\omega)$ dans $I(\tilde{M}(F),\omega)$  se descend en un  homomorphisme de $I(\tilde{G}(F),\omega)_{{\cal O},loc}$ dans $I(\tilde{M}(F),\omega)_{{\cal O}_{\tilde{M}},loc}$ o\`u ${\cal O}_{\tilde{M}}=\tilde{M}(F)\cap {\cal O}$. Il y a deux fa\c{c}ons naturelles de d\'efinir un sous-espace $I_{cusp}(\tilde{G}(F),\omega)_{{\cal O},loc}\subset I(\tilde{G}(F),\omega)_{{\cal O},loc}$: soit comme image par localisation de $I_{cusp}(\tilde{G}(F),\omega)$, soit comme le sous-espace de $I(\tilde{G}(F),\omega)_{{\cal O},loc}$ annul\'e par les homomorphismes $f\mapsto f_{\tilde{M},\omega}$ pour tout $\tilde{M}$ propre. On a
  
  (3) ces deux d\'efinitions co\"{\i}ncident.
 
  Preuve. Supposons $F={\mathbb R}$. La premi\`ere d\'efinition donne \'evidemment un sous-espace de l'espace d\'efini par la seconde. Soit $\varphi_{{\tilde{\cal T}}}\in I(\tilde{G}({\mathbb R}),\omega)$ un \'el\'ement dont l'image par localisation appartient \`a ce dernier espace.   Fixons un \'el\'ement elliptique $\tilde{T}\in {\tilde{\cal T}}$ et un \'el\'ement $\eta\in \tilde{T}({\mathbb R})\cap {\cal O}$.  Comme en 4.13(3), consid\'erons la fonction
  $$(4) \qquad X\mapsto \Delta_{\eta}(X)\varphi_{\tilde{T}}(exp(X)\eta)$$
  au voisinage de $0$ dans  $ \mathfrak{t}^{\theta}({\mathbb R})\cap \mathfrak{g}_{\eta,reg}({\mathbb R})$. Soit $\Omega$ une composante connexe de cet ensemble, contenant $\eta$ dans son adh\'erence. La fonction ci-dessus est $C^{\infty}$ sur  $\Omega$ et on sait  qu'elle se prolonge en une fonction $C^{\infty}$ dans un voisinage de $\Omega$ (cf. [Boua] remarque 3.2). Notons $\phi_{\tilde{T},\Omega}$ un tel prolongement. L'hypoth\`ese de cuspidalit\'e sur $\varphi_{{\tilde{\cal T}}}$ implique que le d\'eveloppement infinit\'esimal au voisinage de $\eta$ de $\phi_{\tilde{T},\Omega}$ ne d\'epend pas de $\Omega$. C'est-\`a-dire que, pour tout $Y\in \mathfrak{U}(T^{\theta,0})$, $(Y\phi_{\tilde{T},\Omega})(\eta)$ est ind\'ependant de $\Omega$.  Consid\'erons le normalisateur de $T^{\theta,0}$ dans $Z_{G}(\eta,{\mathbb R})$ et son quotient fini $W_{\eta}(T^{\theta,0})$ par $T^{\theta,0}({\mathbb R})$. Ce quotient agit sur les fonctions sur $\mathfrak{t}^{\theta}({\mathbb R})$. La fonction $\Delta_{\eta}$ se transforme selon un certain caract\`ere $\chi$ de ce groupe. Parce que les  int\'egrales orbitales sont invariantes par ce groupe, la fonction (4) se transforme selon le m\^eme caract\`ere $\chi$. Donc le d\'eveloppement infinit\'esimal commun des fonctions $\phi_{\tilde{T},\Omega}$ se transforme lui-aussi selon le caract\`ere $\chi$. Fixons $\Omega$ et introduisons la fonction
  $$\phi'_{\tilde{T}}=\vert W_{\eta}(T^{\theta,0})\vert ^{-1}\sum_{w\in W_{\eta}(T^{\theta,0})}\chi(w)^{-1}w\phi_{\tilde{T},\Omega}.$$
  Elle a m\^eme d\'eveloppement infinit\'esimal que nos fonctions $\phi_{\tilde{T},\Omega}$. Il existe une fonction $\varphi'_{\tilde{T}}$ sur $\tilde{T}({\mathbb R})\cap \tilde{G}_{reg}(F)$ v\'erifiant la condition 4.13(2) et telle que la fonction
  $$X\mapsto \Delta_{\eta}(X)\varphi'_{\tilde{T}}(exp(X)\eta)$$
  co\"{\i}ncide avec $\phi'_{\tilde{T}}$ au voisinage de $0$ dans   $ \mathfrak{t}^{\theta}({\mathbb R})\cap \mathfrak{g}_{\eta,reg}({\mathbb R})$. Quitte \`a multiplier cette fonction par une fonction $C^{\infty}$ invariante par conjugaison et \`a support concentr\'e dans un voisinage invariant de $\eta$, on peut supposer que $\varphi'_{\tilde{T}}$ v\'erifie la condition 4.13(3). Donc cette fonction, prolong\'ee par $0$ sur les autres \'el\'ements de ${\tilde{\cal T}}$, appartient \`a $I_{cusp}(\tilde{G}({\mathbb R}),\omega)$. Par construction, $\varphi'_{\tilde{T}}$ a m\^eme d\'eveloppement infinit\'esimal que $\varphi_{\tilde{T}}$ en $\eta$. On fait maintenant varier $\eta$ parmi un ensemble (fini) de repr\'esentants des classes de conjugaison dans ${\cal O}\cap \tilde{T}({\mathbb R})$ et on fait varier $\tilde{T}$ parmi l'ensemble des \'el\'ements elliptiques de ${\tilde{\cal T}}$. Un argument de partition de l'unit\'e nous fournit un \'el\'ement $\varphi'_{{\tilde{\cal T}}}\in I(\tilde{G}({\mathbb R}),\omega)$ qui a m\^eme image que $\varphi_{{\tilde{\cal T}}}$ dans $I(\tilde{G}({\mathbb R}),\omega)_{{\cal O},loc}$. Cela ach\`eve la preuve pour $F={\mathbb R}$. Si $F={\mathbb C}$, il n'y a qu'un seul \'el\'ement dans ${\tilde{\cal T}}$, qui est un espace de Levi minimal. L'espace $I_{cusp}(\tilde{G}({\mathbb C}),\omega)$ est nul sauf si $G$ est un tore, auquel cas $I_{cusp}(\tilde{G}({\mathbb C}),\omega)=I(\tilde{G}({\mathbb C}),\omega)$. Il en est de m\^eme infinit\'esimalement, quelle que soit la d\'efinition. $\square$
  
   Rappelons que, pour nous, un \'el\'ement est elliptique s'il appartient \`a un sous-tore tordu maximal elliptique. On a
  
  (5) supposons ${\cal O}$ form\'e d'\'el\'ements non-elliptiques; alors $D_{g\acute{e}om}({\cal O},\omega) $ est form\'e de combinaisons lin\'eaires de distributions induites \`a partir d'espaces de Levi propres.

  Preuve. L'espace $D_{g\acute{e}om}({\cal O},\omega) $ est engendr\'e par des distributions $\boldsymbol{\gamma}_{\eta,\tilde{T},\Omega,D}$ comme plus haut, o\`u $\tilde{T}\in \tilde{{\cal T}}$ et $\eta\in {\cal O}\cap \tilde{T}(F)$. Notre d\'efinition d'ellipticit\'e implique que $\tilde{T}$ n'est pas elliptique. Il est donc contenu dans un espace de Levi propre $\tilde{M}$. Les m\^emes donn\'ees $\eta$, $\tilde{T}$, $\Omega$, $D$ d\'efinissent une distribution $\boldsymbol{\gamma}_{\tilde{M},\eta,\tilde{T},\Omega,D}\in D_{g\acute{e}om}(\tilde{M}(F),\omega)$ dont $\boldsymbol{\gamma}_{\eta,\tilde{T},\Omega,D}$ est  l'induite. $\square$
  
  D\'ecrivons plus concr\`etement l'espace $D_{g\acute{e}om}({\cal O},\omega) $ dans le cas o\`u $F={\mathbb R}$ et  ${\cal O}$ est une unique classe de conjugaison. Fixons $\eta\in {\cal O}$. 
    Fixons  un ensemble fini  $\tilde{{\cal T}}$ de sous-tores tordus maximaux de $\tilde{G}$ tels que:
 
 $\bullet$ $\eta\in \tilde{T}({\mathbb R})$ pour tout $\tilde{T}\in \tilde{{\cal T}}$;

 $\bullet$ pour tout sous-tore maximal $S$ de $G_{\eta}$, il existe  $\tilde{T}\in \tilde{{\cal T}}$ et il existe $g\in Z_{G}(\eta;{\mathbb R})$ tels que $S=ad_{g}(T^{\theta,0})$.

 Pour tout $\tilde{T}\in \tilde{{\cal T}}$, notons $\underline{\Omega}_{\tilde{T}}$ l'ensemble des composantes connexes de $\mathfrak{t}^{\theta}({\mathbb R})\cap \mathfrak{g}_{\eta,reg}({\mathbb R})$. Pour $f\in C_{c}^{\infty}(\tilde{G}({\mathbb R}))$, pour $\tilde{T}\in \tilde{{\cal T}}$ et $\Omega\in \underline{\Omega}_{\tilde{T}}$, consid\'erons la fonction $\phi_{f,\tilde{T},\Omega}$ sur $\Omega$ d\'efinie par
 $$\phi_{f,\tilde{T},\Omega}(X)= I^{\tilde{G}}(exp(X)\eta,\omega,f).$$
 Elle est nulle si $\omega$ n'est pas trivial sur $T^{\theta}({\mathbb R})$. Comme on l'a dit, Harish-Chandra a prouv\'e que cette fonction se prolongeait en une fonction $C^{\infty}$ dans un voisinage de $\Omega$. Fixons des coordonn\'ees sur $\mathfrak{t}^{\theta}({\mathbb R})$ et notons ${\mathbb C}[[\mathfrak{t}^{\theta}({\mathbb R})]]$ l'espace des s\'eries formelles sur $\mathfrak{t}^{\theta}({\mathbb R})$. On note $\varphi_{f,\tilde{T},\Omega}\in {\mathbb C}[[\mathfrak{t}^{\theta}({\mathbb R})]]$ le d\'eveloppement en s\'erie de la fonction $\phi_{f,\tilde{T},\Omega}$ en $X=0$. On pose $\varphi_{f}=(\varphi_{f,\tilde{T},\Omega})_{\tilde{T}\in \tilde{{\cal T}},\Omega\in \underline{\Omega}_{\tilde{T}}}$. L'espace $I(\tilde{G}({\mathbb R}),\omega)_{{\cal O},loc}$ est celui  de ces familles $\varphi_{f}$ quand $f$ d\'ecrit $C_{c}^{\infty}(\tilde{G}({\mathbb R}))$. C'est un sous-espace de 
 $$(6) \qquad \oplus_{\tilde{T}\in \tilde{{\cal T}}, \Omega\in \underline{\Omega}_{\tilde{T}}}{\mathbb C}[[\mathfrak{t}^{\theta}({\mathbb R})]].$$
  On sait le d\'ecrire. C'est le sous-espace des familles de s\'eries formelles $(\varphi_{\tilde{T},\Omega})_{\tilde{T}\in \tilde{{\cal T}},\Omega\in \underline{\Omega}_{\tilde{T}}}$ qui v\'erifient deux conditions:
 
 (7) soient $\tilde{T}, \tilde{T}'\in \tilde{{\cal T}}$  et $g\in G({\mathbb R})$ tel que $g\eta g^{-1}=\eta$ et $g\tilde{T}g^{-1}=\tilde{T}'$; alors $ad_{g}$ envoie $\underline{\Omega}_{\tilde{T}}$ sur  $\underline{\Omega}_{\tilde{T}'}$    et ${\mathbb C}[\mathfrak{t}^{\theta}({\mathbb R})]$ sur  ${\mathbb C}[\mathfrak{t}^{_{'}\theta}({\mathbb R})]$; pour $\Omega\in \underline{\Omega}_{\tilde{T}}$, on doit avoir $\varphi_{\tilde{T}',ad_{g}(\Omega)}=\omega(g)ad_{g}(\varphi_{\tilde{T},\Omega})$;
 
 (8) soient $\tilde{T}\in \tilde{{\cal T}}$ et $\Omega$, $\Omega'$ deux \'el\'ements adjacents de $\underline{\Omega}_{\tilde{T}}$; alors une condition de saut relie $\varphi_{\tilde{T},\Omega}$, $\varphi_{\tilde{T},\Omega'}$ et $\varphi_{\tilde{T}_{1},\Omega_{1}}$, o\`u $\tilde{T}_{1}$ et $\Omega_{1}$ sont d\'etermin\'es par $\tilde{T}$, $\Omega$, $\Omega'$, $\tilde{T}_{1}$ \'etant plus d\'eploy\'e que $\tilde{T}$ (c'est-\`a-dire que l'on  a $dim(A_{\tilde{T}})<dim( A_{\tilde{T}_{1}})$).
 
 On renvoie \`a [R2] 3.2 pour cette condition de saut. La topologie sur $I(\tilde{G}({\mathbb R}),\omega)_{{\cal O},loc}$ s'identifie \`a celle d\'eduite de la topologie habituelle sur les espaces de s\'eries formelles (un voisinage de $0$ contient les s\'eries qui s'annulent en $0$ \`a un ordre assez grand).
 Pour $\tilde{T}\in \tilde{{\cal T}}$, notons $D[\mathfrak{t}^{\theta}({\mathbb R})]$ l'espace des op\'erateurs diff\'erentiels \`a coefficients constants sur $\mathfrak{t}^{\theta}({\mathbb R})$. Cet espace se plonge naturellement dans le dual de ${\mathbb C}[[\mathfrak{t}^{\theta}({\mathbb R})]]$: on applique un op\'erateur diff\'erentiel \`a une s\'erie formelle et on \'evalue le r\'esultat en $0$. Ainsi 
 $$\oplus_{\tilde{T}\in \tilde{{\cal T}}, \Omega\in \underline{\Omega}_{\tilde{T}}}D[\mathfrak{t}^{\theta}({\mathbb R})]$$
  se plonge dans le dual de l'espace (6). Par restriction, on obtient une application lin\'eaire
  $$\oplus_{\tilde{T}\in \tilde{{\cal T}}, \Omega\in \underline{\Omega}_{\tilde{T}}}D[\mathfrak{t}^{\theta}({\mathbb R})]\to ( I(\tilde{G}({\mathbb R}),\omega)_{{\cal O},loc})^*.$$
  L'espace $D_{g\acute{e}om}({\cal O},\omega)$ est l'image de cette application.
  
  {\bf Remarque.} Arthur donne une description beaucoup plus pr\'ecise en [A3] lemme 1.1.

  \bigskip
  
  \subsection{Filtration de $D_{g\acute{e}om}(\tilde{G}(F),\omega)$}
Fixons des mesures de Haar sur $G(F)$ et sur $M(F)$ pour tout Levi $M$ de $G$. Pour tout entier $n\geq-1$ notons ${\cal F}^nD_{g\acute{e}om}(\tilde{G}(F),\omega)$ le sous-espace de $D_{g\acute{e}om}(\tilde{G}(F),\omega)$ engendr\'e par les distributions induites $(\boldsymbol{\gamma}_{\tilde{M}})^{\tilde{G}}$, o\`u $\tilde{M}$ est un espace de Levi de $\tilde{G}$ tel que $a_{\tilde{M}}=n+1$ et $\boldsymbol{\gamma}_{\tilde{M}}\in D_{g\acute{e}om}(\tilde{M}(F),\omega)$. Ces espaces forment une filtration
$$\{0\}={\cal F}^{a_{\tilde{M}_{0}}}D_{g\acute{e}om}(\tilde{G}(F),\omega)\subset {\cal F}^{a_{\tilde{M}_{0}}-1}D_{g\acute{e}om}(\tilde{G}(F),\omega)\subset...$$
$$...\subset {\cal F}^{a_{\tilde{G}}-1}D_{g\acute{e}om}(\tilde{G}(F),\omega)=D_{g\acute{e}om}(\tilde{G}(F),\omega).$$
Pour une r\'eunion finie ${\cal O}$ de classes de conjugaison semi-simples dans  $ \tilde{G}(F)$, notons ${\cal F}^nD_{g\acute{e}om}({\cal O},\omega)$ le sous-espace de $D_{g\acute{e}om}({\cal O},\omega)$ engendr\'e par distributions induites $(\boldsymbol{\gamma}_{\tilde{M}})^{\tilde{G}}$, o\`u $\tilde{M}$ est un espace de Levi de $\tilde{G}$ tel que $a_{\tilde{M}}=n+1$ et $\boldsymbol{\gamma}_{\tilde{M}}\in D_{g\acute{e}om}({\cal O}\cap \tilde{M}(F),\omega)$.

Rappelons que l'on a d\'efini en 4.2 une filtration $({\cal F}^nI(\tilde{G}(F),\omega))_{n\geq-1}$ de $I(\tilde{G}(F),\omega)$.

\ass{Proposition}{(i) Pour tout entier $n\geq-1$, ${\cal F}^nI(\tilde{G}(F),\omega)$ est l'annulateur de ${\cal F}^nD_{g\acute{e}om}(\tilde{G}(F),\omega)$ dans $I(\tilde{G}(F),\omega)$ et ${\cal F}^nD_{g\acute{e}om}(\tilde{G}(F),\omega)$ est l'annulateur de ${\cal F}^nI(\tilde{G}(F),\omega)$ dans $D_{g\acute{e}om}(\tilde{G}(F),\omega)$.

(ii) Pour toute r\'eunion finie ${\cal O}$ de classes de conjugaison semi-simples dans $ \tilde{G}(F)$ et tout entier $n\geq -1$, on a l'\'egalit\'e
$${\cal F}^nD_{g\acute{e}om}({\cal O},\omega)={\cal F}^nD_{g\acute{e}om}(\tilde{G}(F),\omega)\cap D_{g\acute{e}om}({\cal O},\omega).$$}

Preuve. On aura besoin d'une propri\'et\'e pr\'eliminaire.  Pour tout $n\geq 0$, fixons  un ensemble $\underline{{\cal L}}^n$ de repr\'esentants des classes de conjugaison par $G(F)$ d'espaces de Levi $\tilde{M} $ de $\tilde{G}$ tels que $a_{\tilde{M}}=n$. On consid\`ere l'application
$$(1)\qquad \begin{array}{cccc}p^n:&I(\tilde{G}(F),\omega)&\to&I^n=\oplus_{\tilde{M}\in \underline{{\cal L}}^n}I(\tilde{M}(F),\omega)^{W(\tilde{M})}\\ &f&\mapsto&\oplus_{\tilde{M}\in \underline{{\cal L}}^n}f_{\tilde{M},\omega}.\\ \end{array}$$
Posons 
$$I^n_{cusp}=\oplus_{\tilde{M}\in \underline{{\cal L}}^n}I_{cusp}(\tilde{M}(F),\omega)^{W(\tilde{M})}.$$ 
Par d\'efinition, ${\cal F}^nI(\tilde{G}(F),\omega)$ est l'image r\'eciproque par $p^n$ du sous-espace $I^n_{cusp}$ de $I^n$. On a vu en 4.2 que de l'application $p^n$  se d\'eduisait un isomorphisme
$$(2) \qquad {\cal F}^{n}I(\tilde{G}(F),\omega)/{\cal F}^{n-1}I(\tilde{G}(F),\omega)\simeq I^n_{cusp}.$$
Soit ${\cal O}$ une r\'eunion finie de classes de conjugaison semi-simples dans $\tilde{G}(F)$. On a d\'efini l'espace $I(\tilde{G},\omega)_{{\cal O},loc}$ en 5.1 et 5.2. Pour tout espace de Levi $\tilde{M}$, posons ${\cal O}_{\tilde{M}}={\cal O}\cap \tilde{M}(F)$.     Montrons que

(3) soit $f\in I(\tilde{G}(F),\omega)$; supposons que, pour tout $\tilde{M}\in \underline{{\cal L}}^{n}$, l'image de $f_{\tilde{M},\omega}$ dans $I(\tilde{M}(F),\omega)_{{\cal O}_{\tilde{M}},loc}$ soit nulle; alors il existe $f'\in I(\tilde{G}(F),\omega)$  telle que $p^n(f')=p^n(f)$ et dont l'image dans $I(\tilde{G}(F),\omega)_{{\cal O},loc}$ soit nulle.

 Supposons d'abord $F$ non-archim\'edien. L'hypoth\`ese signifie que $f_{\tilde{M},\omega}\in I(\tilde{M}(F),\omega)_{{\cal O}_{\tilde{M}},0}$ pour tout $\tilde{M}\in \underline{{\cal L}}^{n}$, autrement dit il existe un voisinage $\tilde{V}_{\tilde{M}}$ de ${\cal O}_{\tilde{M}}$ dans $\tilde{M}(F)$, invariant par conjugaison, tel que $f_{\tilde{M},\omega}$ soit nul sur $\tilde{V}_{\tilde{M}}$.  Fixons de tels voisinages. On peut fixer un voisinage $\tilde{V}$ de ${\cal O}$ dans $\tilde{G}(F)$, invariant par conjugaison, tel que $\tilde{V}\cap \tilde{M}(F)\subset \tilde{V}_{\tilde{M}}$ pour tout $\tilde{M}\in \underline{{\cal L}}^{n}$.    On peut supposer $\tilde{V}$ ouvert et ferm\'e.  Alors la fonction $f'=f(1-{\bf 1}_{\tilde{V}})$ r\'epond \`a la question. 

Supposons maintenant $F$ archim\'edien. Si $n=a_{\tilde{G}}$, l'application $p^n$ est l'identit\'e
de $I(\tilde{G}(F),\omega)$ et l'assertion est claire ($f'=f$ r\'epond \`a la question). Supposons $n>a_{\tilde{G}}$ et raisonnons par r\'ecurrence sur $n$. Soit $\tilde{M}\in \underline{{\cal L}}^{n-1}$, consid\'erons l'\'el\'ement $f_{\tilde{M},\omega}\in I(\tilde{M}(F),\omega)^{W(\tilde{M})}$. L'hypoth\`ese implique que son image dans $I(\tilde{M}(F),\omega)_{{\cal O}_{\tilde{M}},loc}$ est cuspidale au sens de la deuxi\`eme d\'efinition de 5.2 (3). Pr\'ecis\'ement, cette relation nous dit qu'il existe $\varphi^{\tilde{M}}\in I_{cusp}(\tilde{M}(F),\omega)$ qui a m\^eme image que $f_{\tilde{M},\omega}$ dans $I(\tilde{M}(F),\omega)_{{\cal O}_{\tilde{M}},loc}$. En moyennant $\varphi^{\tilde{M}}$ sur $W(\tilde{M})$, on peut supposer $\varphi^{\tilde{M}}\in I_{cusp}(\tilde{M}(F),\omega)^{W(\tilde{M})}$. Posons $\varphi=( \varphi^{\tilde{M}})_{\tilde{M}\in \underline{{\cal L}}^{n-1}}$. En appliquant (2) pour $n-1$, on rel\`eve $\varphi$ en un \'el\'ement $f_{0}\in {\cal F}^{n-1}I(\tilde{G}(F),\omega)$. Pour tout $\tilde{M}\in \underline{{\cal L}}^{n-1}$, la fonction $(f-f_{0})_{\tilde{M},\omega}$ est par construction d'image nulle dans $I(\tilde{M}(F),\omega)_{{\cal O}_{\tilde{M}},loc}$. L'hypoth\`ese de r\'ecurrence assure l'existence de $f'\in I(\tilde{G}(F),\omega)$ d'image nulle dans $I(\tilde{G}(F),\omega)$ et telle que $p^{n-1}(f')=p^{n-1}(f-f_{0})$. L'application $p^n$ se factorise par $p^{n-1}$. On a donc aussi $p^n(f')=p^n(f-f_{0})$. Mais $p^n(f_{0})=0$ puisque $f_{0}\in {\cal F}^{n-1}I(\tilde{G}(F),\omega)$. Donc $f'$ r\'epond \`a la question. Cela prouve (3).

Puisque ${\cal F}^nI(\tilde{G}(F),\omega)$ est l'image r\'eciproque par $p^n$ de  $I^n_{cusp}$, (3) entra\^{\i}ne

(4)  soit $f\in {\cal F}^nI(\tilde{G}(F),\omega)$; supposons que, pour tout $\tilde{M}\in \underline{{\cal L}}^{n}$, l'image de $f_{\tilde{M},\omega}$ dans $I(\tilde{M}(F),\omega)_{{\cal O}_{\tilde{M}},loc}$ soit nulle; alors il existe $f'\in{\cal F}^n I(\tilde{G}(F),\omega)$  telle que $p^n(f')=p^n(f)$ et dont l'image dans $I(\tilde{G}(F),\omega)_{{\cal O},loc}$ soit nulle.

Venons-en \`a la preuve de la proposition. Il est clair que l'annulateur de $D_{g\acute{e}om}(\tilde{G}(F),\omega)$ dans $I(\tilde{G}(F),\omega)$ est nul et que l'annulateur de $I(\tilde{G}(F),\omega)$ dans  $D_{g\acute{e}om}(\tilde{G}(F),\omega)$ est nul.   Soient $f\in I(\tilde{G}(F),\omega)$ et $n\geq -1$. Alors $f$ appartient \`a l'annulateur de  
${\cal F}^nD_{g\acute{e}om}(\tilde{G}(F),\omega)$ si et seulement si, pour tout $\tilde{M}\in \underline{{\cal L}}^{n+1}$ et tout $\boldsymbol{\gamma}_{\tilde{M}}\in D_{g\acute{e}om}(\tilde{M}(F),\omega)$, on a $I^{\tilde{G}}((\boldsymbol{\gamma}_{\tilde{M}})^{\tilde{G}},\omega,f)=0$. Cette \'egalit\'e \'equivaut \`a $I^{\tilde{M}}(\boldsymbol{\gamma}_{\tilde{M}},\omega,f_{\tilde{M},\omega})=0$. Comme on vient de le dire, elle est v\'erifi\'ee pour tout $\boldsymbol{\gamma}_{\tilde{M}}$ si et seulement si $f_{\tilde{M},\omega}=0$. Donc $f$ appartient \`a l'annulateur de  
${\cal F}^nD_{g\acute{e}om}(\tilde{G}(F),\omega)$ si et seulement si $f_{\tilde{M}\omega}=0$ pour tout $\tilde{M}\in \underline{{\cal L}}^{n+1}$. Mais c'est la d\'efinition de l'espace ${\cal F}^nI(\tilde{G}(F),\omega)$. Cela prouve la premi\`ere assertion.

Pour tout entier $n\geq-1$, notons $Ann^n$ l'annulateur de ${\cal F}^nI(\tilde{G}(F),\omega)$ dans $D_{g\acute{e}om}(\tilde{G}(F),\omega)$. Fixons une  r\'eunion finie ${\cal O}$ de classes de conjugaison semi-simples  dans $\tilde{G}(F)$.  On va prouver  que

(5)  ${\cal F}^nD_{g\acute{e}om}({\cal O},\omega)=Ann^n\cap D_{g\acute{e}om}({\cal O},\omega)$. 

D'apr\`es ce que l'on a d\'ej\`a d\'emontr\'e, on a 
$${\cal F}^nD_{g\acute{e}om}(\tilde{G}(F),\omega)\subset Ann^n.$$
D'autre part, par d\'efinition, on a
$${\cal F}^nD_{g\acute{e}om}({\cal O},\omega)\subset {\cal F}^nD_{g\acute{e}om}(\tilde{G}(F),\omega).$$
Donc le membre de gauche de (5) est inclus dans celui de droite. On d\'emontre l'inclusion inverse par r\'ecurrence descendante sur $n$. Si $n=a_{\tilde{M}_{0}}$, on a ${\cal F}^{n}I(\tilde{G}(F),\omega)=I(\tilde{G},\omega)$ et $Ann^n=\{0\}$ comme on l'a dit ci-dessus. L'inclusion est \'evidente. Supposons que $n<a_{\tilde{M}_{0}}$ et que l'assertion soit v\'erifi\'ee pour $n+1$. Soit $\boldsymbol{\gamma}\in Ann^n\cap D_{g\acute{e}om}({\cal O},\omega)$. Supposons d'abord $F$ non-archim\'edien.   On a d\'efini en 5.1 l'espace $I(\tilde{G}(F),\omega)_{{\cal O},0}$. C'est  le noyau de l'application $I(\tilde{G}(F),\omega)\to I(\tilde{G}(F),\omega)_{{\cal O},loc}$. La propri\'et\'e (4) entra\^{\i}ne que l'application
$${\cal F}^{n+1}I(\tilde{G}(F),\omega)\cap I(\tilde{G}(F),\omega)_{{\cal O},0}\to \oplus_{\tilde{M}\in \underline{{\cal L}}^{n+1}} I(\tilde{M}(F),\omega)_{{\cal O}_{\tilde{M}},0}\cap I_{cusp}(\tilde{M}(F),\omega)^{W(\tilde{M})}$$
est surjective. Puisque $\boldsymbol{\gamma}\in Ann^n$, la distribution $\boldsymbol{\gamma}$ se factorise en une forme lin\'eaire $\boldsymbol{\gamma}^{n+1}$ sur ${\cal F}^{n+1}I(\tilde{G}(F),\omega)/{\cal F}^{n}I(\tilde{G}(F),\omega)\simeq I^{n+1}_{cusp}$. Puisque $\boldsymbol{\gamma}\in  D_{g\acute{e}om}({\cal O},\omega)$, la surjectivit\'e ci-dessus entra\^{\i}ne que $\boldsymbol{\gamma}^{n+1}$ annule le sous-espace 
$$\oplus_{\tilde{M}\in \underline{{\cal L}}^{n+1}} I(\tilde{M}(F),\omega)_{{\cal O}_{\tilde{M}},0}\cap I_{cusp}(\tilde{M}(F),\omega)^{W(\tilde{M})}\subset I^{n+1}_{cusp}.$$
  On peut donc prolonger $\boldsymbol{\gamma}^{n+1}$ en une forme lin\'eaire sur 
$$\oplus_{\tilde{M}\in \underline{{\cal L}}^{n+1}}\left( I(\tilde{M}(F),\omega)_{{\cal O}_{\tilde{M}},0}+  I_{cusp}(\tilde{M}(F),\omega)^{W(\tilde{M})}\right),$$
nulle sur 
$$\oplus_{\tilde{M}\in \underline{{\cal L}}^{n+1}} I(\tilde{M}(F),\omega)_{{\cal O}_{\tilde{M}},0}.$$
On peut ensuite prolonger cette forme lin\'eaire en une forme lin\'eaire $\oplus_{\tilde{M}\in \underline{{\cal L}}^{n+1}}\boldsymbol{\gamma}_{\tilde{M}}$ sur
$$I^{n+1}=\oplus_{\tilde{M}\in \underline{{\cal L}}^{n+1}} I(\tilde{M}(F),\omega).$$
Pour tout $\tilde{M}$, $\boldsymbol{\gamma}_{\tilde{M}}$ annule $I(\tilde{M}(F),\omega)_{{\cal O}_{\tilde{M}},0}$ donc $\boldsymbol{\gamma}_{\tilde{M}}\in D_{g\acute{e}om}({\cal O}_{\tilde{M}},\omega)$. La distribution
$$\boldsymbol{\gamma}'=\oplus_{\tilde{M}\in \underline{{\cal L}}^{n+1}}(\boldsymbol{\gamma}_{\tilde{M}})^{\tilde{G}}$$
appartient \`a ${\cal F}^nD_{g\acute{e}om}({\cal O},\omega)$ donc annule ${\cal F}^nI(\tilde{G}(F),\omega)$. Puisque $\oplus_{\tilde{M}\in \underline{{\cal L}}^{n+1}}\boldsymbol{\gamma}_{\tilde{M}}$ co\"{\i}ncide par construction avec $\boldsymbol{\gamma}^{n+1}$ sur $I^{n+1}_{cusp}$, $\boldsymbol{\gamma}'$ co\"{\i}ncide avec $\boldsymbol{\gamma}$ sur ${\cal F}^{n+1}I(\tilde{G}(F),\omega)$. Alors $\boldsymbol{\gamma}-\boldsymbol{\gamma}'$ appartient \`a $Ann^{n+1}\cap D_{g\acute{e}om}({\cal O},\omega)$. En appliquant l'hypoth\`ese de r\'ecurrence, on a 
$$\boldsymbol{\gamma}-\boldsymbol{\gamma}'\in {\cal F}^{n+1}D_{g\acute{e}om}({\cal O},\omega)\subset {\cal F}^nD_{g\acute{e}om}({\cal O},\omega).$$
Donc aussi $\boldsymbol{\gamma}\in {\cal F}^nD_{g\acute{e}om}({\cal O},\omega)$. Cela prouve (5) quand $F$ est non-archim\'edien. Supposons maintenant $F$ archim\'edien. Le noyau de l'application $I(\tilde{G}(F),\omega)\to I(\tilde{G}(F),\omega)_{{\cal O},loc}$ est maintenant l'espace $C\ell I(\tilde{G}(F),\omega)_{{\cal O},0}$ d\'efini en 5.2. On peut reprendre le raisonnement en utilisant cet espace \`a la place de $I(\tilde{G}(F),\omega)_{{\cal O},0}$. Il faut v\'erifier que les formes lin\'eaires que l'on construit sont continues. La continuit\'e de $\boldsymbol{\gamma}^{n+1}$ r\'esulte du fait que les espaces ${\cal F}^nI(\tilde{G}(F),\omega)$ sont \'evidemment ferm\'es et que l'isomorphisme (2) est un hom\'eomorphisme  ([R1] th\'eor\`eme 11.2). Il faut pouvoir choisir des $\boldsymbol{\gamma}_{\tilde{M}}$ continus. Pour cela, il suffit de prouver que

(6) pour tout $\tilde{M}\in \underline{{\cal L}}^{n+1}$, le sous-espace 
$$C\ell I(\tilde{M}(F),\omega)_{{\cal O}_{\tilde{M}},0}+ I_{cusp}(\tilde{M}(F),\omega)^{W(\tilde{M})}$$
de $I(\tilde{M}(F),\omega)$ est ferm\'e. 

Le groupe $W(\tilde{M})$ agit sur $I(\tilde{M}(F),\omega)$. On peut d\'ecomposer cet espace en somme de sous-espaces isotypiques pour cette action.  Chacun de ces sous-espaces est ferm\'e et $I(\tilde{M}(F),\omega)$ en est la somme directe topologique. Notons cette d\'ecomposition
$$I(\tilde{M}(F),\omega)=\oplus_{\tau\in W(\tilde{M})^{\vee}}I(\tilde{M}(F),\omega)_{\tau}.$$  
Par d\'efinition de ${\cal O}_{\tilde{M}}$, l'espace $I(\tilde{M}(F),\omega)_{{\cal O}_{\tilde{M}},0}$ est invariant par $W(\tilde{M})$, donc somme directe de ses intersections avec chacun des sous-espaces $I(\tilde{M}(F),\omega)_{\tau}$. Il en r\'esulte que $C\ell I(\tilde{M}(F),\omega)_{{\cal O}_{\tilde{M}},0}$ v\'erifie la m\^eme propri\'et\'e. Notons
$$C\ell I(\tilde{M}(F),\omega)_{{\cal O}_{\tilde{M}},0}=\oplus_{\tau\in W(\tilde{M})^{\vee}}C\ell I(\tilde{M}(F),\omega)_{{\cal O}_{\tilde{M}},0,\tau}$$
la d\'ecomposition obtenue. Remarquons que le sous-espace d'invariants $I(\tilde{M}(F),\omega)^{W(\tilde{M})}$ n'est autre que $I(\tilde{M}(F),\omega)_{{\bf 1}}$, o\`u ${\bf 1}$ est la repr\'esentation triviale de $W(\tilde{M})$. Alors 
$$C\ell I(\tilde{M}(F),\omega)_{{\cal O}_{\tilde{M}},0}+ I_{cusp}(\tilde{M}(F),\omega)^{W(\tilde{M})}$$
est la somme directe de
$$C\ell I(\tilde{M}(F),\omega)_{{\cal O}_{\tilde{M}},0}^{W(\tilde{M})}+ I_{cusp}(\tilde{M}(F),\omega)^{W(\tilde{M})}$$
et des espaces $C\ell I(\tilde{M}(F),\omega)_{{\cal O}_{\tilde{M}},0,\tau}$ pour $\tau\not=1$. Ces derniers \'etant ferm\'es, il suffit de prouver que le premier l'est. Celui-ci 
est l'intersection de $I(\tilde{M}(F),\omega)^{W(\tilde{M})}$ avec
$$C\ell I(\tilde{M}(F),\omega)_{{\cal O}_{\tilde{M}},0}+ I_{cusp}(\tilde{M}(F),\omega).$$
Puisque $I(\tilde{M}(F),\omega)^{W(\tilde{M})}$ est ferm\'e, il suffit de prouver que le sous-espace ci-dessus est ferm\'e. Or la propri\'et\'e 5.2(3) assure que c'est l'image r\'eciproque dans $I(\tilde{M}(F),\omega)$ du sous-espace des \'el\'ements cuspidaux de $I(\tilde{M}(F),\omega)_{{\cal O}_{\tilde{M}},loc}$. Et celui-ci est ferm\'e (d'apr\`es sa seconde d\'efinition, cf 5.2(3)). D'o\`u l'assertion (6).

 Modulo ces propri\'et\'es, la m\^eme preuve que dans le cas non-archim\'edien s'applique. Cela prouve (5) pour tout $F$.

Soit $\boldsymbol{\gamma}\in Ann^n$. Par d\'efinition de $D_{g\acute{e}om}(\tilde{G}(F),\omega)$, il existe une r\'eunion finie ${\cal O}$ de classes de conjugaison semi-simples dans $\tilde{G}(F)$ telle que $\boldsymbol{\gamma}\in D_{g\acute{e}om}({\cal O},\omega)$. En appliquant (5), on obtient
$$\boldsymbol{\gamma}\in {\cal F}^nD_{g\acute{e}om}({\cal O},\omega)\subset {\cal F}^nD_{g\acute{e}om}(\tilde{G}(F),\omega).$$
D'o\`u l'inclusion $Ann^n\subset {\cal F}^nD_{g\acute{e}om}(\tilde{G}(F),\omega)$. On a d\'ej\`a prouv\'e l'inclusion r\'eciproque. D'o\`u l'\'egalit\'e de ces espaces, ce qui est la deuxi\`eme assertion de (i). Gr\^ace \`a cette assertion, le (ii) de l'\'enonc\'e n'est autre que (5). $\square$

     \bigskip
   
   \subsection{Distributions g\'eom\'etriques stables dans le cas non-archim\'edien}
   Supposons $F$ non-archim\'edien et $(G,\tilde{G},{\bf a})$ quasi-d\'eploy\'e et \`a torsion int\'erieure. On note $D^{st}_{g\acute{e}om}(\tilde{G}(F))$ le sous-espace des \'el\'ements de $D_{g\acute{e}om}(\tilde{G}(F))$ qui se factorisent en une forme lin\'eaire sur $SI(\tilde{G}(F))$. Soit ${\cal O}$ une r\'eunion finie de classes de conjugaison stable dans $\tilde{G}(F)$. On note $D^{st}_{g\acute{e}om}({\cal O})=D^{st}_{g\acute{e}om}(\tilde{G}(F))\cap D_{g\acute{e}om}({\cal O})$. Notons $SI(\tilde{G}(F))_{{\cal O},0}$ le sous-espace des \'el\'ements $f\in SI(\tilde{G}(F))$ pour lesquels il existe un voisinage $\tilde{U}$ de ${\cal O}$ tel que $S^{\tilde{G}}(\gamma,f)=0$ pour tout $\gamma\in \tilde{U}$. Posons $SI(\tilde{G}(F))_{{\cal O},loc}=SI(\tilde{G}(F))/SI(\tilde{G}(F))_{{\cal O},0}$. On a encore
   
   (1) $D^{st}_{g\acute{e}om}({\cal O})$ est l'espace des formes lin\'eaires sur $C_{c}^{\infty}(\tilde{G}(F))$ qui se factorisent par la projection $C_{c}^{\infty}(\tilde{G}(F))\to SI(\tilde{G}(F))_{{\cal O},loc}$.
   
   On a aussi
   
   (2) $D^{st}_{g\acute{e}om}(\tilde{G}(F))$ est la somme directe des sous-espaces $D^{st}_{g\acute{e}om}({\cal O})$ quand ${\cal O}$ d\'ecrit les classes de conjugaison stable semi-simples.
   
   Preuve. Soit $\boldsymbol{\delta}\in D^{st}_{g\acute{e}om}(\tilde{G}(F))$. Les parties semi-simples des \'el\'ements de son support restent dans un ensemble fini de classes de conjugaison stable. Notons ${\cal O}_{1},...,{\cal O}_{n}$ ces classes. En utilisant la construction de 4.6, on peut trouver pour chaque $i=1,...,n$ un voisinage ouvert et ferm\'e $\tilde{U}_{i}$ de ${\cal O}_{i}$ de sorte que
   
   - $\tilde{U}_{i}\cap \tilde{U}_{j}=\emptyset$ si $i\not=j$;
   
   - si $\gamma,\gamma'\in \tilde{G}_{reg}(F)$ sont stablement conjugu\'es, alors $\gamma\in \tilde{U}_{i}$ si et seulement si $\gamma'\in\tilde{U}_{i}$.
   
   On note ${\bf 1}_{\tilde{U}_{i}}$ la fonction caract\'eristique de $\tilde{U}_{i}$ et $\boldsymbol{\delta}_{i}$ la distribution $f\mapsto \boldsymbol{\delta}(f{\bf 1}_{\tilde{U}_{i}})$. Elle est encore stable d'apr\`es la seconde condition ci-dessus. Elle appartient clairement \`a $D_{g\acute{e}om}^{st}({\cal O}_{i})$. Enfin, $\boldsymbol{\delta}$ est la somme des $\boldsymbol{\delta}_{i}$ d'apr\`es la premi\`ere condition ci-dessus. $\square$
   
   Soit $\tilde{M}$ un espace de Levi de $\tilde{G}$. L'application d'induction pr\'eserve la stabilit\'e (parce que, si $f\in I(\tilde{G}(F))$ a une image nulle dans $SI(\tilde{G}(F))$, alors l'image de $f_{\tilde{M}}$ dans $SI(\tilde{M}(F))$ est nulle). On a donc un homomorphisme d'induction
    $$\begin{array}{ccc}D^{st}_{g\acute{e}om}(\tilde{M}(F))\otimes Mes(M(F))^*&\to &D^{st}_{g\acute{e}om}(\tilde{G}(F))\otimes Mes(G(F))^*\\ \boldsymbol{\delta}&\mapsto&\boldsymbol{\delta}^{\tilde{G}}\\ \end{array}$$

\bigskip
   
  \subsection{Distributions g\'eom\'etriques stables dans le cas archim\'edien}
  On suppose $F$ archim\'edien et $(G,\tilde{G},{\bf a})$ quasi-d\'eploy\'e et \`a torsion int\'erieure.
   On note $D_{g\acute{e}om}^{st}(\tilde{G}(F))$ le sous-espace des \'el\'ements de $D_{g\acute{e}om}(\tilde{G}(F))$ qui se factorisent en une forme lin\'eaire sur $SI(\tilde{G}(F))$. En adaptant la construction du paragraphe 5.2, on munit $SI(\tilde{G}(F))$ d'une topologie. L'espace $D_{g\acute{e}om}^{st}(\tilde{G}(F))$  s'identifie \`a celui des formes lin\'eaires continues sur cet espace qui sont support\'ees par la r\'eunion d'un nombre fini de classes de conjugaison stable semi-simples. Pour une telle r\'eunion finie ${\cal O}$, on d\'efinit les espaces $D^{st}_{g\acute{e}om}({\cal O})$ et $SI(\tilde{G}(F))_{{\cal O},0}$ comme dans le cas non-archim\'edien. On note $C\ell SI(\tilde{G}(F))_{{\cal O},0}$ sa cl\^oture dans $SI(\tilde{G}(F))$ et le quotient $SI(\tilde{G}(F))_{{\cal O},loc}=SI(\tilde{G}(F))/C\ell SI(\tilde{G}(F))_{{\cal O},0}$. On a comme en 5.2(2)
   
   (1) $D^{st}_{g\acute{e}om}({\cal O})$ s'identifie \`a l'espace des formes lin\'eaires continues sur $SI(\tilde{G}(F))_{{\cal O},loc}$.
   
   La preuve de 5.2(3) s'adapte:
   
   (2) les deux d\'efinitions possibles d'un espace $SI_{cusp}(\tilde{G}(F))_{{\cal O},loc}$ sont \'equivalentes.
   
   Enfin, on a
   
   (3) $D^{st}_{g\acute{e}om}(\tilde{G}(F))$ est la somme directe des sous-espaces $D^{st}_{g\acute{e}om}({\cal O})$, quand ${\cal O}$ d\'ecrit les classes de conjugaison stable semi-simples.
   
   La preuve de 5.4(2) s'adapte, en rempla\c{c}ant les fonctions ${\bf 1}_{\tilde{U}_{i}}$ par des fonctions $C^{\infty}$ convenables.
   
   D\'ecrivons concr\`etement l'espace $D^{st}_{g\acute{e}om}({\cal O})$ dans le cas o\`u $F={\mathbb R}$ et o\`u ${\cal O}$ est une unique classe de conjugaison stable. 
   On doit fixer $\eta\in {\cal O}$ tel que $G_{\eta}$ soit quasi-d\'eploy\'e. On choisit  un ensemble $\tilde{{\cal T}}$ de sous-tores tordus maximaux de $\tilde{G}$ de sorte que
  
   $\bullet$ $\eta\in \tilde{T}({\mathbb R})$ pour tout $\tilde{T}\in \tilde{{\cal T}}$;

 $\bullet$ pour tout sous-tore maximal $S$ de $G_{\eta}$, il existe  $\tilde{T}\in \tilde{{\cal T}}$ et il existe $g\in Z_{G}(\eta)$ tels que $S=ad_{g}(T)$ et l'isomorphisme $ad_{g}:T\to S$ soit d\'efini sur ${\mathbb R}$.

   En rempla\c{c}ant les int\'egrales orbitales par les int\'egrales orbitales stables dans les d\'efinitions de 5.2,  l'espace $SI(\tilde{G}({\mathbb R}))_{{\cal O},loc}$ s'identifie \`a un  sous-espace de l'espace
   $$ \oplus_{\tilde{T}\in \tilde{{\cal T}}, \Omega\in \underline{\Omega}_{\tilde{T}}}{\mathbb C}[[\mathfrak{t}({\mathbb R})]].$$
    Gr\^ace aux r\'esultats de Shelstad, on peut encore le caract\'eriser par des conditions similaires \`a 5.2(7) et (8). On construit de m\^eme une application lin\'eaire
  $$\oplus_{\tilde{T}\in \tilde{{\cal T}},\Omega\in \underline{\Omega}_{\tilde{T}}}D[\mathfrak{t}({\mathbb R})]\to SI(\tilde{G}({\mathbb R}))_{{\cal O},loc}$$
  dont l'image est $D_{g\acute{e}om}^{st}({\cal O})$. 
  
  L'\'ecriture des int\'egrales orbitales stables comme somme d'int\'egrales orbitales fournit une application lin\'eaire surjective
  $$\oplus_{{\cal O}'}I(\tilde{G}({\mathbb R}))_{{\cal O}',loc} \to  SI(\tilde{G}({\mathbb R}))_{{\cal O},loc},$$
  o\`u ${\cal O}'$ d\'ecrit les classes de conjugaison par $G({\mathbb R})$ contenues dans ${\cal O}$. Dualement, on a une application lin\'eaire injective
  $$D_{g\acute{e}om}^{st}({\cal O})\to D_{g\acute{e}om}({\cal O})=\oplus_{{\cal O}'}D_{g\acute{e}om}({\cal O}').$$

   \bigskip
   
   \subsection{Constructions formelles} 
   Le corps $F$ est quelconque et $(G,\tilde{G},{\bf a})$ est quasi-d\'eploy\'e et \`a torsion int\'erieure. On suppose donn\'ee une extension
   $$1\to C_{1}\to G_{1}\to G\to 1$$
   o\`u $C_{1}$ est un tore central induit, une extension compatible
   $$\tilde{G}_{1}\to \tilde{G}$$
   avec $\tilde{G}_{1}$ \`a torsion int\'erieure, et un caract\`ere $\lambda_{1}$ de $C_{1}(F)$.
   
   En adaptant les d\'efinitions des paragraphes pr\'ec\'edents, on d\'efinit les espaces de distributions $D_{g\acute{e}om,\lambda_{1}}(\tilde{G}_{1}(F))$  et $D^{st}_{g\acute{e}om,\lambda_{1}}(\tilde{G}_{1}(F))$. Leurs \'el\'ements sont des formes lin\'eaires respectivement sur $I_{\lambda_{1}}(\tilde{G}_{1}(F))$ et $SI_{\lambda_{1}}(\tilde{G}_{1}(F))$. De m\^eme, pour une r\'eunion finie ${\cal O}$ de classes de conjugaison semi-simples  dans $\tilde{G}(F)$, on d\'efinit des espaces localis\'es que l'on note $I_{\lambda_{1}}(\tilde{G}_{1}(F))_{{\cal O},loc}$ et $D_{g\acute{e}om,\lambda_{1}}(\tilde{G}_{1}(F),{\cal O})$. Si ${\cal O}$ est une r\'eunions finie de  classes de conjugaison stable, on a les variantes $SI_{\lambda_{1}}(\tilde{G}_{1}(F))_{{\cal O},loc}$ et $D^{st}_{g\acute{e}om,\lambda_{1}}(\tilde{G}_{1}(F),{\cal O})$. 
   
   D\'ecrivons concr\`etement $D_{g\acute{e}om,\lambda_{1}}(\tilde{G}_{1}(F),{\cal O})$ quand $F={\mathbb R}$ et ${\cal O}$ est une unique classe de conjugaison. On fixe cette fois $\eta_{1}\in \tilde{G}_{1}({\mathbb R})$ se projetant en un \'el\'ement de ${\cal O}$.  Rempla\c{c}ons $\tilde{G}$ par $\tilde{G}_{1}$ et $\eta$ par $\eta_{1}$ dans les constructions de 5.2 pour d\'efinir un ensemble $\tilde{{\cal T}}_{1}$ et, pour tout $\tilde{T}_{1}\in \tilde{{\cal T}}_{1}$, un ensemble $\underline{\Omega}_{\tilde{T}_{1}}$. Pour $f\in C_{c,\lambda_{1}}^{\infty}(\tilde{G}_{1}({\mathbb R}))$, on d\'efinit la famille 
 $$(\varphi_{f,\tilde{T}_{1},\Omega})_{\tilde{T}_{1}\in \tilde{{\cal T}}_{1},\Omega\in \underline{\Omega}_{\tilde{T}_{1}}}\in  \oplus_{\tilde{T}_{1}\in \tilde{{\cal T}}_{1}, \Omega\in \underline{\Omega}_{\tilde{T}_{1}}}{\mathbb C}[[\mathfrak{t}_{1}({\mathbb R})]]$$
 comme en 5.2. Alors  $I_{\lambda_{1}}(\tilde{G}_{1}({\mathbb R}))_{{\cal O},loc}$ est l'espace de ces familles quand $f$ d\'ecrit $C_{c,\lambda_{1}}^{\infty}(\tilde{G}_{1}({\mathbb R}))$. On peut d\'ecrire cet espace comme celui des familles de s\'eries formelles $(\varphi_{\tilde{T}_{1},\Omega})_{\tilde{T}_{1}\in \tilde{{\cal T}}_{1},\Omega\in \underline{\Omega}_{\tilde{T}_{1}}}$ qui v\'erifient la condition 5.2(8) et les conditions (1) et (2) suivantes.
 
 (1) Soient $\tilde{T}_{1},\tilde{T}'_{1}\in \tilde{{\cal T}}_{1}$, $g_{1}\in G_{1}({\mathbb R})$ et $c\in C_{1}({\mathbb R})$ tels que $g_{1}\eta_{1}g_{1}^{-1}=c\eta_{1}$ et $ad_{g_{1}}(\tilde{T}_{1})=\tilde{T}'_{1}$; pour tout $\Omega\in \underline{\Omega}_{\tilde{T}_{1}}$, on a $\varphi_{\tilde{T}'_{1},ad_{g_{1}}(\Omega)}=\lambda_{1}(c)^{-1}ad_{g_{1}}(\varphi_{\tilde{T}_{1},\Omega})$. 
 
 Remarquons que cette condition est plus forte que  5.2(7). Pour tout $\tilde{T}_{1}\in \tilde{{\cal T}}_{1}$, on peut fixer une d\'ecomposition $\mathfrak{t}_{1}=\mathfrak{t}\oplus \mathfrak{c}_{1}$ de sorte que $\mathfrak{t}$ contienne l'intersection de $\mathfrak{t}_{1}$ avec l'alg\`ebre de Lie du groupe d\'eriv\'e de $G_{1}$. On a alors une application injective
 $${\mathbb C}[[\mathfrak{t}({\mathbb R})]]\otimes {\mathbb C}[[\mathfrak{c}_{1}({\mathbb R})]]\to {\mathbb C}[[\mathfrak{t}_{1}({\mathbb R})]].$$
 En d\'eveloppant en s\'erie formelle le caract\`ere $\lambda_{1}^{-1}$, on obtient un \'el\'ement $\varphi_{\lambda_{1}}$ de  ${\mathbb C}[[\mathfrak{c}_{1}({\mathbb R})]]$. Alors
 
 (2) pour tout $\Omega\in \underline{\Omega}_{1}$, il existe $\varphi_{\tilde{T},\Omega}\in {\mathbb C}[[\mathfrak{t}({\mathbb R})]]$ tel que  $\varphi_{\tilde{T}_{1},\Omega}=\varphi_{\tilde{T},\Omega}\varphi_{\lambda_{1}}$.
 
 De nouveau, on a une application lin\'eaire
 $$\oplus_{\tilde{T}_{1}\in \tilde{{\cal T}}_{1}, \Omega\in \underline{\Omega}_{\tilde{T}_{1}}}D[\mathfrak{t}_{1}({\mathbb R})]\to ( I_{\lambda_{1}}(\tilde{G}_{1}({\mathbb R}))_{{\cal O},loc})^*.$$
  L'espace $D_{g\acute{e}om,\lambda_{1}}(\tilde{G}_{1}({\mathbb R}),{\cal O})$ est l'image de cette application. 
  
  Notons ${\cal O}_{1}$ la classe de conjugaison par $G_{1}({\mathbb R})$ engendr\'ee par $\gamma_{1}$. Il r\'esulte des descriptions ci-dessus que
  $$I_{\lambda_{1}}(\tilde{G}_{1}({\mathbb R}))_{{\cal O},loc}\subset I(\tilde{G}_{1}({\mathbb R}))_{{\cal O}_{1},loc}$$
  et qu'il y a une application lin\'eaire naturelle et surjective
 $$D_{g\acute{e}om}({\cal O}_{1})\to D_{g\acute{e}om,\lambda_{1}}(\tilde{G}_{1}({\mathbb R}),{\cal O}).$$

Supposons toujours $F={\mathbb R}$ et soit ${\cal O}$ une classe de conjugaison stable semi-simple. On suppose $G_{1,\eta_{1}}$ quasi-d\'eploy\'e. Les espaces 
$SI_{\lambda_{1}}(\tilde{G}_{1}({\mathbb R}))_{{\cal O},loc}$ et $D_{g\acute{e}om,\lambda_{1}}^{st}(\tilde{G}_{1}({\mathbb R}),{\cal O})$ se d\'ecrivent comme pr\'ec\'edemment, avec de l\'eg\`eres variantes. On a les m\^emes cons\'equences que ci-dessus, \`a savoir que l'on a l'inclusion
 $$SI_{\lambda_{1}}(\tilde{G}_{1}({\mathbb R}))_{{\cal O},loc}\subset SI(\tilde{G}_{1}({\mathbb R}))_{{\cal O}_{1},loc}$$
  et qu'il y a une application lin\'eaire naturelle et surjective
 $$D_{g\acute{e}om}^{st}({\cal O}_{1})\to D_{g\acute{e}om,\lambda_{1}}^{st}(\tilde{G}_{1}({\mathbb R}),{\cal O}).$$

 On revient \`a un corps de base $F$ quelconque.  Consid\'erons une autre s\'erie de donn\'ees $G_{2},\tilde{G}_{2},C_{2},\lambda_{2}$ v\'erifiant les m\^emes hypoth\`eses. Notons  $G_{12}$ le produit fibr\'e de $G_{1}$ et $G_{2}$ au-dessus de $G$ et $\tilde{G}_{12}$ celui de $\tilde{G}_{1}$ et $\tilde{G}_{2}$ au-dessus de $\tilde{G}$. Consid\'erons un caract\`ere continu $\lambda_{12}$ de $G_{12}(F)$ dont la restriction \`a $C_{1}(F)\times C_{2}(F)$ soit $\lambda_{1}\times \lambda_{2}^{-1}$ et une fonction non nulle $\tilde{\lambda}_{12}$ sur $\tilde{G}_{12}(F)$ telle que $\tilde{\lambda}_{12}(g\gamma)=\lambda_{12}(g)\tilde{\lambda}_{12}(\gamma)$ pour tous $g\in G_{12}(F)$ et $\gamma\in \tilde{G}_{12}(F)$. On a alors un isomorphisme
   $$C_{c,\lambda_{1}}^{\infty}(\tilde{G}_{1}(F))\simeq C_{c,\lambda_{2}}^{\infty}(\tilde{G}_{2}(F))$$
   qui, \`a $f_{1}$ sur $\tilde{G}_{1}(F)$, associe la fonction $f_{2}$ sur $\tilde{G}_{2}(F)$ telle que $f_{2}(\gamma_{2})=f_{1}(\gamma_{1})\tilde{\lambda}_{12}(\gamma_{1},\gamma_{2})$ pour tous $(\gamma_{1},\gamma_{2})\in \tilde{G}_{12}(F)$. Remarquons que, dans le cas archim\'edien, il s'agit d'un hom\'eomorphisme, $\tilde{\lambda}_{12}$ \'etant n\'ecessairement $C^{\infty}$. On voit que l'isomorphisme ci-dessus se dualise en un isomorphisme
   $$D_{g\acute{e}om,\lambda_{2}}(\tilde{G}_{2}(F))\simeq D_{g\acute{e}om,\lambda_{1}}(\tilde{G}_{1}(F))$$
   qui se restreint en un isomorphisme
   $$D^{st}_{g\acute{e}om,\lambda_{2}}(\tilde{G}_{2}(F))\simeq D^{st}_{g\acute{e}om,\lambda_{1}}(\tilde{G}_{1}(F)).$$
   
   Revenons au cas o\`u $(G,\tilde{G},{\bf a})$ est quelconque. Soit ${\bf G}'$ une donn\'ee endoscopique relevante pour $(G,\tilde{G},{\bf a})$. Des constructions ci-dessus se d\'eduisent la d\'efinition de l'espace $D_{g\acute{e}om}({\bf G}')$ et de son sous-espace $D^{st}_{g\acute{e}om}({\bf G}')$. Leurs \'el\'ements sont des formes lin\'eaires sur $I({\bf G}')$, resp. $SI({\bf G}')$, continues dans le cas o\`u $F$ est archim\'edien. Soit ${\cal O}'$ une r\'eunion finie de classes de conjugaison stable semi-simples de $\tilde{G}'(F)$. On d\'efinit comme en 5.4  le sous-espace $D^{st}_{g\acute{e}om}({\bf G}',{\cal O}')\subset D^{st}_{g\acute{e}om}({\bf G}')$, l'espace $SI({\bf G}')_{{\cal O}',0}$, sa cl\^oture $C\ell SI({\bf G}')_{{\cal O}',0}$ dans le cas archim\'edien, et le quotient
   $$SI({\bf G}')_{{\cal O}',loc}=\left\lbrace\begin{array}{cc}SI({\bf G}')/SI({\bf G}')_{{\cal O}',0},& \text{ si }F\text{ est non-archim\'edien }\\SI({\bf G}')/C\ell SI({\bf G}')_{{\cal O}',0},& \text{ si }F\text{ est archim\'edien. }\\ \end{array}\right.$$
   L'espace $D^{st}_{g\acute{e}om}({\bf G}',{\cal O}')$ est celui des formes lin\'eaires sur $C_{c}^{\infty}({\bf G}')$ qui se factorisent en une forme lin\'eaire (continue dans le cas archim\'edien) sur $SI({\bf G}')_{{\cal O}',loc}$.

  \bigskip

\subsection{Transfert de distributions "g\'eom\'etriques" }
  
  Si $F$ est non-archim\'edien ou $F={\mathbb C}$, soit ${\cal O}$ une classe de conjugaison stable semi-simple dans $\tilde{G}(F)$. On a d\'efini l'espace $D_{g\acute{e}om}({\cal O},\omega)\subset D_{g\acute{e}om}(\tilde{G}(F),\omega)$. Si $F={\mathbb R}$, on doit travailler ici avec un $K$-espace $K\tilde{G}$. Soit ${\cal O}$ une classe de conjugaison stable semi-simple dans $K\tilde{G}({\mathbb R})$. On d\'efinit  l'espace $D_{g\acute{e}om}({\cal O},\omega)\subset D_{g\acute{e}om}(\tilde{G}({\mathbb R}),\omega)$, somme directe des $D_{g\acute{e}om}({\cal O}_{p},\omega)$ pour $p\in \Pi$, o\`u ${\cal O}_{p} ={\cal O}\cap \tilde{G}_{p}({\mathbb R})$ (${\cal O}_{p}$ peut \^etre vide). Pour tout ${\bf G}'\in {\cal E}(\tilde{G},{\bf a})$, il  correspond \`a ${\cal O}$ une r\'eunion finie ${\cal O}_{\tilde{G}'}$ de classes de conjugaison stable semi-simples dans $\tilde{G}'(F)$, qui peut \^etre vide. Consid\'erons l'espace
  
  (1) $ \oplus_{{\bf G}'\in {\cal E}(\tilde{G},{\bf a})}D_{g\acute{e}om}^{st}({\cal O}_{\tilde{G}'})\otimes Mes(G'(F))^*$.
  
  Nous allons en d\'efinir diff\'erents sous-espaces. Soient ${\bf G}'\in {\cal E}(\tilde{G},{\bf a})$ et $M'$ un Levi de $G'$. On note ${\cal O}_{\tilde{M}'}$ la r\'eunion des classes de conjugaison stables semi-simples dans $\tilde{M}'(F)$ qui sont incluses dans ${\cal O}_{\tilde{G}'}$. En fixant des donn\'ees suppl\'ementaires $G'_{1}$,...,$\Delta_{1}$, on dispose de l'application
  $$\begin{array}{ccc}SI_{\lambda_{1}}(\tilde{G}'_{1}(F))\otimes Mes(G'(F))&\to&SI_{\lambda_{1}}(\tilde{M}'_{1}(F))\otimes Mes(M'(F))\\ f&\mapsto&f_{\tilde{M}'_{1}}\\ \end{array}$$
  Par dualit\'e, on en d\'eduit un homomorphisme
  $$(2) \qquad \begin{array}{ccc}D^{st}_{g\acute{e}om,\lambda_{1}}(\tilde{M}'_{1}(F) ,{\cal O}_{\tilde{M}'})\otimes Mes(M'(F))^*&\to&D^{st}_{g\acute{e}om,\lambda_{1}}(\tilde{G}'_{1}(F),{\cal O}_{\tilde{G}'})\otimes Mes(G'(F))^*\\ \boldsymbol{\delta}&\mapsto&\boldsymbol{\delta}^{\tilde{G}'}\\ \end{array}$$
  o\`u les espaces de distributions sont d\'efinis de fa\c{c}on \'evidente. Le second espace s'identifie \`a  l'espace d\'ej\`a d\'efini $D^{st}_{g\acute{e}om}({\bf G}',{\cal O}_{\tilde{G}'})$. L'espace $D^{st}_{g\acute{e}om,\lambda_{1}}( \tilde{M}'_{1}(F),{\cal O}_{\tilde{M}'})$ et l'homomorphisme ci-dessus ne sont a priori d\'efinis que modulo le choix de donn\'ees auxiliaires. On v\'erifie toutefois que l'image de cet homomorphisme dans $D^{st}_{g\acute{e}om}({\bf G}',{\cal O}_{\tilde{G}'})\otimes Mes(G'(F))^*$ ne d\'epend pas de ce choix. On note cette image ${\cal I}_{\tilde{M}'}^{\tilde{G}'}({\cal O}_{\tilde{M}'})$.
 Supposons que $M'$ soit relevant. Soit $(\tilde{M},{\bf M}')$ l'\'el\'ement de ${\cal E}_{+}(\tilde{G},{\bf a})$ qui lui est associ\'e par la construction de 3.4. On identifie $\tilde{M}'$ \`a l'espace endoscopique issu de ${\bf M}'$. Remarquons qu'il y a deux fa\c{c}ons de d\'efinir un ensemble ${\cal O}_{\tilde{M}'}$:  soit, comme on l'a fait,  par une suite ${\cal O}\mapsto {\cal O}_{\tilde{G}'} \mapsto {\cal O}_{\tilde{M}'}$, soit par une suite ${\cal O}\mapsto {\cal O}_{\tilde{M}}\mapsto {\cal O}_{\tilde{M}'}$. Les deux proc\'ed\'es donnent le m\^eme r\'esultat. L'espace $D^{st}_{g\acute{e}om,\lambda_{1}}( \tilde{M}'_{1},{\cal O}_{\tilde{M}'})$ s'identifie \`a l'espace  $D^{st}_{g\acute{e}om}({\bf M}',{\cal O}_{\tilde{M}'})$ relatif \`a ${\bf M}'$. Toutefois, l'homomorphisme ci-dessus d\'epend du choix de l'identification. Le groupe $Aut(\tilde{M},{\bf M}')$ agit sur $SI({\bf M}')$. Il  r\'esulte de la d\'efinition de ${\cal O}_{\tilde{M}'}$ que cette action pr\'eserve $SI({\bf M}')_{{\cal O}_{\tilde{M}'},0}$ et sa cl\^oture dans le cas archim\'edien. Donc l'action se descend en une action sur $SI({\bf M}')_{{\cal O}_{\tilde{M}'},loc}$ et il y a aussi une action duale sur $D^{st}_{g\acute{e}om}({\bf M}',{\cal O}_{\tilde{M}'})$. On d\'ecompose  cet espace en la somme du sous-espace des invariants $D^{st}_{g\acute{e}om}({\bf M}',{\cal O}_{\tilde{M}'})^{Aut(\tilde{M},{\bf M}')}$ et de son unique suppl\'ementaire invariant par l'action du groupe. On note $D^{st}_{g\acute{e}om}({\bf M}',{\cal O}_{\tilde{M}'})^{non-inv}$ ce suppl\'ementaire et ${\cal I}_{\tilde{M}'}^{\tilde{G}'}({\cal O}_{\tilde{M}'})^{non-inv}$ son image par l'homomorphisme (2). On v\'erifie que ce dernier espace ne d\'epend pas des choix. Enfin, on v\'erifie que la restriction de (2) au sous-espace des invariants devient un homomorphisme
$$ \begin{array}{ccc}D^{st}_{g\acute{e}om}({\bf M}',{\cal O}_{\tilde{M}'})^{Aut(\tilde{M},{\bf M}')}\otimes Mes(M'(F))^*&\to&D^{st}_{g\acute{e}om}({\bf M}',{\cal O}_{\tilde{G}'})\otimes Mes(G'(F))^*\\ \boldsymbol{\delta}&\mapsto&\boldsymbol{\delta}^{\tilde{G}'}\\ \end{array}$$
qui est ind\'ependant des choix. 

On consid\`ere les sous-espaces suivants de l'espace (1):

(3) les espaces ${\cal I}_{\tilde{M}'}^{\tilde{G}'}({\cal O}_{\tilde{M}'})$, pour un ${\bf G}'\in {\cal E}(\tilde{G},{\bf a})$ et un Levi $M'$ de $G'$ qui n'est pas relevant;

(4) les espaces ${\cal I}_{\tilde{M}'}^{\tilde{G}'}({\cal O}_{\tilde{M}'})^{non-inv}$, pour un ${\bf G}'\in {\cal E}(\tilde{G},{\bf a})$ et un Levi $M'$ de $G'$ qui est relevant;

(5) les espaces images d'un homomorphisme
$$\begin{array}{ccc}D^{st}_{g\acute{e}om}({\bf M}',{\cal O}_{\tilde{M}'})^{Aut(\tilde{M},{\bf M}')}\otimes Mes(M'(F))^*&\to&D^{st}_{g\acute{e}om}({\bf G}',{\cal O}_{\tilde{G}'})\otimes Mes(G'(F))^*\\&&+ D^{st}_{g\acute{e}om}(\underline{{\bf G}}',{\cal O}_{\underline{\tilde{G}}'})\otimes Mes(\underline{G}'(F))^*\\ \boldsymbol{\delta}&\mapsto &\boldsymbol{\delta}^{G'}-\boldsymbol{\delta}^{\underline{G}'}\\ \end{array}$$
pour deux \'el\'ements ${\bf G}',\underline{\bf G}'\in {\cal E}(\tilde{G},{\bf a})$ et un Levi commun $M'$ qui est relevant.
 
 {\bf Remarque.} Pr\'ecis\'ement, dans cette derni\`ere condition, on consid\`ere ${\bf G}',\underline{\bf G}'\in {\cal E}(\tilde{G},{\bf a})$, des Levi $M'$ de $G'$ et $\underline{M}'$ de $\underline{G}'$ et on suppose que l'\'el\'ement $(\tilde{M},{\bf M}')$ de ${\cal E}_{+}(\tilde{G},{\bf a})$ associ\'e \`a ces Levi est le m\^eme. Mais les donn\'ees ${\bf G}'$ et $  \underline{\bf G}'$ peuvent \^etre les m\^emes,  un m\^eme \'el\'ement $(\tilde{M},{\bf M}')$ pouvant \^etre associ\'e \`a deux Levi distincts du m\^eme groupe $G'$. Par exemple, si $G=SO(11)$ et $G'=SO(5)\times SO(7)$, aux Levi $GL(2)\times (GL(1)\times GL(1)\times GL(1))$ et $(GL(1)\times GL(1))\times(GL(2)\times GL(1))$ de $G'$ est associ\'e le m\^eme \'el\'ement de ${\cal E}_{+}(\tilde{G},{\bf a})$. Si les donn\'ees ${\bf G}'$ et $\underline{{\bf G}}'$ sont les m\^emes, les applications $\boldsymbol{\delta}\mapsto \boldsymbol{\delta}^{G'}$ et $\boldsymbol{\delta}\mapsto \boldsymbol{\delta}^{\underline{G}'}$ sont \`a valeurs dans le m\^eme espace mais ne sont pas forc\'ement les m\^emes comme le montre l'exemple ci-dessus (et malgr\'e la notation impr\'ecise qui pourrait le faire croire).
  
   \ass{Proposition}{ Par dualit\'e, le transfert d\'efinit une application lin\'eaire
  $$transfert:\oplus_{{\bf G}'\in {\cal E}(\tilde{G},{\bf a})}D_{g\acute{e}om}^{st}({\bf G}',{\cal O}_{\tilde{G}'})\otimes Mes(G'(F))^*\to D_{g\acute{e}om}({\cal O},\omega)\otimes Mes(G(F))^*.$$
  Elle est surjective.  Son noyau est la somme des sous-espaces d\'ecrits ci-dessus.}
  
  La preuve est donn\'ee dans les deux paragraphes suivants.
  
  \bigskip
  \subsection{Preuve dans le cas non-archim\'edien}
  
 On suppose $F$ non-archim\'edien. Pour simplifier, on fixe des mesures de Haar sur tous les groupes intervenant, ce qui \'elimine les espaces de mesures.  D\'efinissons un espace $I^{\cal E}_{+}(\tilde{G}(F),\omega)_{{\cal O},loc}$. C'est le sous-espace des \'el\'ements $(f_{(\tilde{M},{\bf M}'),loc})\in \oplus_{(\tilde{M},{\bf M}')\in {\cal E}_{+}(\tilde{G},{\bf a})}SI({\bf M}')_{{\cal O}_{\tilde{M}'},loc}$ qui v\'erifient les conditions (1), (2) et (3) de 4.11. Ces conditions conservent un sens pour nos espaces "localis\'es". On note $I^{\cal E}(\tilde{G}(F),\omega)_{{\cal O},loc}$ la projection naturelle de $I^{\cal E}_{+}(\tilde{G}(F),\omega)_{{\cal O},loc}$ dans $\oplus_{{\bf G}'\in {\cal E}(\tilde{G},{\bf a})}SI({\bf G}')_{{\cal O}_{\tilde{G}'},loc}$. Il y a un diagramme naturel de localisation
  $$\begin{array}{ccc}I^{\cal E}_{+}(\tilde{G}(F),\omega)&\to&I^{\cal E}(\tilde{G}(F),\omega)\\ \downarrow&&\downarrow\\I^{\cal E}_{+}(\tilde{G}(F),\omega)_{{\cal O},loc}&\to&I^{\cal E}(\tilde{G}(F),\omega)_{{\cal O},loc}\\ \end{array}$$
  qui est commutatif. Montrons que
  
  (1) les fl\`eches verticales de ce diagramme sont surjectives.
  
  Par d\'efinition, les fl\`eches horizontales le sont. Il suffit donc de prouver que la fl\`eche verticale de gauche l'est. Soit $(f_{(\tilde{M},{\bf M}'),loc})\in I^{\cal E}_{+}(\tilde{G}(F),\omega)_{{\cal O},loc}$. On rel\`eve chaque $f_{(\tilde{M},{\bf M}'),loc}$ en un \'el\'ement $f_{(\tilde{M},{\bf M}')}\in SI({\bf M}')$. On peut remplacer cet \'el\'ement par la moyenne de ses images par l'action de $Aut(\tilde{M},{\bf M}')$. Cela nous permet de supposer que $f_{(\tilde{M},{\bf M}')}$ est invariant par ce groupe.    Soient ${\bf G}'$ et $(\tilde{M},{\bf M}')$ v\'erifiant les hypoth\`eses de la condition (2) de 4.11. Fixons des donn\'ees auxiliaires $G'_{1}$,...,$\Delta_{1}$. Cette condition affirme l'\'egalit\'e $S^{\tilde{G}'_{1}}(\delta_{1},f_{{\bf G}'})=S^{\tilde{M}'_{1}}(\delta_{1},f_{(\tilde{M},{\bf M}')})$ pour tout $\delta_{1}\in \tilde{M}'_{1}(F)$ assez r\'egulier. Elle n'est pas forc\'ement v\'erifi\'ee par les fonctions que l'on vient d'introduire. Mais, parce que  la famille de d\'epart appartient \`a $ I^{\cal E}_{+}(\tilde{G}(F),\omega)_{{\cal O},loc}$, elle l'est si l'image $\delta$ de $\delta_{1}$ dans $\tilde{M}'(F)$ est assez proche de ${\cal O}_{\tilde{M}'}$. Fixons un voisinage $\tilde{V}$ de ${\cal O}$ dans $\tilde{G}(F)$, ouvert et ferm\'e  et tel que $\tilde{V}\cap \tilde{G}_{ss}(F)$ soit invariant par conjugaison stable (un tel voisinage existe, cf. 4.6). De m\^eme que de ${\cal O}$, on a d\'eduit ${\cal O}_{\tilde{M}'}$, de $\tilde{V}$ se d\'eduit un voisinage $\tilde{V}_{\tilde{M}'}$ de ${\cal O}_{\tilde{M}'}$ dans $\tilde{M}'(F)$. Rempla\c{c}ons  chaque fonction $f_{(\tilde{M},{\bf M}')}$ par son produit avec la fonction caract\'eristique de $\tilde{V}_{\tilde{M}'}$. Si $\tilde{V}$ est assez petit, alors l'\'egalit\'e d'int\'egrales orbitales ci-dessus est v\'erifi\'ee pour tout $\delta_{1}$, autrement dit la condition 4.11(2) est satisfaite. Un m\^eme raisonnement s'applique \`a la condition 4.11(3). Donc la famille $(f_{(\tilde{M},{\bf M}')})$ appartient \`a $I^{\cal E}_{+}(\tilde{G}(F),\omega)$. Cela prouve (1).
  
  Il y a un diagramme commutatif naturel de localisation
  $$(2) \qquad \begin{array}{ccc}I(\tilde{G}(F),\omega)&\stackrel{tr}{\to}&\oplus_{{\bf G}'\in {\cal E}(\tilde{G},{\bf a})}SI({\bf G}')\\ \downarrow&&\downarrow\\ I(\tilde{G}(F),\omega)_{{\cal O},loc}&\stackrel{tr_{loc}}{\to}&\oplus_{{\bf G}'\in {\cal E}(\tilde{G},{\bf a})}SI({\bf G}')_{{\cal O}_{{\bf G}'},loc}\\ \end{array}$$
  o\`u $tr$ est le transfert. D'apr\`es la proposition 4.11, l'image de $tr$ est $I^{\cal E}(\tilde{G}(F),\omega)$. Gr\^ace \`a (1), celle de $tr_{loc}$ est donc $I^{\cal E}(\tilde{G}(F),\omega)_{{\cal O},loc}$. Montrons que
  
  (3) l'homomorphisme $tr_{loc}$ est injectif.
  
  Soit $f\in I(\tilde{G}(F),\omega)$ dont l'image dans $I(\tilde{G}(F),\omega)_{{\cal O},loc}$ appartient au noyau de $tr_{loc}$. Les int\'egrales orbitales de $f$ en des \'el\'ements fortement r\'eguliers se calculent par inversion de Fourier \`a partir des int\'egrales orbitales stables des fonctions $f^{{\bf G}'}$ pour ${\bf G}'\in {\cal E}(\tilde{G},{\bf a})$. On a expliqu\'e cela en 4.9(5) pour les \'el\'ements elliptiques   mais cela vaut pour tout \'el\'ement puisque tout \'el\'ement est elliptique dans un espace de Levi convenable. L'hypoth\`ese implique donc que $I^{\tilde{G}}(\gamma,\omega,f)=0$ pour tout $\gamma\in \tilde{G}_{reg}(F)$ assez proche de ${\cal O}$.  Par d\'efinition, cela signifie que l'image de $f$ dans $I(\tilde{G}(F),\omega)_{{\cal O},loc}$ est nulle. Cela prouve (3).
  
  La commutativit\'e du diagramme (2) entra\^{\i}ne que le transfert "dual", restreint \`a  $\oplus_{{\bf G}'\in {\cal E}(\tilde{G},{\bf a})}D_{g\acute{e}om}^{st}({\bf G}',{\cal O}_{\tilde{G}'})$, se factorise par le dual $$tr_{loc}^* :\oplus_{{\bf G}'\in {\cal E}(\tilde{G},{\bf a})}D_{g\acute{e}om}^{st}({\bf G}',{\cal O}_{\tilde{G}'})\to D_{g\acute{e}om}({\cal O},\omega)$$
  de $tr_{loc}$. L'assertion (3) entra\^{\i}ne que $tr_{loc}^*$ est surjective.    Posons pour simplifier $X_{(\tilde{M},{\bf M}')}=D_{g\acute{e}om}^{st}({\bf M}',{\cal O}_{\tilde{M}'})$ pour tout $(\tilde{M},{\bf M}')\in {\cal E}_{+}(\tilde{G},{\bf a})$, $X_{+}=\oplus_{(\tilde{M},{\bf M}')\in {\cal E}_{+}(\tilde{G},{\bf a})}X_{(\tilde{M},{\bf M}')}$, $X=\oplus_{{\bf G}'\in {\cal E}(\tilde{G},{\bf a})}X_{{\bf G}'}$, $Y=\oplus_{(\tilde{M},{\bf M}')\in {\cal E}_{+}(\tilde{G},{\bf a}), \tilde{M}\not=\tilde{G}}X_{(\tilde{M},{\bf M}')}$, $I=I^{\cal E}(\tilde{G}(F),\omega)_{{\cal O},loc}$, $I_{+}=I^{\cal E}_{+}(\tilde{G}(F),\omega)_{{\cal O},loc}$. Le noyau de $tr_{loc}^*$ est l'annulateur de $I $ dans $X$. Puisque $I$ est la projection sur $\oplus_{{\bf G}'\in {\cal E}(\tilde{G},{\bf a})}SI({\bf G}')_{{\cal O}_{{\bf G}'},loc}$ de $I_{+}$, cet annulateur est l'intersection avec $X$ de l'annulateur de $I_{+}$ dans $X_{+}$. L'espace $I_{+}$ est d\'efini par diff\'erentes conditions qui d\'efinissent chacune des sous-espaces. Son annulateur est la somme des annulateurs de ces sous-espaces. La condition 4.11(1) (ou plut\^ot son analogue localis\'ee) fournit l'annulateur
  
  (4) $X_{(\tilde{M},{\bf M}')}^{non-inv}=D^{st}_{g\acute{e}om}({\bf M}',{\cal O}_{\tilde{M}'})^{non-inv}\subset X_{(\tilde{M},{\bf M}')}$ 
  
 \noindent  pour tout $(\tilde{M},{\bf M}')\in {\cal E}_{+}(\tilde{G},{\bf a})$. La condition 4.11(2)  founit pour annulateur l'image de l'application
  $$\begin{array}{ccc}X_{(\tilde{M},{\bf M}')}&\to&X_{(\tilde{M},{\bf M}')}\oplus X_{{\bf G}'}\\ \delta&\mapsto&(\delta,-\delta^{\tilde{G}'})\\ \end{array}$$
  pour $(\tilde{M},{\bf M}')\in {\cal E}_{+}(\tilde{G},{\bf a})$ et ${\bf G}'\in {\cal E}(\tilde{G},{\bf a})$ tel que $M'$ est un Levi propre de $G'$ (pour \^etre correct, il faut choisir des donn\'ees auxiliaires pour d\'efinir l'application ci-dessus). La somme de l'espace (4) avec cette image est aussi la somme de cet espace (4) et des deux espaces suivants:
  
  (5) l'image par $\delta\mapsto \delta^{\tilde{G}'}$ de $X_{(\tilde{M},{\bf M}')}^{non-inv}$; cette image est ${\cal I}_{\tilde{M}'}^{\tilde{G}'}({\cal O}_{\tilde{M}'})^{non-inv}\subset X_{{\bf G}'}$;
  
  (6) l'image de l'application
  $$\begin{array}{ccc}X_{(\tilde{M},{\bf M}')}^{inv}=D^{st}_{g\acute{e}om}({\bf M}',{\cal O}_{\tilde{M}'})^{Aut(\tilde{M},{\bf M}')}&\to&X_{(\tilde{M},{\bf M}')}\oplus X_{{\bf G}'}\\ \delta&\mapsto&(\delta,-\delta^{\tilde{G}'}).\\ \end{array}$$
  
  La condition 4.11(3) fournit pour annulateur l'espace
  
  (7) ${\cal I}_{\tilde{M}'}^{\tilde{G}'}({\cal O}_{\tilde{M}'})\subset X_{{\bf G}'}$,
  
\noindent   pour tout ${\bf G}'\in {\cal E}(\tilde{G},{\bf a})$ et tout Levi $M'$ de $G'$ qui n'est pas relevant. Les espaces (7) sont les m\^emes qu'en  5.7(3). Les espaces (4) pour $\tilde{M}=\tilde{G}$ ou (5) pour $\tilde{M}\not=\tilde{G}$ sont les m\^emes qu'en 5.7(4).  Ces espaces sont inclus dans $X$. Il reste \`a prouver que l'intersection avec $X$ de la somme des espaces (6) et (4) pour $\tilde{M}\not=\tilde{G}$ est  la somme des espaces 5.7(5). Un \'el\'ement de cette intersection est une somme sur $(\tilde{M},{\bf M}')\in {\cal E}_{+}(\tilde{G},{\bf a})$, $\tilde{M}\not=\tilde{G}$, de termes
   $$x_{(\tilde{M},{\bf M}')}=\delta^{non-inv}+\sum_{i=1,...,n}(\delta_{i},-\delta_{i}^{{\bf G}'_{i}}),$$
   o\`u $\delta^{non-inv}\in X_{(\tilde{M},{\bf M}')}^{non-inv}$, $\delta_{i}\in X_{(\tilde{M},{\bf M}')}^{inv}$ pour tout $i$ et o\`u on a not\'e ${\bf G}'_{1}$,...,${\bf G}'_{n}$ les \'el\'ements de ${\cal E}(\tilde{G},{\bf a})$ dont $M'$ est un Levi (ces termes ne sont pas forc\'ement distincts, cf. la remarque suivant 5.7(5)). Fixons $(\tilde{M},{\bf M}')$ et  projetons sur $X_{(\tilde{M},{\bf M}')}$. Cette projection doit \^etre nulle. Cela entra\^{\i}ne que la projection de $x_{(\tilde{M},{\bf M}')}$ est nulle. Avec les notations ci-dessus, on a $\delta^{non-inv}=0$ et $\sum_{i=1,...,n}\delta_{i}=0$. Alors 
   $$x_{(\tilde{M},{\bf M}')}=\sum_{i=1,...,n-1}\left((\delta_{1}+...+\delta_{i})^{{\bf G}'_{i+1}}-(\delta_{1}+...+\delta_{i})^{{\bf G}'_{i}}\right)$$
   qui appartient \`a la somme des espaces 5.7(5). La r\'eciproque est claire. Cela ach\`eve la preuve. $\square$
 
\bigskip
\subsection{Preuve dans le cas archim\'edien}
  On suppose $F={\mathbb R}$ ou ${\mathbb C}$.     Pour unifier les notations, on pose $K\tilde{G}=\tilde{G}$ si $F={\mathbb C}$.   On fixe des mesures de Haar sur tous les groupes qui interviennent. On d\'efinit l'espace $I^{\cal E}_{+}(\tilde{G}(F),\omega)_{{\cal O},loc}\subset \oplus_{(\tilde{M},{\bf M}')\in {\cal E}_{+}(\tilde{G},{\bf a})}SI({\bf M}')_{{\cal O}_{\tilde{M}'}}$ comme dans le cas non-archim\'edien mais on le note plut\^ot $I^{\cal E}_{+}(K\tilde{G}(F),\omega)_{{\cal O},loc}$. On note $I^{\cal E}(K\tilde{G}(F),\omega)_{{\cal O},loc}$ sa projection dans $\oplus_{{\bf G}'\in {\cal E}(\tilde{G},{\bf a})}SI({\bf G}')_{{\cal O}_{\tilde{G}'},loc}$. Remarquons que ces espaces, ainsi que les espaces non localis\'es $I^{\cal E}_{+}(K\tilde{G}(F),\omega)$ et $I^{\cal E}(K\tilde{G}(F),\omega)$,  qui sont d\'efinis comme sous-espaces de certains espaces topologiques, sont ferm\'es dans ceux-ci. On a un diagramme naturel de localisation
 $$\begin{array}{ccc}I^{\cal E}_{+}(K\tilde{G}(F),\omega)&\to&I^{\cal E}(K\tilde{G}(F),\omega)\\ \downarrow&&\downarrow\\I^{\cal E}_{+}(K\tilde{G}(F),\omega)_{{\cal O},loc}&\to&I^{\cal E}(K\tilde{G}(F),\omega)_{{\cal O},loc}\\ \end{array}$$
  qui est commutatif. Montrons que
  
  (1) les fl\`eches verticales de ce diagramme sont surjectives.

   Il suffit de prouver que celle de gauche l'est. On a une filtration sur $I^{\cal E}_{+}(K\tilde{G}(F),\omega)$ dont le gradu\'e est d\'ecrit par 4.12(2). En fait, on a prouv\'e que les inclusions de cette relation \'etaient des \'egalit\'es. Le m\^eme proc\'ed\'e d\'efinit une filtration sur $I^{\cal E}_{+}(K\tilde{G}(F),\omega)_{{\cal O},loc}$ et on a
  $$(2) \qquad Gr\,I^{\cal E}(K\tilde{G}(F),\omega)_{{\cal O},loc}\subset \oplus_{(\tilde{M},{\bf M}')\in{\cal E}_{+}(\tilde{G},{\bf a})}  SI_{cusp}({\bf M}')_{{\cal O}_{\tilde{M}'},loc}^{Aut(\tilde{M},{\bf M}')} .$$
  
  {\bf Remarque.} Cette description est facile \`a condition d'utiliser pour les espaces de droite leur "deuxi\`eme" d\'efinition, cf. 5.4(3). Mais d'apr\`es la propri\'et\'e 5.5(2), la premi\`ere d\'efinition convient aussi bien.
  
  La fl\`eche verticale de gauche est compatible aux filtrations et d\'efinit une fl\`eche 
  $$Gr\,I^{\cal E}(K\tilde{G}(F),\omega)\to Gr\,I^{\cal E}(K\tilde{G}(F),\omega)_{{\cal O},loc}.$$
  D'apr\`es 4.12(2) (qui est une \'egalit\'e) et 5.5(2), $Gr\,I^{\cal E}(K\tilde{G}(F),\omega)$ s'envoie surjectivement sur le membre de droite de (2). Il en r\'esulte que l'homomorphisme ci-dessus entre gradu\'es est surjectif. Donc la fl\`eche verticale de gauche du diagramme est aussi surjective. $\square$
  
  Remarquons que ce raisonnement prouve aussi que (2) est une \'egalit\'e.
  
    Il y a un diagramme naturel de localisation
  $$ \begin{array}{ccc}I(K\tilde{G}(F),\omega)&\stackrel{tr}{\to}&\oplus_{{\bf G}'\in {\cal E}(\tilde{G},{\bf a})}SI({\bf G}')\\ \downarrow&&\downarrow\\ I(K\tilde{G}(F),\omega)_{{\cal O},loc}&\stackrel{tr_{loc}}{\to}&\oplus_{{\bf G}'\in {\cal E}(\tilde{G},{\bf a})}SI({\bf G}')_{{\cal O}_{{\bf G}'},loc}\\ \end{array}$$
  Gr\^ace \`a (1) et \`a la proposition 4.11, l'image de $tr_{loc}$ est $I^{\cal E}(K\tilde{G}(F),\omega)_{{\cal O},loc}$. On a
  
  (3) l'homomorphisme $tr_{loc}$ d\'efinit un hom\'eomorphisme de $I(K\tilde{G}(F))_{{\cal O},loc}$ sur $I^{\cal E}(K\tilde{G}(F))_{{\cal O},loc}$.
  
  Preuve. L'homomorphisme $tr$ se calcule par une formule explicite comme on en a utilis\'e en 4.13. Il r\'esulte de cette formule que $tr$ est continue pourvu que les facteurs de transfert soient des fonctions $C^{\infty}$. Or cela r\'esulte du lemme 2.8. Donc $tr$ est continue. Il en r\'esulte que $tr_{loc}$ l'est aussi. Soit $\varphi_{{\tilde{\cal T}}}\in I(K\tilde{G}(F),\omega)$. Supposons que  son image dans   $I^{\cal E}(K\tilde{G}(F),\omega)_{{\cal O},loc}$ soit nulle. L'\'el\'ement $tr(\varphi_{\tilde{{\cal T}}})$ a un d\'eveloppement infinit\'esimal nul en tout point correspondant \`a un \'el\'ement de ${\cal O}$. Par une formule d'inversion g\'en\'eralisant 4.9(5) au cas non elliptique, la fonction $\varphi_{\tilde{{\cal T}}}$ a elle-m\^eme un d\'eveloppement infinit\'esimal nul en tout \'el\'ement de ${\cal O}$. Donc son image dans $I(K\tilde{G}(F),\omega)_{{\cal O},loc}$ est nulle. Cela prouve que $tr_{loc}$ est injectif. Donc $tr_{loc}$ est une bijection continue de $I(K\tilde{G}(F),\omega)_{{\cal O},loc}$ sur $I^{\cal E}(K\tilde{G}(F),\omega)_{{\cal O},loc}$. Or ces deux espaces sont des espaces de Fr\'echet. Une telle bijection est donc n\'ecessairement ouverte. $\square$

  Gr\^ace \`a (3), l'application duale
 $$tr_{loc}^*: \oplus_{{\bf G}'\in {\cal E}(\tilde{G},{\bf a})}D_{g\acute{e}om}^{st}({\bf G}',{\cal O}_{\tilde{G}'}) \to D_{g\acute{e}om}({\cal O},\omega) $$
 se quotiente en un isomorphisme de l'espace de d\'epart quotient\'e par l'annulateur de $I^{\cal E}(K\tilde{G}(F),\omega)_{{\cal O},loc}$ sur l'espace d'arriv\'ee. Il reste \`a prouver que cet annulateur est la somme des espaces d\'ecrits avant l'\'enonc\'e de la proposition 5.7. Le m\^eme raisonnement que dans le cas non-archim\'edien nous ram\`ene \`a prouver que l'annulateur de $I^{\cal E}_{+}(K\tilde{G}(F),\omega)_{{\cal O},loc}$ est la somme des espaces d\'ecrits en 5.8(4), (5), (6) et (7). Notons  $Ann$ l'annulateur de $I^{\cal E}_{+}(K\tilde{G}(F),\omega)$ et $Ann^?$ la somme de ces espaces . L'espace $I^{\cal E}_{+}(K\tilde{G}(F),\omega)_{{\cal O},loc}$ est intersection  finie de sous-espaces  et  $Ann^?$ n'est autre que la somme des annulateurs de ces sous-espaces. Mais, \`a cause de la topologie, il n'est pas compl\`etement \'evident que l'annulateur de l'intersection soit la somme des annulateurs. On va le prouver.
 
 Consid\'erons d'abord le cas o\`u  $(G,\tilde{G},{\bf a})$ est quasi-d\'eploy\'e et \`a torsion int\'erieure. L'espace $SI(\tilde{G}(F))_{{\cal O},loc}$ est inclus dans $I^{\cal E}(\tilde{G}(F))_{{\cal O},loc}$ (il correspond \`a la donn\'ee maximale ${\bf G}$). Restreinte \`a ce sous-espace, l'inclusion (2), dont on a prouv\'e que c'\'etait une \'egalit\'e, donne une \'egalit\'e
 $$(4)\qquad Gr\,SI(\tilde{G}(F))_{{\cal O},loc}=\oplus_{M\in \underline{\cal L}}SI_{cusp}(\tilde{M}(F))_{{\cal O}_{\tilde{M}},loc}^{W(M)},$$
 o\`u $\underline{\cal L}$ est un ensemble de repr\'esentants des classes de conjugaison  de Levi. L'application naturelle du terme de gauche dans celui de droite est continue. Puisque nos ensembles sont des espaces de Fr\'echet, c'est un hom\'eomorphisme. Pour tout  Levi $M$, notons ${\cal I}_{\tilde{M}}^{\tilde{G}}({\cal O}_{\tilde{M}})$ l'image de l'homomorphisme
 $$\begin{array}{ccc}D^{st}_{g\acute{e}om}({\cal O}_{\tilde{M}})&\to& D^{st}_{g\acute{e}om}({\cal O})\\ \delta&\mapsto&\delta^{\tilde{G}}\\ \end{array}$$
 Notons $ {\cal I}^{\tilde{G}}({\cal O})$ la somme de ces espaces ${\cal I}_{\tilde{M}}^{\tilde{G}}({\cal O}_{\tilde{M}})$ pour $M\not=G$. On peut se limiter aux $M\in \underline{\cal L}$. Le fait que (4) soit un hom\'eomorphisme implique que $ {\cal I}^{\tilde{G}}({\cal O})$ est l'annulateur dans $D^{st}_{g\acute{e}om}({\cal O})$ du sous-espace $SI_{cusp}(\tilde{G})_{{\cal O},loc}$.  Le m\^eme r\'esultat vaut pour tout  Levi $M$. L'action du groupe $W(M)$ pr\'eserve ${\cal I}^{\tilde{M}}({\cal O}_{\tilde{M}})$.  Fixons un suppl\'ementaire $D^{st}_{g\acute{e}om,cusp}({\cal O}_{\tilde{M}})$ de ce sous-espace, invariant par l'action de ce groupe, notons ${\cal I}_{\tilde{M},cusp}^{\tilde{G}}({\cal O}_{\tilde{M}})^{inv}$ son image dans $D^{st}_{g\acute{e}om}({\cal O})$ par l'application ci-dessus. Par dualit\'e, on d\'eduit  de (4) l'\'egalit\'e
 $$D^{st}_{g\acute{e}om}({\cal O})=\oplus_{M\in \underline{\cal L}}{\cal I}_{\tilde{M},cusp}^{\tilde{G}}({\cal O}_{\tilde{M}})^{inv}.$$

 Revenons au cas g\'en\'eral. Ce que l'on vient de dire s'adapte aux espaces $SI({\bf M}')$ pour $(\tilde{M},{\bf M}')\in {\cal E}_{+}(\tilde{G},{\bf a})$, munis cette fois de l'action de $Aut(\tilde{M},{\bf M}')$. En particulier, on fixe un sous-espace $ X_{(\tilde{M},{\bf M}'),cusp}\subset X_{(\tilde{M},{\bf M}')}= D^{st}_{g\acute{e}om}({\bf M}',{\cal O}_{\tilde{M}'})$, qui est un suppl\'ementaire de la somme des espaces induits \`a partir  de Levi propres de $M'$ et qui est invariant par $Aut(\tilde{M},{\bf M}')$. On note $X_{(\tilde{M},{\bf M}'),cusp}^{inv}$ son sous-espace des invariants par ce groupe. Cet espace s'identifie \`a celui des formes lin\'eaires continues sur $SI_{cusp}({\bf M}')_{{\cal O}_{\tilde{M}'},loc}^{Aut(\tilde{M},{\bf M}')}$. Pour les m\^emes raisons que ci-dessus, la bijection (2) est un isomorphisme. Par dualit\'e, on en d\'eduit que le sous-espace
 $$X_{++}=\oplus_{(\tilde{M},{\bf M}')\in {\cal E}_{+}(\tilde{G},{\bf a})}X_{(\tilde{M},{\bf M}'),cusp}^{inv}\subset X_{+}=\oplus_{(\tilde{M},{\bf M}')\in {\cal E}_{+}(\tilde{G},{\bf a})}X_{(\tilde{M},{\bf M}')}$$
 s'identifie par restriction \`a l'espace des formes lin\'eaires continues sur $I^{\cal E}_{+}(K\tilde{G}(F))_{{\cal O},loc}$.  En particulier $Ann\cap X_{++}=\{0\}$. Il est clair que $Ann^?$ est inclus dans $Ann$. Pour prouver que cette inclusion est une \'egalit\'e, il suffit de prouver que $X_{+}=X_{++}+Ann^?$. On d\'emontre par r\'ecurrence descendante sur le corang de $\tilde{M}$ que $X_{(\tilde{M},{\bf M}')}$ est inclus dans $X_{++}+Ann^?$. Fixons $(\tilde{M},{\bf M}')$. L'espace $X_{(\tilde{M},{\bf M}')}$ est somme de $X_{(\tilde{M},{\bf M}'),cusp}^{inv}$, de son suppl\'ementaire $X_{(\tilde{M},{\bf M}'),cusp}^{non-inv}$ conserv\'e par $Aut(\tilde{M},{\bf M}')$ dans $X_{(\tilde{M},{\bf M}'),cusp}$ et des sous-espaces obtenus par induction \`a partir  de Levi propres de $M'$. Le premier espace $X_{(\tilde{M},{\bf M}'),cusp}^{inv}$ est contenu dans $X_{++}$. Le deuxi\`eme $X_{(\tilde{M},{\bf M}'),cusp}^{non-inv}$ est inclus dans $Ann^?$ (5.8(4)). Fixons un  Levi propre $R'\subset M'$ et des donn\'ees auxiliaires pour ${\bf M}'$. Soit $\delta\in D^{st}_{g\acute{e}om,\lambda_{1}}(\tilde{R}'_{1}(F),{\cal O}_{\tilde{R}'})$. On veut prouver que son image $\delta^{\tilde{M}'}$ par induction appartient \`a $X_{++}+Ann^?$. Supposons d'abord $R'$ relevant.   Il lui est associ\'e un \'el\'ement $(\tilde{R},{\bf R}')\in {\cal E}_{+}(\tilde{G},{\bf a})$ et $\delta$ s'identifie \`a un \'el\'ement de $D^{st}_{g\acute{e}om}({\bf R}',{\cal O}_{\tilde{R}'})$. Fixons ${\bf G}'\in {\cal E}(\tilde{G},{\bf a})$ dont un Levi s'identifie \`a $M'$.  En utilisant 5.8 (5) et (6) pour $(\tilde{M},{\bf M}')$ et pour $(\tilde{R},{\bf R}')$, on voit que les deux \'el\'ements $\delta^{\tilde{M}'}-(\delta^{\tilde{M}'})^{\tilde{G}'}$ et $\delta-\delta^{\tilde{G}'}$ appartiennent \`a $Ann^?$. Par transitivit\'e de l'induction, $(\delta^{\tilde{M}'})^{\tilde{G}'}=\delta^{\tilde{G}'}$. Donc $\delta-\delta^{\tilde{M}'}\in X_{++}+Ann^?$. Par hypoth\`ese de r\'ecurrence, $\delta$ appartient \`a $X_{++}+Ann^?$. Donc aussi $\delta^{\tilde{M}'}$. Supposons maintenant $R'$ non relevant. On a de nouveau $\delta^{\tilde{M}'}-\delta^{\tilde{G}'}\in X_{++}+Ann^?$. Mais $\delta^{\tilde{G}'}$ appartient \`a $Ann^?$ (5.8 (7)). Donc $\delta^{\tilde{M}'}\in X_{++}+Ann^?$. Cela ach\`eve la preuve. $\square$
 \bigskip
 
 \subsection{Localisation}
 Fixons un \'el\'ement semi-simple $\eta\in \tilde{G}(F)$ et un voisinage $\mathfrak{u}$ de $0$ dans $\mathfrak{g}_{\eta}(F)$ ayant les m\^emes propri\'et\'es qu'en 4.1. Avec les notations de ce paragraphe (et en r\'etablissant les espaces de mesures), on a d\'efini une application 
 $$desc_{\eta}^{\tilde{G}}:I(\tilde{U},\omega)\otimes Mes(G(F))\to I(U_{\eta},\omega)\otimes Mes(G_{\eta}(F)).$$
 Il s'en d\'eduit une application duale entre espaces de distributions. Pour s'affranchir de l'ensemble $\mathfrak{u}$ qui complique les notations, nous noterons
 $$desc_{\eta}^{\tilde{G},*}:D_{g\acute{e}om}(G_{\eta}(F),\omega)\otimes Mes(G_{\eta}(F))^*\to D_{g\acute{e}om}(\tilde{G}(F),\omega)\otimes Mes(G(F))^*$$
 cette application duale, \'etant entendue qu'elle n'est d\'efinie que pour des distributions dont le support dans $G_{\eta}(F)$ est assez voisin de $1$. Notons ${\cal O}$ la classe de conjugaison de $\eta$ dans $\tilde{G}(F)$. On a d\'efini l'espace $D_{g\acute{e}om}({\cal O},\omega)$. En appliquant la m\^eme d\'efinition en rempla\c{c}ant $\tilde{G}$ par $G_{\eta}$ et ${\cal O}$ par la classe de conjugaison r\'eduite \`a $\{1\}$, on obtient un espace que l'on note plut\^ot $D_{unip}(G_{\eta}(F),\omega)$. L'application ci-dessus se restreint en une application surjective
 $$desc_{\eta}^{\tilde{G},*}:D_{unip}(G_{\eta}(F),\omega)\otimes Mes(G_{\eta}(F))^*\to D_{g\acute{e}om}({\cal O},\omega)\otimes Mes(G(F))^*.$$
 Plus pr\'ecisement, cette application se factorise en
 $$D_{unip}(G_{\eta}(F),\omega)\otimes Mes(G_{\eta}(F))^*\stackrel{p_{\eta}}{\to} D_{unip}(G_{\eta}(F),\omega)^{Z_{G}(\eta;F)}\otimes Mes(G_{\eta}(F))^*$$
 $$\stackrel{desc_{\eta}^{\tilde{G},*}}{\simeq}D_{g\acute{e}om}({\cal O},\omega)\otimes Mes(G(F))^*,$$
 o\`u $p_{\eta}$ est la projection naturelle sur l'espace des invariants (rappelons que l'action naturelle de $Z_{G}(\eta;F)$ tient compte du caract\`ere $\omega$). 
 
 Supposons $(G,\tilde{G},{\bf a})$ quasi-d\'eploy\'e et \`a torsion int\'erieure. On consid\`ere un \'el\'ement semi-simple $\eta\in \tilde{G}(F)$ tel que $G_{\eta}$ soit quasi-d\'eploy\'e. On note ${\cal O}$ sa classe de conjugaison stable et on pose $\Xi_{\eta}=Z_{G}(\eta)/G_{\eta}$. On a de m\^eme une application lin\'eaire
 $$desc_{\eta}^{st,\tilde{G},*}:D_{g\acute{e}om}^{st}(G_{\eta}(F))\otimes Mes(G_{\eta}(F))^*\to D_{g\acute{e}om}^{st}(\tilde{G}(F))\otimes Mes(G(F))^*.$$
 Elle se restreint en une application
  $$D_{unip}^{st}(G_{\eta}(F))\otimes Mes(G_{\eta}(F))^*\stackrel{p^{st}_{\eta}}{\to} D^{st}_{unip}(G_{\eta}(F))^{\Xi_{\eta}^{\Gamma_{F}}}\otimes Mes(G_{\eta}(F))^*$$
 $$\stackrel{desc_{\eta}^{st, \tilde{G},*}}{\simeq}D^{st}_{g\acute{e}om}({\cal O})\otimes Mes(G(F))^*.$$
 {\bf Attention.} L'application $desc_{\eta}^{st,\tilde{G},*}$ n'est pas la restriction de $desc_{\eta}^{\tilde{G},*}$ \`a l'espace des distributions stables. La preuve du lemme 4.8 fournit la relation entre ces deux applications. On a
 $$desc_{\eta}^{st,\tilde{G},*}=\sum_{y\in \dot{{\cal Y}}(\eta)}desc_{\eta[y]}^{\tilde{G},*}\circ transfert_{y},$$
 o\`u $transfert_{y}:D_{g\acute{e}om}^{st}(G_{\eta}(F))\to D_{g\acute{e}om}(G_{\eta[y]}(F))$ est le transfert d\'eduit du torseur int\'erieur $ad_{y}:G_{\eta[y]}\to G_{\eta}$. 
 
 \bigskip
 
 \subsection{Induction et classes de conjugaison stable}
 Soient $\tilde{M}$ un espace de Levi de $\tilde{G}$ et $\eta$ un \'el\'ement semi-simple de $\tilde{M}(F)$. On a d\'efini le groupe $I_{\eta}$ et l'ensemble ${\cal Y}(\eta)$ en 4.6. En rempla\c{c}ant $\tilde{G}$ par $\tilde{M}$, on d\'efinit de m\^eme un groupe et un ensemble que l'on note $I^M_{\eta}$ et    ${\cal Y}^M(\eta)$. Remarquons que
 
 (1)  $I^M_{\eta}=I_{\eta}\cap M$. 
 
 Preuve. On a l'\'egalit\'e $Z(M)^{\theta}=Z(M)^{\theta,0}Z(G)^{\theta}$ et l'inclusion $Z(M)^{\theta,0}
\subset M_{\eta}$. Donc $I^M_{\eta}=Z(M)^{\theta}M_{\eta}= Z(G)^{\theta}M_{\eta}\subset I_{\eta}\cap M$. L'inclusion oppos\'ee provient de l'\'egalit\'e $G_{\eta}\cap M=M_{\eta}$. $\square$ 
 
 Il r\'esulte de (1) que  ${\cal Y}^M(\eta)={\cal Y}(\eta)\cap M$. On en d\'eduit une application naturelle
 $$(2) \qquad I^M_{\eta}\backslash {\cal Y}^M(\eta)/ M(F)\to I_{\eta}\backslash {\cal Y}(\eta)/ G(F).$$
 On note $\underline{{\cal Y}}^M(\eta)$ et $\underline{{\cal Y}}(\eta)$ les ensembles de doubles classes ci-dessus. 

 \ass{Lemme}{L'application (2) est injective. Pour $y\in {\cal Y}(\eta)$, l'image de $y$ dans $  \underline{{\cal Y}}(\eta)$ appartient \`a l'image de cette application si et seulement si le Levi $M_{\eta}$ de $G_{\eta}$ se transf\`ere par le torseur int\'erieur $ad_{y^{-1}}$ en un Levi de $G_{\eta[y]}$. Plus pr\'ecis\'ement, soit $y\in {\cal Y}(\eta)$ dont l'image dans  $  \underline{ {\cal Y}}(\eta)$ n'appartient pas \`a l'image de (2). Soit $T$ un sous-tore maximal de $M_{\eta}$ d\'efini sur $F$. Alors le tore $T$ ne se transf\`ere pas par le torseur int\'erieur $ad_{y^{-1}}$ en un sous-tore maximal de $G_{\eta[y]}$ d\'efini sur $F$.}
 
 Preuve.  Soient $y,y'\in {\cal Y}^M(\eta)$ dont les images dans $ \underline{ {\cal Y}}(\eta)$ sont \'egales. On doit prouver que leurs images dans $\underline{ {\cal Y}}^M(\eta)$ le sont aussi. L'\'el\'ement $(y')^{-1}y$ appartient \`a ${\cal Y}^M(\eta[y'])$. Son image dans 
  $ \underline{ {\cal Y}}(\eta[y'])$ est \'egale \`a celle de $1$. On v\'erifie qu'il suffit de prouver que les images de $(y')^{-1}y$ et de $1$ dans $\underline{ {\cal Y}}^M(\eta[y'])$ sont \'egales. Quitte \`a remplacer $\eta$ par $\eta[y']$, on est ramen\'e au probl\`eme initial avec cette fois $y'=1$. A $y$, on associe le cocycle  $\sigma\mapsto y\sigma(y)^{-1}$ de $\Gamma_{F}$ dans $I^M_{\eta}$. L'hypoth\`ese signifie que ce cocycle, pouss\'e en un cocycle \`a valeurs dans $I_{\eta}$ est un cobord. La conclusion est que ce cocycle lui-m\^eme est un cobord.  Il suffit de prouver que le noyau $K$ de l'application
  $$H^1(\Gamma_{F}; I^M_{\eta})\to H^1(\Gamma_{F};I_{\eta})$$
  est r\'eduit \`a $\{1\}$. Remarquons que, dans le cas o\`u $F$ est archim\'edien, les ensembles ci-dessus ne sont pas des groupes. Le noyau est l'ensemble des \'el\'ements de $H^1(\Gamma_{F}; I^M_{\eta})$ qui s'envoient sur l'\'el\'ement trivial de $H^1(\Gamma_{F};I_{\eta})$. Le centre $Z(I_{\eta})$ de $I_{\eta}$ est \'egal \`a $Z(G)^{\theta}Z(G_{\eta})$ et on a un diagramme commutatif
  $$\begin{array}{ccccc}H^1(\Gamma_{F};Z(I_{\eta}))&\to& H^1(\Gamma_{F};I_{\eta})&\to&H^1(\Gamma_{F};G_{\eta,AD})\\ \parallel&&\uparrow&&\uparrow\\ H^1(\Gamma_{F};Z(I_{\eta}))&\to& H^1(\Gamma_{F};I^M_{\eta})&\to&H^1(\Gamma_{F};M_{\eta,ad})\\ \end{array}$$
  Les suites horizontales sont exactes. Parce que $M_{\eta,ad}$ est un Levi de $G_{\eta,AD}$, la derni\`ere fl\`eche verticale est injective. Il en r\'esulte que $K$ est l'image dans $ H^1(\Gamma_{F};I^M_{\eta})$ du noyau $C$ de l'application $H^1(\Gamma_{F};Z(I_{\eta}))\to H^1(\Gamma_{F};I_{\eta})$. Un \'el\'ement $c\in C$ est un cocycle de la forme $\sigma\mapsto c(\sigma)=x\sigma(x)^{-1}$, o\`u $x$ est un \'el\'ement de $I_{\eta}$ dont l'image $x_{ad}$ dans $G_{\eta,AD}$ appartient \`a $G_{\eta,AD}(F)$. Notons $\pi_{ad}:G_{\eta}\to G_{\eta,AD}$ la projection naturelle. Puisque $M_{\eta,ad}$ est un Levi de $G_{\eta,AD}$, la projection naturelle $M_{\eta,ad}(F)\to G_{\eta,AD}(F)/\pi_{ad}(G_{\eta}(F))$ est surjective. Quitte \`a multiplier $x$ \`a droite par un \'el\'ement de $G_{\eta}(F)$, on peut donc supposer $x_{ad}\in M_{\eta,ad}(F)$. Alors $x\in I^M_{\eta}$ et l'image du cocycle $c$ dans $ H^1(\Gamma_{F};I^M_{\eta})$ est un bord. Cela d\'emontre que  l'image $K$  de $C$ dans $ H^1(\Gamma_{F};I^M_{\eta})$ est r\'eduite \`a $\{1\}$, d'o\`u l'injectivit\'e de l'application (2). 
  
  Pour $y\in {\cal Y}^M(\eta)$, l'image de $M_{\eta}$ par $ad_{y^{-1}}$ est $M_{\eta[y]}$. C'est un  Levi de $G_{\eta[y])}$ (c'est-\`a-dire qu'il est d\'efini sur $F$)  et $M_{\eta}$ se transf\`ere en un tel Levi. Il en r\'esulte plus g\'en\'eralement que, pour $y\in {\cal Y}(\eta)$, si l'image de $y$ dans
$\underline{{\cal Y}}(\eta)$ appartient \`a l'image de l'application (2), le Levi  $M_{\eta}$ de $G_{\eta}$ se transf\`ere par le torseur int\'erieur $ad_{y^{-1}}$ en un Levi de $G_{\eta[y]}$. Soit maintenant $y\in {\cal Y}(\eta)$ et $T$ un sous-tore maximal de $M_{\eta}$. Supposons que $T$ se transf\`ere par $ad_{y^{-1}}$ en un sous-tore maximal de $G_{\eta[y]}$ d\'efini sur $F$. Cela signifie que, quitte \`a multiplier \`a gauche $y$ par un \'el\'ement de $G_{\eta}$, le tore $T_{y}=ad_{y^{-1}}(T)$ est d\'efini sur $F$ et la restriction $ad_{y^{-1}}:T\to T_{y}$ est \'equivariante pour les actions galoisiennes. Il en r\'esulte que $ad_{y^{-1}}$ se restreint en un isomorphisme d\'efini sur $F$ de $A_{T}$ sur $A_{T_{y}}$. Notons $\tilde{R}$ et $\tilde{R}_{y}$ les commutants de $A_{T}$ et $A_{T_{y}}$ dans $\tilde{G}$. Fixons un \'el\'ement $x_{*}$ en position g\'en\'erale dans $X_{*}(A_{T})$. Il d\'etermine un espace parabolique $\tilde{S}\in {\cal P}(\tilde{R})$: $A_{T}$ agit dans $\mathfrak{u}_{S}$  par des caract\`eres $\alpha$ tels que $<\alpha,x_{*}>>0$. A $ad_{y^{-1}}(x_{*})$ est de m\^eme associ\'e un espace parabolique $\tilde{S}_{y}\in {\cal P}(\tilde{R}_{y})$. Alors $ad_{y^{-1}}$ envoie la paire $(\tilde{S},\tilde{R})$ sur $(\tilde{S}_{y},\tilde{R}_{y})$. On sait que deux telles paires d\'efinies sur $F$ qui sont conjugu\'ees par un \'el\'ement de $G(\bar{F})$ le sont aussi par un \'el\'ement de $G(F)$. Quitte \`a multiplier $y$ \`a droite par un \'el\'ement de $G(F)$, on peut donc supposer que les deux paires paraboliques sont \'egales. Cela entra\^{\i}ne $y\in R$. Mais $A_{\tilde{M}}\subset A_{M_{\eta}}\subset A_{T}\subset A_{\tilde{R}}$, donc $R\subset M$ et $y\in {\cal Y}(\eta)\cap M={\cal Y}^M(\eta)$. Cela d\'emontre la derni\`ere assertion de l'\'enonc\'e. Enfin, soit $y\in {\cal Y}(\eta)$, supposons que $M_{\eta}$ se transf\`ere par le torseur int\'erieur $ad_{y^{-1}}$ en un Levi $M_{y}$ de $G_{\eta[y]}$. On choisit un tore maximal $T$ de $M_{\eta}$, d\'efini sur $F$ et elliptique si $F$ est non-archim\'edien, resp. fondamental si $F$ est archim\'edien. Alors $T$ se transf\`ere en un tore maximal d\'efini sur $F$ de $M_{y}$, a fortiori de $G_{\eta[y]}$. D'apr\`es ce que l'on vient de d\'emontrer, l'image de $y$ dans $\underline{{\cal Y}}(\eta)$ appartient \`a l'image de l'application (2). Cela ach\`eve la preuve. $\square$

\bigskip

\subsection{Un r\'esultat de r\'eduction}
On conserve la m\^eme situation. On note ${\cal O}$ la classe de conjugaison stable de $\eta$ dans $\tilde{M}(F)$ et ${\cal O}^{\tilde{G}}$ sa classe de conjugaison stable dans $\tilde{G}(F)$. Remarquons qu'en g\'en\'eral, ${\cal O}$ est plus petit que l'intersection ${\cal O}^{\tilde{G}}\cap \tilde{M}(F)$. Notons $N$ le groupe des $x\in G(F)$ tels que $ad_{x}$ conserve $\tilde{M}$ et ${\cal O}$. Ce groupe agit naturellement sur $D_{g\acute{e}om}({\cal O},\omega)$ via son quotient fini $N/M(F)$. On note $p_{N}$ la projection naturelle sur le sous-espace des invariants par $N$.  Si $(G,\tilde{G},{\bf a})$ est quasi-d\'eploy\'e et \`a torsion int\'erieure, $N$ agit aussi sur $D_{g\acute{e}om}^{st}({\cal O})$. On note $p_{N}^{st}$ la projection sur le sous-espace des invariants

\ass{Lemme}{On suppose $G_{\eta}=M_{\eta}$ et $A_{\tilde{M}}=A_{M_{\eta}}$. 

(i) L'application (2) de 5.11 est bijective.

(ii) La restriction \`a $D_{g\acute{e}om}({\cal O},\omega)\otimes Mes(M(F))^*$ de l'application d'induction de $\tilde{M}$ \`a $\tilde{G}$ se factorise en
$$D_{g\acute{e}om}({\cal O},\omega)\otimes Mes(M(F))^*\stackrel{p_{N}}{\to}D_{g\acute{e}om}({\cal O},\omega)^N\otimes Mes(M(F))^*\simeq D_{g\acute{e}om}({\cal O}^{\tilde{G}},\omega)\otimes Mes(G(F))^*.$$

(iii) Supposons $(G,\tilde{G},{\bf a})$ quasi-d\'eploy\'e et \`a torsion int\'erieure. La restriction \`a $D_{g\acute{e}om}^{st}({\cal O},\omega)\otimes Mes(M(F))^*$ de l'application d'induction de $\tilde{M}$ \`a $\tilde{G}$ se factorise en
$$D_{g\acute{e}om}^{st}({\cal O},\omega)\otimes Mes(M(F))^*\stackrel{p^{st}_{N}}{\to}D_{g\acute{e}om}^{st}({\cal O},\omega)^N\otimes Mes(M(F))^*\simeq D_{g\acute{e}om}^{st}({\cal O}^{\tilde{G}},\omega)\otimes Mes(G(F))^*.$$}

Preuve. L'hypoth\`ese $G_{\eta}=M_{\eta}$ entra\^{\i}ne $I_{\eta}=I^M_{\eta}$. Un \'el\'ement $y\in {\cal Y}(\eta)$ d\'efinit un cocycle $\sigma\mapsto y\sigma(y)^{-1}$ \`a valeurs dans $I^M_{\eta}$ dont l'image dans $H^1(\Gamma_{F}; G)$ est triviale. Mais l'application $H^1(\Gamma_{F};M)\to H^1(\Gamma_{F};G)$ est injective. Donc l'image du cocycle ci-dessus dans $H^1(\Gamma_{F};M)$ est triviale. Cela signifie que l'on peut \'ecrire $y=y'g$, avec $g\in G(F)$ et $y'\in M$. N\'ecessairement, $y'\in {\cal Y}^M(\eta)$, donc l'image dans $\underline{{\cal Y}}(\eta)$ de $y$ appartient \`a l'image de l'application (2) de 5.11. D'o\`u la surjectivit\'e de cette application et sa bijectivit\'e d'apr\`es le lemme pr\'ec\'edent. 

Introduisons le groupe $Z_{\eta}=Z_{G}(\eta)\cap {\cal Y}(\eta)$ et son quotient $\underline{Z}_{\eta}=Z_{\eta}/I_{\eta}$. Le groupe $Z_{\eta}$ agit sur ${\cal Y}(\eta)$ par multiplication \`a gauche. On v\'erifie que l'ensemble de doubles classes
$$\underline{{\cal X}}(\eta)=Z_{\eta}\backslash {\cal Y}(\eta)/G(F)$$
param\`etre les classes de conjugaison par $G(F)$ dans ${\cal O}^{\tilde{G}}$. En rempla\c{c}ant $\tilde{G}$ par $\tilde{M}$, on a de m\^eme un ensemble
$$\underline{{\cal X}}^M(\eta)=Z^M_{\eta}\backslash {\cal Y}^M(\eta)/M(F)$$
qui param\`etre les classes de conjugaison par $M(F)$ dans ${\cal O}$. L'assertion (i) d\'ej\`a prouv\'ee entra\^{\i}ne que l'application naturelle
$$\underline{{\cal X}}^M(\eta)\to \underline{{\cal X}}(\eta)$$
est surjective. On peut donc fixer 
un ensemble de repr\'esentants $\dot{{\cal X}}(\eta)$ de $\underline{{\cal X}}(\eta)$ qui est inclus dans ${\cal Y}^M(\eta)$. Fixons aussi un ensemble de repr\'esentants $\dot{{\cal X}}^M(\eta)$ de $\underline{{\cal X}}(\eta)$. L'application pr\'ec\'edente devient une application surjective
$$q: \dot{{\cal X}}^M(\eta)\to \dot{{\cal X}}(\eta).$$
Pour tout $y\in \dot{{\cal X}}^M(\eta)$, on fixe $z_{y}\in Z_{\eta}$ et $g_{y}\in G(F)$ tels que $y=z_{y}q(y)g_{y}$. 
Remarquons que, pour un \'el\'ement $y$ de l'un ou l'autre de ces ensembles, les \'egalit\'es $G_{\eta}=M_{\eta}$ et $A_{\tilde{M}}=A_{M_{\eta}}$ et le fait que  $y\in M$ entra\^{\i}nent que $G_{\eta[y]}=M_{\eta[y]}$ et $A_{\tilde{M}}=A_{M_{\eta[y]}}$. On pose $D[y]=D_{unip}(M_{\eta[y]},\omega)$ et on note $\zeta_{y}:D[y]\to D[y]^{Z_{G}(\eta[y];F)}$ la projection naturelle. En oubliant pour simplifier les espaces de mesures, la description de 5.10 fournit des isomorphismes
$$D_{g\acute{e}om}({\cal O},\omega)=\oplus_{y\in \dot{{\cal X}}^M(\eta)}D[y]^{Z_{M}(\eta[y];F)},$$
$$D_{g\acute{e}om}({\cal O}^{\tilde{G}},\omega)=\oplus_{y\in \dot{{\cal X}}(\eta)}D[y]^{Z_{G}(\eta[y];F)}.$$
Modulo ces isomorphismes, l'application d'induction se d\'ecrit de la fa\c{c}on suivante. A $(d_{y})_{y\in \dot{{\cal X}}^M(\eta)}\in \oplus_{y\in \dot{{\cal X}}^M(\eta)}D[y]^{Z_{M}(\eta[y];F)}$, elle associe 
$(d'_{y'})_{y'\in \dot{{\cal X}}(\eta)}\in \oplus_{y'\in \dot{{\cal X}}(\eta)}D[y]^{Z_{G}(\eta[y];F)}$, o\`u
$$d'_{y'}=\zeta_{y'}(\sum_{y\in q^{-1}(y')}\omega(g_{y})^{-1}ad_{g_{y}}(d_{y})).$$
On voit que cette application est surjective. D'autre part, l'application d'induction est insensible \`a l'action par conjugaison (tordue par le caract\`ere $\omega$) de tout \'el\'ement de $G(F)$ conservant $\tilde{M}$. Elle se factorise donc par la projection $p_{N}$. Pour obtenir (ii), il reste \`a prouver que l'application d'induction
$$D_{g\acute{e}om}({\cal O},\omega)^N\to D_{g\acute{e}om}({\cal O}^{\tilde{G}},\omega)$$
est injective. A l'aide de la description ci-dessus, cela  r\'esulte de la propri\'et\'e suivante. Soit $(d_{y})_{y\in \dot{{\cal X}}^M(\eta)}\in \oplus_{y\in \dot{{\cal X}}^M(\eta)}D[y]^{Z_{M}(\eta[y];F)}$. Supposons cet \'el\'ement invariant par $N$. Soit $y'\in \dot{{\cal X}}(\eta)$. Alors

(1) l'\'el\'ement $\omega(g_{y})^{-1}ad_{g_{y}}(d_{y})$ est ind\'ependant de $y\in q^{-1}(y')$ et il est invariant par $Z_{G}(\eta[y'];F)$. 

On ne perd rien \`a supposer que $y'=1$ et que $y=1$ appartient \`a $q^{-1}(1)$. Soit $y\in q^{-1}(1)$. Alors $ad_{g_{y}}(\eta[y])=\eta$, donc aussi $ad_{g_{y}}(G_{\eta[y]})=G_{\eta}$. Puisque $g\in G(F)$, $ad_{g_{y}}$ envoie $A_{G_{\eta[y]}}$ sur $A_{G_{\eta}}$. Mais ces deux tores sont \'egaux \`a $A_{\tilde{M}}$. Donc $ad_{g_{y}}$ conserve $A_{\tilde{M}}$ et aussi son commutant $\tilde{M}$. Puisque $ad_{g_{y}}$ envoie $\eta[y]$ sur $\eta$, il conserve la classe de conjugaison stable commune ${\cal O}$ de ces deux \'el\'ements. Donc $g_{y}\in N$. L'hypoth\`ese d'invariance par $N$ entra\^{\i}ne l'\'egalit\'e $\omega(g_{y})^{-1}ad_{g_{y}}(d_{y})=d_{1}$, d'o\`u la premi\`ere assertion de (1). Le m\^eme argument que ci-dessus montre que $Z_{G}(\eta;F)\subset N$. L'hypoth\`ese d'invariance par $N$ entra\^{\i}ne que $d_{1}$ est invariant par $Z_{G}(\eta;F)$. Cela d\'emontre (1) et le (ii) de la proposition. 

Pour le (iii), quitte \`a  changer l'\'el\'ement $\eta$  de ${\cal O}$, on peut supposer $G_{\eta}$ quasi-d\'eploy\'e. La description de 5.10 identifie $D_{g\acute{e}om}^{st}({\cal O}^{\tilde{G}})$ \`a $D^{st}_{unip}(M_{\eta}(F))^{\Xi_{\eta}^{\Gamma_{F}}}$ et $D_{g\acute{e}om}^{st}({\cal O})$ \`a $D^{st}_{unip}(M_{\eta}(F))^{\Xi_{\eta}^{M,\Gamma_{F}}}$. L'application d'induction n'est autre que la projection sur l'espace d'invariants par $\Xi_{\eta}^{\Gamma_{F}}$. Elle est surjective. De nouveau, cette application se factorise par $p_{N}^{st}$ et il reste \`a prouver que cette application d'induction est injective sur $D_{g\acute{e}om}^{st}({\cal O},\omega)^N$. Mais on vient de prouver qu'elle \'etait injective sur l'espace plus gros $D_{g\acute{e}om}({\cal O},\omega)^N$. D'o\`u l'assertion, ce qui ach\`eve la d\'emonstration. $\square$

\bigskip

\subsection{Induction et stabilit\'e}
On suppose $(G,\tilde{G},{\bf a})$ quasi-d\'eploy\'e et \`a torsion int\'erieure. Soient $\tilde{M}$ un espace de Levi de $\tilde{G}$ et $({\cal O}_{j})_{j=1,...,k}$ une famille finie de classes de conjugaison stable semi-simples dans $\tilde{M}(F)$.Rappelons que l'on note $\boldsymbol{\gamma}\mapsto \boldsymbol{\gamma}^{\tilde{G}}$ l'homomorphisme d'induction de $D_{g\acute{e}om}(\tilde{M}(F))\otimes Mes(M(F))^*$ dans $D_{g\acute{e}om}(\tilde{G}(F))\otimes Mes(G(F))^*$.

\ass{Lemme}{Soit $\boldsymbol{\gamma}\in \sum_{j=1,...,k}D_{g\acute{e}om}({\cal O}_{j})\otimes Mes(M(F))^*$. Supposons que $\boldsymbol{\gamma}^{\tilde{G}}$ soit stable. Alors il existe $\boldsymbol{\delta}\in \sum_{j=1,...,k}D^{st}_{g\acute{e}om}({\cal O}_{j})\otimes Mes(M(F))^*$ telle que $\boldsymbol{\delta}^{\tilde{G}}=\boldsymbol{\gamma}^{\tilde{G}}$.}

Preuve. On fixe des mesures de Haar pour se d\'ebarrasser des espaces de mesures. Pour tout $j$, notons ${\cal O}_{j}^{\tilde{G}}$ la classe de conjugaison stable dans $\tilde{G}(F)$ qui contient ${\cal O}_{j}$. On peut regrouper les classes ${\cal O}_{j}$ selon ces classes ${\cal O}_{j}^{\tilde{G}}$. C'est-\`a-dire que l'on peut fixer une famille $({\cal O}'_{l})_{l=1,...,m}$ de classes de conjugaison stable semi-simples dans $\tilde{G}(F)$, distinctes deux-\`a-deux, et une application surjective $q:\{1,...,k\}\to \{1,...,m\}$ de sorte que ${\cal O}_{j}^{\tilde{G}}={\cal O}'_{q(j)}$ pour tout $j=1,...,k$. On peut \'ecrire $\boldsymbol{\gamma}=\sum_{l=1,...,m}\boldsymbol{\gamma}_{l}$, avec $\boldsymbol{\gamma}_{l}\in \sum_{j\in q^{-1}(l)}D_{g\acute{e}om}({\cal O}_{j})$. Alors $\sum_{l=1,...,m}\boldsymbol{\gamma}_{l}^{\tilde{G}}$ est stable. Mais les distributions $\boldsymbol{\gamma}_{l}^{\tilde{G}}$ sont support\'ees par des classes de conjugaison stable distinctes. Il r\'esulte des constructions de 4.6 qu'alors, chaque $\boldsymbol{\gamma}_{l}^{\tilde{G}}$ est stable. Pour r\'esoudre notre probl\`eme, il suffit de trouver pour chaque $l$ une distribution $\boldsymbol{\delta}_{l}\in \sum_{j\in q^{-1}(l)}D^{st}_{g\acute{e}om}({\cal O}_{j})$ telle que $\boldsymbol{\delta}_{l}^{\tilde{G}}=\boldsymbol{\gamma}^{\tilde{G}}$. Cela nous ram\`ene au probl\`eme initial, avec l'hypoth\`ese suppl\'ementaire que chacune des classes ${\cal O}_{j}$ engendre la m\^eme classe de conjugaison stable dans $\tilde{G}(F)$. Nous faisons d\'esormais cette hypoth\`ese et nous posons simplement ${\cal O}^{\tilde{G}}={\cal O}_{j}^{\tilde{G}}$ pour tout $j=1,...,k$. 

On fixe $\eta\in {\cal O}^{\tilde{G}}$ tel que $G_{\eta}$ soit quasi-d\'eploy\'e et on fixe une paire de Borel  \'epingl\'ee ${\cal E}_{\eta}$   de $G_{\eta}$ d\'efinie sur $F$, de paire de Borel $(B_{\eta},T)$. Pour tout $j=1,...,k$, on fixe $\eta_{j}\in {\cal O}_{j}$ tel que $M_{\eta_{j}}$ (donc aussi $G_{\eta_{j}}$) soit quasi-d\'eploy\'e et on fixe une paire de Borel \'epingl\'ee ${\cal E}_{\eta_{j}}$  de $G_{\eta_{j}}$ d\'efinie sur $F$, de paire de Borel  $(B_{\eta_{j}},T_{j})$, de sorte  que $M_{\eta_{j}}$ soit standard. Puisque $\eta_{j}\in {\cal O}^{\tilde{G}}$, on peut fixer $y_{j}\in {\cal Y}(\eta)$ de sorte que $\eta_{j}=\eta[y_{j}]$. L'automorphisme $ad_{y_{j}}$ se restreint en un torseur int\'erieur de $G_{\eta_{j}}$ sur $G_{\eta}$. Quitte \`a multiplier $y_{j}$ \`a gauche par un \'el\'ement de $I_{\eta}=G_{\eta}$, on peut supposer que ce torseur envoie ${\cal E}_{\eta_{j}}$ sur ${\cal E}_{\eta}$. Un tel torseur int\'erieur est alors un isomorphisme d\'efini sur $F$. Il se restreint en un isomorphisme d\'efini sur $F$ de $T_{j}$ sur $T$. Puisque $A_{M}\subset A_{M_{\eta}}\subset T$, le tore  $ad_{y_{j}}(A_{M})$ est d\'efini sur $F$ et l'application $ad_{y_{j}}:A_{M}\to ad_{y_{j}}(A_{M})$ est un isomorphisme d\'efini sur $F$. Notons $\tilde{M}_{j}$ le commutant de $ad_{y_{j}}(A_{M})$ dans $\tilde{G}$.  C'est un espace de Levi de $\tilde{G}$, on a $\eta\in \tilde{M}(F)$ et le groupe $M_{j,\eta}$ est standard pour ${\cal E}_{\eta}$ puisque c'est l'image par $ad_{y_{j}}$ de $M_{\eta_{j}}$. Le m\^eme raisonnement que dans la preuve du lemme 5.11 montre que $y_{j}$ se d\'ecompose en $g_{j}m_{j}$, avec $m_{j}\in M$ et $g_{j}\in G(F)$. On voit que $m_{j}^{-1}$ appartient \`a ${\cal Y}(\eta_{j})$, donc $ad_{m_{j}}(\eta_{j})\in {\cal O}_{j}$. Le groupe $G_{ad_{m_{j}}(\eta_{j})}$ est \'egal \`a $ad_{g_{j}^{-1}}(G_{\eta})$, donc est quasi-d\'eploy\'e. Quitte \`a remplacer $\eta_{j}$ par $ad_{m_{j}}(\eta_{j})$, on peut donc supposer $m_{j}=1$ et $y_{j}=g_{j}\in G(F)$. L'\'el\'ement $g_{j}$ conjugue $\tilde{M}$ en $\tilde{M}_{j}$, $\eta_{j}$ en $\eta$ et la classe ${\cal O}_{j}$ en la classe de conjugaison stable ${\cal O}'_{j}$ de $\eta$ dans $\tilde{M}_{j}(F)$. On peut \'ecrire $\boldsymbol{\gamma}=\sum_{j=1,...,k}\boldsymbol{\gamma}_{j}$, o\`u $\boldsymbol{\gamma}_{j}\in D_{g\acute{e}om}({\cal O}_{j})$. Pour tout $j$, notons $\boldsymbol{\gamma}'_{j}$ l'image de $\boldsymbol{\gamma}_{j}$ par $ad_{g_{j}}$. C'est un \'el\'ement de $D_{g\acute{e}om}({\cal O}'_{j})$. Il est clair que $\boldsymbol{\gamma}_{j}^{\tilde{G}}=\boldsymbol{\gamma}_{j}^{_{'}\tilde{G}}$. Donc $\sum_{j=1,...,k}\boldsymbol{\gamma}_{j}^{_{'}\tilde{G}}$ est stable. Supposons trouv\'ees des distributions stables $\boldsymbol{\delta}_{j}^{_{'}\tilde{G}}\in D^{st}_{g\acute{e}om}({\cal O}'_{j})$ de sorte que $\sum_{j=1,...,k}\boldsymbol{\gamma}_{j}^{_{'}\tilde{G}}=\sum_{j=1,...,k}\boldsymbol{\delta}_{j}^{_{'}\tilde{G}}$. Pour tout $j$, on note alors $\boldsymbol{\delta}_{j}$ l'image de $\boldsymbol{\delta}_{j}'$ par $ad_{g_{j}^{-1}}$. En inversant le calcul ci-dessus, on voit que la distribution $\boldsymbol{\delta}=\sum_{j=1,...,k}\boldsymbol{\delta}_{j}$ r\'esout notre probl\`eme. 

Oubliant notre probl\`eme initial pour simplifier les notations, on est ramen\'e au probl\`eme suivant. On consid\`ere une famille $(\tilde{M}_{j})_{j=1,...,k}$ d'espaces de Levi de $\tilde{G}$ tels que $\eta\in\tilde{M}_{j}(F)$. Pour tout $j$, on note ${\cal O}_{j}$ la classe de conjugaison stable de $\eta$ dans $\tilde{M}_{j}(F)$ et  on consid\`ere une distribution $\boldsymbol{\gamma}_{j}\in  D_{g\acute{e}om}({\cal O}_{j})$. On suppose que $\sum_{j=1,...,k}\boldsymbol{\gamma}_{j}^{\tilde{G}}$ est stable. On veut prouver qu'il existe pour tout $j$ une distribution $\boldsymbol{\delta}_{j}\in D_{g\acute{e}om}^{st}({\cal O}_{j})$ de sorte que $\sum_{j=1,...,k}\boldsymbol{\gamma}_{j}^{\tilde{G}}=\sum_{j=1,...,k}\boldsymbol{\delta}_{j}^{\tilde{G}}$.

Fixons un voisinage $\mathfrak{u}$ de $0$ dans $G_{\eta}(F)$ ayant les m\^emes propri\'et\'es qu'en 4.8. On pose $U_{\eta}=exp(\mathfrak{u})$ et on note $\tilde{U}$ l'ensemble des \'el\'ements de $\tilde{G}(F)$ dont la partie semi-simple est stablement conjugu\'ee \`a un \'el\'ement de $U_{\eta}\eta$. Pour tout $j=1,...,k$, on pose $U_{\eta,j}=U_{\eta}\cap M_{j,\eta}(F)$ et on note $\tilde{U}_{j}$ l'ensemble des \'el\'ements de $\tilde{M}_{j}(F)$ dont la partie semi-simple est stablement conjugu\'ee (dans $\tilde{M}_{j}$) \`a un \'el\'ement de $U_{\eta,j}\eta$. Consid\'erons le diagramme commutatif
$$\begin{array}{ccccccc}I(\tilde{U})&&&\stackrel{res}{\to}&&&\oplus_{j=1,...,k}I(\tilde{U}_{j})\\&\,\searrow \iota&&&&\,\swarrow \underline{\iota}&\\ &&I(\tilde{G}(F))_{{\cal O}^{\tilde{G}},loc}&\stackrel{res_{loc}}{\to}&\oplus_{j=1,...,k}I(\tilde{M}_{j}(F))_{{\cal O}_{j},loc}&&\\ s\downarrow\,&&s_{loc}\downarrow\,&&\underline{s}_{loc}\downarrow\,&&\,\downarrow\underline{s}\\&&SI(\tilde{G}(F))_{{\cal O}^{\tilde{G}},loc}&\stackrel{res^{st}_{loc}}{\to}&\oplus_{j=1,...,k}SI(\tilde{M}_{j}(F))_{{\cal O}_{j},loc}&&\\&\iota^{st}\nearrow\,&&&&\,\nwarrow\underline{\iota}^{st}&\\ SI(\tilde{U})&&&\stackrel{res^{st}}{\to}&&&\oplus_{j=1,...,k}SI(\tilde{U}_{j})\\ \end{array}$$

Les fl\`eches sont les applications naturelles. D\'ecrivons l'espace $I(\tilde{U})$. 
Pour $y\in {\cal Y}(\eta)$, il correspond \`a $\mathfrak{u}$ un voisinage $\mathfrak{u}_{y}$ de $0$ dans $G_{\eta[y]}$. Posons $U_{y}=exp(\mathfrak{u}_{y})$.  On pose $Z_{\eta}=Z_{G}(\eta)\cap {\cal Y}(\eta)$. Comme on l'a vu dans la preuve de 5.12, l'ensemble
$$\underline{{\cal X}}(\eta)=Z_{\eta}\backslash {\cal Y}(\eta)/G(F)$$
param\`etre l'ensemble des classe de conjugaison par $G(F)$ dans ${\cal O}^{\tilde{G}}$. Si on fixe un ensemble de repr\'esentants $\dot{{\cal X}}(\eta)$ de cet ensemble de doubles classes, la th\'eorie de la descente identifie $I(\tilde{U})$ \`a $\oplus_{y\in \dot{{\cal X}}(\eta)}I(U_{y})^{Z_{G}(\eta[y];F)}$. Fixons plut\^ot un ensemble de repr\'esentants $\dot{{\cal Y}}(\eta)$ de l'ensemble de doubles classes 
$$\underline{{\cal Y}}(\eta)=G_{\eta}\backslash {\cal Y}(\eta)/G(F).$$
Alors $I(\tilde{U})$ s'identifie au sous-espace des $(f_{y})_{y\in \dot{{\cal Y}}(\eta)}\in \oplus_{y\in \dot{{\cal X}}(\eta)}I(U_{y})$
qui v\'erifient la condition suivante:

(1) soient $y,y'\in \dot{{\cal Y}}(\eta)$ et $g\in G(F)$ tels que $ad_{g}(\eta[y])=\eta[y']$; alors $f_{y'}=ad_{g}(f_{y})$.

Remarquons que le quotient $Z_{\eta}/G_{\eta}$ est \'egal au groupe $\Xi_{\eta}^{\Gamma_{F}}$ de 4.8. Ce groupe agit sur $\underline{{\cal Y}}(\eta)$ par multiplication \`a gauche. Il s'en d\'eduit une action de ce groupe sur $\dot{{\cal Y}}(\eta)$ que l'on note $(\xi,y)\mapsto \xi\star y$. Le stabilisateur   dans $\Xi_{\eta}^{\Gamma_{F}}$ d'un \'el\'ement $y$ est l'image dans ce groupe de $ad_{y}(Z_{G}(\eta[y];F))$.  Comme on l'a vu en 4.6, le groupe $\Xi_{\eta}^{\Gamma_{F}}$ agit  sur $G_{\eta}$ par automorphismes d\'efinis sur $F$. Rappelons la construction.  Consid\'erons un \'el\'ement $z\in Z_{\eta}$. Quitte \`a multiplier $z$ \`a gauche par un \'el\'ement de $G_{\eta}$, on peut supposer que $ad_{z}$ conserve ${\cal E}_{\eta}$. L'\'el\'ement $z$ est alors bien d\'etermin\'e modulo multiplication \`a gauche par un \'el\'ement de $Z(G_{\eta})$ et on a $z\sigma(z)^{-1}\in Z(G_{\eta})$ pour tout $\sigma\in \Gamma_{F}$.  La restriction de $ad_{z}$ \`a $G_{\eta}$ est un automorphisme   de ce groupe qui est d\'efini sur $F$.   Cet automorphisme ne d\'epend que de l'image de $z$ dans $\Xi^{\Gamma_{F}}_{\eta}$. On note $ad_{\xi}$ l'automorphisme d\'etermin\'e par $\xi\in \Xi_{\eta}^{\Gamma_{F}}$. Posons $\dot{{\cal Y}}^0(\eta)= \dot{{\cal Y}}(\eta)\cap Z_{\eta}G(F)$. Les \'el\'ements de cet ensemble sont les $y\in \dot{{\cal Y}}(\eta)$ tels que $\eta[y]$ est conjugu\'e \`a $\eta$ par un \'el\'ement de $G(F)$. On impose \`a notre syst\`eme de repr\'esentants $\dot{{\cal Y}}(\eta)$ la condition

(2) supposons $y\in \dot{{\cal Y}}^0(\eta)$; alors $y$ est un \'el\'ement de $Z_{\eta}$ tel que $ad_{y}$ conserve ${\cal E}_{\eta}$. 

Il en r\'esulte que, pour un tel \'el\'ement $y$, on a $\eta[y]=\eta$ et, en notant $\xi_{y}$ l'image de $y$ dans $\Xi_{\eta}^{\Gamma_{F}}$, la restriction de $ad_{y}$ \`a $G_{\eta}$ co\"{\i}ncide avec $ad_{\xi_{y}}$. 

On d\'ecrit de fa\c{c}on similaire les espaces $I(\tilde{U}_{j})$ et on impose la m\^eme condition. On ajoute des indices $j$ pour les objets relatifs \`a ces espaces. D'apr\`es le lemme 5.11, il y a pour tout $j$ une injection $q_{j}=\dot{{\cal Y}}_{j}(\eta)\to \dot{{\cal Y}}(\eta)$ de sorte que, pour tout $y\in \dot{{\cal Y}}_{j}(\eta)$, il existe $x_{y}\in G_{\eta}$ et $g_{y}\in G(F)$ tels que $y=x_{y}q_{j}(y)g_{y}$. On fixe de tels \'el\'ements $x_{y}$ et $g_{y}$. Montrons que

(3) soit $y\in \dot{{\cal Y}}_{j}(\eta)$, supposons $q_{j}(y)\in  \dot{{\cal Y}}^0(\eta)$; alors  on peut supposer $x_{y}=1$. 

On a $x_{y}q_{j}(y)\sigma(q_{j}(y))^{-1}\sigma(x_{y})^{-1}=y\sigma(y)^{-1}$ pour tout $\sigma\in \Gamma_{F}$. D'apr\`es (2), le terme $q_{j}(y)\sigma(q_{j}(y))^{-1}$ appartient \`a $Z(G_{\eta})$. Donc, d'une part, il commute \`a $x_{y}$, d'autre part, il appartient \`a $M_{j,\eta}$, a fortiori \`a $M_{j}$. L'\'egalit\'e pr\'ec\'edente entra\^{\i}ne que $x_{y}\sigma(x_{y})^{-1}\in M_{j}$. Puisque c'est aussi un \'el\'ement de $G_{\eta}$, il appartient \`a $M_{j,\eta}$.  On obtient un cocycle $\sigma\mapsto x_{y}\sigma(x_{y})$ \`a valeurs dans $M_{j,\eta}$ qui devient un cobord dans $G_{\eta}$. Puisque $H^1(\Gamma_{F};M_{j,\eta})\to H^1(\Gamma_{F};G_{\eta})$ est injective, il existe $x'\in M_{j,\eta}$ et $g'\in G_{\eta}(F)$ tel que $x_{y}=x'g'$. On a alors $y=x'g'q_{j}(y)g_{y}=x'q_{j}(y)ad_{q_{j}(y)^{-1}}(g')g_{y}$. Puisque  $ad_{q_{j}(y)}$ est un automorphisme d\'efini sur $F$ de $G_{\eta}$, le terme 
$ad_{q_{j}(y)^{-1}}(g')g_{y}$ appartient \`a $G(F)$. On peut remplacer $y$ par $(x')^{-1}y$, $x_{y}$ par $1$ et $g_{y}$ par $ad_{q_{j}(y)^{-1}}(g')g_{y}$. Avec ces nouvelles d\'efinitions, on a $y=q_{j}(y)g_{y}$, ce qui d\'emontre (3).

  L'application $res$ du diagramme se d\'ecrit par
 $$(4) \qquad (f_{y})_{y\in \dot{{\cal Y}}(\eta)}\in I(\tilde{U}) \mapsto (f_{j,y})_{j=1,...,k,y\in \dot{{\cal Y}}_{j}(\eta)}\in \oplus_{j=1,...,k}I(\tilde{U}_{j})$$
  o\`u, pour tout $j$ et tout $y\in \dot{{\cal Y}}_{j}(\eta)$,   $f_{j,y}$ est l'image de  $ad_{g_{y}^{-1}}(f_{q_{j}(y)})$ par l'application $res_{M_{j,\eta[y]}}$. Rappelons que pour tout $y\in {\cal Y}(\eta)$, du torseur int\'erieur $ad_{y}$ se d\'eduit une application $transfert_{y}:I(\tilde{U}_{y})\to SI(\tilde{U}_{\eta})$. Soit $(f_{y})_{y\in \dot{{\cal Y}}(\eta)}\in I(\tilde{U})$. Pour tout $y\in \dot{{\cal Y}}(\eta)$ et tout  $\xi\in \Xi_{\eta}^{\Gamma_{F}}$, on a l'\'egalit\'e 
  
  $$(5) \qquad transfert_{\xi\star y}(f_{\xi\star y})=ad_{\xi}(transfert_{y}(f_{y})).$$

  A ce point, nous allons  s\'eparer les cas $F$ non-archim\'edien  et $F$ archim\'edien.  
 
 \bigskip
  
  \subsection{Suite de la preuve, cas $F$ non-archim\'edien}

On suppose $F$ non-archim\'edien.   On va prouver

(1) soit $f\in I(\tilde{U})$; supposons que l'image de $f$ dans $\oplus_{j=1,...,k}SI(\tilde{M}_{j}(F))_{{\cal O}_{j},loc}$ est nulle; alors il existe $f'\in I(\tilde{U})$ qui a m\^eme image que $f$ dans  $\oplus_{j=1,...,k}I(\tilde{M}_{j}(F))_{{\cal O}_{j},loc}$ et dont l'image dans $SI(\tilde{U}) $ est nulle.

Soit $f=(f_{y})_{y\in \dot{{\cal Y}}(\eta)}\in I(\tilde{U})$. On note $(f_{j,y})_{j=1,...,k, y\in \dot{{\cal Y}}_{j}(\eta)}$ son image dans $\oplus_{j=1,...,k}I(\tilde{U}_{j})$, cf. 5.13 (4). Supposons  
que l'image de $f$ dans $\oplus_{j=1,...,k}SI(\tilde{M}_{j}(F))_{{\cal O}_{j},loc}$ est nulle. Posons $\phi=\sum_{y\in \dot{{\cal Y}}(\eta)}transfert_{y}(f_{y})$. C'est un \'el\'ement de $SI(U_{\eta})$.  Montrons que

(2) pour tout $j=1,...,k$, l'image  $ \phi_{M_{j,\eta}}$ de $\phi$ dans $SI(U_{\eta,j})$ est nulle au voisinage de $0$.

Soit $j\in \{1,...,k\}$. Posons $\phi_{j}=\sum_{y\in \dot{{\cal Y}}_{j}(\eta)}transfert_{y}(f_{j,y})$. C'est un \'el\'ement de $SI(U_{\eta, j})$. D'apr\`es la description de 4.8, dire que l'image de $f$ dans $SI(\tilde{M}_{j}(F))_{{\cal O}_{j},loc}$ est nulle revient \`a dire que $\phi_{j}$ est nulle au  voisinage de $0$. Il suffit donc de prouver que $\phi_{j}= \phi_{M_{j,\eta}}$. Par commutation du transfert \`a la restriction, on voit que, pour tout $y\in \dot{{\cal Y}}_{j}(\eta)$, on a $transfert_{y}(f_{j,y})=(transfert_{q_{j}(y)}(f_{q_{j}(y)}))_{M_{j,\eta}}$.  D'autre part,  pour $y\in \dot{{\cal Y}}(\eta)$ qui n'appartient pas \`a l'image de $q_{j}$, aucun sous-tore maximal de $M_{j,\eta}$ ne se transf\`ere \`a $G_{\eta[y]}$, cf. lemme 5.11. Il en r\'esulte que $ (transfert_{y}(f_{y}))_{M_{j,\eta}}=0$. Cela d\'emontre l'\'egalit\'e $\phi_{j}=\phi_{M_{j,\eta}}$ et (2).

Quitte \`a multiplier $f$ par la fonction caract\'eristique d'un voisinage ouvert et ferm\'e de ${\cal O}^{\tilde{G}}$ invariant par conjugaison stable (c'est-\`a-dire tel qu'en 4.6) et assez petit, ce qui ne change pas l'image de $f$ dans $\oplus_{j=1,...,k}I(\tilde{M}_{j}(F))_{{\cal O}_{j},loc}$, on peut donc supposer que $\phi_{M_{j,\eta}}=0$. On dispose d'une action de $\Xi_{\eta}^{\Gamma_{F}}$ sur  $G_{\eta}(F)$, donc aussi sur $I(G_{\eta}(F))$ et $SI(G_{\eta}(F))$. On a aussi une action de $G_{\eta,AD}(F)$. Les deux actions se combinent en une action du produit semi-direct  $H_{\eta}=G_{\eta,AD}(F)\rtimes \Xi_{\eta}^{\Gamma_{F}}$. 
On sait que $\phi$ est invariant par $\Xi_{\eta}^{\Gamma_{F}}$, cf. lemme 4.8. On retrouve d'ailleurs ce r\'esultat en utilisant 5.13(5).   D'autre part,  les classes de conjugaison stable dans $G_{\eta}(F)$ d'\'el\'ements fortement r\'eguliers sont invariantes par l'action de $G_{\eta,AD}(F)$. Il en r\'esulte que $U_{\eta}$ est invariant par   $G_{\eta,AD}(F)$  et que l'action de ce groupe $G_{\eta,AD}(F)$ sur $SI(G_{\eta}(F))$ est triviale. Donc $\phi$ est invariant par $H_{\eta}$. Cela entra\^{\i}ne que $\phi_{ad_{h}(M_{j,\eta})}=0$ pour tout $j$ et tout $h\in H_{\eta}$.  L'action de $H_{\eta}$
 sur $I(G_{\eta}(F))$ se factorise par l'action d'un groupe fini puisque l'image de  $G_{\eta}(F)$ dans $G_{\eta,AD}(F)$ agit trivialement. Il en r\'esulte que l'on peut relever $\phi$ en un \'el\'ement $\varphi\in I(U_{\eta})$ qui est invariant par $H_{\eta}$. Cet \'el\'ement v\'erifie: l'image de $\varphi_{ad_{h}(M_{j,\eta})}$ dans $SI(ad_{h}(M_{j,\eta}(F)))$ est nulle pour tout $j=1,...,k$ et tout $h\in H_{\eta}$. Pour la m\^eme raison que ci-dessus, l'ensemble des Levi intervenant dans cette relation est fini modulo conjugaison par $G_{\eta}(F)$. On peut donc appliquer le 4.16: il existe $\varphi_{0}\in I^{inst}(G_{\eta}(F))$ tel que $\varphi_{0,ad_{h}(M_{j,\eta})}=\varphi_{ad_{h}(M_{j,\eta})}$ pour tous $j,h$. On peut moyenner $\varphi_{0}$ sous l'action de $H_{\eta}$ et supposer $\varphi_{0}$ invariant par ce groupe. On peut aussi remplacer $\varphi_{0}$ par son produit avec la fonction caract\'eristique de $U_{\eta}$ et supposer $\varphi_{0}\in I(U_{\eta})$. 
 
 Notons $N$ le nombre d'\'el\'ements de $  \dot{{\cal Y}}^0(\eta)$. D\'efinissons une famille $f'=(f'_{y})_{y\in \dot{{\cal Y}}(\eta)}\in \oplus_{y\in \dot{{\cal Y}}(\eta)}I(U_{y})$ par
 $f'_{y}=f_{y}$ pour $y\not\in  \dot{{\cal Y}}^0(\eta)$ et $f'_{y}=f_{y}+\frac{1}{N}(\varphi_{0}-\varphi)$ pour $y\in  \dot{{\cal Y}}^0(\eta)$. Remarquons qu'en vertu de l'hypoth\`ese 5.13(2), on a $\eta[y]=\eta$ et $U_{y}=U_{\eta}$ pour $y\in \dot{{\cal Y}}^0(\eta)$. Nos fonctions appartiennent bien \`a l'espace indiqu\'e. Montrons que
 
 (3) la famille $f'$ appartient \`a $I(\tilde{U})$.
 
 On doit v\'erifier la condition 5.13(1). Soient $y,y'\in \dot{{\cal Y}}(\eta)$ et $g\in G(F)$ tels que $ad_{g}(\eta[y])=\eta[y']$. Ces conditions entra\^{\i}nent que $y\in  \dot{{\cal Y}}^0(\eta)$ si et seulement si $y'\in \dot{{\cal Y}}^0(\eta)$. Supposons d'abord que $y,y'\not\in \dot{{\cal Y}}^0(\eta)$. Alors la condition $ad_{g}(f'_{y})=f'_{y'}$ r\'esulte de la condition initiale $ad_{g}(f_{y})=f_{y'}$. Supposons maintenant $y,y'\in \dot{{\cal Y}}^0(\eta)$. Dans ce cas $\eta[y]=\eta[y']$, donc $g\in Z_{G}(\eta;F)$. En vertu de la condition initiale $ad_{g}(f_{y})=f_{y'}$, il nous suffit de prouver que $\varphi$ et $\varphi_{0}$ sont invariantes par $ad_{g}$. Puisque ces fonctions sont invariantes par $H_{\eta}$, il suffit de prouver qu'il existe $h\in H_{\eta}$ tel que $ad_{g}=ad_{h}$. Or $Z_{G}(\eta;F)\subset Z_{\eta}$. On peut donc trouver $x\in G_{\eta}$ et $z\in Z_{\eta}$ de sorte que $g=xz$ et $ad_{z}$ conserve ${\cal E}_{\eta}$. On a $ad_{g}=ad_{x}ad_{z}$. On a $ad_{z}=ad_{\xi}$,  o\`u $\xi$ est l'image de $z$ dans $\Xi_{\eta}^{\Gamma_{F}}$. Puisque $ad_{g}$ et $ad_{\xi}$ sont d\'efinis sur $F$, $ad_{x}$ aussi, ce qui implique que l'image de $x$ dans $G_{\eta,AD}$ appartient \`a $G_{\eta,AD}(F)$. On a bien d\'ecompos\'e $ad_{g}$ en produit de l'action d'un \'el\'ement de $G_{\eta,AD}(F)$ et d'un \'el\'ement de $\Xi_{\eta}^{\Gamma_{F}}$. Cela prouve (3).

On a

(4) l'image de $f'$ dans  $SI(\tilde{U})$ est nulle. 

Posons $\phi'=\sum_{y\in \dot{{\cal Y}}(\eta)}transfert_{y}(f'_{y})$. C'est un \'el\'ement de $SI(U_{\eta})$. En vertu de 4.8, il s'agit de prouver que $\phi'=0$. Par d\'efinition, 
$$\phi'=\phi+\frac{1}{N}\sum_{y\in \dot{{\cal Y}}^0(\eta)}transfert_{y}(\varphi_{0}-\varphi).$$
Rappelons que l'image de $\varphi$ dans $SI(U_{\eta})$ est $\phi$. Pour $y\in \dot{{\cal Y}}^0(\eta)$, l'image de $transfert_{y}(\varphi)$ est $\xi_{y}(\phi)$, qui est \'egale \`a $\phi$ puisque $\phi$ est invariant par $\Xi_{\eta}^{\Gamma_{F}}$. L'image de $\varphi_{0}$ dans $SI(U_{\eta})$ est nulle, et celle de $transfert_{y}(\varphi_{0})$ est l'image de la pr\'ec\'edente par $\xi_{y}$, donc est nulle. L'\'egalit\'e ci-dessus entra\^{\i}ne $\phi'=0$, d'o\`u (4).

Montrons que

(5) pour tout $j=1,...,k$, $f$ et $f'$ ont m\^eme image dans $ I(\tilde{M}_{j}(F))_{{\cal O}_{j},loc}$.

Par 5.13(4), la famille $f'$ d\'efinit une famille $(f'_{j,y})_{j=1,...,k, y\in \dot{{\cal Y}}_{j}(\eta)}$. On doit prouver que $f_{j,y}=f'_{j,y}$ pour tous $j$, $y$. Fixons $j$ et $y$. Alors $f_{j,y}$ et $f'_{j,y}$ sont les images de $ad_{g_{y}^{-1}}(f_{q_{j}(y)})$ et $ad_{g_{y}^{-1}}(f'_{q_{j}(y)})$ par $res_{M_{j,\eta[y]}}$. Si $q_{j}(y)\not\in \dot{{\cal Y}}^0(\eta)$, on a $f_{q_{j}(y)}=f'_{q_{j}(y)}$, d'o\`u l'\'egalit\'e cherch\'ee. Supposons $q_{j}(y)\in \dot{{\cal Y}}^0(\eta)$. En vertu de la d\'efinition de $f'_{q_{j}(y)}$, il suffit de prouver que les images de $ad_{g_{y}^{-1}}(\varphi)$ et de $ad_{g_{y}^{-1}}(\varphi_{0})$ par $res_{M_{j,\eta[y]}}$ sont \'egales. Cela \'equivaut \`a $\varphi_{ad_{g_{y}}(M_{j,\eta[y]})}=\varphi_{0,ad_{g_{y}}(M_{j,\eta[y]})}$. Posons $z=q_{j}(y)$. D'apr\`es 5.13(3), on a $y=zg_{y}$. Donc $ad_{g_{y}}(M_{j,\eta[y]})=ad_{z^{-1}}ad_{y}(M_{j,\eta[y]})=ad_{z^{-1}}(M_{j,\eta})$, puisque $y\in M_{j}$. D'o\`u $ad_{g_{y}}(M_{j,\eta[y]})=ad_{\xi}(M_{j,\eta})$, o\`u $\xi$ est l'image de $z^{-1}$ dans $\Xi_{\eta}^{\Gamma_{F}}$. La d\'efinition de $\varphi_{0}$ entra\^{\i}ne  que $\varphi_{ad_{\xi}(M_{j,\eta})}=\varphi_{0,ad_{\xi}(M_{j,\eta})}$, ce qui prouve (5).

 D'apr\`es (3), (4) et (5), on a prouv\'e (1). Prouvons maintenant le lemme 5.13. Pour tout $j=1,...,k$, soit $\boldsymbol{\gamma}_{j}\in D_{g\acute{e}om}({\cal O}_{j})$. On suppose que $\sum_{j=1,...,k}\boldsymbol{\gamma}_{j}^{\tilde{G}}$ est stable. On s'est ramen\'e \`a trouver pour tout $j$ une distribution $\boldsymbol{\delta}_{j}\in D^{st}_{g\acute{e}om}({\cal O}_{j})$ de sorte que $\sum_{j=1,...,k}\boldsymbol{\gamma}_{j}^{\tilde{G}}=\sum_{j=1,...,k}\boldsymbol{\delta}_{j}^{\tilde{G}}$. L'\'el\'ement  $\oplus_{j=1,...,k} \boldsymbol{\gamma}_{j}$ est une forme lin\'eaire sur $\oplus_{j=1,...,k}I(\tilde{M}_{j}(F))_{{\cal O}_{j},loc}$. L'\'el\'ement $\oplus_{j=1,...,k}\boldsymbol{\delta}_{j}$ cherch\'e est une forme lin\'eaire sur $\oplus_{j=1,...,k}SI(\tilde{M}_{j}(F))_{{\cal O}_{j},loc}$. On peut la consid\'erer comme une forme lin\'eaire sur $\oplus_{j=1,...,k}I(\tilde{M}_{j}(F))_{{\cal O}_{j},loc}$ nulle sur le noyau de $\underline{s}_{loc}$, avec la notation du diagramme de 5.13. La condition d'\'egalit\'e des induites revient \`a ce que ces deux formes lin\'eaires co\"{\i}ncident sur l'image $Im$ de l'application $res_{loc}$. La condition n\'ecessaire et suffisante pour qu'il existe une solution est que $\oplus_{j=1,...,k} \boldsymbol{\gamma}_{j}$  annule $Im\cap Ker(\underline{s}_{loc})$. Un \'el\'ement de $Im\cap Ker(\underline{s}_{loc})$ est l'image d'un $f\in I(\tilde{U})$ tel que  l'image de $f$ dans $\oplus_{j=1,...,k}SI(\tilde{M}_{j}(F))_{{\cal O}_{j},loc}$ est nulle. D'apr\`es (1), on peut supposer que l'image de $f$ dans $SI(\tilde{U})$ est nulle. Par ailleurs,  la valeur de $\oplus_{j=1,...,k} \boldsymbol{\gamma}_{j}$ sur  l'image de $f$  est \'egale \`a celle de $\oplus_{j=1,...,k} \boldsymbol{\gamma}_{j}^{\tilde{G}}$ sur $f$. Celle-ci est nulle puisque cette distribution est stable. Cela ach\`eve la d\'emonstration. 
 
   \bigskip
  
  \subsection{Suite de la preuve, cas $F$ archim\'edien}
  Le probl\`eme pour $F={\mathbb C}$ se ram\`ene au m\^eme probl\`eme pour $F={\mathbb R}$ en rempla\c{c}ant chaque groupe et chaque espace par l'objet sur ${\mathbb R}$ obtenu par restriction des scalaires. On suppose donc $F={\mathbb R}$. Tous les ensembles du diagramme de 5.13 sont des espaces de Fr\'echet et toutes les applications sont continues. Les applications $s$, $\iota$, $\iota^{st}$ et $s_{loc}$ sont surjectives. Il en est de m\^eme de $\underline{s}$, $\underline{\iota}$, $\underline{\iota}^{st}$ et $\underline{s}_{loc}$. Montrons que
  
  (1) les images de $res$ et $res^{st}$ sont ferm\'ees.
  
    On a d\'ecrit $\oplus_{j=1,...,k}I(\tilde{U}_{j})$ comme un espace de familles  $(f_{j,y})_{j=1,...,k,y\in \dot{{\cal Y}}_{j}(\eta)}$ o\`u $f_{j,y}\in I(U_{j,y})$ pour tous $j,y$.   
   On va montrer que l'image de $res$ s'identifie au sous-espace  des  familles $(f_{j,y})_{j=1,...,k,y\in \dot{{\cal Y}}_{j}(\eta)} $ qui v\'erifient la condition suivante:

  (2) soient $j,j'\in \{1,...,k\}$, $y\in \dot{{\cal Y}}_{j}(\eta)$, $y'\in \dot{{\cal Y}}_{j'}(\eta)$, $R_{y}$ un Levi de $M_{j,\eta[y]}$, $R'_{y'}$ un Levi de $M_{j',\eta[y']}$ et $g\in G({\mathbb R})$ tel que $ad_{g}(\eta[y])=\eta[y']$ et $ad_{g}(R_{y})=R'_{y'}$; alors $f_{j',y',R'_{y'}}=ad_{g}(f_{j,y,R_{y}})$. 
  
  La condition est n\'ecessaire. En effet, soit $x\in R_{y}({\mathbb R})$ en position g\'en\'erale. Si notre collection $(f_{j,y})_{j=1,...,k,y\in \dot{{\cal Y}}_{j}(\eta)}$ provient de $f\in I(\tilde{U})$, on a
  $$I^{R'_{y'}}(ad_{g}(x),f_{j',y',R'_{y'}})=I^{M_{j',\eta[y']}}(ad_{g}(x),f_{j',y'})= I^{\tilde{M}_{j}'}(exp(ad_{g}(x))\eta[y'],f_{\tilde{M}_{j'}})$$
  $$= I^{\tilde{G}}(exp(ad_{g}(x))\eta[y'],f)=I^{\tilde{G}}(ad_{g}(exp(x)\eta[y]),f)=I^{\tilde{G}}(exp(x)\eta[y],f)$$
  $$=I^{\tilde{M}_{j}}(exp(x)\eta[y],f_{\tilde{M}_{j}})=I^{M_{j,\eta[y]}}(x,f_{j,y})=I^{R_{y}}(x,f_{j,y,R_{y}}).$$
Inversement, supposons (2) v\'erifi\'ee. Pour tout $y\in \dot{{\cal Y}}(\eta)$, consid\'erons l'ensemble des triplets $(j,y',g)$ tels que $j\in \{1,...,k\}$, $y'\in \dot{{\cal Y}}_{j}(\eta)$, $g\in G({\mathbb R})$ tels que $ad_{g}(\eta[y])=\eta[y']$. Le groupe $G_{\eta[y]}({\mathbb R})$ agit sur cet ensemble par multiplication de $g$ \`a droite. L'ensemble des orbites est fini. Fixons un ensemble de repr\'esentants $\dot{{\cal G}}_{y}$ de cet ensemble d'orbites. A tout \'el\'ement $\underline{g}=(j,y',g)\in \dot{{\cal G}}_{y}$ sont associ\'es un Levi $L_{\underline{g}}=ad_{g^{-1}}(M_{j,\eta[y']})$ de  $G_{\eta[y]}$ et une fonction $f_{\underline{g}}=ad_{g^{-1}}(f_{j,y'})\in I(L_{\underline{g}}({\mathbb R}))$. La condition (2) assure que ces familles de Levi et de fonctions v\'erifient la condition du lemme 4.3. On peut donc fixer une fonction $\phi_{y}\in I(G_{\eta[y]}({\mathbb R}))$ de sorte que $(\phi_{y})_{L_{\underline{g}}}=f_{\underline{g}}$ pour tout $\underline{g}\in \dot{{\cal G}}_{y}$.  Puisque chaque $f_{j,y'}$ est \`a support dans $U_{j,y'}$, il est plus ou moins clair que l'on peut fixer une fonction $\alpha$   sur $G_{\eta[y]}({\mathbb R})$,  qui est $C^{\infty}$ et invariante par conjugaison, dont le support est contenu dans $U_{y}$, de sorte que $f_{\underline{g}}=\alpha f_{\underline{g}}$ pour tout $\underline{g}\in \dot{{\cal G}}_{y}$. On peut aussi bien remplacer $\phi_{y}$ par $\alpha\phi_{y}$ et supposer $\phi_{y}\in I(U_{y})$. 
Consid\'erons l'ensemble des couples $(y',g)$ tels que $y'\in \dot{{\cal Y}}(\eta)$ et $g\in G({\mathbb R})$ tels que $ad_{g}(\eta[y])=\eta[y']$. De nouveau, le groupe $G_{\eta[y]}({\mathbb R})$ agit sur cet ensemble par multiplication de $g$ \`a droite. On fixe un ensemble $\dot{{\cal H}}_{y}$ de repr\'esentants de l'ensemble d'orbites. Il r\'esulte de (2) que, pour tout $(y',g)\in \dot{{\cal H}}_{y}$, la fonction $ad_{g^{-1}}(\phi_{y'})$ v\'erifie la m\^eme condition que $\phi_{y}$. On pose
$$f_{y}=\vert \dot{{\cal H}}_{y}\vert ^{-1}\sum_{(y',g)\in \dot{{\cal H}}_{y}}ad_{g^{-1}}(\phi_{y'}).$$
On voit  que la famille $(f_{y})_{y\in \dot{{\cal Y}}(\eta)}$  v\'erifie la condition (1) de 5.13. Elle s'identifie donc \`a un \'el\'ement de $I(\tilde{U})$. On voit que son image par $res$ est la famille 
 $(f_{j,y})_{j=1,...,k,y\in \dot{{\cal Y}}_{j}(\eta)}$ de d\'epart. Cela prouve (2).
 
 Cette relation (2) d\'ecrit l'image de $res$ par des conditions qui sont ferm\'ees. Il en r\'esulte que cette image est ferm\'ee. Une preuve similaire s'applique \`a l'application $res^{st}$. D'o\`u (1).

Montrons que (1) vaut aussi pour les applications localis\'ees, c'est-\`a-dire

(3) les images de $res_{loc}$ et $res_{loc}^{st}$ sont ferm\'ees.

L'espace $I(\tilde{G}({\mathbb R}))_{{\cal O}^{\tilde{G}},loc}$ s'identifie \`a celui des familles $(f_{y})_{y\in \dot{{\cal Y}}(\eta)}$ telles que

- on  a  $f_{y}\in I(G_{\eta[y]}({\mathbb R}))_{unip,loc}$ pour tout $y\in \dot{{\cal Y}}(\eta)$, o\`u l'indice $unip$ signifie le localis\'e relatif \`a la classe de conjugaison $\{1\}$ de $G_{\eta[y]}({\mathbb R})$;
 
-  soient $y,y'\in \dot{{\cal Y}}(\eta)$ et $g\in G({\mathbb R})$ tels que $ad_{g}(\eta[y))=\eta[y']$; alors $f_{y'}=ad_{g}(f_{y})$.
 
 On d\'ecrit de fa\c{c}on analogue l'espace $\oplus_{j=1,...,k}I(\tilde{M}_{j}({\mathbb R}))_{{\cal O}_{j},loc}$. Il est facile de reprendre la preuve de (1) et de montrer  que l'image de $res_{loc}$ est form\'e des familles $(f_{j,y})_{j=1,...,k,y\in \dot{{\cal Y}}_{j}(\eta)}\in \oplus_{j=1,...,k}I(\tilde{M}_{j}({\mathbb R}))_{{\cal O}_{j},loc}$ qui v\'erifient la condition (2) ci-dessus. On laisse cette preuve au lecteur. De nouveau, ces conditions sont ferm\'ees, ce qui prouve que l'image de $res_{loc}$ est ferm\'ee. Une preuve similaire s'applique \`a $res^{st}_{loc}$. D'o\`u (3).

De la commutativit\'e du diagramme de 5.13 r\'esulte que l'image de $Ker(s)$ par $res_{loc}\circ \iota$ est incluse dans $Im(res_{loc})\cap Ker(\underline{s}_{loc})$. 
On va prouver

(4) l'image de $Ker(s)$ par $res_{loc}\circ\iota$ est dense dans $Im(res_{loc})\cap Ker(\underline{s}_{loc})$.

On a $Im(res_{loc})=Im(res_{loc}\circ \iota)$. Soit $f\in I(\tilde{U})$, supposons $res_{loc}\circ \iota(f)\in Ker(\underline{s}_{loc})$. Soit $V_{1}$ un voisinage de $0$ dans $\oplus_{j=1,...,k}I(\tilde{M}_{j}(F))_{{\cal O}_{j},loc}$. Puisque $res_{loc}\circ \iota$ est continue, on peut fixer un voisinage $V_{2}$ de $0$ dans $I(\tilde{U})$ tel que $res_{loc}\circ \iota(V_{2})\subset V_{1}$. L'application $\underline{s}\circ res=res^{st}\circ s$ est d'image ferm\'ee d'apr\`es (1). Puisqu'il s'agit d'une application continue entre espaces de Fr\'echet, elle est ouverte sur son image. Il existe donc un voisinage $V_{3}$ de $0$ dans $\oplus_{j=1,...,k}SI(\tilde{U}_{j})$ tel que $V_{3}\cap Im(\underline{s}\circ res)\subset \underline{s}\circ res(V_{2})$. Fixons un tore maximal $T$ de $G_{\eta}$ et munissons $\mathfrak{t}({\mathbb C})$ d'une forme hermitienne d\'efinie positive invariante par le groupe de Weyl absolu de $T$ dans $G$. Si l'on suppose $\mathfrak{u}$ assez petit, tout \'el\'ement $\gamma\in \tilde{U}$ est conjugu\'e par un \'el\'ement de $G({\mathbb C})$ \`a un \'el\'ement $exp(X)\eta$ avec $X\in \mathfrak{t}({\mathbb C})$ proche de $0$. La norme $\vert X\vert $ est bien d\'etermin\'ee. Soit $b$ une fonction $C^{\infty}$ sur ${\mathbb R}$ qui vaut $1$ dans un voisinage de $0$ et est nulle sur $[1,+\infty[$. Pour tout entier $n\geq1$, on d\'efinit une fonction $B_{n}$ sur $\tilde{U}$ par $B_{n}(\gamma)=b(n\vert X\vert^2 )$ avec la notation pr\'ec\'edente. Elle vaut $1$ dans un voisinage de ${\cal O}^{\tilde{G}}$ et sa restriction aux \'el\'ements fortement r\'eguliers est invariante par conjugaison stable. On a

(5) $lim_{n\to \infty}\underline{s}\circ res(fB_{n})=0$.

En effet, fixons $j=1,...,k$ et un sous-tore maximal de $M_{j,\eta}$ d\'efini sur $F$. Pour simplifier la notation, on peut aussi bien supposer que c'est le tore $T$ pr\'ec\'edent. D\'efinissons des  fonctions $\psi$ et $\psi_{n}$  sur $\mathfrak{t}({\mathbb R})$ par $\psi(X)=S^{\tilde{G}}(exp(X)\eta,f)$ et $\psi_{n}(X)=S^{\tilde{G}}(exp(X)\eta,fB_{n})$.  
 Soit $D$ un op\'erateur diff\'erentiel sur $\mathfrak{t}$  \`a coefficients constants. On doit prouver que
$$lim_{n\to \infty}sup_{X\in \mathfrak{t}({\mathbb R})}\vert D\psi_{n}(X)\vert =0.$$
On a $\psi_{n}(X)=\psi(X)b(n\vert X\vert ^2)$. On voit que $D\psi_{n}(X)$ est combinaison lin\'eaire de termes $n^kD_{1}\psi(X)(D_{2}b)(n\vert X\vert^2 )P(X)$, avec des op\'erateurs diff\'erentiels $D_{1}$ et $D_{2}$ \`a coefficients constants et un polyn\^ome $P$, les coefficients de cette combinaison lin\'eaire ne d\'ependant pas de $n$ (les termes $n^k$ et $P(X)$ proviennent par d\'erivation de $n\vert X\vert ^2$). L'hypoth\`ese sur $f$ est que $\underline{\iota}^{st}\circ\underline{s}\circ res(f)=0$. Cela implique que toutes les d\'eriv\'ees de $\psi$ sont nulles en $0$. Le d\'eveloppement d'Euler-Mac-Laurin entra\^{\i}ne que l'on a pour tout $m\in {\mathbb N}$ une majoration $\vert D_{1} \psi(X)\vert \leq C_{m}\vert X\vert ^m$. Le terme ci-dessus est donc major\'e par 
$$C_{2k+2}n^k\vert X\vert ^{2k+2}\vert P(X)(D_{2}b)(n\vert X\vert ^2)\vert.$$
Le terme $(D_{2}b)(n\vert X\vert ^2)$ est major\'e uniform\'ement et sa non-nullit\'e  implique  $\vert X\vert^2 \leq n^{-1}$.  A fortiori, $\vert X\vert ^2\leq 1$ et $\vert P(X)\vert $ est uniform\'ement major\'e dans ce domaine. Le terme $n^kD_{1}\psi(X)(D_{2}b)(n\vert X\vert^2 )P(X)$ est donc major\'e par $Cn^{-1}$ pour une constante $C$ convenable. Cela prouve (5). 

Pour $n$ assez grand, on a donc $\underline{s}\circ res(fB_{n})\in V_{3}$. On peut alors choisir une fonction $f_{n}\in V_{2}$ de sorte que $\underline{s}\circ res(fB_{n}-f_{n})=0$. On peut alors reprendre la d\'emonstration du cas non-archim\'edien en l'appliquant \`a $fB_{n}-f_{n}$. On a l'analogue de 5.14(1), \`a savoir qu'il existe $f'\in I(\tilde{U})$ qui a m\^eme image que $fB_{n}-f_{n}$ dans $\oplus_{j=1,...,k}I(\tilde{M}_{j}({\mathbb R}))_{{\cal O}_{j},loc}$  et dont l'image dans $SI(\tilde{U})$ est nulle. Cette derni\`ere condition signifie que $f'$ appartient \`a $Ker(s)$. La premi\`ere condition signifie que $res_{loc}\circ \iota(f')=res_{loc}\circ \iota(fB_{n}-f_{n})$. Puisque $B_{n}$ vaut $1$ dans un voisinage de ${\cal O}^{\tilde{G}}$, on a $res_{loc}\circ \iota(fB_{n})=res_{loc}\circ\iota(f)$. On a aussi $res_{loc}\circ\iota(f_{n})\in V_{1}$. Cela prouve qu'il existe un \'el\'ement $f'\in Ker(s)$ tel que $res_{loc}\circ \iota(f-f')\in V_{1}$. D'o\`u la densit\'e affirm\'ee par (4).

Prouvons maintenant le lemme 5.13. Pour tout $j=1,...,k$, soit $\boldsymbol{\gamma}_{j}\in D_{g\acute{e}om}({\cal O}_{j})$. On suppose que $\sum_{j=1,...,k}\boldsymbol{\gamma}_{j}^{\tilde{G}}$ est stable. On s'est ramen\'e \`a trouver pour tout $j$ une distribution $\boldsymbol{\delta}_{j}\in D^{st}_{g\acute{e}om}({\cal O}_{j})$ de sorte que $\sum_{j=1,...,k}\boldsymbol{\gamma}_{j}^{\tilde{G}}=\sum_{j=1,...,k}\boldsymbol{\delta}_{j}^{\tilde{G}}$.  Posons pour simplifier $\boldsymbol{\gamma}= \oplus_{j=1,...,k} \boldsymbol{\gamma}_{j}$. C'est une forme lin\'eaire continue sur $\oplus_{j=1,...,k}I(\tilde{M}_{j}(F))_{{\cal O}_{j},loc}$.  Comme dans le cas non-archim\'edien, la stabilit\'e de $\sum_{j=1,...,k}\boldsymbol{\gamma}_{j}^{\tilde{G}}$  implique que $\boldsymbol{\gamma}$ est nulle sur l'image de $Ker(s)$ par $res_{loc}\circ \iota$. D'apr\`es (4) et puisque cette forme lin\'eaire est continue, elle est nulle sur $Im(res_{loc})\cap Ker(\underline{s}_{loc})$. Les espaces intervenant ici sont ferm\'es d'apr\`es (3). Donc $\boldsymbol{\gamma}$ se descend en une forme lin\'eaire continue sur $Im(res_{loc})/(Im(res_{loc})\cap Ker(\underline{s}_{loc}))$. L'application $\underline{s}_{loc}$ se quotiente en une bijection continue de cet espace sur $Im(res_{loc}^{st})$. D'apr\`es (3) et parce que nos espaces sont de Fr\'echet, cette bijection est un hom\'eomorphisme. On obtient qu'il existe une forme lin\'eaire continue $\boldsymbol{\delta}'$ sur  $Im(res_{loc}^{st})$ telle que $\boldsymbol{\delta}'\circ \underline{s}_{loc}$ co\"{\i}ncide avec $\boldsymbol{\gamma}$ sur $Im(res_{loc})$. Toujours d'apr\`es (3), on peut prolonger $\boldsymbol{\delta}'$ en une forme lin\'eaire continue $\boldsymbol{\delta}=\oplus_{j=1,...,k}\boldsymbol{\delta}_{j}\in \oplus_{j=1,...,k}D^{st}_{g\acute{e}om}({\cal O}_{j})$. La condition pr\'ec\'edente signifie que $\sum_{j=1,...,k}\boldsymbol{\gamma}_{j}^{\tilde{G}}=\sum_{j=1,...,k}\boldsymbol{\delta}_{j}^{\tilde{G}}$. 
 Cela ach\`eve la d\'emonstration. $\square$

 \bigskip

   \section{Le cas non ramifi\'e}
 
 \subsection{La situation non ramifi\'ee}
  Les donn\'ees sont les m\^emes qu'en 1.5. On suppose
 
 (1) $F$ est local non archim\'edien;
 
 (2) $G$ est non ramifi\'e (quasi-d\'eploy\'e sur $F$ et d\'eploy\'e sur une extension non ramifi\'ee);
 
 (3) ${\bf a}$ est non ramifi\'e (si on note ${\mathbb F}_{q}$ le corps r\'esiduel de $F$ et $\Gamma_F^{nr}=Gal(\bar{{\mathbb F}}_{q}/{\mathbb F}_{q})$, ${\bf a}$ provient par inflation d'un \'el\'ement de $H^1(\Gamma_F^{nr},Z(\hat{G}))$);
 
 (4) $\tilde{G}(F)$ poss\`ede un sous-espace hypersp\'ecial.
 
 Expliquons cette derni\`ere condition. Soit $K\subset G(F)$ un sous-groupe compact hypersp\'ecial (il en existe d'apr\`es (2)). Le normalisateur $Norm_{\tilde{G}(F)}(K)=\{\gamma\in \tilde{G}(F); ad_{\gamma}(K)=K\}$ peut \^etre vide. Sinon, c'est un espace principal homog\`ene sous $Z(G;F)K$ et on appelle sous-espace hypersp\'ecial une classe $\tilde{K}=\gamma K=K\gamma$ pour un \'el\'ement $\gamma\in Norm_{\tilde{G}(F)}(K)$. On dit que $\tilde{G}(F)$ poss\`ede un sous-espace hypersp\'ecial s'il existe $K$ tel que $Norm_{\tilde{G}(F)}(K)$ ne soit pas vide.
 
 {\bf Remarque.} L'hypoth\`ese que $G$ est non ramifi\'e n'implique pas l'existence d'un sous-espace hypersp\'ecial. Par exemple, pour un entier $n\geq1$ et un \'el\'ement $d\in F^{\times}$, consid\'erons $G=SL(n)$ et $\tilde{G}=\{g\in GL(n); det(g)=d\}$. On v\'erifie que $\tilde{G}(F)$ poss\`ede un sous-espace hypersp\'ecial si et seulement si la valuation de $d$ est divisible par $n$.
 \bigskip

 On fixe un couple $(K,\tilde{K})$ comme ci-dessus.
 
 Dans certains cas (en particulier pour les applications globales), on peut imposer une hypoth\`ese suppl\'ementaire, \`a savoir
 
 (Hyp) la caract\'eristique r\'esiduelle $p$ de $F$ est grande, plus pr\'ecis\'ement $p>N(G)e_{F}+1$, o\`u $N(G)$ est  l'entier d\'ependant de $G$ d\'efini en [W1] 4.3 et $e_{F}$ est l'indice de ramification de $F/{\mathbb Q}_{p}$.
 
 Nous ne l'imposons pas ici.

 \bigskip
  \subsection{Donn\'ees endoscopiques non ramifi\'ees}
 On note $I_{F}\subset W_{F}$ le groupe d'inertie. Soit ${\bf G}'=(G',{\cal G}',\tilde{s})$ une donn\'ee endoscopique de $(G,\tilde{G},{\bf a})$. On dit qu'elle est non ramifi\'ee si $I_{F}\subset {\cal G}'$. Cela entra\^{\i}ne:
 
(1)  $G'$ est non ramifi\'e.
 
 Preuve. Pour $w\in I_{F}$, soit $g_{w}=(g(w),w)\in {\cal G}'$ qui agit par $w_{G'}$ sur $\hat{G}'$. Puisque $w\in {\cal G}'$, on a aussi $g(w)\in {\cal G}'$. Puisque ${\cal G}'\cap \hat{G}=\hat{G}'$, on a $g(w)\in \hat{G}'$. On a $w_{G'}=ad_{g(w)}\circ w_{G}=ad_{g(w)}$, car $w_{G}=1$ ($G$ est non ramifi\'e). Donc $w_{G'}$ est un automorphisme int\'erieur de $\hat{G}'$. Il conserve par d\'efinition une paire de Borel \'epingl\'ee, c'est donc l'identit\'e. $\square$
 
    On suppose d\'esormais ${\bf G}'$ non ramifi\'ee.
 
 Rappelons que, si ${\cal E}$ est une paire de Borel \'epingl\'ee de $G$ d\'efinie sur $F$, la th\'eorie de Bruhat-Tits lui associe un sch\'ema en groupes ${\cal K}$ d\'efini sur l'anneau des entiers $\mathfrak{o}$ de $F$, et ${\cal K}(\mathfrak{o})$ est un sous-groupe compact hypersp\'ecial de $G(F)$. R\'eciproquement, tout tel sous-groupe est construit ainsi. Fixons donc   une paire de Borel \'epingl\'ee de $G$ d\'efinie sur $F$ dont est issu le groupe $K$ d\'ej\`a fix\'e. On peut la noter  ${\cal E}^*=(B^*,T^*,(E^*_{\alpha})_{\alpha\in \Delta})$. Notons $F^{nr}$ l'extension non ramifi\'ee maximale de $F$ et $\mathfrak{o}^{nr}$ son anneau d'entiers. Montrons que
 
 (2) l'ensemble $Z(\tilde{G},{\cal E}^*)(F^{nr})\cap T^*(\mathfrak{o}^{nr})\tilde{K}$ n'est pas vide.
 
 Soit $\gamma\in \tilde{K}$. La paire $ad_{\gamma}({\cal E}^*)$ est une paire de Borel \'epingl\'ee d\'efinie sur $F$. Deux telles paires sont conjugu\'ees sous le groupe adjoint $G_{AD}(F)$. Soit donc $x\in G_{AD}(F)$ tel que $ad_{x}\circ ad_{\gamma}({\cal E}^*)={\cal E}^*$. L'automorphisme $ad_{x}\circ ad_{\gamma}$ est d\'efini sur $F$. Puisqu'il conserve ${\cal E}^*$, il conserve aussi le sous-groupe hypersp\'ecial associ\'e \`a ${\cal E}$: $ad_{x}\circ ad_{\gamma}(K)=K$. Puisque $ad_{\gamma}(K)=K$, on a donc $ad_{x}(K)=K$. Cela entra\^{\i}ne que $x$ appartient au sous-groupe hypersp\'ecial $K_{AD}$ de $G_{AD}(F)$ associ\'e \`a la paire de Borel \'epingl\'ee $(B^*_{ad},T^*_{ad},(E^*_{\alpha})_{\alpha\in \Delta})$ d\'eduite de ${\cal E}^*$. On montre que l'application produit 
 $$T^*_{ad}(\mathfrak{o})\times K\to K_{AD}$$
 et l'application naturelle
 $$T^*(\mathfrak{o}^{nr})\to T^*_{ad}(\mathfrak{o}^{nr})$$
 sont toutes deux surjectives. Donc $x$ est l'image dans $G_{AD}(F)$ d'un produit $tk$, avec $t\in T(\mathfrak{o}^{nr})$ et $k\in K$. Puisque $ad_{x}\circ ad_{\gamma}({\cal E}^*)={\cal E}^*$, on a $tk\gamma\in Z(\tilde{G},{\cal E}^*)$. On a aussi $tk\gamma\in T(\mathfrak{o}^{nr})\tilde{K}$. Cela prouve (2). $\square$
 
 Fixons un \'el\'ement $e\in {\cal Z}(\tilde{G})(F^{nr})$, image d'un \'el\'ement de l'intersection $Z(\tilde{G},{\cal E}^*)(F^{nr})\cap T(\mathfrak{o}^{nr})\tilde{K}$, soit $e'$ son image dans ${\cal Z}(\tilde{G}')(F^{nr})$. Fixons une paire de Borel \'epingl\'ee ${\cal E}'=(B',T',(E'_{\alpha})_{\alpha\in \Delta'})$ de $G'$ d\'efinie sur $F$. Soit $K'$ le sous-groupe compact hypersp\'ecial de $G'(F)$ qui s'en d\'eduit. Pour $\sigma\in \Gamma_F$, soit $z'(\sigma)\in Z(G')$ tel que $e'=z'(\sigma)\sigma(e')$. Par construction, le cocycle $z'$ est non ramifi\'e et prend ses valeurs dans $T'(\mathfrak{o}^{nr})$. Or ce groupe est cohomologiquement trivial (cela r\'esulte du th\'eor\`eme de Lang). On peut choisir $t'\in T'(\mathfrak{o}^{nr})$ tel que $z'(\sigma)=\sigma(t')t^{_{'}-1}$. Alors $t'e'\in \tilde{G}'(F)$ et il est clair que $t'e'\in Norm_{\tilde{G}'(F)}(K')$. L'ensemble $\tilde{K}'=K't'e'$ est un sous-espace hypersp\'ecial de $\tilde{G}'(F)$. On voit qu'il ne d\'epend pas des choix de $e$ et $t'$. La classe de conjugaison par $G'_{AD}(F)$ du couple $(K',\tilde{K}')$ ne d\'epend pas des choix des paires de Borel \'epingl\'ees. Elle d\'epend par contre du couple $(K,\tilde{K})$ que l'on a fix\'e. 
 
 Ainsi l'espace $\tilde{K}$ que l'on a fix\'e d\'etermine un espace analogue $\tilde{K}'$ pour $\tilde{G}'(F)$, \`a conjugaison pr\`es par $G'_{AD}(F)$. Dans les raisonnements par r\'ecurrence, et dans ce qui suit, $\tilde{G}'(F)$ sera suppos\'e muni d'un tel ensemble $\tilde{K}'$ issu de $\tilde{K}$.

 \ass{Lemme}{La donn\'ee ${\bf G}'$ est relevante.}
 
 Preuve. Notons $\theta^*$ l'automorphisme $ad_{e}$ pour tout \'el\'ement $e\in Z(\tilde{G},{\cal E}^*)$. Il est d\'efini sur $F$. Introduisons le groupe $G_{1}=G^{\theta^*,0}$. A ${\cal E}^*$ est associ\'e une paire de Borel \'epingl\'ee ${\cal E}_{1}=(B_{1},T_{1},(E_{\alpha_{1}})_{\alpha_{1}\in \Delta_{1}})$ de $G_{1}$. On a $B_{1}=B^*\cap G_{1}$, $T_{1}=T^*\cap G_{1}$, $\Delta_{1}$ est  l'image de $\Delta$ par restriction \`a $T_{1}$. Pour $\alpha_{1}\in \Delta_{1}$, $E_{\alpha_{1}}$ est la somme des $E_{\alpha}$ pour $\alpha\in \Delta$ de restriction $\alpha_{1}$ (ces $\alpha$ forment une seule orbite pour l'action du groupe engendr\'e par $\theta^*$). De la paire ${\cal E}_{1}$ est issu un sous-groupe compact hypersp\'ecial $K_{1}$ de $G_{1}(F)$. Il r\'esulte des constructions de Bruhat et Tits que $K_{1}\subset K$. Des paires ${\cal E}^*$ et ${\cal E}'$ est issu un homomorphisme $\xi_{T^*,T'}:T^*\to T'$. Il existe un cocycle $\omega_{G'}:\Gamma_{F}\to W^{\theta^*}$ tel que $\sigma(\xi_{T^*,T'})=\xi_{T^*,T'}\circ \omega_{G'}(\sigma)$. Il est \'evidemment non ramifi\'e. Choisissons un \'el\'ement de Frobenius $\phi\in \Gamma_{F}$. Introduisons la section de Springer $n_{1}:W^{\theta^*}\to G_{1}$. Comme on l'a dit, \`a $K_{1}$ est associ\'e un sch\'ema en groupes ${\cal K}_{1}$ sur $\mathfrak{o}$ tel que ${\cal K}_{1}\times_{\mathfrak{o}}F=G_{1}$ et ${\cal K}_{1}(\mathfrak{o})=K_{1}$. Il r\'esulte des constructions que $n_{1}$ prend ses valeurs dans ${\cal K}_{1}(\mathfrak{o}^{nr})$. Posons $x=n_{1}(\omega_{G'}(\phi))$. On v\'erifie que $x$ appartient \`a un sous-groupe fini de ${\cal K}_{1}(\mathfrak{o}^{nr})$ invariant par $\Gamma_{F}$ (le groupe engendr\'e par l'image de $n_{1}$ et les \'el\'ements d'ordre $2$ de $T_{1}(\mathfrak{o}^{nr})$ convient). Appliquant par exemple [W1] 4.2(1), on voit qu'il existe $k\in {\cal K}_{1}(\mathfrak{o}^{nr})$ tel que $x=k\phi(k)^{-1}$. Posons ${\cal E}=ad_{k^{-1}}({\cal E}^*)$, notons $(B,T)$ la paire de Borel sous-jacente \`a ${\cal E}$. L'homomorphisme $\xi_{T,T'}$ d\'eduit de cette paire et de ${\cal E}'$ est $\xi_{T^*,T'}\circ ad_{k}$. D'apr\`es les constructions, il est \'equivariant pour les actions galoisiennes. Fixons $e\in Z(\tilde{G},{\cal E}^*)(F^{nr})\cap T^*(\mathfrak{o}^{nr})\tilde{K}$. Pour $\sigma\in \Gamma_{F}$, soit $z(\sigma)\in Z(G)$ tel que $e=z(\sigma)\sigma(e)$. Alors $z$ est un cocycle non ramifi\'e \`a valeurs dans $Z(G)\cap T^*(\mathfrak{o}^{nr})$. Mais $Z(G)\cap T^*(\mathfrak{o}^{nr})=Z(G)\cap T(\mathfrak{o}^{nr})$. Le groupe $T(\mathfrak{o}^{nr})$ \'etant cohomologiquement trivial, on peut choisir $\tau\in T(\mathfrak{o}^{nr})$ tel que $z(\sigma)=\sigma(\tau)\tau^{-1}$. Alors $\tau e\in \tilde{G}(F)$. Puisque $T(\mathfrak{o}^{nr})$ et $T^*(\mathfrak{o}_{nr})$ sont tous deux inclus dans ${\cal K}(\mathfrak{o}^{nr})$ (o\`u ${\cal K}$ est le sch\'ema en groupes associ\'e \`a $K$), on a m\^eme $\tau e\in \tilde{K}$. Puisque $k\in G_{1}$, la paire ${\cal E}$ est fix\'ee par $\theta^*=ad_{e}$. Il en r\'esulte que $ad_{\tau e}$ conserve $(B,T)$. Soit maintenant $t\in T(\mathfrak{o})$, posons $\gamma=t\tau e$. Notons $e'$ l'image de $e\in Z(\tilde{G},{\cal E})$ dans ${\cal Z}(\tilde{G}')$, posons $t'=\xi_{T,T'}(t\tau)$ et $\delta=t'e'$. Il est clair que $\delta\in \tilde{G}'(F)$ et que $(\delta,B',T',B,T,\gamma)$ est un diagramme. Si $t$ est en position g\'en\'erale, $\gamma$ est fortement r\'egulier, donc $(\delta,\gamma)\in {\cal D}({\bf G}')$. $\square$

   \bigskip
 
 \subsection{Facteur de transfert}
 Soit ${\bf G}'=(G',{\cal G}',\tilde{s})$ une donn\'ee endoscopique non ramifi\'ee de $(G,\tilde{G},{\bf a})$. Consid\'erons des donn\'ees auxiliaires $G'_{1}$, $\tilde{G}'_{1}$, $C_{1}$, $\hat{\xi}_{1}$. On dit qu'elles sont non ramifi\'ees si $G'_{1}$ est non ramifi\'e et le plongement $\hat{\xi}_{1}:{\cal G}'\to {^LG}'_{1}$ est l'identit\'e sur $I_{F}$. De telles donn\'ees existent. En fait
 
 (1) on peut choisir $G'_{1}=G'$.
 
 Preuve. On normalise l'action galoisienne sur $\hat{G}$ et $\hat{G}'$ en fixant des paires de Borel \'epingl\'ees de ces groupes et en imposant que les actions conservent ces paires. Choisissons un Frobenius $\phi\in W_{F}$ et un \'el\'ement $g_{\phi}\in {\cal G}'$ agissant comme $\phi_{G'}$ sur $\hat{G}'$. Alors ${\cal G}'$ est le produit semi-direct $(\hat{G}'\times I_{F})\rtimes g_{\phi}^{{\mathbb Z}}$. On d\'efinit une application $\hat{\xi}_{1}:{\cal G}'\to {^LG}'$ par $\hat{\xi}_{1}((x,w)g_{\phi}^n)=(x,w\phi^n)$ pour $x\in \hat{G}'$, $w\in I_{F}$, $n\in{\mathbb Z}$. C'est un isomorphisme. $\square$

   Supposons   les donn\'ees auxiliaires non ramifi\'ees. De $K'$ se d\'eduit un sous-groupe compact hypersp\'ecial $K'_{1}$ de $G'_{1}(F)$. Choisissons un \'el\'ement $\delta_{1,0}\in \tilde{G}'_{1}(F)$ dont l'image $\delta_{0}$ dans $\tilde{G}'(F)$ appartient \`a $ \tilde{K}'$. Alors $\tilde{K}'_{1}=K'_{1}\delta_{1,0}$ est un sous-espace hypersp\'ecial de $\tilde{G}'_{1}(F)$. Ce sous-espace \'etant fix\'e, nous allons d\'efinir un facteur de transfert $\Delta_{1}$ sur ${\cal D}_{1}$. 
   
     On fixe $g_{\phi}=(g(\phi),\phi)\in {\cal G}'$ comme dans la preuve de (1) ci-dessus et un \'el\'ement $g_{sc}(\phi)\in \hat{G}_{SC}$ dont l'image dans $\hat{G}_{AD}$ est la m\^eme que celle de $g(\phi)$. Il existe un unique cocycle $w\mapsto g(w)$ de $W_{F}$ dans $\hat{G}$ qui est non ramifi\'e et tel que $g(\phi)$ soit l'\'el\'ement que l'on vient de fixer. De m\^eme, il existe un unique cocycle $w\mapsto g_{sc}(w)$ de $W_{F}$ dans $\hat{G}_{SC}$ qui est non ramifi\'e et tel que $g_{sc}(\phi)$ soit l'\'el\'ement que l'on vient de fixer. Soit $w\mapsto z(w)$ le cocycle de $W_{F}$ dans $Z(\hat{G})$ tel que $g(w)=z(w)\pi(g_{sc}(w))$. On a \'evidemment $(g(w),w)\in {\cal G}'$ pour tout $w\in W_{F}$ et on pose $\hat{\xi}_{1}(g(w),w)=(\zeta_{1}(w),w)$. L'application $\zeta_{1}$ est un cocycle de $W_{F}$ dans $Z(\hat{G}'_{1})$. Les cocycles $z$ et $\zeta_{1}$  d\'eterminent des caract\`eres $\lambda_{z}$ de $G(F)$ et $\lambda_{\zeta_{1}}$ de $G'_{1}(F)$. Parce que les cocycles sont non ramifi\'es, $\lambda_{z}$ est trivial sur $K$ et $\lambda_{\zeta_{1}}$ est trivial sur $K'_{1}$ ([W1] 4.1(1)). Il existe donc une unique application $\tilde{\lambda}_{z}:\tilde{G}(F)\to {\mathbb C}^{\times}$ qui vaut $1$ sur $\tilde{K}$ et v\'erifie $\tilde{\lambda}_{z}(g\gamma)=\lambda_{z}(g)\tilde{\lambda}_{z}(\gamma)$ pour tous $g\in G(F)$ et $\gamma\in \tilde{G}(F)$. De m\^eme, il existe une unique application  $\tilde{\lambda}_{\zeta_{1}}:\tilde{G}'_{1}(F)\to {\mathbb C}^{\times}$ qui vaut $1$ sur $\tilde{K}'_{1}$ et v\'erifie $\tilde{\lambda}_{\zeta_{1}}(g_{1}\delta_{1})=\lambda_{\zeta_{1}}(g_{1})\tilde{\lambda}_{\zeta_{1}}(\delta_{1})$ pour tous $g_{1}\in G'_{1}(F)$ et $\delta_{1}\in \tilde{G}'_{1}(F)$.
   
  On fixe comme en 6.2 une  paire de Borel \'epingl\'ee ${\cal E}^*$ de $G$, d\'efinie sur $F$,   dont le groupe $K$  est issu. Soit $(\delta_{1},\gamma)\in {\cal D}_{1}$. On fixe un diagramme $(\delta,B',T',B,T,\gamma)$ et on utilise les constructions de 2.2.  En particulier, on compl\`ete $(B,T)$ en une paire de Borel \'epingl\'ee ${\cal E}$. On fixe $g\in G_{SC}$ tel que $ad_{g}({\cal E})={\cal E}^*$. On choisit pour cocha\^{\i}ne $u_{{\cal E}}$ l'application $u_{{\cal E}}(\sigma)=g^{-1}\sigma(g)$. On fixe $e\in Z(\tilde{G},{\cal E})$. Comme en 2.2, on d\'efinit une cocha\^{\i}ne $V_{T}:\Gamma_{F}\to T_{sc}$ par
   $$V_{T}(\sigma)=r_{T}(\sigma)n_{{\cal E}}(\omega_{T}(\sigma))u_{{\cal E}}(\sigma).$$
   La cocha\^{\i}ne $V_{T}$ est un cocycle. On \'ecrit $\gamma=\nu e$, avec   $\nu\in T$. On note $\nu_{ad}$ l'image de $\nu$ dans $T_{ad}$. Alors le couple $(V_{T},\nu_{ad})$ appartient \`a $Z^{1,0}(\Gamma_{F};T_{sc}\stackrel{1-\theta}{\to}T_{ad})$. On d\'efinit une cochaine $t_{T,sc}:W_{F}\to \hat{T}_{sc}$ par la m\^eme formule qu'en 2.2:
   $$t_{T,sc}(w)=\hat{r}_{T}(w)\hat{n}(\omega_{T}(w))g_{sc}(w)^{-1}\hat{n}_{G'}(\omega_{T,G'}(w))^{-1}\hat{r}_{T,G'}(w)^{-1}.$$
   C'est un cocycle. On note $s_{ad}$ l'image de $s$ dans $\hat{T}_{ad}$ (rappelons que $\tilde{s}=s\hat{\theta}$). Le couple $(t_{T,sc},s_{ad})$ appartient \`a $Z^{1,0}(W_{F};\hat{T}_{sc}\stackrel{1-\hat{\theta}}{\to}\hat{T}_{ad})$. On dispose du produit
   $$<.,.>:H^{1,0}(\Gamma_{F};T_{sc}\stackrel{1-\theta}{\to}T_{ad})\times H^{1,0}(W_{F};\hat{T}_{sc}\stackrel{1-\hat{\theta}}{\to}\hat{T}_{ad})\to {\mathbb C}^{\times}.$$
   On pose
   $$\Delta_{imp}(\delta_{1},\gamma)=\tilde{\lambda}_{\zeta_{1}}(\delta_{1})^{-1}\tilde{\lambda}_{z}(\gamma)<(V_{T},\nu_{ad}),(t_{T,sc},s_{ad})>^{-1}$$
   et
   $$\Delta_{1}(\delta_{1},\gamma)=\Delta_{II}(\delta,\gamma)\Delta_{imp}(\delta_{1},\gamma).$$
   
\ass{Lemme}{(i) Le facteur $\Delta_{1}$ ne d\'epend que des choix des sous-espaces hypersp\'eciaux $\tilde{K}$ et $\tilde{K}'_{1}$, c'est-\`a-dire qu'il ne d\'epend d'aucune autre donn\'ee auxiliaire.

(ii) Pour $(\delta_{1},\gamma), (\underline{\delta}_{1},\underline{\gamma})\in {\cal D}_{1}$, on a l'\'egalit\'e
$$\boldsymbol{\Delta}_{1}(\delta_{1},\gamma;   \underline{\delta}_{1},\underline{\gamma})=\Delta_{1}(\delta_{1},\gamma)\Delta_{1}(\underline{\delta}_{1},\underline{\gamma})^{-1}.$$}
   
  Preuve. On commence par d\'emontrer (ii), sous la r\'eserve que les choix de donn\'ees auxiliaires pour les deux paires $(\delta_{1},\gamma)$ et $ (\underline{\delta}_{1},\underline{\gamma})$ soient coh\'erents. Dans les constructions de 2.2 intervient un \'el\'ement $r\in G_{SC}$ tel que $ad_{r}({\cal E})=\underline{{\cal E}}$. Puisqu'on a choisi $g\in G_{SC}$ tel que $ad_{g}({\cal E})={\cal E}^*$ et de m\^eme $\underline{g}\in G_{SC}$ tel que $ad_{\underline{g}}(\underline{{\cal E}})={\cal E}^*$, on peut choisir et on choisit $r=\underline{g}^{-1}g$. Il est clair que le cocycle $V$ de 2.2 est l'image de $(V_{T},V_{\underline{T}}^{-1})$ par l'homomorphisme naturel $T_{sc}\times \underline{T}_{sc}\to U$. En utilisant la compatibilit\'e des produits aux deux diagrammes duaux
  $$\begin{array}{ccc}T_{sc}\times \underline{T}_{sc}&\stackrel{1-\theta}{\to}& S_{1}\\ \downarrow&&\downarrow\\ U&\stackrel{1-\theta}{\to}&S_{1}\\ \end{array}$$
et
   $$\begin{array}{ccc} \hat{S}_{1}&\stackrel{1-\hat{\theta}}{\to}& \hat{T}_{ad}\times \hat{\underline{T}}_{ad}\\ \uparrow&&\uparrow\\ \hat{S}_{1}&\stackrel{1-\hat{\theta}}{\to}&\hat{U},\\ \end{array}$$
on voit que 
$$\Delta_{imp}(\delta_{1},\gamma;   \underline{\delta}_{1},\underline{\gamma})=<((V_{T},V_{\underline{T}}^{-1}),\boldsymbol{\nu}_{1}),(\hat{V}_{1},(s_{ad},s_{ad}))>^{-1},$$
le produit \'etant celui sur
$$H^{1,0}(\Gamma_{F};T_{sc}\times\underline{T}_{sc}\stackrel{1-\theta}{\to}S_{1})\times H^{1,0}(W_{F};\hat{S}_{1}\stackrel{1-\hat{\theta}}{\to}\hat{T}_{ad}\times \hat{T}_{ad}).$$
Le cocycle $\hat{V}_{1}$ est le produit des deux cocycles suivants:

- l'image $\hat{V}_{sc}$ de $(t_{T,sc},t_{\underline{T},sc})$ par l'homomorphisme naturel  $\hat{q}:\hat{T}_{sc}\times \hat{\underline{T}}_{sc}\to \hat{S}_{1}$ qui, \`a $(t_{sc},\underline{t}_{sc})$, associe $\hat{q}(t_{sc},\underline{t}_{sc})=(j(t_{sc}),j(\underline{t}_{sc}), t_{sc}\underline{t}_{sc}^{-1})$;

- le cocycle $w\mapsto Z(w)=((\zeta_{1}(w),z(w)^{-1}),(\zeta_{1}(w),z(w)^{-1}),1)\in \hat{S}_{1}$.

Et le cocycle $(\hat{V}_{1},(s_{ad},s_{ad}))$ est le produit des deux cocycles $(\hat{V}_{sc},(s_{ad},s_{ad}))$ et du cocycle $(Z,1)$. On en d\'eduit l'\'egalit\'e
$$(2) \qquad \boldsymbol{\Delta}_{imp}(\delta_{1},\gamma;   \underline{\delta}_{1},\underline{\gamma})=<((V_{T},V_{\underline{T}}^{-1}),\boldsymbol{\nu}_{1}),(\hat{V}_{sc},(s_{ad},s_{ad}))>^{-1}<(V_{T},V_{\underline{T}}^{-1}),\boldsymbol{\nu}_{1}),(Z,1)>^{-1}.$$
En utilisant de nouveau une compatibilit\'e des produits, le premier terme est \'egal \`a
$$(3) \qquad <((V_{T},V_{\underline{T}}^{-1}),q(\boldsymbol{\nu}_{1})),((t_{T,sc},t_{\underline{T},sc}),(s_{ad},s_{ad}))>^{-1},$$
o\`u $q:S_{1}\to T_{ad}\times \underline{T}_{ad}$ est dual de l'homomorphisme $\hat{q}$ d\'efini ci-dessus. On voit que $q(\boldsymbol{\nu}_{1})=(\nu_{ad},\underline{\nu}_{ad})$. Le produit ci-dessus est maintenant celui sur
$$H^{1,0}(\Gamma_{F};T_{sc}\times \underline{T}_{sc}\stackrel{1-\theta}{\to}T_{ad}\times \underline{T}_{ad})\times H^{1,0}(W_{F};\hat{T}_{sc}\times \hat{\underline{T}}_{sc}\stackrel{1-\hat{\theta}}{\to}\hat{T}_{ad}\times \hat{\underline{T}}_{ad}).$$
Ces espaces comme ce produit se scindent selon les termes provenant de $T$ et ceux provenant de $\underline{T}$. Le produit (3) est alors \'egal \`a
$$(4) \qquad <(V_{T},\nu_{ad}),(t_{T,sc},s_{ad})>^{-1}<(V_{\underline{T}},\underline{\nu}_{ad}),(t_{\underline{T},sc},s_{ad})>.$$
Introduisons le tore $\hat{R}$ form\'e des $(t,\underline{t},t_{sc})\in \hat{T}\times\hat{\underline{T}}\times \hat{T}_{sc}$ tels que $j(t_{sc})=t\underline{t}^{-1}$ et le tore $\hat{R}_{1}$ form\'e des  $(t,\underline{t},t_{sc})\in \hat{T}'_{1}\times\hat{\underline{T}}'_{1}\times \hat{T}_{sc}^{\hat{\theta}}$ tels que $j(t_{sc})=t\underline{t}^{-1}$. On a des diagrammes commutatifs
$$\begin{array}{ccc}\hat{R}&\stackrel{\hat{\pi}}{\to}&\hat{T}_{ad}\times \hat{\underline{T}}_{ad}\\ \hat{\rho} \downarrow \,\,&&1-\hat{\theta}\,\downarrow\,\,\\ \hat{S}_{1}&\stackrel{1-\hat{\theta}}{\to}&\hat{T}_{ad}\times \hat{\underline{T}}_{ad},\\ \end{array}$$
   $$\begin{array}{ccc}\hat{R}_{1}&\to& 1\\ \hat{\rho}_{1} \downarrow\,\,&&\downarrow\\ \hat{S}_{1}&\stackrel{1-\hat{\theta}}{\to}&\hat{T}_{ad}\times \hat{\underline{T}}_{ad},\\ \end{array}$$
   o\`u $\hat{\pi}$, $\hat{\rho}$ et $\hat{\rho}_{1}$ sont les homomorphismes naturels.
   On introduit aussi les tores duaux $R$ et $R_{1}$ et les homomorphismes $\pi:T_{sc}\times\underline{T}_{sc}\to R$, $\rho:S_{1}\to R$ et $\rho_{1}:S_{1}\to R_{1}$ duaux de $\hat{\pi}$, $\hat{\rho}$ et $\hat{\rho}_{1}$.
 Le cocycle $Z$ est le produit des images des deux cocycles suivants:
 
 - l'inverse du  cocycle $\boldsymbol{z}:w\mapsto (z(w),z(w),1)\in \hat{R}$;
 
 - le cocycle $\boldsymbol{\zeta}_{1}:w\mapsto (\zeta_{1}(w),\zeta_{1}(w),1)\in \hat{R}_{1}$.
 
 On utilise la compatibilit\'e des produits aux diagrammes ci-dessus et la relation [KS1] A.3.13 (o\`u le signe dispara\^{\i}t d'apr\`es [KS2] 4.3). On voit que le deuxi\`eme terme de (2) est \'egal \`a
 $$(5) \qquad <(((1-\theta)(V_{T}),(1-\theta)(V_{\underline{T}}^{-1})),\rho(\boldsymbol{\nu}_{1})),({\bf z},1)><\rho_{1}(\boldsymbol{\nu}_{1}),\boldsymbol{\zeta}_{1}>^{-1},$$
 le premier produit \'etant celui sur
 $$H^{1,0}(\Gamma_{F};T_{sc}\times\underline{T}_{sc}\stackrel{\pi}{\to}R)\times H^{1,0}(W_{F};\hat{R}\stackrel{\hat{\pi}}{\to}\hat{T}_{ad}\times\hat{\underline{T}}_{ad})$$
  et le second celui sur
 $$H^{1,0}(\Gamma_{F};R_{1})\times H^{1,0}(W_{F};\hat{R}_{1}).$$
 On a l'\'egalit\'e $R_{1}=(T'_{1}\times \underline{T}'_{1})/diag_{-}(Z(G'_{1},G))$, o\`u $Z(G'_{1},G)$ est le sous-groupe des \'el\'ements de $Z(G'_{1})$ dont l'image dans $Z(G')$ appartient \`a l'image naturelle de $Z(G)$ (ou encore, c'est la projection dans $Z(G'_{1})$ du groupe $\mathfrak{Z}_{1}$ de 2.2). Le tore $R_{1}$ est un sous-tore maximal du groupe $(G'_{1}\times G'_{1})/diag_{-}(Z(G'_{1},G))$.  L'\'el\'ement $\rho_{1}(\boldsymbol{\nu}_{1})$ est \'egal \`a l'image dans $R_{1}$ de $(\mu_{1},\underline{\mu}_{1}^{-1})$. Son image dans $(G'_{1}\times G'_{1})/diag_{-}(Z(G'_{1},G))$ est celle de $(x\underline{\mu}_{1},\underline{\mu}_{1}^{-1})$, o\`u $x\in G'_{1}(F)$ est l'\'el\'ement tel que $x\delta_{1}=\underline{\delta}_{1}$. Le calcul de la preuve du lemme 2.5 montre que le produit de cet \'el\'ement avec $\boldsymbol{\zeta}_{1}$ vaut $\lambda_{\zeta_{1}}(x)$.   En appliquant les d\'efinitions, on obtient
 $$(6) \qquad <\rho_{1}(\boldsymbol{\nu}_{1}),\boldsymbol{\zeta}_{1}>^{-1}=\tilde{\lambda}_{\zeta_{1}}(\delta_{1})^{-1}\tilde{\lambda}_{\zeta_{1}}(\underline{\delta}_{1}).$$
 On a l'\'egalit\'e $R=(T\times\underline{T})/diag_{-}(Z(G))$. C'est un sous-tore maximal du groupe $G^{\flat}=(G\times G)/diag_{-}(Z(G))$. On a l'\'egalit\'e $G^{\flat}_{SC}=G_{SC}\times G_{SC}$ et $T_{sc}\times \underline{T}_{sc}$ est l'image r\'eciproque de $R$ dans $G^{\flat}_{SC}$. On se retrouve dans la situation de 2.4. C'est-\`a-dire que ${\bf z}$ est un cocycle \`a valeurs dans $Z(\hat{G}^{\flat})$ qui d\'etermine un caract\`ere $\lambda_{{\bf z}}$ du groupe $G^{\flat}(F)$. Si $(((1-\theta)(V_{T}),(1-\theta)(V_{\underline{T}}^{-1})),\rho(\boldsymbol{\nu}_{1}))$ est l'image de $y^{\flat}\in G^{\flat}(F)$ par l'homomorphisme surjectif
 $$G^{\flat}(F)\to H^{1,0}(\Gamma_{F};T_{sc}\times\underline{T}_{sc}\stackrel{\pi}{\to}R),$$
 on a l'\'egalit\'e
$$(7) \qquad <(((1-\theta)(V_{T}^{-1}),(1-\theta)(V_{\underline{T}})),\rho(\boldsymbol{\nu}_{1})),({\bf z},1)>=\lambda_{{\bf z}}(y^{\flat}).$$
 Il reste \`a calculer un \'el\'ement $y^{\flat}$ v\'erifiant la propri\'et\'e ci-dessus. Introduisons l'\'el\'ement $e^*=geg^{-1}\in Z(\tilde{G},{\cal E}^*)$. Remarquons que, d'apr\`es nos choix, on a aussi $e^*=\underline{g}\underline{e}\underline{g}^{-1}$. Ecrivons $\gamma=ye^*$, $\underline{\gamma}=\underline{y}e^*$ avec $y, \underline{y}\in G$. Puisque ${\cal E}^*$ est d\'efini sur $F$, on a $\sigma(e^*)\in Z(G)e^*$ pour tout $\sigma\in \Gamma_{F}$. Il en r\'esulte que l'image de $(y,\underline{y}^{-1})$ dans $G^{\flat}$ appartient \`a $G^{\flat}(F)$. Montrons que
 
 (8) on peut choisir pour $y^{\flat}$ l'image de $(y,\underline{y}^{-1})$ dans $G^{\flat}(F)$.

   D\'ecomposons $\nu$ en $\pi(\nu_{sc})\nu_{Z}$, avec $\nu_{Z}\in Z(G)$ et $\nu_{sc}\in T_{sc}$. On a $\gamma=\nu e=\nu g^{-1}e^*g=\nu g^{-1}ad_{e^*}(g)e^*$. Donc $y=\pi(y_{sc})\nu_{Z}$, avec $y_{sc}=\nu_{sc}g^{-1}ad_{e^*}(g)$. On d\'efinit le cocycle $\tau:\Gamma_{F}\to Z(G_{SC})$ par $\tau(\sigma)=y_{sc}\sigma(y_{sc})^{-1}$. En utilisant des notations analogues pour l'\'el\'ement $\underline{\gamma}$, le calcul  de 2.4 montre que l'image de $(y,\underline{y}^{-1})$ dans  $H^{1,0}(\Gamma_{F};T_{sc}\times\underline{T}_{sc}\stackrel{\pi}{\to}R)$ est le cocycle $((\tau,\underline{\tau}^{-1}),(\nu_{Z},\underline{\nu}_{Z}^{-1}))$. On doit montrer que celui-ci est cohomologue \`a $(((1-\theta)(V_{T}),(1-\theta)(V_{\underline{T}}^{-1})),\rho(\boldsymbol{\nu}_{1}))$. Tout d'abord, on a l'\'egalit\'e $\rho(\boldsymbol{\nu}_{1})=(\nu,\underline{\nu}^{-1})$. Donc $(((1-\theta)(V_{T}),(1-\theta)(V_{\underline{T}}^{-1})),\rho(\boldsymbol{\nu}_{1}))$ est cohomologue \`a $((\tau',(\underline{\tau}')^{-1}),(\nu_{Z},\underline{\nu}_{Z}^{-1}))$, o\`u $\tau'(\sigma)=\nu_{sc}(1-\theta)(V_{T}(\sigma))\sigma(\nu_{sc})^{-1}$. Rappelons que le $\theta$ de cette relation est plus pr\'ecis\'ement $ad_{e}$, c'est-\`a-dire $ad_{g}^{-1}\circ ad_{e^*}\circ ad_{g}$. En reprenant la d\'efinition de $V_{T}$ et en se rappelant que les termes $r_{T}(\sigma)$ et $n_{{\cal E}}(\omega_{T}(\sigma))$ sont fixes par $ad_{e}$, on obtient
 $$\tau'(\sigma)=\nu_{sc}ad_{e}(u_{{\cal E}}(\sigma)^{-1})u_{{\cal E}}(\sigma)
\sigma(\nu_{sc})^{-1}$$
$$= \nu_{sc} \,ad_{g}^{-1}\circ ad_{e^*}\circ ad_{g}(\sigma(g)^{-1}g)g^{-1}\sigma(g)\sigma(\nu_{sc})^{-1}$$
$$=\nu_{sc}g^{-1}ad_{e^*}(g\sigma(g)^{-1})\sigma(g\nu_{sc}^{-1}).$$
L'automorphisme $ad_{e^*}$ est d\'efini sur $F$. D'o\`u
$$\tau'(\sigma)=\nu_{sc}g^{-1}ad_{e^*}(g)\sigma(ad_{e^*}(g^{-1})g\nu_{sc}^{-1})=y_{sc}\sigma(y_{sc})^{-1}=\tau(\sigma).$$
Un calcul analogue vaut pour $\underline{\tau}'$, ce qui d\'emontre (8).

 On peut  donc appliquer (7) en prenant pour $y^{\flat}$ l'image de $(y,\underline{y}^{-1})$. Un calcul analogue \`a celui de la preuve du lemme 2.5 montre que $\lambda_{{\bf z}}(y^{\flat})=\lambda_{z}(x)$, o\`u $x$ est l'\'el\'ement de $G(F)$ tel que $y=x\underline{y}$ ou encore $\gamma=x\underline{\gamma}$. D'o\`u
$$(9) <(((1-\theta)(V_{T}),(1-\theta)(V_{\underline{T}}^{-1})),\rho(\boldsymbol{\nu}_{1})),({\bf z},1)>=\tilde{\lambda}_{z}(\gamma)\tilde{\lambda}_{z}(\underline{\gamma})^{-1}.$$
   
Rassemblons nos calculs. Le facteur $\boldsymbol{\Delta}_{imp}(\delta_{1},\gamma;\underline{\delta}_{1},\underline{\gamma})$ est le produit des termes (4), (6) et (9). Autrement dit
$$\boldsymbol{\Delta}_{imp}(\delta_{1},\gamma;\underline{\delta}_{1},\underline{\gamma})=\Delta_{imp}(\delta_{1},\gamma)\Delta_{imp}(\underline{\delta}_{1},\underline{\gamma})^{-1}.$$
Cela d\'emontre le (ii) de l'\'enonc\'e.
 
 Prouvons maintenant l'assertion (i). Les donn\'ees auxiliaires pour une paire $(\delta_{1},\gamma)$ sont
 
 (10) le diagramme $(\delta,B',T',B,T,\gamma)$, la paire de Borel \'epingl\'ee ${\cal E}$, l'\'el\'ement $g\in G_{SC}$, les $a$-data et les $\chi$-data;
 
 (11) la paire de Borel \'epingl\'ee ${\cal E}^*$, les paires de Borel \'epingl\'ees des groupes duaux, les termes $g(\phi)$ et $g_{sc}(\phi)$;
 
 (12) l'\'el\'ement $e\in Z(\tilde{G},{\cal E})$.
 
 On voit tout de suite que le choix de $e$ n'influe pas: ce terme ne sert qu'\`a d\'efinir $ad_{e}$ et $\nu$. L'automorphisme $ad_{e}$ ne d\'epend pas du choix de $e$. Le terme $\nu$ en d\'epend, mais il n'intervient que via $\nu_{ad}$ qui, lui, n'en d\'epend pas. Quand on consid\`ere deux couples $(\delta_{1},\gamma)$ et $(\underline{\delta}_{1},\underline{\gamma})$, faire des choix coh\'erents signifie que l'on prend les m\^emes objets (11) pour les deux couples (il y a aussi une condition portant sur les termes $e$ et $\underline{e}$, mais on peut l'oublier d'apr\`es ce que l'on vient de dire). Il n'y a aucune condition de coh\'erence portant sur les objets (10). Puisque $\boldsymbol{\Delta}_{1}(\delta_{1},\gamma;\underline{\delta}_{1},\underline{\gamma})$ ne d\'epend d'aucun choix et puisque $\Delta_{1}(\underline{\delta}_{1},\underline{\gamma})$ ne d\'epend pas des objets (10) relatifs au couple $(\delta_{1},\gamma)$, on  d\'eduit de notre preuve (partielle) de (ii) que $\Delta_{1}(\delta_{1},\gamma)$ ne d\'epend pas des objets (10) et qu'il ne d\'epend des objets (11) que par multiplication par un scalaire. Il nous suffit donc de prouver que pour un couple particulier $(\delta_{1},\gamma)$, le facteur $\Delta_{1}(\delta_{1},\gamma)$ ne d\'epend pas des objets (11). On choisit l'une des paires $(\delta,\gamma)$ que l'on a construites dans la preuve du lemme 6.2. L'\'el\'ement $\delta$ appartient \`a l'espace $\tilde{K}'$.  On v\'erifie facilement que l'application $\tilde{K}_{1}'\to \tilde{K}'$ est surjective. On rel\`eve $\delta$ en un \'el\'ement $\delta_{1}\in \tilde{K}'_{1}$. On choisit pour diagramme et pour \'el\'ement $g$ le diagramme et l'\'el\'ement $k$ que l'on a construits dans cette preuve. Les tores $T$ et $T'$ sont non ramifi\'es.  On peut supposer que $\chi_{\alpha}$ est trivial pour un \'el\'ement $\alpha\in \Sigma(T)_{res,ind}$ appartenant \`a une orbite asym\'etrique et est non ramifi\'e pour un $\alpha$ appartenant \`a une orbite sym\'etrique. Cette derni\`ere condition d\'etermine $\chi_{\alpha}$: on a $\chi_{\alpha}(x)=(-1)^{val_{F_{\alpha}}(x)}$ pour $x\in F_{\alpha}$, o\`u $val_{F_{\alpha}}$ est la valuation usuelle de $F_{\alpha}$. On peut aussi supposer que les $a$-data $a_{\alpha}$ sont des unit\'es de $F_{\alpha}$. Il r\'esulte alors des constructions que $(V_{T},\nu_{ad})$ appartient \`a 
 $$H^{1,0}(\Gamma_{F}/\Gamma_{F^{nr}};T_{sc}(\mathfrak{o}^{nr})\stackrel{1-\theta}{\to}T_{ad}(\mathfrak{o}^{nr})).$$
 Par ailleurs, $(t_{T,sc},s_{ad})$ appartient \`a
 $$H^{1,0}(W_{F}/W_{F^{nr}};\hat{T}_{sc}\stackrel{1-\hat{\theta}}{\to}\hat{T}_{ad}).$$
 Or la restriction de la dualit\'e de Kottwitz-Shelstad  au produit des deux groupes ci-dessus est  triviale. Donc 
 $$<(V_{T},\nu_{ad}),(t_{T,sc},s_{ad})>=1.$$
 Puisque $\delta_{1}\in \tilde{K}'_{1}$ et $\gamma\in \tilde{K}$, on a $\tilde{\lambda}_{\zeta_{1}}(\delta_{1})=\tilde{\lambda}_{z}(\gamma)=1$. D'o\`u $\Delta_{imp}(\delta_{1},\gamma)=1$ et $\Delta_{1}(\delta_{1},\gamma)=\Delta_{II}(\delta,\gamma)$. Ce terme ne d\'ependant pas des donn\'ees (11), cela ach\`eve la d\'emonstration. $\square$

 Dans [W1], on a donn\'e une autre fa\c{c}on de normaliser le facteur de transfert, sous  l'hypoth\`ese (Hyp) de 6.1. On a

(13) sous  l'hypoth\`ese (Hyp), le facteur de [W1] co\"{\i}ncide avec celui ci-dessus.

Le facteur de [W1] est caract\'eris\'e par le fait que, pous  $(\delta_{1},\gamma)$ appartenant \`a un certain sous-ensemble ${\cal D}_{1,nr}\subset {\cal D}_{1}$, on a $\Delta_{1}(\delta_{1},\gamma)=\Delta_{II}(\delta,\gamma)$. Or, parmi les couples que  l'on a consid\'er\'e \`a la fin de la d\'emonstration ci-dessus, il y en a qui appartiennent \`a ${\cal D}_{1,nr}$. On a prouv\'e que notre pr\'esent facteur v\'erifiait l'\'egalit\'e ci-dessus pour ces couples-l\`a. Cela conclut. $\square$

   \bigskip
   
   \subsection{Le lemme fondamental}
   On suppose ${\bf G}'$ non ramifi\'e. Consid\'erons des donn\'ees auxiliaires $G'_{1}$,..., $\hat{\xi}_{1}$ non ramifi\'ees. On fixe comme dans le paragraphe pr\'ec\'edent un sous-espace hypersp\'ecial $\tilde{K}'_{1}\subset \tilde{G}'_{1}(F)$.       Notons ${\bf 1}_{\tilde{K}}$ la fonction caract\'eristique de $\tilde{K}$ et ${\bf 1}_{\tilde{K}'_{1},\lambda_{1}}$ l'\'el\'ement de $C_{c,\lambda_{1}}^{\infty}(\tilde{G}'_{1}(F))$ qui est \`a support dans $C_{1}(F)\tilde{K}'_{1}$ et vaut $1$ sur $\tilde{K}'_{1}$. On utilise le facteur de transfert normalis\'e de 6.3 pour d\'efinir la notion de transfert. Gr\^ace \`a Ngo Bao Chau, on a:
 
 \ass{Th\'eor\`eme (lemme fondamental pour les unit\'es)}{${\bf 1}_{\tilde{K}'_{1},\lambda_{1}}$ est un transfert de ${\bf 1}_{\tilde{K}}$.}
 
 Notons ${\cal H}$, resp. ${\cal H}'_{1}$, l'alg\`ebre des fonctions sur $G(F)$, resp. $G'_{1}(F)$, \`a support compact et biinvariantes par $K$, resp. $K'_{1}$. Notons $\phi\in  W_{F}$ un \'el\'ement de Frobenius et $\hat{{\cal H}}$, resp. $\hat{{\cal H}}'$, resp. $\hat{{\cal H}}'_{1}$, l'alg\`ebre des fonctions polynomiales sur $\hat{G}\rtimes \phi\subset {^LG}$, resp. ${\cal G}\cap (\hat{G}\rtimes \phi)$, resp. $\hat{G}'_{1}\rtimes \phi\subset {^LG}'_{1}$, invariantes par conjugaison par $\hat{G}$, resp. $\hat{G'}$, resp. $\hat{G}'_{1}$. On a un diagramme
 $$\begin{array}{cccc}{\cal H}&\stackrel{{\rm Satake}}{\simeq}&\hat{{\cal H}}&\\ &&\downarrow&{\rm restriction}\\ && \hat{{\cal H}}'&\\ &&\uparrow&{\rm restriction}\\ {\cal H}'_{1}&\stackrel{{\rm Satake}}{\simeq}&\hat{{\cal H}}'_{1}&\\ \end{array}$$
 D'autre part, ${\cal H}$ agit par convolution \`a droite et \`a gauche sur $C_{c}^{\infty}(\tilde{G}(F))$ et ${\cal H}'_{1}$ agit par convolution \`a droite et \`a gauche sur $C_{c,\lambda_{1}}^{\infty}(\tilde{G}'_{1}(F))$. On peut peut-\^etre \'enoncer un lemme fondamental sous la forme suivante.
 
 \ass{Conjecture}{Soient $h\in {\cal H}$ et $h'_{1}\in {\cal H}'_{1}$. On suppose que $h$ et $h'_{1}$ ont m\^eme image dans $\hat{{\cal H}}'$. Alors $h'_{1}*{\bf 1}_{\tilde{K}'_{1},\lambda_{1}}={\bf 1}_{\tilde{K}'_{1},\lambda_{1}}*h'_{1}$ est un transfert de $h*{\bf 1}_{\tilde{K}}$ comme de ${\bf 1}_{\tilde{K}}*(\omega^{-1}h)$.}

 Ces \'enonc\'es se traduisent ais\'ement selon le formalisme introduit en 2.5. A l'aide du facteur de transfert normalis\'e, on identifie $C_{c,\lambda_{1}}^{\infty}(\tilde{G}'_{1}(F))$ \`a $C_{c}^{\infty}({\bf G}')$. Notons ${\bf 1}_{\tilde{K}',{\bf G}'}$ l'image de ${\bf 1}_{\tilde{K}'_{1},\lambda_{1}}$ dans ce dernier espace. On v\'erifie qu'elle ne d\'epend pas des donn\'ees auxiliaires choisies. Le th\'eor\`eme signifie que cet \'el\'ement est un transfert de ${\bf 1}_{\tilde{K}}$. De m\^eme, on peut introduire une alg\`ebre ${\cal H}'$ limite inductive des alg\`ebres ${\cal H}'_{1}$ quand $G'_{1}$,..., $\Delta_{1}$ parcourent toutes les donn\'ees auxiliaires non ramifi\'ees. Elle s'identifie, mais de fa\c{c}on non canonique, \`a l'alg\`ebre des fonctions sur $G'(F)$
\`a support compact et biinvariantes par $K'$. L'isomorphisme de Satake identifie ${\cal H}'$ \`a $\hat{{\cal H}}'$. L'alg\`ebre ${\cal H}'$ agit sur $C_{c}^{\infty}({\bf G}')$ et la conjecture ci-dessus se
r\'ecrit imm\'ediatement en termes de cette action.

\bigskip

{\bf Bibliographie}

[A1] J. Arthur: {\it On the transfer of distributions: weighted orbital integrals}, Duke Math. J. 99 (1999), p. 209-283

[A2] ---------: {\it On local character relations}, Selecta math. 2 (1996), p. 501-579

[A3] ---------: { \it Germ expansions for real groups}

[Bor] A. Borel: {\it Automorphic $L$-functions}, in {\it Automorphic forms, representations and $L$-functions}, Proc. of Symposia in Pure Math XXXIII, part 2, A. Borel et W. Casselman ed., AMS 1979

[Boua] A. Bouaziz: {\it Int\'egrales orbitales sur les groupes de Lie r\'eductifs}, Ann. Sc. ENS 27 (1994), p. 573-609

[Bour] N. Bourbaki: {\it Groupes et alg\`ebres de Lie, chapitres 4,5 et 6}, Hermann 1968

[K1] R. Kottwitz: {\it Rational conjugacy classes in reductive groups}, Duke Math. J. 49 (1982), p. 785-806

[K2] -------------: { \it Stable trace formula: elliptic singular terms}, Math. Annalen 275 (1986), p. 365-399

[KS1] R. Kottwitz, D. Shelstad: {\it Foundations of twisted endoscopy}, Ast\'erisque 255 (1999)

[KS2] -------------------------------: {\it On splitting invariants and sign conventions in endoscopic transfer}, pr\'epublication 2012

[Lab1] J.-P. Labesse: {\it Cohomologie, stabilisation et changement de base}, Ast\'erisque 257 (1999)

[Lab2] ------------------: {\it Stable twisted trace formula: elliptic terms}, Journal of the Inst. of Math. Jussieu 3 (2004), p. 473-530

[LS] R. P. Langlands, D. Shelstad: {\it On the definition of transfer factors}, Math. Annalen 278 (1987), p. 219-271

[R1] D. Renard: {\it Int\'egrales orbitales tordues sur les groupes de Lie r\'eductifs r\'eels}, J. of Funct. Analysis 145 (1997), p. 374-454

[R2] -------------: {\it Endoscopy for real reductive groups}, in {\it On the stabilization of the trace formula}, L. Clozel, M. Harris, J.-P. Labesse, B.-C. Ng\^o \'ed., International Press, 2011

[S1] D. Shelstad: {\it On geometric transfer in real twisted endoscopy}, pr\'epublication 2011

[S2] ----------------: {\it  Characters and inner forms of a quasi-split group over ${\mathbb R}$}, Compositio Math. 39 (1979), p.11-45

[W1] J.-L. Waldspurger: {\it L'endoscopie tordue n'est pas si tordue}, Memoirs AMS 908 (2008)

[W2] ------------------------: {\it A propos du lemme fondamental pond\'er\'e tordu}, Math. Ann. 343 (2009), p. 103-174

[W3] -----------------------: {\it La formule des traces locale tordue}, pr\'epublication 2012

\bigskip

Institut de Math\'ematiques de Jussieu, CNRS

2 place Jussieu

75005 Paris

e-mail: waldspur@math.jussieu.fr

\end{document}